\documentclass[a4paper,11pt]{article}
\usepackage{amsmath, amsfonts, amssymb}
\usepackage[english, russian]{babel}
\usepackage {color}

\definecolor{mycolor1}{rgb}{0.020,0.388,0.757}
\definecolor{mycolor2}{rgb}{0.208,0.208,0.208}

\definecolor{mycolor1}{rgb}{0.020,0.388,0.757}
\definecolor{mycolor2}{rgb}{0.208,0.208,0.208}

\begin{document}
\begin{center}
\Large{PDE-Systems associated with the hypergeometric functions in three variables and their  particular solutions near the origin}
\end{center}

\bigskip

\begin{center}
M.Ruzhansky${}^{1,a}$, A.Hasanov${}^{1,2,b}$,  T.G.Ergashev${}^{1,3,c}$
\end{center}

\bigskip

${}^{\left.1\right)}$ Ghent University, Department of Mathematics, Analysis, Logic and Discrete Math., Belgium;

${}^{\left.2\right)}$ Romanovskii Institute of Mathematics, Uzbek Academy of Sciences, Tashkent, Uzbekistan;

$^{\left.3\right)}$ TIIAME National Research University, Tashkent, Uzbekistan;

\bigskip

E-mails: 

${}^{\left.a\right)}$michael.ruzhansky@ugent.be; 

 ${}^{\left.b\right)}$anvarhasanov@yahoo.com;

${}^{\left.c\right)}${ergashev.tukhtasin@gmail.com}

\bigskip

\textbf{Abstract.} The great success of the theory of hypergeometric series in one variable has stimulated the development of a corresponding theory in two and more variables. Horn has investigated the convergence of 34 (14 complete and 20 confluent) hypergeometric series of two variables and established the systems of partial differential equations which they satisfy.   At present, 600 (of which 205 are complete and 395 are confluent) hypergeometric functions of three second-order variables are known. The present work is devoted to the composition of systems of partial differential equations satisfied by 600 confluent hypergeometric functions of three variables. In addition, a particular solutions (if such solutions exist) of some systems of differential equations have been found near the origin.

\bigskip

\textit{Keywords:} Gauss hypergeometric function, Appell functions, Horn functions, Humbert functions, confluent functions, hypergeometric functions of three variables, systems of partial differential equations of hypergeometric type.

\bigskip

MSC 2020: 33C05, 33C15, 33C65, 33C70.

\bigskip

\section{Historical introduction}

A function
\begin{equation}
\label{hyper}
F(a,b;c;x)=\sum\limits_{n=0}^\infty\frac{(a)_n(b)_n}{(c)_n}\frac{x^n}{n!},\,\,|x|<1,
\end{equation}
is known as the Gaussian hypergeometric function, where
$(\nu)_n$  is a Pochhammer symbol defined by
$$(\nu)_0:=1, \,\,(\nu)_n:=\nu(\nu+1)...(\nu+n-1),\,\, n\in N: =\{1,2,3,...\}, $$ for which the following equality is true
\begin{equation}
\label{poch1}
(\nu)_n = \frac{\Gamma(\nu+n)}{\Gamma(\nu)}, \,\,\, n \in \{0\}\cup N;
\end{equation}
$\Gamma(\nu) $  is a famous Euler gamma-function.

In his work \textit{Arithmetica Infinitorum} [42], the Oxford professor John Wallis (1616--1703) first used the term "hypergeometric" (from the Greek $\acute{\nu}\pi\varepsilon\rho$, above or beyond) to denote any series which was beyond the ordinary geometric series
\[
1+x+x^2+x^3+...
\]

During the next one hundred and fifty years many other mathematicians studied similar series, notably the Swiss L.Euler (1707 -- 1783) who gave amongst many other results [8].

In 1812, C.F.Gauss (1777 -- 1855) delivered his famous thesis [13] before the Royal Society in G\"{o}ttingen. In it, this  brilliant mathematician defined the modern infinite series of (1) above and introduced the notation $F(a,b;c;x)$ for it. He also proved his famous summation theorem
\[
F(a,b;c;1)=\frac{\Gamma(c)\Gamma(c-a-b)}{\Gamma(c-a)\Gamma(c-b)},
\]
and he gave many relations between two or more of these series.

The next major advance was made in 1936 by E.E.Kummer (1810 -- 1893), who first used the term "hypergeometric" for series of the  (1) only.  He showed that the differential equation
\begin{equation}  \label{eqgauss1}
x(1-x)\frac{d^2u}{dx^2}+\left[c-(a+b+1)x\right]\frac{du}{dx}-abu=0
\end{equation}
is satisfied by the function $F(a,b;c;x)$, and has in all twenty-four solutions in terms of similar Gauss functions
 [22]. In 1857  G.F.B.Riemann (1826 -- 1866 )[31] extended this theory by the introduction of his $P$ functions, which in a way are generalizations of the Gaussian $F(a,b;c;x)$.

 Riemann also discussed the general theory of the transformation of the variable in a differential equation and this theory was applied to Kummer's work by J.Thomae [40] who, in 1879, worked out in detail the relationships between Kummer's twenty-four solutions.

 The great success of the theory of hypergeometric series in one variable has stimulated the development of a corresponding theory in two and more variables. Appell    [1]  has defined, in 1880, four series $F_1$ -- $F_4$ (cf. equations  (11) -- (14) \textit{infra}) which are all analogous to Gauss' $F(a,b;c;x)$.   Picard  [29] has pointed out that one of these series is intimately related to a function studied by Pochhammer  [30] in 1870,  and  Picard  [29] and Goursat   [14] also constructed  a theory of Appell's series which is analogous to Riemann's theory [31] of Gauss' hypergeometric series . P.Humbert
  [20] has studied confluent hypergeometric series in two variables. An exposition of the results of the French school together with references to the original literature are to be found in the monograph by Appell and Kamp\'{e} de F\'{e}riet
 [2], which is the standard work on the subject. This work also contains an extensive bibliography of all relevant papers up to 1926.

A great amount of work was also done by Jakob Horn  [18], which he published in a long series of papers extending over fifty years.

 Horn defined thirty types of double series, including series with suffixes of the type $m-n$ as well as the normal $m+n$ suffixes. and he investigated the relationships between them [5, Chap. 5].

The concept of a double hypergeometric series can be extended to triple, quadruple or multiple sums, in general, though the results become progressively more complicated. Such series were first studied by Lauricella  [27] in 1893, whose name they carry.  The theory of general multiple series was investigated more fully by Appell and Kamp\'{e} de F\'{e}riet
 [2,Chap. 5].

The theory of hypergeometric functions of three variables developed relatively slowly. In 1985 H.M. Srivastava Per W. Karlsson [39] established that there exist 205 {complete} triple series of the second order. A. Hasanov and M. Ruzhansky constructed integral representations [16] and also composed 205 systems of partial differential equations  [17] satisfied by complete hypergeometric functions of three variables, and wrote out explicit linearly independent solutions of these systems at the origin, if such solutions exist. Relatively few works are devoted to the study of confluent hypergeometric functions of three variables. Let us note, for example, the work  [21], in which 38 confluent hypergeometric functions of variables are defined and corresponding systems of differential equations of hypergeometric type are composed, and in the work  [9] the author defines confluent forms of the known Lauricella hypergeometric functions of three variables. Some triple confluent series appeared in solving applied problems (see, for example, [15]). In the work  [7] a list of
395 \textit{confluent} series is compiled, which are limit forms for 205 complete series.

The main objective of this work is to establish systems of partial differential equations corresponding to 205 complete and 395 confluent hypergeometric functions of three variables.

In Section 2 we briefly recall the Gauss hypergeometric equation, in Section 3, using the example of the Kummer hypergeometric function, we show that any confluent hypergeometric function is the limit form of the complete hypergeometric function.

In Sections
4, 5 and 6 a hypergeometric functions of two second-order variables (Horn's List), Horn's theorem of convergence  and the corresponding systems of partial differential equations are given with corrections based on the results of [6].

In Section
7 a general definition of a hypergeometric function of three variables  is given.

Sections 8 and
9 are devoted to constructing systems of partial differential equations for 205 total and 395 confluent hypergeometric functions of three variables, respectively. In addition, we determine particular solutions of the constructed systems of equations near the origin, if such solutions exist.

The last section 10 presents a method for determining particular solutions of systems of partial differential equations near the origin.

\section{The Gauss equation}   \label{S2}

The differential equation
\begin{equation}  \label{eqgauss}
x(1-x)\frac{d^2u}{dx^2}+\left[c-(a+b+1)x\right]\frac{du}{dx}-abu=0
\end{equation}
is called the Gauss equation or the hypergeometric equation. In the region  $|x|<1$, one solution is
\[u_1=F(a,b;c;x). \] This can be verified by direct differentiation of the series  (1), and substitution in the above differential equation
(4).
But an alternative form of writing this equation is
\begin{equation}  \label{may10}
\frac{d}{dx}\left(x\frac{d}{dx}+c-1\right)u = \left(x\frac{d}{dx}+a\right)\left(x\frac{d}{dx}+b\right)u
\end{equation}
and this leads to an elegant proof, for
$$
 \left(x\frac{d}{dx}+a\right)u=\sum\limits_{n=0}^\infty\frac{(a)_n(b)_n}{(c)_nn!}(n+a)x^n.
$$
Hence
$$
 \left(x\frac{d}{dx}+a\right)\left(x\frac{d}{dx}+b\right)u=\sum\limits_{n=0}^\infty\frac{(a)_{n+1}(b)_{n+1}}{(c)_nn!}x^n.
$$
Similarly,
$$
 \left(x\frac{d}{dx}+c-1\right)u=\sum\limits_{n=1}^\infty\frac{(a)_{n}(b)_{n}}{(c)_{n-1}n!}x^n.
$$
Hence
$$
 \frac{d}{dx}\left(x\frac{d}{dx}+c-1\right)u=\sum\limits_{n=1}^\infty\frac{(a)_{n}(b)_{n}}{(c)_{n-1}n!}nx^{n-1},
$$
or
$$
 \frac{d}{dx}\left(x\frac{d}{dx}+c-1\right)u=\sum\limits_{n=0}^\infty\frac{(a)_{n+1}(b)_{n+1}}{(c)_{n}n!}x^{n}.
$$

The Gauss equation   (4) can be rewritten
$$
\frac{d^2u}{dx^2}+\left[\frac{c}{x(1-x)}-\frac{a+b+1}{1-x}\right]\frac{du}{dx}-\frac{ab}{x(1-x)}u=0,
$$
from which 0 and 1 are seen to be regular singularities. If we write  $1/x$ for $x$, we find that infinity is also a regular singularity of the Gauss equation [43, Section 10.3].

In the notation of operators, where
$\Delta \equiv x \dfrac{d}{dx}$, the Gauss equation (5) can also be written
$$
\Delta(\Delta +c-1)u=x(\Delta+a)(\Delta +b)u.
$$

\section{Confluent hypergeometric functions}  \label{S3}

If we put $z=\dfrac{x}{b}$ in Gauss' hypergeometric series
$$
F(a,b;c;z)=1+\frac{ab}{1\cdot c}z+\frac{a(a+1)b(b+1)}{1\cdot 2 \cdot c(c+1)}z^2+...
$$
in which we assume that neither  $a$ nor $c$   is zero or a negative integer we obtain a power series in  $x$ whose radius of convergence is    $|b|$  and which defines an analytic function with singularities at  $x=0$, $b$  and $\infty$. As $b \to \infty$, the limiting case will define an entire function whose singularity at  $x=\infty$ is a \textit{confluence} of two singularities of    $F\left(a,b; c; \dfrac{x}{b}\right)$. In this manner we are lead to Kummer's series [23]
\begin{equation} \label{kummerseries}
\Phi(a,c;x)=1+\frac{a}{c}\frac{x}{1!} +\frac{a(a+1)}{c(c+1)}\frac{x^2}{2!}+...
\end{equation}

The series (6) satisfies the differential equation
\begin{equation}  \label{kummer}
x\frac{d^2u}{dx^2}+(c-x)\frac{du}{dx}-au=0.
\end{equation}

The substitution
$$
u=x^{-c/2}e^{x/2}w, \,\,\, a=\frac{1}{2}-k+\mu, \,\,\,c=1+2\mu
$$
reduces (7) to Whittaker's standard form
$$
\frac{d^2w}{dx^2}+ \left(-\frac{1}{4}+\frac{k}{x}+\frac{\frac{1}{4}-\mu^2}{x^2} \right)w=0.
$$
Either of the two equations will be called a \textit{confluent hypergeometric equation}, and any solution of either of them a \textit{confluent hypergeometric function};
 $a$ and $c$ ( or $k$ and $\mu$) will be called the  \textit{parameters}, $x$ the   \textit{variable}.

\section{Hypergeometric series in two variables}  \label{S4}

 Horn [18] gave the following general definition: the double power series

  \begin{equation} \label{defin}
\sum\limits_{m,n=0}^\infty A(m,n)x^my^n
\end{equation}
is a hypergeometric series if the two quotients
  \begin{equation} \label{quot}
f(m,n)=\frac{A(m+1,n)}{A(m,n)},\,\,\,\,\,\,\,g(m,n)=\frac{A(m,n+1)}{A(m,n)}
\end{equation}
are rational functions of  $m$ and $n$.

Horn puts
\begin{equation} \label{quot1}
f(m,n)=\frac{F(m,n)}{F'(m,n)},\,\,\,\,\,\,\,g(m,n)=\frac{G(m,n)}{G'(m,n)},
\end{equation}
where  $F$, $F'$, $G$, $G'$   are polynomials in  $m$, $n$, of respective degrees  $p$, $p'$, $q$, $q'$.  $F'$ is assumed to have a factor  $m+1$, and $G'$ a factor  $n+1$;  $F$ and $F'$  have no common factor except, possibly, $m+1$; and $G$ and  $G'$  no common factor except possibly $n+1$. The highest of the four numbers  $p$, $p'$, $q$, $q'$, is the order of the hypergeometric series. Horn investigated in particular hypergeometric series of \textit{order two} and found that, apart from certain series which are either expressible in terms of one variable or are products of two hypergeometric series, each in one variable, there are essentially 34 distinct convergent series of order two. He distinguished double series into complete and confluent series, that is,  according to Horn's definition, a second-order double hypergeometric series is called \textit{complete} series if $p = p' = q = q'=2$, otherwise \textit{confluent} series.

 There are  14 \textit{complete} series for which  $p=p'=q=q'=2:$

Appell hypergeometric series [1]:
\begin{equation}
\label{i01}
F_{1} \left( {a,b,b';c;x,y} \right) = \sum\limits_{m,n = 0}^{\infty}
{\frac{{(a)_{m + n} (b)_{m} \left(b'\right)_{n}}} {(c)_{m+n}}}\frac{x^{m}}{m!}\frac{y^{n}}{n!},\,\,\,\,\,\,\,\,\,\,\,\,\,\,\,\,\,\,\,\,\,\,\,\,\,\,\,\,\,\,\,\,\,\,\,\,\,\,\,\,\,\,
\end{equation}
\begin{equation}
\label{i02}
F_{2} \left( {a,b,b';c,c';x,y} \right) = {\sum\limits_{m,n = 0}^{\infty}
{{\frac{{({a})_{m+n} ({b})_{m} \left(b'\right)_{n}}} {(c)_m\left(c'\right)_n}}}}\frac{x^{m}}{m!}\frac{y^{n}}{n!},\,\,\,\,\,\,\,\,\,\,\,\,\,\,\,\,\,\,\,\,\,\,\,\,\,\,\,\,\,\,\,\,\,\,
\end{equation}
\begin{equation}
\label{i03}
F_{3} \left( {a,a',b,b';c;x,y} \right) = {\sum\limits_{m,n = 0}^{\infty}
{{\frac{{({a})_{m} \left(a'\right)_{n}(b)_m\left(b'\right)_n}} {{(c)_{m+n}}}}}}\frac{x^{m}}{m!}\frac{y^{n}}{n!},\,\,\,\,\,\,\,\,\,\,\,\,\,\,\,\,\,\,\,\,\,\,\,\,\,
\end{equation}
\begin{equation}
\label{i04}
F_{4} \left( {a,b;c,c';x,y} \right) = {\sum\limits_{m,n = 0}^{\infty}
{{\frac{{({a})_{m + n} ({b})_{m + n} }} {{(c)_m\left(c'\right)_{n} }}}}} \frac{x^{m}}{m!}\frac{y^{n}}{n!},\,\,\,\,\,\,\,\,\,\,\,\,\,\,\,\,\,\,\,\,\,\,\,\,\,\,\,\,\,\,\,\,\,\,\,\,\,\,\,\,\,\,\,\,\,\,\,
\end{equation}

Horn hypergeometric series [19]:

\begin{equation}
\label{g1}
G_{1} \left( {\alpha,\beta,\beta';x,y} \right) = {\sum\limits_{m,n = 0}^{\infty}
{{{{(\alpha)_{m + n} (\beta)_{n - m} (\beta')_{m - n}}} }}}\frac{x^{m}}{m!}\frac{y^{n}}{n!},\,\,\,\,\,\,\,\,\,\,\,\,\,\,\,\,\,\,\,\,\,\,\,\,\,\,\,\,\,
\end{equation}

\begin{equation}
\label{g2}
G_{2} \left( {\alpha,\alpha',\beta,\beta';x,y} \right) = {\sum\limits_{m,n = 0}^{\infty}
{{{{\left( {\alpha} \right)_{m} \left( {\alpha'} \right)_{n} (\beta)_{n - m} (\beta')_{m -
n}}} }}}\frac{x^{m}}{m!}\frac{y^{n}}{n!},\,\,\,\,\,\,\,\,\,
\end{equation}

\begin{equation}
\label{g3}
G_{3} \left( {\alpha,\alpha';x,y} \right) = {\sum\limits_{m,n = 0}^{\infty}
{{{{\left( {\alpha} \right)_{2n - m} (\alpha')_{2m - n}}} }}}\frac{x^{m}}{m!}\frac{y^{n}}{n!},\,\,\,\,\,\,\,\,\,\,\,\,\,\,\,\,\,\,\,\,\,\,\,\,\,\,\,\,\,\,\,\,\,\,\,\,\,\,\,\,\,\,\,\,\,\,\,
\end{equation}

\begin{equation}
\label{h1}
H_{1} \left( {\alpha,\beta,\gamma;\delta;x,y} \right) = {\sum\limits_{m,n = 0}^{\infty}
{{\frac{{\left( {\alpha} \right)_{m - n} \left( {\beta} \right)_{m + n} \left( {\gamma}
\right)_{n}}} {{\left( {\delta} \right)_{m} }}}}} \frac{x^{m}}{m!}\frac{y^{n}}{n!},
\,\,\,\,\,\,\,\,\,\,\,\,\,\,\,\,\,\,\,\,\,\,\,\,\,\,\,
\end{equation}

\begin{equation}
\label{h2}
H_{2} \left( {\alpha,\beta,\gamma,\delta;\varepsilon; x,y} \right) = {\sum\limits_{m,n = 0}^{\infty}
{{\frac{{\left( {\alpha} \right)_{m - n} \left( {\beta} \right)_{m} \left( {\gamma}
\right)_{n} \left( {\delta} \right)_{n}}} {{\left( {\varepsilon} \right)_{m} }}}}
}\frac{x^{m}}{m!}\frac{y^{n}}{n!},\,\,\,\,\,\,\,\,\,\,\,\,\,\,\,\,
\end{equation}

\begin{equation}
\label{h3}
H_{3} \left( {\alpha,\beta;\gamma; x,y} \right) = {\sum\limits_{m,n = 0}^{\infty}
{{\frac{{\left( {\alpha} \right)_{2m + n} \left( {\beta} \right)_{n}}} {{\left( {\gamma}
\right)_{m + n} }}}}} \frac{x^{m}}{m!}\frac{y^{n}}{n!},\,\,\,\,\,\,\,\,\,\,\,\,\,\,\,\,\,\,\,\,\,\,\,\,\,\,\,\,\,\,\,\,\,\,\,\,\,\,\,\,\,\,\,\,\,\,\,
\end{equation}

\begin{equation}
\label{h4}
H_{4} \left( {\alpha,\beta;\gamma,\delta; x,y} \right) = {\sum\limits_{m,n = 0}^{\infty}
{{\frac{{\left( {\alpha} \right)_{2m + n} \left( {\beta} \right)_{n}}} {{\left( {\gamma}
\right)_{m} \left( {\delta} \right)_{n} }}}}} \frac{x^{m}}{m!}\frac{y^{n}}{n!},\,\,\,\,\,\,\,\,\,\,\,\,\,\,\,\,\,\,\,\,\,\,\,\,\,\,\,\,\,\,\,\,\,\,\,\,
\end{equation}

\begin{equation}
\label{h5}
H_{5} \left( {\alpha,\beta;\gamma; x,y} \right) = {\sum\limits_{m,n = 0}^{\infty}
{{\frac{{\left( {\alpha} \right)_{2m + n} \left( {\beta} \right)_{n - m}}} {{\left(
{\gamma} \right)_{n} }}}}} \frac{x^{m}}{m!}\frac{y^{n}}{n!},\,\,\,\,\,\,\,\,\,\,\,\,\,\,\,\,\,\,\,\,\,\,\,\,\,\,\,\,\,\,\,\,\,\,\,\,\,\,\,\,\,\,
\end{equation}

\begin{equation}
\label{h6}
H_{6} \left( {\alpha,\beta,\gamma; x,y} \right) = \sum\limits_{m,n = 0}^{\infty}
{{{\left( {\alpha} \right)_{2m - n} \left( {\beta} \right)_{n - m} \left( {\gamma}
\right)_{n}}}} \frac{x^{m}}{m!}\frac{y^{n}}{n!},\,\,\,\,\,\,\,\,\,\,\,\,\,\,\,\,\,\,\,\,\,\,\,\,\,\,\,\,
\end{equation}

\begin{equation}
\label{h7}
H_{7} \left( {\alpha,\beta,\gamma;\delta;x,y} \right) = {\sum\limits_{m,n = 0}^{\infty}
{{\frac{{\left( {\alpha} \right)_{2m - n} \left( {\beta} \right)_{n} \left( {\gamma}
\right)_{n}}} {{\left( {\delta} \right)_{m} }}}}} \frac{x^{m}}{m!}\frac{y^{n}}{n!},\,\,\,\,\,\,\,\,\,\,\,\,\,\,\,\,\,\,\,\,\,\,\,\,\,\,\,\,\,\,
\end{equation}
and there are 20 \textit{confluent} series which are limiting forms of the complete ones and for which $p\leq p'=2$ , $q\leq q'=2$ and $p$, $q$ not both $=2$:

Humbert confluent hypergeometric series [20]:
\begin{equation}
\label{rf1}
\Phi _{1} \left( {\alpha,\beta; \gamma;x,y} \right) = {\sum\limits_{m,n = 0}^{\infty}
{{\frac{{\left( {\alpha} \right)_{m+n} \left( {\beta} \right)_{m} }} {{(\gamma)_{m+n}}}}}} \frac{x^{m}}{m!}\frac{y^{n}}{n!},\,\,\,|x|<1,\,\,\,\,\,\,\,\,\,\,\,\,\,\,\,\,\,\,\,\,\,\,\,\,\,\,\,\,\,\,\,\,\,\,\,\,
\end{equation}

\begin{equation}
\label{rf2}
\Phi _{2} \left( {\beta, \beta';\gamma;x,y} \right) = {\sum\limits_{m,n = 0}^{\infty}
{{\frac{{\left( {\beta} \right)_{m} \left( {\beta'} \right)_{n}}} {{(\gamma)_{m+n}}}}}
}\frac{x^{m}}{m!}\frac{y^{n}}{n!},\,\,\,\,\,\,\,\,\,\,\,\,\,\,\,\,\,\,\,\,\,\,\,\,\,\,\,\,\,\,\,\,\,\,\,\,\,\,\,\,\,\,\,\,\,\,\,\,\,\,\,\,\,\,\,\,\,\,\,\,
\end{equation}

\begin{equation}
\label{rf3}
{\Phi}_{3} \left( {\beta;\gamma;x,y} \right) = {\sum\limits_{m,n = 0}^{\infty}
{{\frac{{\left( {\beta} \right)_{m} }} {{\left(
{\gamma} \right)_{m+n} }}}}} \frac{x^{m}}{m!}\frac{y^{n}}{n!},\,\,\,\,\,\,\,\,\,\,\,\,\,\,\,\,\,\,\,\,\,\,\,\,\,\,\,\,\,\,\,\,\,\,\,\,\,\,\,\,\,\,\,\,\,\,\,\,\,\,\,\,\,\,\,\,\,\,\,\,\,\,\,\,
\end{equation}

\begin{equation}
\label{ph1}
{\Psi}_{1} \left( {\alpha,\beta;\gamma,\gamma';x,y} \right) = {\sum\limits_{m,n = 0}^{\infty}
{{\frac{{\left( {\alpha} \right)_{m + n} \left( {\beta} \right)_{m} }} {{(\gamma)_m ({\gamma'})_{n} }}}}} \frac{x^{m}}{m!}\frac{y^{n}}{n!},\,\,\,|x|<1,\,\,\,\,\,\,\,\,\,\,\,\,\,\,\,\,\,\,\,\,\,\,
\end{equation}

\begin{equation}
\label{ph2}
{\Psi}_{2} \left( {\alpha;\gamma,\gamma';x,y} \right) = {\sum\limits_{m,n = 0}^{\infty}
{{\frac{{\left( {\alpha} \right)_{m + n} }} {{(\gamma)_m\left( {\gamma'}
\right)_{n}}}}}} \frac{x^{m}}{m!}\frac{y^{n}}{n!},
\,\,\,\,\,\,\,\,\,\,\,\,\,\,\,\,\,\,\,\,\,\,\,\,\,\,\,\,\,\,\,\,\,\,\,\,\,\,\,\,\,\,\,\,\,\,\,\,\,
\end{equation}

\begin{equation}
\label{kh1}
{\Xi}_{1} \left( {\alpha,\alpha',\beta;\gamma;x,y} \right) = {\sum\limits_{m,n = 0}^{\infty}
{{\frac{{\left( {\alpha} \right)_{m} \left( {\alpha'} \right)_{n}(\beta)_m}} {{\left( {\gamma}
\right)_{m+n}}}}}} \frac{x^{m}}{m!}\frac{y^{n}}{n!},\,\,\,|x|<1,\,\,\,\,\,\,\,\,\,\,\,\,\,\,\,\,
\end{equation}

\begin{equation}
\label{kh2}
{\Xi}_{2} \left( {\alpha,\beta;\gamma;x,y} \right) = {\sum\limits_{m,n = 0}^{\infty}
{{\frac{{\left( {\alpha} \right)_{m}(\beta)_m}} {{\left( {\gamma} \right)_{m+n}}}}}
}\frac{x^{m}}{m!}\frac{y^{n}}{n!},\,\,\,|x|<1,\,\,\,\,\,\,\,\,\,\,\,\,\,\,\,\,\,\,\,\,\,\,\,\,\,\,\,\,\,\,\,\,
\end{equation}

Horn confluent hypergeometric series [19]:

\begin{equation}
\label{rg1}
\Gamma _{1} \left( {\alpha,\beta,\beta';x,y} \right) = {\sum\limits_{m,n = 0}^{\infty}
{{{{\left( {\alpha} \right)_{m} \left( {\beta} \right)_{n - m} \left( {\beta'}
\right)_{m - n}}} }}} \frac{x^{m}}{m!}\frac{y^{n}}{n!},\,\,\,|x|<1,\,\,\,\,\,\,\,
\end{equation}

\begin{equation}
\label{rg2}
\Gamma _{2} \left( {\beta,\beta';x,y} \right) = {\sum\limits_{m,n = 0}^{\infty}
{{{{\left( {\beta} \right)_{n - m} \left( {\beta'} \right)_{m - n}}} }}
}\frac{x^{m}}{m!}\frac{y^{n}}{n!},\,\,\,\,\,\,\,\,\,\,\,\,\,\,\,\,\,\,\,\,\,\,\,\,\,\,\,\,\,\,\,\,\,\,\,\,\,\,\,\,\,\,\,\,\,\,\,\,\,\,
\end{equation}

\begin{equation}
\label{rh1}
{\rm H}_{1} \left( {\alpha,\beta;\delta;x,y} \right) = {\sum\limits_{m,n = 0}^{\infty}
{{\frac{{\left( {\alpha} \right)_{m - n} \left( {\beta} \right)_{m + n}}} {{\left(
{\delta} \right)_{m}}}}}} \frac{x^{m}}{m!}\frac{y^{n}}{n!},\,\,\,|x|<1,\,\,\,\,\,\,\,\,\,\,\,\,\,\,\,\,\,\,\,\,\,\,
\end{equation}

\begin{equation}
\label{rh2}
{\rm H}_{2} \left( {\alpha,\beta,\gamma;\delta;x,y} \right) = {\sum\limits_{m,n = 0}^{\infty}
{{\frac{{\left( {\alpha} \right)_{m - n} \left( {\beta} \right)_{m} \left( {\gamma}
\right)_{n}}} {{\left( {\delta} \right)_{m}}}}}} \frac{x^{m}}{m!}\frac{y^{n}}{n!},\,\,\,|x|<1,\,\,\,\,\,\,\,\,\,\,
\end{equation}

\begin{equation}
\label{rh3}
{\rm H}_{3} \left( {\alpha,\beta;\delta;x,y} \right) = {\sum\limits_{m,n = 0}^{\infty}
{{\frac{{\left( {\alpha} \right)_{m - n} \left( {\beta} \right)_{m}}} {{\left( {\delta}
\right)_{m} }}}}} \frac{x^{m}}{m!}\frac{y^{n}}{n!},\,\,\,|x|<1,\,\,\,\,\,\,\,\,\,\,\,\,\,\,\,\,\,\,\,\,\,\,\,\,\,\,\,
\end{equation}

\begin{equation}
\label{rh4}
{\rm H}_{4} \left( {\alpha,\gamma;\delta;x,y} \right) = {\sum\limits_{m,n = 0}^{\infty}
{{\frac{{\left( {\alpha} \right)_{m - n} \left( {\gamma} \right)_{n}}} {{\left( {\delta}
\right)_{m}}}}}} \frac{x^{m}}{m!}\frac{y^{n}}{n!},
\,\,\,\,\,\,\,\,\,\,\,\,\,\,\,\,\,\,\,\,\,\,\,\,\,\,\,\,\,\,\,\,\,\,\,\,\,\,\,\,\,\,\,\,\,\,\,\,\,\,\,\,
\end{equation}

\begin{equation}
\label{rh5}
{\rm H}_{5} \left( {\alpha;\delta;x,y} \right) = {\sum\limits_{m,n = 0}^{\infty}
{{\frac{{\left( {\alpha} \right)_{m - n}}} {{\left( {\delta} \right)_{m}}}}}
}\frac{x^{m}}{m!}\frac{y^{n}}{n!},\,\,\,\,\,\,\,\,\,\,\,\,\,\,\,\,\,\,\,\,\,\,\,\,\,\,\,\,\,\,\,\,\,\,\,\,\,\,\,\,\,\,\,\,\,\,\,\,\,\,\,\,\,\,\,\,\,\,\,\,\,\,\,\,
\end{equation}

\begin{equation}
\label{rh6}
{\rm H}_{6} \left( {\alpha;\gamma;x,y} \right) = {\sum\limits_{m,n = 0}^{\infty}
{{\frac{{\left( {\alpha} \right)_{2m + n}}} {{\left( {\gamma} \right)_{m + n}}}}}
}\frac{x^{m}}{m!}\frac{y^{n}}{n!},\,\,\,|x|<\frac{1}{4},\,\,\,\,\,\,\,\,\,\,\,\,\,\,\,\,\,\,\,\,\,\,\,\,\,\,\,\,\,\,\,
\end{equation}

\begin{equation}
\label{rh7}
{\rm H}_{7} \left( {\alpha;\gamma,\delta;x,y} \right) = {\sum\limits_{m,n = 0}^{\infty}
{{\frac{{\left( {\alpha} \right)_{2m + n}}} {{\left( {\gamma} \right)_{m} \left( {\delta}
\right)_{n}}}}}}\frac{x^{m}}{m!}\frac{y^{n}}{n!} , |x|<\frac{1}{4},\,\,\,\,\,\,\,\,\,\,\,\,\,\,\,\,\,\,\,\,\,\,\,\,
\end{equation}

\begin{equation}
\label{rh8}
{\rm H}_{8} \left( {\alpha,\beta;x,y} \right) = {\sum\limits_{m,n = 0}^{\infty}
{{{{\left( {\alpha} \right)_{2m - n} \left( {\beta} \right)_{n - m}}} }}
}\frac{x^{m}}{m!}\frac{y^{n}}{n!},\,\,\,|x|<\frac{1}{4},\,\,\,\,\,\,\,\,\,\,\,\,\,\,\,\,\,\,\,\,
\end{equation}

\begin{equation}
\label{rh9}
{\rm H}_{9} \left( {\alpha,\beta;\delta;x,y} \right) = {\sum\limits_{m,n = 0}^{\infty}
{{\frac{{\left( {\alpha} \right)_{2m - n} \left( {\beta} \right)_{n}}} {{\left( {\delta}
\right)_{m}}}}}} \frac{x^{m}}{m!}\frac{y^{n}}{n!},\,\,\,|x|<\frac{1}{4},\,\,\,\,\,\,\,\,\,\,\,\,\,\,\,\,\,\,\,\,\,\,
\end{equation}

\begin{equation}
\label{rh10}
{\rm H}_{10} \left( {\alpha;\delta;x,y} \right) = {\sum\limits_{m,n = 0}^{\infty}
{{\frac{{\left( {\alpha} \right)_{2m - n}}} {{\left( {\delta} \right)_{m}
}}}}} \frac{x^{m}}{m!}\frac{y^{n}}{n!},\,\,\,|x|<\frac{1}{4},\,\,\,\,\,\,\,\,\,\,\,\,\,\,\,\,\,\,\,\,\,\,\,\,\,\,\,\,\,\,\,\,\,\,
\end{equation}

\begin{equation}
\label{rh1111}
{\rm H}_{11} \left( {\alpha,\beta,\gamma;\delta;x,y} \right) = {\sum\limits_{m,n = 0}^{\infty}
{{\frac{{\left( {\alpha} \right)_{m - n} \left( {\beta} \right)_{n} \left( {\gamma}
\right)_{n}}} {{\left( {\delta} \right)_{m}}}}}} \frac{x^{m}}{m!}\frac{y^{n}}{n!},\,\,\,|y|<1.\,\,\,
\end{equation}

Horn studied the convergence of hypergeometric series in two variables. The domains of the hypergeometric functions of Appell and Horn, defined in (11)--(24), can be found, for example, in [5, p. 227--229]. In the case of confluent series either $\Phi$ or $\Psi$  vanishes identically, the region of convergence simplifies considerably, and any inequalities which may be necessary to secure convergence are recorded in
(25) to  (44).

\section{Horn's theorem on convergence}  \label{convergence}

The problem of convergence hardly exists for single hypergeometric series. The region of convergence of the Gaussian series (1) is obtained immediately by appearing to d'Alambert's ratio test. A generalization of this test was given by Horn [18] whose theorem on the convergence of double and triple hypergeometric series is considered in greater detail in Sections 2.2 and 5.1 of the monograph by Srivastava and Karlsson [39]. 

In the present work we use the convergence results of [39]. 

For future we introduce a following notations [39, p. 58]:

\begin{equation}
\Phi_1(\xi)=\frac{2A(\xi)+1}{3[A(\xi)+1]^2}, \,\,\,\,\,\,A(\xi)=\sqrt{1+3\xi}
\end{equation}

\begin{equation}
\Phi_2(\xi)=\frac{2B(\xi)-1}{3[B(\xi)-1]^2}, \,\,\,\,\,\,B(\xi)=\sqrt{1-3\xi}
\end{equation}

\begin{equation}
\Psi_1(\xi)=\frac{2[2-a(\xi)]^2}{9[a(\xi)-1]}, \,\,\,\,\,\,a(\xi)=\sqrt{1+12\xi}
\end{equation}

\begin{equation}
\Psi_2(\xi)=\frac{2[2+b(\xi)]^2}{9[b(\xi)+1]}, \,\,\,\,\,\,b(\xi)=\sqrt{1-12\xi}
\end{equation}

\begin{equation}
\Theta_1(\xi)=\frac{[1+3\alpha(\xi)]\,[\alpha(\xi)-1]}{12[\alpha(\xi)+1]^2}, \,\,\,\,\,\,\alpha(\xi)=\sqrt{1+\frac{8}{9\xi}}
\end{equation}

\begin{equation}
\Theta_2(\xi)=\frac{[1-3\beta(\xi)]\,[\beta(\xi)+1]}{12[\beta(\xi)-1]^2}, \,\,\,\,\,\,\beta(\xi)=\sqrt{1-\frac{8}{9\xi}}.
\end{equation}

\section{Systems of partial differential equations associated with the hypergeometric functions in two variables} \label{SS}

The series
\begin{equation}
\sum\limits_{m,n=0}^\infty A(m,n)x^my^n,
\end{equation}
where

$$
\frac{{{A({m + 1,n})}}}{{{A({m,n})}}} = \frac{{F\left( {m,n} \right)}}{{F'\left( {m,n} \right)}},\,\,\,\frac{{{A({m,n + 1})}}}{{{A({m,n})}}} = \frac{{G\left( {m,n} \right)}}{{G'\left( {m,n} \right)}},
$$
and  $F\left( {m,n} \right),F'\left( {m,n} \right),\,\,\,G\left( {m,n} \right),G'\left( {m,n} \right)$ , are polynomials as in  (10),  satisfies a system of linear partial differential equations which can be written in terms of the differential operators

$$
 {\delta} = x\frac{\partial }{{\partial x}},\,\,\,{\delta'} = y\frac{\partial }{{\partial y}}
 $$
as
$$
\left\{ {\begin{array}{*{20}{c}}
  {\left[ {F'\left( {{\delta},{\delta'}} \right){x^{ - 1}} - F\left( {{\delta},{\delta'}} \right)} \right]u = 0}, \\
  {}\\
  {\left[ {G'\left( {{\delta},{\delta'}} \right){y^{ - 1}} - G\left( {{\delta},{\delta'}} \right)} \right]u = 0}. \\
 \end{array}} \right.
$$

In the book [5, Section 5.9]  all 34 systems of differential equations of hypergeometric type are given, which are satisfied by the second-order hypergeometric functions in two variables.     In 1991, Volkodavov and Bystrova  [41]  drew attention to the typos in the system corresponding to the confluent hypergeometric function $\rm{H}_3$, defined by equality (36).

 A detailed discussion showed that in Section 5.9 of  [5] in 26 cases the hypergeometric function from Horn's List actually satisfies the corresponding system of differential equations, and in 8 cases the systems of differential equations are erroneously given, which are allegedly satisfied by the functions $H_4$, $H_5$, $H_7$, $\Gamma_1$, ${\rm{H}}_2$, ${\rm{H}}_3$, ${\rm{H}}_5$ and ${\rm{H}}_7$ from Horn's List. It is possible that these circumstances have been a brake on the applications of these functions.

 We note only the work [6], in which the definitions of hypergeometric functions of two variables are clarified, the truth of all systems of differential equations of hypergeometric type is rechecked according to the original sources, and in cases where they are absent, the process of compiling systems of differential equations with partial derivatives is discussed in detail.

Marichev  [25, 26, 27] also pointed out possible typos in the  Sections  5.7, 5.8, and 5.11 of the book [5],  devoted to hypergeometric functions of two variables.

In the following list of partial differential equations $z$ is the unknown function of $x$ and $y$,
$$
p=\frac{\partial z}{\partial x},\,\,\,\,\,q=\frac{\partial z}{\partial y},\,\,\,\,\,r=\frac{\partial^2 z}{\partial x^2},\,\,\,\,\,s=\frac{\partial^2 z}{\partial x \partial y},\,\,\,\,\,t=\frac{\partial^2 z}{\partial y^2}.\,\,\,\,\,\,\,\,\,\,\,\,\,\,\,\,\,\,\,\,\,\,\,\,
$$
\begin{equation} \label{difeq09}
\left.\begin{array}{*{20}c}
 {x(1-x)r+(1-x)ys+\left[\gamma-\left(\alpha+\beta+1\right)x\right]p-\beta yq-\alpha \beta z=0} \hfill \\
 {y(1-y)t+x(1-y)s+\left[\gamma-(\alpha+\beta'+1)y\right]q-\beta'xp-\alpha \beta'z=0} \hfill \\
\end{array}\right\}F_1;\,\,\,\,\,\,\,\,\,\,\,\,\,\,\,\,\,\,\,\,\,\,\,\,\,\,\,\,\,\,\,\,\,\,\,\,\,\,\,\,\,\,\,\,\,\,\,\,\,\,\,\,\,\,\,\,\,\,\,\,\,\,\,\,
\end{equation}
\begin{equation} \label{difeq10}
\left.\begin{array}{*{20}c}
 {x(1-x)r-xys+\left[\gamma-\left(\alpha+\beta+1\right)x\right]p-\beta yq-\alpha \beta z=0} \hfill \\
 {y(1-y)t-xys+\left[\gamma'-(\alpha+\beta'+1)y\right]q-\beta'xp-\alpha \beta'z=0} \hfill \\
\end{array}\right\}F_2,\,\,\,\,\,\,\,\,\,\,\,\,\,\,\,\,\,\,\,\,\,\,\,\,\,\,\,\,\,\,\,\,\,\,\,\,\,\,\,\,\,\,\,\,\,\,\,\,\,\,\,\,\,\,\,\,\,\,\,\,\,\,\,\,\,\,\,\,\,\,
\end{equation}

particular solutions:

$
{z_1} = F_2\left( \alpha, \beta, \beta'; \gamma, \gamma'; x,y\right),
$

$
{z_2} = {x^{1 - \gamma}}F_2\left( 1-\gamma+\alpha, 1-\gamma+ \beta, \beta'; 2-\gamma, \gamma'; x,y\right),
$

$
z_3= y^{1-\gamma'}F_2\left( 1-\gamma'+\alpha, \beta, 1-\gamma'+\beta'; \gamma, 2-\gamma'; x,y\right),
$

$
z_4= x^{1-\gamma}y^{1-\gamma'}F_2\left( 2-\gamma-\gamma'+\alpha, 1-\gamma+\beta, 1-\gamma'+\beta'; 2-\gamma, 2-\gamma'; x,y\right);
$

\begin{equation} \label{difeq11}
\left.\begin{array}{*{20}c}
 {x(1-x)r+ys+\left[\gamma-\left(\alpha+\beta+1\right)x\right]p-\alpha \beta z=0} \hfill \\
 {y(1-y)t+xs+\left[\gamma-(\alpha'+\beta'+1)y\right]q-\alpha' \beta'z=0} \hfill \\
\end{array}\right\}F_3;\,\,\,\,\,\,\,\,\,\,\,\,\,\,\,\,\,\,\,\,\,\,\,\,\,\,\,\,\,\,\,\,\,\,\,\,\,\,\,\,\,\,\,\,\,\,\,\,\,\,\,\,\,\,\,\,\,\,\,\,\,\,\,\,\,\,\,\,\,\,\,\,\,\,\,\,\,\,\,\,\,\,\,\,\,\,\,\,\,\,\,\,\,\,\,\,\,
\end{equation}
\begin{equation} \label{difeq12}
\left.\begin{array}{*{20}c}
 {x(1-x)r-y^2t-2xys+\left[\gamma-\left(\alpha+\beta+1\right)x\right]p-(\alpha+\beta+1) yq-\alpha \beta z=0} \hfill \\
 {y(1-y)t-x^2r-2xys+\left[\gamma'-(\alpha+\beta+1)y\right]q-(\alpha+\beta+1)xp-\alpha \beta z=0} \hfill \\
\end{array}\right\}F_4,\,\,\,\,\,\,\,\,\,\,\,\,\,\,\,\,\,\,\,\,\,\,\,\,\,\,\,\,\,\,\,\,\,\,\,\,
\end{equation}

particular solutions:

$
{z_1} = F_4\left( \alpha, \beta; \gamma, \gamma'; x,y\right),
$

$
{z_2} = {x^{1 - \gamma}}F_4\left( 1-\gamma+\alpha, 1-\gamma+ \beta; 2-\gamma, \gamma'; x,y\right),
$

$
z_3= y^{1-\gamma'}F_4\left( 1-\gamma'+\alpha, 1-\gamma'+\beta; \gamma, 2-\gamma'; x,y\right),
$

$
z_4= x^{1-\gamma}y^{1-\gamma'}F_4\left( 2-\gamma-\gamma'+\alpha, 2-\gamma-\gamma'+\beta; 2-\gamma, 2-\gamma'; x,y\right);
$

\begin{equation} \label{difeq13}
\left.\begin{array}{*{20}c}
 {x(1+x)r-y^2t-ys+\left[1-\beta+\left(\alpha+\beta'+1\right)x\right]p+(\beta'-\alpha-1) yq+\alpha \beta' z=0} \hfill \\
 {y(1+y)t-x^2r-xs+\left[1-\beta'+(\alpha+\beta+1)y\right]q+(\beta-\alpha-1)xp+\alpha \beta z=0} \hfill \\
\end{array}\right\}G_1;\,\,\,\,\,\,\,\,\,\,\,\,\,\,\,\,\,\,\,\,\,\,\,\,
\end{equation}
\begin{equation} \label{difeq14}
\left.\begin{array}{*{20}c}
 {x(1+x)r-(1+x)ys+\left[1-\beta+\left(\alpha+\beta'+1\right)x\right]p-\alpha yq+\alpha \beta' z=0} \hfill \\
 {y(1+y)t-x(1+y)s+\left[1-\beta'+(\alpha'+\beta+1)y\right]q-\alpha'xp+\alpha' \beta z=0} \hfill \\
\end{array}\right\}G_2;\,\,\,\,\,\,\,\,\,\,\,\,\,\,\,\,\,\,\,\,\,\,\,\,\,\,\,\,\,\,\,\,\,\,\,\,\,\,\,\,\,\,\,\,\,\,\,\,\,
\end{equation}
\begin{equation} \label{difeq15}
\left.\begin{array}{*{20}c}
 {x(1+4x)r-2(1+2x)ys+y^2t+\left[1-\alpha+2\left(2\alpha'+3\right)x\right]p-2\alpha'y q+\alpha'\left(\alpha'+1\right) z=0} \hfill \\
 {y(1+4y)t-2x(1+2y)s+x^2r+\left[1-\alpha'+2(2\alpha+3)y\right]q-2\alpha xp+\alpha\left( \alpha+1\right) z=0} \hfill \\
\end{array}\right\}G_3;\,\,\,\,\,\,\,
\end{equation}
\begin{equation} \label{difeq16}
\left.\begin{array}{*{20}c}
 {x(1-x)r+y^2t+\left[\delta-\left(\alpha+\beta+1\right)x\right]p+(1-\alpha+\beta) yq-\alpha \beta z=0} \hfill \\
 {y(1+y)t-x(1-y)s+\left[1-\alpha+(\beta+\gamma+1)y\right]q+\gamma xp+ \beta \gamma z=0} \hfill \\
\end{array}\right\}H_1,\,\,\,\,\,\,\,\,\,\,\,\,\,\,\,\,\,\,\,\,\,\,\,\,\,\,\,\,\,\,\,\,\,\,\,\,\,\,\,\,\,\,\,\,\,\,\,\,\,\,\,\,\,\,
\end{equation}

particular solutions:

$
{z_1} = H_1\left( \alpha, \beta,  \gamma; \delta;  x,y\right),
$

$
{z_2} = {x^{1 - \delta}}H_1\left( 1-\delta+\alpha, 1-\delta+\beta, \gamma; 2-\delta; x,y\right);
$

\begin{equation} \label{difeq17}
\left.\begin{array}{*{20}c}
 {x(1-x)r+xys+\left[\varepsilon-\left(\alpha+\beta+1\right)x\right]p+\beta yq-\alpha \beta z=0} \hfill \\
 {y(1+y)t-xs+\left[1-\alpha+(\gamma+\delta+1)y\right]q+\gamma \delta z=0} \hfill \\
\end{array}\right\}H_2,\,\,\,\,\,\,\,\,\,\,\,\,\,\,\,\,\,\,\,\,\,\,\,\,\,\,\,\,\,\,\,\,\,\,\,\,\,\,\,\,\,\,\,\,\,\,\,\,\,\,\,\,\,\,\,\,\,\,\,\,\,\,\,\,\,\,\,\,\,\,\,\,\,\,\,\,\,\,\,\,
\end{equation}

particular solutions:

$
{z_1} = H_2\left( \alpha, \beta,  \gamma, \delta; \varepsilon; x,y\right),
$

$
{z_2} = {x^{1 - \varepsilon}}H_2\left( 1-\varepsilon+\alpha, 1-\varepsilon+\beta, \gamma, \delta; 2-\varepsilon; x,y\right);
$

\begin{equation} \label{difeq18}
\left.\begin{array}{*{20}c}
 {x(1-4x)r+(1-4x)ys-y^2t+\left[\gamma-2\left(2\alpha +3\right)x\right]p-2(\alpha+1)y q-\alpha\left(\alpha+1\right) z=0} \hfill \\
 {y(1-y)t+x(1-2y)s+\left[\gamma-(\alpha+\beta+1)y\right]q-2\beta xp-\alpha\beta z=0} \hfill \\
\end{array}\right\}H_3;\,\,\,\,\,\,\,\,\,\,\,\,
\end{equation}
\begin{equation} \label{difeq19}
\left.\begin{array}{*{20}c}
 {x(1-4x)r-4xys-y^2t+\left[\gamma-2\left(2\alpha+3\right)x\right]p-2(\alpha+1)y q-\alpha(\alpha+1) z=0} \hfill \\
 {y(1-y)t-2xys+\left[\delta-(\alpha+\beta)y\right]q-2\beta xp-\alpha\beta z=0} \hfill \\
\end{array}\right\}H_4,\,\,\,\,\,\,\,\,\,\,\,\,\,\,\,\,\,\,\,\,\,\,\,\,\,\,\,\,
\end{equation}

particular solutions:

$
{z_1} = H_4\left( \alpha, \beta; \gamma, \delta; x,y\right),
$

$
{z_2} = {x^{1 - \gamma}}H_4\left( 2-2\gamma+\alpha, \beta; 2-\gamma, \delta; x,y\right),
$

$
z_3= y^{1-\delta}H_4\left( 1-\delta+\alpha, 1-\delta+\beta; \gamma, 2-\delta; x,y\right),
$

$
z_4= x^{1-\gamma}y^{1-\delta}H_4\left( 3-2\gamma-\delta+\alpha, 1-\delta+\beta; 2-\gamma, 2-\delta; x,y\right);
$

\begin{equation} \label{difeq20}
\left.\begin{array}{*{20}c}
 {x(1+4x)r-(1-4x)ys+y^2t+\left[1-\beta+2\left(2\alpha+3\right)x\right]p+2(\alpha+1)y q+\alpha(\alpha+1) z=0} \hfill \\
 {y(1-y)t-xys+2x^2r+\left[\gamma-(\alpha+\beta+1)y\right]q+(2+\alpha-2\beta) xp-\alpha\beta z=0} \hfill \\
\end{array}\right\}H_5,
\end{equation}

particular solutions:

$
{z_1} = H_5\left( \alpha, \beta; \gamma; x,y\right),
$

$
{z_2} = {y^{1 - \gamma}}H_5\left( 1-\gamma+\alpha, 1-\gamma+ \beta; 2-\gamma; x,y\right),
$

\begin{equation} \label{difeq21}
\left.\begin{array}{*{20}c}
 {x(1+4x)r-(1+4x)ys+y^2t+\left[1-\beta+2(2\alpha+3)x\right]p-2\alpha y q+\alpha(\alpha+1) z=0} \hfill \\
 {y(1+y)t-x(2+y)s+\left[1-\alpha+(\beta+\gamma+1)y\right]q-\gamma xp+\beta \gamma z=0} \hfill \\
\end{array}\right\}H_6;\,\,\,\,\,\,\,\,\,\,\,\,\,\,\,\,\,
\end{equation}

\begin{equation} \label{difeq22}
\left.\begin{array}{*{20}c}
 {x(1-4x)r+4xys-y^2t+\left[\delta-2(2\alpha+3)x\right]p+2\alpha y q-\alpha(\alpha+1) z=0} \hfill \\
 {y(1+y)t-2xs+\left[1-\alpha+(\beta+\gamma+1)y\right]q+\beta \gamma z=0} \hfill \\
\end{array}\right\}H_7,\,\,\,\,\,\,\,\,\,\,\,\,\,\,\,\,\,\,\,\,\,\,\,\,\,\,\,\,\,\,\,\,\,\,\,\,\,\,\,\,\,\,\,
\end{equation}

particular solutions:

$
{z_1} = H_7\left( \alpha, \beta, \gamma; \delta; x,y\right),
$

$
{z_2} = {x^{1 - \gamma}}H_7\left( 2-2\gamma+\alpha, \beta, \gamma; 2-\delta; x,y\right),
$

\begin{equation} \label{difeq23}
\left.\begin{array}{*{20}c}
 {x(1-x)r+(1-x)ys+\left[\gamma-\left(\alpha+\beta+1\right)x\right]p-\beta yq-\alpha \beta z=0} \hfill \\
 {yt+xs+(\gamma-y)q- xp-\alpha z=0} \hfill \\
\end{array}\right\}\Phi_1;\,\,\,\,\,\,\,\,\,\,\,\,\,\,\,\,\,\,\,\,\,\,\,\,\,\,\,\,\,\,\,\,\,\,\,\,\,\,\,\,\,\,\,\,\,\,\,\,\,\,\,\,\,\,\,\,\,\,\,\,\,\,\,
\end{equation}
\begin{equation} \label{difeq24}
\left.\begin{array}{*{20}c}
 {xr+ys+(\gamma-x)p- \beta z=0} \hfill \\
 {yt+xs+(\gamma-y)q -\beta' z=0} \hfill \\
\end{array}\right\}\Phi_2;\,\,\,\,\,\,\,\,\,\,\,\,\,\,\,\,\,\,\,\,\,\,\,\,\,\,\,\,\,\,\,\,\,\,\,\,\,\,\,\,\,\,\,\,\,\,\,\,\,\,\,\,\,\,\,\,\,\,\,\,\,\,\,\,\,\,\,\,\,\,\,\,\,\,\,\,\,\,\,\,\,\,\,\,\,\,\,\,\,\,\,\,\,\,\,\,\,\,\,\,\,\,\,\,\,\,\,\,\,\,\,\,\,\,\,\,\,\,\,\,\,\,\,\,\,\,\,\,\,\,\,\,\,\,\,\,\,\,\,\,\,\,\,\,\,\,\,\,\,\,
\end{equation}
\begin{equation} \label{difeq25}
\left.\begin{array}{*{20}c}
 {xr+ys+(\gamma-x)p- \beta z=0} \hfill \\
 {yt+xs+\gamma q - z=0} \hfill \\
\end{array}\right\}\Phi_3;\,\,\,\,\,\,\,\,\,\,\,\,\,\,\,\,\,\,\,\,\,\,\,\,\,\,\,\,\,\,\,\,\,\,\,\,\,\,\,\,\,\,\,\,\,\,\,\,\,\,\,\,\,\,\,\,\,\,\,\,\,\,\,\,\,\,\,\,\,\,\,\,\,\,\,\,\,\,\,\,\,\,\,\,\,\,\,\,\,\,\,\,\,\,\,\,\,\,\,\,\,\,\,\,\,\,\,\,\,\,\,\,\,\,\,\,\,\,\,\,\,\,\,\,\,\,\,\,\,\,\,\,\,\,\,\,\,\,\,\,\,\,\,\,\,\,\,\,\,\,
\end{equation}
\begin{equation} \label{difeq26}
\left.\begin{array}{*{20}c}
 {x(1-x)r-xys+\left[\gamma-\left(\alpha+\beta+1\right)x\right]p-\beta yq-\alpha \beta z=0} \hfill \\
 {yt+\left(\gamma'-y\right)q- xp-\alpha z=0} \hfill \\
\end{array}\right\}\Psi_1,\,\,\,\,\,\,\,\,\,\,\,\,\,\,\,\,\,\,\,\,\,\,\,\,\,\,\,\,\,\,\,\,\,\,\,\,\,\,\,\,\,\,\,\,\,\,\,\,\,\,\,\,\,\,\,\,\,\,\,\,\,\,\,\,\,\,\,\,\,\,\,\,\,\,\,\,\,\,
\end{equation}

particular solutions:

$
{z_1} =\Psi_1\left( \alpha, \beta; \gamma, \gamma'; x,y\right),
$

$
{z_2} = {x^{1 - \gamma}}\Psi_1\left( 1-\gamma+\alpha, 1-\gamma+ \beta; 2-\gamma, \gamma'; x,y\right),
$

$
z_3= y^{1-\gamma'}\Psi_1\left( 1-\gamma'+\alpha, \beta; \gamma, 2-\gamma'; x,y\right),
$

$
z_4= x^{1-\gamma}y^{1-\gamma'}\Psi_1\left( 2-\gamma-\gamma'+\alpha, 1-\gamma +\beta; 2-\gamma, 2-\gamma'; x,y\right);
$

\begin{equation} \label{difeq27}
\left.\begin{array}{*{20}c}
 {xr+(\gamma-x)p- yq- \alpha z=0} \hfill \\
 {yt+\left(\gamma'-y\right)q -xp -\alpha z=0} \hfill \\
\end{array}\right\}\Psi_2,\,\,\,\,\,\,\,\,\,\,\,\,\,\,\,\,\,\,\,\,\,\,\,\,\,\,\,\,\,\,\,\,\,\,\,\,\,\,\,\,\,\,\,\,\,\,\,\,\,\,\,\,\,\,\,\,\,\,\,\,\,\,\,\,\,\,\,\,\,\,\,\,\,\,\,\,\,\,\,\,\,\,\,\,\,\,\,\,\,\,\,\,\,\,\,\,\,\,\,\,\,\,\,\,\,\,\,\,\,\,\,\,\,\,\,\,\,\,\,\,\,\,\,\,\,\,\,\,\,\,\,\,\,\,\,\,\,\,\,\,\,\,\,\,\,\,\,\,\,
\end{equation}

particular solutions:

$
{z_1} =\Psi_2\left( \alpha; \gamma, \gamma'; x,y\right),
$

$
{z_2} = {x^{1 - \gamma}}\Psi_2\left( 1-\gamma+\alpha; 2-\gamma, \gamma'; x,y\right),
$

$
z_3= y^{1-\gamma'}\Psi_2\left( 1-\gamma'+\alpha; \gamma, 2-\gamma'; x,y\right),
$

$
z_4= x^{1-\gamma}y^{1-\gamma'}\Psi_2\left( 2-\gamma-\gamma'+\alpha; 2-\gamma, 2-\gamma'; x,y\right);
$

\begin{equation} \label{difeq28}
\left.\begin{array}{*{20}c}
 {x(1-x)r+ys+\left[\gamma-\left(\alpha+\beta+1\right)x\right]p-\alpha \beta z=0} \hfill \\
 {yt+xs+\left(\gamma-y\right)q-\alpha' z=0} \hfill \\
\end{array}\right\}\Xi_1;\,\,\,\,\,\,\,\,\,\,\,\,\,\,\,\,\,\,\,\,\,\,\,\,\,\,\,\,\,\,\,\,\,\,\,\,\,\,\,\,\,\,\,\,\,\,\,\,\,\,\,\,\,\,\,\,\,\,\,\,\,\,\,\,\,\,\,\,\,\,\,\,\,\,\,\,\,\,\,\,\,\,\,\,\,\,\,\,\,\,\,\,\,\,\,\,\,\,
\end{equation}
\begin{equation} \label{difeq29}
\left.\begin{array}{*{20}c}
 {x(1-x)r+ys+\left[\gamma-\left(\alpha+\beta+1\right)x\right]p-\alpha \beta z=0} \hfill \\
 {yt+xs+\gamma q- z=0} \hfill \\
\end{array}\right\}\Xi_2;\,\,\,\,\,\,\,\,\,\,\,\,\,\,\,\,\,\,\,\,\,\,\,\,\,\,\,\,\,\,\,\,\,\,\,\,\,\,\,\,\,\,\,\,\,\,\,\,\,\,\,\,\,\,\,\,\,\,\,\,\,\,\,\,\,\,\,\,\,\,\,\,\,\,\,\,\,\,\,\,\,\,\,\,\,\,\,\,\,\,\,\,\,\,\,\,\,\,
\end{equation}
\begin{equation} \label{difeq30}
\left.\begin{array}{*{20}c}
 {x(1+x)r-(x+1)ys+\left[1-\beta+\left(\alpha+\beta'+1\right)x\right]p-\alpha yq+\alpha \beta' z=0} \hfill \\
 {yt-xs+\left(1-\beta'+y\right)q-xp+ \beta z=0} \hfill \\
\end{array}\right\}\Gamma_1;\,\,\,\,\,\,\,\,\,\,\,\,\,\,\,\,\,\,\,\,\,\,\,\,\,\,\,\,\,\,\,\,\,\,\,\,\,\,\,\,\,\,\,\,\,\,\,\,\,
\end{equation}
\begin{equation} \label{difeq31}
\left.\begin{array}{*{20}c}
 {xr-ys+(1-\beta+x)p- yq+ \beta' z=0} \hfill \\
 {yt-xs+(1-\beta'+y)q -xp+\beta z=0} \hfill \\
\end{array}\right\}\Gamma_2;\,\,\,\,\,\,\,\,\,\,\,\,\,\,\,\,\,\,\,\,\,\,\,\,\,\,\,\,\,\,\,\,\,\,\,\,\,\,\,\,\,\,\,\,\,\,\,\,\,\,\,\,\,\,\,\,\,\,\,\,\,\,\,\,\,\,\,\,\,\,\,\,\,\,\,\,\,\,\,\,\,\,\,\,\,\,\,\,\,\,\,\,\,\,\,\,\,\,\,\,\,\,\,\,\,\,\,\,\,\,\,\,\,\,\,\,\,\,\,\,\,\,\,\,\,\,\,
\end{equation}
\begin{equation} \label{difeq32}
\left.\begin{array}{*{20}c}
 {x(1-x)r+y^2t+\left[\delta-\left(\alpha+\beta+1\right)x\right]p+(\beta-\alpha+1) yq-\alpha \beta z=0} \hfill \\
 {yt-xs+\left[1-\alpha+y\right]q+xp+ \beta z=0} \hfill \\
\end{array}\right\}\rm{H}_1,\,\,\,\,\,\,\,\,\,\,\,\,\,\,\,\,\,\,\,\,\,\,\,\,\,\,\,\,\,\,\,\,\,\,\,\,\,\,\,\,\,\,\,\,\,\,\,\,\,\,\,\,\,\,
\end{equation}

particular solutions:

$
{z_1} = {\rm{H}_1}\left( \alpha, \beta; \delta; x,y\right),
$

$
{z_2} = {x^{1 - \delta}}{\rm{H}_1}\left( 1-\delta+\alpha, 1-\delta+\beta; 2-\delta; x,y\right);
$

\begin{equation} \label{difeq33}
\left.\begin{array}{*{20}c}
 {x(1-x)r+xys+\left[\delta-\left(\alpha+\beta+1\right)x\right]p+\beta yq-\alpha \beta z=0} \hfill \\
 {yt-xs+\left[1-\alpha+y\right]q+ \gamma z=0} \hfill \\
\end{array}\right\}\rm{H}_2,\,\,\,\,\,\,\,\,\,\,\,\,\,\,\,\,\,\,\,\,\,\,\,\,\,\,\,\,\,\,\,\,\,\,\,\,\,\,\,\,\,\,\,\,\,\,\,\,\,\,\,\,\,\,\,\,\,\,\,\,\,\,\,\,\,\,\,\,\,\,\,\,\,\,\,\,\,\,\,
\end{equation}

particular solutions:

$
{z_1} = {\rm{H}_2}\left( \alpha, \beta, \gamma; \delta; x,y\right),
$

$
{z_2} = x^{1 - \delta}{\rm{H}_2}\left( 1-\delta+\alpha, 1-\delta+\beta, \gamma; 2-\delta; x,y\right);
$

\begin{equation} \label{difeq34}
\left.\begin{array}{*{20}c}
 {x(1-x)r+xys+\left[\delta-\left(\alpha+\beta+1\right)x\right]p+\beta yq-\alpha \beta z=0} \hfill \\
 {yt-xs+(1-\alpha)q+ z=0} \hfill \\
\end{array}\right\}\rm{H}_3,\,\,\,\,\,\,\,\,\,\,\,\,\,\,\,\,\,\,\,\,\,\,\,\,\,\,\,\,\,\,\,\,\,\,\,\,\,\,\,\,\,\,\,\,\,\,\,\,\,\,\,\,\,\,\,\,\,\,\,\,\,\,\,\,\,\,\,\,\,\,\,\,\,\,\,\,\,\,
\end{equation}

particular solutions:

$
{z_1} = {\rm{H}_3}\left( \alpha, \beta; \delta; x,y\right),
$

$
{z_2} = {x^{1 - \delta}}{\rm{H}_3}\left( 1-\delta+\alpha, 1-\delta+\beta; 2-\delta; x,y\right);
$

\begin{equation} \label{difeq35}
\left.\begin{array}{*{20}c}
 {xr+(\delta-x)p+yq- \alpha z=0} \hfill \\
 {yt-xs+(1-\alpha+y)q +\gamma z=0} \hfill \\
\end{array}\right\}\rm{H}_4,
\,\,\,\,\,\,\,\,\,\,\,\,\,\,\,\,\,\,\,\,\,\,\,\,\,\,\,\,\,\,\,\,\,\,\,\,\,\,\,\,\,\,\,\,\,\,\,\,\,\,\,\,\,\,\,\,\,\,\,\,\,\,\,\,\,\,\,\,\,\,\,\,\,\,\,\,\,\,\,\,\,\,\,\,\,\,\,\,\,\,\,\,\,\,\,\,\,\,\,\,\,\,\,\,\,\,\,\,\,\,\,\,\,\,\,\,\,\,\,\,\,\,\,\,\,\,\,\,\,\,\,\,\,\,\,\,\,\,\,\,\,\,
\end{equation}

particular solutions:

$
{z_1} = {\rm{H}_4}\left( \alpha,  \gamma; \delta; x,y\right),
$

$
{z_2} = x^{1 - \delta}{\rm{H}_4}\left( 1-\delta+\alpha,  \gamma; 2-\delta; x,y\right);
$

\begin{equation} \label{difeq36}
\left.\begin{array}{*{20}c}
 {xr+(\delta-x)p+yq- \alpha z=0} \hfill \\
 {yt-xs+(1-\alpha)q + z=0} \hfill \\
\end{array}\right\}\rm{H}_5;
\,\,\,\,\,\,\,\,\,\,\,\,\,\,\,\,\,\,\,\,\,\,\,\,\,\,\,\,\,\,\,\,\,\,\,\,\,\,\,\,\,\,\,\,\,\,\,\,\,\,\,\,\,\,\,\,\,\,\,\,\,\,\,\,\,\,\,\,\,\,\,\,\,\,\,\,\,\,\,\,\,\,\,\,\,\,\,\,\,\,\,\,\,\,\,\,\,\,\,\,\,\,\,\,\,\,\,\,\,\,\,\,\,\,\,\,\,\,\,\,\,\,\,\,\,\,\,\,\,\,\,\,\,\,\,\,\,\,\,\,\,\,\,\,\,\,\,\,\,\,\,\,
\end{equation}
\begin{equation} \label{difeq37}
\left.\begin{array}{*{20}c}
 {x(1-4x)r+(1-4x)ys-y^2t+\left[\gamma-2(2\alpha+3)x\right]p-2(\alpha+1) y q-\alpha(\alpha+1) z=0} \hfill \\
 {yt+xs+(\gamma-y)q-2 xp - \alpha z=0} \hfill \\
\end{array}\right\}\rm{H}_6,\,\,\,\,\,\,\,\,\,\,\,\,\,
\,\,
\end{equation}
\begin{equation} \label{difeq38}
\left.\begin{array}{*{20}c}
 {x(1-4x)r-4xys-y^2t+\left[\gamma-2(2\alpha+3)x\right]p-2(\alpha+1) y q-\alpha(\alpha+1) z=0} \hfill \\
 {yt+(\delta-y)q-2 xp - \alpha z=0} \hfill \\
\end{array}\right\}\rm{H}_7,\,\,\,\,\,\,\,\,\,\,\,\,\,\,\,\,\,\,\,\,\,\,\,\,\,
\,\,\,\,
\,
\end{equation}

particular solutions:

$
{z_1} = \rm{H}_7\left( \alpha; \gamma, \delta; x,y\right),
$

$
{z_2} = {x^{1 - \gamma}}\rm{H}_7\left( 2-2\gamma+\alpha; 2-\gamma, \delta; x,y\right),
$

$
z_3= y^{1-\delta}\rm{H}_7\left( 1-\delta+\alpha; \gamma, 2-\delta; x,y\right),
$

$
z_4= x^{1-\gamma}y^{1-\delta}\rm{H}_7\left( 3-2\gamma-\delta+\alpha; 2-\gamma, 2-\delta; x,y\right);
$

\begin{equation} \label{difeq39}
\left.\begin{array}{*{20}c}
 {x(1+4x)r-(1+4x)ys+y^2t+\left[1-\beta+2(2\alpha+3)x\right]p-2\alpha y q+\alpha(\alpha+1) z=0} \hfill \\
 {yt-2xs+(1-\alpha+y)q- xp+\beta z=0} \hfill \\
\end{array}\right\}\rm{H}_8;\,\,\,\,\,\,\,\,\,\,\,\,\,\,\,
\,\,\,\,\,\,
\end{equation}
\begin{equation} \label{difeq40}
\left.\begin{array}{*{20}c}
 {x(1-4x)r+4xys-y^2t+\left[\delta-2(2\alpha+3)x\right]p+2\alpha y q-\alpha(\alpha+1) z=0} \hfill \\
 {yt-2xs+(1-\alpha+y)q +\beta z=0} \hfill \\
\end{array}\right\}\rm{H}_9,\,\,\,\,\,\,\,\,\,\,\,\,
\,\,\,\,\,\,\,\,\,\,\,\,\,\,\,\,\,\,\,
\,\,\,\,\,\,\,\,\,\,\,\,\,\,\,
\,
\end{equation}

particular solutions:

$
{z_1} =\rm{H}_9\left( \alpha, \beta; \delta;  x,y\right),
$

$
{z_2} = {x^{1 - \delta}}\rm{H}_9\left( 2-2\delta+\alpha, \beta; 2-\delta; x,y\right),
$

\begin{equation} \label{difeq41}
\left.\begin{array}{*{20}l}
 {x(1-4x)r+4xys-y^2t+\left[\delta-2(2\alpha+3)x\right]p+2\alpha y q-\alpha(\alpha+1) z=0} \hfill \\
 {yt-2xs+(1-\alpha)q + z=0} \hfill \\
\end{array}\right\}\rm{H}_{10},\,\,\,\,\,\,\,\,\,\,\,\,\,\,\,
\,\,\,\,\,\,\,\,\,\,\,\,\,\,\,\,\,\,\,
\,\,\,\,\,\,\,\,
\,
\end{equation}

particular solutions:

$
{z_1} =\rm{H}_{10}\left( \alpha; \delta;  x,y\right),
$

$
{z_2} = {x^{1 - \delta}}\rm{H}_{10}\left( 2-2\delta+\alpha; 2-\delta; x,y\right),
$

\begin{equation} \label{difeq42}
\left.\begin{array}{*{20}c}
 {xr+\left(\delta-x\right)p+ yq-\alpha z=0} \hfill \\
 {y(1+y)t-xs+\left[1-\alpha+(\beta+\gamma+1)y\right]q+ \beta \gamma z=0} \hfill \\
\end{array}\right\}\rm{H}_{11},\,\,\,\,\,\,\,\,\,\,\,\,\,\,\,
\,\,\,\,\,\,\,\,\,\,\,\,\,\,\,\,\,\,\,
\,\,\,\,\,\,\,\,\,\,\,\,\,\,\,\,\,\,\,
\,\,\,\,\,\,\,\,\,\,\,\,\,\,\,\,\,\,\,
\,\,\,\,\,\,\,\,\,\,\,\,\,\,\,\,
\end{equation}

particular solutions:

$
{z_1} = \rm{H}_{11}\left( \alpha, \beta,  \gamma; \delta;  x,y\right),
$

$
{z_2} = {x^{1 - \delta}}\rm{H}_{11}\left( 1-\delta+\alpha, \beta, \gamma; 2-\delta; x,y\right).
$

Note that the systems (52) -- (85), except for the systems (59), (62), (63), (65), (73), (76), (77), (79), (81),     are borrowed from [5], and the systems, related to the hypergeometric functions $H_1$,\, $H_4$, \,$H_5$,\, $H_7$,\, $\Gamma_1$,\, ${\rm{H}}_2$,\, ${\rm{H}}_3$,\,${\rm{H}}_5$, and ${\rm{H}}_7$  are compiled in [6].

\section{General definition of the triple hypergeometric series} \label{S5}

Following Horn [18], we define a hypergeometric series in three variables: the triple power series
\begin{equation} \label{def33}
\sum\limits_{m,n,p=0}^\infty A(m,n,p)x^my^nz^p
\end{equation}
is a hypergeometric series, if the three quotients
\begin{equation}
\label{otn}
\frac{A(m+1,n,p)}{A(m,n,p)}=f(m,n,p),\,\,\,\,\,
\frac{A(m,n+1,p)}{A(m,n,p)}=g(m,n,p),\,\,\,\,\,
\frac{A(m,n,p+1)}{A(m,n,p)}=h(m,n,p)
\end{equation}
are rational functions of $m$, $n$ and $p$.

We put
\begin{equation}
\label{otn53}
f(m,n,p)=\frac{F(m,n,p)}{F'(m,n,p)},\,\,\,\,\,g(m,n,p)=\frac{G(m,n,p)}{G'(m,n,p)},\,\,\,\,\,h(m,n,p)=\frac{H(m,n,p)}{H'(m,n,p)},
\end{equation}
where $F$, $F'$, $G$, $G'$, $H$, $H'$ are polynomials in   $m$, $n$,  $p$ , of respective degrees   $p$, $p'$, $q$, $q'$, $h$, $h'$.   $F'$ is assumed to have a factor  $m+1$, and $G'$  a factor  $n+1$; $F$    and  $F'$  have no common factor except, possibly,  $m+1$;    $G$   and $G'$ no common factor except $n+1$, and  $H$ and $H'$  no common factor except possibly  $p+1$. The highest of the six numbers    $p$, $p'$, $q$, $q'$, $h$, $h'$, is the \textit{order} of the hypergeometric series    (96). H.M.Srivastava  and  Per W. Karlsson [39] investigated hypergeometric series of order two and found that, apart from certain series which are either expressible in terms of one variable, there are essentially 205 distinct complete  convergent series of order two for which
$p=p'=q=q'=h=h'=2$. A. Hasanov and M. Ruzhanskii constructed integral representations [16]  for 205 hypergeometric functions, and also composed systems of partial differential equations [17] satisfied by 205 complete hypergeometric functions of three variables, and wrote out explicit linearly independent solutions of these systems at the origin, if such solutions exist. In the paper [7], 395 \textit{confluent} series were defined, which are limit forms for complete series and for which $p \leq p'=2$, $q \leq q' =2$, $h \leq h' =2$, and $p$, $q$ and $h$ cannot simultaneously be equal to two.

\section{Systems of partial differential equations, associated with the complete hypergeometric functions in three variables} \label{S7}

In this section regions of convergence in the absolute octant are given for the triple Gaussian series as listed in [39] without any permutations of variables. For brevity, the co-ordinates are written $(r,s,t)$  instead of $(|x|, |y|, |z|)$, and the contracted notation for point sets introduced in Section 5. The symbols
\[
\Phi_1, \Phi_2, \Psi_1, \Psi_2, \Theta_1, \Theta_2; \,\,A, B, a, b, \alpha, \beta
\]
denote the auxiliary functions and associated square roots defined by (56) to (61).  Other auxiliary variables and functions are explained in conjunction with the series to which they apply. Occasionally, the definitions involve  double signs ($\pm$ and $\mp$); in such cases it will be understood that all upper signs are taken together and all lower signs are taken together, such that we then have exactly two expressions. The regions of convergence of the double and triple hypergeometric series  are discussed in detail in [39].

Below, we present systems of partial differential equations that are satisfied by the hypergeometric functions defined in  [39].
\\
The series
$$
\sum\limits_{m,n,p = 0}^\infty  {{A_{m,n,p}}{x^m}{y^n}{z^p}},
$$
where
$$
\frac{{{A_{m + 1,n,p}}}}{{{A_{m,n,p}}}} = \frac{{F\left( {m,n,p} \right)}}{{F'\left( {m,n,p} \right)}},\,\,\,\frac{{{A_{m,n + 1,p}}}}{{{A_{m,n,p}}}} = \frac{{G\left( {m,n,p} \right)}}{{G'\left( {m,n,p} \right)}},\,\,\,\frac{{{A_{m,n,p + 1}}}}{{{A_{m,n,p}}}} = \frac{{H\left( {m,n,p} \right)}}{{H'\left( {m,n,p} \right)}}
$$
and
$F\left( {m,n,p} \right), F'\left( {m,n,p} \right), G\left( {m,n,p} \right), G'\left( {m,n,p} \right), H\left( {m,n,p} \right), H'\left( {m,n,p} \right)$ are polynomials as in the previous section, satisfies a system of linear partial differential equations which can be written in terms of differential operators
\begin{equation} \label{oper}
 {\delta _x} = x\frac{\partial }{{\partial x}},\,\,\,{\delta _y} = y\frac{\partial }{{\partial y}},\,\,\,{\delta _z} = z\frac{\partial }{{\partial z}}
 \end{equation}
as
\begin{equation}  \label{system}
\left\{ {\begin{array}{*{20}{l}}
  {\left[ {F'\left( {{\delta _x},{\delta _y},{\delta _z}} \right){x^{ - 1}} - F\left( {{\delta _x},{\delta _y},{\delta _z}} \right)} \right]u = 0,} \\
  {\left[ {G'\left( {{\delta _x},{\delta _y},{\delta _z}} \right){y^{ - 1}} - G\left( {{\delta _x},{\delta _y},{\delta _z}} \right)} \right]u = 0,} \\
  {\left[ {H'\left( {{\delta _x},{\delta _y},{\delta _z}} \right){z^{ - 1}} - H\left( {{\delta _x},{\delta _y},{\delta _z}} \right)} \right]u = 0.}
\end{array}} \right.
\end{equation}

In what follows we shall restrict ourselves to hypergeometric functions of the second order, in which case we find three partial differential equations of second order. The three equations are certainly compatible (since the hypergeometric series satisfies three them), and from the general theory of such systems it follows that they have at most eight, and possibly less, linearly independent solutions in common.

The difficulties in dealing with these systems of partial differential equations have two sources. One is the unsatisfactory state of the general analytic theory of systems of partial differential equations; in particular our very scant knowledge of the behavior of solutions in the neighborhood of points at which more than two singular curves of the system intersect, or at which two singular curves are at contact. The second difficulty is the large number of apparently distinct systems.

For an example of constructing a system of equations, see Section 9.

\bigskip 

In this Section we consider a complete hypergeometric functions in three variables.

\begin{equation} \label{eqf1}
  {F_{1b}}\left( {{a_1},{a_2},{a_3},{a_4},{b_1},{b_2};c;x,y,z} \right)
  =  \sum\limits_{m,n,p = 0}^\infty  {} \frac{{{{\left( {{a_1}} \right)}_n}{{\left( {{a_2}} \right)}_n}{{\left( {{a_3}} \right)}_p}{{\left( {{a_4}} \right)}_p}{{\left( {{b_1}} \right)}_{m - n}}{{\left( {{b_2}} \right)}_{m - p}}}}{{{{\left( c \right)}_m}}}\frac{x^m}{m!}\frac{y^n}{n!}\frac{z^p}{p!},
  \end{equation}

region of convergence:
$$
\left\{{s < 1,\,\,\,t < 1,\,\,\,r < \min \left\{ {1,\frac{{1 - s}}{s},\frac{{1 - t}}{t},\frac{{\left( {1 - s} \right)\left( {1 - t} \right)}}{{st}}} \right\}}\right\}.
$$

System of partial differential equations:

$\left\{
{\begin{array}{*{20}{l}}
  x\left( {1 - x} \right){u_{xx}} + xy{u_{xy}} + xz{u_{xz}} - yz{u_{yz}}  + \left[ {c - \left( {{b_1} + {b_2} + 1} \right)x} \right]{u_x} + {b_2}y{u_y} + {b_1}z{u_z} - {b_1}{b_2}u = 0,  \\
  {y\left( {1 + y} \right){u_{yy}} - x{u_{xy}} + \left[ {1 - {b_1} + \left( {{a_1} + {a_2} + 1} \right)y} \right]{u_y} + {a_1}{a_2}u = 0,} \\
  {z\left( {1 + z} \right){u_{zz}} - x{u_{xz}} + \left[ {1 - {b_2} + \left( {{a_3} + {a_4} + 1} \right)z} \right]{u_z} + {a_3}{a_4}u = 0,}
\end{array}} \right.
$

where $u\equiv {F_{1b}}\left( {{a_1},{a_2},{a_3},{a_4},{b_1},{b_2};c;x,y,z} \right)$.

Particular solutions:

$
{u_1} = {F_{1b}}\left( {{a_1},{a_2},{a_3},{a_4},{b_1},{b_2};c;x,y,z} \right),
$

$
{u_2} = {x^{1 - c}}{F_{1b}}\left( {{a_1},{a_2},{a_3},{a_4},1 - c + {b_1},1 - c + {b_2};2 - c;x,y,z} \right).
$

\bigskip 

\begin{equation} \label{eqf1c}
{F_{1c}}\left( {{a_1},{a_2},{a_3},{a_4},{b_1},{b_2};c;x,y,z} \right)\hfill \\
  = \sum\limits_{m,n,p = 0}^\infty  {} \frac{{{{\left( {{a_1}} \right)}_m}{{\left( {{a_2}} \right)}_n}{{\left( {{a_3}} \right)}_p}{{\left( {{a_4}} \right)}_p}{{\left( {{b_1}} \right)}_{m - n}}{{\left( {{b_2}} \right)}_{n - p}}}}{{{{\left( c \right)}_m}}}\frac{x^m}{m!}\frac{y^n}{n!}\frac{z^p}{p!},
\end{equation}

region of convergence:
$$
\left\{r < 1,\,\,\,s < \frac{1}{{1 + r}},\,\,\,t < \frac{1}{1 + s + rs}\right\}.
$$

System of partial differential equations:

$
\left\{
{\begin{array}{*{20}{l}}
  {x\left( {1 - x} \right){u_{xx}} + xy{u_{xy}} + \left[ {c - \left( {{a_1} + {b_1} + 1} \right)x} \right]{u_x} + {a_1}y{u_y} - {a_1}{b_1}u = 0,} \\
  y\left( {1 + y} \right){u_{yy}} - x{u_{xy}} - yz{u_{yz}} + \left[ {1 - {b_1} + \left( {{a_2} + {b_2} + 1} \right)y} \right]{u_y} - {a_2}z{u_z} + {a_2}{b_2}u = 0,\\
  {z\left( {1 + z} \right){u_{zz}} - y{u_{yz}} + \left[ {1 - {b_2} + \left( {{a_3} + {a_4} + 1} \right)z} \right]{u_z} + {a_3}{a_4}u = 0,}
\end{array}} \right.
$

where $u\equiv {F_{1c}}\left( {{a_1},{a_2},{a_3},{a_4},{b_1},{b_2};c;x,y,z} \right).   $

Particular solutions:

$
{u_1} = {F_{1c}}\left({{a_1},{a_2},{a_3},{a_4},{b_1},{b_2};c;x,y,z}\right),
$

$
{u_2} = {x^{1 - c}}{F_{1c}}\left( {1 - c + {a_1},{a_2},{a_3},{a_4},1 - c + {b_1},{b_2};2 - c;x,y,z} \right).
$

\bigskip

\begin{equation} \label{eqf1c}
 {F_{1d}}\left( {{a_1},{a_2},{a_3},{b_1},{b_2},{b_3};x,y,z} \right) \hfill \\
  = \sum\limits_{m,n,p = 0}^\infty  {} {{{{\left( {{a_1}} \right)}_m}{{\left( {{a_2}} \right)}_n}{{\left( {{a_3}} \right)}_p}{{\left( {{b_1}} \right)}_{m - n}}{{\left( {{b_2}} \right)}_{n - p}}{{\left( {{b_3}} \right)}_{p - m}}}}\frac{x^m}{m!}\frac{y^n}{n!}\frac{z^p}{p!},
\end{equation}

region of convergence:
$$\left\{ t\left( {1 + s} \right) < 1,\,\,\,r\left( {1 + t} \right) < 1,\,\,\,s\left( {1 + r} \right) < 1\right\}.
$$

System of partial differential equations:

$ \left\{ {\begin{array}{*{20}{l}}
  x\left( {1 + x} \right){u_{xx}} - xy{u_{xy}} - z{u_{xz}} + \left[ {1 - {b_3} + \left( {{a_1} + {b_1} + 1} \right)x} \right]{u_x}
   - {a_1}y{u_y} + {a_1}{b_1}u = 0, \\
  y\left( {1 + y} \right){u_{yy}} - x{u_{xy}} - yz{u_{yz}} + \left[ {1 - {b_1} + \left( {{a_2} + {b_2} + 1} \right)y} \right]{u_y}
   - {a_2}z{u_z} + {a_2}{b_2}u = 0,  \\
  z\left( {1 + z} \right){u_{zz}} - xz{u_{xz}} - y{u_{yz}} + \left[ {1 - {b_2} + \left( {{a_3} + {b_3} + 1} \right)z} \right]{u_z}
   - {a_3}x{u_x} + {a_3}{b_3}u = 0,
\end{array}} \right.
$

where $u\equiv  {F_{1d}}\left( {{a_1},{a_2},{a_3},{b_1},{b_2},{b_3};x,y,z} \right)$.

\bigskip 

\begin{equation} \label{eqf1c}
{F_{1e}}\left( {{a_1},{a_2},{a_3},{b_1},{b_2},{b_3};x,y,z} \right) \hfill \\
  = \sum\limits_{m,n,p = 0}^\infty  {} {{{{\left( {{a_1}} \right)}_m}{{\left( {{a_2}} \right)}_p}{{\left( {{a_3}} \right)}_p}{{\left( {{b_1}} \right)}_{m - n}}{{\left( {{b_2}} \right)}_{n - p}}{{\left( {{b_3}} \right)}_{n - m}}}}\frac{x^m}{m!}\frac{y^n}{n!}\frac{z^p}{p!},
\end{equation}

region of convergence:
$$\left\{
{r < 1,\,\,\,s < 1,\,\,\,t < \frac{1}{{1 + s}}}\right\}.
$$

System of partial differential equations:

$
\left\{ {\begin{array}{*{20}{l}}
  x(1 + x){u_{xx}} - ( {1 + x})y{u_{xy}} + \left[ {1 - {b_3} + \left( {{a_1} + {b_1} + 1} \right)x} \right]{u_x}
   - {a_1}y{u_y} + {a_1}{b_1}u = 0, \\
  y\left( {1 + y} \right){u_{yy}} - x\left( {1 + y} \right){u_{xy}} + xz{u_{xz}} - yz{u_{yz}} \\\,\,\,\,\,\,\,\,\,\,\,- {b_2}x{u_x}
   + \left[ {1 - {b_1} + \left( {{b_2} + {b_3} + 1} \right)y} \right]{u_y} - {b_3}z{u_z} + {b_2}{b_3}u = 0, \\
  {z\left( {1 + z} \right){u_{zz}} - y{u_{yz}} + \left[ {1 - {b_2} + \left( {{a_2} + {a_3} + 1} \right)z} \right]{u_z} + {a_2}{a_3}u = 0,}
\end{array}} \right.
$

where $u\equiv  {F_{1e}}\left( {{a_1},{a_2},{a_3},{b_1},{b_2},{b_3};x,y,z} \right)$.

\bigskip

\begin{equation} \label{eqf1c}
{F_{2b}}\left( {{a_1},{a_2},{a_3},{a_4},{a_5},b;c;x,y,z} \right) \hfill \\
  = \sum\limits_{m,n,p = 0}^\infty  {} \frac{{{{\left( {{a_1}} \right)}_m}{{\left( {{a_2}} \right)}_n}{{\left( {{a_3}} \right)}_n}{{\left( {{a_4}} \right)}_p}{{\left( {{a_5}} \right)}_p}{{\left( b \right)}_{m - p}}}}{{{{\left( c \right)}_{m + n}}}}\frac{x^m}{m!}\frac{y^n}{n!}\frac{z^p}{p!},
\end{equation}

region of convergence:
$$\left\{
{r < 1,\,\,\,s < 1,\,\,\,t < \frac{1}{{1 + r}}}\right\}.
$$

System of partial differential equations:

$
\left\{ {\begin{array}{*{20}{l}}
  x\left( {1 - x} \right){u_{xx}} + y{u_{xy}} + xz{u_{xz}} + \left[ {c - \left( {{a_1} + b + 1} \right)x} \right]{u_x}  + {a_1}z{u_z} - {a_1}bu = 0, \\
  y\left( {1 - y} \right){u_{yy}} + x{u_{xy}} + \left[ {c - \left( {{a_2} + {a_3} + 1} \right)y} \right]{u_y}
   - {a_2}{a_3}u = 0,  \\
  z\left( {1 + z} \right){u_{zz}} - x{u_{xz}} + \left[ {1 - b + \left( {{a_4} + {a_5} + 1} \right)z} \right]{u_z}
   + {a_4}{a_5}u = 0,
\end{array}} \right.
$

where $u\equiv  {F_{2b}}\left( {{a_1},{a_2},{a_3},{a_4},{a_5},b;c;x,y,z} \right)$.

\bigskip

\begin{equation} \label{eqf1c}
{F_{2c}}\left( {{a_1},{a_2},{a_3},{a_4},{a_5},b;c;x,y,z} \right) \hfill \\
  = \sum\limits_{m,n,p = 0}^\infty  {} \frac{{{{\left( {{a_1}} \right)}_m}{{\left( {{a_2}} \right)}_m}{{\left( {{a_3}} \right)}_n}{{\left( {{a_4}} \right)}_n}{{\left( {{a_5}} \right)}_p}{{\left( b \right)}_{p - m - n}}}}{{{{\left( c \right)}_p}}}\frac{x^m}{m!}\frac{y^n}{n!}\frac{z^p}{p!}, \end{equation}

first appearance of this function in the literature, and old notation:  [4], \,$H_{3,1}$,

region of convergence:
$$\left\{
{t < 1,\,\,\,\max \left[ {r,s} \right] < \frac{1}{{1 + t}}}\right\}.
$$

System of partial differential equations:

$
\left\{
{\begin{array}{*{20}{l}}
  x\left( {1 + x} \right){u_{xx}} + y{u_{xy}} - z{u_{xz}} + \left[ {1 - b + \left( {{a_1} + {a_2} + 1} \right)x} \right]{u_x}
   + {a_1}{a_2}u = 0, \\
  y\left( {1 + y} \right){u_{yy}} - z{u_{yz}} + x{u_{xy}} + \left[ {1 - b + \left( {{a_3} + {a_4} + 1} \right)y} \right]{u_y}
   + {a_3}{a_4}u = 0, \\
  z\left( {1 - z} \right){u_{zz}} + xz{u_{xz}} + yz{u_{yz}} + \left[ {c - \left( {{a_5} + b + 1} \right)z} \right]{u_z}
   + {a_5}x{u_x} + {a_5}y{u_y} - {a_5}bu = 0,
\end{array}} \right.
$

where $u\equiv    {F_{2c}}\left( {{a_1},{a_2},{a_3},{a_4},{a_5},b;c;x,y,z} \right)$.

Particular solutions:

$
{u_1} = {F_{2c}}\left( {{a_1},{a_2},{a_3},{a_4},{a_5},b;c;x,y,z} \right),
$

$
{u_2} = {z^{1 - c}}{F_{2c}}\left( {{a_1},{a_2},{a_3},{a_4},1 - c + {a_5},1 - c + b;2 - c;x,y,z} \right).
$

\bigskip 

\begin{equation} \label{eqf1c}
{F_{2d}}\left( {{a_1},{a_2},{a_3},{a_4},{b_1},{b_2};x,y,z} \right) \hfill \\
  = \sum\limits_{m,n,p = 0}^\infty  {} {{{{\left( {{a_1}} \right)}_m}{{\left( {{a_2}} \right)}_n}{{\left( {{a_3}} \right)}_n}{{\left( {{a_4}} \right)}_p}{{\left( {{b_1}} \right)}_{p - m - n}}{{\left( {{b_2}} \right)}_{m - p}}}}\frac{x^m}{m!}\frac{y^n}{n!}\frac{z^p}{p!},
\end{equation}

first appearance of this function in the literature, and old notation  [3], \,$G_{D}$,

region of convergence:
$$\left\{
{r < 1,\,\,\,t < 1,\,\,\,s < \frac{1}{{1 + t}}}\right\}.
$$

System of partial differential equations:

$
\left\{
{\begin{array}{*{20}{l}}
  x\left( {1 + x} \right){u_{xx}} + y{u_{xy}} - \left( {1 + x} \right)z{u_{xz}}
   + \left[ {1 - {b_1} + \left( {{a_1} + {b_2} + 1} \right)x} \right]{u_x} - {a_1}z{u_z} + {a_1}{b_2}u = 0, \\
  y\left( {1 + y} \right){u_{yy}} + x{u_{xy}} - z{u_{yz}}
   + \left[ {1 - {b_1} + \left( {{a_2} + {a_3} + 1} \right)y} \right]{u_y} + {a_2}{a_3}u = 0, \\
  z\left( {1 + z} \right){u_{zz}} - x\left( {1 + z} \right){u_{xz}} - yz{u_{yz}}\\
   \,\,\,\,\,\,\,\,\,\,\,+ \left[ {1 - {b_2} + \left( {{a_4} + {b_1} + 1} \right)z} \right]{u_z} - {a_4}x{u_x} - {a_4}y{u_y} + {a_4}{b_1}u = 0,
\end{array}} \right.
$

where $u\equiv \,\,  {F_{2d}}\left( {{a_1},{a_2},{a_3},{a_4},{b_1},{b_2};x,y,z} \right)$.

\bigskip 

\begin{equation} \label{eqf1c}
{F_{3a}}\left( {{a_1},{a_2},{a_3},{a_4},{a_5},{a_6};c;x,y,z} \right) \hfill \\
  = \sum\limits_{m,n,p = 0}^\infty  {} \frac{{{{\left( {{a_1}} \right)}_m}{{\left( {{a_2}} \right)}_m}{{\left( {{a_3}} \right)}_n}{{\left( {{a_4}} \right)}_n}{{\left( {{a_5}} \right)}_p}{{\left( {{a_6}} \right)}_p}}}{{{{\left( c \right)}_{m + n + p}}}}{x^m}\frac{x^m}{m!}\frac{y^n}{n!}\frac{z^p}{p!},
\end{equation}

first appearance of this function in the literature, and old notation: [24], $F_B^{(3)}$, $ F_2,$

region of convergence:
$$\left\{
{r < 1,\,\,\,s < 1,\,\,\,t < 1}\right\}.
$$

System of partial differential equations:

$
\left\{
{\begin{array}{*{20}{l}}
  {x\left( {1 - x} \right){u_{xx}} + y{u_{xy}} + z{u_{xz}} + \left[ {c - \left( {{a_1} + {a_2} + 1} \right)x} \right]{u_x} - {a_1}{a_2}u = 0,} \\
  {y\left( {1 - y} \right){u_{yy}} + x{u_{xy}} + z{u_{yz}} + \left[ {c - \left( {{a_3} + {a_4} + 1} \right)y} \right]{u_y} - {a_3}{a_4}u = 0,} \\
  {z\left( {1 - z} \right){u_{zz}} + x{u_{xz}} + y{u_{yz}} + \left[ {c - \left( {{a_5} + {a_6} + 1} \right)z} \right]{u_z} - {a_5}{a_6}u = 0,}
\end{array}} \right.
$

where $u\equiv \,\,   {F_{3a}}\left( {{a_1},{a_2},{a_3},{a_4},{a_5},{a_6};c;x,y,z} \right) $.

\bigskip 

\begin{equation} \label{eqf1c}
{F_{4b}}\left( {{a_1},{a_2},{a_3},{a_4},b;{c_1},{c_2};x,y,z} \right) \hfill \\
  = \sum\limits_{m,n,p = 0}^\infty  {} \frac{{{{\left( {{a_1}} \right)}_m}{{\left( {{a_2}} \right)}_n}{{\left( {{a_3}} \right)}_p}{{\left( {{a_4}} \right)}_p}{{\left( b \right)}_{m + n - p}}}}{{{{\left( {{c_1}} \right)}_m}{{\left( {{c_2}} \right)}_n}}}\frac{x^m}{m!}\frac{y^n}{n!}\frac{z^p}{p!},
\end{equation}
first appearance of this function in the literature, and old notation: [4],  $H_{3,2},$

region of convergence:
$$\left\{
{r + s < 1,\,\,\,t < \frac{1}{{1 + r + s}}}\right\} .
$$

System of partial differential equations:

$
\left\{
{\begin{array}{*{20}{l}}
  x\left( {1 - x} \right){u_{xx}} - xy{u_{xy}} + xz{u_{xz}} + \left[ {{c_1} - \left( {{a_1} + b + 1} \right)x} \right]{u_x}
   - {a_1}y{u_y} + {a_1}z{u_z} - {a_1}bu = 0, \\
  y\left( {1 - y} \right){u_{yy}} - xy{u_{xy}} + yz{u_{yz}} + \left[ {{c_2} - \left( {{a_2} + b + 1} \right)y} \right]{u_y}
  - {a_2}x{u_x} + {a_2}z{u_z} - {a_2}bu = 0, \\
  z\left( {1 + z} \right){u_{zz}} - x{u_{xz}} - y{u_{yz}} + \left[ {1 - b + \left( {{a_3} + {a_4} + 1} \right)z} \right]{u_z}
  + {a_3}{a_4}u = 0,
\end{array}} \right.
$

where $u\equiv \,\,   {F_{4b}}\left( {{a_1},{a_2},{a_3},{a_4},b;{c_1},{c_2};x,y,z} \right)$.

Particular solutions:

$
{u_1} = {F_{4b}}\left( {{a_1},{a_2},{a_3},{a_4},b;{c_1},{c_2};x,y,z} \right),
$

$
{u_2} = {x^{1 - {c_1}}}{F_{4b}}\left( {1 - {c_1} + {a_1},{a_2},{a_3},{a_4},1 - {c_1} + b;2 - {c_1},{c_2};x,y,z} \right),
$

$
{u_3} = {y^{1 - {c_2}}}{F_{4b}}\left( {{a_1},1 - {c_2} + {a_2},{a_3},{a_4},b + 1 - {c_2};{c_1},2 - {c_2};x,y,z} \right),
$

$
  {u_4} = {x^{1 - {c_1}}}{y^{1 - {c_2}}}{F_{4b}}\left( {1 - {c_1} + {a_1},1 - {c_2} + {a_2},{a_3},{a_4},2 + b - {c_1} - {c_2};2 - {c_1},2 - {c_2};x,y,z} \right).
$

\bigskip 

\begin{equation} \label{eqf1c}
 {F_{4c}}\left( {{a_1},{a_2},{a_3},{a_4},b;{c_1},{c_2};x,y,z} \right) \hfill \\
  = \sum\limits_{m,n,p = 0}^\infty  {} \frac{{{{\left( {{a_1}} \right)}_{m + n}}{{\left( {{a_2}} \right)}_n}{{\left( {{a_3}} \right)}_p}{{\left( {{a_4}} \right)}_p}{{\left( b \right)}_{m - p}}}}{{{{\left( {{c_1}} \right)}_m}{{\left( {{c_2}} \right)}_n}}}\frac{x^m}{m!}\frac{y^n}{n!}\frac{z^p}{p!},
\end{equation}

first appearance of this function in the literature: [34],

region of convergence:
$$\left\{
{r + s < 1,\,\,\,t < \frac{{1 - s}}{{1 + r - s}}}\right\}.
$$

System of partial differential equations:

$
\left\{
{\begin{array}{*{20}{l}}
  x\left( {1 - x} \right){u_{xx}} - xy{u_{xy}} + xz{u_{xz}} + yz{u_{yz}}
   + \left[ {{c_1} - \left( {{a_1} + b + 1} \right)x} \right]{u_x} \hfill
   - by{u_y} + {a_1}z{u_z} - {a_1}bu = 0, \\
  y\left( {1 - y} \right){u_{yy}} - xy{u_{xy}} + \left[ {{c_2} - \left( {{a_1} + {a_2} + 1} \right)y} \right]{u_y} - {a_2}x{u_x} - {a_1}{a_2}u = 0 , \\
  {z\left( {1 + z} \right){u_{zz}} - x{u_{xz}} + \left[ {1 - b + \left( {{a_3} + {a_4} + 1} \right)z} \right]{u_z} + {a_3}{a_4}u = 0,}
\end{array}} \right.
$

where $u\equiv \,\,   {F_{4c}}\left( {{a_1},{a_2},{a_3},{a_4},b;{c_1},{c_2};x,y,z} \right)$.

Particular solutions:

$
{u_1} = {F_{4c}}\left( {{a_1},{a_2},{a_3},{a_4},b;{c_1},{c_2};x,y,z} \right),
$

$
{u_2} = {x^{1 - {c_1}}}{F_{4c}}\left( {1 - {c_1} + {a_1},{a_2},{a_3},{a_4},1 - {c_1} + b;2 - {c_1},{c_2};x,y,z} \right),
$

$
{u_3} = {y^{1 - {c_2}}}{F_{4c}}\left( {1 - {c_2} + {a_1},1 - {c_2} + {a_2},{a_3},{a_4},b;{c_1},2 - {c_2};x,y,z} \right),
$

$
{u_4} = {x^{1 - {c_1}}}{y^{1 - {c_2}}}{F_{4c}}\left( {2 - {c_1} - {c_2} + {a_1},1 - {c_2} + {a_2},{a_3},{a_4},1 - {c_1} + b;2 - {c_1},2 - {c_2};x,y,z} \right).
$

\bigskip 

\begin{equation} \label{eqf1c}
{F_{4d}}\left( {{a_1},{a_2},{a_3},{a_4},b;{c_1},{c_2};x,y,z} \right) \hfill \\
  = \sum\limits_{m,n,p = 0}^\infty  {} \frac{{{{\left( {{a_1}} \right)}_{m + n}}{{\left( {{a_2}} \right)}_m}{{\left( {{a_3}} \right)}_n}{{\left( {{a_4}} \right)}_p}{{\left( b \right)}_{p - m}}}}{{{{\left( {{c_1}} \right)}_n}{{\left( {{c_2}} \right)}_p}}}\frac{x^m}{m!}\frac{y^n}{n!}\frac{z^p}{p!},
\end{equation}

region of convergence:
$$\left\{
{s < 1,\,\,\,t < 1,\,\,\,r < \frac{{1 - s}}{{1 + t}}}\right\}.
$$

System of partial differential equations:

$
\left\{
{\begin{array}{*{20}{l}}
  x\left( {1 + x} \right){u_{xx}} + xy{u_{xy}} - z{u_{xz}} + \left[ {1 - b + \left( {{a_1} + {a_2} + 1} \right)x} \right]{u_x}
  + {a_2}y{u_y} + {a_1}{a_2}u = 0, \\
  y\left( {1 - y} \right){u_{yy}} - xy{u_{xy}} + \left[ {{c_1} - \left( {{a_1} + {a_3} + 1} \right)y} \right]{u_y}
   - {a_3}x{u_x} - {a_1}{a_3}u = 0, \\
  z\left( {1 - z} \right){u_{zz}} + xz{u_{xz}} + \left[ {{c_2} - \left( {{a_4} + b + 1} \right)z} \right]{u_z}
  + {a_4}x{u_x} - {a_4}bu = 0,
\end{array}} \right.
$

where $u\equiv \,\,   {F_{4d}}\left( {{a_1},{a_2},{a_3},{a_4},b;{c_1},{c_2};x,y,z} \right)$.

Particular solutions:

$
{u_1} = {F_{4d}}\left( {{a_1},{a_2},{a_3},{a_4},b;{c_1},{c_2};x,y,z} \right),
$

$
{u_2} = {y^{1 - {c_1}}}{F_{4d}}\left( {1 - {c_1} + {a_1},{a_2},1 - {c_1} + {a_3},{a_4},b;{2-c_1},{c_2};x,y,z} \right),
$

$
{u_3} = {z^{1 - {c_2}}}{F_{4d}}\left( {{a_1},{a_2},{a_3},1 - {c_2} + {a_4},1 - {c_2} + b;{c_1},2 - {c_2};x,y,z} \right),
$

$
{u_4} = {y^{1 - {c_1}}}{z^{1 - {c_2}}}{F_{4d}}\left( {1 - {c_1} + {a_1},{a_2},1 - {c_1} + {a_3},1 - {c_2} + {a_4},1 - {c_2} + b;2 - {c_1},2 - {c_2};x,y,z} \right).
$

\bigskip 

\begin{equation} \label{eqf1c}
{F_{4e}}\left( {{a_1},{a_2},{a_3},{b_1},{b_2};c;x,y,z} \right)
   = \sum\limits_{m,n,p = 0}^\infty  {} \frac{{{{\left( {{a_1}} \right)}_n}{{\left( {{a_2}} \right)}_p}{{\left( {{a_3}} \right)}_p}{{\left( {{b_1}} \right)}_{m + n - p}}{{\left( {{b_2}} \right)}_{m - n}}}}{{{{\left( c \right)}_m}}}\frac{x^m}{m!}\frac{y^n}{n!}\frac{z^p}{p!},
\end{equation}

region of convergence:
$$\left\{
{r < 1,\,\,\,s + 2\sqrt {rs}  < 1,\,\,\,t < \min \left[ {\frac{1}{{1 + r}},\frac{1}{{1 + s + 2\sqrt {rs} }}} \right]}\right\}.
$$

System of partial differential equations:

$
\left\{
{\begin{array}{*{20}{l}}
  x\left( {1 - x} \right){u_{xx}} + xz{u_{xz}} - yz{u_{yz}} + {y^2}{u_{yy}} \\
 \,\,\,\,\,\,\,\,\,\, + \left[ {c - \left( {{b_1} + {b_2} + 1} \right)x} \right]{u_x} + \left( {{b_1} - {b_2} + 1} \right)y{u_y} + {b_2}z{u_z} - {b_1}{b_2}u = 0,  \\
  y\left( {1 + y} \right){u_{yy}} - x\left( {1 - y} \right){u_{xy}} - zy{u_{yz}}\\
 \,\,\,\,\,\,\,\,\,\, + \left[ {1 - {b_2} + \left( {{a_1} + {b_1} + 1} \right)y} \right]{u_y}
   + {a_1}x{u_x} - {a_1}z{u_z} + {a_1}{b_1}u = 0,  \\
  z\left( {1 + z} \right){u_{zz}} - x{u_{xz}} - y{u_{yz}}
   + \left[ {1 - {b_1} + \left( {{a_2} + {a_3} + 1} \right)z} \right]{u_z} + {a_2}{a_3}u = 0,
\end{array}} \right.
$

where $u\equiv \,\,   {F_{4e}}\left( {{a_1},{a_2},{a_3},{b_1},{b_2};c;x,y,z} \right)$.

Particular solutions:

$
{u_1} = {F_{4e}}\left( {{a_1},{a_2},{a_3},{b_1},{b_2};c;x,y,z} \right),
$

$
{u_2} = {x^{1 - c}}{F_{4e}}\left( {{a_1},{a_2},{a_3},1 - c + {b_1},1 - c + {b_2};2 - c;x,y,z} \right).
$

\bigskip 

\begin{equation} \label{eqf1c}
{F_{4f}}\left( {{a_1},{a_2},{a_3},{b_1},{b_2};c;x,y,z} \right) \hfill \\
  = \sum\limits_{m,n,p = 0}^\infty  {} \frac{{{{\left( {{a_1}} \right)}_m}{{\left( {{a_2}} \right)}_n}{{\left( {{a_3}} \right)}_p}{{\left( {{b_1}} \right)}_{m + n - p}}{{\left( {{b_2}} \right)}_{p - n}}}}{{{{\left( c \right)}_m}}}\frac{x^m}{m!}\frac{y^n}{n!}\frac{z^p}{p!},
\end{equation}

first appearance of this function in the literature, and old notation: [10], $D_{(3)}^{1,2}$,

region of convergence:
$$\left\{
{r + s < 1,\,\,t < \frac{1}{{1 + r}}}\right\}.
$$

System of partial differential equations:

$\left\{
{\begin{array}{*{20}{l}}
  x\left( {1 - x} \right){u_{xx}} - xy{u_{xy}} + xz{u_{xz}} + \left[ {c - \left( {{a_1} + {b_1} + 1} \right)x} \right]{u_x}
   - {a_1}y{u_y} + {a_1}z{u_z} - {a_1}{b_1}u = 0 ,  \\
  y\left( {1 + y} \right){u_{yy}} + xy{u_{xy}} - ({1 + y})z{u_{yz}} + \left[ {1 - {b_2} + \left( {{a_2} + {b_1} + 1} \right)y} \right]{u_y}\\
  \,\,\,\,\,\,\,\,\,\, + {a_2}x{u_x} - {a_2}z{u_z} + {a_2}{b_1}u = 0, \\
  z\left( {1 + z} \right){u_{zz}} - x{u_{xz}} - y({1 + z}){u_{yz}} + \left[ {1 - {b_1} + \left( {{a_3} + {b_2} + 1} \right)z} \right]{u_z}
  - {a_3}y{u_y} + {a_3}{b_2}u = 0,
\end{array}} \right.
$

where $u\equiv \,\,   {F_{4f}}\left( {{a_1},{a_2},{a_3},{b_1},{b_2};c;x,y,z} \right)$.

Particular solutions:

$
{u_1} = {F_{4f}}\left( {{a_1},{a_2},{a_3},{b_1},{b_2};c;x,y,z} \right),
$

$
{u_2} = {x^{1 - c}}{F_{4f}}\left( {1 - c + {a_1},{a_2},{a_3},1 - c + {b_1},{b_2};2-c;x,y,z} \right).
$

\bigskip

\begin{equation} \label{eqf1c}
 {F_{4g}}\left( {{a_1},{a_2},{a_3},{b_1},{b_2};c;x,y,z} \right) \hfill \\
  = \sum\limits_{m,n,p = 0}^\infty  {} \frac{{{{\left( {{a_1}} \right)}_{m + n}}{{\left( {{a_2}} \right)}_p}{{\left( {{a_3}} \right)}_p}{{\left( {{b_1}} \right)}_{m - p}}{{\left( {{b_2}} \right)}_{n - m}}}}{{{{\left( c \right)}_n}}}\frac{x^m}{m!}\frac{y^n}{n!}\frac{z^p}{p!},
\end{equation}

region of convergence:
$$\left\{
{r < 1,\,\,\,t < \frac{1}{{1 + r}},\,\,\,\sqrt s  < \min \left\{ {1,\frac{{1 - r}}{{2\sqrt r }},\frac{1}{2}\sqrt {\frac{{1 - t}}{{rt}}}  - \frac{1}{2}\sqrt {\frac{{rt}}{{1 - t}}} } \right\}}\right\}.
$$

System of partial differential equations:

$\left\{
{\begin{array}{*{20}{l}}
  x( {1 + x}){u_{xx}} - ({1 - x})y{u_{xy}} - xz{u_{xz}} - yz{u_{yz}}\hfill \\
   \,\,\,\,\,\,\,\,\,\,+ \left[ {1 - {b_2} + \left( {{a_1} + {b_1} + 1} \right)x} \right]{u_x}    + {b_1}y{u_y} - {a_1}z{u_z} + {a_1}{b_1}u = 0 ,  \\
  y\left( {1 - y} \right){u_{yy}} + {x^2}{u_{xx}} + \left[ {c - \left( {{a_1} + {b_2} + 1} \right)y} \right]{u_y}
   + \left( {{a_1} - {b_2} + 1} \right)x{u_x} - {a_1}{b_2}u = 0,  \\
  {z\left( {1 + z} \right){u_{zz}} - x{u_{xz}} + \left[ {1 - {b_1} + \left( {{a_2} + {a_3} + 1} \right)z} \right]{u_z} + {a_2}{a_3}u = 0,}
\end{array}} \right.
$

where $u\equiv \,\,   {F_{4g}}\left( {{a_1},{a_2},{a_3},{b_1},{b_2};c;x,y,z} \right)$.

Particular solutions:

$
{u_1} = {F_{4g}}\left( {{a_1},{a_2},{a_3},{b_1},{b_2};c;x,y,z} \right),
$

$
{u_2} = {y^{1 - c}}{F_{4g}}\left( {1 - c + {a_1},{a_2},{a_3},{b_1},1 - c + {b_2};2 - c;x,y,z} \right).
$

\bigskip

\begin{equation} \label{eqf1c}
 {F_{4h}}\left( {{a_1},{a_2},{a_3},{b_1},{b_2};c;x,y,z} \right)
  = \sum\limits_{m,n,p = 0}^\infty  {} \frac{{{{\left( {{a_1}} \right)}_{m + n}}{{\left( {{a_2}} \right)}_n}{{\left( {{a_3}} \right)}_p}{{\left( {{b_1}} \right)}_{m - p}}{{\left( {{b_2}} \right)}_{p - m}}}}{{{{\left( c \right)}_n}}}\frac{x^m}{m!}\frac{y^n}{n!}\frac{z^p}{p!},
\end{equation}

region of convergence:
$$\left\{
{r + s < 1,\,\,\,t < 1}\right\} .
$$

System of partial differential equations:

$\left\{
{\begin{array}{*{20}{l}}
  x( {1 + x}){u_{xx}} + xy{u_{xy}} - ( {1 + x})z{u_{xz}} - yz{u_{yz}} \hfill \\
  \,\,\,\,\,\,\,\,\,\, + \left[ {1 - {b_2} + \left( {{a_1} + {b_1} + 1} \right)x} \right]{u_x} + {b_1}y{u_y} - {a_1}z{u_z} + {a_1}{b_1}u = 0,  \\
  y( {1 - y}){u_{yy}} - xy{u_{xy}} + \left[ {c - \left( {{a_1} + {a_2} + 1} \right)y} \right]{u_y} - {a_2}x{u_x} - {a_1}{a_2}u = 0,  \\
  z({1 + z}){u_{zz}} - x({1 + z}){u_{xz}}
  + \left[ {1 - {b_1} + \left( {{a_3} + {b_2} + 1} \right)z} \right]{u_z} - {a_3}x{u_x} + {a_3}{b_2}u = 0,
\end{array}} \right.
$

where $u\equiv \,\,   {F_{4h}}\left( {{a_1},{a_2},{a_3},{b_1},{b_2};c;x,y,z} \right)$.

Particular solutions:

$
{u_1} = {F_{4h}}\left( {{a_1},{a_2},{a_3},{b_1},{b_2};c;x,y,z} \right),
$

$
{u_2} = {y^{1 - c}}{F_{4h}}\left( {1 - c + {a_1},1 - c + {a_2},{a_3},{b_1},{b_2};2 - c;x,y,z} \right).
$

\bigskip 

\begin{equation} \label{eqf1c}
{F_{4i}}\left( {{a_1},{a_2},{a_3},{b_1},{b_2};c;x,y,z} \right)
   = \sum\limits_{m,n,p = 0}^\infty  {} \frac{{{{\left( {{a_1}} \right)}_{m + n}}{{\left( {{a_2}} \right)}_n}{{\left( {{a_3}} \right)}_p}{{\left( {{b_1}} \right)}_{p - m}}{{\left( {{b_2}} \right)}_{m - n}}}}{{{{\left( c \right)}_p}}}\frac{x^m}{m!}\frac{y^n}{n!}\frac{z^p}{p!},
\end{equation}

region of convergence:
$$\left\{
{t < 1,\,\,\,r < \frac{1}{{1 + t}},\,\,\,\sqrt s  < \sqrt {1 + r + rt}  - \sqrt {r + rt}}\right\}.
$$

System of partial differential equations:

$\left\{
{\begin{array}{*{20}{l}}
  x( {1 + x}){u_{xx}} - {y^2}{u_{yy}} - z{u_{xz}}
   + \left[ {1 - {b_1} + \left( {{a_1} + {b_2} + 1} \right)x} \right]{u_x} - \left( {{a_1} - {b_2} + 1} \right)y{u_y} + {a_1}{b_2}u = 0,  \\
  y\left( {1 + y} \right){u_{yy}} - x\left( {1 - y} \right){u_{xy}}
   + \left[ {1 - {b_2} + \left( {{a_1} + {a_2} + 1} \right)y} \right]{u_y} + x{a_2}{u_x} + {a_1}{a_2}u = 0 , \\
  z\left( {1 - z} \right){u_{zz}} + xz{u_{xz}}
   + \left[ {c - \left( {{a_3} + {b_1} + 1} \right)z} \right]{u_z} + {a_3}x{u_x} - {a_3}{b_1}u = 0,
\end{array}} \right.
$

where $u\equiv \,\,   {F_{4i}}\left( {{a_1},{a_2},{a_3},{b_1},{b_2};c;x,y,z} \right)$.

Particular solutions:

$
{u_1} = {F_{4i}}\left( {{a_1},{a_2},{a_3},{b_1},{b_2};c;x,y,z} \right),
$

$
{u_2} = {z^{1 - c}}{F_{4i}}\left( {{a_1},{a_2},1 - c + {a_3},1 - c + {b_1},{b_2};2 - c;x,y,z} \right).
$

\bigskip 

\begin{equation} \label{eqf1c}
{F_{4j}}\left( {{a_1},{a_2},{a_3},{b_1},{b_2};c;x,y,z} \right)
  = \sum\limits_{m,n,p = 0}^\infty  {} \frac{{{{\left( {{a_1}} \right)}_{m + n}}{{\left( {{a_2}} \right)}_m}{{\left( {{a_3}} \right)}_n}{{\left( {{b_1}} \right)}_{p - m}}{{\left( {{b_2}} \right)}_{p - n}}}}{{{{\left( c \right)}_p}}}\frac{x^m}{m!}\frac{y^n}{n!}\frac{z^p}{p!},
\end{equation}

region of convergence:
$$\left\{
{r + s < 1,\,\,\,\sqrt t  < \min \left\{ {1,\sqrt {\frac{{1 - s}}{s}} ,\sqrt {\frac{{1 - r}}{r}} ,\frac{{1 - r - s}}{{2\sqrt {rs} }}} \right\}}\right\}.
$$

System of partial differential equations:

$
\left\{
{\begin{array}{*{20}{l}}
  x\left( {1 + x} \right){u_{xx}} + xy{u_{xy}} - z{u_{xz}}
   + \left[ {1 - {b_1} + \left( {{a_1} + {a_2} + 1} \right)x} \right]{u_x} + {a_2}y{u_y} + {a_1}{a_2}u = 0, \\
  y\left( {1 + y} \right){u_{yy}} + xy{u_{xy}} - z{u_{yz}}
   + \left[ {1 - {b_2} + \left( {{a_1} + {a_3} + 1} \right)y} \right]{u_y} + {a_3}x{u_x} + {a_1}{a_3}u = 0, \\
  z\left( {1 - z} \right){u_{zz}}- xy{u_{xy}} + xz{u_{xz}} + yz{u_{yz}}
   + \left[ {c - \left( {{b_1} + {b_2} + 1} \right)z} \right]{u_z} + {b_2}x{u_x} + {b_1}y{u_y} - {b_1}{b_2}u = 0,
\end{array}} \right.
$

where $u\equiv \,\,   {F_{4j}}\left( {{a_1},{a_2},{a_3},{b_1},{b_2};c;x,y,z} \right) $.

Particular solutions:

$
{u_1} = {F_{4j}}\left( {{a_1},{a_2},{a_3},{b_1},{b_2};c;x,y,z} \right),
$

$
{u_2} = {z^{1 - c}}{F_{4j}}\left( {{a_1},{a_2},{a_3},1 - c + {b_1},1 - c + {b_2};2 - c;x,y,z} \right).
$

\bigskip

\begin{equation} \label{eqf1c}
{F_{4k}}\left( {{a_1},{a_2},{b_1},{b_2},{b_3};x,y,z} \right)
  = \sum\limits_{m,n,p = 0}^\infty  {} {{{{\left( {{a_1}} \right)}_p}{{\left( {{a_2}} \right)}_p}{{\left( {{b_1}} \right)}_{m + n - p}}{{\left( {{b_2}} \right)}_{m - n}}{{\left( {{b_3}} \right)}_{n - m}}}}\frac{x^m}{m!}\frac{y^n}{n!}\frac{z^p}{p!},
\end{equation}

region of convergence:
$$\left\{
{r + s < 1,\,\,t < \frac{1}{{1 + r + s}}} \right\}.
$$

System of partial differential equations:

$
\left\{ {\begin{array}{*{20}{l}}
  x\left( {1 + x} \right){u_{xx}} - y{u_{xy}} - xz{u_{xz}} + yz{u_{yz}} - {y^2}{u_{yy}}\hfill \\\,\,\,\,\,\,\,\,\,
   + \left[ {1 - {b_3} + \left( {{b_2} + {b_1} + 1} \right)x} \right]{u_x} - \left( {{b_1} - {b_2} + 1} \right)y{u_y} - {b_2}z{u_z} + {b_1}{b_2}u = 0,  \\
  y\left( {1 + y} \right){u_{yy}} - x{u_{xy}} - yz{u_{yz}} + xz{u_{xz}} - {x^2}{u_{xx}} \hfill \\\,\,\,\,\,\,\,\,\,
   + \left[ {1 - {b_2} + \left( {{b_3} + {b_1} + 1} \right)y} \right]{u_y} - \left( {{b_1} - {b_3} + 1} \right)x{u_x} - {b_3}z{u_z} + {b_1}{b_3}u = 0,  \\
  {z\left( {1 + z} \right){u_{zz}} - x{u_{xz}} - y{u_{yz}} + \left[ {1 - {b_1} + \left( {{a_1} + {a_2} + 1} \right)z} \right]{u_z} + {a_1}{a_2}u = 0,}
\end{array}} \right.
$

where $u\equiv \,\,   {F_{4k}}\left( {{a_1},{a_2},{b_1},{b_2},{b_3};x,y,z} \right)$.

\bigskip

\begin{equation} \label{eqf1c}
{F_{4l}}\left( {{a_1},{a_2},{b_1},{b_2},{b_3};x,y,z} \right) \hfill \\
  = \sum\limits_{m,n,p = 0}^\infty  {} {{{{\left( {{a_1}} \right)}_n}{{\left( {{a_2}} \right)}_p}{{\left( {{b_1}} \right)}_{m + n - p}}{{\left( {{b_2}} \right)}_{m - n}}{{\left( {{b_3}} \right)}_{p - m}}}}\frac{x^m}{m!}\frac{y^n}{n!}\frac{z^p}{p!},
\end{equation}

region of convergence:
$$\left\{
{s < 1,\,\,\,t < \frac{1}{{1 + s}},\,\,\,r < \min \left\{ {1,\frac{{{{\left( {1 - s} \right)}^2}}}{{4s}},\frac{{1 - t - st}}{{s{t^2}}}} \right\}}\right\}.
$$

System of partial differential equations:

$
\left\{ {\begin{array}{*{20}{l}}
  x\left( {1 + x} \right){u_{xx}} - ( {1 + x})z{u_{xz}} + yz{u_{yz}} - {y^2}{u_{yy}} \hfill \\
   \,\,\,\,\,\,\,\,\,\,\,\,+ \left[ {1 - {b_3} + \left( {{b_1} + {b_2} + 1} \right)x} \right]{u_x} - \left( {{b_1} - {b_2} + 1} \right)y{u_y} - {b_2}z{u_z} + {b_1}{b_2}u = 0,   \\
  y\left( {1 + y} \right){u_{yy}} - x\left( {1 - y} \right){u_{xy}} - yz{u_{yz}}\\
 \,\,\,\,\,\,\,\,\,\,  + \left[ {1 - {b_2} + \left( {{a_1} + {b_1} + 1} \right)y} \right]{u_y} + {a_1}x{u_x} - {a_1}z{u_z} + {a_1}{b_1}u = 0,  \\
  z\left( {1 + z} \right){u_{zz}} - x( {1 + z} ){u_{xz}} - y{u_{yz}} + \left[ {1 - {b_1} + \left( {{a_2} + {b_3} + 1} \right)z} \right]{u_z}
   - {a_2}x{u_x} + {a_2}{b_3}u = 0,
\end{array}} \right.
$

where $u\equiv \,\,   {F_{4l}}\left( {{a_1},{a_2},{b_1},{b_2},{b_3};x,y,z} \right)$.

\bigskip 

\begin{equation} \label{eqf1c}
 {F_{4m}}\left( {{a_1},{a_2},{b_1},{b_2},{b_3};x,y,z} \right)
  = \sum\limits_{m,n,p = 0}^\infty  {} {{{{\left( {{a_1}} \right)}_m}{{\left( {{a_2}} \right)}_n}{{\left( {{b_1}} \right)}_{m + n - p}}{{\left( {{b_2}} \right)}_{p - m}}{{\left( {{b_3}} \right)}_{p - n}}}}\frac{x^m}{m!}\frac{y^n}{n!}\frac{z^p}{p!},
\end{equation}

region of convergence:
$$\left\{
{r + s < 1,\,\,\,t < \min \left\{ {1,\frac{{1 - r - s}}{{rs}}} \right\}}\right\}.
$$

System of partial differential equations:

$
\left\{ {\begin{array}{*{20}{l}}
  x\left( {1 + x} \right){u_{xx}} + xy{u_{xy}} - \left( {1 + x} \right)z{u_{xz}}\\
  \,\,\,\,\,\,\,\,\,\, + \left[ {1 - {b_2} + \left( {{a_1} + {b_1} + 1} \right)x} \right]{u_x} + {a_1}y{u_y} - {a_1}z{u_z} + {a_1}{b_1}u = 0,  \\
  y\left( {1 + y} \right){u_{yy}} + xy{u_{xy}} - \left( {1 + y} \right)z{u_{yz}}\\
  \,\,\,\,\,\,\,\,\,\, + \left[ {1 - {b_3} + \left( {{a_2} + {b_1} + 1} \right)y} \right]{u_y} + {a_2}x{u_x} - {a_2}z{u_z} + {a_2}{b_1}u = 0,   \\
  z\left( {1 + z} \right){u_{zz}} - x\left( {1 + z} \right){u_{xz}} - y\left( {1 + z} \right){u_{yz}} + xy{u_{xy}} \hfill \\
 \,\,\,\,\,\,\,\,\, + \left[ {1 - {b_1} + \left( {{b_2} + {b_3} + 1} \right)z} \right]{u_z} - {b_3}x{u_x} - {b_2}y{u_y} + {b_2}{b_3}u = 0,
\end{array}} \right.
$

where $u\equiv \,\,   {F_{4m}}\left( {{a_1},{a_2},{b_1},{b_2},{b_3};x,y,z} \right)$.

\bigskip

\begin{equation} \label{eqf1c}
{F_{4n}}\left( {{a_1},{a_2},{a_3},{b_1},{b_2};c;x,y,z} \right)
  = \sum\limits_{m,n,p = 0}^\infty  {} \frac{{{{\left( {{a_1}} \right)}_{m + n}}{{\left( {{a_2}} \right)}_n}{{\left( {{a_3}} \right)}_p}{{\left( {{b_1}} \right)}_{m - p}}{{\left( {{b_2}} \right)}_{p - n}}}}{{{{\left( c \right)}_m}}}\frac{x^m}{m!}\frac{y^n}{n!}\frac{z^p}{p!},
\end{equation}

region of convergence:
$$\left\{
{t < 1,\,\,\,s < \frac{1}{{1 + t}},\,\,\,r < \min \left\{ {1 - s,\,\,\,\frac{{1 - t}}{t},\,\,\,\frac{{{{\left( {1 - s - st} \right)}^2}}}{{4st}}} \right\}}\right\}.
$$

System of partial differential equations:

$\left\{
{\begin{array}{*{20}{l}}
  x\left( {1 - x} \right){u_{xx}} - xy{u_{xy}} + xz{u_{xz}} + yz{u_{yz}}\\
\,\,\,\,\,\,\,\,\,\,  + \left[ {c - \left( {{a_1} + {b_1} + 1} \right)x} \right]{u_x} - {b_1}y{u_y} + {a_1}z{u_z} - {a_1}{b_1}u = 0,   \\
  y\left( {1 + y} \right){u_{yy}} + xy{u_{xy}} - z{u_{yz}}
   + \left[ {1 - {b_2} + \left( {{a_1} + {a_2} + 1} \right)y} \right]{u_y} + {a_2}x{u_x} + {a_1}{a_2}u = 0,   \\
  z\left( {1 + z} \right){u_{zz}} - x{u_{xz}} - yz{u_{yz}}
  + \left[ {1 - {b_1} + \left( {{b_2} + {a_3} + 1} \right)z} \right]{u_z} - {a_3}y{u_y} + {a_3}{b_2}u = 0,
\end{array}} \right.
$

where $u\equiv \,\,   {F_{4n}}\left( {{a_1},{a_2},{a_3},{b_1},{b_2};c;x,y,z} \right)$.

Particular solutions:

$
{u_1} = {F_{4n}}\left( {{a_1},{a_2},{a_3},{b_1},{b_2};c;x,y,z} \right),
$

$
{u_2} = {x^{1 - c}}{F_{4n}}\left( {1 - c + {a_1},{a_2},{a_3},1 - c + {b_1},{b_2};2 - c;x,y,z} \right).
$

\bigskip 

\begin{equation} \label{eqf1c}
{F_{4o}}\left( {{a_1},{a_2},{b_1},{b_2},{b_3};x,y,z} \right)\hfill \\
   = \sum\limits_{m,n,p = 0}^\infty  {} {{{{\left( {{a_1}} \right)}_{m + n}}{{\left( {{a_2}} \right)}_m}{{\left( {{b_1}} \right)}_{n - p}}{{\left( {{b_2}} \right)}_{p - n}}{{\left( {{b_3}} \right)}_{p - m}}}}\frac{x^m}{m!}\frac{y^n}{n!}\frac{z^p}{p!},
\end{equation}

region of convergence:
$$\left\{
{s < 1,\,\,t < 1,\,\,\,r < \frac{{1 - s}}{{1 + t}}}\right\}.
$$

System of partial differential equations:

$
\left\{
{\begin{array}{*{20}{l}}
  x\left( {1 + x} \right){u_{xx}} + xy{u_{xy}} - z{u_{xz}}
   + \left[ {1 - {b_3} + \left( {{a_1} + {a_2} + 1} \right)x} \right]{u_x} + {a_2}y{u_y} + {a_1}{a_2}u = 0, \\
  y\left( {1 + y} \right){u_{yy}} + xy{u_{xy}} - \left( {1 + y} \right)z{u_{yz}} - xz{u_{xz}} \\\,\,\,\,\,\,\,\,\,\,
  + \left[ {1 - {b_2} + \left( {{a_1} + {b_1} + 1} \right)y} \right]{u_y} + {b_1}x{u_x} - {a_1}z{u_z} + {a_1}{b_1}u = 0,  \\
  z\left( {1 + z} \right){u_{zz}} + xy{u_{xy}} - xz{u_{xz}} - y\left( {1 + z} \right){u_{yz}} \\\,\,\,\,\,\,\,\,\,\,
   + \left[ {1 - {b_1} + \left( {{b_3} + {b_2} + 1} \right)z} \right]{u_z} - {b_2}x{u_x} - {b_3}y{u_y} + {b_2}{b_3}u = 0,
\end{array}} \right.
$

where $u\equiv \,\,   {F_{4o}}\left( {{a_1},{a_2},{b_1},{b_2},{b_3};x,y,z} \right)$.

\bigskip 

\begin{equation} \label{eqf1c}
{F_{4p}}\left( {{a_1},{a_2},{b_1},{b_2},{b_3};x,y,z} \right) \hfill \\
  = \sum\limits_{m,n,p = 0}^\infty  {} {{{{\left( {{a_1}} \right)}_{m + n}}{{\left( {{a_2}} \right)}_p}{{\left( {{b_1}} \right)}_{m - n}}{{\left( {{b_2}} \right)}_{n - p}}{{\left( {{b_3}} \right)}_{p - m}}}}\frac{x^m}{m!}\frac{y^n}{n!}\frac{z^p}{p!},
\end{equation}

region of convergence:
$$\left\{
{t < 1,\,\,\,\,r < \frac{1}{{1 + t}}}\right.,
$$
$$
{\left.s < \min \left\{ {{{\left( {\sqrt {1 + r}  - \sqrt r } \right)}^2},\,\,\,\frac{{\left( {1 + t} \right)\left( {1 - r - rt} \right)}}{t},\,\,\,\frac{{\left( {1 - t} \right)\left( {1 - r + rt} \right)}}{t}} \right\}\right\}}.
$$

System of partial differential equations:

$
\left\{
{\begin{array}{*{20}{l}}
  x\left( {1 + x} \right){u_{xx}} - z{u_{xz}} - {y^2}{u_{yy}} + \left[ {1 - {b_3} + \left( {{a_1} + {b_1} + 1} \right)x} \right]{u_x}
   - \left( {{a_1} - {b_1} + 1} \right)y{u_y} + {a_1}{b_1}u = 0,
\\
  y\left( {1 + y} \right){u_{yy}} - x\left( {1 - y} \right){u_{xy}} - xz{u_{xz}} - yz{u_{yz}} \\\,\,\,\,\,\,\,\,\,\,
   + \left[ {1 - {b_1} + \left( {{a_1} + {b_2} + 1} \right)y} \right]{u_y}  + {b_2}x{u_x} - {a_1}z{u_z} + {a_1}{b_2}u = 0,   \\
  z\left( {1 + z} \right){u_{zz}} - xz{u_{xz}} - y{u_{yz}} + \left[ {1 - {b_2} + \left( {{a_2} + {b_3} + 1} \right)z} \right]{u_z}
  - {a_2}x{u_x} + {a_2}{b_3}u = 0,
\end{array}} \right.
$

where $u\equiv \,\,   {F_{4p}}\left( {{a_1},{a_2},{b_1},{b_2},{b_3};x,y,z} \right)$.

\bigskip

\begin{equation} \label{eqf1c}
{F_{5b}}\left( {{a_1},{a_2},{a_3},{a_4},b;c;x,y,z} \right) \hfill \\
  = \sum\limits_{m,n,p = 0}^\infty  {} \frac{{{{\left( {{a_1}} \right)}_m}{{\left( {{a_2}} \right)}_n}{{\left( {{a_3}} \right)}_p}{{\left( {{a_4}} \right)}_p}{{\left( b \right)}_{m + n - p}}}}{{{{\left( c \right)}_{m + n}}}}\frac{x^m}{m!}\frac{y^n}{n!}\frac{z^p}{p!},
\end{equation}

first appearance of this function in the literature, and old notation: [10], $D_{(3)}^{1,3}$,

region of convergence:
$$\left\{
{r < 1,\,\,\,s < 1,\,\,t < \frac{1}{{1 + \max \left\{ {r,s} \right\}}}}\right\}.
$$

System of partial differential equations:

$\left\{
{\begin{array}{*{20}{l}}
  x\left( {1 - x} \right){u_{xx}} + \left( {1 - x} \right)y{u_{xy}} + xz{u_{xz}}
  + \left[ {c - \left( {{a_1} + b + 1} \right)x} \right]{u_x} - {a_1}y{u_y} + {a_1}z{u_z} - {a_1}bu = 0, \\
  y\left( {1 - y} \right){u_{yy}} + x\left( {1 - y} \right){u_{xy}} + yz{u_{yz}}
   + \left[ {c - \left( {{a_2} + b + 1} \right)y} \right]{u_y} - {a_2}x{u_x} + {a_2}z{u_z} - {a_2}bu = 0,  \\
  z\left( {1 + z} \right){u_{zz}} - x{u_{xz}} - y{u_{yz}} + \left[ {1 - b + \left( {{a_3} + {a_4} + 1} \right)z} \right]{u_z}
  + {a_3}{a_4}u = 0,
\end{array}} \right.
$

where $u\equiv \,\,   {F_{5b}}\left( {{a_1},{a_2},{a_3},{a_4},b;c;x,y,z} \right)$.

\bigskip 

\begin{equation} \label{eqf1c}
{F_{5c}}\left( {{a_1},{a_2},{a_3},{a_4},b;c;x,y,z} \right)  = \sum\limits_{m,n,p = 0}^\infty  {} \frac{{{{\left( {{a_1}} \right)}_{m + n}}{{\left( {{a_2}} \right)}_m}{{\left( {{a_3}} \right)}_p}{{\left( {{a_4}} \right)}_p}{{\left( b \right)}_{n - p}}}}{{{{\left( c \right)}_{m + n}}}}\frac{x^m}{m!}\frac{y^n}{n!}\frac{z^p}{p!},
\end{equation}

region of convergence:
$$\left\{
{r < 1,\,\,\,s < 1,\,\,\,t < \frac{1}{{1 + s}}}\right\}.
$$

System of partial differential equations:

$\left\{
{\begin{array}{*{20}{l}}
  x\left( {1 - x} \right){u_{xx}} + \left( {1 - x} \right)y{u_{xy}} + \left[ {c - \left( {{a_1} + {a_2} + 1} \right)x} \right]{u_x}
   - {a_2}y{u_y} - {a_1}{a_2}u = 0,\\
  y\left( {1 - y} \right){u_{yy}} + x\left( {1 - y} \right){u_{xy}} + xz{u_{xz}} + yz{u_{yz}}\\
   \,\,\,\,\,\,\,\,\,\,+ \left[ {c - \left( {{a_1} + b + 1} \right)y} \right]{u_y}    - bx{u_x} + {a_1}z{u_z} - {a_1}bu = 0, \\
  {z\left( {1 + z} \right){u_{zz}} - y{u_{yz}} + \left[ {1 - b + \left( {{a_3} + {a_4} + 1} \right)z} \right]{u_z} + {a_3}{a_4}u = 0,}
\end{array}} \right.
$

where $u\equiv \,\,   {F_{5c}}\left( {{a_1},{a_2},{a_3},{a_4},b;c;x,y,z} \right)$.

\bigskip

\begin{equation} \label{eqf1c}
{F_{5d}}\left( {{a_1},{a_2},{a_3},{a_4},b;c;x,y,z} \right) \hfill \\
  = \sum\limits_{m,n,p = 0}^\infty  {} \frac{{{{\left( {{a_1}} \right)}_{m + n}}{{\left( {{a_2}} \right)}_m}{{\left( {{a_3}} \right)}_n}{{\left( {{a_4}} \right)}_p}{{\left( b \right)}_{p - m - n}}}}{{{{\left( c \right)}_p}}}\frac{x^m}{m!}\frac{y^n}{n!}\frac{z^p}{p!},
\end{equation}

first appearance of this function in the literature, and old notation: [10], $D_{(3)}^{2,3}$,

region of convergence:
$$\left\{
t < 1,\,\,\,\max \left\{ {r,s} \right\} < \frac{1}{{1 + t}}\right\} .
$$

System of partial differential equations:

$
\left\{ {\begin{array}{*{20}{l}}
  x\left( {1 + x} \right){u_{xx}} + \left( {1 + x} \right)y{u_{xy}} - z{u_{xz}}
  + \left[ {1 - b + \left( {{a_1} + {a_2} + 1} \right)x} \right]{u_x} +  {a_2}y{u_y} + {a_1}{a_2}u = 0, \\
  y\left( {1 + y} \right){u_{yy}} + x\left( {1 + y} \right){u_{xy}} - z{u_{yz}}
  + \left[ {1 - b + \left( {{a_1} + {a_3} + 1} \right)y} \right]{u_y} + {a_3}x{u_x} + {a_1}{a_3}u = 0,  \\
  z\left( {1 - z} \right){u_{zz}} + xz{u_{xz}} + yz{u_{yz}} + \left[ {c - \left( {{a_4} + b + 1} \right)z} \right]{u_z}
  + {a_4}x{u_x} + {a_4}y{u_y} - {a_4}bu = 0,
\end{array}} \right.
$

where $u\equiv \,\,   {F_{5d}}\left( {{a_1},{a_2},{a_3},{a_4},b;c;x,y,z} \right)$.

Particular solutions:

$
{u_1} = {F_{5d}}\left( {{a_1},{a_2},{a_3},{a_4},b;c;x,y,z} \right),
$

$
{u_2} = {z^{1 - c}}{F_{5d}}\left( {{a_1},{a_2},{a_3},1 - c + {a_4},1 - c + b;2 - c;x,y,z} \right).
$

\bigskip 

\begin{equation} \label{eqf1c}
{F_{5e}}\left( {{a_1},{a_2},{a_3},{b_1},{b_2};x,y,z} \right) \hfill \\
  = \sum\limits_{m,n,p = 0}^\infty  {} {{{{\left( {{a_1}} \right)}_m}{{\left( {{a_2}} \right)}_n}{{\left( {{a_3}} \right)}_p}{{\left( {{b_1}} \right)}_{m + n - p}}{{\left( {{b_2}} \right)}_{p - m - n}}}}\frac{x^m}{m!}\frac{y^n}{n!}\frac{z^p}{p!},
\end{equation}

first appearance of this function in the literature, and old notation: [28], $G_{B}$,

region of convergence:
$$\left\{
{r < 1,\,\,\,s < 1,\,\,\,t < 1}\right\}.
$$

System of partial differential equations:

$
\left\{
{\begin{array}{*{20}{l}}
  x\left( {1 + x} \right){u_{xx}} + \left( {1 + x} \right)y{u_{xy}} - \left( {1 + x} \right)z{u_{xz}} \\\,\,\,\,\,\,\,\,\,
   + \left[ {1 - {b_2} + \left( {{a_1} + {b_1} + 1} \right)x} \right]{u_x}  + {a_1}y{u_y} - {a_1}z{u_z} + {a_1}{b_1}u = 0,   \\
  y\left( {1 + y} \right){u_{yy}} + x\left( {1 + y} \right){u_{xy}} - \left( {1 + y} \right)z{u_{yz}}\\ \,\,\,\,\,\,\,\,\,
  + \left[ {1 - {b_2} + \left( {{a_2} + {b_1} + 1} \right)y} \right]{u_y} + {a_2}x{u_x}- {a_2}z{u_z}  + {a_2}{b_1}u = 0,   \\
  z\left( {1 + z} \right){u_{zz}} - x\left( {1 + z} \right){u_{xz}} - y\left( {1 + z} \right){u_{yz}}\\ \,\,\,\,\,\,\,\,\,
   + \left[ {1 - {b_1} + \left( {{a_3} + {b_2} + 1} \right)z} \right]{u_z} - {a_3}x{u_x} - {a_3}y{u_y} + {a_3}{b_2}u = 0,
\end{array}} \right.
$

where $u\equiv \,\,   {F_{5e}}\left( {{a_1},{a_2},{a_3},{b_1},{b_2};x,y,z} \right)$.

\bigskip 

\begin{equation} \label{eqf1c}
{F_{5f}}\left( {{a_1},{a_2},{a_3},{b_1},{b_2};x,y,z} \right) \hfill \\
  = \sum\limits_{m,n,p = 0}^\infty  {} {{{{\left( {{a_1}} \right)}_{m + n}}{{\left( {{a_2}} \right)}_m}{{\left( {{a_3}} \right)}_p}{{\left( {{b_1}} \right)}_{n - p}}{{\left( {{b_2}} \right)}_{p - m - n}}}}\frac{x^m}{m!}\frac{y^n}{n!}\frac{z^p}{p!},
\end{equation}

region of convergence:
$$\left\{
{t < 1,\,\,\,r < \frac{1}{{1 + t}},\,\,\,s < \min \left\{ {1 ,\,\frac{{1 - r-rt}}{t}\frac{}{}} \right\}}\right\}.
$$

System of partial differential equations:

$\left\{
{\begin{array}{*{20}{l}}
  x\left( {1 + x} \right){u_{xx}} + \left( {1 + x} \right)y{u_{xy}} - z{u_{xz}}
   + \left[ {1 - {b_2} + \left( {{a_1} + {a_2} + 1} \right)x} \right]{u_x}  + {a_2}y{u_y} + {a_1}{a_2}u = 0,   \\
  y\left( {1 + y} \right){u_{yy}} + x\left( {1 + y} \right){u_{xy}} - xz{u_{xz}} - \left( {1 + y} \right)z{u_{yz}} \\\,\,\,\,\,\,\,\,\,
  + \left[ {1 - {b_2} + \left( {{a_1} + {b_1} + 1} \right)y} \right]{u_y}  + {b_1}x{u_x} - {a_1}z{u_z} + {a_1}{b_1}u = 0,   \\
 z\left( {1 + z} \right){u_{zz}} - xz{u_{xz}} - y\left( {1 + z} \right){u_{yz}}\\
  \,\,\,\,\,\,\,\,\,\,+ \left[ {1 - {b_1} + \left( {{a_3} + {b_2} + 1} \right)z} \right]{u_z}  - {a_3}x{u_x} - {a_3}y{u_y} + {a_3}{b_2}u = 0,
\end{array}} \right.
$

where $u\equiv \,\,   {F_{5f}}\left( {{a_1},{a_2},{a_3},{b_1},{b_2};x,y,z} \right)$.

\bigskip

\begin{equation} \label{eqf1c}
{F_{6a}}\left( {{a_1},{a_2},{a_3},{a_4},{a_5};{c_1},{c_2};x,y,z} \right)
   = \sum\limits_{m,n,p = 0}^\infty  {} \frac{{{{\left( {{a_1}} \right)}_{m + n}}{{\left( {{a_2}} \right)}_m}{{\left( {{a_3}} \right)}_n}{{\left( {{a_4}} \right)}_p}{{\left( {{a_5}} \right)}_p}}}{{{{\left( {{c_1}} \right)}_m}{{\left( {{c_2}} \right)}_{n + p}}}}\frac{x^m}{m!}\frac{y^n}{n!}\frac{z^p}{p!},
\end{equation}

first appearance of this function in the literature, and old notation: [24], $F_{6}$,\,\,[33], $F_N$,

region of convergence:
$$\left\{
{r + s < 1,\,\,\,t < 1}\right\}.
$$

System of partial differential equations:

$\left\{
{\begin{array}{*{20}{l}}
  x\left( {1 - x} \right){u_{xx}} - xy{u_{xy}} + \left[ {{c_1} - \left( {{a_1} + {a_2} + 1} \right)x} \right]{u_x}
   - {a_2}y{u_y} - {a_1}{a_2}u = 0, \\
  y\left( {1 - y} \right){u_{yy}} - xy{u_{xy}} + z{u_{yz}} + \left[ {{c_2} - \left( {{a_1} + {a_3} + 1} \right)y} \right]{u_y}
  - {a_3}x{u_x} - {a_1}{a_3}u = 0, \\
  {z\left( {1 - z} \right){u_{zz}} + y{u_{yz}} + \left[ {{c_2} - \left( {{a_4} + {a_5} + 1} \right)z} \right]{u_z} - {a_4}{a_5}u = 0,}
\end{array}} \right.
$

where $u\equiv \,\,   {F_{6a}}\left( {{a_1},{a_2},{a_3},{a_4},{a_5};{c_1},{c_2};x,y,z} \right)$.

Particular solutions:

$
{u_1} = {F_{6a}}\left( {{a_1},{a_2},{a_3},{a_4},{a_5};{c_1},{c_2};x,y,z} \right),
$

$
{u_2} = {x^{1 - {c_1}}}{F_{6a}}\left( {1 - {c_1} + {a_1},1 - {c_1} + {a_2},{a_3},{a_4},{a_5};2 - {c_1},{c_2};x,y,z} \right).
$

\bigskip

\begin{equation} \label{eqf1c}
{F_{6b}}\left( {{a_1},{a_2},{a_3},{a_4},b;c;x,y,z} \right)\hfill \\
   = \sum\limits_{m,n,p = 0}^\infty  {} \frac{{{{\left( {{a_1}} \right)}_{m + n}}{{\left( {{a_2}} \right)}_n}{{\left( {{a_3}} \right)}_p}{{\left( {{a_4}} \right)}_p}{{\left( b \right)}_{m - n - p}}}}{{{{\left( c \right)}_m}}}\frac{x^m}{m!}\frac{y^n}{n!}\frac{z^p}{p!},
\end{equation}

region of convergence:
$$\left\{
{s < 1,\,\,\,t < 1,\,\,\,r < \min \left\{ {1,\,\,\,\frac{{{{\left( {1 - s} \right)}^2}}}{{4s}},\,\,\,\frac{{1 - t}}{t}} \right\}}\right\}.
$$

System of partial differential equations:

$\left\{
{\begin{array}{*{20}{l}}
  x\left( {1 - x} \right){u_{xx}} + xz{u_{xz}} + yz{u_{yz}} + {y^2}{u_{yy}}\\
 \,\,\,\,\,\,\,\,\,\, + \left[ {c - \left( {{a_1} + b + 1} \right)x} \right]{u_x}    + \left( {{a_1} - b + 1} \right)y{u_y} + {a_1}z{u_z} - {a_1}bu = 0,   \\
  y\left( {1 + y} \right){u_{yy}} - x\left( {1 - y} \right){u_{xy}}
  + z{u_{yz}} + \left[ {1 - b + \left( {{a_1} + {a_2} + 1} \right)y} \right]{u_y}
   + {a_2}x{u_x} + {a_1}{a_2}u = 0,  \\
  z\left( {1 + z} \right){u_{zz}} - x{u_{xz}} + y{u_{yz}} + \left[ {1 - b + \left( {{a_3} + {a_4} + 1} \right)z} \right]{u_z}
   + {a_3}{a_4}u = 0,
\end{array}} \right.
$

where $u\equiv \,\,   {F_{6b}}\left( {{a_1},{a_2},{a_3},{a_4},b;c;x,y,z} \right)$.

Particular solutions:

$
{u_1} = {F_{6b}}\left( {{a_1},{a_2},{a_3},{a_4},b;c;x,y,z} \right),
$

$
{u_2} = {x^{1 - c}}{F_{6b}}\left( {1 - c + {a_1},{a_2},{a_3},{a_4},1 - c + b;2 - c;x,y,z} \right).
$

\bigskip

\begin{equation} \label{eqf1c}
{F_{6c}}\left( {{a_1},{a_2},{a_3},{a_4},b;c;x,y,z} \right) \hfill \\
  = \sum\limits_{m,n,p = 0}^\infty  {} \frac{{{{\left( {{a_1}} \right)}_{m + n}}{{\left( {{a_2}} \right)}_m}{{\left( {{a_3}} \right)}_p}{{\left( {{a_4}} \right)}_p}{{\left( b \right)}_{n - m}}}}{{{{\left( c \right)}_{n + p}}}}\frac{x^m}{m!}\frac{y^n}{n!}\frac{z^p}{p!},
\end{equation}

region of convergence:
$$\left\{
{s < 1,\,\,\,r + 2\sqrt {rs}  < 1,\,\,\,t < 1}\right\}.
$$

System of partial differential equations:

$
\left\{
{\begin{array}{*{20}{l}}
  x\left( {1 + x} \right){u_{xx}} - \left( {1 - x} \right)y{u_{xy}} - \left[ {1-b + \left( {{a_1} + {a_2} + 1} \right)x} \right]{u_x}
   + {a_2}y{u_y} + {a_1}{a_2}u = 0, \\
  y\left( {1 - y} \right){u_{yy}} + z{u_{yz}} + {x^2}{u_{xx}} + \left[ {c - \left( {{a_1} + b + 1} \right)y} \right]{u_y}
  + \left( {{a_1} - b + 1} \right)x{u_x} - {a_1}bu = 0, \\
  {z\left( {1 - z} \right){u_{zz}} + y{u_{yz}} + \left[ {c - \left( {{a_4} + {a_3} + 1} \right)z} \right]{u_z} - {a_3}{a_4}u = 0,}
\end{array}} \right.
$

where $u\equiv \,\,   {F_{6c}}\left( {{a_1},{a_2},{a_3},{a_4},b;c;x,y,z} \right)$.

\bigskip 

\begin{equation} \label{eqf1c}
{F_{6d}}\left( {{a_1},{a_2},{a_3},{a_4},b;c;x,y,z} \right) \hfill \\
  = \sum\limits_{m,n,p = 0}^\infty  {} \frac{{{{\left( {{a_1}} \right)}_{m + n}}{{\left( {{a_2}} \right)}_m}{{\left( {{a_3}} \right)}_n}{{\left( {{a_4}} \right)}_p}{{\left( b \right)}_{p - m}}}}{{{{\left( c \right)}_{n + p}}}}\frac{x^m}{m!}\frac{y^n}{n!}\frac{z^p}{p!},
\end{equation}

region of convergence:
$$\left\{
{t < 1,\,\,\,r < \frac{1}{{1 + t}},\,\,\,s < 1 - r}\right\}.
$$

System of partial differential equations:

$\left\{
{\begin{array}{*{20}{l}}
  x\left( {1 + x} \right){u_{xx}} + xy{u_{xy}} - z{u_{xz}} + \left[ {1 - b + \left( {{a_1} + {a_2} + 1} \right)x} \right]{u_x}
   + {a_2}y{u_y} + {a_1}{a_2}u = 0, \\
  y\left( {1 - y} \right){u_{yy}} - xy{u_{xy}} + z{u_{yz}} + \left[ {c - \left( {{a_1} + {a_3} + 1} \right)y} \right]{u_y}
  - {a_3}x{u_x} - {a_1}{a_3}u = 0, \\
  z\left( {1 - z} \right){u_{zz}} + xz{u_{xz}} + y{u_{yz}} + \left[ {c - \left( {{a_4} + b + 1} \right)z} \right]{u_z}
  + {a_4}x{u_x} - {a_4}bu = 0,
\end{array}} \right.
$

where $u\equiv \,\,   {F_{6d}}\left( {{a_1},{a_2},{a_3},{a_4},b;c;x,y,z} \right)$.

\bigskip

\begin{equation} \label{eqf1c}
{F_{6e}}\left( {{a_1},{a_2},{a_3},{b_1},{b_2};x,y,z} \right) \hfill \\
  = \sum\limits_{m,n,p = 0}^\infty  {} {{{{\left( {{a_1}} \right)}_{m + n}}{{\left( {{a_2}} \right)}_p}{{\left( {{a_3}} \right)}_p}{{\left( {{b_1}} \right)}_{n - m}}{{\left( {{b_2}} \right)}_{m - n - p}}}}\frac{x^m}{m!}\frac{y^n}{n!}\frac{z^p}{p!},
\end{equation}

first appearance of this function in the literature, and old notation: [3], $G_C$,

region of convergence:
$$\left\{
{r + s < 1,\,\,\,t < \frac{{1 + 2s + \sqrt {1 - 4rs} }}{{2\left( {1 + r + s} \right)}}}\right\}.
$$

System of partial differential equations:

$
\left\{
{\begin{array}{*{20}{l}}
  x\left( {1 + x} \right){u_{xx}} - y{u_{xy}} - xz{u_{xz}} - {y^2}{u_{yy}} - yz{u_{yz}} \hfill \\
   \,\,\,\,\,\,\,\,\,+ \left[ {1 - {b_1} + \left( {{a_1} + {b_2} + 1} \right)x} \right]{u_x} - \left( {{a_1} - {b_2} + 1} \right)y{u_y} - {a_1}z{u_z} + {a_1}{b_2}u = 0, \\
  y\left( {1 + y} \right){u_{yy}} - x{u_{xy}} + z{u_{yz}} - {x^2}{u_{xx}} \hfill \\
  \,\,\,\,\,\,\,\,\, + \left[ {1 - {b_2} + \left( {{a_1} + {b_1} + 1} \right)y} \right]{u_y} - \left( {{a_1} - {b_1} + 1} \right)x{u_x} + {a_1}{b_1}u = 0,   \\
  {z\left( {1 + z} \right){u_{zz}} - x{u_{xz}} + y{u_{yz}} + \left[ {1 - {b_2} + \left( {{a_2} + {a_3} + 1} \right)z} \right]{u_z} + {a_2}{a_3}u = 0,}
\end{array}} \right.
$

where $u\equiv \,\,   {F_{6e}}\left( {{a_1},{a_2},{a_3},{b_1},{b_2};x,y,z} \right)$.

\bigskip 

\begin{equation} \label{eqf1c}
{F_{6f}}\left( {{a_1},{a_2},{a_3},{b_1},{b_2};x,y,z} \right)\hfill \\
   = \sum\limits_{m,n,p = 0}^\infty  {} {{{{\left( {{a_1}} \right)}_{m + n}}{{\left( {{a_2}} \right)}_n}{{\left( {{a_3}} \right)}_p}{{\left( {{b_1}} \right)}_{p - m}}{{\left( {{b_2}} \right)}_{m - n - p}}}}\frac{x^m}{m!}\frac{y^n}{n!}\frac{z^p}{p!},
\end{equation}

region of convergence:
$$\left\{
 {r < 1,\,\,\,s + 2\sqrt {rs}  < 1,\,\,\,t < \min \left[ {1,\frac{{1 - s + \sqrt {{{\left( {1 - s} \right)}^2} - 4rs} }}{{2r}}} \right]}\right\}.
$$

System of partial differential equations:

$\left\{ {\begin{array}{*{20}{l}}
  x\left( {1 + x} \right){u_{xx}} - \left( {1 + x} \right)z{u_{xz}} - {y^2}{u_{yy}} - yz{u_{yz}} \hfill \\
 \,\,\,\,\,\,\,\,\,+ \left[ {1 - {b_1} + \left( {{a_1} + {b_2} + 1} \right)x} \right]{u_x}
 - \left( {{a_1} - {b_2} + 1} \right)y{u_y} - {a_1}z{u_z} + {a_1}{b_2}u = 0, \hfill \\
  y\left( {1 + y} \right){u_{yy}} - x\left( {1 - y} \right){u_{xy}} + z{u_{yz}}
   + \left[ {1 - {b_2} + \left( {{a_1} + {a_2} + 1} \right)y} \right]{u_y} + {a_2}x{u_x} + {a_1}{a_2}u = 0, \hfill \\
  z\left( {1 + z} \right){u_{zz}} - x\left( {1 + z} \right){u_{xz}} + y{u_{yz}}
   + \left[ {1 - {b_2} + \left( {{b_1} + {a_3} + 1} \right)z} \right]{u_z} - {a_3}x{u_x} + {a_3}{b_1}u = 0,
\end{array}} \right.
$

where $u\equiv \,\,   {F_{6f}}\left( {{a_1},{a_2},{a_3},{b_1},{b_2};x,y,z} \right)$.

\bigskip 

\begin{equation} \label{eqf1c}
{F_{7a}}\left( {{a_1},{a_2},{a_3},{a_4},{a_5};c;x,y,z} \right)
  = \sum\limits_{m,n,p = 0}^\infty  {} \frac{{{{\left( {{a_1}} \right)}_{m + n}}{{\left( {{a_2}} \right)}_m}{{\left( {{a_3}} \right)}_n}{{\left( {{a_4}} \right)}_p}{{\left( {{a_5}} \right)}_p}}}{{{{\left( c \right)}_{m + n + p}}}}\frac{x^m}{m!}\frac{y^n}{n!}\frac{z^p}{p!},
\end{equation}

first appearance of this function in the literature, and old notation: [24], $F_{7}$,\,\,[33], $F_S$,

region of convergence:
$$\left\{
{r < 1,\,\,\,s < 1,\,\,\,t < 1}\right\}.
$$

System of partial differential equations:

$
\left\{
{\begin{array}{*{20}{l}}
  x\left( {1 - x} \right){u_{xx}} + \left( {1 - x} \right)y{u_{xy}} + z{u_{xz}} + \left[ {c - \left( {{a_1} + {a_2} + 1} \right)x} \right]{u_x}
   - {a_2}y{u_y} - {a_1}{a_2}u = 0, \\
  y\left( {1 - y} \right){u_{yy}} + x\left( {1 - y} \right){u_{xy}} + z{u_{yz}} + \left[ {c - \left( {{a_1} + {a_3} + 1} \right)y} \right]{u_y}
   - {a_3}x{u_x} - {a_1}{a_3}u = 0, \\
  z\left( {1 - z} \right){u_{zz}} + x{u_{xz}} + y{u_{yz}} + \left[ {c - \left( {{a_4} + {a_5} + 1} \right)z} \right]{u_z}
  - {a_4}{a_5}u = 0,
\end{array}} \right.
$

where $u\equiv \,\,   {F_{7a}}\left( {{a_1},{a_2},{a_3},{a_4},{a_5};c;x,y,z} \right) $.

\bigskip

\begin{equation} \label{eqf1c}
{F_{8b}}\left( {{a_1},{a_2},{a_3},b;{c_1},{c_2};x,y,z} \right) \hfill \\
  = \sum\limits_{m,n,p = 0}^\infty  {} \frac{{{{\left( {{a_1}} \right)}_{m + n}}{{\left( {{a_2}} \right)}_p}{{\left( {{a_3}} \right)}_p}{{\left( b \right)}_{m + n - p}}}}{{{{\left( {{c_1}} \right)}_m}{{\left( {{c_2}} \right)}_n}}}\frac{x^m}{m!}\frac{y^n}{n!}\frac{z^p}{p!},
\end{equation}

first appearance of this function in the literature, and old notation: [32],

region of convergence:
$$\left\{
{\sqrt r  + \sqrt s  < 1,\,\,\,t < \frac{1}{{1 + {{\left( {\sqrt r  + \sqrt s } \right)}^2}}}}\right\}.
$$

System of partial differential equations:

$\left\{
{\begin{array}{*{20}{l}}
  x\left( {1 - x} \right){u_{xx}} - 2xy{u_{xy}} + xz{u_{xz}} - {y^2}{u_{yy}} + yz{u_{yz}}\hfill \\
  \,\,\,\,\,\,\,\,\, + \left[ {{c_1} - \left( {{a_1} + b + 1} \right)x} \right]{u_x}  - \left( {{a_1} + b + 1} \right)y{u_y}
   + {a_1}z{u_z} - {a_1}bu = 0,   \\
  y\left( {1 - y} \right){u_{yy}} - 2xy{u_{xy}} + yz{u_{yz}} - {x^2}{u_{xx}} + xz{u_{xz}} \hfill \\
   \,\,\,\,\,\,\,\,\, + \left[ {{c_2} - \left( {{a_1} + b + 1} \right)y} \right]{u_y} - \left( {{a_1} + b + 1} \right)x{u_x}
   + {a_1}z{u_z} - {a_1}bu = 0, \\
  z\left( {1 + z} \right){u_{zz}} - x{u_{xz}} - y{u_{yz}}
   + \left[ {1 - b + \left( {{a_2} + {a_3} + 1} \right)z} \right]{u_z} + {a_2}{a_3}u = 0,
\end{array}} \right.
$

where $u\equiv \,\,   {F_{8b}}\left( {{a_1},{a_2},{a_3},b;{c_1},{c_2};x,y,z} \right)$.

Particular solutions:

$
{u_1} = {F_{8b}}\left( {{a_1},{a_2},{a_3},b;{c_1},{c_2};x,y,z} \right),
$

$
{u_2} = {x^{1 - {c_1}}}{F_{8b}}\left( {1 - {c_1} + {a_1},{a_2},{a_3},1 - {c_1} + b;2 - {c_1},{c_2};x,y,z} \right),
$

$
{u_3} = {y^{1 - {c_2}}}{F_{8b}}\left( {1 - {c_2} + {a_1},{a_2},{a_3},1 - {c_2} + b;{c_1},2 - {c_2};x,y,z} \right),
$

$
{u_4} = {x^{1 - {c_1}}}{y^{1 - {c_2}}}{F_{8b}}\left( {2 - {c_1} - {c_2} + {a_1},{a_2},{a_3},2 - {c_1} - {c_2} + b;2 - {c_1},2 - {c_2};x,y,z} \right).
$

\bigskip 

\begin{equation} \label{eqf1c}
{F_{8c}}\left( {{a_1},{a_2},{a_3},b;{c_1},{c_2};x,y,z} \right) \hfill \\
  = \sum\limits_{m,n,p = 0}^\infty  {} \frac{{{{\left( {{a_1}} \right)}_{m + n}}{{\left( {{a_2}} \right)}_{m + n}}{{\left( {{a_3}} \right)}_p}{{\left( b \right)}_{p - m}}}}{{{{\left( {{c_1}} \right)}_n}{{\left( {{c_2}} \right)}_p}}}\frac{x^m}{m!}\frac{y^n}{n!}\frac{z^p}{p!},
\end{equation}

region of convergence:
$$\left\{
{t < 1,\,\,\,\sqrt s  + \sqrt {r + rt}  < 1}\right\}.
$$

System of partial differential equations:

$
\left\{
{\begin{array}{*{20}{l}}
  x\left( {1 + x} \right){u_{xx}} + 2xy{u_{xy}} - z{u_{xz}} + {y^2}{u_{yy}}\\
\,\,\,\,\,\,\,\,\,\,  + \left[ {1 - b + \left( {{a_1} + {a_2} + 1} \right)x} \right]{u_x}   + \left( {{a_1} + {a_2} + 1} \right)y{u_y} + {a_1}{a_2}u = 0, \\
  y\left( {1 - y} \right){u_{yy}} - 2xy{u_{xy}} - {x^2}{u_{xx}}
   + \left[ {{c_1} - \left( {{a_1} + {a_2} + 1} \right)y} \right]{u_y}  - \left( {{a_1} + {a_2} + 1} \right)x{u_x} - {a_1}{a_2}u = 0,   \\
  z\left( {1 - z} \right){u_{zz}} + xz{u_{xz}} + \left[ {{c_2} - \left( {{a_3} + b + 1} \right)z} \right]{u_z}
  + {a_3}x{u_x} - {a_3}bu = 0,
\end{array}} \right.
$

where $u\equiv \,\,   {F_{8c}}\left( {{a_1},{a_2},{a_3},b;{c_1},{c_2};x,y,z} \right)$.

Particular solutions:

$
{u_1} = {F_{8c}}\left( {{a_1},{a_2},{a_3},b;{c_1},{c_2};x,y,z} \right),
$

$
{u_2} = {y^{1 - {c_1}}}{F_{8c}}\left( {1 - {c_1} + {a_1},1 - {c_1} + {a_2},{a_3},b;2 - {c_1},{c_2};x,y,z} \right),
$

$
{u_3} = {z^{1 - {c_2}}}{F_{8c}}\left( {{a_1},{a_2},{c_2} + {a_3},{c_2} + b;{c_1},2 - {c_2};x,y,z} \right),
$

$
{u_4} = {y^{1 - {c_1}}}{z^{1 - {c_2}}} {F_{8c}}\left( {1 - {c_1} + {a_1},1 - {c_1} + {a_2},1 - {c_2} + {a_3},1 - {c_2} + b;2 - {c_1},2 - {c_2};x,y,z} \right).
$

\bigskip 

\begin{equation} \label{eqf1c}
{F_{8d}}\left( {{a_1},{a_2},{b_1},{b_2};c;x,y,z} \right) \hfill \\
  = \sum\limits_{m,n,p = 0}^\infty  {} \frac{{{{\left( {{a_1}} \right)}_{m + n}}{{\left( {{a_2}} \right)}_p}{{\left( {{b_1}} \right)}_{m + n - p}}{{\left( {{b_2}} \right)}_{p - m}}}}{{{{\left( c \right)}_n}}}\frac{x^m}{m!}\frac{y^n}{n!}\frac{z^p}{p!},
\end{equation}

region of convergence:
$$\left\{
 {s < 1,\,\,\,t < \frac{1}{{1 + s}},\,\,\,r < \min \left\{ {{{\left( {1 - \sqrt s } \right)}^2},\frac{{1 - t - st}}{{t\left( {1 - t} \right)}}} \right\}}\right\}.
$$

System of partial differential equations:

$
\left\{ {\begin{array}{*{20}{l}}
  x\left( {1 + x} \right){u_{xx}} + 2xy{u_{xy}} - \left( {1 + x} \right)z{u_{xz}} - yz{u_{yz}} + {y^2}{u_{yy}} \hfill \\
    \,\,\,\,\,\,\,\,\,+ \left[ {1 - {b_2} + \left( {{a_1} + {b_1} + 1} \right)x} \right]{u_x} + \left( {{a_1} + {b_1} + 1} \right)y{u_y} - {a_1}z{u_z} + {a_1}{b_1}u = 0,  \\
  y\left( {1 - y} \right){u_{yy}} - 2xy{u_{xy}} + xz{u_{xz}} + yz{u_{yz}} - {x^2}{u_{xx}} \hfill \\
  \,\,\,\,\,\,\,\,\, + \left[ {c - \left( {{a_1} + {b_1} + 1} \right)y} \right]{u_y} - \left( {{a_1} + {b_1} + 1} \right)x{u_x} + {a_1}z{u_z} - {a_1}{b_1}u = 0,   \\
  z\left( {1 + z} \right){u_{zz}} - x\left( {1 + z} \right){u_{xz}} - y{u_{yz}}
   + \left[ {1 - {b_1} + \left( {{a_2} + {b_2} + 1} \right)z} \right]{u_z} - {a_2}x{u_x} + {a_2}{b_2}u = 0,
\end{array}} \right.
$

where $u\equiv \,\,   {F_{8d}}\left( {{a_1},{a_2},{b_1},{b_2};c;x,y,z} \right)$.

Particular solutions:

$
{u_1} = {F_{8d}}\left( {{a_1},{a_2},{b_1},{b_2};c;x,y,z} \right),
$

$
{u_2} = {y^{1 - c}}{F_{8d}}\left( {1 - c + {a_1},{a_2},1 - c + {b_1},{b_2};2 - c;x,y,z} \right).
$

\bigskip 

\begin{equation} \label{eqf1c}
{F_{8e}}\left( {{a_1},{a_2},{b_1},{b_2};c;x,y,z} \right) \hfill \\
  = \sum\limits_{m,n,p = 0}^\infty  {} \frac{{{{\left( {{a_1}} \right)}_{m + n}}{{\left( {{a_2}} \right)}_{m + n}}{{\left( {{b_1}} \right)}_{p - m}}{{\left( {{b_2}} \right)}_{p - n}}}}{{{{\left( c \right)}_p}}}\frac{x^m}{m!}\frac{y^n}{n!}\frac{z^p}{p!},
\end{equation}

region of convergence:
\[
\begin{gathered}
  \left\{ {\sqrt r  + \sqrt s  < 1,\,\,\,t < \min \left\{ {1,U_{rs}^ + \left( {{w_1}} \right),U_{rs}^ - \left( {{w_{rs}}} \right),U_{sr}^ - \left( {{w_{sr}}} \right)} \right\}} \right\}, \hfill \\
  P_{rs}^ \pm \left( w \right) = w \pm \frac{1}{2}\sqrt {{w^2} \pm 4r}  + \frac{1}{2}\sqrt {{w^2} \pm 4s}  - 1, \hfill \\
  {w_1}: {\rm{the\,positive\,root\,of}}\,P_{rs}^ + \left( w \right) = 0, \hfill \\
  {w_{rs}}: {\rm{the\,greater\,root\,in}}\,\left( {2\sqrt r ,\infty } \right)\,\,{\rm{of}}\,P_{rs}^ + \left( w \right) = 0, \hfill \\
  U_{rs}^ \pm \left( {{w}} \right) = \frac{{{w^2}}}{{4rs}}\left( {w \pm \sqrt {{w^2} \pm 4r} } \right)\left( {w + \sqrt {{w^2} + 4s} } \right). \hfill \\
\end{gathered}
\]

System of partial differential equations:

$
\left\{
{\begin{array}{*{20}{l}}
  x\left( {1 + x} \right){u_{xx}} + 2xy{u_{xy}} - z{u_{xz}} + {y^2}{u_{yy}}  + \left[ {1 - {b_1} + \left( {{a_1} + {a_2} + 1} \right)x} \right]{u_x}  \\ \,\,\,\,\,\,\,\,\,\,+ \left( {{a_1} + {a_2} + 1} \right)y{u_y} + {a_1}{a_2}u = 0,  \\
  y\left( {1 + y} \right){u_{yy}} + 2xy{u_{xy}} - z{u_{yz}} + {x^2}{u_{xx}} + \left[ {1 - {b_2} + \left( {{a_1} + {a_2} + 1} \right)y} \right]{u_y} \\\,\,\,\,\,\,\,\,\,\,+ \left( {{a_1} + {a_2} + 1} \right)x{u_x} + {a_1}{a_2}u = 0,  \\
 z\left( {1 - z} \right){u_{zz}} - xy{u_{xy}} + xz{u_{xz}} + yz{u_{yz}}
   + \left[ {c - \left( {{b_1} + {b_2} + 1} \right)z} \right]{u_z} + {b_2}x{u_x} + {b_1}y{u_y} - {b_1}{b_2}u = 0,
\end{array}} \right.
$

where $u\equiv \,\,   {F_{8e}}\left( {{a_1},{a_2},{b_1},{b_2};c;x,y,z} \right)$.

Particular solutions:

$
{u_1} = {F_{8e}}\left( {{a_1},{a_2},{b_1},{b_2};c;x,y,z} \right),
$

$
{u_2} = {z^{1 - c}}{F_{8e}}\left( {{a_1},{a_2},1 - c + {b_1},1 - c + {b_2};2 - c;x,y,z} \right).
$

\bigskip 

\begin{equation} \label{eqf1c}
{F_{8f}}\left( {a,{b_1},{b_2},{b_3};c;x,y,z} \right) \hfill \\
  = \sum\limits_{m,n,p = 0}^\infty  {} {{{{\left( a \right)}_{m + n}}{{\left( {{b_1}} \right)}_{m + n - p}}{{\left( {{b_2}} \right)}_{p - m}}{{\left( {{b_3}} \right)}_{p - n}}}}\frac{x^m}{m!}\frac{y^n}{n!}\frac{z^p}{p!},
\end{equation}

region of convergence:
$$\left\{
{\sqrt r  + \sqrt s  < 1,\,\,\,t < \min \left\{ {1,\frac{{{{\left( {1 - r - s} \right)}^2} - 4rs}}{{4rs}}} \right\}}\right\}.
$$

System of partial differential equations:

$\left\{ {\begin{array}{*{20}{l}}
  x\left( {1 + x} \right){u_{xx}} + 2xy{u_{xy}} - \left( {1 + x} \right)z{u_{xz}} + {y^2}{u_{yy}} - yz{u_{yz}} \hfill \\
  \,\,\,\,\,\,\,\,\, + \left[ {1 - {b_2} + \left( {a + {b_1} + 1} \right)x} \right]{u_x}   + \left( {a + {b_1} + 1} \right)y{u_y} - az{u_z} + a{b_1}u = 0,  \\
  y\left( {1 + y} \right){u_{yy}} + 2xy{u_{xy}} + x^2{u_{xx}} - xz{u_{xz}} - \left( {1 + y} \right)z{u_{yz}} \hfill \\
  \,\,\,\,\,\,\,\,\, + \left[ {1 - {b_3} + \left( {a + {b_1} + 1} \right)y} \right]{u_y}  + \left( {a + {b_1} + 1} \right)x{u_x} - az{u_z} + a{b_1}u = 0 , \\
  z\left( {1 + z} \right){u_{zz}} - x\left( {1 + z} \right){u_{xz}} + xy{u_{xy}} - y\left( {1 + z} \right){u_{yz}}\hfill \\
  \,\,\,\,\,\,\,\,\,  + \left[ {1 - {b_1} + \left( {{b_2} + {b_3} + 1} \right)z} \right]{u_z}   - {b_3}x{u_x} - {b_2}y{u_y} + {b_2}{b_3}u = 0,
\end{array}} \right.
$

where $u\equiv \,\,   {F_{8f}}\left( {a,{b_1},{b_2},{b_3};c;x,y,z} \right)$.

\bigskip 

\begin{equation} \label{eqf1c}
{F_{9a}}\left( {{a_1},{a_2},{a_3},{a_4};{c_1},{c_2};x,y,z} \right)
  = \sum\limits_{m,n,p = 0}^\infty  {} \frac{{{{\left( {{a_1}} \right)}_{m + n}}{{\left( {{a_2}} \right)}_{m + n}}{{\left( {{a_3}} \right)}_p}{{\left( {{a_4}} \right)}_p}}}{{{{\left( {{c_1}} \right)}_m}{{\left( {{c_2}} \right)}_{n + p}}}}\frac{x^m}{m!}\frac{y^n}{n!}\frac{z^p}{p!},
\end{equation}

first appearance of this function in the literature, and old notation: [24], $F_{10}$,\,\,[33], $F_R$,

region of convergence:
$$\left\{
{\sqrt r  + \sqrt s  < 1,\,\,\,t < 1}\right\}.
$$

System of partial differential equations:

$
\left\{
{\begin{array}{*{20}{l}}
  x\left( {1 - x} \right){u_{xx}} - 2xy{u_{xy}} - {y^2}{u_{yy}} + \left[ {{c_1} - \left( {{a_1} + {a_2} + 1} \right)x} \right]{u_x}
  - \left( {{a_1} + {a_2} + 1} \right)y{u_y} - {a_1}{a_2}u = 0,  \\
  y\left( {1 - y} \right){u_{yy}} - 2xy{u_{xy}} - {x^2}{u_{xx}} + z{u_{yz}}  + \left[ {{c_2} - \left( {{a_1} + {a_2} + 1} \right)y} \right]{u_y}  \\\,\,\,\,\,\,\,\,\,\, - \left( {{a_1} + {a_2} + 1} \right)x{u_x} - {a_1}{a_2}u = 0, \\
  {z ( {1 - z}){u_{zz}} + y{u_{yz}} + \left[ {{c_2} - \left( {{a_3} + {a_4} + 1} \right)z} \right]{u_z} - {a_3}{a_4}u = 0,}
\end{array}} \right.
$

where $u\equiv \,\,   {F_{9a}}\left( {{a_1},{a_2},{a_3},{a_4};{c_1},{c_2};x,y,z} \right)$.

Particular solutions:

$
{u_1} = {F_{9a}}\left( {{a_1},{a_2},{a_3},{a_4};{c_1},{c_2};x,y,z} \right),
$

$
{u_2} = {x^{1 - {c_1}}}{F_{9a}}\left( {1 - {c_1} + {a_1},1 - {c_1} + {a_2},{a_3},{a_4};2 - {c_1},{c_2};x,y,z} \right).
$

\bigskip

\begin{equation} \label{eqf1c}
{F_{9b}}\left( {{a_1},{a_2},{a_3},b;c;x,y,z} \right) \hfill \\
  = \sum\limits_{m,n,p = 0}^\infty  {} \frac{{{{\left( {{a_1}} \right)}_{m + n}}{{\left( {{a_2}} \right)}_{m + n}}{{\left( {{a_3}} \right)}_p}{{\left( b \right)}_{p - m}}}}{{{{\left( c \right)}_{n + p}}}}\frac{x^m}{m!}\frac{y^n}{n!}\frac{z^p}{p!},
\end{equation}

region of convergence:
$$\left\{
 {\sqrt r  + \sqrt s  < 1,\,\,\,t < \min \left\{ {1,\frac{{1 - r - s + \sqrt {{{\left( {1 - r - s} \right)}^2} - 4rs} }}{{2r}}} \right\}}\right\}.
$$

System of partial differential equations:

$
\left\{
{\begin{array}{*{20}{l}}
  x\left( {1 + x} \right){u_{xx}} + 2xy{u_{xy}} - z{u_{xz}} + {y^2}{u_{yy}} + \left[ {1 - b + \left( {{a_1} + {a_2} + 1} \right)x} \right]{u_x}\\\,\,\,\,\,\,\,\,\,\,
  + \left( {{a_1} + {a_2} + 1} \right)y{u_y} + {a_1}{a_2}u = 0 ,  \\
  y\left( {1 - y} \right){u_{yy}} - 2xy{u_{xy}} + z{u_{yz}} - {x^2}{u_{xx}} + \left[ {c - \left( {{a_1} + {a_2} + 1} \right)y} \right]{u_y}\\\,\,\,\,\,\,\,\,\,\,
   - \left( {{a_1} + {a_2} + 1} \right)x{u_x} - {a_1}{a_2}u = 0 ,  \\
  {z\left( {1 - z} \right){u_{zz}} + xz{u_{xz}} + y{u_{yz}} + \left[ {c - \left( {{a_3} + b + 1} \right)z} \right]{u_z} + {a_3}x{u_x} - {a_3}bu = 0,}
\end{array}} \right.
$

where $u\equiv \,\,   {F_{9b}}\left( {{a_1},{a_2},{a_3},b;c;x,y,z} \right)$.

\bigskip

\begin{equation} \label{eqf1c}
{F_{10a}}\left( {{a_1},{a_2},{a_3},{a_4};{c_1},{c_2},{c_3};x,y,z} \right) \hfill \\
  = \sum\limits_{m,n,p = 0}^\infty  {} \frac{{{{\left( {{a_1}} \right)}_{m + n}}{{\left( {{a_2}} \right)}_{n + p}}{{\left( {{a_3}} \right)}_m}{{\left( {{a_4}} \right)}_p}}}{{{{\left( {{c_1}} \right)}_m}{{\left( {{c_2}} \right)}_n}{{\left( {{c_3}} \right)}_p}}}\frac{x^m}{m!}\frac{y^n}{n!}\frac{z^p}{p!},
\end{equation}

first appearance of this function in the literature, and old notation: [24], $F_{3}$,\,\,[33], $F_K$,

region of convergence:
$$\left\{
{r < 1,\,\,\,t < 1,\,\,\,s < \left( {1 - r} \right)\left( {1 - t} \right)}\right\}.
$$

System of partial differential equations:

$\left\{
{\begin{array}{*{20}{l}}
  {x\left( {1 - x} \right){u_{xx}} - xy{u_{xy}} + \left[ {{c_1} - \left( {{a_1} + {a_3} + 1} \right)x} \right]{u_x} - {a_3}y{u_y} - {a_1}{a_3}u = 0,} \\
  y\left( {1 - y} \right){u_{yy}} - xy{u_{xy}} - xz{u_{xz}} - yz{u_{yz}} + \left[ {{c_2} - \left( {{a_1} + {a_2} + 1} \right)y} \right]{u_y}\\\,\,\,\,\,\,\,\,\,\,
  - {a_2}x{u_x} - {a_1}z{u_z} - {a_1}{a_2}u = 0,  \\
  {z\left( {1 - z} \right){u_{zz}} - yz{u_{yz}} + \left[ {{c_3} - \left( {{a_2} + {a_4} + 1} \right)z} \right]{u_z} - {a_4}y{u_y} - {a_2}{a_4}u = 0,}
\end{array}} \right.
$

where $u\equiv \,\,   {F_{10a}}\left( {{a_1},{a_2},{a_3},{a_4};{c_1},{c_2},{c_3};x,y,z} \right)$.

Particular solutions:

$
{u_1} = {F_{10a}}\left( {{a_1},{a_2},{a_3},{a_4};{c_1},{c_2},{c_3};x,y,z} \right),
$

$
{u_2} = {x^{1 - {c_1}}}{F_{10a}}\left( {1 - {c_1} + {a_1},{a_2},1 - {c_1} + {a_3},{a_4};2 - {c_1},{c_2},{c_3};x,y,z} \right),
$

$
{u_3} = {y^{1 - {c_2}}}{F_{10a}}\left( {1 - {c_2} + {a_1},1 - {c_2} + {a_2},{a_3},{a_4};2 - {c_1},{c_2},{c_3};x,y,z} \right),
$

$
{u_4} = {z^{1 - {c_3}}}{F_{10a}}\left( {{a_1},1 - {c_3} + {a_2},{a_3},1 - {c_3} + {a_4};{c_1},{c_2},2 - {c_3};x,y,z} \right),
$

$
{u_5} = {x^{1 - {c_1}}}{y^{1 - {c_2}}}{F_{10a}}\left( {2 - {c_1} - {c_2} + {a_1},1 - {c_2} + {a_2},1 - {c_1} + {a_3},{a_4};2 - {c_1},2 - {c_2},{c_3};x,y,z} \right),
$

$
{u_6} = {y^{1 - {c_2}}}{z^{1 - {c_3}}}
{F_{10a}}\left( {1 - {c_2} + {a_1},1 - {c_2} + 1 - {c_3} + {a_2},{a_3},1 - {c_3} + {a_4};{c_1},2 - {c_2},2 - {c_3};x,y,z} \right),
$

$
{u_7} = {x^{1 - {c_1}}}{z^{1 - {c_3}}}{F_{10a}}\left( {1 - {c_1} + {a_1},1 - {c_2} + {a_2},1 - {c_1} + {a_3},1 - {c_2} + {a_4};2 - {c_1},{c_2},2 - {c_3};x,y,z} \right),
$

$
\begin{aligned}
  {u_8}& = {x^{1 - {c_1}}}{y^{1 - {c_2}}}{z^{1 - {c_3}}} \times\\& \times {F_{10a}}\left( {2 - {c_1} - {c_2} + {a_1},2 - {c_2} - {c_3} + {a_2},1 - {c_1} + {a_3},1 - {c_3} + {a_4};}
 {2 - {c_1},2 - {c_2},2 - {c_3};x,y,z} \right).
 \end{aligned}
$

\bigskip

\begin{equation} \label{eqf1c}
{F_{10b}}\left( {{a_1},{a_2},{a_3},b;{c_1},{c_2};x,y,z} \right) \hfill \\
  = \sum\limits_{m,n,p = 0}^\infty  {} \frac{{{{\left( {{a_1}} \right)}_{n + p}}{{\left( {{a_2}} \right)}_m}{{\left( {{a_3}} \right)}_p}{{\left( b \right)}_{m + n - p}}}}{{{{\left( {{c_1}} \right)}_m}{{\left( {{c_2}} \right)}_n}}}\frac{x^m}{m!}\frac{y^n}{n!}\frac{z^p}{p!},
\end{equation}

region of convergence:
$$\left\{
{s < 1,\,\,\,t + 2\sqrt {st}  < 1,\,\,\,r < \min \left\{ {1 - s,\frac{{1 - t - 2\sqrt {st} }}{t}} \right\}}\right\}.
$$

System of partial differential equations:

$
\left\{
{\begin{array}{*{20}{l}}
  x\left( {1 - x} \right){u_{xx}} - xy{u_{xy}} + xz{u_{xz}} + \left[ {{c_1} - \left( {{a_2} + b + 1} \right)x} \right]{u_x}
   - {a_2}y{u_y} + {a_2}z{u_z} - {a_2}bu = 0, \\
  y\left( {1 - y} \right){u_{yy}} - xy{u_{xy}} - xz{u_{xz}} + {z^2}{u_{zz}} + \left[ {{c_2} - \left( {{a_1} + b + 1} \right)y} \right]{u_y}\\\,\,\,\,\,\,\,\,\,\,
   - {a_1}x{u_x} + \left( {{a_1} - b + 1} \right)z{u_z} - {a_1}bu = 0,  \\
  z\left( {1 + z} \right){u_{zz}} - x{u_{xz}} - y\left( {1 - z} \right){u_{yz}} + \left[ {1 - b + \left( {{a_1} + {a_3} + 1} \right)z} \right]{u_z}
   + {a_3}y{u_y} + {a_1}{a_3}u = 0,
\end{array}} \right.
$

where $u\equiv \,\,   {F_{10b}}\left( {{a_1},{a_2},{a_3},b;{c_1},{c_2};x,y,z} \right)$.

Particular solutions:

$
{u_1} = {F_{10b}}\left( {{a_1},{a_2},{a_3},b;{c_1},{c_2};x,y,z} \right),
$

$
{u_2} = {x^{1 - {c_1}}}{F_{10b}}\left( {{a_1},1 - {c_1} + {a_2},{a_3},1 - {c_1} + b;2 - {c_1},{c_2};x,y,z} \right),
$

$
{u_3} = {y^{1 - {c_2}}}{F_{10b}}\left( {1 - {c_2} + {a_1},{a_2},{a_3},1 - {c_2} + b;{c_1},2 - {c_2};x,y,z} \right),
$

$
  {u_4} = {x^{1 - {c_1}}}{y^{1 - {c_2}}}  {F_{10b}}\left( {1 - {c_2} + {a_1},1 - {c_1} + {a_2},{a_3},2 - {c_1} - {c_2} + b;2 - {c_1},2 - {c_2};x,y,z} \right).
$

\bigskip

\begin{equation} \label{eqf1c}
{F_{10c}}\left( {{a_1},{a_2},{a_3},b;{c_1},{c_2};x,y,z} \right) \hfill \\
  = \sum\limits_{m,n,p = 0}^\infty  {} \frac{{{{\left( {{a_1}} \right)}_{m + n}}{{\left( {{a_2}} \right)}_{n + p}}{{\left( {{a_3}} \right)}_p}{{\left( b \right)}_{m - p}}}}{{{{\left( {{c_1}} \right)}_m}{{\left( {{c_2}} \right)}_n}}}\frac{x^m}{m!}\frac{y^n}{n!}\frac{z^p}{p!},
\end{equation}

region of convergence:
$$\left\{
{r < 1,\,\,\,t < \frac{1}{{1 + r}},\,\,\,s < \min \left\{ {1 - r,{{\left( {\sqrt {1 - t}  - \sqrt {rt} } \right)}^2}} \right\}}\right\}.
$$

System of partial differential equations:

$
\left\{
{\begin{array}{*{20}{l}}
  x\left( {1 - x} \right){u_{xx}} - xy{u_{xy}} + xz{u_{xz}} + yz{u_{yz}} + \left[ {{c_1} - \left( {{a_1} + b + 1} \right)x} \right]{u_x}
   - by{u_y} + {a_1}z{u_z} - {a_1}bu = 0,  \\
  y\left( {1 - y} \right){u_{yy}} - xy{u_{xy}} - xz{u_{xz}} - yz{u_{yz}} + \left[ {{c_2} - \left( {{a_1} + {a_2} + 1} \right)y} \right]{u_y}\\\,\,\,\,\,\,\,\,\,\,
  - {a_2}x{u_x} - {a_1}z{u_z} - {a_1}{a_2}u = 0, \\
  z\left( {1 + z} \right){u_{zz}} - x{u_{xz}} + yz{u_{yz}} + \left[ {1 - b + \left( {{a_2} + {a_3} + 1} \right)z} \right]{u_z}
   + {a_3}y{u_y} + {a_2}{a_3}u = 0,
\end{array}} \right.
$

where $u\equiv \,\,   {F_{10c}}\left( {{a_1},{a_2},{a_3},b;{c_1},{c_2};x,y,z} \right) $.

Particular solutions:

$
{u_1} = {F_{10c}}\left( {{a_1},{a_2},{a_3},b;{c_1},{c_2};x,y,z} \right),
$

${u_2} = {x^{1 - {c_1}}}{F_{10c}}\left( {1 - {c_1} + {a_1},{a_2},{a_3},1 - {c_1} + b;2 - {c_1},{c_2};x,y,z} \right),
$

$
{u_3} = {y^{1 - {c_2}}}{F_{10c}}\left( {1 - {c_2} + {a_1},1 - {c_2} + {a_2},{a_3},b;{c_1},2 - {c_2};x,y,z} \right),
$

$
  {u_4} = {x^{1 - {c_1}}}{y^{1 - {c_2}}} {F_{10c}}\left( {2 - {c_1} - {c_2} + {a_1},1 - {c_2} + {a_2},{a_3},1 - {c_1} + b;2 - {c_1},2 - {c_2};x,y,z} \right).
  $

\bigskip

\begin{equation} \label{eqf1c}
{F_{10d}}\left( {{a_1},{a_2},{a_3},b;{c_1},{c_2};x,y,z} \right) \hfill \\
  = \sum\limits_{m,n,p = 0}^\infty  {} \frac{{{{\left( {{a_1}} \right)}_{m + n}}{{\left( {{a_2}} \right)}_{n + p}}{{\left( {{a_3}} \right)}_p}{{\left( b \right)}_{m - n}}}}{{{{\left( {{c_1}} \right)}_m}{{\left( {{c_2}} \right)}_p}}}\frac{x^m}{m!}\frac{y^n}{n!}\frac{z^p}{p!},
\end{equation}

region of convergence:
$$\left\{
{r < 1,\,\,\,t < 1,\,\,\,\sqrt s  < \min \left\{ {1 - r,\left( {\sqrt {1 + r}  - \sqrt r } \right)\sqrt {1 - t} } \right\}}\right\}.
$$

System of partial differential equations:

$
\left\{
{\begin{array}{*{20}{l}}
  x\left( {1 - x} \right){u_{xx}} + {y^2}{u_{yy}}
  + \left[ {{c_1} - \left( {{a_1} + b + 1} \right)x} \right]{u_x} + \left( {{a_1} - b + 1} \right)y{u_y} - {a_1}bu = 0,  \\
  y\left( {1 + y} \right){u_{yy}} - x\left( {1 - y} \right){u_{xy}} + xz{u_{xz}} + yz{u_{yz}}\\ \,\,\,\,\,\,\,\,\,
   + \left[ {1 - b + \left( {{a_1} + {a_2} + 1} \right)y} \right]{u_y}  + {a_2}x{u_x} + {a_1}z{u_z} + {a_1}{a_2}u = 0,   \\
  z\left( {1 - z} \right){u_{zz}} - yz{u_{yz}} + \left[ {{c_2} - \left( {{a_2} + {a_3} + 1} \right)z} \right]{u_z}
  - {a_3}y{u_y} - {a_2}{a_3}u = 0,
\end{array}} \right.
$

where $u\equiv \,\,   {F_{10d}}\left( {{a_1},{a_2},{a_3},b;{c_1},{c_2};x,y,z} \right) $.

Particular solutions:

$
{u_1} = {F_{10d}}\left( {{a_1},{a_2},{a_3},b;{c_1},{c_2};x,y,z} \right),
$

$
{u_2} = {x^{1 - {c_1}}}{F_{10d}}\left( {1 - {c_1} + {a_1},{a_2},{a_3},1 - {c_1} + b;2 - {c_1},{c_2};x,y,z} \right),
$

$
{u_3} = {z^{1 - {c_2}}}{F_{10d}}\left( {{a_1},1 - {c_2} + {a_2},1 - {c_2} + {a_3},b;{c_1},2 - {c_2};x,y,z} \right),
$

$
  {u_4} = {x^{1 - {c_1}}}{z^{1 - {c_2}}}  {F_{10d}}\left( {1 - {c_1} + {a_1},1 - {c_2} + {a_2},1 - {c_2} + {a_3},1 - {c_1} + b;2 - {c_1},2 - {c_2};x,y,z} \right).
  $

\bigskip

\begin{equation} \label{eqf1c}
{F_{10e}}\left( {{a_1},{a_2},{b_1},{b_2};c;x,y,z} \right) \hfill \\
  = \sum\limits_{m,n,p = 0}^\infty  {} \frac{{{{\left( {{a_1}} \right)}_m}{{\left( {{a_2}} \right)}_p}{{\left( {{b_1}} \right)}_{m + n - p}}{{\left( {{b_2}} \right)}_{n + p - m}}}}{{{{\left( c \right)}_n}}}\frac{x^m}{m!}\frac{y^n}{n!}\frac{z^p}{p!},
\end{equation}

region of convergence:
$$\left\{
{r < 1,\,\,\,t < 1,\,\,\,s < \min \left\{ {1,\frac{{{{\left( {1 - r} \right)}^2}}}{{4r}},\frac{{{{\left( {1 - t} \right)}^2}}}{{4t}}} \right\}}\right\}.
$$

System of partial differential equations:

$
\left\{ {\begin{array}{*{20}{l}}
  x\left( {1 + x} \right){u_{xx}} - ({1 - x})y{u_{xy}} - ({1 + x})z{u_{xz}} \hfill \\
  \,\,\,\,\,\,\,\,\, + \left[ {1 - {b_2} + \left( {{a_1} + {b_1} + 1} \right)x} \right]{u_x}   + {a_1}y{u_y} - {a_1}z{u_z} + {a_1}{b_1}u = 0,   \\
  y\left( {1 - y} \right){u_{yy}} - 2xz{u_{xz}} + {x^2}{u_{xx}} + {z^2}{u_{zz}} + \left[ {c - \left( {{b_1} + {b_2} + 1} \right)y} \right]{u_y}  \hfill \\
  \,\,\,\,\,\,\,\,\, + \left( {{b_1} - {b_2} + 1} \right)x{u_x} + \left( {{b_2} - {b_1} + 1} \right)z{u_z} - {b_1}{b_2}u = 0,   \\
  z\left( {1 + z} \right){u_{zz}} - x\left( {1 + z} \right){u_{xz}} - y\left( {1 - z} \right){u_{yz}} \hfill \\
 \,\,\,\,\,\,\,\,\, + \left[ {1 - {b_1} + \left( {{a_2} + {b_2} + 1} \right)z} \right]{u_z}  - {a_2}x{u_x} + {a_2}y{u_y} + {a_2}{b_2}u = 0, \hfill \\
\end{array}} \right.
$

where $u\equiv \,\,   {F_{10e}}\left( {{a_1},{a_2},{b_1},{b_2};c;x,y,z} \right)$.

Particular solutions:

$
{u_1} = {F_{10e}}\left( {{a_1},{a_2},{b_1},{b_2};c;x,y,z} \right),
$

$
{u_2} = {y^{1 - c}}{F_{10e}}\left( {{a_1},{a_2},1 - c + {b_1},1 - c + {b_2};2 - c;x,y,z} \right).
$

\bigskip

\begin{equation} \label{eqf1c}
{F_{10f}}\left( {{a_1},{a_2},{b_1},{b_2};c;x,y,z} \right)\hfill \\
   = \sum\limits_{m,n,p = 0}^\infty  {} \frac{{{{\left( {{a_1}} \right)}_{n + p}}{{\left( {{a_2}} \right)}_m}{{\left( {{b_1}} \right)}_{m + n - p}}{{\left( {{b_2}} \right)}_{p - m}}}}{{{{\left( c \right)}_n}}}\frac{x^m}{m!}\frac{y^n}{n!}\frac{z^p}{p!},
\end{equation}

region of convergence:
$$\left\{
{r + s < 1,\,\,\,t < \min \left\{ {{{\left( {\sqrt {1 + s}  - \sqrt s } \right)}^2},\frac{{1 - r - s}}{{r\left( {1 - r} \right)}}} \right\}}\right\}.
$$

System of partial differential equations:

$\left\{
{\begin{array}{*{20}{l}}
  x\left( {1 + x} \right){u_{xx}} + xy{u_{xy}} - ( {1 + x})z{u_{xz}} + \left[ {1 - {b_2} + \left( {{a_2} + {b_1} + 1} \right)x} \right]{u_x}\\\,\,\,\,\,\,\,\,\,\,
  + {a_2}y{u_y} - {a_2}z{u_z} + {a_2}{b_1}u = 0,   \\
  y\left( {1 - y} \right){u_{yy}} - xy{u_{xy}} - xz{u_{xz}} + {z^2}{u_{zz}} + \left[ {c - \left( {{a_1} + {b_1} + 1} \right)y} \right]{u_y}  \hfill \\
   \,\,\,\,\,\,\,\,\, - {a_1}x{u_x} + \left( {{a_1} - {b_1} + 1} \right)z{u_z} - {a_1}{b_1}u = 0,  \\
  \begin{gathered}
  z\left( {1 + z} \right){u_{zz}} - x\left( {1 + z} \right){u_{xz}} - xy{u_{xy}} - y\left( {1 - z} \right){u_{yz}} \hfill \\
   \,\,\,\,\,\,\,\,\,+ \left[ {1 - {b_1} + \left( {{a_1} + {b_2} + 1} \right)z} \right]{u_z}  - {a_1}x{u_x} + {b_2}y{u_y} + {a_1}{b_2}u = 0, \hfill \\
\end{gathered}
\end{array}} \right.
$

where $u\equiv \,\,   {F_{10f}}\left( {{a_1},{a_2},{b_1},{b_2};c;x,y,z} \right)$.

Particular solutions:

$
{u_1} = {F_{10f}}\left( {{a_1},{a_2},{b_1},{b_2};c;x,y,z} \right),
$

$
{u_2} = {y^{1 - c}}{F_{10f}}\left( {1 - c + {a_1},{a_2},1 - c + {b_1},{b_2};2 - c;x,y,z} \right).
$

\bigskip

\begin{equation} \label{eqf1c}
{F_{10g}}\left( {{a_1},{a_2},{b_1},{b_2};c;x,y,z} \right)\hfill \\
   = \sum\limits_{m,n,p = 0}^\infty  {} \frac{{{{\left( {{a_1}} \right)}_{n + p}}{{\left( {{a_2}} \right)}_p}{{\left( {{b_1}} \right)}_{m + n - p}}{{\left( {{b_2}} \right)}_{m - n}}}}{{{{\left( c \right)}_m}}}\frac{x^m}{m!}\frac{y^n}{n!}\frac{z^p}{p!},
\end{equation}

region of convergence:

$$
 \left\{ {s < 1,\,\,\,t + 2\sqrt {st}  < 1,} \,\,
   {r < \min \left\{ {1,\frac{{{{\left( {1 - s} \right)}^2}}}{{4s}},\frac{{1 - t}}{t}{\Psi _1}\left( {\frac{{st}}{{{{\left( {1 - t} \right)}^2}}}} \right),\frac{{1 - t}}{t}{\Psi _2}\left( {\frac{{st}}{{{{\left( {1 - t} \right)}^2}}}} \right)} \right\}}\right\}.
$$

System of partial differential equations:

$
\left\{ {\begin{array}{*{20}{l}}
  x\left( {1 - x} \right){u_{xx}} + xz{u_{xz}} - yz{u_{yz}} + {y^2}{u_{yy}} + \left[ {c - \left( {{b_1} + {b_2} + 1} \right)x} \right]{u_x}  \hfill \\
    \,\,\,\,\,\,\,\,\,+ \left( {{b_1} - {b_2} + 1} \right)y{u_y} + {b_2}z{u_z} - {b_1}{b_2}u = 0,  \\
  y\left( {1 + y} \right){u_{yy}} - x\left( {1 - y} \right){u_{xy}} + xz{u_{xz}} - {z^2}{u_{zz}} + \left[ {1 - {b_2} + \left( {{a_1} + {b_1} + 1} \right)y} \right]{u_y}  \hfill \\
   \,\,\,\,\,\,\,\,\, + {a_1}x{u_x} - \left( {{a_1} - {b_1} + 1} \right)z{u_z} + {a_1}{b_1}u = 0, \\
  z\left( {1 + z} \right){u_{zz}} - x{u_{xz}} - y( {1 - z}){u_{yz}} + \left[ {1 - {b_1} + \left( {{a_1} + {a_2} + 1} \right)z} \right]{u_z}
  + {a_2}y{u_y} + {a_1}{a_2}u = 0,
\end{array}} \right.
$

where $u\equiv \,\,   {F_{10g}}\left( {{a_1},{a_2},{b_1},{b_2};c;x,y,z} \right)$.

Particular solutions:

$
{u_1} = {F_{10g}}\left( {{a_1},{a_2},{b_1},{b_2};c;x,y,z} \right),
$

$
{u_2} = {x^{1 - c}}{F_{10g}}\left( {{a_1},{a_2},1 - c + {b_1},1 - c + {b_2};2 - c;x,y,z} \right).
$

\bigskip

\begin{equation} \label{eqf1c}
{F_{10h}}\left( {{a_1},{a_2},{b_1},{b_2};c;x,y,z} \right)\hfill \\
   = \sum\limits_{m,n,p = 0}^\infty  {} \frac{{{{\left( {{a_1}} \right)}_{m + n}}{{\left( {{a_2}} \right)}_{n + p}}{{\left( {{b_1}} \right)}_{m - p}}{{\left( {{b_2}} \right)}_{p - m}}}}{{{{\left( c \right)}_n}}}\frac{x^m}{m!}\frac{y^n}{n!}\frac{z^p}{p!},
\end{equation}

region of convergence:
$$\left\{
{r < 1,\,\,\,t < 1,\,\,\,s < \left( {1 - r} \right)\left( {1 - t} \right)}\right\}.
$$

System of partial differential equations:

$
\left\{ {\begin{array}{*{20}{l}}
  x\left( {1 + x} \right){u_{xx}} + xy{u_{xy}} - ( {1 + x})z{u_{xz}} - yz{u_{yz}}\hfill \\
    \,\,\,\,\,\,\,\,\,+ \left[ {1 - {b_2} + \left( {{a_1} + {b_1} + 1} \right)x} \right]{u_x} + {b_1}y{u_y}  - {a_1}z{u_z} + {a_1}{b_1}u = 0,   \\
  y\left( {1 - y} \right){u_{yy}} - xy{u_{xy}} - xz{u_{xz}} - yz{u_{yz}}
  + \left[ {c - \left( {{a_1} + {a_2} + 1} \right)y} \right]{u_y} - {a_2}x{u_x}  - {a_1}z{u_z} - {a_1}{a_2}u = 0,  \\
  z\left( {1 + z} \right){u_{zz}} - xy{u_{xy}} - x\left( {1 + z} \right){u_{xz}} + yz{u_{yz}} \hfill \\
  \,\,\,\,\,\,\,\,\, + \left[ {1 - {b_1} + \left( {{a_2} + {b_2} + 1} \right)z} \right]{u_z} - {a_2}x{u_x}  + {b_2}y{u_y} + {a_2}{b_2}u = 0,
\end{array}} \right.
$

where $u\equiv \,\,   {F_{10h}}\left( {{a_1},{a_2},{b_1},{b_2};c;x,y,z} \right)$.

Particular solutions:

$
{u_1} = {F_{10h}}\left( {{a_1},{a_2},{b_1},{b_2};c;x,y,z} \right),
$

$
{u_2} = {y^{1 - c}}{F_{10h}}\left( {1 - c + {a_1},1 - c + {a_2},{b_1},{b_2};2 - c;x,y,z} \right).
$

\bigskip

\begin{equation} \label{eqf1c}
{F_{10i}}\left( {{a_1},{a_2},{b_1},{b_2};c;x,y,z} \right) \hfill \\
  = \sum\limits_{m,n,p = 0}^\infty  {} \frac{{{{\left( {{a_1}} \right)}_{m + n}}{{\left( {{a_2}} \right)}_{n + p}}{{\left( {{b_1}} \right)}_{m - p}}{{\left( {{b_2}} \right)}_{p - n}}}}{{{{\left( c \right)}_m}}}\frac{x^m}{m!}\frac{y^n}{n!}\frac{z^p}{p!},
\end{equation}

region of convergence:
$$
\begin{gathered}
    \left\{ {t < 1,\,\,\,s + 2\sqrt {st}  < 1,\,\,\,r < \min \left\{ {U\left( {{w_1}} \right),U\left( {{w_2}} \right),U\left( {{w_3}} \right)} \right\}}\right\} , \hfill \\
  {P_{st}}\left( w \right) = 2{w^3} - \left( {1 + 2s + 4st} \right){w^2} + {s^2}, \hfill \\
  {w_1}: {\rm{the\,root\,in}}\,\left( {s,1} \right)\,\,{\rm{of}}\,\,{P_{st}}\left( w \right) = 0, \hfill \\
  {w_2}:\rm{the\,smaller\,\,root\,in}\,\left( {s,1} \right)\,\,\rm{of}\,\,{P_{ - s, - t}}\left( w \right) = 0, \hfill \\
  {w_3}: {\rm{the\,smaller\,\,root\,in}}\,\left( {s,1} \right)\,\,{\rm{of}}\,\,{P_{s, - t}}\left( w \right) = 0, \hfill \\
  U\left( w \right) = \frac{{\left( {{w^2} - {s^2}} \right)\left( {1 - {w^2}} \right)}}{{4st{w^2}}}. \hfill \\
\end{gathered}
$$

System of partial differential equations:

$
\left\{ {\begin{array}{*{20}{l}}
  x\left( {1 - x} \right){u_{xx}} - xy{u_{xy}} + xz{u_{xz}} + yz{u_{yz}} + \left[ {c - \left( {{a_1} + {b_1} + 1} \right)x} \right]{u_x}
   - {b_1}y{u_y} + {a_1}z{u_z} - {a_1}{b_1}u = 0,   \\
  y\left( {1 + y} \right){u_{yy}} + xy{u_{xy}} + xz{u_{xz}} - ({1 - y})z{u_{yz}}\hfill \\
   \,\,\,\,\,\,\,\,\, + \left[ {1 - {b_2} + \left( {{a_1} + {a_2} + 1} \right)y} \right]{u_y}  + {a_2}x{u_x} + {a_1}z{u_z} + {a_1}{a_2}u = 0,   \\
  z\left( {1 + z} \right){u_{zz}} - x{u_{xz}}  - {y^2}{u_{yy}}
   + \left[ {1 - {b_1} + \left( {{a_2} + {b_2} + 1} \right)z} \right]{u_z} - \left( {{a_2} - {b_2} + 1} \right)y{u_y} + {a_2}{b_2}u = 0,
\end{array}} \right.
$

where $u\equiv \,\,   {F_{10i}}\left( {{a_1},{a_2},{b_1},{b_2};c;x,y,z} \right)$.

Particular solutions:

$
{u_1} = {F_{10i}}\left( {{a_1},{a_2},{b_1},{b_2};c;x,y,z} \right),
$

$
{u_2} = {x^{1 - c}}{F_{10i}}\left( {1 - c + {a_1},{a_2},1 - c + {b_1},{b_2};2 - c;x,y,z} \right).
$

\bigskip

\begin{equation} \label{eqf1c}
{F_{10j}}\left( {a,{b_1},{b_2},{b_3};x,y,z} \right) \hfill \\
  = \sum\limits_{m,n,p = 0}^\infty  {} {{{{\left( a \right)}_p}{{\left( {{b_1}} \right)}_{m + n - p}}{{\left( {{b_2}} \right)}_{n + p - m}}{{\left( {{b_3}} \right)}_{m - n}}}}\frac{x^m}{m!}\frac{y^n}{n!}\frac{z^p}{p!},
\end{equation}

region of convergence:
$$\left\{
{s < 1,\,\,\,t + 2\sqrt {st}  < 1,\,\,\,r < \min \left\{ {1 - s,\frac{{{{\left( {1 - t} \right)}^2} - 4st}}{{4s{t^2}}}} \right\}}\right\}.
$$

System of partial differential equations:

$
\left\{ {\begin{array}{*{20}{l}}
  x\left( {1 + x} \right){u_{xx}} - y{u_{xy}} - \left( {1 + x} \right)z{u_{xz}} - {y^2}{u_{yy}} + yz{u_{yz}} \\
  \,\,\,\,\,\,\,\,\, + \left[ {1 - {b_2} + \left( {{b_1} + {b_3} + 1} \right)x} \right]{u_x} - \left( {{b_1} - {b_3} + 1} \right)y{u_y}
   - {b_3}z{u_z} + {b_1}{b_3}u = 0,  \\
  y\left( {1 + y} \right){u_{yy}} - x{u_{xy}} + 2xz{u_{xz}} - {x^2}{u_{xx}} - {z^2}{u_{zz}} \\
  \,\,\,\,\,\,\,\,\, + \left[ {1 - {b_3} + \left( {{b_1} + {b_2} + 1} \right)y} \right]{u_y}  - \left( {{b_1} - {b_2} + 1} \right)x{u_x}
   - \left( {{b_2} - {b_1} + 1} \right)z{u_z} + {b_1}{b_2}u = 0,   \\
  z( {1 + z}){u_{zz}} - x( {1 + z}){u_{xz}} - y( {1 - z}){u_{yz}}\\\,\,\,\,\,\,\,\,\,\,
  + \left[ {1 - {b_1} + \left( {a + {b_2} + 1} \right)z} \right]{u_z}  - ax{u_x} + ay{u_y} + a{b_2}u = 0,
\end{array}} \right.
$

where $u\equiv \,\,   {F_{10j}}\left( {a,{b_1},{b_2},{b_3};x,y,z} \right)$.

\bigskip

\begin{equation} \label{eqf1c}
{F_{10k}}\left( {a,{b_1},{b_2},{b_3};x,y,z} \right) \hfill \\
  = \sum\limits_{m,n,p = 0}^\infty  {} {{{{\left( a \right)}_{n + p}}{{\left( {{b_1}} \right)}_{m + n - p}}{{\left( {{b_2}} \right)}_{m - n}}{{\left( {{b_3}} \right)}_{p - m}}}}\frac{x^m}{m!}\frac{y^n}{n!}\frac{z^p}{p!},
\end{equation}

region of convergence:
$$
\begin{gathered}
    \left\{ {s < 1,\,\,\,t + 2\sqrt {st}  < 1,}
   \,\,{r < \min \left\{ {1,\frac{{{{\left( {1 - s} \right)}^2}}}{{s + st}}{\Phi _1}\left( { - \frac{{t\left( {1 - s} \right)}}{{{{\left( {1 + t} \right)}^2}}}} \right),\frac{{{{\left( {1 + s} \right)}^2}}}{{s + st}}{\Phi _2}\left( {\frac{{t\left( {1 + s} \right)}}{{{{\left( {1 + t} \right)}^2}}}} \right)} \right\}} \right\}. \hfill \\
\end{gathered}
$$

System of partial differential equations:

$
\left\{ {\begin{array}{*{20}{l}}
  x\left( {1 + x} \right){u_{xx}} - ( {1 + x})z{u_{xz}} - {y^2}{u_{yy}} + yz{u_{yz}} \\
 \,\,\,\,\,\,\,\,\,  + \left[ {1 - {b_3} + \left( {{b_1} + {b_2} + 1} \right)x} \right]{u_x}
  - \left( {{b_1} - {b_2} + 1} \right)y{u_y} - {b_2}z{u_z} + {b_1}{b_2}u = 0,   \\
  y\left( {1 + y} \right){u_{yy}} - x ( {1 - y}){u_{xy}} - {z^2}{u_{zz}} + xz{u_{xz}} \\
  \,\,\,\,\,\,\,\,\, + \left[ {1 - {b_2} + \left( {a + {b_1} + 1} \right)y} \right]{u_y}
  - \left( {a - {b_1} + 1} \right)z{u_z} + ax{u_x} + a{b_1}u = 0,   \\
  z( {1 + z}){u_{zz}} - x( {1 + z}){u_{xz}} - y( {1 - z}){u_{yz}} - xy{u_{xy}} \\ \,\,\,\,\,\,\,\,\,
  + \left[ {1 - {b_1} + \left( {a + {b_3} + 1} \right)z} \right]{u_z} - ax{u_x} + {b_3}y{u_y} + a{b_3}u = 0,
\end{array}} \right.
$

where $u\equiv \,\,   {F_{10k}}\left( {a,{b_1},{b_2},{b_3};x,y,z} \right) $.

\bigskip

\begin{equation} \label{eqf1c}
{F_{10l}}\left( {{a_1},{a_2},{b_1},{b_2};c;x,y,z} \right) \hfill \\
  = \sum\limits_{m,n,p = 0}^\infty  {} \frac{{{{\left( {{a_1}} \right)}_{n + p}}{{\left( {{a_2}} \right)}_m}{{\left( {{b_1}} \right)}_{m + n - p}}{{\left( {{b_2}} \right)}_{p - n}}}}{{{{\left( c \right)}_m}}}\frac{x^m}{m!}\frac{y^n}{n!}\frac{z^p}{p!},
\end{equation}

region of convergence:
$$
\left\{ {s + t < 1,\,\,\,r < \frac{{\sqrt {1 - 4st}  - \left| {1 - 2t} \right|}}{{2t}}} \right\}.
$$

System of partial differential equations:

$
\left\{ {\begin{array}{*{20}{l}}
  x ( {1 - x}){u_{xx}} - xy{u_{xy}} + xz{u_{xz}} + \left[ {c - \left( {{a_2} + {b_1} + 1} \right)x} \right]{u_x}
   - {a_2}y{u_y} + {a_2}z{u_z} - {a_2}{b_1}u = 0,  \\
  y( {1 + y}){u_{yy}} + xy{u_{xy}} + xz{u_{xz}} - z{u_{yz}} - {z^2}{u_{zz}}
   + \left[ {1 - {b_2} + \left( {{a_1} + {b_1} + 1} \right)y} \right]{u_y} \hfill \\
   \,\,\,\,\,\,\,\,\, + {a_1}x{u_x} - \left( {{a_1} - {b_1} + 1} \right)z{u_z} + {a_1}{b_1}u = 0,   \\
  z( {1 + z}){u_{zz}} - x{u_{xz}} - {y^2}{u_{yy}} - y{u_{yz}}\\\,\,\,\,\,\,\,\,\,\,
  + \left[ {1 - {b_1} + \left( {{a_1} + {b_2} + 1} \right)z} \right]{u_z}  - \left( {{a_1} - {b_2} + 1} \right)y{u_y} + {a_1}{b_2}u = 0,
\end{array}} \right.
$

where $u\equiv \,\,   {F_{10l}}\left( {{a_1},{a_2},{b_1},{b_2};c;x,y,z} \right)$.

Particular solutions:

$
{u_1} = {F_{10l}}\left( {{a_1},{a_2},{b_1},{b_2};c;x,y,z} \right),
$

$
{u_2} = {x^{1 - c}}{F_{10l}}\left( {{a_1},1 - c + {a_2},1 - c + {b_1},{b_2};2 - c;x,y,z} \right).
$

\bigskip

\begin{equation}
{F_{11a}}\left( {{a_1},{a_2},{a_3},{a_4};{c_1},{c_2};x,y,z} \right)\hfill \\
  = \sum\limits_{m,n,p = 0}^\infty  {} \frac{{{{\left( {{a_1}} \right)}_{m + n}}{{\left( {{a_2}} \right)}_{n + p}}{{\left( {{a_3}} \right)}_m}{{\left( {{a_4}} \right)}_p}}}{{{{\left( {{c_1}} \right)}_{m + n}}{{\left( {{c_2}} \right)}_p}}}\frac{x^m}{m!}\frac{y^n}{n!}\frac{z^p}{p!},
\end{equation}

first appearance of this function in the literature, and old notation: [24], $F_{11}$,\,\,[33], $F_M$,

region of convergence:
$$
\left\{ {r < 1,\,\,\,s + t < 1} \right\}.
$$

System of partial differential equations:

$
\left\{ {\begin{array}{*{20}{l}}
  x( {1 - x}){u_{xx}} + ( {1 - x})y{u_{xy}} + \left[ {{c_1} - \left( {{a_1} + {a_3} + 1} \right)x} \right]{u_x}
   - {a_3}y{u_y} - {a_1}{a_3}u = 0, \\
  y( {1 - y}){u_{yy}} + x( {1 - y}){u_{xy}} - xz{u_{xz}} - yz{u_{yz}}\\\,\,\,\,\,\,\,\,\,\,
  + \left[ {{c_1} - \left( {{a_1} + {a_2} + 1} \right)y} \right]{u_y} - {a_2}x{u_x} - {a_1}z{u_z} - {a_1}{a_2}u = 0,\\
  z( {1 - z}){u_{zz}} - yz{u_{yz}} + \left[ {{c_2} - \left( {{a_2} + {a_4} + 1} \right)z} \right]{u_z}
  - {a_4}y{u_y} - {a_2}{a_4}u = 0,
\end{array}} \right.
$

where $u\equiv \,\,   {F_{11a}}\left( {{a_1},{a_2},{a_3},{a_4};{c_1},{c_2};x,y,z} \right)$.

Particular solutions:

$
{u_1} = {F_{11a}}\left( {{a_1},{a_2},{a_3},{a_4};{c_1},{c_2};x,y,z} \right),
$

$
{u_2} = {z^{1 - c_2}}{F_{11a}}\left( {{a_1},1 - c_2 + {a_2},{a_3},1 - c_2 + {a_4};{c_1},2 - {c_2};x,y,z} \right).
$

\bigskip

\begin{equation}
{F_{11b}}\left( {{a_1},{a_2},{a_3},b;c;x,y,z} \right) \hfill \\
  = \sum\limits_{m,n,p = 0}^\infty  {} \frac{{{{\left( {{a_1}} \right)}_{n + p}}{{\left( {{a_2}} \right)}_m}{{\left( {{a_3}} \right)}_p}{{\left( b \right)}_{m + n - p}}}}{{{{\left( c \right)}_{m + n}}}}\frac{x^m}{m!}\frac{y^n}{n!}\frac{z^p}{p!},
\end{equation}

region of convergence:
$$
 \left\{ {s < 1,\,\,\,t + 2\sqrt {st}  < 1,\,\,\,r < \min \left\{ {1,\frac{1-t+{\sqrt {{{\left( {1 - t} \right)}^2} - 4st} }}{{2t}}} \right\}} \right\}.
$$

System of partial differential equations:

$
\left\{ {\begin{array}{*{20}{l}}
  x( {1 - x}){u_{xx}} + ( {1 - x})y{u_{xy}} + xz{u_{xz}}
  + \left[ {c - \left( {{a_2} + b + 1} \right)x} \right]{u_x} - {a_2}y{u_y}   + {a_2}z{u_z} - {a_2}bu = 0,  \\
  y( {1 - y}){u_{yy}} + x( {1 - y}){u_{xy}} - xz{u_{xz}} + {z^2}{u_{zz}} \\ \,\,\,\,\,\,\,\,\,
   + \left[ {c - \left( {{a_1} + b + 1} \right)y} \right]{u_y} - {a_1}x{u_x}  + \left( {{a_1} - b + 1} \right)z{u_z} - {a_1}bu = 0,  \\
  z( {1 + z}){u_{zz}} - x{u_{xz}}  -y( {1 - z}){u_{yz}} + {a_3}y{u_y}
   + \left[ {1 - b + \left( {{a_1} + {a_3} + 1} \right)z} \right]{u_z}    + {a_1}{a_3}u = 0,
\end{array}} \right.
$

where $u\equiv \,\,   {F_{11b}}\left( {{a_1},{a_2},{a_3},b;c;x,y,z} \right)$.

\bigskip

\begin{equation}
{F_{11c}}\left( {{a_1},{a_2},{a_3},b;c;x,y,z} \right) \hfill \\
  = \sum\limits_{m,n,p = 0}^\infty  {} \frac{{{{\left( {{a_1}} \right)}_{m + n}}{{\left( {{a_2}} \right)}_{n + p}}{{\left( {{a_3}} \right)}_p}{{\left( b \right)}_{m - p}}}}{{{{\left( c \right)}_{m + n}}}}\frac{x^m}{m!}\frac{y^n}{n!}\frac{z^p}{p!},
\end{equation}

region of convergence:
$$
\left\{ {r < 1,\,\,\,s < 1,\,\,\,t < \frac{{1 - s}}{{1 + r}}} \right\}.
$$

System of partial differential equations:

$
\left\{ {\begin{array}{*{20}{l}}
  x( {1 - x}){u_{xx}} + ( {1 - x})y{u_{xy}} + xz{u_{xz}} + yz{u_{yz}}\\\,\,\,\,\,\,\,\,\,\,
   + \left[ {c - \left( {{a_1} + b + 1} \right)x} \right]{u_x} - by{u_y}   + {a_1}z{u_z} - {a_1}bu = 0, \\
  y( {1 - y}){u_{yy}} + x( {1 - y}){u_{xy}} - xz{u_{xz}} - yz{u_{yz}}\\\,\,\,\,\,\,\,\,\,\,
   + \left[ {c - \left( {{a_1} + {a_2} + 1} \right)y} \right]{u_y} - {a_2}x{u_x}  - {a_1}z{u_z} - {a_1}{a_2}u = 0,   \\
  z( {1 + z}){u_{zz}} - x{u_{xz}} + yz{u_{yz}}
   + \left[ {1 - b + \left( {{a_2} + {a_3} + 1} \right)z} \right]{u_z} + {a_3}y{u_y} + {a_2}{a_3}u = 0,
\end{array}} \right.
$

where $u\equiv \,\,   {F_{11c}}\left( {{a_1},{a_2},{a_3},b;c;x,y,z} \right)$.

\bigskip

\begin{equation}
{F_{11d}}\left( {{a_1},{a_2},{a_3},b;c;x,y,z} \right)\hfill \\
   = \sum\limits_{m,n,p = 0}^\infty  {} \frac{{{{\left( {{a_1}} \right)}_{m + n}}{{\left( {{a_2}} \right)}_{n + p}}{{\left( {{a_3}} \right)}_m}{{\left( b \right)}_{p - m - n}}}}{{{{\left( c \right)}_p}}}{x^m}\frac{x^m}{m!}\frac{y^n}{n!}\frac{z^p}{p!},
\end{equation}

region of convergence:
$$
\left\{ {r < 1,\,\,\,s + 2\sqrt {st}  < 1,\,\,\,r <  {\frac{{1 + s + \sqrt {{{\left( {1 - s} \right)}^2} - 4st} }}{{2(1+t)}}} } \right\}.
$$

System of partial differential equations:

$
\left\{ {\begin{array}{*{20}{l}}
  x( {1 + x}){u_{xx}} + \left( {1 + x} \right)y{u_{xy}} - z{u_{xz}}
   + \left[ {1 - b + \left( {{a_1} + {a_3} + 1} \right)x} \right]{u_x}  + {a_3}y{u_y} + {a_1}{a_3}u = 0,  \\
  y( {1 + y}){u_{yy}} + x( {1 + y}){u_{xy}} + xz{u_{xz}} - ( {1 - y})z{u_{yz}} \hfill \\ \,\,\,\,\,\,\,\,\,
   + \left[ {1 - b + \left( {{a_1} + {a_2} + 1} \right)y} \right]{u_y}  + {a_2}x{u_x} + {a_1}z{u_z} + {a_1}{a_2}u = 0 ,  \\
  z( {1 - z}){u_{zz}} + xy{u_{xy}} + xz{u_{xz}} + {y^2}{u_{yy}}\\\,\,\,\,\,\,\,\,\,\,
  + \left[ {c - \left( {{a_2} + b + 1} \right)z} \right]{u_z} + {a_2}x{u_x} + \left( {{a_2} - b + 1} \right)y{u_y} - {a_2}bu = 0,
\end{array}} \right.
$

where $u\equiv \,\,   {F_{11d}}\left( {{a_1},{a_2},{a_3},b;c;x,y,z} \right)$.

Particular solutions:

$
{u_1} = {F_{11d}}\left( {{a_1},{a_2},{a_3},b;c;x,y,z} \right),
$

$
{u_2} = {z^{1 - c}}{F_{11d}}\left( {{a_1},1 - c + {a_2},{a_3},1 - c + b;2 - c;x,y,z} \right).
$

\bigskip

\begin{equation}
{F_{11e}}\left( {{a_1},{a_2},{b_1},{b_2};x,y,z} \right) \hfill \\
  = \sum\limits_{m,n,p = 0}^\infty  {} {{{{\left( {{a_1}} \right)}_{n + p}}{{\left( {{a_2}} \right)}_m}{{\left( {{b_1}} \right)}_{m + n - p}}{{\left( {{b_2}} \right)}_{p - m - n}}}}\frac{x^m}{m!}\frac{y^n}{n!}\frac{z^p}{p!},
\end{equation}

first appearance of this function in the literature, and old notation: [28], $G_{A}$,

region of convergence:
$$
 \left\{ {r < 1,\,\,\,s + t < 1} \right\}.
$$

System of partial differential equations:

$
\left\{ {\begin{array}{*{20}{l}}
  x( {1 + x}){u_{xx}} + (1+x)y{u_{xy}} - (1+x)z{u_{xz}}\\\,\,\,\,\,\,\,\,\,\,
   + \left[ {1 - {b_2} + \left( {{a_2} + {b_1} + 1} \right)x} \right]{u_x}  + {a_2}y{u_y} - {a_2}z{u_z} + {a_2}{b_1}u = 0, \\
  y\left( {1 + y} \right){u_{yy}} + x\left( {1 + y} \right){u_{xy}} + xz{u_{xz}} - z{u_{yz}} - {z^2}{u_{zz}}\hfill \\ \,\,\,\,\,\,\,\,\,
   + \left[ {1 - {b_2} + \left( {{a_1} + {b_1} + 1} \right)y} \right]{u_y}  + {a_1}x{u_x} -
    \left( {{a_1} -b_1 +1} \right)z{u_z} + {a_1}{b_1}u = 0, \\
  z\left( {1 + z} \right){u_{zz}} - xy{u_{xy}} - x\left( {1 + z} \right){u_{xz}} - y{u_{yz}} - {y^2}{u_{yy}} \hfill \\ \,\,\,\,\,\,\,\,\,
  + \left[ {1 - {b_1} + \left( {{a_1} + {b_2} + 1} \right)z} \right]{u_z}   - {a_1}x{u_x} - \left( {{a_1} - {b_2} + 1} \right)y{u_y} + {a_1}{b_2}u = 0,
\end{array}} \right.
$

where $u\equiv \,\,   {F_{11e}}\left( {{a_1},{a_2},{b_1},{b_2};x,y,z} \right)$.

\bigskip

\begin{equation}
{F_{11f}}\left( {{a_1},{a_2},{b_1},{b_2};x,y,z} \right) \hfill \\
  = \sum\limits_{m,n,p = 0}^\infty  {} {{{{\left( {{a_1}} \right)}_{m + n}}{{\left( {{a_2}} \right)}_{n + p}}{{\left( {{b_1}} \right)}_{p - m - n}}{{\left( {{b_2}} \right)}_{m - p}}}}\frac{x^m}{m!}\frac{y^n}{n!}\frac{z^p}{p!},
\end{equation}

first appearance of this function in the literature, and old notation: [38], $G_{C}$,

region of convergence:
$$
 \left\{ {r < 1,\,\,\,s + 2\sqrt {st}  < 1,\,\,\,r < \min \left\{ {1,\frac{{1 + s - 2\sqrt {s + st} }}{t}} \right\}} \right\}.
$$

System of partial differential equations:

$
\left\{ {\begin{array}{*{20}{l}}
  x\left( {1 + x} \right){u_{xx}} + \left( {1 + x} \right)y{u_{xy}} - \left( {1 + x} \right)z{u_{xz}} - yz{u_{yz}} \hfill \\
   \,\,\,\,\,\,\,\,\, + \left[ {1 - {b_1} + \left( {{a_1} + {b_2} + 1} \right)x} \right]{u_x}   + {b_2}y{u_y} - {a_1}z{u_z} + {a_1}{b_2}u = 0,  \\
  y\left( {1 + y} \right){u_{yy}} + x\left( {1 + y} \right){u_{xy}} + xz{u_{xz}} - \left( {1 - y} \right)z{u_{yz}} \hfill \\
   \,\,\,\,\,\,\,\,\, + \left[ {1 - {b_1} + \left( {{a_1} + {a_2} + 1} \right)y} \right]{u_y}  + {a_2}x{u_x} + {a_1}z{u_z} + {a_1}{a_2}u = 0,\\
  z\left( {1 + z} \right){u_{zz}} - xy{u_{xy}} - x\left( {1 + z} \right){u_{xz}} - {y^2}{u_{yy}} \hfill \\ \,\,\,\,\,\,\,\,\,
   + \left[ {1 - {b_2} + \left( {{a_2} + {b_1} + 1} \right)z} \right]{u_z} - {a_2}x{u_x}
   - \left( {{a_2} - {b_1} + 1} \right)y{u_y} + {a_2}{b_1}u = 0,
\end{array}} \right.
$

where $u\equiv \,\,   {F_{11f}}\left( {{a_1},{a_2},{b_1},{b_2};x,y,z} \right)$.

\bigskip

\begin{equation}
{F_{12a}}\left( {{a_1},{a_2},{a_3},{a_4};{c_1},{c_2};x,y,z} \right) \hfill \\
  = \sum\limits_{m,n,p = 0}^\infty  {} \frac{{{{\left( {{a_1}} \right)}_{m + n}}{{\left( {{a_2}} \right)}_{n + p}}{{\left( {{a_3}} \right)}_m}{{\left( {{a_4}} \right)}_p}}}{{{{\left( {{c_1}} \right)}_{m + p}}{{\left( {{c_2}} \right)}_n}}}\frac{x^m}{m!}\frac{y^n}{n!}\frac{z^p}{p!},
\end{equation}

first appearance of this function in the literature, and old notation: [24], $F_{12}$,\,\,[33], $F_P$,

region of convergence:
$$
 \left\{ {r + s < 1,\,\,\,s + t < 1} \right\}.
$$

System of partial differential equations:

$
\left\{ {\begin{array}{*{20}{l}}
  x\left( {1 - x} \right){u_{xx}} - xy{u_{xy}} + z{u_{xz}} + \left[ {{c_1} - \left( {{a_1} + {a_3} + 1} \right)x} \right]{u_x}
  - {a_3}y{u_y} - {a_1}{a_3}u = 0, \\
  y\left( {1 - y} \right){u_{yy}} - xy{u_{xy}} - xz{u_{xz}} - yz{u_{yz}}\\\,\,\,\,\,\,\,\,\,\,
  + \left[ {{c_2} - \left( {{a_1} + {a_2} + 1} \right)y} \right]{u_y} - {a_2}x{u_x} - {a_1}z{u_z}   - {a_1}{a_2}u = 0, \\
  z\left( {1 - z} \right){u_{zz}} + x{u_{xz}} - yz{u_{yz}} + \left[ {{c_1} - \left( {{a_2} + {a_4} + 1} \right)z} \right]{u_z}
   - {a_4}y{u_y} - {a_2}{a_4}u = 0,
\end{array}} \right.
$

where $u\equiv \,\,   {F_{12a}}\left( {{a_1},{a_2},{a_3},{a_4};{c_1},{c_2};x,y,z} \right)$.

Particular solutions:

$
{u_1} = {F_{12a}}\left( {{a_1},{a_2},{a_3},{a_4};{c_1},{c_2};x,y,z} \right),
$

$
{u_2} = {y^{1 - {c_2}}}{F_{12a}}\left( {1 - {c_2} + {a_1},1 - {c_2} + {a_2},{a_3},{a_4};{c_1},2 - {c_2};x,y,z} \right).
$

\bigskip

\begin{equation}
 {F_{12b}}\left( {{a_1},{a_2},{a_3},b;c;x,y,z} \right) \hfill \\
  = \sum\limits_{m,n,p = 0}^\infty  {} \frac{{{{\left( {{a_1}} \right)}_{m + n}}{{\left( {{a_2}} \right)}_{n + p}}{{\left( {{a_3}} \right)}_p}{{\left( b \right)}_{m - n}}}}{{{{\left( c \right)}_{m + p}}}}\frac{x^m}{m!}\frac{y^n}{n!}\frac{z^p}{p!},
\end{equation}

region of convergence:
$$
  \left\{ {r < 1,\,\,\,s + 2\sqrt {rs}  < 1,\,\,\,t <  {\frac{{1 - s + 2r + \sqrt {{{\left( {1 - s} \right)}^2} - 4rs} }}{{2\left( {1 + r} \right)}}} } \right\}.
$$

System of partial differential equations:

$
\left\{ {\begin{array}{*{20}{l}}
  x\left( {1 - x} \right){u_{xx}} + z{u_{xz}} + {y^2}{u_{yy}} + \left[ {c - \left( {{a_1} + b + 1} \right)x} \right]{u_x}
   + \left( {{a_1} - b + 1} \right)y{u_y} - {a_1}bu = 0, \\
  y\left( {1 + y} \right){u_{yy}} - x\left( {1 - y} \right){u_{xy}} + xz{u_{xz}} + yz{u_{yz}} \hfill \\
  \,\,\,\,\,\,\,\,\,  + \left[ {1 - b + \left( {{a_1} + {a_2} + 1} \right)y} \right]{u_y} + {a_2}x{u_x} + {a_1}z{u_z} + {a_1}{a_2}u = 0, \\
  z\left( {1 - z} \right){u_{zz}} + x{u_{xz}} - yz{u_{yz}} + \left[ {c - \left( {{a_2} + {a_3} + 1} \right)z} \right]{u_z}
  - {a_3}y{u_y} - {a_2}{a_3}u = 0,
\end{array}} \right.
$

where $u\equiv \,\,   {F_{12b}}\left( {{a_1},{a_2},{a_3},b;c;x,y,z} \right)$.

\bigskip

\begin{equation}
{F_{13a}}\left( {{a_1},{a_2},{a_3},{a_4};c;x,y,z} \right) \hfill \\
  = \sum\limits_{m,n,p = 0}^\infty  {} \frac{{{{\left( {{a_1}} \right)}_{m + n}}{{\left( {{a_2}} \right)}_{n + p}}{{\left( {{a_3}} \right)}_m}{{\left( {{a_4}} \right)}_p}}}{{{{\left( c \right)}_{m + n + p}}}}\frac{x^m}{m!}\frac{y^n}{n!}\frac{z^p}{p!},
\end{equation}

first appearance of this function in the literature, and old notation: [24], $F_{13}$,\,\,[33], $F_T$,

region of convergence:
$$
\left\{ {r < 1,\,\,\,s < 1,\,\,\,t < 1} \right\}.
$$

System of partial differential equations:

$
\left\{ {\begin{array}{*{20}{l}}
  x\left( {1 - x} \right){u_{xx}} + \left( {1 - x} \right)y{u_{xy}} + z{u_{xz}} + \left[ {c - \left( {{a_1} + {a_3} + 1} \right)x} \right]{u_x}
   - {a_3}y{u_y} - {a_1}{a_3}u = 0, \hfill \\
  y\left( {1 - y} \right){u_{yy}} + x\left( {1 - y} \right){u_{xy}} - xz{u_{xz}} + \left( {1 - y} \right)z{u_{yz}}\\  \,\,\,\,\,\,\,\,\,
   + \left[ {c - \left( {{a_1} + {a_2} + 1} \right)y} \right]{u_y}  - {a_2}x{u_x} - {a_1}z{u_z} - {a_1}{a_2}u = 0,  \\
  z\left( {1 - z} \right){u_{zz}} + x{u_{xz}} + y\left( {1 - z} \right){u_{yz}} + \left[ {c - \left( {{a_2} + {a_4} + 1} \right)z} \right]{u_z}
   - {a_4}y{u_y} - {a_2}{a_4}u = 0,
\end{array}} \right.
$

where $u\equiv \,\,   {F_{13a}}\left( {{a_1},{a_2},{a_3},{a_4};c;x,y,z} \right)$.

\bigskip

\begin{equation}
 {F_{14a}}\left( {{a_1},{a_2},{a_3};{c_1},{c_2},{c_3};x,y,z} \right) \hfill \\
  = \sum\limits_{m,n,p = 0}^\infty  {} \frac{{{{\left( {{a_1}} \right)}_{m + n}}{{\left( {{a_2}} \right)}_{n + p}}{{\left( {{a_3}} \right)}_{p + m}}}}{{{{\left( {{c_1}} \right)}_m}{{\left( {{c_2}} \right)}_n}{{\left( {{c_3}} \right)}_p}}}\frac{x^m}{m!}\frac{y^n}{n!}\frac{z^p}{p!},
\end{equation}

first appearance of this function in the literature, and old notation: [36], $H_{B}$,

region of convergence:
$$
\left\{ {r + s + t + 2\sqrt {rst}  < 1} \right\}.
$$

System of partial differential equations:

$
\left\{ {\begin{array}{*{20}{l}}
  x\left( {1 - x} \right){u_{xx}} - xy{u_{xy}} - xz{u_{xz}} - yz{u_{yz}}\\\,\,\,\,\,\,\,\,\,\,
   + \left[ {{c_1} - \left( {{a_1} + {a_3} + 1} \right)x} \right]{u_x} - {a_3}y{u_y} - {a_1}z{u_z}   - {a_1}{a_3}u = 0,   \\
  y\left( {1 - y} \right){u_{yy}} - xy{u_{xy}} - xz{u_{xz}} - yz{u_{yz}}\\\,\,\,\,\,\,\,\,\,\,
   + \left[ {{c_2} - \left( {{a_1} + {a_2} + 1} \right)y} \right]{u_y} - {a_2}x{u_x} - {a_1}z{u_z}   - {a_1}{a_2}u = 0,  \\
  z\left( {1 - z} \right){u_{zz}} - xy{u_{xy}} - xz{u_{xz}} - yz{u_{yz}}\\\,\,\,\,\,\,\,\,\,\,
   + \left[ {{c_3} - \left( {{a_2} + {a_3} + 1} \right)z} \right]{u_z} - {a_2}x{u_x} - {a_3}y{u_y}   - {a_2}{a_3}u = 0,
\end{array}} \right.
$

where $u\equiv \,\,   {F_{14a}}\left( {{a_1},{a_2},{a_3};{c_1},{c_2},{c_3};x,y,z} \right)$.

Particular solutions:

$
{u_1} = {F_{14a}}\left( {{a_1},{a_2},{a_3};{c_1},{c_2},{c_3};x,y,z} \right),
$

$
{u_2} = {x^{1 - {c_1}}}{F_{14a}}\left( {1 - {c_1} + {a_1},{a_2},1 - {c_1} + {a_3};2 - {c_1},{c_2},{c_3};x,y,z} \right),
$

$
{u_3} = {y^{1 - {c_2}}}{F_{14a}}\left( {1 - {c_2} + {a_1},1 - {c_2} + {a_2},{a_3};{c_1},2 - {c_2},{c_3};x,y,z} \right),
$

$
{u_4} = {z^{1 - {c_3}}}{F_{14a}}\left( {{a_1},1 - {c_3} + {a_2},1 - {c_3} + {a_3};{c_1},{c_2},2 - {c_3};x,y,z} \right),
$

$
{u_5} = {x^{1 - {c_1}}}{y^{1 - {c_2}}}{F_{14a}}\left( {2 - {c_1} - {c_2} + {a_1},1 - {c_2} + {a_2},1 - {c_1} + {a_3};2 - {c_1},2 - {c_2},{c_3};x,y,z} \right),
$

$
{u_6} = {y^{1 - {c_2}}}{z^{1 - {c_3}}}{F_{14a}}\left( {1 - {c_2} + {a_1},2 - {c_2} - {c_3} + {a_2},1 - {c_3} + {a_3};{c_1},{c_2},{c_3};x,y,z} \right),
$

$
{u_7} = {x^{1 - {c_1}}}{z^{1 - {c_3}}}{F_{14a}}\left( {1 - {c_1} + {a_1},1 - {c_3} + {a_2},2 - {c_1} - {c_3} + {a_3};2 - {c_1},{c_2},2 - {c_3};x,y,z} \right),
$

$
  {u_8} = {x^{1 - {c_1}}}{y^{1 - {c_2}}}{z^{1 - {c_3}}}\times$
    
    $\,\,\,\,\,\,\,\,\,\,\times{F_{14a}}\left( {2 - {c_1} - {c_2} + {a_1},2 - {c_2} - {c_3} + {a_2},2 - {c_1} - {c_3} + {a_3};}  {2 - {c_1},2 - {c_2},2 - {c_3};x,y,z} \right).
  $

\bigskip

\begin{equation}
{F_{14b}}\left( {{a_1},{a_2},b;{c_1},{c_2};x,y,z} \right)\hfill \\
   = \sum\limits_{m,n,p = 0}^\infty  {} \frac{{{{\left( {{a_1}} \right)}_{n + p}}{{\left( {{a_2}} \right)}_{m + p}}{{\left( b \right)}_{m + n - p}}}}{{{{\left( {{c_1}} \right)}_m}{{\left( {{c_2}} \right)}_n}}}\frac{x^m}{m!}\frac{y^n}{n!}\frac{z^p}{p!},
\end{equation}

region of convergence:
$$
\begin{gathered}
  \left\{ {r + s < 1,\,\,\,t < \min \left\{ {U\left( {{w_1}} \right),U\left( {{w_2}} \right)} \right\}} \right\}, \hfill \\
  P\left( w \right) =  - w + \frac{1}{2}r + \frac{1}{2}s + \frac{1}{2}\sqrt {{r^2} + 4wr}  + \frac{1}{2}\sqrt {{s^2} + 4ws} , \hfill \\
  {w_1}:{\rm{the\,positive\,root\,of}}\,P\left( w \right) = -1, \hfill \\
  {w_2}:{\rm{the\,smaller\,positive\,root\,of}}\,\,P\left( w \right) = 1, \hfill \\
  U\left( w \right) = \frac{{16w}}{{{{\left( {\sqrt r  + \sqrt {r + 4w} } \right)}^2}{{\left( {\sqrt s  + \sqrt {s + 4w} } \right)}^2}}}. \hfill \\
\end{gathered}
$$

System of partial differential equations:

$
\left\{ {\begin{array}{*{20}{l}}
  x\left( {1 - x} \right){u_{xx}} - xy{u_{xy}} - yz{u_{yz}} + {z^2}{u_{zz}}  + \left[ {{c_1} - \left( {{a_2} + b + 1} \right)x} \right]{u_x}  - {a_2}y{u_y}\\\,\,\,\,\,\,\,\,\,\, + \left( {{a_2} - b + 1} \right)z{u_z} - {a_2}bu = 0 , \\
  y\left( {1 - y} \right){u_{yy}} - xy{u_{xy}} - xz{u_{xz}} + {z^2}{u_{zz}}   + \left[ {{c_2} - \left( {{a_1} + b + 1} \right)y} \right]{u_y}  - {a_1}x{u_x} \\\,\,\,\,\,\,\,\,\,\, + \left( {{a_1} - b + 1} \right)z{u_z} - {a_1}bu = 0 , \\
  z\left( {1 + z} \right){u_{zz}} + xy{u_{xy}} - x\left( {1 - z} \right){u_{xz}} - y\left( {1 - z} \right){u_{yz}} \\ \,\,\,\,\,\,\,\,\,
  + \left[ {1 - b + \left( {{a_1} + {a_2} + 1} \right)z} \right]{u_z}  + {a_1}x{u_x} + {a_2}y{u_y} + {a_1}{a_2}u = 0,
\end{array}} \right.
$

where $u\equiv \,\,   {F_{14b}}\left( {{a_1},{a_2},b;{c_1},{c_2};x,y,z} \right)$.

Particular solutions:

$
{u_1} = {F_{14b}}\left( {{a_1},{a_2},b;{c_1},{c_2};x,y,z} \right),
$

$
{u_2} = {x^{1 - {c_1}}}{F_{14b}}\left( {{a_1},1 - {c_1} + {a_2},1 - {c_1} + b;2 - {c_1},{c_2};x,y,z} \right),
$

$
{u_3} = {y^{1 - {c_2}}}{F_{14b}}\left( {1 - {c_2} + {a_1},{a_2},1 - {c_2} + b;{c_1},2 - {c_2};x,y,z} \right),
$

$
{u_4} = {x^{1 - {c_1}}}{y^{1 - {c_2}}}{F_{14b}}\left( {1 - {c_2} + {a_1},1 - {c_1} + {a_2},2 - {c_1} - {c_2} + b;2 - {c_1},2 - {c_2};x,y,z} \right).
$

\bigskip

\begin{equation}
{F_{14c}}\left( {{a},{b_1},{b_2};c;x,y,z} \right) \hfill \\
  = \sum\limits_{m,n,p = 0}^\infty  {} \frac{{{{\left( {{a}} \right)}_{m + p}}{{\left( {{b_1}} \right)}_{m + n - p}}{{\left( {{b_2}} \right)}_{n + p - m}}}}{{{{\left( c \right)}_n}}}\frac{x^m}{m!}\frac{y^n}{n!}\frac{z^p}{p!},
\end{equation}

region of convergence:
$$
\left\{ {r + t < 1,\,\,\,s < \min \left\{ {1,\frac{{{{\left( {1 - r - t} \right)}^2}}}{{2\left( {r + t - 4rt + \left| {r - t} \right|\sqrt {1 - 4rt} } \right)}}} \right\}} \right\}.
$$

System of partial differential equations:

$
\left\{ {\begin{array}{*{20}{l}}
  x\left( {1 + x} \right){u_{xx}} - \left( {1 - x} \right)y{u_{xy}} - z{u_{xz}} + yz{u_{yz}} - {z^2}{u_{zz}}\hfill \\
 \,\,\,\,\,\,\,\,\,  + \left[ {1 - {b_2} + \left( {{a} + {b_1} + 1} \right)x} \right]{u_x}    + {a}y{u_y}
   - \left( {{a} - {b_1} + 1} \right)z{u_z} + {a}{b_1}u = 0, \\
  y\left( {1 - y} \right){u_{yy}} - 2xz{u_{xz}} + {x^2}{u_{xx}} + {z^2}{u_{zz}} \hfill \\
 \,\,\,\,\,\,\,\,\,  + \left[ {c - \left( {{b_1} + {b_2} + 1} \right)y} \right]{u_y}
   + \left( {{b_1} - {b_2} + 1} \right)x{u_x} + \left( {{b_2} - {b_1} + 1} \right)z{u_z} - {b_1}{b_2}u = 0,  \\
  z\left( {1 + z} \right){u_{zz}} + xy{u_{xy}} - x{u_{xz}} - y\left( {1 - z} \right){u_{yz}} - {x^2}{u_{xx}}\hfill \\
 \,\,\,\,\,\,\,\,\,  + \left[ {1 - {b_1} + \left( {{a} + {b_2} + 1} \right)z} \right]{u_z}
   - \left( {{a} - {b_2} + 1} \right)x{u_x} + {a}y{u_y} + {a}{b_2}u = 0, \
\end{array}} \right.
$

where $u\equiv \,\,   {F_{14c}}\left( {{a},{b_1},{b_2};c;x,y,z} \right)$.

Particular solutions:

$
{u_1} = {F_{14c}}\left( {{a},{b_1},{b_2};c;x,y,z} \right),
$

$
{u_2} = {y^{1 - c}}{F_{14c}}\left( {{a},1 - c + {b_1},1 - c + {b_2};2 - c;x,y,z} \right).
$

\bigskip

\begin{equation}
{F_{14d}}\left( {{b_1},{b_2},{b_3};x,y,z} \right) \hfill \\
  = \sum\limits_{m,n,p = 0}^\infty  {} {{{{\left( {{b_1}} \right)}_{m + n - p}}{{\left( {{b_2}} \right)}_{n + p - m}}{{\left( {{b_3}} \right)}_{p + m - n}}}}\frac{x^m}{m!}\frac{y^n}{n!}\frac{z^p}{p!},
\end{equation}

region of convergence:
$$
 \left\{ {r + s + t - 4rst < 1} \right\}.
$$

System of partial differential equations:

$
\left\{ {\begin{array}{*{20}{l}}
  x\left( {1 + x} \right){u_{xx}} - y{u_{xy}} - z{u_{xz}} + 2yz{u_{yz}} - {y^2}{u_{yy}} - {z^2}{u_{zz}} \\
 \,\,\,\,\,\,\,\,\,  + \left[ {1 - {b_2} + \left( {{b_1} + {b_3} + 1} \right)x} \right]{u_x}
   - \left( {{b_1} - {b_3} + 1} \right)y{u_y} + \left( {{b_1} - {b_3} - 1} \right)z{u_z} + {b_1}{b_3}u = 0 , \\
  y\left( {1 + y} \right){u_{yy}} - x{u_{xy}} - z{u_{yz}} + 2xz{u_{xz}} - {x^2}{u_{xx}} - {z^2}{u_{zz}} \\
 \,\,\,\,\,\,\,\,\,  + \left[ {1 - {b_3} + \left( {{b_1} + {b_2} + 1} \right)y} \right]{u_y}
   - \left( {{b_1} - {b_2} + 1} \right)x{u_x} + \left( {{b_1} - {b_2} - 1} \right)z{u_z} + {b_1}{b_2}u = 0,
 \\
  z\left( {1 + z} \right){u_{zz}} + 2xy{u_{xy}} - x{u_{xz}} - y{u_{yz}} - {x^2}{u_{xx}} - {y^2}{u_{yy}} \\
  \,\,\,\,\,\,\,\,\, + \left[ {1 - {b_1} + \left( {{b_2} + {b_3} + 1} \right)z} \right]{u_z}
   + \left( {{b_2} - {b_3} - 1} \right)x{u_x} - \left( {{b_2} - {b_3} + 1} \right)y{u_y} + {b_2}{b_3}u = 0,
\end{array}} \right.
$

where $u\equiv \,\,   {F_{14d}}\left( {{b_1},{b_2},{b_3};x,y,z} \right)$.

\bigskip

\begin{equation}
{F_{15a}}\left( {{a_1},{a_2},{a_3};{c_1},{c_2};x,y,z} \right)\hfill \\
   = \sum\limits_{m,n,p = 0}^\infty  {} \frac{{{{\left( {{a_1}} \right)}_{m + n}}{{\left( {{a_2}} \right)}_{n + p}}{{\left( {{a_3}} \right)}_{p + m}}}}{{{{\left( {{c_1}} \right)}_{m + n}}{{\left( {{c_2}} \right)}_p}}}\frac{x^m}{m!}\frac{y^n}{n!}\frac{z^p}{p!},
\end{equation}

first appearance of this function in the literature, and old notation: [36], $H_{A}$,

region of convergence:
$$
 \left\{ {r < 1,\,\,\,s < 1,\,\,\,t < \left( {1 - r} \right)\left( {1 - s} \right)} \right\}.
$$

System of partial differential equations:

$
\left\{ {\begin{array}{*{20}{l}}
  \begin{gathered}
  x\left( {1 - x} \right){u_{xx}} + \left( {1 - x} \right)y{u_{xy}} - xz{u_{xz}} - yz{u_{yz}} \hfill\\\,\,\,\,\,\,\,\,\,\, + \left[ {{c_1} - \left( {{a_1} + {a_3} + 1} \right)x} \right]{u_x}   - {a_3}y{u_y} - {a_1}z{u_z} - {a_1}{a_3}u = 0 ,\hfill \\
\end{gathered}  \\
  \begin{gathered}
  y\left( {1 - y} \right){u_{yy}} + x\left( {1 - y} \right){u_{xy}} - xz{u_{xz}} - yz{u_{yz}} \hfill \\\,\,\,\,\,\,\,\,\,\, + \left[ {{c_1} - \left( {{a_1} + {a_2} + 1} \right)y} \right]{u_y}  - {a_2}x{u_x} - {a_1}z{u_z} - {a_1}{a_2}u = 0, \hfill \\
\end{gathered}  \\
  \begin{gathered}
  z\left( {1 - z} \right){u_{zz}} - xy{u_{xy}} - xz{u_{xz}} - yz{u_{yz}} \hfill \\ \,\,\,\,\,\,\,\,\,\, + \left[ {{c_2} - \left( {{a_2} + {a_3} + 1} \right)z} \right]{u_z}   - {a_2}x{u_x} - {a_3}y{u_y} - {a_2}{a_3}u = 0, \hfill \\
\end{gathered}
\end{array}} \right.
$

where $u\equiv \,\,   {F_{15a}}\left( {{a_1},{a_2},{a_3};{c_1},{c_2};x,y,z} \right)$.

Particular solutions:

$
{u_1} = {F_{15a}}\left( {{a_1},{a_2},{a_3};{c_1},{c_2};x,y,z} \right),
$

$
{u_2} = {z^{1 - {c_2}}}{F_{15a}}\left( {{a_1},1 - {c_2} + {a_2},1 - {c_2} + {a_3};{c_1},2 - {c_2};x,y,z} \right).
$

\bigskip

\begin{equation}
{F_{15b}}\left( {{a_1},{a_2},b;c;x,y,z} \right)  = \sum\limits_{m,n,p = 0}^\infty  {} \frac{{{{\left( {{a_1}} \right)}_{n + p}}{{\left( {{a_2}} \right)}_{m + p}}{{\left( b \right)}_{m + n - p}}}}{{{{\left( c \right)}_{m + n}}}}\frac{x^m}{m!}\frac{y^n}{n!}\frac{z^p}{p!},
\end{equation}

region of convergence:
$$
\begin{gathered}
  \left\{ {r < 1,\,\,\,s < 1,\,\,\,t < \min \left\{ {{U^+ }\left( {{w_1}} \right), - {U^ - }\left( {{w_2}} \right)} \right\}} \right\}, \hfill \\
  {P^ \pm }\left( w \right) = \left( {1 \pm r} \right)\left( {1 \pm s} \right){w^3} - \left( {3rs \pm r \pm s} \right) + 2rs, \hfill \\
  {w_1}: {\rm{the\,root\,in}}\,\left( {\sqrt {\max \left( {\frac{r}{{1 + r}},\frac{s}{{1 + s}}} \right)} ,1} \right)\,\,{\rm{of}}\,\,{P^ + }\left( w \right) = 0, \hfill \\
  {w_2}: {\rm{the\,greater\,root\,in}}\,\,\left( {1,\infty } \right)\,\,{\rm{of}}\,{P^- }\left( w \right) = 0, \hfill \\
  {U^ \pm }\left( w \right) = \frac{{1 - w}}{{{w^2}}}\left[ {w \pm r\left( {w - 1} \right)} \right]\left[ {w \pm s\left( {w - 1} \right)} \right]. \hfill \\
\end{gathered}
$$

System of partial differential equations:

$
\left\{ {\begin{array}{*{20}{l}}
  \begin{gathered}
  x\left( {1 - x} \right){u_{xx}} + \left( {1 - x} \right)y{u_{xy}} - yz{u_{yz}} + {z^2}{u_{zz}} + \hfill \\
 \,\,\,\,\,\,\,\,\, \left[ {c - \left( {{a_2} + b + 1} \right)x} \right]{u_x}
   - {a_2}y{u_y} + \left( {{a_2} - b + 1} \right)z{u_z} - {a_2}bu = 0, \hfill \\
\end{gathered}  \\
  \begin{gathered}
  y\left( {1 - y} \right){u_{yy}} + x\left( {1 - y} \right){u_{xy}} - xz{u_{xz}} + {z^2}{u_{zz}} + \hfill \\
 \,\,\,\,\,\,\,\,\, \left[ {c - \left( {{a_1} + b + 1} \right)y} \right]{u_y}
   - {a_1}x{u_x} + \left( {{a_1} - b + 1} \right)z{u_z} - {a_1}bu = 0, \hfill \\
\end{gathered}  \\
  \begin{gathered}
  z\left( {1 + z} \right){u_{zz}} + xy{u_{xy}} - x\left( {1 - z} \right){u_{xz}} - y\left( {1 - z} \right){u_{yz}} \hfill \\
 \,\,\,\,\,\,\,\,\,
   + \left[ {1 - b + \left( {{a_1} + {a_2} + 1} \right)z} \right]{u_z} + {a_1}x{u_x} + {a_2}y{u_y} + {a_1}{a_2}u = 0, \hfill \\
\end{gathered}
\end{array}} \right.
$

where $u\equiv \,\,   {F_{15b}}\left( {{a_1},{a_2},b;c;x,y,z} \right)$.

\bigskip

\begin{equation}
{F_{16a}}\left( {{a_1},{a_2},{a_3};c;x,y,z} \right)\hfill \\
   = \sum\limits_{m,n,p = 0}^\infty  {} \frac{{{{\left( {{a_1}} \right)}_{m + n}}{{\left( {{a_2}} \right)}_{n + p}}{{\left( {{a_3}} \right)}_{p + m}}}}{{{{\left( c \right)}_{m + n + p}}}}\frac{x^m}{m!}\frac{y^n}{n!}\frac{z^p}{p!},
\end{equation}

first appearance of this function in the literature, and old notation: [37], $H_{C}$,

region of convergence:
$$
\left\{ {r < 1,\,\,\,s < 1,\,\,\,t < 1,\,\,\,r + s + t - 2\sqrt {\left( {1 - r} \right)\left( {1 - s} \right)\left( {1 - t} \right)}  < 2} \right\}.
$$

System of partial differential equations:

$
\left\{ {\begin{array}{*{20}{l}}
  \begin{gathered}
  x\left( {1 - x} \right){u_{xx}} + \left( {1 - x} \right)y{u_{xy}} + \left( {1 - x} \right)z{u_{xz}} - yz{u_{yz}} \hfill \\
  \,\,\,\,\,\,\,\,\,     + \left[ {c - \left( {{a_1} + {a_3} + 1} \right)x} \right]{u_x}   - {a_3}y{u_y} - {a_1}z{u_z} - {a_1}{a_3}u = 0, \hfill \\
\end{gathered}  \\
  \begin{gathered}
  y\left( {1 - y} \right){u_{yy}} + x\left( {1 - y} \right){u_{xy}} - xz{u_{xz}} + \left( {1 - y} \right)z{u_{yz}} \hfill \\
 \,\,\,\,\,\,\,\,\,   + \left[ {c - \left( {{a_1} + {a_2} + 1} \right)y} \right]{u_y}   - {a_2}x{u_x} - {a_1}z{u_z} - {a_1}{a_2}u = 0 ,\hfill \\
\end{gathered}  \\
  \begin{gathered}
  z\left( {1 - z} \right){u_{zz}} - xy{u_{xy}} + x\left( {1 - z} \right){u_{xz}} + y\left( {1 - z} \right){u_{yz}} \hfill \\
\,\,\,\,\,\,\,\,\,    + \left[ {c - \left( {{a_2} + {a_3} + 1} \right)z} \right]{u_z}   - {a_2}x{u_x} - {a_3}y{u_y} - {a_2}{a_3}u = 0, \hfill \\
\end{gathered}
\end{array}} \right.
$

where $u\equiv \,\,   {F_{16a}}\left( {{a_1},{a_2},{a_3};c;x,y,z} \right)$.

\bigskip

\begin{equation}
 {F_{17a}}\left( {{a_1},{a_2},{a_3},{a_4};{c_1},{c_2},{c_3};x,y,z} \right) \hfill \\
  = \sum\limits_{m,n,p = 0}^\infty  {} \frac{{{{\left( {{a_1}} \right)}_{m + n + p}}{{\left( {{a_2}} \right)}_m}{{\left( {{a_3}} \right)}_n}{{\left( {{a_4}} \right)}_p}}}{{{{\left( {{c_1}} \right)}_m}{{\left( {{c_2}} \right)}_n}{{\left( {{c_3}} \right)}_p}}}\frac{x^m}{m!}\frac{y^n}{n!}\frac{z^p}{p!},
\end{equation}

first appearance of this function in the literature, and old notation: [24], $F_{A}^{(3)}$, $F_1$,

region of convergence:
$$
\left\{ {r + s + t < 1} \right\}.
$$

System of partial differential equations:

$
\left\{ {\begin{array}{*{20}{l}}
  x\left( {1 - x} \right){u_{xx}} - xy{u_{xy}} - xz{u_{xz}} + \left[ {{c_1} - \left( {{a_1} + {a_2} + 1} \right)x} \right]{u_x}
   - {a_2}y{u_y} - {a_2}z{u_z} - {a_1}{a_2}u = 0, \\
  y\left( {1 - y} \right){u_{yy}} - xy{u_{xy}} - yz{u_{yz}} + \left[ {{c_2} - \left( {{a_1} + {a_3} + 1} \right)y} \right]{u_y}
   - {a_3}x{u_x} - {a_3}z{u_z} - {a_1}{a_3}u = 0. \\
  z\left( {1 - z} \right){u_{zz}} - xz{u_{xz}} - yz{u_{yz}} + \left[ {{c_3} - \left( {{a_1} + {a_4} + 1} \right)z} \right]{u_z}
   - {a_4}x{u_x} - {a_4}y{u_y} - {a_1}{a_4}u = 0,
\end{array}} \right.
$

where $u\equiv \,\,   {F_{17a}}\left( {{a_1},{a_2},{a_3},{a_4};{c_1},{c_2},{c_3};x,y,z} \right)$.

Particular solutions:

$
{u_1} = {F_{17a}}\left( {{a_1},{a_2},{a_3},{a_4};{c_1},{c_2},{c_3};x,y,z} \right),
$

$
{u_2} = {x^{1 - {c_1}}}{F_{17a}}\left( {1 - {c_1} + {a_1},1 - {c_1} + {a_2},{a_3},{a_4};2 - {c_1},{c_2},{c_3};x,y,z} \right),
$

$
{u_3} = {y^{1 - {c_2}}}{F_{17a}}\left( {1 - {c_2} + {a_1},{a_2},1 - {c_2} + {a_3},{a_4};{c_1},2 - {c_2},{c_3};x,y,z} \right),
$

$
{u_4} = {z^{1 - {c_3}}}{F_{17a}}\left( {{a_1},{a_2},{a_3},{a_4};{c_1},{c_2},{c_3};x,y,z} \right),
$

$
{u_5} = {x^{1 - {c_1}}}{y^{1 - {c_2}}}{F_{17a}}\left( {2 - {c_1} - {c_2} + {a_1},1 - {c_1} + {a_2},1 - {c_2} + {a_3},{a_4};2 - {c_1},2 - {c_2},{c_3};x,y,z} \right),
$

$
{u_6} = {y^{1 - {c_2}}}{z^{1 - {c_3}}}{F_{17a}}\left( {2 - {c_2} - {c_3} + {a_1},{a_2},1 - {c_2} + {a_3},1 - {c_3} + {a_4};{c_1},2 - {c_2},2 - {c_3};x,y,z} \right),
$

$
{u_7} = {x^{1 - {c_1}}}{z^{1 - {c_3}}}{F_{17a}}\left( {2 - {c_1} - {c_3} + {a_1},1 - {c_1} + {a_2},{a_3},1 - {c_3} + {a_4};2 - {c_1},{c_2},2 - {c_3};x,y,z} \right),
$

$
  {u_8} = {x^{1 - {c_1}}}{y^{1 - {c_2}}}{z^{1 - {c_3}}}\times $

  $\,\,\,\,\,\,\,\,\,\times  {F_{17a}}\left( {3 - {c_1} - {c_2} - {c_3} + {a_1},1 - {c_1} + {a_2},1 - {c_2} + {a_3},1 - {c_3} + {a_4};}
   {2 - {c_1},2 - {c_2},2 - {c_3};x,y,z} \right).
   $

\bigskip

\begin{equation}
{F_{17b}}\left( {{a_1},{a_2},{a_3},b;{c_1},{c_2};x,y,z} \right) \hfill \\
  = \sum\limits_{m,n,p = 0}^\infty  {} \frac{{{{\left( {{a_1}} \right)}_{m + n + p}}{{\left( {{a_2}} \right)}_n}{{\left( {{a_3}} \right)}_p}{{\left( b \right)}_{m - n}}}}{{{{\left( {{c_1}} \right)}_m}{{\left( {{c_2}} \right)}_p}}}\frac{x^m}{m!}\frac{y^n}{n!}\frac{z^p}{p!},
\end{equation}

region of convergence:
$$
\left\{ {r + t < 1,\,\,\,\sqrt s  < \sqrt {1 + r - t}  - \sqrt r } \right\}.
$$

System of partial differential equations:

$
\left\{ {\begin{array}{*{20}{l}}
  \begin{gathered}
  x\left( {1 - x} \right){u_{xx}} - xz{u_{xz}} + {y^2}{u_{yy}} + yz{u_{yz}}\hfill\\\,\,\,\,\,\,\,\,\,\,
   + \left[ {{c_1} - \left( {{a_1} + b + 1} \right)x} \right]{u_x} + \left( {{a_1} - b + 1} \right)y{u_y}   - bz{u_z} - {a_1}bu = 0 ,\hfill \\
\end{gathered}  \\
  \begin{gathered}
  y\left( {1 + y} \right){u_{yy}} - x\left( {1 - y} \right){u_{xy}} + yz{u_{yz}}\hfill\\\,\,\,\,\,\,\,\,\,\,
   + \left[ {1 - b + \left( {{a_1} + {a_2} + 1} \right)y} \right]{u_y} + {a_2}x{u_x}   + {a_2}z{u_z} + {a_1}{a_2}u = 0 , \hfill \\
\end{gathered}  \\
  \begin{gathered}
  z\left( {1 - z} \right){u_{zz}} - xz{u_{xz}} - yz{u_{yz}} + \left[ {{c_2} - \left( {{a_1} + {a_3} + 1} \right)z} \right]{u_z}
   - {a_3}x{u_x} - {a_3}y{u_y} - {a_1}{a_3}u = 0, \hfill \\
\end{gathered}
\end{array}} \right.
$

where $u\equiv \,\,   {F_{17b}}\left( {{a_1},{a_2},{a_3},b;{c_1},{c_2};x,y,z} \right)$.

Particular solutions:

$
{u_1} = {F_{17b}}\left( {{a_1},{a_2},{a_3},b;{c_1},{c_2};x,y,z} \right),
$

$
{u_2} = {x^{1 - {c_1}}}{F_{17b}}\left( {1 - {c_1} + {a_1},{a_2},{a_3},1 - {c_1} + b;2 - {c_1},{c_2};x,y,z} \right),
$

$
{u_3} = {z^{1 - {c_2}}}{F_{17b}}\left( {1 - {c_2} + {a_1},{a_2},1 - {c_2} + {a_3},b;{c_1},2 - {c_2};x,y,z} \right),
$

$
{u_4} = {x^{1 - {c_1}}}{z^{1 - {c_2}}}{F_{17b}}\left( {2 - {c_1} - {c_2} + {a_1},{a_2},1 - {c_2} + {a_3},1 - {c_1} + b;2 - {c_1},2 - {c_2};x,y,z} \right).
$

\bigskip

\begin{equation}
{F_{17c}}\left( {{a_1},{a_2},{b_1},{b_2};c;x,y,z} \right)  = \sum\limits_{m,n,p = 0}^\infty  {} \frac{{{{\left( {{a_1}} \right)}_{m + n + p}}{{\left( {{a_2}} \right)}_p}{{\left( {{b_1}} \right)}_{m - n}}{{\left( {{b_2}} \right)}_{n - m}}}}{{{{\left( c \right)}_p}}}\frac{x^m}{m!}\frac{y^n}{n!}\frac{z^p}{p!},
\end{equation}

region of convergence:
$$
\left\{ {r + s + t < 1} \right\}.
$$

System of partial differential equations:

$
\left\{ {\begin{array}{*{20}{l}}
  x\left( {1 + x} \right){u_{xx}} - y{u_{xy}} + xz{u_{xz}} - {y^2}{u_{yy}} - yz{u_{yz}} \\
 \,\,\,\,\,\,\,\,\,  + \left[ {1 - {b_2} + \left( {{a_1} + {b_1} + 1} \right)x} \right]{u_x}
  - \left( {{a_1} - {b_1} + 1} \right)y{u_y} + {b_1}z{u_z} + {a_1}{b_1}u = 0,   \\
  y\left( {1 + y} \right){u_{yy}} - x{u_{xy}} - {x^2}{u_{xx}} - xz{u_{xz}} + yz{u_{yz}} \\
\,\,\,\,\,\,\,\,\,   + \left[ {1 - {b_1} + \left( {{a_1} + {b_2} + 1} \right)y} \right]{u_y}
  - \left( {{a_1} - {b_2} + 1} \right)x{u_x} + {b_2}z{u_z} + {a_1}{b_2}u = 0, \\
  z\left( {1 - z} \right){u_{zz}} - xz{u_{xz}} - yz{u_{yz}} + \left[ {c - \left( {{a_1} + {a_2} + 1} \right)z} \right]{u_z}
  - {a_2}x{u_x}    - {a_2}y{u_y} - {a_1}{a_2}u = 0,
\end{array}} \right.
$

where $u\equiv \,\,   {F_{17c}}\left( {{a_1},{a_2},{b_1},{b_2};c;x,y,z} \right)$.

Particular solutions:

$
{u_1} = {F_{17c}}\left( {{a_1},{a_2},{b_1},{b_1};c;x,y,z} \right),
$

$
{u_2} = {z^{1 - c}}{F_{17c}}\left( {1 - c + {a_1},1 - c + {a_2},{b_1},{b_1};2 - c;x,y,z} \right).
$

\bigskip

\begin{equation}
 {F_{17d}}\left( {{a_1},{a_2},{b_1},{b_2};c;x,y,z} \right)\hfill \\
   = \sum\limits_{m,n,p = 0}^\infty  {} \frac{{{{\left( {{a_1}} \right)}_{m + n + p}}{{\left( {{a_2}} \right)}_p}{{\left( {{b_1}} \right)}_{m - n}}{{\left( {{b_2}} \right)}_{n - p}}}}{{{{\left( c \right)}_m}}}{x^m}\frac{x^m}{m!}\frac{y^n}{n!}\frac{z^p}{p!},
\end{equation}

region of convergence:
$$
\begin{gathered}
    \left\{ {s < 1,\,\,\,t + 2\sqrt {st}  < 1,}  {r < \min \left\{ {\frac{{{{\left( {1 - s} \right)}^2}}}{{4s}},\left( {1 - t} \right){\Psi _1}\left( {\frac{{st}}{{{{\left( {1 - t} \right)}^2}}}} \right),\left( {1 - t} \right){\Psi _2}\left( {\frac{{st}}{{{{\left( {1 - t} \right)}^2}}}} \right)} \right\}} \right\}. \hfill \\
\end{gathered}
$$

System of partial differential equations:

$
\left\{ {\begin{array}{*{20}{l}}
  \begin{gathered}
  x\left( {1 - x} \right){u_{xx}} - xz{u_{xz}} + {y^2}{u_{yy}} + yz{u_{yz}} \hfill \\
\,\,\,\,\,\,\,\,\,   + \left[ {c - \left( {{a_1} + {b_1} + 1} \right)x} \right]{u_x}
   + \left( {{a_1} - {b_1} + 1} \right)y{u_y} - {b_1}z{u_z} - {a_1}{b_1}u = 0, \hfill \\
\end{gathered}  \\
  \begin{gathered}
  y\left( {1 + y} \right){u_{yy}} - x\left( {1 - y} \right){u_{xy}} - xz{u_{xz}} - {z^2}{u_{zz}} \hfill \\
\,\,\,\,\,\,\,\,\,   + \left[ {1 - {b_1} + \left( {{a_1} + {b_2} + 1} \right)y} \right]{u_y}
   + {b_2}x{u_x} - \left( {{a_1} - {b_2} + 1} \right)z{u_z} + {a_1}{b_2}u = 0, \hfill \\
\end{gathered}  \\
  \begin{gathered}
  z\left( {1 + z} \right){u_{zz}} + xz{u_{xz}} - y\left( {1 - z} \right){u_{yz}}\hfill \\ \,\,\,\,\,\,\,\,\,\,
  + \left[ {1 - {b_2} + \left( {{a_1} + {a_2} + 1} \right)z} \right]{u_z}
   + {a_2}x{u_x} + {a_2}y{u_y} + {a_1}{a_2}u = 0, \hfill \\
\end{gathered}
\end{array}} \right.
$

where $u\equiv \,\,   {F_{17d}}\left( {{a_1},{a_2},{b_1},{b_2};c;x,y,z} \right)$.

Particular solutions:

$
{u_1} = {F_{17d}}\left( {{a_1},{a_2},{b_1},{b_2};c;x,y,z} \right),
$

$
{u_2} = {x^{1 - c}}{F_{17d}}\left( {1 - c + {a_1},{a_2},1 - c + {b_1},{b_2};2 - c;x,y,z} \right).
$

\bigskip

\begin{equation}
 {F_{17e}}\left( {{a},{b_1},{b_2},{b_3};x,y,z} \right) \hfill \\
  = \sum\limits_{m,n,p = 0}^\infty  {} {{{{\left( {{a}} \right)}_{m + n + p}}{{\left( {{b_1}} \right)}_{m - n}}{{\left( {{b_2}} \right)}_{n - p}}{{\left( {{b_3}} \right)}_{p - m}}}}\frac{x^m}{m!}\frac{y^n}{n!}\frac{z^p}{p!},
\end{equation}

region of convergence:
$$
\begin{gathered}
  {{E_{rst}} \cap {E_{str}} \cap {E_{trs}}}, \hfill \\
  {E_{rst}} = \left\{ {r < 1,\,\,\,s + 2\sqrt {rs}  < 1,\,\,\,t < U\left( {r,s} \right)} \right\}, \hfill \\
  {P_{rs}}\left( w \right) = 2{w^3} - \left( {1 - r - s} \right){w^2} - {r^2}s, \hfill \\
  {w_{rs}}:\,\,{\rm{the\,root\,in}}\,\left( {\sqrt {rs} ,\infty } \right)\,{\rm{of}}\,{P_{rs}}\left( w \right) = 0, \hfill \\
  U\left( {r,s} \right) = \frac{{{{\left( {r + {w_{rs}}} \right)}^2}\left( {w_{rs}^2 - rs} \right)}}{{rw_{rs}^2}}. \hfill \\
\end{gathered}
$$

System of partial differential equations:

$
\left\{ {\begin{array}{*{20}{l}}
  \begin{gathered}
  x\left( {1 + x} \right){u_{xx}} - \left( {1 - x} \right)z{u_{xz}} - {y^2}{u_{yy}} - yz{u_{yz}} \hfill\\
\,\,\,\,\,\,\,\,\,   + \left[ {1 - {b_3} + \left( {{a} + {b_1} + 1} \right)x} \right]{u_x}
   - \left( {{a} - {b_1} + 1} \right)y{u_y} + {b_1}z{u_z} + {a}{b_1}u = 0, \hfill \\
\end{gathered}  \\
  \begin{gathered}
  y\left( {1 + y} \right){u_{yy}} - x\left( {1 - y} \right){u_{xy}} - xz{u_{xz}} - {z^2}{u_{zz}}\hfill\\
 \,\,\,\,\,\,\,\,\,   + \left[ {1 - {b_1} + \left( {{a} + {b_2} + 1} \right)y} \right]{u_y}
   + {b_2}x{u_x} - \left( {{a} - {b_2} + 1} \right)z{u_z} + {a}{b_2}u = 0, \hfill \\
\end{gathered}  \\
  \begin{gathered}
  z\left( {1 + z} \right){u_{zz}} - xy{u_{xy}} - {x^2}{u_{xx}} - y\left( {1 - z} \right){u_{yz}}\hfill\\
 \,\,\,\,\,\,\,\,\,   + \left[ {1 - {b_2} + \left( {{a} + {b_3} + 1} \right)z} \right]{u_z}
   - \left( {{a} - {b_3} + 1} \right)x{u_x} + {b_3}y{u_y} + {a}{b_3}u = 0, \hfill \\
\end{gathered}
\end{array}} \right.
$

where $u\equiv \,\,    {F_{17e}}\left( {{a},{b_1},{b_2},{b_3};x,y,z} \right)$.

\bigskip

\begin{equation}
{F_{18a}}\left( {{a_1},{a_2},{a_3},{a_4};{c_1},{c_2};x,y,z} \right)
  = \sum\limits_{m,n,p = 0}^\infty  {} \frac{{{{\left( {{a_1}} \right)}_{m + n + p}}{{\left( {{a_2}} \right)}_m}{{\left( {{a_3}} \right)}_n}{{\left( {{a_4}} \right)}_p}}}{{{{\left( {{c_1}} \right)}_{m + n}}{{\left( {{c_2}} \right)}_p}}}\frac{x^m}{m!}\frac{y^n}{n!}\frac{z^p}{p!},
\end{equation}

first appearance of this function in the literature, and old notation: [24], $F_{8}$,\,\,[33], $F_G$,

region of convergence:
$$
\left\{ {r + t < 1,\,\,\,s + t < 1} \right\}.
$$

System of partial differential equations:

$
\left\{ {\begin{array}{*{20}{l}}
  \begin{gathered}
  x\left( {1 - x} \right){u_{xx}} + \left( {1 - x} \right)y{u_{xy}} - xz{u_{xz}}
   + \left[ {{c_1} - \left( {{a_1} + {a_2} + 1} \right)x} \right]{u_x} - {a_2}y{u_y} - {a_2}z{u_z} - {a_1}{a_2}u = 0 ,\hfill \\
\end{gathered}  \\
  \begin{gathered}
  y\left( {1 - y} \right){u_{yy}} + x\left( {1 - y} \right){u_{xy}} - yz{u_{yz}}
   + \left[ {{c_1} - \left( {{a_1} + {a_3} + 1} \right)y} \right]{u_y} - {a_3}x{u_x} - {a_3}z{u_z} - {a_1}{a_3}u = 0, \hfill \\
\end{gathered}  \\
  \begin{gathered}
  z\left( {1 - z} \right){u_{zz}} - xz{u_{xz}} - yz{u_{yz}}
   + \left[ {{c_2} - \left( {{a_1} + {a_4} + 1} \right)z} \right]{u_z} - {a_4}x{u_x} - {a_4}y{u_y} - {a_1}{a_4}u = 0, \hfill \\
\end{gathered}
\end{array}} \right.
$

where $u\equiv \,\,   {F_{18a}}\left( {{a_1},{a_2},{a_3},{a_4};{c_1},{c_2};x,y,z} \right)$.

Particular solutions:

$
{u_1} = {F_{18a}}\left( {{a_1},{a_2},{a_3},{a_4};{c_1},{c_2};x,y,z} \right),
$

$
{u_2} = {z^{1 - {c_2}}}{F_{18a}}\left( {1 - {c_2} + {a_1},{a_2},{a_3},1 - {c_2} + {a_4};{c_1},2 - {c_2};x,y,z} \right).
$

\bigskip

\begin{equation}
{F_{18b}}\left( {{a_1},{a_2},{a_3},b;c;x,y,z} \right) \hfill \\
  = \sum\limits_{m,n,p = 0}^\infty  {} \frac{{{{\left( {{a_1}} \right)}_{m + n + p}}{{\left( {{a_2}} \right)}_n}{{\left( {{a_3}} \right)}_p}{{\left( b \right)}_{m - p}}}}{{{{\left( c \right)}_{m + n}}}}\frac{x^m}{m!}\frac{y^n}{n!}\frac{z^p}{p!},
\end{equation}

region of convergence:
$$
 \left\{ {r < 1,\,\,\,t + 2\sqrt {rt}  < 1,\,\,\,s < \frac{1}{2}\left( {1 - t} \right) + \frac{1}{2}\sqrt {{{\left( {1 - t} \right)}^2} - 4rt} } \right\}.
$$

System of partial differential equations:

$
\left\{ {\begin{array}{*{20}{l}}
  \begin{gathered}
  x\left( {1 - x} \right){u_{xx}} + \left( {1 - x} \right)y{u_{xy}} + yz{u_{yz}} + {z^2}{u_{zz}} \hfill \\ \,\,\,\,\,\,\,\,\,
  + \left[ {c - \left( {{a_1} + b + 1} \right)x} \right]{u_x} - by{u_y}  + \left( {{a_1} - b + 1} \right)z{u_z} - {a_1}bu = 0, \hfill \\
\end{gathered}  \\
  \begin{gathered}
  y\left( {1 - y} \right){u_{yy}} + x\left( {1 - y} \right){u_{xy}} - yz{u_{yz}}  \hfill \\ \,\,\,\,\,\,\,\,\,\, + \left[ {c - \left( {{a_1} + {a_2} + 1} \right)y} \right]{u_y}
  - {a_2}x{u_x}     - {a_2}z{u_z} - {a_1}{a_2}u = 0, \hfill \\
\end{gathered}  \\
  \begin{gathered}
  z\left( {1 + z} \right){u_{zz}} - x\left( {1 - z} \right){u_{xz}} + yz{u_{yz}} \hfill \\ \,\,\,\,\,\,\,\,\,\,
   + \left[ {1 - b + \left( {{a_1} + {a_3} + 1} \right)z} \right]{u_z} + {a_3}x{u_x} +
   {a_3}y{u_y} + {a_1}{a_3}u = 0,\hfill \\
\end{gathered}
\end{array}} \right.
$

where $u\equiv \,\,   {F_{18b}}\left( {{a_1},{a_2},{a_3},b;c;x,y,z} \right)$.

\bigskip

\begin{equation}
{F_{18c}}\left( {{a_1},{a_2},{a_3},b;c;x,y,z} \right)\hfill \\
   = \sum\limits_{m,n,p = 0}^\infty  {} \frac{{{{\left( {{a_1}} \right)}_{m + n + p}}{{\left( {{a_2}} \right)}_m}{{\left( {{a_3}} \right)}_n}{{\left( b \right)}_{p - m - n}}}}{{{{\left( c \right)}_p}}}\frac{x^m}{m!}\frac{y^n}{n!}\frac{z^p}{p!},
\end{equation}

region of convergence:
$$
\left\{ {r < 1,\,\,\,s < 1,\,\,\,t < \min \left\{ {1,\frac{{{{\left( {1 - \max \left( {r,s} \right)} \right)}^2}}}{{4\max \left( {r,s} \right)}}} \right\}} \right\}.
$$

System of partial differential equations:

$
\left\{ {\begin{array}{*{20}{l}}
  \begin{gathered}
  x\left( {1 + x} \right){u_{xx}} + \left( {1 + x} \right)y{u_{xy}} - \left( {1 - x} \right)z{u_{xz}} \hfill \\ \,\,\,\,\,\,\,\,\,\,
   + \left[ {1 - b + \left( {{a_1} + {a_2} + 1} \right)x} \right]{u_x} + {a_2}y{u_y}     + {a_2}z{u_z} + {a_1}{a_2}u = 0, \hfill \\
\end{gathered}  \\
  \begin{gathered}
  y\left( {1 + y} \right){u_{yy}} + x\left( {1 + y} \right){u_{xy}} - \left( {1 - y} \right)z{u_{yz}}  \hfill \\ \,\,\,\,\,\,\,\,\,\,
   + \left[ {1 - b + \left( {{a_1} + {a_3} + 1} \right)y} \right]{u_y} + {a_3}x{u_x}   + {a_3}z{u_z} + {a_1}{a_3}u = 0, \hfill \\
\end{gathered}  \\
  \begin{gathered}
  z\left( {1 - z} \right){u_{zz}} + 2xy{u_{xy}} + {x^2}{u_{xx}} + {y^2}{u_{yy}} + \left[ {c - \left( {{a_1} + b + 1} \right)z} \right]{u_z}\hfill \\
 \,\,\,\,\,\,\,\,\,   + \left( {{a_1} - b + 1} \right)x{u_x}
   + \left( {{a_1} - b + 1} \right)y{u_y} - {a_1}bu = 0, \hfill \\
\end{gathered}
\end{array}} \right.
$

where $u\equiv \,\,   {F_{18c}}\left( {{a_1},{a_2},{a_3},b;c;x,y,z} \right)$.

Particular solutions:

$
{u_1} = {F_{18c}}\left( {{a_1},{a_2},{a_3},b;c;x,y,z} \right),
$

$
{u_2} = {z^{1 - c}}{F_{18c}}\left( {1 - c + {a_1},{a_2},{a_3},1 - c + b;2 - c;x,y,z} \right).
$

\bigskip

\begin{equation}
{F_{18d}}\left( {{a_1},{a_2},{b_1},{b_2};x,y,z} \right) \hfill \\
  = \sum\limits_{m,n,p = 0}^\infty  {} {{{{\left( {{a_1}} \right)}_{m + n + p}}{{\left( {{a_2}} \right)}_n}{{\left( {{b_1}} \right)}_{m - p}}{{\left( {{b_2}} \right)}_{p - m - n}}}}\frac{x^m}{m!}\frac{y^n}{n!}\frac{z^p}{p!},
\end{equation}

region of convergence:
$$
 \left\{ {t < 1,\,\,\,s + 2\sqrt {st}  < 1,\,\,\,r < \min \left\{ {1 - s,\frac{{{{\left( {1 - s} \right)}^2} - 4st}}{{4t}}} \right\}} \right\}.
$$

System of partial differential equations:

$
\left\{ {\begin{array}{*{20}{l}}
  \begin{gathered}
  x\left( {1 + x} \right){u_{xx}} + \left( {1 + x} \right)y{u_{xy}} - z{u_{xz}} - yz{u_{yz}} - {z^2}{u_{zz}} \hfill \hfill \\
  \,\,\,\,\,\,\,\,\,  + \left[ {1 - {b_2} + \left( {{a_1} + {b_1} + 1} \right)x} \right]{u_x} + {b_1}y{u_y}
    - \left( {{a_1} - {b_1} + 1} \right)z{u_z} + {a_1}{b_1}u = 0, \hfill \\
\end{gathered}  \\
  \begin{gathered}
  y\left( {1 + y} \right){u_{yy}} + x\left( {1 + y} \right){u_{xy}} - \left( {1 - y} \right)z{u_{yz}} \hfill\\
  \,\,\,\,\,\,\,\,\,
   + \left[ {1 - {b_2} + \left( {{a_1} + {a_2} + 1} \right)y} \right]{u_y} + {a_2}x{u_x} + {a_2}z{u_z} + {a_1}{a_2}u = 0 ,\hfill \\
\end{gathered}  \\
  \begin{gathered}
  z\left( {1 + z} \right){u_{zz}} - 2xy{u_{xy}} - x{u_{xz}} - {x^2}{u_{xx}} - {y^2}{u_{yy}} \hfill \\
 \,\,\,\,\,\,\,\,\,   + \left[ {1 - {b_1} + \left( {{a_1} + {b_2} + 1} \right)z} \right]{u_z} - \left( {{a_1} - {b_2} + 1} \right)x{u_x}
   - \left( {{a_1} - {b_2} + 1} \right)y{u_y} + {a_1}{b_2}u = 0, \hfill \\
\end{gathered}
\end{array}} \right.
$

where $u\equiv \,\,   {F_{18d}}\left( {{a_1},{a_2},{b_1},{b_2};x,y,z} \right)$.

\bigskip

\begin{equation}
{F_{19a}}\left( {{a_1},{a_2},{a_3},{a_4};c;x,y,z} \right) \hfill \\
  = \sum\limits_{m,n,p = 0}^\infty  {} \frac{{{{\left( {{a_1}} \right)}_{m + n + p}}{{\left( {{a_2}} \right)}_m}{{\left( {{a_3}} \right)}_n}{{\left( {{a_4}} \right)}_p}}}{{{{\left( c \right)}_{m + n + p}}}}\frac{x^m}{m!}\frac{y^n}{n!}\frac{z^p}{p!},
\end{equation}

first appearance of this function in the literature, and old notation: [35]; [24], $F_D^{(3)}$,  $F_{9}$,

region of convergence:
$$
\left\{ {r < 1,\,\,\,s < 1,\,\,\,t < 1} \right\}.
$$

System of partial differential equations:

$
\left\{ {\begin{array}{*{20}{l}}
  \begin{gathered}
  x\left( {1 - x} \right){u_{xx}} + \left( {1 - x} \right)y{u_{xy}} + \left( {1 - x} \right)z{u_{xz}}\hfill\\\,\,\,\,\,\,\,\,\,\,
   + \left[ {c - \left( {{a_1} + {a_2} + 1} \right)x} \right]{u_x} - {a_2}y{u_y} - {a_2}z{u_z} - {a_1}{a_2}u = 0 ,\hfill \\
\end{gathered}  \\
  \begin{gathered}
  y\left( {1 - y} \right){u_{yy}} + x\left( {1 - y} \right){u_{xy}} + \left( {1 - y} \right)z{u_{yz}}\hfill\\ \,\,\,\,\,\,\,\,\,\,
   + \left[ {c - \left( {{a_1} + {a_3} + 1} \right)y} \right]{u_y} - {a_3}x{u_x} - {a_3}z{u_z} - {a_1}{a_3}u = 0, \hfill \\
\end{gathered}  \\
  \begin{gathered}
  z\left( {1 - z} \right){u_{zz}} + x\left( {1 - z} \right){u_{xz}} + y\left( {1 - z} \right){u_{yz}}\hfill\\\,\,\,\,\,\,\,\,\,\,
   + \left[ {c - \left( {{a_1} + {a_4} + 1} \right)z} \right]{u_z} - {a_4}x{u_x} - {a_4}y{u_y} - {a_1}{a_4}u = 0, \hfill \\
\end{gathered}
\end{array}} \right.
$

where $u\equiv \,\,   {F_{19a}}\left( {{a_1},{a_2},{a_3},{a_4};c;x,y,z} \right)$.

\bigskip

\begin{equation}
{F_{20a}}\left( {{a_1},{a_2},{a_3};{c_1},{c_2},{c_3};x,y,z} \right) \hfill \\
  = \sum\limits_{m,n,p = 0}^\infty  {} \frac{{{{\left( {{a_1}} \right)}_{m + n + p}}{{\left( {{a_2}} \right)}_{m + n}}{{\left( {{a_3}} \right)}_p}}}{{{{\left( {{c_1}} \right)}_m}{{\left( {{c_2}} \right)}_n}{{\left( {{c_3}} \right)}_p}}}\frac{x^m}{m!}\frac{y^n}{n!}\frac{z^p}{p!},
\end{equation}

first appearance of this function in the literature, and old notation: [24], $F_{4}$;\,\,[33], $F_E$,

region of convergence:
$$
\left\{ {{{\left( {\sqrt r  + \sqrt s } \right)}^2} + t < 1} \right\}.
$$

System of partial differential equations:

$
\left\{ {\begin{array}{*{20}{l}}
  \begin{gathered}
  x\left( {1 - x} \right){u_{xx}} - 2xy{u_{xy}} - xz{u_{xz}} - {y^2}{u_{yy}} - yz{u_{yz}} \hfill \\
 \,\,\,\,\,\,\,\,\,   + \left[ {{c_1} - \left( {{a_1} + {a_2} + 1} \right)x} \right]{u_x} - \left( {{a_1} + {a_2} + 1} \right)y{u_y}
   - {a_2}z{u_z} - {a_1}{a_2}u = 0, \hfill \\
\end{gathered}  \\
  \begin{gathered}
  y\left( {1 - y} \right){u_{yy}} - {x^2}{u_{xx}} - 2xy{u_{xy}} - xz{u_{xz}} - yz{u_{yz}} \hfill \\
 \,\,\,\,\,\,\,\,\,   + \left[ {{c_2} - \left( {{a_1} + {a_2} + 1} \right)y} \right]{u_y} - \left( {{a_1} + {a_2} + 1} \right)x{u_x}
   - {a_2}z{u_z} - {a_1}{a_2}u = 0, \hfill \\
\end{gathered}  \\
  \begin{gathered}
  z\left( {1 - z} \right){u_{zz}} - xz{u_{xz}} - yz{u_{yz}}
   + \left[ {{c_3} - \left( {{a_1} + {a_3} + 1} \right)z} \right]{u_z} - {a_3}x{u_x} - {a_3}y{u_y} - {a_1}{a_3}u = 0, \hfill \\
\end{gathered}
\end{array}} \right.
$

where $u\equiv \,\,   {F_{20a}}\left( {{a_1},{a_2},{a_3};{c_1},{c_2},{c_3};x,y,z} \right)$.

Particular solutions:

$
{u_1} = {F_{20a}}\left( {{a_1},{a_2},{a_3};{c_1},{c_2},{c_3};x,y,z} \right),
$

$
{u_2} = {x^{1 - {c_1}}}{F_{20a}}\left( {1 - {c_1} + {a_1},1 - {c_1} + {a_2},{a_3};{c_1},{c_2},{c_3};x,y,z} \right),
$

$
{u_3} = {y^{1 - {c_2}}}{F_{20a}}\left( {1 - {c_2} + {a_1},1 - {c_2} + {a_2},{a_3};{c_1},2 - {c_2},{c_3};x,y,z} \right),
$

$
{u_4} = {z^{1 - {c_3}}}{F_{20a}}\left( {1 - {c_3} + {a_1},{a_2},1 - {c_3} + {a_3};{c_1},{c_2},2 - {c_3};x,y,z} \right),
$

$
{u_5} = {x^{1 - {c_1}}}{y^{1 - {c_2}}}{F_{20a}}\left( {2 - {c_1} - {c_2} + {a_1},2 - {c_1} - {c_2} + {a_2},{a_3};2 - {c_1},2 - {c_2},{c_3};x,y,z} \right),
$

$
{u_6} = {y^{1 - {c_2}}}{z^{1 - {c_3}}}{F_{20a}}\left( {2 - {c_2} - {c_3} + {a_1},1 - {c_2} + {a_2},1 - {c_3} + {a_3};{c_1},2 - {c_2},2 - {c_3};x,y,z} \right),
$

$
{u_7} = {x^{1 - {c_1}}}{z^{1 - {c_3}}}{F_{20a}}\left( {2 - {c_1} - {c_3} + {a_1},1 - {c_1} + {a_2},1 - {c_3} + {a_3};2 - {c_1},{c_2},2 - {c_3};x,y,z} \right),
$

$
    {u_8} = {x^{1 - {c_1}}}{y^{1 - {c_2}}}{z^{1 - {c_3}}}\times $
     
     $ \,\,\,\,\,\,\,\,\,\times{F_{20a}}\left( {3 - {c_1} - {c_2} - {c_3} + {a_1},2 - {c_2} - {c_1} + {a_2},1 - {c_3} + {a_3};}  {2 - {c_1},2 - {c_2},2 - {c_3};x,y,z} \right).
    $

\bigskip

\begin{equation}
{F_{20b}}\left( {{a_1},{a_2},b;{c_1},{c_2};x,y,z} \right) = \sum\limits_{m,n,p = 0}^\infty  {} \frac{{{{\left( {{a_1}} \right)}_{m + n + p}}{{\left( {{a_2}} \right)}_p}{{\left( b \right)}_{m + n - p}}}}{{{{\left( {{c_1}} \right)}_m}{{\left( {{c_2}} \right)}_n}}}\frac{x^m}{m!}\frac{y^n}{n!}\frac{z^p}{p!},
\end{equation}

region of convergence:
$$
 \left\{ r < 1,\,\,\sqrt{r}+\sqrt{s}< \min\left\{1, \frac{1-t}{2\sqrt{t}}\right\} \right\}.
$$

System of partial differential equations:

$
\left\{ {\begin{array}{*{20}{l}}
  \begin{gathered}
  x\left( {1 - x} \right){u_{xx}} - 2xy{u_{xy}} - {y^2}{u_{yy}} + {z^2}{u_{zz}} \hfill \\
 \,\,\,\,\,\,\,\,\,  + \left[ {{c_1} - \left( {{a_1} + b + 1} \right)x} \right]{u_x}
   - \left( {{a_1} + b + 1} \right)y{u_y} + \left( {{a_1} - b + 1} \right)z{u_z} - {a_1}bu = 0, \hfill \\
\end{gathered}  \\
  \begin{gathered}
  y\left( {1 - y} \right){u_{yy}} - 2xy{u_{xy}} - {x^2}{u_{xx}} + {z^2}{u_{zz}}\hfill \\
 \,\,\,\,\,\,\,\,\,   + \left[ {{c_2} - \left( {{a_1} + b + 1} \right)y} \right]{u_y}
   - \left( {{a_1} + b + 1} \right)x{u_x} + \left( {{a_1} - b + 1} \right)z{u_z} - {a_1}bu = 0, \hfill \\
\end{gathered}  \\
  \begin{gathered}
  z\left( {1 + z} \right){u_{zz}} - x\left( {1 - z} \right){u_{xz}} - y\left( {1 - z} \right){u_{yz}} \hfill\\ \,\,\,\,\,\,\,\,\,  + \left[ {1 - b + \left( {{a_1} + {a_2} + 1} \right)z} \right]{u_z}     + {a_2}x{u_x} + {a_2}y{u_y} + {a_1}{a_2}u = 0, \hfill \\
\end{gathered}
\end{array}} \right.
$

where $u\equiv \,\,   {F_{20b}}\left( {{a_1},{a_2},b;{c_1},{c_2};x,y,z} \right)$.

Particular solutions:

$
{u_1} = {F_{20b}}\left( {{a_1},{a_2},b;{c_1},{c_2};x,y,z} \right),
$

$
{u_2} = {x^{1 - {c_1}}}{F_{20b}}\left( {{a_1},{a_2},b;{c_1},{c_2};x,y,z} \right),
$

$
{u_3} = {y^{1 - {c_2}}}{F_{20b}}\left( {1 - {c_2} + {a_1},{a_2},1 - {c_2} + b;{c_1},2 - {c_2};x,y,z} \right),
$

$
{u_4} = {x^{1 - {c_1}}}{y^{1 - {c_2}}}{F_{20b}}\left( {2 - {c_1} - {c_2} + {a_1},{a_2},2 - {c_1} - {c_2} + b;2 - {c_1},2 - {c_2};x,y,z} \right).
$

\bigskip

\begin{equation}
{F_{20c}}\left( {{a_1},{a_2},b;{c_1},{c_2};x,y,z} \right) = \sum\limits_{m,n,p = 0}^\infty  {} \frac{{{{\left( {{a_1}} \right)}_{m + n + p}}{{\left( {{a_2}} \right)}_{m + n}}{{\left( b \right)}_{p - m}}}}{{{{\left( {{c_1}} \right)}_n}{{\left( {{c_2}} \right)}_p}}}\frac{x^m}{m!}\frac{y^n}{n!}\frac{z^p}{p!},
\end{equation}

region of convergence:
$$
\begin{gathered}
  \left\{ {\sqrt r  + \sqrt s  < 1,\,\,\,t < \min \left\{ {U\left( {{w_1}} \right),U\left( {{w_2}} \right)} \right\}} \right\}, \hfill \\
  {P^ \pm }\left( w \right) = {w^3} - 2\left( {2 + s \pm 2r} \right)sw + 2\left( {1 \mp r} \right){s^2}, \hfill \\
  {w_1}:\,\,{\rm{the\,root\,in}}\,\left( {s,2\sqrt s  - s} \right)\,{\rm{of}}\,{P^ + }\left( w \right) = 0, \hfill \\
  {w_2}:\,\,{\rm{the\,smaller\,root\,in}}\,\left( {s,2\sqrt s  - s} \right)\,{\rm{of}}\,{P^ - }\left( w \right) = 0, \hfill \\
  U\left( w \right) = \frac{{{{\left[ {4s - {{\left( {s + w} \right)}^2}} \right]}^2}\left( {w - s} \right)}}{{16r{s^2}\left( {w + s} \right)}}. \hfill \\
\end{gathered}
$$

System of partial differential equations:

$
\left\{ {\begin{array}{*{20}{l}}
  \begin{gathered}
  x\left( {1 + x} \right){u_{xx}} + 2xy{u_{xy}} - \left( {1 - x} \right)z{u_{xz}} + {y^2}{u_{yy}} + yz{u_{yz}} \hfill \\
\,\,\,\,\,\,\,\,\,    + \left[ {1 - b + \left( {{a_1} + {a_2} + 1} \right)x} \right]{u_x} + \left( {{a_1} + {a_2} + 1} \right)y{u_y}
   + {a_2}z{u_z} + {a_1}{a_2}u = 0, \hfill \\
\end{gathered}  \\
  \begin{gathered}
  y\left( {1 - y} \right){u_{yy}} - {x^2}{u_{xx}} - 2xy{u_{xy}} - xz{u_{xz}} - yz{u_{yz}} \hfill \\
  \,\,\,\,\,\,\,\,\,  + \left[ {{c_1} - \left( {{a_1} + {a_2} + 1} \right)y} \right]{u_y} - \left( {{a_1} + {a_2} + 1} \right)x{u_x} - {a_2}z{u_z} - {a_1}{a_2}u = 0, \hfill \\
\end{gathered}  \\
  \begin{gathered}
  z\left( {1 - z} \right){u_{zz}} + {x^2}{u_{xx}} + xy{u_{xy}} - yz{u_{yz}} \hfill \\\,\,\,\,\,\,\,\,\,
  + \left[ {{c_2} - \left( {{a_1} + b + 1} \right)z} \right]{u_z} + \left( {{a_1} - b + 1} \right)x{u_x}
  - by{u_y} - {a_1}bu = 0, \hfill \\
\end{gathered}
\end{array}} \right.
$

where $u\equiv \,\,   {F_{20c}}\left( {{a_1},{a_2},b;{c_1},{c_2};x,y,z} \right)$.

Particular solutions:

$
{u_1} = {F_{20c}}\left( {{a_1},{a_2},b;{c_1},{c_2};x,y,z} \right),
$

$
{u_2} = {y^{1 - {c_1}}}{F_{20c}}\left( {1 - {c_1} + {a_1},1 - {c_1} + {a_2},b;2 - {c_1},{c_2};x,y,z} \right),
$

$
{u_3} = {z^{1 - {c_2}}}{F_{20c}}\left( {1 - {c_2} + {a_1},{a_2},1 - {c_2} + b;{c_1},2 - {c_2};x,y,z} \right),
$

$
{u_4} = {y^{1 - {c_1}}}{z^{1 - {c_2}}}{F_{20c}}\left( {2 - {c_1} - {c_2} + {a_1},1 - {c_1} + {a_2},1 - {c_2} + b;2 - {c_1},2 - {c_2};x,y,z} \right).
$

\bigskip

\begin{equation}
{F_{20d}}\left( {a,{b_1},{b_2};c;x,y,z} \right) = \sum\limits_{m,n,p = 0}^\infty  {} \frac{{{{\left( a \right)}_{m + n + p}}{{\left( {{b_1}} \right)}_{m + n - p}}{{\left( {{b_2}} \right)}_{p - m}}}}{{{{\left( c \right)}_n}}}\frac{x^m}{m!}\frac{y^n}{n!}\frac{z^p}{p!},
\end{equation}

region of convergence:
$$
\begin{gathered}
  \left\{ {r + t < 1,\,\,\,s < \min \left\{ {U\left( {{w_1}} \right), - U\left( {{w_2}} \right)} \right\}} \right\}, \hfill \\
  P\left( w \right) = r{w^3} - \left( {1 + t} \right)w + 2t, \hfill \\
  {w_1}:\,\,{\rm{the\,greater\,positive\,root\,of}}\,\,P\left( w \right) = 0, \hfill \\
  {w_2}:\,\,{\rm{the\,smaller\,positive\,root\,\,of}}\,P\left( w \right) = 0, \hfill \\
  U\left( w \right) = \left( {1 - \frac{{2t}}{w}} \right){\left( {1 - \frac{1}{w}} \right)^2}. \hfill \\
\end{gathered}
$$

System of partial differential equations:

$
\left\{ {\begin{array}{*{20}{l}}
  \begin{gathered}
  x\left( {1 + x} \right){u_{xx}} + 2xy{u_{xy}} - z{u_{xz}} + {y^2}{u_{yy}} - {z^2}{u_{zz}} \hfill \\
 \,\,\,\,\,\,\,\,\,   + \left[ {1 - {b_2} + \left( {a + {b_1} + 1} \right)x} \right]{u_x} + \left( {a + {b_1} + 1} \right)y{u_y}
    - \left( {a - {b_1} + 1} \right)z{u_z} + a{b_1}u = 0, \hfill \\
\end{gathered}  \\
  \begin{gathered}
  y\left( {1 - y} \right){u_{yy}} - 2xy{u_{xy}} - {x^2}{u_{xx}} + {z^2}{u_{zz}} \hfill \\
 \,\,\,\,\,\,\,\,\,   + \left[ {c - \left( {a + {b_1} + 1} \right)y} \right]{u_y} - \left( {a + {b_1} + 1} \right)x{u_x}
   + \left( {a - {b_1} + 1} \right)z{u_z} - a{b_1}u = 0, \hfill \\
\end{gathered}  \\
  \begin{gathered}
  z\left( {1 + z} \right){u_{zz}} - xy{u_{xy}} - x{u_{xz}} - y\left( {1 - z} \right){u_{yz}} - {x^2}{u_{xx}} \hfill \\
 \,\,\,\,\,\,\,\,\,   + \left[ {1 - {b_1} + \left( {a + {b_2} + 1} \right)z} \right]{u_z} - \left( {a - {b_2} + 1} \right)x{u_x}
    + {b_2}y{u_y} + a{b_2}u = 0, \hfill \\
\end{gathered}
\end{array}} \right.
$

where $u\equiv \,\,   {F_{20d}}\left( {a,{b_1},{b_2};c;x,y,z} \right)$.

Particular solutions:

$
{u_1} = {F_{20d}}\left( {a,{b_1},{b_2};c;x,y,z} \right),
$

$
{u_2} = {y^{1 - c}}{F_{20d}}\left( {1 - c + a,1 - c + {b_1},{b_2};2 - c;x,y,z} \right).
$

\bigskip

\begin{equation}
{F_{21a}}\left( {{a_1},{a_2},{a_3};{c_1},{c_2};x,y,z} \right)
  = \sum\limits_{m,n,p = 0}^\infty  {} \frac{{{{\left( {{a_1}} \right)}_{m + n + p}}{{\left( {{a_2}} \right)}_{m + n}}{{\left( {{a_3}} \right)}_p}}}{{{{\left( {{c_1}} \right)}_m}{{\left( {{c_2}} \right)}_{n + p}}}}\frac{x^m}{m!}\frac{y^n}{n!}\frac{z^p}{p!},
\end{equation}

first appearance of this function in the literature, and old notation: [24], $F_{14}$;\,\,[33], $F_F$,

region of convergence:
$$
 \left\{ {\sqrt r  + \sqrt s  < 1,\,\,\,t < \frac{1}{2}\left( {1 + s - r + \sqrt {{{\left( {1 + s - r} \right)}^2} - 4s} } \right)} \right\}.
$$

System of partial differential equations:

$
\left\{ {\begin{array}{*{20}{l}}
  \begin{gathered}
  x\left( {1 - x} \right){u_{xx}} - 2xy{u_{xy}} - xz{u_{xz}} - {y^2}{u_{yy}} - yz{u_{yz}} \hfill \\
  \,\,\,\,\,\,\,\,\,  + \left[ {{c_1} - \left( {{a_1} + {a_2} + 1} \right)x} \right]{u_x} - \left( {{a_1} + {a_2} + 1} \right)y{u_y}
   - {a_2}z{u_z} - {a_1}{a_2}u = 0, \hfill \\
\end{gathered}  \\
  \begin{gathered}
  y\left( {1 - y} \right){u_{yy}} - 2xy{u_{xy}} - {x^2}{u_{xx}} - xz{u_{xz}} + \left( {1 - y} \right)z{u_{yz}} \hfill \\
  \,\,\,\,\,\,\,\,\,  + \left[ {{c_2} - \left( {{a_1} + {a_2} + 1} \right)y} \right]{u_y} - \left( {{a_1} + {a_2} + 1} \right)x{u_x}
   - {a_2}z{u_z} - {a_1}{a_2}u = 0, \hfill \\
\end{gathered}  \\
  \begin{gathered}
  z\left( {1 - z} \right){u_{zz}} - xz{u_{xz}} + y\left( {1 - z} \right){u_{yz}}
   + \left[ {{c_2} - \left( {{a_1} + {a_3} + 1} \right)z} \right]{u_z} - {a_3}x{u_x} - {a_3}y{u_y} - {a_1}{a_3}u = 0, \hfill \\
\end{gathered}
\end{array}} \right.
$

where $u\equiv \,\,   {F_{21a}}\left( {{a_1},{a_2},{a_3};{c_1},{c_2};x,y,z} \right)$.

Particular solutions:

$
{u_1} = {F_{21a}}\left( {{a_1},{a_2},{a_3};{c_1},{c_2};x,y,z} \right),
$

$
{u_2} = {F_{21a}}\left( {1 - {c_1} + {a_1},1 - {c_1} + {a_2},{a_3};2 - {c_1},{c_2};x,y,z} \right).
$

\bigskip

\begin{equation}
{F_{21b}}\left( {{a_1},{a_2},{a_3};c;x,y,z} \right)
  = \sum\limits_{m,n,p = 0}^\infty  \frac{{{{\left( {{a_1}} \right)}_{m + n + p}}{{\left( {{a_2}} \right)}_{m + n}}{{\left( {{a_3}} \right)}_{p - m}}}}{{{{\left( c \right)}_{n + p}}}}\frac{x^m}{m!}\frac{y^n}{n!}\frac{z^p}{p!},
\end{equation}

region of convergence:
$$
 \left\{ {t < 1,\,\,\,r + 2\sqrt {rt}  < 1,\,\,\,s < 1 + r - 2\sqrt {r + rt} } \right\}.
$$

System of partial differential equations:

$
\left\{ \begin{array}{*{20}{l}}
  \begin{gathered}
  x\left( {1 + x} \right){u_{xx}} + 2xy{u_{xy}} - \left( {1 - x} \right)z{u_{xz}} + {y^2}{u_{yy}} + yz{u_{yz}} \hfill \\
 \,\,\,\,\,\,\,\,\,   + \left[ {1 - {a_3} + \left( {{a_1} + {a_2} + 1} \right)x} \right]{u_x} + \left( {{a_1} + {a_2} + 1} \right)y{u_y}
    + {a_2}z{u_z} + {a_1}{a_2}u = 0, \hfill \\
\end{gathered}  \\
  \begin{gathered}
  y\left( {1 - y} \right){u_{yy}} - 2xy{u_{xy}} + \left( {1 - y} \right)z{u_{yz}} - {x^2}{u_{xx}} - xz{u_{xz}} \hfill \\
 \,\,\,\,\,\,\,\,\,   + \left[ {c - \left( {{a_1} + {a_2} + 1} \right)y} \right]{u_y} - \left( {{a_1} + {a_2} + 1} \right)x{u_x}
   - {a_2}z{u_z} - {a_1}{a_2}u = 0, \hfill \\
\end{gathered}  \\
  \begin{gathered}
  z\left( {1 - z} \right){u_{zz}} + xy{u_{xy}} + y\left( {1 - z} \right){u_{yz}} + {x^2}{u_{xx}} \hfill \\
 \,\,\,\,\,\,\,\,\,   + \left[ {c - \left( {{a_1} + {a_3} + 1} \right)z} \right]{u_z} + \left( {{a_1} - {a_3} + 1} \right)x{u_x}
    - {a_3}y{u_y} - {a_1}{a_3}u = 0, \hfill \\
\end{gathered}
\end{array} \right.
$

where $u\equiv \,\,   {F_{21b}}\left( {{a_1},{a_2},{a_3};c;x,y,z} \right)$.

\bigskip

\begin{equation}
{F_{22a}}\left( {{a_1},{a_2};{c_1},{c_2},{c_3};x,y,z} \right)  = \sum\limits_{m,n,p = 0}^\infty  {} \frac{{{{\left( {{a_1}} \right)}_{m + n + p}}{{\left( {{a_2}} \right)}_{m + n + p}}}}{{{{\left( {{c_1}} \right)}_m}{{\left( {{c_2}} \right)}_n}{{\left( {{c_3}} \right)}_p}}}\frac{x^m}{m!}\frac{y^n}{n!}\frac{z^p}{p!},
\end{equation}

first appearance of this function in the literature, and old notation: [24], $F_{C}^{(3)}$, $F_5$,

region of convergence:
$$
\left\{ {\sqrt r  + \sqrt s  + \sqrt t  < 1} \right\}.
$$

System of partial differential equations:

$
\left\{ {\begin{array}{*{20}{l}}
  \begin{gathered}
  x\left( {1 - x} \right){u_{xx}} - 2xy{u_{xy}} - 2xz{u_{xz}} - 2yz{u_{yz}} - {y^2}{u_{yy}} - {z^2}{u_{zz}} \hfill \\
  \,\,\,\,\,\,\,\,\,  + \left[ {{c_1} - \left( {{a_1} + {a_2} + 1} \right)x} \right]{u_x} - \left( {{a_1} + {a_2} + 1} \right)\left( {y{u_y} + z{u_z}} \right)
    - {a_1}{a_2}u = 0, \hfill \\
\end{gathered}  \\
  \begin{gathered}
  y\left( {1 - y} \right){u_{yy}} - 2xy{u_{xy}} - 2xz{u_{xz}} - 2yz{u_{yz}} - {x^2}{u_{xx}} - {z^2}{u_{zz}} \hfill \\
 \,\,\,\,\,\,\,\,\,   + \left[ {{c_2} - \left( {{a_1} + {a_2} + 1} \right)y} \right]{u_y} - \left( {{a_1} + {a_2} + 1} \right)\left( {x{u_x} + z{u_z}} \right)
    - {a_1}{a_2}u = 0 \hfill \\
\end{gathered}  \\
  \begin{gathered}
  z\left( {1 - z} \right){u_{zz}} - 2xy{u_{xy}} - 2xz{u_{xz}} - 2yz{u_{yz}} - {x^2}{u_{xx}} - {y^2}{u_{yy}} \hfill \\
 \,\,\,\,\,\,\,\,\,   + \left[ {{c_3} - \left( {{a_1} + {a_2} + 1} \right)z} \right]{u_z} - \left( {{a_1} + {a_2} + 1} \right)\left( {y{u_y} + x{u_x}} \right)
   - {a_1}{a_2}u = 0, \hfill \\
\end{gathered}
\end{array}} \right.
$

where $u\equiv \,\,   {F_{22a}}\left( {{a_1},{a_2};{c_1},{c_2},{c_3};x,y,z} \right) $.

Particular solutions:

$
{u_1} = {F_{22a}}\left( {{a_1},{a_2};{c_1},{c_2},{c_3};x,y,z} \right),
$

$
{u_2} = {x^{1 - {c_1}}}{F_{22a}}\left( {1 - {c_1} + {a_1},1 - {c_1} + {a_2};2 - {c_1},{c_2},{c_3};x,y,z} \right),
$

$
{u_3} = {y^{1 - {c_2}}}{F_{22a}}\left( {1 - {c_2} + {a_1},1 - {c_2} + {a_2};{c_1},2 - {c_2},{c_3};x,y,z} \right),
$

$
{u_4} = {z^{1 - {c_3}}}{F_{22a}}\left( {1 - {c_3} + {a_1},1 - {c_3} + {a_2};{c_1},{c_2},2 - {c_3};x,y,z} \right),
$

$
{u_5} = {x^{1 - {c_1}}}{y^{1 - {c_2}}}{F_{22a}}\left( {2 - {c_1} - {c_2} + {a_1},2 - {c_1} - {c_2} + {a_2};{c_1},{c_2},{c_3};x,y,z} \right),
$

$
{u_6} = {y^{1 - {c_2}}}{z^{1 - {c_3}}}{F_{22a}}\left( {2 - {c_2} - {c_3} + {a_1},2 - {c_2} - {c_3} + {a_2};{c_1},2 - {c_2},2 - {c_3};x,y,z} \right),
$

$
{u_7} = {x^{1 - {c_1}}}{z^{1 - {c_3}}}{F_{22a}}\left( {2 - {c_1} - {c_3} + {a_1},2 - {c_1} - {c_3} + {a_2};2 - {c_1},{c_2},2 - {c_3};x,y,z} \right),
$

$
  {u_8} = {x^{1 - {c_1}}}{y^{1 - {c_2}}}{z^{1 - {c_3}}}\times$
    
    $ \,\,\,\,\,\,\,\,\,\times{F_{22a}}\left( {3 - {c_1} - {c_2} - {c_3} + {a_1},3 - {c_1} - {c_2} - {c_3} + {a_2};2 - {c_1},2 - {c_2},2 - {c_3};x,y,z} \right). $

\bigskip

\begin{equation}
{F_{23b}}\left( {{a_1},{a_2},{a_3},{b_1},{b_2};c;x,y,z} \right)\hfill \\
   = \sum\limits_{m,n,p = 0}^\infty  {} \frac{{{{\left( {{a_1}} \right)}_m}{{\left( {{a_2}} \right)}_n}{{\left( {{a_3}} \right)}_n}{{\left( {{b_1}} \right)}_{m - n}}{{\left( {{b_2}} \right)}_{2p - m}}}}{{{{\left( c \right)}_p}}}\frac{x^m}{m!}\frac{y^n}{n!}\frac{z^p}{p!},
\end{equation}

region of convergence:
$$
\left\{ {t < \frac{1}{4},\,\,\,r < \frac{1}{{1 + 2\sqrt t }},\,\,\,s < \frac{1}{{1 + r + 2r\sqrt t }}} \right\}.
$$

System of partial differential equations:

$
\left\{ {\begin{array}{*{20}{l}}
  \begin{gathered}
  x\left( {1 + x} \right){u_{xx}} - xy{u_{xy}} - 2z{u_{xz}}
   + \left[ {1 - {b_2} + \left( {{a_1} + {b_1} + 1} \right)x} \right]{u_x} - {a_1}y{u_y} + {a_1}{b_1}u = 0, \hfill \\
\end{gathered}  \\
  {y\left( {1 + y} \right){u_{yy}} - x{u_{xy}} + \left[ {1 - {b_1} + \left( {{a_2} + {a_3} + 1} \right)y} \right]{u_y} + {a_2}{a_3}u = 0,} \\
  \begin{gathered}
  z\left( {1 - 4z} \right){u_{zz}} + 4xz{u_{xz}} - {x^2}{u_{xx}}
   + \left[ {c - 2\left( {2{b_2} + 3} \right)z} \right]{u_z} + 2{b_2}x{u_x} - {b_2}\left( {1 + {b_2}} \right)u = 0, \hfill \\
\end{gathered}
\end{array}} \right.
$

where $u\equiv \,\,   {F_{23b}}\left( {{a_1},{a_2},{a_3},{b_1},{b_2};c;x,y,z} \right)$.

Particular solutions:

$
{u_1} = {F_{23b}}\left( {{a_1},{a_2},{a_3},{b_1},{b_2};c;x,y,z} \right),
$

$
{u_2} = {z^{1 - c}}{F_{23b}}\left( {{a_1},{a_2},{a_3},{b_1},2 - 2c + {b_2};2 - c;x,y,z} \right).
$

\bigskip

\begin{equation}
{F_{23c}}\left( {{a_1},{a_2},{a_3},{b_1},{b_2};c;x,y,z} \right)\hfill \\
   = \sum\limits_{m,n,p = 0}^\infty  {} \frac{{{{\left( {{a_1}} \right)}_m}{{\left( {{a_2}} \right)}_m}{{\left( {{a_3}} \right)}_n}{{\left( {{b_1}} \right)}_{n - p}}{{\left( {{b_2}} \right)}_{2p - m}}}}{{{{\left( c \right)}_n}}}\frac{x^m}{m!}\frac{y^n}{n!}\frac{z^p}{p!},
\end{equation}

region of convergence:
$$
  \left\{ {s < 1,\,\,\,t < \frac{1}{{4\left( {1 + s} \right)}},\,\,\,r < \frac{1}{{1 + 2\sqrt {t + st} }}} \right\}.
$$

System of partial differential equations:

$
\left\{ {\begin{array}{*{20}{l}}
  x\left( {1 + x} \right){u_{xx}} - 2z{u_{xz}} + \left[ {1 - {b_2} + \left( {{a_1} + {a_2} + 1} \right)x} \right]{u_x}
  + {a_1}{a_2}u = 0, \\
  y\left( {1 - y} \right){u_{yy}} + yz{u_{yz}} + \left[ {c - \left( {{a_3} + {b_1} + 1} \right)y} \right]{u_y} + {a_3}z{u_z}
   - {a_3}{b_1}u = 0, \\
  \begin{gathered}
  z\left( {1 + 4z} \right){u_{zz}} - 4xz{u_{xz}} + {x^2}{u_{xx}} - y{u_{yz}} \hfill\\\,\,\,\,\,\,\,\,\,
   + \left[ {1 - {b_1} + 2\left( {2{b_2} + 3} \right)z} \right]{u_z} - 2{b_2}x{u_x} + {b_2}\left( {1 + {b_2}} \right)u = 0, \hfill \\
\end{gathered}
\end{array}} \right.
$

where $u\equiv \,\,   {F_{23c}}\left( {{a_1},{a_2},{a_3},{b_1},{b_2};c;x,y,z} \right)$.

Particular solutions:

$
{u_1} = {F_{23c}}\left( {{a_1},{a_2},{a_3},{b_1},{b_2};c;x,y,z} \right),
$

$
{u_2} = {y^{1 - c}}{F_{23c}}\left( {{a_1},{a_2},1 - c + {a_3},1 - c + {b_1},{b_2};2 - c;x,y,z} \right).
$

\bigskip

\begin{equation}
{F_{23d}}\left( {{a_1},{a_2},{b_1},{b_2},{b_3};x,y,z} \right)\hfill \\
   = \sum\limits_{m,n,p = 0}^\infty  {} {{{{\left( {{a_1}} \right)}_n}{{\left( {{a_2}} \right)}_n}{{\left( {{b_1}} \right)}_{m - n}}{{\left( {{b_2}} \right)}_{m - p}}{{\left( {{b_3}} \right)}_{2p - m}}}}\frac{x^m}{m!}\frac{y^n}{n!}\frac{z^p}{p!},
\end{equation}

region of convergence:
$$
\left\{ {r < 1,\,\,\,s < \frac{1}{{1 + r}},\,\,\,t < \min \left\{ {\frac{1}{4},\frac{{1 - r}}{{{r^2}}},\frac{{\left( {1 - s} \right)\left( {1 - s - rs} \right)}}{{{r^2}{s^2}}}} \right\}} \right\}.
$$

System of partial differential equations:

$
\left\{ {\begin{array}{*{20}{l}}
  \begin{gathered}
  x\left( {1 + x} \right){u_{xx}} - xy{u_{xy}} - \left( {2 + x} \right)z{u_{xz}} + yz{u_{yz}} \hfill\\
  \,\,\,\,\,\,\,\,\, + \left[ {1 - {b_3} + \left( {{b_1} + {b_2} + 1} \right)x} \right]{u_x} - {b_2}y{u_y} - {b_1}z{u_z}
   + {b_1}{b_2}u = 0, \hfill \\
\end{gathered}  \\
  y\left( {1 + y} \right){u_{yy}} - x{u_{xy}} + \left[ {1 - {b_1} + \left( {{a_1} + {a_2} + 1} \right)y} \right]{u_y}
   + {a_1}{a_2}u = 0, \\
  \begin{gathered}
  z\left( {1 + 4z} \right){u_{zz}} - x\left( {1 + 4z} \right){u_{xz}} + {x^2}{u_{xx}}\hfill \\  \,\,\,\,\,\,\,\,\,
   + \left[ {1 - {b_2} + 2\left( {2{b_3} + 3} \right)z} \right]{u_z} - 2{b_3}x{u_x} + {b_3}\left( {{b_3} + 1} \right)u = 0, \hfill \\
\end{gathered}
\end{array}} \right.
$

where $u\equiv \,\,   {F_{23d}}\left( {{a_1},{a_2},{b_1},{b_2},{b_3};x,y,z} \right)$.

\bigskip

\begin{equation}
{F_{23e}}\left( {{a_1},{a_2},{b_1},{b_2},{b_3};x,y,z} \right) \hfill \\
  = \sum\limits_{m,n,p = 0}^\infty  {} {{{{\left( {{a_1}} \right)}_m}{{\left( {{a_2}} \right)}_n}{{\left( {{b_1}} \right)}_{m - n}}{{\left( {{b_2}} \right)}_{n - p}}{{\left( {{b_3}} \right)}_{2p - m}}}}\frac{x^m}{m!}\frac{y^n}{n!}\frac{z^p}{p!},
\end{equation}

region of convergence:
$$
 \left\{ {r < 1,\,\,\,s < \frac{1}{{1 + r}},\,\,\,t < \min \left\{ {\frac{1}{{4\left( {1 + s} \right)}},\frac{{{{\left( {1 - r} \right)}^2}}}{{4{r^2}}},\frac{{1 - s - rs}}{{{r^2}{s^2}}}} \right\}} \right\}.
$$

System of partial differential equations:

$
\left\{ {\begin{array}{*{20}{l}}
  x\left( {1 + x} \right){u_{xx}} - xy{u_{xy}} - 2z{u_{xz}} + \left[ {1 - {b_3} + \left( {{a_1} + {b_1} + 1} \right)x} \right]{u_x}
  - {a_1}y{u_y} + {a_1}{b_1}u = 0, \\
  y\left( {1 + y} \right){u_{yy}} - x{u_{xy}} - yz{u_{yz}} + \left[ {1 - {b_1} + \left( {{a_2} + {b_2} + 1} \right)y} \right]{u_y}
   - {a_2}z{u_z} + {a_2}{b_2}u = 0, \\
  \begin{gathered}
  z\left( {1 + 4z} \right){u_{zz}} - 4xz{u_{xz}} - y{u_{yz}} + {x^2}{u_{xx}}\hfill\\  \,\,\,\,\,\,\,\,\,
  + \left[ {1 - {b_2} + 2\left( {2{b_3} + 3} \right)z} \right]{u_z} - 2{b_3}x{u_x} + {b_3}\left( {{b_3} + 1} \right)u = 0, \hfill \\
\end{gathered}
\end{array}} \right.
$

where $u\equiv \,\,   {F_{23e}}\left( {{a_1},{a_2},{b_1},{b_2},{b_3};x,y,z} \right)$.

\bigskip

\begin{equation}
{F_{24b}}\left( {{a_1},{a_2},{a_3},{a_4},b;c;x,y,z} \right) \hfill \\
  = \sum\limits_{m,n,p = 0}^\infty  {} \frac{{{{\left( {{a_1}} \right)}_m}{{\left( {{a_2}} \right)}_m}{{\left( {{a_3}} \right)}_n}{{\left( {{a_4}} \right)}_n}{{\left( b \right)}_{2p - m - n}}}}{{{{\left( c \right)}_p}}}{x^m}\frac{x^m}{m!}\frac{y^n}{n!}\frac{z^p}{p!},
\end{equation}

region of convergence:
$$
\left\{ {t < \frac{1}{4},\,\,\,\max \left\{ {r,s} \right\} < \frac{1}{{1 + 2\sqrt t }}} \right\}.
$$

System of partial differential equations:

$
\left\{ {\begin{array}{*{20}{l}}
  x\left( {1 + x} \right){u_{xx}} + y{u_{xy}} - 2z{u_{xz}} + \left[ {1 - b + \left( {{a_1} + {a_2} + 1} \right)x} \right]{u_x}
   + {a_1}{a_2}u = 0, \\
  y\left( {1 + y} \right){u_{yy}} + x{u_{xy}} - 2z{u_{yz}} + \left[ {1 - b + \left( {{a_3} + {a_4} + 1} \right)y} \right]{u_y}
   + {a_3}{a_4}u = 0, \\
  \begin{gathered}
  z\left( {1 - 4z} \right){u_{zz}} - 2xy{u_{xy}} + 4xz{u_{xz}} + 4yz{u_{yz}} - {x^2}{u_{xx}} - {y^2}{u_{yy}} \hfill \\
  \,\,\,\,\,\,\,\,\,  + \left[ {c - 2\left( {2b + 3} \right)z} \right]{u_z} + 2bx{u_x} + 2by{u_y} - b\left( {b + 1} \right)u = 0, \hfill \\
\end{gathered}
\end{array}} \right.
$

where $u\equiv \,\,   {F_{24b}}\left( {{a_1},{a_2},{a_3},{a_4},b;c;x,y,z} \right)$.

Particular solutions:

$
{u_1} = {F_{24b}}\left( {{a_1},{a_2},{a_3},{a_4},b;c;x,y,z} \right),
$

$
{u_2} = {z^{1 - c}}{F_{24b}}\left( {{a_1},{a_2},{a_3},{a_4},2 - 2c + b;2 - c;x,y,z} \right).
$

\bigskip

\begin{equation}
{F_{24c}}\left( {{a_1},{a_2},{a_3},{b_1},{b_2};x,y,z} \right)\hfill \\
   = \sum\limits_{m,n,p = 0}^\infty  {} {{{{\left( {{a_1}} \right)}_m}{{\left( {{a_2}} \right)}_n}{{\left( {{a_3}} \right)}_n}{{\left( {{b_1}} \right)}_{m - p}}{{\left( {{b_2}} \right)}_{2p - m - n}}}}\frac{x^m}{m!}\frac{y^n}{n!}\frac{z^p}{p!},
\end{equation}

region of convergence:
$$
  \left\{ {t < \frac{1}{4},\,\,\,r < \frac{{\sqrt {1 + 4t}  - 1}}{{2t}},\,\,\,s < \frac{1}{{1 + 2\sqrt t }}} \right\}.
$$

System of partial differential equations:

$
\left\{ {\begin{array}{*{20}{l}}
  \begin{gathered}
  x\left( {1 + x} \right){u_{xx}} + y{u_{xy}} - \left( {2 + x} \right)z{u_{xz}}
   + \left[ {1 - {b_2} + \left( {{a_1} + {b_1} + 1} \right)x} \right]{u_x} - {a_1}z{u_z} + {a_1}{b_1}u = 0, \hfill \\
\end{gathered}  \\
  {y\left( {1 + y} \right){u_{yy}} + x{u_{xy}} - 2z{u_{yz}} + \left[ {1 - {b_2} + \left( {{a_2} + {a_3} + 1} \right)y} \right]{u_y} + {a_2}{a_3}u = 0,} \\
  \begin{gathered}
  z\left( {1 + 4z} \right){u_{zz}} + 2xy{u_{xy}} - x\left( {1 + 4z} \right){u_{xz}} - 4yz{u_{yz}} + {x^2}{u_{xx}} + {y^2}{u_{yy}} + \hfill \\
  \,\,\,\,\,\,\,\,\,  + \left[ {1 - {b_1} + 2\left( {2{b_2} + 3} \right)z} \right]{u_z} - 2{b_2}x{u_x} - 2{b_2}y{u_y} + {b_2}\left( {{b_2} + 1} \right)u = 0, \hfill \\
\end{gathered}
\end{array}} \right.
$

where $u\equiv \,\,   {F_{24c}}\left( {{a_1},{a_2},{a_3},{b_1},{b_2};x,y,z} \right)$.

\bigskip

\begin{equation}
{F_{25b}}\left( {{a_1},{a_2},{a_3},{a_4},b;c;x,y,z} \right)\hfill \\
   = \sum\limits_{m,n,p = 0}^\infty  {} \frac{{{{\left( {{a_1}} \right)}_m}{{\left( {{a_2}} \right)}_m}{{\left( {{a_3}} \right)}_n}{{\left( {{a_4}} \right)}_n}{{\left( b \right)}_{2p - m}}}}{{{{\left( c \right)}_{n + p}}}}\frac{x^m}{m!}\frac{y^n}{n!}\frac{z^p}{p!},
\end{equation}

region of convergence:
$$
 \left\{ {s < 1,\,\,\,t < \frac{1}{4},\,\,\,r < \frac{1}{{1 + 2\sqrt t }}} \right\}.
$$

System of partial differential equations:

$
\left\{ {\begin{array}{*{20}{l}}
  x\left( {1 + x} \right){u_{xx}} - 2z{u_{xz}} + \left[ {1 - b + \left( {{a_1} + {a_2} + 1} \right)x} \right]{u_x}
   + {a_1}{a_2}u = 0, \\
  {y\left( {1 - y} \right){u_{yy}} + z{u_{yz}} + \left[ {c - \left( {{a_3} + {a_4} + \,1} \right)y} \right]{u_y} - {a_3}{a_4}u = 0,} \\
  \begin{gathered}
  z\left( {1 - 4z} \right){u_{zz}} + 4xz{u_{xz}} + y{u_{yz}} - {x^2}{u_{xx}}
   + \left[ {c - 2\left( {2b + 3} \right)z} \right]{u_z} + 2bx{u_x} - b\left( {b + 1} \right)u = 0, \hfill \\
\end{gathered}
\end{array}} \right.
$

where $u\equiv \,\,   {F_{25b}}\left( {{a_1},{a_2},{a_3},{a_4},b;c;x,y,z} \right)$.

\bigskip

\begin{equation}
{F_{25c}}\left( {{a_1},{a_2},{a_3},{b_1},{b_2};x,y,z} \right) \hfill \\
  = \sum\limits_{m,n,p = 0}^\infty  {} {{{{\left( {{a_1}} \right)}_m}{{\left( {{a_2}} \right)}_n}{{\left( {{a_3}} \right)}_n}{{\left( {{b_1}} \right)}_{2p - m}}{{\left( {{b_2}} \right)}_{m - n - p}}}}\frac{x^m}{m!}\frac{y^n}{n!}\frac{z^p}{p!},
\end{equation}

region of convergence:
$$
\left\{ {r < 1,\,\,\,s < \frac{1}{{1 + r}},\,\,\,t < \min \left\{ {\frac{1}{4},\frac{{1 - r}}{{{r^2}}},\frac{{1 - s - rs}}{{{r^2}s}}} \right\}} \right\}.
$$

System of partial differential equations:

$
\left\{ {\begin{array}{*{20}{l}}
  \begin{gathered}
  x\left( {1 + x} \right){u_{xx}} - xy{u_{xy}} - \left( {2 + x} \right)z{u_{xz}}\hfill\\ \,\,\,\,\,\,\,\,\,
   + \left[ {1 - {b_1} + \left( {{a_1} + {b_2} + 1} \right)x} \right]{u_x} - {a_1}y{u_y} - {a_1}z{u_z} + {a_1}{b_2}u = 0, \hfill \\
\end{gathered}  \\
  y\left( {1 + y} \right){u_{yy}} - x{u_{xy}} + z{u_{yz}} + \left[ {1 - {b_2} + \left( {{a_2} + {a_3} + 1} \right)y} \right]{u_y}
   + {a_2}{a_3}u = 0, \\
  \begin{gathered}
  z\left( {1 + 4z} \right){u_{zz}} - x\left( {1 + 4z} \right){u_{xz}} + y{u_{yz}} + {x^2}{u_{xx}}\hfill\\ \,\,\,\,\,\,\,\,\,
 + \left[ {1 - {b_2} + 2\left( {2{b_1} + 3} \right)z} \right]{u_z} - 2{b_1}x{u_x} + {b_1}\left( {{b_1} + 1} \right)u = 0, \hfill \\
\end{gathered}
\end{array}} \right.
$

where $u\equiv \,\,   {F_{25c}}\left( {{a_1},{a_2},{a_3},{b_1},{b_2};x,y,z} \right)$.

\bigskip

\begin{equation}
{F_{26b}}\left( {{a_1},{a_2},{a_3},b;{c_1},{c_2};x,y,z} \right)  = \sum\limits_{m,n,p = 0}^\infty  {} \frac{{{{\left( {{a_1}} \right)}_{m + n}}{{\left( {{a_2}} \right)}_m}{{\left( {{a_3}} \right)}_n}{{\left( b \right)}_{2p - m}}}}{{{{\left( {{c_1}} \right)}_n}{{\left( {{c_2}} \right)}_p}}}\frac{x^m}{m!}\frac{y^n}{n!}\frac{z^p}{p!},
\end{equation}

region of convergence:
$$
\left\{ {s < 1,\,\,\,t < \frac{1}{4},\,\,\,r < \frac{{1 - s}}{{1 + 2\sqrt t }}} \right\}.
$$

System of partial differential equations:

$
\left\{ {\begin{array}{*{20}{l}}
  \begin{gathered}
  x\left( {1 + x} \right){u_{xx}} + xy{u_{xy}} - 2z{u_{xz}}
  + \left[ {1 - b + \left( {{a_1} + {a_2} + 1} \right)x} \right]{u_x} + {a_2}y{u_y} + {a_1}{a_2}u = 0, \hfill \\
\end{gathered}  \\
  y\left( {1 - y} \right){u_{yy}} - xy{u_{xy}} + \left[ {{c_1} - \left( {{a_1} + {a_3} + 1} \right)y} \right]{u_y} - {a_3}x{u_x}
   - {a_1}{a_3}u = 0, \\
  \begin{gathered}
  z\left( {1 - 4z} \right){u_{zz}} + 4xz{u_{xz}} - {x^2}{u_{xx}}
  + \left[ {{c_2} - 2\left( {2b + 3} \right)z} \right]{u_z} + 2bx{u_x} - b\left( {b + 1} \right)u = 0, \hfill \\
\end{gathered}
\end{array}} \right.
$

where $u\equiv \,\,   {F_{26b}}\left( {{a_1},{a_2},{a_3},b;{c_1},{c_2};x,y,z} \right)$.

Particular solutions:

$
{u_1} = {F_{26b}}\left( {{a_1},{a_2},{a_3},b;{c_1},{c_2};x,y,z} \right),
$

$
{u_2} = {y^{1 - {c_1}}}{F_{26b}}\left( {1 - {c_1} + {a_1},{a_2},1 - {c_1} + {a_3},b;2 - {c_1},{c_2};x,y,z} \right),
$

$
{u_3} = {z^{1 - {c_2}}}{F_{26b}}\left( {{a_1},{a_2},{a_3},2 - 2{c_2} + b;{c_1},2 - {c_2};x,y,z} \right),
$

$
{u_4} = {y^{1 - {c_1}}}{z^{1 - {c_2}}}{F_{26b}}\left( {1 - {c_1} + {a_1},{a_2},1 - {c_1} + {a_3},2 - 2{c_2} + b;2 - {c_1},2 - {c_2};x,y,z} \right).
$

\bigskip

\begin{equation}
{F_{26c}}\left( {{a_1},{a_2},{b_1},{b_2};c;x,y,z} \right) \hfill \\
  = \sum\limits_{m,n,p = 0}^\infty  {} \frac{{{{\left( {{a_1}} \right)}_m}{{\left( {{a_2}} \right)}_n}{{\left( {{b_1}} \right)}_{2p - m}}{{\left( {{b_2}} \right)}_{m + n - p}}}}{{{{\left( c \right)}_n}}}\frac{x^m}{m!}\frac{y^n}{n!}\frac{z^p}{p!},
\end{equation}

region of convergence:
$$
  \left\{ {r + s < 1,\,\,\,t < \min \left\{ {\frac{1}{{4\left( {1 + s} \right)}},\frac{{1 - r - s}}{{{r^2}}}} \right\}} \right\}.
$$

System of partial differential equations:

$
\left\{ {\begin{array}{*{20}{l}}
  \begin{gathered}
  x\left( {1 + x} \right){u_{xx}} + xy{u_{xy}} - \left( {2 + x} \right)z{u_{xz}}\hfill\\  \,\,\,\,\,\,\,\,\,
   + \left[ {1 - {b_1} + \left( {{a_1} + {b_2} + 1} \right)x} \right]{u_x} + {a_1}y{u_y} - {a_1}z{u_z} + {a_1}{b_2}u = 0, \hfill \\
\end{gathered}  \\
  \begin{gathered}
  y\left( {1 - y} \right){u_{yy}} - xy{u_{xy}} + yz{u_{yz}}
   + \left[ {c - \left( {{a_2} + {b_2} + 1} \right)y} \right]{u_y} - {a_2}x{u_x} + {a_2}z{u_z} - {a_2}{b_2}u = 0, \hfill \\
\end{gathered}  \\
  \begin{gathered}
  z\left( {1 + 4z} \right){u_{zz}} - x\left( {1 + 4z} \right){u_{xz}} - y{u_{yz}} + {x^2}{u_{xx}}\hfill \\  \,\,\,\,\,\,\,\,\,
   + \left[ {1 - {b_2} + 2\left( {2{b_1} + 3} \right)z} \right]{u_z} - 2{b_1}x{u_x} + {b_1}\left( {{b_1} + 1} \right)u = 0, \hfill \\
\end{gathered}
\end{array}} \right.
$

where $u\equiv \,\,   {F_{26c}}\left( {{a_1},{a_2},{b_1},{b_2};c;x,y,z} \right)$.

Particular solutions:

$
{u_1} = {F_{26c}}\left( {{a_1},{a_2},{b_1},{b_2};c;x,y,z} \right),
$

$
{u_2} = {y^{1 - c}}{F_{26c}}\left( {{a_1},1 - c + {a_2},{b_1},1 - c + {b_2};2 - c;x,y,z} \right).
$

\bigskip

\begin{equation}
{F_{26d}}\left( {{a_1},{a_2},{b_1},{b_2};c;x,y,z} \right) \hfill \\
  = \sum\limits_{m,n,p = 0}^\infty  {} \frac{{{{\left( {{a_1}} \right)}_{m + n}}{{\left( {{a_2}} \right)}_n}{{\left( {{b_1}} \right)}_{2p - m}}{{\left( {{b_2}} \right)}_{m - p}}}}{{{{\left( c \right)}_n}}}\frac{x^m}{m!}\frac{y^n}{n!}\frac{z^p}{p!},
\end{equation}

region of convergence:
$$
 \left\{ {s < 1,\,\,\,t < \frac{1}{4},\,\,\,r < \frac{{\left( {1 - s} \right)\left( {\sqrt {1 + 4t}  - 1} \right)}}{{2t}}} \right\}.
$$

System of partial differential equations:

$
\left\{ {\begin{array}{*{20}{l}}
  \begin{gathered}
  x\left( {1 + x} \right){u_{xx}} + xy{u_{xy}} - \left( {2 + x} \right)z{u_{xz}} - yz{u_{yz}} \hfill \\
 \,\,\,\,\,\,\,\,\,   + \left[ {1 - {b_1} + \left( {{a_1} + {b_2} + 1} \right)x} \right]{u_x} + {b_2}y{u_y} - {a_1}z{u_z} + {a_1}{b_2}u = 0, \hfill \\
\end{gathered}  \\
  y\left( {1 - y} \right){u_{yy}} - xy{u_{xy}} + \left[ {c - \left( {{a_1} + {a_2} + 1} \right)y} \right]{u_y} - {a_2}x{u_x}
   - {a_1}{a_2}u = 0, \\
  \begin{gathered}
  z\left( {1 + 4z} \right){u_{zz}} - x\left( {1 + 4z} \right){u_{xz}} + {x^2}{u_{xx}}\hfill\\ \,\,\,\,\,\,\,\,\,
   + \left[ {1 - {b_2} + 2\left( {2{b_1} + 3} \right)z} \right]{u_z} - 2{b_1}x{u_x} + {b_1}\left( {{b_1} + 1} \right)u = 0, \hfill \\
\end{gathered}
\end{array}} \right.
$

where $u\equiv \,\,   {F_{26d}}\left( {{a_1},{a_2},{b_1},{b_2};c;x,y,z} \right)$.

Particular solutions:

$
{u_1} = {F_{26d}}\left( {{a_1},{a_2},{b_1},{b_2};c;x,y,z} \right),
$

$
{u_2} = {y^{1 - c}}{F_{26d}}\left( {1 - c + {a_1},1 - c + {a_2},{b_1},{b_2};2 - c;x,y,z} \right).
$

\bigskip

\begin{equation}
{F_{26e}}\left( {{a_1},{a_2},{b_1},{b_2};c;x,y,z} \right) \hfill \\
  = \sum\limits_{m,n,p = 0}^\infty  {} \frac{{{{\left( {{a_1}} \right)}_{m + n}}{{\left( {{a_2}} \right)}_m}{{\left( {{b_1}} \right)}_{2p - m}}{{\left( {{b_2}} \right)}_{n - p}}}}{{{{\left( c \right)}_n}}}\frac{x^m}{m!}\frac{y^n}{n!}\frac{z^p}{p!},
\end{equation}

region of convergence:
$$
    \left\{ {t < \frac{1}{4},\,\,\,r < \frac{1}{{1 + 2\sqrt t }},} \,\,\, {s < \min \left\{ {\frac{{1 - 4t}}{{4t}},\left( {1 - r} \right){\Psi _1}\left( {\frac{{{r^2}t}}{{{{\left( {1 - r} \right)}^2}}}} \right),\left( {1 - r} \right){\Psi _2}\left( {\frac{{{r^2}t}}{{{{\left( {1 - r} \right)}^2}}}} \right)} \right\}} \right\}.
$$

System of partial differential equations:

$
\left\{ {\begin{array}{*{20}{l}}
  \begin{gathered}
  x\left( {1 + x} \right){u_{xx}} + xy{u_{xy}} - 2z{u_{xz}}
   + \left[ {1 - {b_1} + \left( {{a_1} + {a_2} + 1} \right)x} \right]{u_x} + {a_2}y{u_y} + {a_1}{a_2}u = 0, \hfill \\
\end{gathered}  \\
  \begin{gathered}
  y\left( {1 - y} \right){u_{yy}} - xy{u_{xy}} + xz{u_{xz}} + yz{u_{yz}}\hfill\\ \,\,\,\,\,\,\,\,\,
   + \left[ {c - \left( {{a_1} + {b_2} + 1} \right)y} \right]{u_y} - {b_2}x{u_x} + {a_1}z{u_z} - {a_1}{b_2}u = 0 ,\hfill \\
\end{gathered}  \\
  \begin{gathered}
  z\left( {1 + 4z} \right){u_{zz}} - 4xz{u_{xz}} - y{u_{yz}} + {x^2}{u_{xx}}\hfill\\ \,\,\,\,\,\,\,\,\,
   + \left[ {1 - {b_2} + 2\left( {2{b_1} + 3} \right)z} \right]{u_z} - 2{b_1}x{u_x} + {b_1}\left( {{b_1} + 1} \right)u = 0, \hfill \\
\end{gathered}
\end{array}} \right.
$

where $u\equiv \,\,   {F_{26e}}\left( {{a_1},{a_2},{b_1},{b_2};c;x,y,z} \right)$.

Particular solutions:

$
{u_1} = {F_{26e}}\left( {{a_1},{a_2},{b_1},{b_2};c;x,y,z} \right),
$

$
{u_2} = {y^{1 - c}}{F_{26e}}\left( {1 - c + {a_1},{a_2},{b_1},1 - c + {b_2};2 - c;x,y,z}\right).
$

\bigskip

\begin{equation}
{F_{26f}}\left( {{a_1},{a_2},{b_1},{b_2};c;x,y,z} \right) \hfill \\
  = \sum\limits_{m,n,p = 0}^\infty  {} \frac{{{{\left( {{a_1}} \right)}_{m + n}}{{\left( {{a_2}} \right)}_n}{{\left( {{b_1}} \right)}_{2p - m}}{{\left( {{b_2}} \right)}_{m - n}}}}{{{{\left( c \right)}_p}}}\frac{x^m}{m!}\frac{y^n}{n!}\frac{z^p}{p!},
\end{equation}

region of convergence:
$$
\left\{ {t < \frac{1}{4},\,\,\,r < \frac{1}{{1 + 2\sqrt t }},\,\,\,s < \sqrt {1 + r\left( {1 + 2\sqrt t } \right)}  - \sqrt {r\left( {1 + 2\sqrt t } \right)} } \right\}.
$$

System of partial differential equations:

$\left\{ {\begin{array}{*{20}{l}}
  \begin{gathered}
  x\left( {1 + x} \right){u_{xx}} - 2z{u_{xz}} - {y^2}{u_{yy}}
   + \left[ {1 - {b_1} + \left( {{a_1} + {b_2} + 1} \right)x} \right]{u_x} - \left( {{a_1} - {b_2} + 1} \right)y{u_y}
   + {a_1}{b_2}u = 0, \hfill \\
\end{gathered}  \\
  y\left( {1 + y} \right){u_{yy}} - x\left( {1 - y} \right){u_{xy}} + \left[ {1 - {b_2} + \left( {{a_1} + {a_2} + 1} \right)y} \right]{u_y}
  + {a_2}x{u_x} + {a_1}{a_2}u = 0, \\
  \begin{gathered}
  z\left( {1 - 4z} \right){u_{zz}} + 4xz{u_{xz}} - {x^2}{u_{xx}}
   + \left[ {c - 2\left( {2{b_1} + 3} \right)z} \right]{u_z} + 2{b_1}x{u_x} - {b_1}\left( {{b_1} + 1} \right)u = 0, \hfill \\
\end{gathered}
\end{array}} \right.
$

where $u\equiv \,\,   {F_{26f}}\left( {{a_1},{a_2},{b_1},{b_2};c;x,y,z} \right)$.

Particular solutions:

$
{u_1} = {F_{26f}}\left( {{a_1},{a_2},{b_1},{b_2};c;x,y,z} \right),
$

$
{u_2} = {z^{1 - c}}{F_{26f}}\left( {{a_1},{a_2},2 - 2c + {b_1},{b_2};2 - c;x,y,z} \right).
$

\bigskip

\begin{equation}
{F_{26g}}\left( {{a},{b_1},{b_2},{b_3};x,y,z} \right)\hfill \\
   = \sum\limits_{m,n,p = 0}^\infty  {} {{{{\left( {{a}} \right)}_n}{{\left( {{b_1}} \right)}_{2p - m}}{{\left( {{b_2}} \right)}_{m - n}}{{\left( {{b_3}} \right)}_{m + n - p}}}}\frac{x^m}{m!}\frac{y^n}{n!}\frac{z^p}{p!},
\end{equation}

region of convergence:
$$
\begin{gathered}
  \left\{ {r < 1,\,\,\,s + 2\sqrt {rs}  < 1,} \,\,\, {t < \min \left\{ {\frac{{1 - r}}{{{r^2}}},\frac{1}{{1 + s}}{\Phi _1}\left( {\frac{{rs}}{{{{\left( {1 + s} \right)}^2}}}} \right),\frac{1}{{1 - s}}{\Phi _2}\left( {\frac{{rs}}{{{{\left( {1 - s} \right)}^2}}}} \right)} \right\}} \right\}. \hfill \\
\end{gathered}
$$

System of partial differential equations:

$
\left\{ {\begin{array}{*{20}{l}}
  \begin{gathered}
  x\left( {1 + x} \right){u_{xx}} + yz{u_{yz}} - \left( {2 + x} \right)z{u_{xz}} - {y^2}{u_{yy}} \hfill \\ \,\,\,\,\,\,\,\,\,
   + \left[ {1 - {b_1} + \left( {{b_2} + {b_3} + 1} \right)x} \right]{u_x} + \left( {{b_2} - {b_3} - 1} \right)y{u_y}
    - {b_2}z{u_z} + {b_2}{b_3}u = 0 ,\hfill \\
\end{gathered}  \\
  \begin{gathered}
  y\left( {1 + y} \right){u_{yy}} - x\left( {1 - y} \right){u_{xy}} - yz{u_{yz}}
   + \left[ {1 - {b_2} + \left( {{a} + {b_3} + 1} \right)y} \right]{u_y} + {a}x{u_x} - {a}z{u_z} + {a}{b_3}u = 0 ,\hfill \\
\end{gathered}  \\
  \begin{gathered}
  z\left( {1 + 4z} \right){u_{zz}} - x\left( {1 + 4z} \right){u_{xz}} - y{u_{yz}} + {x^2}{u_{xx}}\hfill\\ \,\,\,\,\,\,\,\,\,
   + \left[ {1 - {b_3} + 2\left( {2{b_1} + 3} \right)z} \right]{u_z} - 2{b_1}x{u_x} + {b_1}\left( {1 + {b_1}} \right)u = 0, \hfill \\
\end{gathered}
\end{array}} \right.
$

where $u\equiv \,\,   {F_{26g}}\left( {{a},{b_1},{b_2},{b_3};x,y,z} \right)$.

\bigskip

\begin{equation}
{F_{26h}}\left( {{a_1},{b_1},{b_2},{b_3};x,y,z} \right) \hfill \\
  = \sum\limits_{m,n,p = 0}^\infty  {} {{{{\left( {{a_1}} \right)}_{m + n}}{{\left( {{b_1}} \right)}_{2p - m}}{{\left( {{b_2}} \right)}_{m - n}}{{\left( {{b_3}} \right)}_{n - p}}}}\frac{x^m}{m!}\frac{y^n}{n!}\frac{z^p}{p!},
\end{equation}

region of convergence:
$$
\begin{gathered}
  \left\{ {r < 1,\,\,\,s + 2\sqrt {rs}  < 1,\,\,\,t < \min \left\{ {{U^ + }\left( {{w_1}} \right), - {U^ + }\left( {{w_2}} \right),{U^ - }\left( {{w_3}} \right)} \right\}} \right\}, \hfill \\
  {P_{rs}}\left( w \right) = 4{w^3} - 3\left( {2 + s} \right){w^2} + 2\left( {1 + s + rs} \right)w - rs, \hfill \\
  {w_1}:\,\,{\rm{the\,\,root\,in}}\,\,\left( {0,\frac{1}{2}} \right)\,{\rm{of}}\,{P_{rs}}\left( w \right) = 0, \hfill \\
  {w_1}:\,\,{\rm{the\,\,root\,in}}\,\,\left( {\frac{1}{2},1} \right)\,{\rm{of}}\,{P_{r, - s}}\left( w \right) = 0, \hfill \\
  {w_3}:\,\,{\rm{the\,\,root\,in}}\,\,\left( { - r + \sqrt {{r^2} + r} ,\frac{1}{2}} \right)\,{\rm{of}}\,{P_{ - r, - s}}\left( w \right) = 0, \hfill \\
  {U^ \pm }\left( w \right) = \frac{{w\left( {1 - w} \right)\left( {{w^2} \mp 2rw \pm r} \right)}}{{2\left( {1 - 2w} \right){r^2}s}}. \hfill \\
\end{gathered}
$$

System of partial differential equations:

$
\left\{ {\begin{array}{*{20}{l}}
  \begin{gathered}
  x\left( {1 + x} \right){u_{xx}} - 2z{u_{xz}} - {y^2}{u_{yy}}
   + \left[ {1 - {b_1} + \left( {{a_1} + {b_2} + 1} \right)x} \right]{u_x} - \left( {{a_1} - {b_2} + 1} \right)y{u_y} + {a_1}{b_2}u = 0 ,\hfill \\
\end{gathered}  \\
  \begin{gathered}
  y\left( {1 + y} \right){u_{yy}} - x\left( {1 - y} \right){u_{xy}} - xz{u_{xz}} - yz{u_{yz}} \hfill \\
 \,\,\,\,\,\,\,\,\,   + \left[ {1 - {b_2} + \left( {{a_1} + {b_3} + 1} \right)y} \right]{u_y} + {b_3}x{u_x} - {a_1}z{u_z} + {a_1}{b_3}u = 0 ,\hfill \\
\end{gathered}  \\
  \begin{gathered}
  z\left( {1 + 4z} \right){u_{zz}} - 4xz{u_{xz}} - y{u_{yz}} + {x^2}{u_{xx}}\hfill\\ \,\,\,\,\,\,\,\,\,
   + \left[ {1 - {b_3} + 2\left( {2{b_1} + 3} \right)z} \right]{u_z} - 2{b_1}x{u_x} + {b_1}\left( {1 + {b_1}} \right)u = 0, \hfill \\
\end{gathered}
\end{array}} \right.
$

where $u\equiv \,\,   {F_{26h}}\left( {{a_1},{b_1},{b_2},{b_3};x,y,z} \right)$.

\bigskip

\begin{equation}
{F_{27b}}\left( {{a_1},{a_2},{a_3},b;c;x,y,z} \right) \hfill \\
  = \sum\limits_{m,n,p = 0}^\infty  {} \frac{{{{\left( {{a_1}} \right)}_{m + n}}{{\left( {{a_2}} \right)}_m}{{\left( {{a_3}} \right)}_n}{{\left( b \right)}_{2p - m - n}}}}{{{{\left( c \right)}_p}}}\frac{x^m}{m!}\frac{y^n}{n!}\frac{z^p}{p!},
\end{equation}

region of convergence:
$$
 \left\{ {t < \frac{1}{4},\,\,\,\max \left(r,s\right) < \frac{1}{{1 + 2\sqrt t }}} \right\}.
$$

System of partial differential equations:

$
\left\{ {\begin{array}{*{20}{l}}
  \begin{gathered}
  x\left( {1 + x} \right){u_{xx}} + \left( {1 + x} \right)y{u_{xy}} - 2z{u_{xz}}
   + \left[ {1 - b + \left( {{a_1} + {a_2} + 1} \right)x} \right]{u_x} + {a_2}y{u_y} + {a_1}{a_2}u = 0, \hfill \\
\end{gathered}  \\
  \begin{gathered}
  y\left( {1 + y} \right){u_{yy}} + x\left( {1 + y} \right){u_{xy}} - 2z{u_{yz}}
   + \left[ {1 - b + \left( {{a_1} + {a_3} + 1} \right)y} \right]{u_y} + {a_3}x{u_x} + {a_1}{a_3}u = 0, \hfill \\
\end{gathered}  \\
  \begin{gathered}
  z\left( {1 - 4z} \right){u_{zz}} - 2xy{u_{xy}} + 4xz{u_{xz}}
  \hfill \\\,\,\,\,\,\,\,\,\,
  + 4yz{u_{yz}} - {x^2}{u_{xx}} - {y^2}{u_{yy}}
   + \left[ {c - 2\left( {2b + 3} \right)z} \right]{u_z} + 2bx{u_x} + 2by{u_y} - b\left( {b + 1} \right)u = 0, \hfill \\
\end{gathered}
\end{array}} \right.
$

where $u\equiv \,\,   {F_{27b}}\left( {{a_1},{a_2},{a_3},b;c;x,y,z} \right)$.

Particular solutions:

$
{u_1} = {F_{27b}}\left( {{a_1},{a_2},{a_3},b;c;x,y,z} \right),
$

$
{u_2} = {z^{1 - c}}{F_{27b}}\left( {{a_1},{a_2},{a_3},2 - 2c + b;2 - c;x,y,z} \right).
$

\bigskip

\begin{equation}
{F_{27c}}\left( {{a_1},{a_2},{b_1},{b_2};x,y,z} \right)\hfill \\
   = \sum\limits_{m,n,p = 0}^\infty  {} {{{{\left( {{a_1}} \right)}_m}{{\left( {{a_2}} \right)}_n}{{\left( {{b_1}} \right)}_{2p - n - m}}{{\left( {{b_2}} \right)}_{m + n - p}}}}\frac{x^m}{m!}\frac{y^n}{n!}\frac{z^p}{p!},
\end{equation}

region of convergence:
$$
\left\{ {r < 1,\,\,\,s < 1,\,\,\,t < \min \left\{ {\frac{1}{4},\frac{{1 - \max \left( {r,s} \right)}}{{{{\left( {\max \left( {r,s} \right)} \right)}^2}}}} \right\}} \right\}.
$$

System of partial differential equations:

$
\left\{ {\begin{array}{*{20}{l}}
  \begin{gathered}
  x\left( {1 + x} \right){u_{xx}} + \left( {1 + x} \right)y{u_{xy}} - \left( {2 + x} \right)z{u_{xz}}\hfill \\ \,\,\,\,\,\,\,\,\,
   + \left[ {1 - {b_1} + \left( {{a_1} + {b_2} + 1} \right)x} \right]{u_x} + {a_1}y{u_y} - {a_1}z{u_z} + {a_1}{b_2}u = 0, \hfill \\
\end{gathered}  \\
  \begin{gathered}
  y\left( {1 + y} \right){u_{yy}} + x\left( {1 + y} \right){u_{xy}} - \left( {2 + y} \right)z{u_{yz}}\hfill \\  \,\,\,\,\,\,\,\,\,
   + \left[ {1 - {b_1} + \left( {{a_2} + {b_2} + 1} \right)y} \right]{u_y} + {a_2}x{u_x} - {a_2}z{u_z} + {a_2}{b_2}u = 0, \hfill \\
\end{gathered}  \\
  \begin{gathered}
  z\left( {1 + 4z} \right){u_{zz}} + 2xy{u_{xy}} - x\left( {1 + 4z} \right){u_{xz}} - y\left( {1 + 4z} \right){u_{yz}}\hfill\\
 \,\,\,\,\,\,\,\,\,   + {x^2}{u_{xx}} + {y^2}{u_{yy}} +  \left[ {1 - {b_2} + 2\left( {2{b_1} + 3} \right)z} \right]{u_z} - 2{b_1}x{u_x}
    - 2{b_1}y{u_y} + {b_1}\left( {{b_1} + 1} \right)u = 0, \hfill \\
\end{gathered}
\end{array}} \right.
$

where $u\equiv \,\,   {F_{27c}}\left( {{a_1},{a_2},{b_1},{b_2};x,y,z} \right)$.

\bigskip

\begin{equation}
{F_{27d}}\left( {{a_1},{a_2},{b_1},{b_2};x,y,z} \right) \hfill \\
  = \sum\limits_{m,n,p = 0}^\infty  {} {{{{\left( {{a_1}} \right)}_{m + n}}{{\left( {{a_2}} \right)}_n}{{\left( {{b_1}} \right)}_{2p - n - m}}{{\left( {{b_2}} \right)}_{m - p}}}}\frac{x^m}{m!}\frac{y^n}{n!}\frac{z^p}{p!},
\end{equation}

region of convergence:
$$
 \left\{ {t < \frac{1}{4},\,\,\,s < \frac{1}{{1 + 2\sqrt t }},\,\,\,r < \min \left\{ {\frac{{\sqrt {1 + 4t}  - 1}}{{2t}},\frac{{{{\left( {1 - s} \right)}^2} - 4{s^2}t}}{{4st}}} \right\}} \right\}.
$$

System of partial differential equations:

$
\left\{ {\begin{array}{*{20}{l}}
  \begin{gathered}
  x\left( {1 + x} \right){u_{xx}} + \left( {1 + x} \right)y{u_{xy}} - \left( {2 + x} \right)z{u_{xz}} - yz{u_{yz}} \hfill \\ \,\,\,\,\,\,\,\,\,
   + \left[ {1 - {b_1} + \left( {{a_1} + {b_2} + 1} \right)x} \right]{u_x} + {b_2}y{u_y} - {a_1}z{u_z} + {a_1}{b_2}u = 0, \hfill \\
\end{gathered}  \\
  \begin{gathered}
  y\left( {1 + y} \right){u_{yy}} + x\left( {1 + y} \right){u_{xy}} - 2z{u_{yz}}
   + \left[ {1 - {b_1} + \left( {{a_1} + {a_2} + 1} \right)y} \right]{u_y} + {a_2}x{u_x} + {a_1}{a_2}u = 0, \hfill \\
\end{gathered}  \\
  \begin{gathered}
  z\left( {1 + 4z} \right){u_{zz}} + 2xy{u_{xy}} - x\left( {1 + 4z} \right){u_{xz}} - 4yz{u_{yz}}
   + {x^2}{u_{xx}} + {y^2}{u_{yy}} \hfill \\\,\,\,\,\,\,\,\,\,  +  \left[ {1 - {b_2} + 2\left( {2{b_1} + 3} \right)z} \right]{u_z} - 2{b_1}x{u_x}
    - 2{b_1}y{u_y} + {b_1}\left( {{b_1} + 1} \right)u = 0, \hfill \\
\end{gathered}
\end{array}} \right.
$

where $u\equiv \,\,   {F_{27d}}\left( {{a_1},{a_2},{b_1},{b_2};x,y,z} \right)$.

\bigskip

\begin{equation}
{F_{28b}}\left( {{a_1},{a_2},{a_3},b;c;x,y,z} \right) \hfill \\
  = \sum\limits_{m,n,p = 0}^\infty  {} \frac{{{{\left( {{a_1}} \right)}_{m + n}}{{\left( {{a_2}} \right)}_m}{{\left( {{a_3}} \right)}_n}{{\left( b \right)}_{2p - m}}}}{{{{\left( c \right)}_{n + p}}}}\frac{x^m}{m!}\frac{y^n}{n!}\frac{z^p}{p!},
\end{equation}

region of convergence:
$$
\left\{ {t < \frac{1}{4},\,\,\,r < \frac{1}{{1 + 2\sqrt t }},\,\,\,s < \frac{1}{2}\left( {1 - r} \right) + \frac{1}{2}\sqrt {{{\left( {1 - r} \right)}^2} - 4{r^2}t} } \right\}.
$$

System of partial differential equations:

$
\left\{ {\begin{array}{*{20}{l}}
  x\left( {1 + x} \right){u_{xx}} + xy{u_{xy}} - 2z{u_{xz}} + \left[ {1 - b + \left( {{a_1} + {a_2} + 1} \right)x} \right]{u_x}
   + {a_2}y{u_y} + {a_1}{a_2}u = 0, \\
  y\left( {1 - y} \right){u_{yy}} - xy{u_{xy}} + z{u_{yz}} + \left[ {c - \left( {{a_1} + {a_3} + 1} \right)y} \right]{u_y}
   - {a_3}x{u_x} - {a_1}{a_3}u = 0, \\
  z\left( {1 - 4z} \right){u_{zz}} + 4xz{u_{xz}} + y{u_{yz}} - {x^2}{u_{xx}}
   + \left[ {c - 2\left( {2b + 3} \right)z} \right]{u_z} + 2bx{u_x} - b\left( {b + 1} \right)u = 0,
\end{array}} \right.
$

where $u\equiv \,\,   {F_{28b}}\left( {{a_1},{a_2},{a_3},b;c;x,y,z} \right)$.

\bigskip

\begin{equation}
{F_{28c}}\left( {{a_1},{a_2},{b_1},{b_2};x,y,z} \right) \hfill \\
  = \sum\limits_{m,n,p = 0}^\infty  {} {{{{\left( {{a_1}} \right)}_{m + n}}{{\left( {{a_2}} \right)}_n}{{\left( {{b_1}} \right)}_{2p - m}}{{\left( {{b_2}} \right)}_{m - n - p}}}}\frac{x^m}{m!}\frac{y^n}{n!}\frac{z^p}{p!},
\end{equation}

region of convergence:
$$
\left\{ r < 1,\,\,\,s + 2\sqrt {rs}  < 1,\,\,\,t < \min \left\{\frac{1}{4},\frac{{1 - r}}{{{r^2}}},\frac{s}{{1 - s}}{\Phi _2}\left( {\frac{{rs}}{{{{\left( {1 - s} \right)}^2}}}} \right) \right\} \right\}.
$$

System of partial differential equations:

$
\left\{ {\begin{array}{*{20}{l}}
  \begin{gathered}
  x\left( {1 + x} \right){u_{xx}} - \left( {2 + x} \right)z{u_{xz}} - yz{u_{yz}} - {y^2}{u_{yy}} \hfill \\
  \,\,\,\,\,\,\,\,\,  + \left[ {1 - {b_1} + \left( {{a_1} + {b_2} + 1} \right)x} \right]{u_x} - \left( {{a_1} - {b_2} + 1} \right)y{u_y}
   - {a_1}z{u_z} + {a_1}{b_2}u = 0, \hfill \\
\end{gathered}  \\
  \begin{gathered}
  y\left( {1 + y} \right){u_{yy}} - x\left( {1 - y} \right){u_{xy}} + z{u_{yz}}
   + \left[ {1 - {b_2} + \left( {{a_1} + {a_2} + 1} \right)y} \right]{u_y} + {a_2}x{u_x} + {a_1}{a_2}u = 0, \hfill \\
\end{gathered}  \\
  \begin{gathered}
  z\left( {1 + 4z} \right){u_{zz}} - x\left( {1 + 4z} \right){u_{xz}} + y{u_{yz}} + {x^2}{u_{xx}}\hfill \\  \,\,\,\,\,\,\,\,\,
   + \left[ {1 - {b_2} + 2\left( {2{b_1} + 3} \right)z} \right]{u_z} - 2{b_1}x{u_x} + {b_1}\left( {{b_1} + 1} \right)u = 0, \hfill \\
\end{gathered}
\end{array}} \right.
$

where $u\equiv \,\,   {F_{28c}}\left( {{a_1},{a_2},{b_1},{b_2};x,y,z} \right)$.

\bigskip

\begin{equation}
{F_{29b}}\left( {{a_1},{a_2},{a_3},b;{c_1},{c_2};x,y,z} \right) \hfill \\
  = \sum\limits_{m,n,p = 0}^\infty  {} \frac{{{{\left( {{a_1}} \right)}_n}{{\left( {{a_2}} \right)}_p}{{\left( {{a_3}} \right)}_p}{{\left( b \right)}_{2m + n - p}}}}{{{{\left( {{c_1}} \right)}_m}{{\left( {{c_2}} \right)}_n}}}\frac{x^m}{m!}\frac{y^n}{n!}\frac{z^p}{p!},
\end{equation}

region of convergence:
$$
 \left\{ {s + 2\sqrt r  < 1,\,\,\,t < \frac{1}{{1 + s + 2\sqrt r }}} \right\}.
$$

System of partial differential equations:

$
\left\{ {\begin{array}{*{20}{l}}
  \begin{gathered}
  x\left( {1 - 4x} \right){u_{xx}} - 4xy{u_{xy}} + 4xz{u_{xz}} + 2yz{u_{yz}} - {y^2}{u_{yy}} - {z^2}{u_{zz}} \hfill \\
 \,\,\,\,\,\,\,\,\,   + \left[ {{c_1} - 2\left( {2b + 3} \right)x} \right]{u_x} - 2\left( {b + 1} \right)y{u_y} + 2bz{u_z}
    - b\left( {b + 1} \right)u = 0, \hfill \\
\end{gathered}  \\
  \begin{gathered}
  y\left( {1 - y} \right){u_{yy}} - 2xy{u_{xy}} + yz{u_{yz}}
   + \left[ {{c_2} - \left( {{a_1} + b + 1} \right)y} \right]{u_y} - 2{a_1}x{u_x} + {a_1}z{u_z} - {a_1}bu = 0, \hfill \\
\end{gathered}  \\
  z\left( {1 + z} \right){u_{zz}} - 2x{u_{xz}} - y{u_{yz}} + \left[ {1 - b + \left( {{a_2} + {a_3} + 1} \right)z} \right]{u_z}
   + {a_2}{a_3}u = 0,
\end{array}} \right.
$

where $u\equiv \,\,   {F_{29b}}\left( {{a_1},{a_2},{a_3},b;{c_1},{c_2};x,y,z} \right) $.

Particular solutions:

$
{u_1} = {F_{29b}}\left( {{a_1},{a_2},{a_3},b;{c_1},{c_2};x,y,z} \right),
$

$
{u_2} = {x^{1 - {c_1}}}{F_{29b}}\left( {{a_1},{a_2},{a_3},2 - 2{c_1} + b;2 - {c_1},{c_2};x,y,z} \right),
$

$
{u_3} = {y^{1 - {c_2}}}{F_{29b}}\left( {1 - {c_2} + {a_1},{a_2},{a_3},1 - {c_2} + b;{c_1},2 - {c_2};x,y,z} \right),
$

$
{u_4} = {x^{1 - {c_1}}}{y^{1 - {c_2}}}{F_{29b}}\left( {1 - {c_2} + {a_1},{a_2},{a_3},3 - 2{c_1} - {c_2} + b;{2-c_1},{2-c_2};x,y,z} \right).
$

\bigskip

\begin{equation}
{F_{29c}}\left( {{a_1},{a_2},{a_3},b;{c_1},{c_2};x,y,z} \right) \hfill \\
  = \sum\limits_{m,n,p = 0}^\infty  {} \frac{{{{\left( {{a_1}} \right)}_{2m + n}}{{\left( {{a_2}} \right)}_p}{{\left( {{a_3}} \right)}_p}{{\left( b \right)}_{n - p}}}}{{{{\left( {{c_1}} \right)}_m}{{\left( {{c_2}} \right)}_n}}}\frac{x^m}{m!}\frac{y^n}{n!}\frac{z^p}{p!},
\end{equation}

region of convergence:
$$
 \left\{ {s + 2\sqrt r  < 1,\,\,\,t < \frac{{1 - 2\sqrt r }}{{1 + s - 2\sqrt r }}} \right\}.
$$

System of partial differential equations:

$
\left\{ {\begin{array}{*{20}{l}}
  \begin{gathered}
  x\left( {1 - 4x} \right){u_{xx}} - 4xy{u_{xy}} - {y^2}{u_{yy}}\hfill\\ \,\,\,\,\,\,\,\,\,
   + \left[ {{c_1} - 2\left( {2{a_1} + 3} \right)x} \right]{u_x} - 2\left( {{a_1} + 1} \right)y{u_y} - {a_1}\left( {{a_1} + 1} \right)u = 0, \hfill \\
\end{gathered}  \\
  \begin{gathered}
  y\left( {1 - y} \right){u_{yy}} - 2xy{u_{xy}} + 2xz{u_{xz}} + yz{u_{yz}}\hfill\\ \,\,\,\,\,\,\,\,\,
  + \left[ {{c_2} - \left( {{a_1} + b + 1} \right)y} \right]{u_y} - 2bx{u_x} + {a_1}z{u_z} - {a_1}bu = 0, \hfill \\
\end{gathered}  \\
  {z\left( {1 + z} \right){u_{zz}} - y{u_{yz}} + \left[ {1 - b + \left( {{a_2} + {a_3} + 1} \right)z} \right]{u_z} + {a_2}{a_3}u = 0,}
\end{array}} \right.
$

where $u\equiv \,\,   {F_{29c}}\left( {{a_1},{a_2},{a_3},b;{c_1},{c_2};x,y,z} \right) $.

Particular solutions:

$
{u_1} = {F_{29c}}\left( {{a_1},{a_2},{a_3},b;{c_1},{c_2};x,y,z} \right),
$

$
{u_2} = {x^{1 - {c_1}}}{F_{29c}}\left( {2 - 2{c_1} + {a_1},{a_2},{a_3},b;2 - {c_1},{c_2};x,y,z} \right),
$

$
{u_3} = {y^{1 - {c_2}}}{F_{29c}}\left( {1 - {c_2} + {a_1},{a_2},{a_3},1 - {c_2} + b;{c_1},2 - {c_2};x,y,z} \right),
$

$
{u_4} = {x^{1 - {c_1}}}{y^{1 - {c_2}}}{F_{29c}}\left( {3 - 2{c_1} - {c_2} + {a_1},{a_2},{a_3},1 - {c_2} + b;2 - {c_1},2 - {c_2};x,y,z} \right).
$

\bigskip

\begin{equation}
{F_{29d}}\left( {{a_1},{a_2},{a_3},b;{c_1},{c_2};x,y,z} \right) \hfill \\
  = \sum\limits_{m,n,p = 0}^\infty  {} \frac{{{{\left( {{a_1}} \right)}_{2m + n}}{{\left( {{a_2}} \right)}_n}{{\left( {{a_3}} \right)}_p}{{\left( b \right)}_{p - n}}}}{{{{\left( {{c_1}} \right)}_m}{{\left( {{c_2}} \right)}_p}}}\frac{x^m}{m!}\frac{y^n}{n!}\frac{z^p}{p!},
\end{equation}

region of convergence:
$$
 \left\{ {r < \frac{1}{4},\,\,\,t < 1,\,\,\,s < \frac{{1 - 2\sqrt t }}{{1 + t}}} \right\}.
$$

System of partial differential equations:

$
\left\{ {\begin{array}{*{20}{l}}
  \begin{gathered}
  x\left( {1 - 4x} \right){u_{xx}} - 4xy{u_{xy}} - {y^2}{u_{yy}}  + \left[ {{c_1} - 2\left( {2{a_1} + 3} \right)x} \right]{u_x}
  - 2\left( {{a_1} + 1} \right)y{u_y} - {a_1}\left( {{a_1} + 1} \right)u = 0, \hfill \\
\end{gathered}  \\
  \begin{gathered}
  y\left( {1 + y} \right){u_{yy}} + 2xy{u_{xy}} - z{u_{yz}}
   + \left[ {1 - b + \left( {{a_1} + {a_2} + 1} \right)y} \right]{u_y} + 2{a_2}x{u_x} + {a_1}{a_2}u = 0, \hfill \\
\end{gathered}  \\
  z\left( {1 - z} \right){u_{zz}} + yz{u_{yz}} + \left[ {{c_2} - \left( {{a_3} + b + 1} \right)z} \right]{u_z}
   + {a_3}y{u_y} - {a_3}bu = 0,
\end{array}} \right.
$

where $u\equiv \,\,   {F_{29d}}\left( {{a_1},{a_2},{a_3},b;{c_1},{c_2};x,y,z} \right)$.

Particular solutions:

$
{u_1} = {F_{29d}}\left( {{a_1},{a_2},{a_3},b;{c_1},{c_2};x,y,z} \right),
$

$
{u_2} = {x^{1 - {c_1}}}{F_{29d}}\left( {{a_1},{a_2},{a_3},b;{c_1},{c_2};x,y,z} \right),
$

$
{u_3} = {z^{1 - {c_2}}}{F_{29d}}\left( 2-2c+{{a_1},{a_2},1 - {c_2} + {a_3},1 - {c_2} + b;{c_1},2 - {c_2};x,y,z} \right),
$

$
{u_4} = {x^{1 - {c_1}}}{z^{1 - {c_2}}}{F_{29d}}\left( {2 - 2{c_1} + {a_1},{a_2},1 - {c_2} + {a_3},1 - {c_2} + b;2 - {c_1},2 - {c_2};x,y,z} \right).
$

\bigskip

\begin{equation}
 {F_{29e}}\left( {{a_1},{a_2},{a_3},b;{c_1},{c_2};x,y,z} \right)\hfill \\
   = \sum\limits_{m,n,p = 0}^\infty  {} \frac{{{{\left( {{a_1}} \right)}_{2m + n}}{{\left( {{a_2}} \right)}_n}{{\left( {{a_3}} \right)}_p}{{\left( b \right)}_{p - m}}}}{{{{\left( {{c_1}} \right)}_n}{{\left( {{c_2}} \right)}_p}}}\frac{x^m}{m!}\frac{y^n}{n!}\frac{z^p}{p!},
\end{equation}

region of convergence:
$$
\left\{ {s < 1,\,\,\,t < 1,\,\,\,r < \frac{{{{\left( {1 - s} \right)}^2}}}{{4\left( {1 + t} \right)}}} \right\}.
$$

System of partial differential equations:

$
\left\{ {\begin{array}{*{20}{l}}
  \begin{gathered}
  x\left( {1 + 4x} \right){u_{xx}} + 4xy{u_{xy}} + {y^2}{u_{yy}} - z{u_{xz}}\\hfill \\ \,\,\,\,\,\,\,\,\,
   + \left[ {1 - b + 2\left( {2{a_1} + 3} \right)x} \right]{u_x} + 2\left( {{a_1} + 1} \right)y{u_y}
    + {a_1}\left( {{a_1} + 1} \right)u = 0, \hfill \\
\end{gathered}  \\
  y\left( {1 - y} \right){u_{yy}} - 2xy{u_{xy}} + \left[ {{c_1} - \left( {{a_1} + {a_2} + 1} \right)y} \right]{u_y}- 2{a_2}x{u_x}
  - {a_1}{a_2}u = 0, \\
  z\left( {1 - z} \right){u_{zz}} + xz{u_{xz}} + \left[ {{c_2} - \left( {{a_3} + b + 1} \right)z} \right]{u_z} + {a_3}x{u_x}
   - {a_3}bu = 0,
\end{array}} \right.
$

where $u\equiv \,\,    {F_{29e}}\left( {{a_1},{a_2},{a_3},b;{c_1},{c_2};x,y,z} \right)$.

Particular solutions:

$
{u_1} = {F_{29e}}\left( {{a_1},{a_2},{a_3},b;{c_1},{c_2};x,y,z} \right),
$

$
{u_2} = {y^{1 - {c_1}}}{F_{29e}}\left( {1 - {c_1} + {a_1},1 - {c_1} + {a_2},{a_3},b;2 - {c_1},{c_2};x,y,z} \right),
$

$
{u_3} = {z^{1 - {c_2}}}{F_{29e}}\left( {{a_1},{a_2},1 - {c_2} + {a_3},1 - {c_2} + b;{c_1},2 - {c_2};x,y,z} \right),
$

$
{u_4} = {y^{1 - {c_1}}}{z^{1 - {c_2}}}{F_{29e}}\left( {1 - {c_1} + {a_1},1 - {c_1} + {a_2},1 - {c_2} + {a_3},1 - {c_2} + b;2 - {c_1},2 - {c_2};x,y,z} \right).
$

\bigskip

\begin{equation}
{F_{29f}}\left( {{a_1},{a_2},{b_1},{b_2};c;x,y,z} \right)\hfill \\
   = \sum\limits_{m,n,p = 0}^\infty  {} \frac{{{{\left( {{a_1}} \right)}_p}{{\left( {{a_2}} \right)}_p}{{\left( {{b_1}} \right)}_{2m + n - p}}{{\left( {{b_2}} \right)}_{n - m}}}}{{{{\left( c \right)}_n}}}\frac{x^m}{m!}\frac{y^n}{n!}\frac{z^p}{p!},
\end{equation}

region of convergence:
$$
\begin{gathered}
    \left\{ {r < \frac{1}{4},\,\,\,t < \frac{1}{{1 + 2\sqrt r }},} \,\,\, s < \min \left\{{\Psi_1}\left( r \right),{\Psi _2}\left( r \right),\frac{{1 - t}}{t}{\Psi _1}\left( {\frac{{r{t^2}}}{{{{\left( {1 - t} \right)}^2}}}} \right),\frac{{1 - t}}{t}{\Psi_2}\left( {\frac{{r{t^2}}}{{{{\left( {1 - t} \right)}^2}}}} \right)  \right\}\right\}. \hfill \\
\end{gathered}
$$

System of partial differential equations:

$
\left\{ {\begin{array}{*{20}{l}}
  \begin{gathered}
  x\left( {1 + 4x} \right){u_{xx}} - \left( {1 - 4x} \right)y{u_{xy}} - 4xz{u_{xz}} - 2yz{u_{yz}} + {y^2}{u_{yy}}\hfill \\
  \,\,\,\,\,\,\,\,\,  + {z^2}{u_{zz}}  + \left[ {1 - {b_2} + 2\left( {2{b_1} + 3} \right)x} \right]{u_x} + 2\left( {{b_1} + 1} \right)y{u_y} - 2{b_1}z{u_z}
   + {b_1}\left( {{b_1} + 1} \right)u = 0, \hfill \\
\end{gathered}  \\
  \begin{gathered}
  y\left( {1 - y} \right){u_{yy}} - xy{u_{xy}} - xz{u_{xz}} + yz{u_{yz}} + 2{x^2}{u_{xx}} \hfill \\ \,\,\,\,\,\,\,\,\,
 + \left[ {c - \left( {{b_1} + {b_2} + 1} \right)y} \right]{u_y} + \left( {{b_1} - 2{b_2} + 2} \right)x{u_x}
  + {b_2}z{u_z} - {b_1}{b_2}u = 0 ,\hfill \\
\end{gathered}  \\
  z\left( {1 + z} \right){u_{zz}} - 2x{u_{xz}} - y{u_{yz}} + \left[ {1 - {b_1} + \left( {{a_1} + {a_2} + 1} \right)z} \right]{u_z}
   + {a_1}{a_2}u = 0,
\end{array}} \right.
$

where $u\equiv \,\,   {F_{29f}}\left( {{a_1},{a_2},{b_1},{b_2};c;x,y,z} \right)$.

Particular solutions:

$
{u_1} = {F_{29f}}\left( {{a_1},{a_2},{b_1},{b_2};c;x,y,z} \right),
$

$
{u_2} = {y^{1 - c}}{F_{29f}}\left( {{a_1},{a_2},1 - c + {b_1},1 - c + {b_2};2 - c;x,y,z} \right).
$

\bigskip

\begin{equation}
{F_{29g}}\left( {{a_1},{a_2},{b_1},{b_2};c;x,y,z} \right) \hfill \\
  = \sum\limits_{m,n,p = 0}^\infty  {} \frac{{{{\left( {{a_1}} \right)}_n}{{\left( {{a_2}} \right)}_p}{{\left( {{b_1}} \right)}_{2m + n - p}}{{\left( {{b_2}} \right)}_{p - m}}}}{{{{\left( c \right)}_n}}}\frac{x^m}{m!}\frac{y^n}{n!}\frac{z^p}{p!},
\end{equation}

region of convergence:
$$
\left\{ s < 1,\,\,\,t < \frac{1}{{1 + s}},\,\,\,r < \min \left\{\frac{(1-s)^2}{4},\frac{1-t-st}{{{t^2}}} \right\}\right\}.
$$

System of partial differential equations:

$
\left\{ {\begin{array}{*{20}{l}}
  \begin{gathered}
  x\left( {1 + 4x} \right){u_{xx}} + 4xy{u_{xy}} - \left( {1 + 4x} \right)z{u_{xz}} - 2yz{u_{yz}} + {y^2}{u_{yy}}\hfill \\
 \,\,\,\,\,\,\,\,\,   + {z^2}{u_{zz}}  + \left[ {1 - {b_2} + 2\left( {2{b_1} + 3} \right)x} \right]{u_x} + 2\left( {{b_1} + 1} \right)y{u_y}
    - 2{b_1}z{u_z} + {b_1}\left( {{b_1} + 1} \right)u = 0, \hfill \\
\end{gathered}  \\
  \begin{gathered}
  y\left( {1 - y} \right){u_{yy}} - 2xy{u_{xy}} + yz{u_{yz}}
   + \left[ {c - \left( {{a_1} + {b_1} + 1} \right)y} \right]{u_y} - 2{a_1}x{u_x} + {a_1}z{u_z} - {a_1}{b_1}u = 0 ,\hfill \\
\end{gathered}  \\
  \begin{gathered}
  z\left( {1 + z} \right){u_{zz}} - x\left( {2 + z} \right){u_{xz}} - y{u_{yz}}
   + \left[ {1 - {b_1} + \left( {{a_2} + {b_2} + 1} \right)z} \right]{u_z} - {a_2}x{u_x} + {a_2}{b_2}u = 0, \hfill \\
\end{gathered}
\end{array}} \right.
$

where $u\equiv \,\,   {F_{29g}}\left( {{a_1},{a_2},{b_1},{b_2};c;x,y,z} \right)$.

Particular solutions:

$
{u_1} = {F_{29g}}\left( {{a_1},{a_2},{b_1},{b_2};c;x,y,z} \right),
$

$
{u_2} = {y^{1 - c}}{F_{29g}}\left( {1 - c + {a_1},{a_2},1 - c + {b_1},{b_2};2 - c;x,y,z} \right).
$

\bigskip

\begin{equation}
{F_{29h}}\left( {{a_1},{a_2},{b_1},{b_2};c;x,y,z} \right)\hfill \\
   = \sum\limits_{m,n,p = 0}^\infty  {} \frac{{{{\left( {{a_1}} \right)}_n}{{\left( {{a_2}} \right)}_p}{{\left( {{b_1}} \right)}_{2m + n - p}}{{\left( {{b_2}} \right)}_{p - n}}}}{{{{\left( c \right)}_m}}}\frac{x^m}{m!}\frac{y^n}{n!}\frac{z^p}{p!},
\end{equation}

region of convergence:
$$
 \left\{ {s + 2\sqrt r  < 1,\,\,\,t < \frac{1}{{1 + 2\sqrt r }}} \right\}.
$$

System of partial differential equations:

$
\left\{ {\begin{array}{*{20}{l}}
  \begin{gathered}
  x\left( {1 - 4x} \right){u_{xx}} - 4xy{u_{xy}} + 4xz{u_{xz}} + 2yz{u_{yz}} - {y^2}{u_{yy}} - {z^2}{u_{zz}} \hfill \\
  \,\,\,\,\,\,\,\,\,  + \left[ {c - 2\left( {2{b_1} + 3} \right)x} \right]{u_x} - \left( {{b_1} + {b_1} + 2} \right)y{u_y} + 2{b_1}z{u_z}
   - {b_1}\left( {{b_1} + 1} \right)u = 0, \hfill \\
\end{gathered}  \\
  \begin{gathered}
  y\left( {1 + y} \right){u_{yy}} + 2xy{u_{xy}} - \left( {1 + y} \right)z{u_{yz}}\hfill \\  \,\,\,\,\,\,\,\,\,
   + \left[ {1 - {b_2} + \left( {{a_1} + {b_1} + 1} \right)y} \right]{u_y} + 2{a_1}x{u_x} - {a_1}z{u_z} + {a_1}{b_1}u = 0, \hfill \\
\end{gathered}  \\
  \begin{gathered}
  z\left( {1 + z} \right){u_{zz}} - 2x{u_{xz}} - y\left( {1 + z} \right){u_{yz}}
   + \left[ {1 - {b_1} + \left( {{a_2} + {b_2} + 1} \right)z} \right]{u_z} - {a_2}y{u_y} + {a_2}{b_2}u = 0, \hfill \\
\end{gathered}
\end{array}} \right.
$

where $u\equiv \,\,   {F_{29h}}\left( {{a_1},{a_2},{b_1},{b_2};c;x,y,z} \right)$.

Particular solutions:

$
{u_1} = {F_{29h}}\left( {{a_1},{a_2},{b_1},{b_2};c;x,y,z} \right),
$

$
{u_2} = {x^{1 - c}}{F_{29h}}\left( {{a_1},{a_2},2 - 2c + {b_1},{b_2};2 - c;x,y,z} \right).
$

\bigskip

\begin{equation}
{F_{29i}}\left( {{a_1},{a_2},{b_1},{b_2};c;x,y,z} \right) \hfill \\
  = \sum\limits_{m,n,p = 0}^\infty  {} \frac{{{{\left( {{a_1}} \right)}_{2m + n}}{{\left( {{a_2}} \right)}_p}{{\left( {{b_1}} \right)}_{n - p}}{{\left( {{b_2}} \right)}_{p - m}}}}{{{{\left( c \right)}_n}}}\frac{x^m}{m!}\frac{y^n}{n!}\frac{z^p}{p!},
\end{equation}

region of convergence:
$$
\begin{gathered}
   \left\{ {s < 1,\,\,\,t < \frac{1}{{1 + s}},\,} \,\,\, r < \min \left\{\frac{{{{\left( {1 - s} \right)}^2}}}{4},\frac{1}{{1 + t}}{\Theta _1}\left( {\frac{{st}}{{1 + t}}} \right),\frac{1}{{1 - t}}{\Theta _1}\left( {\frac{{st}}{{1 - t}}} \right) \right\} \right\}. \hfill \\
\end{gathered}
$$

System of partial differential equations:

$
\left\{ {\begin{array}{*{20}{l}}
  \begin{gathered}
  x\left( {1 + 4x} \right){u_{xx}} + 4xy{u_{xy}} - z{u_{xz}} + {y^2}{u_{yy}} \hfill \\
 \,\,\,\,\,\,\,\,\,   + \left[ {1 - {b_2} + 2\left( {2{a_1} + 3} \right)x} \right]{u_x} + 2\left( {{a_1} + 1} \right)y{u_y}
    + {a_1}\left( {{a_1} + 1} \right)u = 0, \hfill \\
\end{gathered}  \\
  \begin{gathered}
  y\left( {1 - y} \right){u_{yy}} - 2xy{u_{xy}} + 2xz{u_{xz}} + yz{u_{yz}}\hfill\\  \,\,\,\,\,\,\,\,\,
   + \left[ {c - \left( {{a_1} + {b_1} + 1} \right)y} \right]{u_y} - 2{b_1}x{u_x} + {a_1}z{u_z} - {a_1}{b_1}u = 0, \hfill \\
\end{gathered}  \\
  z\left( {1 + z} \right){u_{zz}} - xz{u_{xz}} - y{u_{yz}} + \left[ {1 - {b_1} + \left( {{a_2} + {b_2} + 1} \right)z} \right]{u_z}
   - {a_2}x{u_x} + {a_2}{b_2}u = 0,
\end{array}} \right.
$

where $u\equiv \,\,   {F_{29i}}\left( {{a_1},{a_2},{b_1},{b_2};c;x,y,z} \right)$.

Particular solutions:

$
{u_1} = {F_{29i}}\left( {{a_1},{a_2},{b_1},{b_2};c;x,y,z} \right),
$

$
{u_2} = {y^{1 - c}}{F_{29i}}\left( {1 - c + {a_1},{a_2},1 - c + {b_1},{b_2};2 - c;x,y,z} \right).
$

\bigskip

\begin{equation}
{F_{29j}}\left( {{a_1},{a_2},{b_1},{b_2};c;x,y,z} \right)\hfill \\
   = \sum\limits_{m,n,p = 0}^\infty  {} \frac{{{{\left( {{a_1}} \right)}_{2m + n}}{{\left( {{a_2}} \right)}_p}{{\left( {{b_1}} \right)}_{n - p}}{{\left( {{b_2}} \right)}_{p - n}}}}{{{{\left( c \right)}_m}}}\frac{x^m}{m!}\frac{y^n}{n!}\frac{z^p}{p!},
\end{equation}

region of convergence:
$$
 \left\{ {s + 2\sqrt r  < 1,\,\,\,t < 1} \right\}.
$$

System of partial differential equations:

$
\left\{ {\begin{array}{*{20}{l}}
  \begin{gathered}
  x\left( {1 - 4x} \right){u_{xx}} - 4xy{u_{xy}} - {y^2}{u_{yy}}
   + \left[ {c - 2\left( {2{a_1} + 3} \right)x} \right]{u_x} - 2\left( {{a_1} + 1} \right)y{u_y} - {a_1}\left( {{a_1} + 1} \right)u = 0, \hfill \\
\end{gathered}  \\
  \begin{gathered}
  y\left( {1 + y} \right){u_{yy}} + 2xy{u_{xy}} - 2xz{u_{xz}} - \left( {1 + y} \right)z{u_{yz}} \hfill \\\,\,\,\,\,\,\,\,\,
   + \left[ {1 - {b_2} + \left( {{a_1} + {b_1} + 1} \right)y} \right]{u_y} + 2{b_1}x{u_x} - {a_1}z{u_z} + {a_1}{b_1}u = 0, \hfill \\
\end{gathered}  \\
  z\left( {1 + z} \right){u_{zz}} - y\left( {1 + z} \right){u_{yz}} + \left[ {1 - {b_1} + \left( {{a_2} + {b_2} + 1} \right)z} \right]{u_z}
   - {a_2}y{u_y} + {a_2}{b_2}u = 0,
\end{array}} \right.
$

where $u\equiv \,\,   {F_{29j}}\left( {{a_1},{a_2},{b_1},{b_2};c;x,y,z} \right)$.

Particular solutions:

$
{u_1} = {F_{29j}}\left( {{a_1},{a_2},{b_1},{b_2};c;x,y,z} \right),
$

$
{u_2} = {x^{1 - c}}{F_{29j}}\left( {2 - 2c + {a_1},{a_2},{b_1},{b_2};2 - c;x,y,z} \right).
$

\bigskip

\begin{equation}
{F_{29k}}\left( {{a_1},{a_2},{b_1},{b_2};c;x,y,z} \right)\hfill \\
   = \sum\limits_{m,n,p = 0}^\infty  {} \frac{{{{\left( {{a_1}} \right)}_{2m + n}}{{\left( {{a_2}} \right)}_p}{{\left( {{b_1}} \right)}_{n - m}}{{\left( {{b_2}} \right)}_{p - n}}}}{{{{\left( c \right)}_p}}}\frac{x^m}{m!}\frac{y^n}{n!}\frac{z^p}{p!},
\end{equation}

region of convergence:
$$
 \left\{ r < \frac{1}{4},\,\,\,t < 1,\,\,\,s < \frac{1}{{1 + t}}\min \left\{{\Psi _1}\left( r \right),{\Psi _2}\left( r \right) \right\} \right\}.
$$

System of partial differential equations:

$
\left\{ {\begin{array}{*{20}{l}}
  x\left( {1 + 4x} \right){u_{xx}} - \left( {1 - 4x} \right)y{u_{xy}} + {y^2}{u_{yy}}\hfill\\  \,\,\,\,\,\,\,\,\,
   + \left[ {1 - {b_1} + 2\left( {2{a_1} + 3} \right)x} \right]{u_x} + 2\left( {{a_1} + 1} \right)y{u_y}
    + {a_1}\left( {{a_1} + 1} \right)u = 0,  \\
  y\left( {1 + y} \right){u_{yy}} + xy{u_{xy}} - z{u_{yz}} - 2{x^2}{u_{xx}}\hfill \\ \,\,\,\,\,\,\,\,\,
   + \left[ {1 - {b_2} + \left( {{a_1} + {b_1} + 1} \right)y} \right]{u_y} - \left( {{a_1} - 2{b_1} + 2} \right)x{u_x}
    + {a_1}{b_1}u = 0,   \\
  z\left( {1 - z} \right){u_{zz}} + yz{u_{yz}} + \left[ {c - \left( {{a_2} + {b_2} + 1} \right)z} \right]{u_z} + {a_2}y{u_y}
   - {a_2}{b_2}u = 0,
\end{array}} \right.
$

where $u\equiv \,\,   {F_{29k}}\left( {{a_1},{a_2},{b_1},{b_2};c;x,y,z} \right)$.

Particular solutions:

$
{u_1} = {F_{29k}}\left( {{a_1},{a_2},{b_1},{b_2};c;x,y,z} \right),
$

$
{u_2} = {z^{1 - c}}{F_{29k}}\left( {{a_1},1 - c + {a_2},{b_1},1 - c + {b_2};2 - c;x,y,z} \right).
$

\bigskip

\begin{equation}
{F_{29l}}\left( {{a_1},{a_2},{b_1},{b_2};c;x,y,z} \right) \hfill \\
  = \sum\limits_{m,n,p = 0}^\infty  {} \frac{{{{\left( {{a_1}} \right)}_{2m + n}}{{\left( {{a_2}} \right)}_n}{{\left( {{b_1}} \right)}_{p - n}}{{\left( {{b_2}} \right)}_{p - m}}}}{{{{\left( c \right)}_p}}}\frac{x^m}{m!}\frac{y^n}{n!}\frac{z^p}{p!},
\end{equation}

region of convergence:
$$
   \left\{ {s + 2\sqrt r  < 1,} \,\,\, t < \min \left\{1,\frac{{1 - 4r}}{{4r}},\frac{{1 - s}}{s}{\Psi _1}\left( {\frac{r}{{{{\left( {1 - s} \right)}^2}}}} \right),\frac{{1 - s}}{s}{\Psi _2}\left( \frac{r}{{{{\left( {1 - s} \right)}^2}}} \right)  \right\}\right\}.
$$

System of partial differential equations:

$
\left\{ {\begin{array}{*{20}{l}}
  \begin{gathered}
  x\left( {1 + 4x} \right){u_{xx}} + 4xy{u_{xy}} - z{u_{xz}} + {y^2}{u_{yy}} \hfill \\ \,\,\,\,\,\,\,\,\,
   + \left[ {1 - {b_2} + 2\left( {2{a_1} + 3} \right)x} \right]{u_x} + 2\left( {{a_1} + 1} \right)y{u_y}
    + {a_1}\left( {{a_1} + 1} \right)u = 0, \hfill \\
\end{gathered}  \\
  y\left( {1 + y} \right){u_{yy}} + 2xy{u_{xy}} - z{u_{yz}} + \left[ {1 - {b_1} + \left( {{a_1} + {a_2} + 1} \right)y} \right]{u_y}
   + 2{a_2}x{u_x}   + {a_1}{a_2}u = 0, \hfill \\
  \begin{gathered}
  z\left( {1 - z} \right){u_{zz}} - xy{u_{xy}} + xz{u_{xz}} + yz{u_{yz}}
  + \left[ {c - \left( {{b_1} + {b_2} + 1} \right)z} \right]{u_z} + {b_1}x{u_x} + {b_2}y{u_y} - {b_1}{b_2}u = 0, \hfill \\
\end{gathered}
\end{array}} \right.
$

where $u\equiv \,\,   {F_{29l}}\left( {{a_1},{a_2},{b_1},{b_2};c;x,y,z} \right)$.

Particular solutions:

$
{u_1} = {F_{29l}}\left( {{a_1},{a_2},{b_1},{b_2};c;x,y,z} \right),
$

$
{u_2} = {z^{1 - c}}{F_{29l}}\left( {{a_1},{a_2},1 - c + {b_1},1 - c + {b_2};2 - c;x,y,z} \right).
$

\bigskip

\begin{equation}
{F_{29m}}\left( {a,{b_1},{b_2},{b_3};x,y,z} \right) \hfill \\
  = \sum\limits_{m,n,p = 0}^\infty  {} {{{{\left( a \right)}_p}{{\left( {{b_1}} \right)}_{2m + n - p}}{{\left( {{b_2}} \right)}_{n - m}}{{\left( {{b_3}} \right)}_{p - n}}}}\frac{x^m}{m!}\frac{y^n}{n!}\frac{z^p}{p!},
\end{equation}

region of convergence:
$$
 \left\{ r < \frac{1}{4},\,\,\,t < \frac{1}{{1 + 2\sqrt r }},\,\,\,s < \min \left\{{\Psi _1}\left( r \right),{\Psi _2}\left( r \right),\frac{{{{\left( {1 - t} \right)}^2} - 4r{t^2}}}{{4r{t^3}}} \right\} \right\}.
$$

System of partial differential equations:

$
\left\{ {\begin{array}{*{20}{l}}
  \begin{gathered}
  x\left( {1 + 4x} \right){u_{xx}} - \left( {1 - 4x} \right)y{u_{xy}} - 4xz{u_{xz}} - 2yz{u_{yz}} \hfill \\
 \,\,\,\,\,\,\,\,\,  + {y^2}{u_{yy}} + {z^2}{u_{zz}}   + \left[ {1 - {b_2} + 2\left( {2{b_1} + 3} \right)x} \right]{u_x} + 2\left( {{b_1} + 1} \right)y{u_y}
  - 2{b_1}z{u_z} + {b_1}\left( {{b_1} + 1} \right)u = 0, \hfill \\
\end{gathered}  \\
  \begin{gathered}
  y\left( {1 + y} \right){u_{yy}} + xy{u_{xy}} + xz{u_{xz}} - \left( {1 + y} \right)z{u_{yz}} - 2{x^2}{u_{xx}} \hfill \\
  \,\,\,\,\,\,\,\,\,  + \left[ {1 - {b_3} + \left( {{b_1} + {b_2} + 1} \right)y} \right]{u_y} - \left( {{b_1} - 2{b_2} + 2} \right)x{u_x}
    - {b_2}z{u_z} + {b_1}{b_2}u = 0, \hfill \\
\end{gathered}  \\
  z\left( {1 + z} \right){u_{zz}} - 2x{u_{xz}} - y\left( {1 + z} \right){u_{yz}}
   + \left[ {1 - {b_1} + \left( {a + {b_3} + 1} \right)z} \right]{u_z} - ay{u_y} + a{b_3}u = 0,
\end{array}} \right.
$

where $u\equiv \,\,   {F_{29m}}\left( {a,{b_1},{b_2},{b_3};x,y,z} \right)$.

\bigskip

\begin{equation}
{F_{29n}}\left( {a,{b_1},{b_2},{b_3};x,y,z} \right) \hfill \\
  = \sum\limits_{m,n,p = 0}^\infty  {} {{{{\left( a \right)}_n}{{\left( {{b_1}} \right)}_{2m + n - p}}{{\left( {{b_2}} \right)}_{p - m}}{{\left( {{b_3}} \right)}_{p - n}}}}\frac{x^m}{m!}\frac{y^n}{n!}\frac{z^p}{p!},
\end{equation}

region of convergence:
$$
\left\{ {s + 2\sqrt r  < 1,\,\,\,t < \min \left\{ {\frac{{\sqrt {1 + 4r}  - 1}}{{2r}},\frac{{{{\left( {1 - s} \right)}^2} - 4r}}{{4rs}}} \right\}} \right\}.
$$

System of partial differential equations:

$
\left\{ {\begin{array}{*{20}{l}}
  \begin{gathered}
  x\left( {1 + 4x} \right){u_{xx}} + 4xy{u_{xy}} - \left( {1 + 4x} \right)z{u_{xz}} - 2yz{u_{yz}}+ {y^2}{u_{yy}}\hfill \\
 \,\,\,\,\,\,\,\,\,    + {z^2}{u_{zz}}    + \left[ {1 - {b_2} + 2\left( {2{b_1} + 3} \right)x} \right]{u_x} + 2\left( {{b_1} + 1} \right)y{u_y}
   - 2{b_1}z{u_z} + {b_1}\left( {{b_1} + 1} \right)u = 0, \hfill \\
\end{gathered}  \\
  \begin{gathered}
  y\left( {1 + y} \right){u_{yy}} + 2xy{u_{xy}} - \left( {1 + y} \right)z{u_{yz}}\hfill\\\,\,\,\,\,\,\,\,\,
   + \left[ {1 - {b_3} + \left( {a + {b_1} + 1} \right)y} \right]{u_y} + 2ax{u_x} - az{u_z} + a{b_1}u = 0, \hfill \\
\end{gathered}  \\
  \begin{gathered}
  z\left( {1 + z} \right){u_{zz}} + xy{u_{xy}} - x\left( {2 + z} \right){u_{xz}} - y\left( {1 + z} \right){u_{yz}}
   \hfill \\
  \,\,\,\,\,\,\,\,\,  + \left[ {1 - {b_1} + \left( {{b_2} + {b_3} + 1} \right)z} \right]{u_z} - {b_3}x{u_x} - {b_2}y{u_y} + {b_2}{b_3}u = 0, \hfill \\
\end{gathered}
\end{array}} \right.
$

where $u\equiv \,\,   {F_{29n}}\left( {a,{b_1},{b_2},{b_3};x,y,z} \right)$.

\bigskip

\begin{equation}
{F_{29o}}\left( {a,{b_1},{b_2},{b_3};x,y,z} \right)  = \sum\limits_{m,n,p = 0}^\infty  {} {{{{\left( a \right)}_{2m + n}}{{\left( {{b_1}} \right)}_{n - p}}{{\left( {{b_2}} \right)}_{p - m}}{{\left( {{b_3}} \right)}_{p - n}}}}\frac{x^m}{m!}\frac{y^n}{n!}\frac{z^p}{p!},
\end{equation}

region of convergence:
$$
\left\{ {s < 1,\,\,\,t < 1,\,\,\,r < \frac{{{{\left( {1 - s} \right)}^2}}}{{4\left( {1 + t} \right)}}} \right\}.
$$

System of partial differential equations:

$
\left\{ {\begin{array}{*{20}{l}}
  \begin{gathered}
  x\left( {1 + 4x} \right){u_{xx}} + 4xy{u_{xy}} - z{u_{xz}} + {y^2}{u_{yy}}\hfill\\\,\,\,\,\,\,\,\,\,
   + \left[ {1 - {b_2} + 2\left( {2a + 3} \right)x} \right]{u_x} + 2\left( {a + 1} \right)y{u_y} + a\left( {a + 1} \right)u = 0 ,\hfill \\
\end{gathered}  \\
  \begin{gathered}
  y\left( {1 + y} \right){u_{yy}} + 2xy{u_{xy}} - 2xz{u_{xz}} - \left( {1 + y} \right)z{u_{yz}}
  \hfill \\  \,\,\,\,\,\,\,\,\, + \left[ {1 - {b_3} + \left( {a + {b_1} + 1} \right)y} \right]{u_y} + 2{b_1}x{u_x} - az{u_z} + a{b_1}u = 0, \hfill \\
\end{gathered}  \\
  \begin{gathered}
  z\left( {1 + z} \right){u_{zz}} + xy{u_{xy}} - xz{u_{xz}} - y\left( {1 + z} \right){u_{yz}} \hfill \\\,\,\,\,\,\,\,\,\,
   + \left[ {1 - {b_1} + \left( {{b_2} + {b_3} + 1} \right)z} \right]{u_z} - {b_3}x{u_x} - {b_2}y{u_y} + {b_2}{b_3}u = 0, \hfill \\
\end{gathered}
\end{array}} \right.
$

where $u\equiv \,\,   {F_{29o}}\left( {a,{b_1},{b_2},{b_3};x,y,z} \right)$.

\bigskip

\begin{equation}
 {F_{30b}}\left( {{a_1},{a_2},{a_3},b;c;x,y,z} \right)  = \sum\limits_{m,n,p = 0}^\infty  {} \frac{{{{\left( {{a_1}} \right)}_n}{{\left( {{a_2}} \right)}_p}{{\left( {{a_3}} \right)}_p}{{\left( b \right)}_{2m + n - p}}}}{{{{\left( c \right)}_{m + n}}}}\frac{x^m}{m!}\frac{y^n}{n!}\frac{z^p}{p!},
\end{equation}

region of convergence:
$$
\begin{gathered}
    \left\{ {r < \frac{1}{4},\,}  \,\,{t < \frac{1}{{1 + 2\sqrt r }},\,\,\,s < \frac{1}{2}\min \left\{ 1+{\sqrt {1 - 4r} ,\frac{{1 - t + \sqrt {{{\left( {1 - t} \right)}^2} - 4r{t^2}} }}{t}} \right\}} \right\}. \hfill \\
\end{gathered}
$$

System of partial differential equations:

$
\left\{ {\begin{array}{*{20}{l}}
  \begin{gathered}
  x\left( {1 - 4x} \right){u_{xx}} + \left( {1 - 4x} \right)y{u_{xy}} + 4xz{u_{xz}} + 2yz{u_{yz}} - {y^2}{u_{yy}} - {z^2}{u_{zz}} \hfill \\
\,\,\,\,\,\,\,\,\,    + \left[ {c - 2\left( {2b + 3} \right)x} \right]{u_x} - 2\left( {b + 1} \right)y{u_y} + 2bz{u_z} - b\left( {b + 1} \right)u = 0, \hfill \\
\end{gathered}  \\
  \begin{gathered}
  y\left( {1 - y} \right){u_{yy}} + x\left( {1 - 2y} \right){u_{xy}} + yz{u_{yz}}   + \left[ {c - \left( {{a_1} + b + 1} \right)y} \right]{u_y} - 2{a_1}x{u_x} + {a_1}z{u_z} - {a_1}bu = 0, \hfill \\
\end{gathered}  \\
  {z\left( {1 + z} \right){u_{zz}} - 2x{u_{xz}} - y{u_{yz}} + \left[ {1 - b + \left( {{a_2} + {a_3} + 1} \right)z} \right]{u_z} + {a_2}{a_3}u = 0,}
\end{array}} \right.
$

where $u\equiv \,\,   {F_{30b}}\left( {{a_1},{a_2},{a_3},b;c;x,y,z} \right)$.

\bigskip

\begin{equation}
{F_{30c}}\left( {{a_1},{a_2},{a_3},b;c;x,y,z} \right) \hfill \\
  = \sum\limits_{m,n,p = 0}^\infty  {} \frac{{{{\left( {{a_1}} \right)}_{2m + n}}{{\left( {{a_2}} \right)}_p}{{\left( {{a_3}} \right)}_p}{{\left( b \right)}_{n - p}}}}{{{{\left( c \right)}_{m + n}}}}\frac{x^m}{m!}\frac{y^n}{n!}\frac{z^p}{p!},
\end{equation}

region of convergence:
$$
\left\{ {r < \frac{1}{4},\,\,\,t < 1,\,\,\,s < \left( {\frac{1}{2} + \frac{1}{2}\sqrt {1 - 4r} } \right)\min \left\{ {1,\frac{{1 - t}}{t}} \right\}} \right\}.
$$

System of partial differential equations:

$
\left\{ {\begin{array}{*{20}{l}}
  \begin{gathered}
  x\left( {1 - 4x} \right){u_{xx}} + \left( {1 - 4x} \right)y{u_{xy}} - {y^2}{u_{yy}}\hfill\\\,\,\,\,\,\,\,\,\,
   + \left[ {c - 2\left( {2{a_1} + 3} \right)x} \right]{u_x} - 2\left( {{a_1} + 1} \right)y{u_y} - {a_1}\left( {{a_1} + 1} \right)u = 0, \hfill \\
\end{gathered}  \\
  \begin{gathered}
  y\left( {1 - y} \right){u_{yy}} + x\left( {1 - 2y} \right){u_{xy}} + 2xz{u_{xz}} + yz{u_{yz}}\hfill\\\,\,\,\,\,\,\,\,\,
   + \left[ {c - \left( {{a_1} + b + 1} \right)y} \right]{u_y} - 2bx{u_x} + {a_1}z{u_z} - {a_1}bu = 0, \hfill \\
\end{gathered}  \\
  {z\left( {1 + z} \right){u_{zz}} - y{u_{yz}} + \left[ {1 - b + \left( {{a_2} + {a_3} + 1} \right)z} \right]{u_z} + {a_2}{a_3}u = 0,}
\end{array}} \right.
$

where $u\equiv \,\,   {F_{30c}}\left( {{a_1},{a_2},{a_3},b;c;x,y,z} \right)$.

\bigskip

\begin{equation}
{F_{30d}}\left( {{a_1},{a_2},{a_3},b;c;x,y,z} \right) \hfill \\
    = \sum\limits_{m,n,p = 0}^\infty  {} \frac{{{{\left( {{a_1}} \right)}_{2m + n}}{{\left( {{a_2}} \right)}_n}{{\left( {{a_3}} \right)}_p}{{\left( b \right)}_{p - m - n}}}}{{{{\left( c \right)}_p}}}\frac{x^m}{m!}\frac{y^n}{n!}\frac{z^p}{p!},
\end{equation}

region of convergence:
$$
 \left\{ {t < 1,\,\,\,r < \frac{1}{{4\left( {1 + t} \right)}},\,\,\,s < \frac{{1 + \sqrt {1 - 4r\left( {1 + t} \right)} }}{{2\left( {1 + t} \right)}}} \right\}.
$$

System of partial differential equations:

$
\left\{ {\begin{array}{*{20}{l}}
  \begin{gathered}
  x\left( {1 + 4x} \right){u_{xx}} + \left( {1 + 4x} \right)y{u_{xy}} - z{u_{xz}} + {y^2}{u_{yy}} \hfill \\ \,\,\,\,\,\,\,\,\,
   + \left[ {1 - b + 2\left( {2{a_1} + 3} \right)x} \right]{u_x} + 2\left( {{a_1} + 1} \right)y{u_y} + {a_1}\left( {{a_1} + 1} \right)u = 0, \hfill \\
\end{gathered}  \\
  \begin{gathered}
  y\left( {1 + y} \right){u_{yy}} + x\left( {1 + 2y} \right){u_{xy}} - z{u_{yz}}
   + \left[ {1 - b + \left( {{a_1} + {a_2} + 1} \right)} \right]{u_y} + 2{a_2}x{u_x} + {a_1}{a_2}u = 0, \hfill \\
\end{gathered}  \\
  z\left( {1 - z} \right){u_{zz}} + xz{u_{xz}} + yz{u_{yz}} + \left[ {c - \left( {{a_3} + b + 1} \right)z} \right]{u_z}
   + {a_3}x{u_x} + {a_3}y{u_y} - {a_3}bu = 0,
\end{array}} \right.
$

where $u\equiv \,\,   {F_{30d}}\left( {{a_1},{a_2},{a_3},b;c;x,y,z} \right)$.

Particular solutions:

$
{u_1} = {F_{30d}}\left( {{a_1},{a_2},{a_3},b;c;x,y,z} \right),
$

$
{u_2} = {z^{1 - c}}{F_{30d}}\left( {{a_1},{a_2},1 - c + {a_3},1 - c + b;2 - c;x,y,z} \right).
$

\bigskip

\begin{equation}
{F_{30e}}\left( {{a_1},{a_2},{b_1},{b_2};x,y,z} \right) \hfill \\
  = \sum\limits_{m,n,p = 0}^\infty  {} {{{{\left( {{a_1}} \right)}_n}{{\left( {{a_2}} \right)}_p}{{\left( {{b_1}} \right)}_{p - m - n}}{{\left( {{b_2}} \right)}_{2m + n - p}}}}\frac{x^m}{m!}\frac{y^n}{n!}\frac{z^p}{p!},
\end{equation}

region of convergence:
$$
 \left\{ {r < \frac{1}{4},\,\,\,s < \frac{1}{2} + \frac{1}{2}\sqrt {1 - 4r} ,\,\,\,t < \frac{{\sqrt {1 + 4r}  - 1}}{{2r}}} \right\}.
$$

System of partial differential equations:

$
\left\{ {\begin{array}{*{20}{l}}
  \begin{gathered}
  x\left( {1 + 4x} \right){u_{xx}} + \left( {1 + 4x} \right)y{u_{xy}} - \left( {1 + 4x} \right)z{u_{xz}} - 2yz{u_{yz}}+ {y^2}{u_{yy}} \hfill \\
  \,\,\,\,\,\,\,\,\,  + {z^2}{u_{zz}}  + \left[ {1 - {b_1} + 2\left( {2{b_2} + 3} \right)x} \right]{u_x} + 2\left( {{b_2} + 1} \right)y{u_y}
    - 2{b_2}z{u_z} + {b_2}\left( {{b_2} + 1} \right)u = 0, \hfill \\
\end{gathered}  \\
  \begin{gathered}
  y\left( {1 + y} \right){u_{yy}} + x\left( {1 + 2y} \right){u_{xy}} - \left( {1 + y} \right)z{u_{yz}}
 \hfill \\ \,\,\,\,\,\,\,\,\,   + \left[ {1 - {b_1} + \left( {{a_1} + {b_2} + 1} \right)y} \right]{u_y} + 2{a_1}x{u_x} - {a_1}z{u_z} + {a_1}{b_2}u = 0, \hfill \\
\end{gathered}  \\
  \begin{gathered}
  z\left( {1 + z} \right){u_{zz}} - x\left( {2 + z} \right){u_{xz}} - y\left( {1 + z} \right){u_{yz}} \hfill\\\,\,\,\,\,\,\,\,\,
   + \left[ {1 - {b_2} + \left( {{a_2} + {b_1} + 1} \right)z} \right]{u_z} - {a_2}x{u_x} - {a_2}y{u_y} + {a_2}{b_1}u = 0, \hfill \\
\end{gathered}
\end{array}} \right.
$

where $u\equiv \,\,   {F_{30e}}\left( {{a_1},{a_2},{b_1},{b_2};x,y,z} \right) $.

\bigskip

\begin{equation}
{F_{30f}}\left( {{a_1},{a_2},{b_1},{b_2};x,y,z} \right) \hfill \\
  = \sum\limits_{m,n,p = 0}^\infty  {} {{{{\left( {{a_1}} \right)}_{2m + n}}{{\left( {{a_2}} \right)}_p}{{\left( {{b_1}} \right)}_{p - m - n}}{{\left( {{b_2}} \right)}_{n - p}}}}\frac{x^m}{m!}\frac{y^n}{n!}\frac{z^p}{p!},
\end{equation}

region of convergence:
$$
\left\{ {t < 1,\,\,\,r < \frac{1}{{4\left( {1 + t} \right)}},\,\,\,s < \min \left\{ {\frac{1}{2} + \frac{1}{2}\sqrt {1 - 4r} ,\frac{{1 - 2\sqrt {r + rt} }}{t}} \right\}} \right\}.
$$

System of partial differential equations:

$
\left\{ {\begin{array}{*{20}{l}}
  \begin{gathered}
  x\left( {1 + 4x} \right){u_{xx}} + \left( {1 + 4x} \right)y{u_{xy}} - z{u_{xz}} + {y^2}{u_{yy}} \hfill \\
  \,\,\,\,\,\,\,\,\,  + \left[ {1 - {b_1} + 2\left( {2{a_1} + 3} \right)x} \right]{u_x} + 2\left( {{a_1} + 1} \right)y{u_y} + {a_1}\left( {{a_1} + 1} \right)u = 0, \hfill \\
\end{gathered}  \\
  \begin{gathered}
  y\left( {1 + y} \right){u_{yy}} + x\left( {1 + 2y} \right){u_{xy}} - 2xz{u_{xz}} - z\left( {1 + y} \right){u_{yz}} \hfill \\
\,\,\,\,\,\,\,\,\,    + \left[ {1 - {b_1} + \left( {{a_1} + {b_2} + 1} \right)y} \right]{u_y} + 2{b_2}x{u_x} - {a_1}z{u_z} + {a_1}{b_2}u = 0 \hfill \\
\end{gathered}  \\
  \begin{gathered}
  z\left( {1 + z} \right){u_{zz}} - xz{u_{xz}} - y\left( {1 + z} \right){u_{yz}}\hfill\\ \,\,\,\,\,\,\,\,\,
  + \left[ {1 - {b_2} + \left( {{a_2} + {b_1} + 1} \right)z} \right]{u_z} - {a_2}x{u_x} - {a_2}y{u_y} + {a_2}{b_1}u = 0, \hfill \\
\end{gathered}
\end{array}} \right.
$

where $u\equiv \,\,   {F_{30f}}\left( {{a_1},{a_2},{b_1},{b_2};x,y,z} \right)$.

\bigskip

\begin{equation}
{F_{31a}}\left( {{a_1},{a_2},{a_3},{a_4};{c_1},{c_2};x,y,z} \right) \hfill \\
  = \sum\limits_{m,n,p = 0}^\infty  {} \frac{{{{\left( {{a_1}} \right)}_{2m + n}}{{\left( {{a_2}} \right)}_n}{{\left( {{a_3}} \right)}_p}{{\left( {{a_4}} \right)}_p}}}{{{{\left( {{c_1}} \right)}_{m + p}}{{\left( {{c_2}} \right)}_n}}}{x^m}\frac{x^m}{m!}\frac{y^n}{n!}\frac{z^p}{p!},
\end{equation}

first appearance of this function in the literature, and old notation: [12], $X_{20}$,

region of convergence:
$$
\left\{ {s + 2\sqrt r  < 1,\,\,\,t < 1} \right\}.
$$

System of partial differential equations:

$
\left\{ {\begin{array}{*{20}{l}}
  \begin{gathered}
  x\left( {1 - 4x} \right){u_{xx}} - 4xy{u_{xy}} + z{u_{xz}} - {y^2}{u_{yy}}\hfill\\\,\,\,\,\,\,\,\,\,
   + \left[ {{c_1} - 2\left( {2{a_1} + 3} \right)x} \right]{u_x} - 2\left( {{a_1} + 1} \right)y{u_y} - {a_1}\left( {{a_1} + 1} \right)u = 0, \hfill \\
\end{gathered}  \\
  \begin{gathered}
  y\left( {1 - y} \right){u_{yy}} - 2xy{u_{xy}}
   + \left[ {{c_2} - \left( {{a_1} + {a_2} + 1} \right)y} \right]{u_y} - 2{a_2}x{u_x} - {a_1}{a_2}u = 0, \hfill \\
\end{gathered}  \\
  {z\left( {1 - z} \right){u_{zz}} + x{u_{xz}} + \left[ {{c_1} - \left( {{a_3} + {a_4} + 1} \right)z} \right]{u_z} - {a_3}{a_4}u = 0,}
\end{array}} \right.
$

where $u\equiv \,\,   {F_{31a}}\left( {{a_1},{a_2},{a_3},{a_4};{c_1},{c_2};x,y,z} \right)$.

Particular solutions:

$
{u_1} = {F_{31a}}\left( {{a_1},{a_2},{a_3},{a_4};{c_1},{c_2};x,y,z} \right),
$

$
{u_2} = {y^{1 - {c_2}}}{F_{31a}}\left( {1 - {c_2} + {a_1},1 - {c_2} + {a_2},{a_3},{a_4};{c_1},2 - {c_2};x,y,z} \right).
$

\bigskip

\begin{equation}
{F_{31b}}\left( {{a_1},{a_2},{a_3},b;c;x,y,z} \right) \hfill \\
  = \sum\limits_{m,n,p = 0}^\infty  {} \frac{{{{\left( {{a_1}} \right)}_{2m + n}}{{\left( {{a_2}} \right)}_n}{{\left( {{a_3}} \right)}_p}{{\left( b \right)}_{p - n}}}}{{{{\left( c \right)}_{m + p}}}}\frac{x^m}{m!}\frac{y^n}{n!}\frac{z^p}{p!},
\end{equation}

region of convergence:
$$
  \left\{ {s + 2\sqrt r  < 1,\,\,\,t < \min \left\{ {1,\frac{{1 - s + \sqrt {{{\left( {1 - s} \right)}^2} - 4r} }}{{2s}}} \right\}} \right\}.
$$

System of partial differential equations:

$
\left\{ {\begin{array}{*{20}{l}}
  \begin{gathered}
  x\left( {1 - 4x} \right){u_{xx}} - 4xy{u_{xy}} + z{u_{xz}} - {y^2}{u_{yy}}\hfill\\\,\,\,\,\,\,\,\,\,
   + \left[ {c - 2\left( {2{a_1} + 3} \right)x} \right]{u_x} - 2\left( {{a_1} + 1} \right)y{u_y} - {a_1}\left( {{a_1} + 1} \right)u = 0, \hfill \\
\end{gathered}  \\
  \begin{gathered}
  y\left( {1 + y} \right){u_{yy}} + 2xy{u_{xy}} - z{u_{yz}}
   - \left[ {b - 1 - \left( {{a_1} + {a_2} + 1} \right)y} \right]{u_y} + 2{a_2}x{u_x} + {a_1}{a_2}u = 0, \hfill \\
\end{gathered}  \\
  z\left( {1 - z} \right){u_{zz}} + x{u_{xz}} + yz{u_{yz}} + \left[ {c - \left( {{a_3} + b + 1} \right)z} \right]{u_z}
   + {a_3}y{u_y} - {a_3}bu = 0,
\end{array}} \right.
$

where $u\equiv \,\,   {F_{31b}}\left( {{a_1},{a_2},{a_3},b;c;x,y,z} \right)$.

\bigskip

\begin{equation}
{F_{31c}}\left( {{a_1},{a_2},{a_3},b;c;x,y,z} \right) \hfill \\
  = \sum\limits_{m,n,p = 0}^\infty  {} \frac{{{{\left( {{a_1}} \right)}_{2m + n}}{{\left( {{a_2}} \right)}_p}{{\left( {{a_3}} \right)}_p}{{\left( b \right)}_{n - m - p}}}}{{{{\left( c \right)}_n}}}\frac{x^m}{m!}\frac{y^n}{n!}\frac{z^p}{p!},
\end{equation}

region of convergence:
$$
\left\{ {s < 1,\,\,\,t < \frac{1}{{1 + s}},\,\,\,r < \min \left\{ {{\Theta _1}\left( s \right),{\Theta _2}\left( s \right),\frac{1}{{1 - t}}{\Theta _2}\left( {\frac{{st}}{{1 - t}}} \right)} \right\}} \right\}.
$$

System of partial differential equations:

$
\left\{ {\begin{array}{*{20}{l}}
  \begin{gathered}
  x\left( {1 + 4x} \right){u_{xx}} - \left( {1 - 4x} \right)y{u_{xy}} + z{u_{xz}} + {y^2}{u_{yy}} \hfill \\\,\,\,\,\,\,\,\,\,
   + \left[ {1 - b + 2\left( {2{a_1} + 3} \right)x} \right]{u_x} + 2\left( {{a_1} + 1} \right)y{u_y}
   + {a_1}\left( {{a_1} + 1} \right)u = 0, \hfill \\
\end{gathered}  \\
  \begin{gathered}
  y\left( {1 - y} \right){u_{yy}} - xy{u_{xy}} + 2xz{u_{xz}} + yz{u_{yz}} + 2{x^2}{u_{xx}} \hfill \\\,\,\,\,\,\,\,\,\,
   + \left[ {c - \left( {{a_1} + b + 1} \right)y} \right]{u_y} + \left( {{a_1} - 2b + 2} \right)x{u_x}
   + {a_1}z{u_z} - {a_1}bu = 0, \hfill \\
\end{gathered}  \\
  z\left( {1 + z} \right){u_{zz}} - y{u_{yz}} + x{u_{xz}} + \left[ {1 - b + \left( {{a_2} + {a_3} + 1} \right)z} \right]{u_z}
   + {a_2}{a_3}u = 0,
\end{array}} \right.
$

where $u\equiv \,\,   {F_{31c}}\left( {{a_1},{a_2},{a_3},b;c;x,y,z} \right)$.

Particular solutions:

$
{u_1} = {F_{31c}}\left( {{a_1},{a_2},{a_3},b;c;x,y,z} \right),
$

$
{u_2} = {y^{1 - c}}{F_{31c}}\left( {1 - c + {a_1},{a_2},{a_3},1 - c + b;2 - c;x,y,z} \right).
$

\bigskip

\begin{equation}
{F_{31d}}\left( {{a_1},{a_2},{b_1},{b_2};x,y,z} \right) \hfill \\
  = \sum\limits_{m,n,p = 0}^\infty  {} {{{{\left( {{a_1}} \right)}_{2m + n}}{{\left( {{a_2}} \right)}_p}{{\left( {{b_1}} \right)}_{n - m - p}}{{\left( {{b_2}} \right)}_{p - n}}}}\frac{x^m}{m!}\frac{y^n}{n!}\frac{z^p}{p!},
\end{equation}

region of convergence:
$$
 \left\{ {s < 1,\,\,\,t \leq \frac{3}{4}\left[ {\alpha(s) - 1} \right],\,\,t<1, \,\,\,r < \min \left\{ {{\Theta _1}\left( s \right),{\Theta _2}\left( s \right)} \right\}} \right\},
$$
$$
 \cup\left\{ {\frac{1}{5}<s < 1,\,\,\frac{3}{4}\left[ {\alpha(s) - 1} \right]<t<1, \,\,\,r < \min \left\{ {{\Theta _2}\left( s \right)}, \,\,\frac{t(1-st)^2}{4(1+t)} \right\}} \right\}.
$$

System of partial differential equations:

$
\left\{ {\begin{array}{*{20}{l}}
  \begin{gathered}
  x\left( {1 + 4x} \right){u_{xx}} - \left( {1 - 4x} \right)y{u_{xy}} + z{u_{xz}} + {y^2}{u_{yy}} \hfill \\
 \,\,\,\,\,\,\,\,\,   + \left[ {1 - {b_1} + 2\left( {2{a_1} + 3} \right)x} \right]{u_x} + 2\left( {{a_1} + 1} \right)y{u_y}
    + {a_1}\left( {{a_1} + 1} \right)u = 0, \hfill \\
\end{gathered}  \\
  \begin{gathered}
  y\left( {1 + y} \right){u_{yy}} + xy{u_{xy}} - 2xz{u_{xz}} - \left( {1 + y} \right)z{u_{yz}} - 2{x^2}{u_{xx}}   \hfill \\
 \,\,\,\,\,\,\,\,\,   + \left[ {1 - {b_2} + \left( {{a_1} + {b_1} + 1} \right)y} \right]{u_y} - \left( {{a_1} - 2{b_1}}+2 \right)x{u_x}
    - {a_1}z{u_z} + {a_1}{b_1}u = 0, \hfill \\
\end{gathered}  \\
  z\left( {1 + z} \right){u_{zz}} + x{u_{xz}} - y\left( {1 + z} \right){u_{yz}}
  + \left[ {1 - {b_1} + \left( {{a_2} + {b_2} + 1} \right)z} \right]{u_z} - {a_2}y{u_y} + {a_2}{b_2}u = 0,\hfill \\
\end{array}} \right.
$

where $u\equiv \,\,   {F_{31d}}\left( {{a_1},{a_2},{b_1},{b_2};x,y,z} \right)$.

\bigskip

\begin{equation}
{F_{32a}}\left( {{a_1},{a_2},{a_3},{a_4};{c_1},{c_2};x,y,z} \right) \hfill \\
  = \sum\limits_{m,n,p = 0}^\infty  {} \frac{{{{\left( {{a_1}} \right)}_{2m + n}}{{\left( {{a_2}} \right)}_n}{{\left( {{a_3}} \right)}_p}{{\left( {{a_4}} \right)}_p}}}{{{{\left( {{c_1}} \right)}_{n + p}}{{\left( {{c_2}} \right)}_m}}}\frac{x^m}{m!}\frac{y^n}{n!}\frac{z^p}{p!},
\end{equation}

first appearance of this function in the literature, and old notation:
[12], $X_{19}$,

region of convergence:
$$
\left\{ {s + 2\sqrt r  < 1,\,\,\,t < 1} \right\}.
$$

System of partial differential equations:

$
\left\{ {\begin{array}{*{20}{l}}
  \begin{gathered}
  x\left( {1 - 4x} \right){u_{xx}} - 4xy{u_{xy}} - {y^2}{u_{yy}}
   + \left[ {{c_2} - 2\left( {2{a_1} + 3} \right)x} \right]{u_x} - 2\left( {{a_1} + 1} \right)y{u_y} - {a_1}\left( {{a_1} + 1} \right)u = 0, \hfill \\
\end{gathered}  \\
  \begin{gathered}
  y\left( {1 - y} \right){u_{yy}} - 2xy{u_{xy}} + z{u_{yz}}
   + \left[ {{c_1} - \left( {{a_1} + {a_2} + 1} \right)y} \right]{u_y} - 2{a_2}x{u_x} - {a_1}{a_2}u = 0, \hfill \\
\end{gathered}  \\
  {z\left( {1 - z} \right){u_{zz}} + y{u_{yz}} + \left[ {{c_1} - \left( {{a_3} + {a_4} + 1} \right)z} \right]{u_z} - {a_3}{a_4}u = 0,}
\end{array}} \right.
$

where $u\equiv \,\,   {F_{32a}}\left( {{a_1},{a_2},{a_3},{a_4};{c_1},{c_2};x,y,z} \right) $.

Particular solutions:

$
{u_1} = {F_{32a}}\left( {{a_1},{a_2},{a_3},{a_4};{c_1},{c_2};x,y,z} \right),
$

$
{u_2} = {x^{1 - {c_2}}}{F_{32a}}\left( {2 - 2{c_2} + {a_1},{a_2},{a_3},{a_4};{c_1},2 - {c_2};x,y,z} \right).
$

\bigskip

\begin{equation}
{F_{32b}}\left( {{a_1},{a_2},{a_3},b;c;x,y,z} \right)\hfill \\
   = \sum\limits_{m,n,p = 0}^\infty  {} \frac{{{{\left( {{a_1}} \right)}_{2m + n}}{{\left( {{a_2}} \right)}_p}{{\left( {{a_3}} \right)}_p}{{\left( b \right)}_{n - m}}}}{{{{\left( c \right)}_{n + p}}}}\frac{x^m}{m!}\frac{y^n}{n!}\frac{z^p}{p!},
\end{equation}

region of convergence:
$$
 \left\{ {r < \frac{1}{4},\,\,\,t < 1,\,\,\,s < \min \left\{ {{\Psi _1}\left( r \right),{\Psi _2}\left( r \right)} \right\}} \right\}.
$$

System of partial differential equations:

$
\left\{ {\begin{array}{*{20}{l}}
  \begin{gathered}
  x\left( {1 + 4x} \right){u_{xx}} - \left( {1 - 4x} \right)y{u_{xy}} + {y^2}{u_{yy}} \hfill\\\,\,\,\,\,\,\,\,\,
   + \left[ {1 - b + 2\left( {2{a_1} + 3} \right)x} \right]{u_x} + 2\left( {{a_1} + 1} \right)y{u_y}
   + {a_1}\left( {{a_1} + 1} \right)u = 0, \hfill \\
\end{gathered}  \\
  \begin{gathered}
  y\left( {1 - y} \right){u_{yy}} - xy{u_{xy}} + z{u_{yz}} + 2{x^2}{u_{xx}}\hfill\\ \,\,\,\,\,\,\,\,\,
   + \left[ {c - \left( {{a_1} + b + 1} \right)y} \right]{u_y} + {a_1}x{u_x} - \left( a_1- 2b+2 \right)x{u_x} - {a_1}bu = 0, \hfill \\
\end{gathered}  \\
  {z\left( {1 - z} \right){u_{zz}} + y{u_{yz}} + \left[ {c - \left( {{a_2} + {a_3} + 1} \right)z} \right]{u_z} - {a_2}{a_3}u = 0,}
\end{array}} \right.
$

where $u\equiv \,\,   {F_{32b}}\left( {{a_1},{a_2},{a_3},b;c;x,y,z} \right)$.

\bigskip

\begin{equation}
{F_{32c}}\left( {{a_1},{a_2},{a_3},b;c;x,y,z} \right) \hfill \\
  = \sum\limits_{m,n,p = 0}^\infty  {} \frac{{{{\left( {{a_1}} \right)}_{2m + n}}{{\left( {{a_2}} \right)}_n}{{\left( {{a_3}} \right)}_p}{{\left( b \right)}_{p - m}}}}{{{{\left( c \right)}_{n + p}}}}\frac{x^m}{m!}\frac{y^n}{n!}\frac{z^p}{p!},
\end{equation}

region of convergence:
$$
\left\{r < \frac{1}{4},\,\,\,s \leq \frac{1}{2} - 2r,\,\,\,t < \min \left\{ {1,\frac{{1 - 4r}}{{4r}}} \right\} \right\}
$$
$$
\cup \left\{ r < \frac{1}{4},\,\,\,\frac{1}{2}-2r<s<1-2\sqrt{r},\,\,\,t < \min \left\{ 1,\frac{s}{2r}\left(1-s+\sqrt{(1-s)^2-4r}\right) \right\}\right\}.
$$

System of partial differential equations:

$
\left\{ {\begin{array}{*{20}{l}}
  \begin{gathered}
  x\left( {1 + 4x} \right){u_{xx}} + 4xy{u_{xy}} - z{u_{xz}} + {y^2}{u_{yy}}\hfill\\\,\,\,\,\,\,\,\,\,
   + \left[ {1 - b + 2\left( {2{a_1} + 3} \right)x} \right]{u_x} + 2\left( {{a_1} + 1} \right)y{u_y}
   + {a_1}\left( {1 + {a_1}} \right)u = 0, \hfill \\
\end{gathered}  \\
  \begin{gathered}
  y\left( {1 - y} \right){u_{yy}} - 2xy{u_{xy}} + z{u_{yz}}
   + \left[ {c - \left( {{a_1} + {a_2} + 1} \right)y} \right]{u_y} - 2{a_2}x{u_x} - {a_1}{a_2}u = 0, \hfill \\
\end{gathered}  \\
  z\left( {1 - z} \right){u_{zz}} + xz{u_{xz}} + y{u_{yz}} + \left[ {c - \left( {{a_3} + b + 1} \right)z} \right]{u_z}
   + {a_3}x{u_x} - {a_3}bu = 0,
\end{array}} \right.
$

where $u\equiv \,\,   {F_{32c}}\left( {{a_1},{a_2},{a_3},b;c;x,y,z} \right)$.

\bigskip

\begin{equation}
{F_{33a}}\left( {{a_1},{a_2},{a_3},{a_4};c;x,y,z} \right) \hfill \\
  = \sum\limits_{m,n,p = 0}^\infty  {} \frac{{{{\left( {{a_1}} \right)}_{2m + n}}{{\left( {{a_2}} \right)}_n}{{\left( {{a_3}} \right)}_p}{{\left( {{a_4}} \right)}_p}}}{{{{\left( c \right)}_{m + n + p}}}}\frac{x^m}{m!}\frac{y^n}{n!}\frac{z^p}{p!},
\end{equation}

first appearance of this function in the literature, and old notation:
[12], $X_{18}$,

region of convergence:
$$
 \left\{ {r < \frac{1}{4},\,\,\,t < 1,\,\,\,s < \frac{1}{2} + \frac{1}{2}\sqrt {1 - 4r} } \right\}.
$$

System of partial differential equations:

$
\left\{ {\begin{array}{*{20}{l}}
  \begin{gathered}
  x\left( {1 - 4x} \right){u_{xx}} + \left( {1 - 4x} \right)y{u_{xy}} + z{u_{xz}} - {y^2}{u_{yy}} \hfill \\\,\,\,\,\,\,\,\,\,
   + \left[ {c - 2\left( {2{a_1} + 3} \right)x} \right]{u_x} - 2\left( {{a_1} + 1} \right)y{u_y} - {a_1}\left( {1 + {a_1}} \right)u = 0, \hfill \\
\end{gathered}  \\
  \begin{gathered}
  y\left( {1 - y} \right){u_{yy}} + x\left( {1 - 2y} \right){u_{xy}} + z{u_{yz}}
   + \left[ {c - \left( {{a_1} + {a_2} + 1} \right)y} \right]{u_y} - 2{a_2}x{u_x} - {a_1}{a_2}u = 0, \hfill \\
\end{gathered}  \\
  z\left( {1 - z} \right){u_{zz}} + x{u_{xz}} + y{u_{yz}} + \left[ {c - \left( {{a_3} + {a_4} + 1} \right)z} \right]{u_z}
  - {a_3}{a_4}u = 0,
\end{array}} \right.
$

where $u\equiv \,\,   {F_{33a}}\left( {{a_1},{a_2},{a_3},{a_4};c;x,y,z} \right) $.

\bigskip

\begin{equation}
{F_{34b}}\left( {{a_1},{a_2},b;{c_1},{c_2};x,y,z} \right) \hfill \\
  = \sum\limits_{m,n,p = 0}^\infty  {} \frac{{{{\left( {{a_1}} \right)}_{m + n}}{{\left( {{a_2}} \right)}_{m + n}}{{\left( b \right)}_{2p - m}}}}{{{{\left( {{c_1}} \right)}_n}{{\left( {{c_2}} \right)}_p}}}\frac{x^m}{m!}\frac{y^n}{n!}\frac{z^p}{p!},
\end{equation}

region of convergence:
$$
\left\{ {r < \frac{1}{4},\,\,\,\sqrt s  + \sqrt {r + 2r\sqrt t }  < 1} \right\}.
$$

System of partial differential equations:

$
\left\{ {\begin{array}{*{20}{l}}
  \begin{gathered}
  x\left( {1 + x} \right){u_{xx}} + 2xy{u_{xy}} - 2z{u_{xz}} + {y^2}{u_{yy}}\hfill \\\,\,\,\,\,\,\,\,\,
   + \left[ {1 - b + \left( {{a_1} + {a_2} + 1} \right)x} \right]{u_x} + \left( {{a_1} + {a_2} + 1} \right)y{u_y} + {a_1}{a_2}u = 0, \hfill \\
\end{gathered}  \\
  \begin{gathered}
  y\left( {1 - y} \right){u_{yy}} - 2xy{u_{xy}} - {x^2}{u_{xx}}
   + \left[ {{c_1} - \left( {{a_1} + {a_2} + 1} \right)y} \right]{u_y} - \left( {{a_1} + {a_2} + 1} \right)x{u_x} - {a_1}{a_2}u = 0, \hfill \\
\end{gathered}  \\
  z\left( {1 - 4z} \right){u_{zz}} + 4xz{u_{xz}} - {x^2}{u_{xx}} + \left[ {{c_2} - 2\left( {2b + 3} \right)z} \right]{u_z}
  + 2bx{u_x} - b\left( {1 + b} \right)u = 0,
\end{array}} \right.
$

where $u\equiv \,\,   {F_{34b}}\left( {{a_1},{a_2},b;{c_1},{c_2};x,y,z} \right)$.

Particular solutions:

$
{u_1} = {F_{34b}}\left( {{a_1},{a_2},b;{c_1},{c_2};x,y,z} \right),
$

$
{u_2} = {y^{1 - {c_1}}}{F_{34b}}\left( {1 - {c_1} + {a_1},1 - {c_1} + {a_2},b;2 - {c_1},{c_2};x,y,z} \right),
$

$
{u_3} = {z^{1 - {c_2}}}{F_{34b}}\left( {{a_1},{a_2},2 - 2{c_2} + b;{c_1},2 - {c_2};x,y,z} \right),
$

$
{u_4} = {y^{1 - {c_1}}}{z^{1 - {c_2}}}{F_{34b}}\left( {1 - {c_1} + {a_1},1 - {c_1} + {a_2},2 - 2{c_2} + b;2 - {c_1},2 - {c_2};x,y,z} \right).
$

\bigskip

\begin{equation}
{F_{34c}}\left( {a,{b_1},{b_2};c;x,y,z} \right)  = \sum\limits_{m,n,p = 0}^\infty  {} \frac{{{{\left( a \right)}_{m + n}}{{\left( {{b_1}} \right)}_{2p - m}}{{\left( {{b_2}} \right)}_{m + n - p}}}}{{{{\left( c \right)}_n}}}\frac{x^m}{m!}\frac{y^n}{n!}\frac{z^p}{p!},
\end{equation}

region of convergence:
$$
\begin{gathered}
    \left\{ {\sqrt r  + \sqrt s  < 1,\,\,\,t < \min \left\{ {{U^ + }\left( {{w_1}} \right),{U^ - }\left( {{w_2}} \right)} \right\}} \right\}, \hfill \\
  {P^ \pm }\left( w \right) = {w^3} \pm \left( {2 \mp s - r} \right)sw \mp 2{s^2}, \hfill \\
  {w_1}:\,\,{\rm{the\,\,root\,\,in}}\,\,\left( {s,\infty } \right)\,\,{\rm{of}}\,\,{P^ + }\left( w \right) = 0, \hfill \\
  {w_2}:\,\,{\rm{the\,\,root\,\,in}}\,\,\left( {\sqrt s ,\infty } \right)\,\,{\rm{of}}\,\,{P^ - }\left( w \right) = 0, \hfill \\
  {U^ \pm }\left( w \right) = \frac{{{{\left( {w - s} \right)}^2}\left( {{w^2} \pm s} \right)}}{{{r^2}s{w^2}}}. \hfill \\
\end{gathered}
$$

System of partial differential equations:

$
\left\{ {\begin{array}{*{20}{l}}
  \begin{gathered}
  x\left( {1 + x} \right){u_{xx}} + 2xy{u_{xy}} - \left( {2 + x} \right)z{u_{xz}} - yz{u_{yz}} + {y^2}{u_{yy}} \hfill \\
 \,\,\,\,\,\,\,\,\,   + \left[ {1 - {b_1} + \left( {a + {b_2} + 1} \right)x} \right]{u_x} + \left( {a + {b_2} + 1} \right)y{u_y}
    - az{u_z} + a{b_2}u = 0, \hfill \\
\end{gathered}  \\
  \begin{gathered}
  y\left( {1 - y} \right){u_{yy}} - 2xy{u_{xy}} + xz{u_{xz}} + yz{u_{yz}} - {x^2}{u_{xx}} \hfill \\
  \,\,\,\,\,\,\,\,\,  + \left[ {c - \left( {a + {b_2} + 1} \right)y} \right]{u_y} - \left( {a + {b_2} + 1} \right)x{u_x} + az{u_z} - a{b_2}u = 0, \hfill \\
\end{gathered}  \\
  \begin{gathered}
  z\left( {1 + 4z} \right){u_{zz}} - x\left( {1 + 4z} \right){u_{xz}} - y{u_{yz}} + {x^2}{u_{xx}}\hfill\\\,\,\,\,\,\,\,\,\,
  + \left[ {1 - {b_2} + 2\left( {2{b_1} + 3} \right)z} \right]{u_z} - 2{b_1}x{u_x} + {b_1}\left( {1 + {b_1}} \right)u = 0, \hfill \\
\end{gathered}
\end{array}} \right.
$

where $u\equiv \,\,   {F_{34c}}\left( {a,{b_1},{b_2};c;x,y,z} \right)$.

Particular solutions:

$
{u_1} = {F_{34c}}\left( {a,{b_1},{b_2};c;x,y,z} \right),
$

$
{u_2} = {y^{1 - c}}{F_{34c}}\left( {1 - c + a,{b_1},1 - c + {b_2};2 - c;x,y,z} \right).
$

\bigskip

\begin{equation}
{F_{35b}}\left( {{a_1},{a_2},b;c;x,y,z} \right)  = \sum\limits_{m,n,p = 0}^\infty  {} \frac{{{{\left( {{a_1}} \right)}_{m + n}}{{\left( {{a_2}} \right)}_{m + n}}{{\left( b \right)}_{2p - m}}}}{{{{\left( c \right)}_{n + p}}}}\frac{x^m}{m!}\frac{y^n}{n!}\frac{z^p}{p!},
\end{equation}

region of convergence:
$$
\begin{gathered}
  \left\{ {\sqrt r  + \sqrt s  < 1,\,\,\,t < \min \left\{ {1,\frac{{s\left( {1 - w_1^2} \right){{\left( {1 - {w_1}\sqrt s } \right)}^2}}}{{rw_1^2}}} \right\}} \right\}, \hfill \\
  P\left( w \right) = s{w^3} - \left( {1 + 2s - r} \right)w + 2\sqrt s , \hfill \\
{w_1}:\,\,{\rm{the\,\,root\,\,in}}\,\,\left( {\sqrt s ,1} \right)\,\,{\rm{of}}\,\,P\left( w \right) = 0. \hfill \\
\end{gathered}
$$

System of partial differential equations:

$
\left\{ {\begin{array}{*{20}{l}}
  \begin{gathered}
  x\left( {1 +x} \right){u_{xx}} + 2xy{u_{xy}} -2z{u_{xz}}  + {y^2}{u_{yy}} \hfill\\ \,\,\,\,\,\,\,\,\, + \left[1-b+\left(a_1+a_2+1\right)x\right]u_x+\left(a_1+a_2+1\right)yu_y+a_1a_2u = 0, \hfill \\
\end{gathered}  \\
  \begin{gathered}
  y\left( {1 - y} \right){u_{yy}}-x^2u_{xx}- 2xy{u_{xy}} + z{u_{yz}}\hfill\\ \,\,\,\,\,\,\,\,\,
  - (a_1+a_2+1)x{u_x}  + \left[ {c - \left( {{a_1} + a_2 + 1} \right)y} \right]{u_y}
    - {a_1}a_2u = 0, \hfill \\
\end{gathered}  \\
  \begin{gathered}
  z\left( {1 -4z} \right){u_{zz}} -x^2u_{xx}+ 4xz{u_{xz}} + y{u_{yz}}
  +2bxu_x + \left[ {c -\left( 4b+6 \right)z} \right]{u_z}  - b(b+1)u = 0, \hfill \\
\end{gathered}
\end{array}} \right.
$

where $u\equiv \,\,   {F_{35b}}\left( {{a_1},{a_2},b;c;x,y,z} \right)$.

\bigskip

\begin{equation}
{F_{36a}}\left( {{a_1},{a_2},{a_3};{c_1},{c_2},{c_3};x,y,z} \right) \hfill \\
  = \sum\limits_{m,n,p = 0}^\infty  {} \frac{{{{\left( {{a_1}} \right)}_{2m + n}}{{\left( {{a_2}} \right)}_{n + p}}{{\left( {{a_3}} \right)}_p}}}{{{{\left( {{c_1}} \right)}_m}{{\left( {{c_2}} \right)}_n}{{\left( {{c_3}} \right)}_p}}}\frac{x^m}{m!}\frac{y^n}{n!}\frac{z^p}{p!},
\end{equation}

first appearance of this function in the literature, and old notation:
[12], $X_{17}$,

region of convergence:
$$
  \left\{ {r < \frac{1}{4},\,\,\,t < 1,\,\,\,s < \left( {1 - 2\sqrt r } \right)\left( {1 - t} \right)} \right\}.
$$

System of partial differential equations:

$
\left\{ {\begin{array}{*{20}{l}}
  \begin{gathered}
  x\left( {1 - 4x} \right){u_{xx}} - 4xy{u_{xy}} - {y^2}{u_{yy}}\hfill\\ \,\,\,\,\,\,\,\,\,
   + \left[ {{c_1} - 2\left( {2{a_1} + 3} \right)x} \right]{u_x} - 2\left( {1 + {a_1}} \right)y{u_y} - {a_1}\left( {1 + {a_1}} \right)u = 0, \hfill \\
\end{gathered}  \\
  \begin{gathered}
  y\left( {1 - y} \right){u_{yy}} - 2xy{u_{xy}} - 2xz{u_{xz}} - yz{u_{yz}}\hfill\\\,\,\,\,\,\,\,\,\,
  + \left[ {{c_2} - \left( {{a_1} + {a_2} + 1} \right)y} \right]{u_y} - 2{a_2}x{u_x} - {a_1}z{u_z} - {a_1}{a_2}u = 0, \hfill \\
\end{gathered}  \\
  \begin{gathered}
  z\left( {1 - z} \right){u_{zz}} - yz{u_{yz}}
  + \left[ {{c_3} - \left( {{a_2} + {a_3} + 1} \right)z} \right]{u_z} - {a_3}y{u_y} - {a_2}{a_3}u = 0, \hfill \\
\end{gathered}
\end{array}} \right.
$

where $u\equiv \,\,   {F_{36a}}\left( {{a_1},{a_2},{a_3};{c_1},{c_2},{c_3};x,y,z} \right)$.

Particular solutions:

$
{u_1} = {F_{36a}}\left( {{a_1},{a_2},{a_3};{c_1},{c_2},{c_3};x,y,z} \right),
$

$
{u_2} = {x^{1 - {c_1}}}{F_{36a}}\left( {2 - 2{c_1} + {a_1},{a_2},{a_3};2 - {c_1},{c_2},{c_3};x,y,z} \right),
$

$
{u_3} = {y^{1 - {c_2}}}{F_{36a}}\left( {1 - {c_2} + {a_1},1 - {c_2} + {a_2},{a_3};{c_1},2 - {c_2},{c_3};x,y,z} \right),
$

$
{u_4} = {z^{1 - {c_3}}}{F_{36a}}\left( {{a_1},1 - {c_3} + {a_2},1 - {c_3} + {a_3};{c_1},{c_2},2 - {c_3};x,y,z} \right),
$

$
{u_5} = {x^{1 - {c_1}}}{y^{1 - {c_2}}}{F_{36a}}\left( {3 - 2{c_1} - {c_2} + {a_1},1 - {c_2} + {a_2},{a_3};2 - {c_1},2 - {c_2},{c_3};x,y,z} \right),
$

$
{u_6} = {y^{1 - {c_2}}}{z^{1 - {c_3}}}{F_{36a}}\left( {1 - {c_2} + {a_1},2 - {c_2} - {c_3} + {a_2},1 - {c_3} + {a_3};{c_1},2 - {c_2},2 - {c_3};x,y,z} \right),
$

$
{u_7} = {x^{1 - {c_1}}}{z^{1 - {c_3}}}{F_{36a}}\left( {2 - 2{c_1} + {a_1},1 - {c_3} + {a_2},1 - {c_3} + {a_3};2 - {c_1},{c_2},2 - {c_3};x,y,z} \right),
$

$
  {u_8} = {x^{1 - {c_1}}}{y^{1 - {c_2}}}{z^{1 - {c_3}}}\times$
   
   $\,\,\,\,\,\,\,\,\,\times{F_{36a}}\left( {3 - 2{c_1} - {c_2} + {a_1},2 - {c_2} - {c_3} + {a_2},1 - {c_3} + {a_3};2 - {c_1},2 - {c_2},2 - {c_3};x,y,z} \right).
  $

\bigskip

\begin{equation}
{F_{36b}}\left( {{a_1},{a_2},b;{c_1},{c_2};x,y,z} \right)  = \sum\limits_{m,n,p = 0}^\infty  {} \frac{{{{\left( {{a_1}} \right)}_{n + p}}{{\left( {{a_2}} \right)}_p}{{\left( b \right)}_{2m + n - p}}}}{{{{\left( {{c_1}} \right)}_m}{{\left( {{c_2}} \right)}_n}}}\frac{x^m}{m!}\frac{y^n}{n!}\frac{z^p}{p!},
\end{equation}

region of convergence:
$$
 \left\{ {s < 1,\,\,\,t+2\sqrt{st} <1, \,\,\,\sqrt{r} < \frac{1}{2}\min \left\{ 1-s, \,\,\, \frac{1-t-2\sqrt{st}}{t}\right\}} \right\}.
$$

System of partial differential equations:

$
\left\{ {\begin{array}{*{20}{l}}
  \begin{gathered}
  x\left( {1 - 4x} \right){u_{xx}} - 4xy{u_{xy}} + 4xz{u_{xz}} + 2yz{u_{yz}} - {y^2}{u_{yy}} - {z^2}{u_{zz}} \hfill \\
 \,\,\,\,\,\,\,\,\,  + \left[ {{c_1} - 2\left( {2b + 3} \right)x} \right]{u_x} - 2\left( {b + 1} \right)y{u_y} + 2bz{u_z} - b\left( {1 + b} \right)u = 0, \hfill \\
\end{gathered}  \\
  \begin{gathered}
  y\left( {1 - y} \right){u_{yy}} - 2xy{u_{xy}} - 2xz{u_{xz}} + {z^2}{u_{zz}} \hfill \\
 \,\,\,\,\,\,\,\,\,   + \left[ {{c_2} - \left( {{a_1} + b + 1} \right)y} \right]{u_y} - 2{a_1}x{u_x} + \left( {{a_1} - b + 1} \right)z{u_z} - {a_1}bu = 0, \hfill \\
\end{gathered}  \\
  \begin{gathered}
  z\left( {1 + z} \right){u_{zz}} - 2x{u_{xz}} - y\left( {1 - z} \right){u_{yz}}
   + \left[ {1 - b + \left( {{a_1} + {a_2} + 1} \right)z} \right]{u_z} + {a_2}y{u_y} + {a_1}{a_2}u = 0, \hfill \\
\end{gathered}
\end{array}} \right.
$

where $u\equiv \,\,   {F_{36b}}\left( {{a_1},{a_2},b;{c_1},{c_2};x,y,z} \right)$.

Particular solutions:

$
{u_1} = {F_{36b}}\left( {{a_1},{a_2},b;{c_1},{c_2};x,y,z} \right),
$

$
{u_2} = {x^{1 - {c_1}}}{F_{36b}}\left( {{a_1},{a_2},2 - 2{c_1} + b;2 - {c_1},{c_2};x,y,z} \right),
$

$
{u_3} = {y^{1 - {c_2}}}{F_{36b}}\left( {1 - {c_2} + {a_1},{a_2},1 - {c_2} + b;{c_1},2 - {c_2};x,y,z} \right),
$

$
{u_4} = {x^{1 - {c_1}}}{y^{1 - {c_2}}}{F_{36b}}\left( {1 - {c_2} + {a_1},{a_2},3 - 2{c_1} - {c_2} + b;2 - {c_1},2 - {c_2};x,y,z} \right).
$

\bigskip

\begin{equation}
{F_{36c}}\left( {{a_1},{a_2},b;{c_1},{c_2};x,y,z} \right)  = \sum\limits_{m,n,p = 0}^\infty  {} \frac{{{{\left( {{a_1}} \right)}_{2m + n}}{{\left( {{a_2}} \right)}_p}{{\left( b \right)}_{n + p - m}}}}{{{{\left( {{c_1}} \right)}_n}{{\left( {{c_2}} \right)}_p}}}\frac{x^m}{m!}\frac{y^n}{n!}\frac{z^p}{p!},
\end{equation}

region of convergence:
$$
\left\{ {s + t < 1,\,\,\,r < \min \left\{ {\frac{1}{{1 + t}}{\Theta _1}\left( {\frac{s}{{1 + t}}} \right),\frac{1}{{1 - t}}{\Theta _2}\left( {\frac{s}{{1 - t}}} \right)} \right\}} \right\}.
$$

System of partial differential equations:

$
\left\{ {\begin{array}{*{20}{l}}
  \begin{gathered}
  x\left( {1 + 4x} \right){u_{xx}} - \left( {1 - 4x} \right)y{u_{xy}} - z{u_{xz}} + {y^2}{u_{yy}} \hfill \\
\,\,\,\,\,\,\,\,\,    + \left[ {1 - b + 2\left( {2{a_1} + 3} \right)x} \right]{u_x} + 2\left( {{a_1} + 1} \right)y{u_y} + {a_1}\left( {1 + {a_1}} \right)u = 0, \hfill \\
\end{gathered}  \\
  \begin{gathered}
  y\left( {1 - y} \right){u_{yy}} - xy{u_{xy}} - 2xz{u_{xz}} - yz{u_{yz}} + 2{x^2}{u_{xx}} \hfill \\
  \,\,\,\,\,\,\,\,\,  + \left[ {{c_1} - \left( {{a_1} + b + 1} \right)y} \right]{u_y} + \left( {{a_1} - 2b + 2} \right)x{u_x} - {a_1}z{u_z} - {a_1}bu = 0, \hfill \\
\end{gathered}  \\
  \begin{gathered}
  z\left( {1 - z} \right){u_{zz}} + xz{u_{xz}} - yz{u_{yz}}
   + \left[ {{c_2} - \left( {{a_2} + b + 1} \right)z} \right]{u_z} - {a_2}y{u_y} + {a_2}x{u_x} - {a_2}bu = 0, \hfill \\
\end{gathered}
\end{array}} \right.
$

where $u\equiv \,\,   {F_{36c}}\left( {{a_1},{a_2},b;{c_1},{c_2};x,y,z} \right)$.

Particular solutions:

$
{u_1} = {F_{36c}}\left( {{a_1},{a_2},b;{c_1},{c_2};x,y,z} \right),
$

$
{u_2} = {y^{1 - {c_1}}}{F_{36c}}\left( {1 - {c_1} + {a_1},{a_2},1 - {c_1} + b;2 - {c_1},{c_2};x,y,z} \right),
$

$
{u_3} = {z^{1 - {c_2}}}{F_{36c}}\left( {{a_1},1 - {c_2} + {a_2},1 - {c_2} + b;{c_1},2 - {c_2};x,y,z} \right),
$

$
{u_4} = {y^{1 - {c_1}}}{z^{1 - {c_2}}}{F_{36c}}\left( {1 - {c_1} + {a_1},1 - {c_2} + {a_2},2 - {c_1} - {c_2} + b;2 - {c_1},2 - {c_2};x,y,z} \right).
$

\bigskip

\begin{equation}
{F_{36d}}\left( {{a_1},{a_2},b;{c_1},{c_2};x,y,z} \right)  = \sum\limits_{m,n,p = 0}^\infty  {} \frac{{{{\left( {{a_1}} \right)}_{2m + n}}{{\left( {{a_2}} \right)}_{n + p}}{{\left( b \right)}_{p - m}}}}{{{{\left( {{c_1}} \right)}_n}{{\left( {{c_2}} \right)}_p}}}\frac{x^m}{m!}\frac{y^n}{n!}\frac{z^p}{p!},
\end{equation}

region of convergence:
$$
\begin{gathered}
   \left\{ {s + t < 1,} \right. \hfill \\
  \left. {r < \min \left[ {\frac{{{{\left( {1 - s} \right)}^2}}}{{1 + t}}{\Theta _1}\left( {\frac{{st}}{{\left( {1 - s} \right)\left( {1 + t} \right)}}} \right),\frac{{{{\left( {1 - s} \right)}^2}}}{{1 - t}}{\Theta _2}\left( {\frac{{st}}{{\left( {1 - s} \right)\left( {1 - t} \right)}}} \right)} \right]} \right\}. \hfill \\
\end{gathered}
$$

System of partial differential equations:

$
\left\{ {\begin{array}{*{20}{l}}
  \begin{gathered}
  x\left( {1 + 4x} \right){u_{xx}} + 4xy{u_{xy}} - z{u_{xz}} + {y^2}{u_{yy}} \hfill \\
  \,\,\,\,\,\,\,\,\,  + \left[ {1 - b + 2\left( {2{a_1} + 3} \right)x} \right]{u_x} + 2\left( {{a_1} + 1} \right)y{u_y} + {a_1}\left( {1 + {a_1}} \right)u = 0, \hfill \\
\end{gathered}  \\
  \begin{gathered}
  y\left( {1 - y} \right){u_{yy}} - 2xy{u_{xy}} - 2xz{u_{xz}} - yz{u_{yz}}
   + \left[ {{c_1} - \left( {{a_1} + {a_2} + 1} \right)y} \right]{u_y} - 2{a_2}x{u_x} - {a_1}z{u_z} - {a_1}{a_2}u = 0, \hfill \\
\end{gathered}  \\
  \begin{gathered}
  z\left( {1 - z} \right){u_{zz}} + xy{u_{xy}} + xz{u_{xz}} - yz{u_{yz}}
   + \left[ {{c_2} - \left( {{a_2} + b + 1} \right)z} \right]{u_z} + {a_2}x{u_x} - by{u_y} - {a_2}bu = 0, \hfill \\
\end{gathered}
\end{array}} \right.
$

where $u\equiv \,\,   {F_{36d}}\left( {{a_1},{a_2},b;{c_1},{c_2};x,y,z} \right) $.

Particular solutions:

$
{u_1} = {F_{36d}}\left( {{a_1},{a_2},b;{c_1},{c_2};x,y,z} \right),
$

$
{u_2} = {y^{1 - {c_1}}}{F_{36d}}\left( {1 - {c_1} + {a_1},1 - {c_1} + {a_2},b;2 - {c_1},{c_2};x,y,z} \right),
$

$
{u_3} = {z^{1 - {c_2}}}{F_{36d}}\left( {{a_1},1 - {c_2} + {a_2},1 - {c_2} + b;{c_1},2 - {c_2};x,y,z} \right),
$

$
{u_4} = {y^{1 - {c_1}}}{z^{1 - {c_2}}}{F_{36d}}\left( {1 - {c_1} + {a_1},2 - {c_1} - {c_2} + {a_2},1 - {c_2} + b;2 - {c_1},2 - {c_2};x,y,z} \right).
$

\bigskip

\begin{equation}
{F_{36e}}\left( {{a_1},{a_2},b;{c_1},{c_2};x,y,z} \right) = \sum\limits_{m,n,p = 0}^\infty  {} \frac{{{{\left( {{a_1}} \right)}_{2m + n}}{{\left( {{a_2}} \right)}_{n + p}}{{\left( b \right)}_{p - n}}}}{{{{\left( {{c_1}} \right)}_m}{{\left( {{c_2}} \right)}_p}}}\frac{x^m}{m!}\frac{y^n}{n!}\frac{z^p}{p!},
\end{equation}

region of convergence:
$$
 \left\{ {r < \frac{1}{4},\,\,\,t < 1,\,\,\,\sqrt s  < \sqrt {1 - 2\sqrt r } \left( {\sqrt {1 + t}  - \sqrt t } \right)} \right\}.
$$

System of partial differential equations:

$
\left\{ {\begin{array}{*{20}{l}}
  \begin{gathered}
  x\left( {1 - 4x} \right){u_{xx}} - 4xy{u_{xy}} - {y^2}{u_{yy}}
   + \,\left[ {{c_1} - 2\left( {2{a_1} + 3} \right)x} \right]{u_x} - 2\left( {{a_1} + 1} \right)y{u_y} - {a_1}\left( {1 + {a_1}} \right)u = 0, \hfill \\
\end{gathered}  \\
  \begin{gathered}
  y\left( {1 + y} \right){u_{yy}} + 2xy{u_{xy}} + 2xz{u_{xz}} - \left( {1 - y} \right)z{u_{yz}} \hfill \\
  \,\,\,\,\,\,\,\,\,  + \left[ {1 - b + \left( {{a_1} + {a_2} + 1} \right)y} \right]{u_y} + 2{a_2}x{u_x} + {a_1}z{u_z} + {a_1}{a_2}u = 0, \hfill \\
\end{gathered}  \\
  \begin{gathered}
  z\left( {1 - z} \right){u_{zz}} + {y^2}{u_{yy}}
   + \left[ {{c_2} - \left( {{a_2} + b + 1} \right)z} \right]{u_z} + \left( {{a_2} - b + 1} \right)y{u_y} - {a_2}bu = 0, \hfill \\
\end{gathered}
\end{array}} \right.
$

where $u\equiv \,\,   {F_{36e}}\left( {{a_1},{a_2},b;{c_1},{c_2};x,y,z} \right)$.

Particular solutions:

$
{u_1} = {F_{36e}}\left( {{a_1},{a_2},b;{c_1},{c_2};x,y,z} \right),
$

$
{u_2} = {x^{1 - {c_1}}}{F_{36e}}\left( {2 - 2{c_1} + {a_1},{a_2},b;2 - {c_1},{c_2};x,y,z} \right),
$

$
{u_3} = {z^{1 - {c_3}}}{F_{36e}}\left( {{a_1},1 - {c_3} + {a_2},1 - {c_3} + b;{c_1},2 - {c_2};x,y,z} \right),
$

$
{u_4} = {x^{1 - {c_1}}}{z^{1 - {c_3}}}{F_{36e}}\left( {2 - 2{c_1} + {a_1},1 - {c_3} + {a_2},1 - {c_3} + b;2 - {c_1},2 - {c_2};x,y,z} \right).
$

\bigskip

\begin{equation}
{F_{36f}}\left( {a,{b_1},{b_2};c;x,y,z} \right)  = \sum\limits_{m,n,p = 0}^\infty  {} \frac{{{{\left( a \right)}_p}{{\left( {{b_1}} \right)}_{2m + n - p}}{{\left( {{b_2}} \right)}_{n + p - m}}}}{{{{\left( c \right)}_n}}}\frac{x^m}{m!}\frac{y^n}{n!}\frac{z^p}{p!},
\end{equation}

region of convergence:
$$
 \left\{ {r < \frac{1}{4},\,\,\,t < \frac{{\sqrt {1 + 4t}  - 1}}{{2r}},\,\,\,s < \min \left\{ {{\Psi _1}\left( r \right),{\Psi _2}\left( r \right),\frac{{{{\left( {1 - t - r{t^2}} \right)}^2}}}{{4t}}} \right\}} \right\}.
$$

System of partial differential equations:

$
\left\{ {\begin{array}{*{20}{l}}
  \begin{gathered}
  x\left( {1 + 4x} \right){u_{xx}} - \left( {1 - 4x} \right)y{u_{xy}} - \left( {1 + 4x} \right)z{u_{xz}} - 2yz{u_{yz}}
   + {y^2}{u_{yy}} + {z^2}{u_{zz}} \hfill \\ \,\,\,\,\,\,\,\,\,+ \left[ {1 - {b_2} + 2\left( {2{b_1} + 3} \right)x} \right]{u_x} + 2\left( {{b_1} + 1} \right)y{u_y}
    - 2{b_1}z{u_z} + {b_1}\left( {1 + {b_1}} \right)u = 0, \hfill \\
\end{gathered}  \\
  \begin{gathered}
  y\left( {1 - y} \right){u_{yy}} - xy{u_{xy}} - 3xz{u_{xz}} + 2{x^2}{u_{xx}} + {z^2}{u_{zz}} \hfill \\
 \,\,\,\,\,\,\,\,\,  + \left[ {c - \left( {{b_1} + {b_2} + 1} \right)y} \right]{u_y} + \left( {{b_1} - 2{b_2} + 2} \right)x{u_x}
    + \left( {{b_2} - {b_1} + 1} \right)z{u_z} - {b_1}{b_2}u = 0, \hfill \\
\end{gathered}  \\
  \begin{gathered}
  z\left( {1 + z} \right){u_{zz}} - x(2+z){u_{xz}} - y\left( {1 - z} \right){u_{yz}} - ax{u_x} + ay{u_y}\hfill\\
  \,\,\,\,\,\,\,\,\, + \left[ {1 - {b_1} + \left( {a + {b_2} + 1} \right)z} \right]{u_z}  + a{b_2}u = 0, \hfill \\
\end{gathered}
\end{array}} \right.
$

where $u\equiv \,\,   {F_{36f}}\left( {a,{b_1},{b_2};c;x,y,z} \right)$.

Particular solutions:

$
{u_1} = {F_{36f}}\left( {a,{b_1},{b_2};c;x,y,z} \right),
$

$
{u_2} = {y^{1 - c}}{F_{36f}}\left( {a,1 - c + {b_1},1 - c + {b_2};2 - c;x,y,z} \right).
$

\bigskip

\begin{equation}
{F_{36g}}\left( {a,{b_1},{b_2};c;x,y,z} \right) = \sum\limits_{m,n,p = 0}^\infty  {} \frac{{{{\left( a \right)}_{n + p}}{{\left( {{b_1}} \right)}_{2m + n - p}}{{\left( {{b_2}} \right)}_{p - m}}}}{{{{\left( c \right)}_n}}}\frac{x^m}{m!}\frac{y^n}{n!}\frac{z^p}{p!},
\end{equation}

region of convergence:
$$
\begin{gathered}
    \left\{ {r < \frac{1}{4},\,\,\,t + r{t^2} < 1,\,\,\,s < \min \left\{ {{U^ + }\left( {{w_1}} \right),{U^ - }\left( {{w_2}} \right)} \right\}} \right\}, \hfill \\
  {P^ \pm }\left( w \right) = r{w^3} \pm \left( {1 + t} \right)w \mp 2t, \hfill \\
{w_1}:\,\,{\rm{the\,\,root\,\,in}}\,\,\left( {t,\infty } \right)\,\,{\rm{of}}\,\,{P^ + }\left( w \right) = 0, \hfill \\
{w_2}:\,\,{\rm{the\,\,root\,\,in}}\,\,\left( {\frac{1}{{\sqrt r }},\infty } \right)\,\,{\rm{of}}\,\,{P^ - }\left( w \right) = 0, \hfill \\
  {U^ \pm }\left( w \right) = {\left( {1 - \frac{t}{w}} \right)^2}\frac{{r{w^2} \pm 1}}{t}. \hfill \\
\end{gathered}
$$

System of partial differential equations:

$
\left\{ {\begin{array}{*{20}{l}}
  \begin{gathered}
  x\left( {1 + 4x} \right){u_{xx}} + 4xy{u_{xy}} - \left( {1 + 4x} \right)z{u_{xz}} - 2yz{u_{yz}}
   + {y^2}{u_{yy}} + {z^2}{u_{zz}} \hfill \\\,\,\,\,\,\,\,\,\, + \left[ {1 - {b_2} + 2\left( {2{b_1} + 3} \right)x} \right]{u_x} + 2\left( {{b_1} + 1} \right)y{u_y}
    - 2{b_1}z{u_z} + {b_1}\left( {1 + {b_1}} \right)u = 0, \hfill \\
\end{gathered}  \\
  \begin{gathered}
  y\left( {1 - y} \right){u_{yy}} - 2xy{u_{xy}} - 2xz{u_{xz}} + {z^2}{u_{zz}}\hfill \\ \,\,\,\,\,\,\,\,\,
   + \left[ {c - \left( {a + {b_1} + 1} \right)y} \right]{u_y} - 2ax{u_x}
   + \left( {a - {b_1} + 1} \right)z{u_z} - a{b_1}u = 0, \hfill \\
\end{gathered}  \\
  \begin{gathered}
  z\left( {1 + z} \right){u_{zz}} - xy{u_{xy}} - x\left( {2 + z} \right){u_{xz}} - y\left( {1 - z} \right){u_{yz}}\hfill \\ \,\,\,\,\,\,\,\,\,
   + \left[ {1 - {b_1} + \left( {{b_2} + a + 1} \right)z} \right]{u_z} - ax{u_x} + {b_2}y{u_y} + a{b_2}u = 0, \hfill \\
\end{gathered}
\end{array}} \right.
$

where $u\equiv \,\,   {F_{36g}}\left( {a,{b_1},{b_2};c;x,y,z} \right)$.

Particular solutions:

$
{u_1} = {F_{36g}}\left( {a,{b_1},{b_2};c;x,y,z} \right),
$

$
{u_2} = {y^{1 - c}}{F_{36g}}\left( {1 - c + a,1 - c + {b_1},{b_2};2 - c;x,y,z} \right).
$

\bigskip

\begin{equation}
{F_{36h}}\left( {a,{b_1},{b_2};c;x,y,z} \right)  = \sum\limits_{m,n,p = 0}^\infty  {} \frac{{{{\left( a \right)}_{n + p}}{{\left( {{b_1}} \right)}_{2m + n - p}}{{\left( {{b_2}} \right)}_{p - n}}}}{{{{\left( c \right)}_m}}}\frac{x^m}{m!}\frac{y^n}{n!}\frac{z^p}{p!},
\end{equation}

region of convergence:
$$
 \left\{ {s + t < 1,\,\,\,\sqrt r  < \frac{{\sqrt {1 - 4st}  - \left| {1 - 2t} \right|}}{{4t}}} \right\}.
$$

System of partial differential equations:

$
\left\{ {\begin{array}{*{20}{l}}
  \begin{gathered}
  x\left( {1 - 4x} \right){u_{xx}} - 4xy{u_{xy}} + 4xz{u_{xz}} + 2yz{u_{yz}} - {y^2}{u_{yy}} - {z^2}{u_{zz}} \hfill \\\,\,\,\,\,\,\,\,\,
   + \left[ {c - 2\left( {2{b_1} + 3} \right)x} \right]{u_x} - 2\left( {{b_1} + 1} \right)y{u_y} + 2{b_1}z{u_z} - {b_1}\left( {1 + {b_1}} \right)u = 0, \hfill \\
\end{gathered}  \\
  \begin{gathered}
  y\left( {1 + y} \right){u_{yy}} + 2xy{u_{xy}} + 2xz{u_{xz}} - z{u_{yz}} - {z^2}{u_{zz}} \hfill \\ \,\,\,\,\,\,\,\,\,
  + \left[ {1 - {b_2} + \left( {a + {b_1} + 1} \right)y} \right]{u_y} + 2ax{u_x} - \left( {a - {b_1} + 1} \right)z{u_z} + a{b_1}u = 0, \hfill \\
\end{gathered}  \\
  \begin{gathered}
  z\left( {1 + z} \right){u_{zz}} - 2x{u_{xz}} - y{u_{yz}} - {y^2}{u_{yy}}\hfill\\\,\,\,\,\,\,\,\,\,
   + \left[ {1 - {b_1} + \left( {a + {b_2} + 1} \right)z} \right]{u_z} + \left( {a - {b_2} + 1} \right)y{u_y} + a{b_2}u = 0, \hfill \\
\end{gathered}
\end{array}} \right.
$

where $u\equiv \,\,   {F_{36h}}\left( {a,{b_1},{b_2};c;x,y,z} \right)$.

Particular solutions:

$
{u_1} = {F_{36h}}\left( {a,{b_1},{b_2};c;x,y,z} \right),
$

$
{u_2} = {x^{1 - c}}{F_{36h}}\left( {a,2 - 2c + {b_1},{b_2};2 - c;x,y,z} \right).
$

\bigskip

\begin{equation}
{F_{36i}}\left( {a,{b_1},{b_2};c;x,y,z} \right)  = \sum\limits_{m,n,p = 0}^\infty  {} \frac{{{{\left( a \right)}_{2m + n}}{{\left( {{b_1}} \right)}_{n + p - m}}{{\left( {{b_2}} \right)}_{p - n}}}}{{{{\left( c \right)}_p}}}\frac{x^m}{m!}\frac{y^n}{n!}\frac{z^p}{p!},
\end{equation}

region of convergence:
$$
\begin{gathered}
   \left\{ {r < \frac{1}{4},\,\,\,s < \min \left\{ {{\Psi _1}\left( r \right),{\Psi _2}\left( r \right)} \right\},}  {\sqrt t  < \min \left\{ {1,{U^ + }\left( {{w_1}} \right),{U^ - }\left( {{w_2}} \right),{U^ + }\left( {{w_3}} \right)} \right\}} \right\}, \hfill \\
  {P_{rs}}\left( w \right) = 10{w^3} - 9{w^2} + 2\left( {1 - r} \right)w + r\left( {1 - s} \right), \hfill \\
{w_1}:\,\,{\rm{the\,\,root\,\,in}}\,\,\left( {\frac{1}{6}\left[ {1 + a\left( r \right)} \right],\frac{1}{2}} \right)\,\,{\rm{of}}\,\,{P_{rs}}\left( w \right) = 0, \hfill \\
{w_2}:\,\,{\rm{the\,\,smaller\,\,root\,\,in}}\,\,\left( {0,\frac{1}{2}} \right)\,\,{\rm{of}}\,\,{P_{ - r, - s}}\left( w \right) = 0, \hfill \\
{w_3}:\,\,{\rm{the\,\,greater\,\,root\,\,in}}\,\,\left( {\frac{1}{6}\left[ {1 + a\left( r \right)} \right],\frac{1}{2}} \right)\,\,{\rm{of}}\,\,{P_{r, - s}}\left( w \right) = 0, \hfill \\
  {U^ \pm }\left( w \right) = \frac{{\left( {3{w^2} - w \mp r} \right)\sqrt {1 - 2w} }}{{r\sqrt s }}. \hfill \\
\end{gathered}
$$

System of partial differential equations:

$
\left\{ {\begin{array}{*{20}{l}}
  \begin{gathered}
  x\left( {1 + 4x} \right){u_{xx}} - \left( {1 - 4x} \right)y{u_{xy}} - z{u_{xz}} + {y^2}{u_{yy}} \hfill \\
 \,\,\,\,\,\,\,\,\,  + \left[ {1 - {b_1} + 2\left( {2a + 3} \right)x} \right]{u_x} + 2\left( {a + 1} \right)y{u_y} + a\left( {1 + a} \right)u = 0, \hfill \\
\end{gathered}  \\
  \begin{gathered}
  y\left( {1 + y} \right){u_{yy}} + xy{u_{xy}} + 2xz{u_{xz}} - \left( {1 - y} \right)z{u_{yz}} - 2{x^2}{u_{xx}} \hfill \\
 \,\,\,\,\,\,\,\,\,  + \left[ {1 - {b_2} + \left( {a + {b_1} + 1} \right)y} \right]{u_y} - \left( {a - 2{b_1} + 2} \right)x{u_x} + az{u_z} + a{b_1}u = 0, \hfill \\
\end{gathered}  \\
  \begin{gathered}
  z\left( {1 - z} \right){u_{zz}} - xy{u_{xy}} + xz{u_{xz}} + {y^2}{u_{yy}} \hfill \\\,\,\,\,\,\,\,\,\,
   + \left[ {c - \left( {{b_1} + {b_2} + 1} \right)z} \right]{u_z} + {b_2}x{u_x} + \left( {{b_1} - {b_2} + 1} \right)y{u_y} - {b_1}{b_2}u = 0, \hfill \\
\end{gathered}
\end{array}} \right.
$

where $u\equiv \,\,   {F_{36i}}\left( {a,{b_1},{b_2};c;x,y,z} \right)$.

Particular solutions:

$
{u_1} = {F_{36i}}\left( {a,{b_1},{b_2};c;x,y,z} \right),
$

$
{u_2} = {z^{1 - c}}{F_{36i}}\left( {a,1 - c + {b_1},1 - c + {b_2};2 - c;x,y,z} \right).
$

\bigskip

\begin{equation}
{F_{36j}}\left( {{b_1},{b_2},{b_3};x,y,z} \right)  = \sum\limits_{m,n,p = 0}^\infty  {} {{{{\left( {{b_1}} \right)}_{2m + n - p}}{{\left( {{b_2}} \right)}_{n + p - m}}{{\left( {{b_3}} \right)}_{p - n}}}}\frac{x^m}{m!}\frac{y^n}{n!}\frac{z^p}{p!},
\end{equation}

region of convergence:
$$
\begin{gathered}
    \left\{ {r < \frac{1}{4},\,\,\,t < \frac{2}{{1 + \sqrt {1 + 4r} }},}  {s < \left( {1 + 2rt} \right)\min \left\{ {\left( {1 + t + r{t^2}} \right){\Psi _1}\left( \xi  \right),\left( {1 - t - r{t^2}} \right){\Psi _2}\left( \eta  \right)} \right\}} \right\}, \hfill \\
  \xi  = \frac{{r\left( {1 + t + r{t^2}} \right)}}{{{{\left( {1 + 2rt} \right)}^2}}},\,\,\,\eta  = \frac{{r\left( {1 - t - r{t^2}} \right)}}{{{{\left( {1 + 2rt} \right)}^2}}}. \hfill \\
\end{gathered}
$$

System of partial differential equations:

$
\left\{ {\begin{array}{*{20}{l}}
  \begin{gathered}
  x\left( {1 + 4x} \right){u_{xx}} - \left( {1 - 4x} \right)y{u_{xy}} - \left( {1 + 4x} \right)z{u_{xz}} - 2yz{u_{yz}}
   + {y^2}{u_{yy}} + {z^2}{u_{zz}} \hfill \\ \,\,\,\,\,\,\,\,\, + \left[ {1 - {b_2} + 2\left( {2{b_1} + 3} \right)x} \right]{u_x} + 2\left( {{b_1} + 1} \right)y{u_y}
    - 2{b_1}z{u_z} + {b_1}\left( {1 + {b_1}} \right)u = 0, \hfill \\
\end{gathered}  \\
  \begin{gathered}
  y\left( {1 + y} \right){u_{yy}} + xy{u_{xy}} + 3xz{u_{xz}} - z{u_{yz}} - 2{x^2}{u_{xx}} - {z^2}{u_{zz}} \hfill \\
 \,\,\,\,\,\,\,\,\,  + \left[ {1 - {b_3} + \left( {{b_1} + {b_2} + 1} \right)y} \right]{u_y} - \left( {{b_1} - 2{b_2} + 2} \right)x{u_x}
    + \left( {{b_1} - {b_2} - 1} \right)z{u_z} + {b_1}{b_2}u = 0, \hfill \\
\end{gathered}  \\
  \begin{gathered}
  z\left( {1 + z} \right){u_{zz}} + xy{u_{xy}} - x\left( {2 + z} \right){u_{xz}} - y{u_{yz}} - {y^2}{u_{yy}} \hfill \\
\,\,\,\,\,\,\,\,\,   + \left[ {1 - {b_1} + \left( {{b_2} + {b_3} + 1} \right)z} \right]{u_z} - {b_3}x{u_x}
   - \left( {{b_2} - {b_3} + 1} \right)y{u_y} + {b_2}{b_3}u = 0, \hfill \\
\end{gathered}
\end{array}} \right.
$

where $u\equiv \,\,   {F_{36j}}\left( {{b_1},{b_2},{b_3};x,y,z} \right)$.

\bigskip

\begin{equation}
{F_{37a}}\left( {{a_1},{a_2},{a_2};{c_1},{c_2};x,y,z} \right)  = \sum\limits_{m,n,p = 0}^\infty  {} \frac{{{{\left( {{a_1}} \right)}_{2m + n}}{{\left( {{a_2}} \right)}_{n + p}}{{\left( {{a_3}} \right)}_p}}}{{{{\left( {{c_1}} \right)}_{m + n}}{{\left( {{c_2}} \right)}_p}}}\frac{x^m}{m!}\frac{y^n}{n!}\frac{z^p}{p!},
\end{equation}

first appearance of this function in the literature, and old notation:
[12], $X_{14}$,

region of convergence:
$$
\left\{ {r < \frac{1}{4},\,\,\,t < 1,\,\,\,s < \left( {1 - t} \right)\left( {\frac{1}{2} + \frac{1}{2}\sqrt {1 - 4r} } \right)} \right\}.
$$

System of partial differential equations:

$
\left\{ \begin{array}{*{20}{l}}
  \begin{gathered}
  x\left( {1 - 4x} \right){u_{xx}} + \left( {1 - 4x} \right)y{u_{xy}} - {y^2}{u_{yy}}\hfill\\\,\,\,\,\,\,\,\,\,
   + \left[ {{c_1} - 2\left( {2{a_1} + 3} \right)x} \right]{u_x} - 2\left( {{a_1} + 1} \right)y{u_y} - {a_1}\left( {1 + {a_1}} \right)u = 0, \hfill \\
\end{gathered}  \\
  \begin{gathered}
  y\left( {1 - y} \right){u_{yy}} + x\left( {1 - 2y} \right){u_{xy}} - 2xz{u_{xz}} - yz{u_{yz}} \hfill \\ \,\,\,\,\,\,\,\,\,
   + \left[ {{c_1} - \left( {{a_1} + {a_2} + 1} \right)y} \right]{u_y} - 2{a_2}x{u_x} - {a_1}z{u_z} - {a_1}{a_2}u = 0, \hfill \\
\end{gathered}  \\
  \begin{gathered}
  z\left( {1 - z} \right){u_{zz}} - yz{u_{yz}}
   + \left[ {{c_2} - \left( {{a_2} + {a_3} + 1} \right)z} \right]{u_z} - {a_3}y{u_y} - {a_2}{a_3}u = 0, \hfill \\
\end{gathered}
\end{array} \right.
$

where $u\equiv \,\,   {F_{37a}}\left( {{a_1},{a_2},{a_2};{c_1},{c_2};x,y,z} \right)$.

Particular solutions:

$
{u_1} = {F_{37a}}\left( {{a_1},{a_2},{a_2};{c_1},{c_2};x,y,z} \right),
$

$
{u_2} = {z^{1 - {c_2}}}{F_{37a}}\left( {{a_1},1 - {c_2} + {a_2},1 - {c_2} + {a_3};{c_1},2 - {c_2};x,y,z} \right).
$

\bigskip

\begin{equation}
{F_{37b}}\left( {{a_1},{a_2},b;c;x,y,z} \right)  = \sum\limits_{m,n,p = 0}^\infty  {} \frac{{{{\left( {{a_1}} \right)}_{n + p}}{{\left( {{a_2}} \right)}_p}{{\left( b \right)}_{2m + n - p}}}}{{{{\left( c \right)}_{m + n}}}}\frac{x^m}{m!}\frac{y^n}{n!}\frac{z^p}{p!},
\end{equation}

region of convergence:
$$
 \left\{ {r < \frac{1}{4},\,\,\,t < \frac{1}{{1 + 2\sqrt r }},\,\,\,s < \min \left\{ {\frac{1}{2} + \frac{1}{2}\sqrt {1 - 4r} ,\frac{{{{\left( {1 - t} \right)}^2} - 4r{t^2}}}{{4t}}} \right\}} \right\}.
$$

System of partial differential equations:

$
\left\{ {\begin{array}{*{20}{l}}
  \begin{gathered}
  x\left( {1 - 4x} \right){u_{xx}} + \left( {1 - 4x} \right)y{u_{xy}} + 4xz{u_{xz}} + 2yz{u_{yz}} - {y^2}{u_{yy}} - {z^2}{u_{zz}} \hfill \\
  \,\,\,\,\,\,\,\,\, + \left[ {c - 2\left( {2b + 3} \right)x} \right]{u_x} - 2\left( {b + 1} \right)y{u_y} + 2bz{u_z} - b\left( {1 + b} \right)u = 0, \hfill \\
\end{gathered}  \\
  \begin{gathered}
  y\left( {1 - y} \right){u_{yy}} + x\left( {1 - 2y} \right){u_{xy}} - 2xz{u_{xz}} + {z^2}{u_{zz}} \hfill \\
 \,\,\,\,\,\,\,\,\,  + \left[ {c - \left( {{a_1} + b + 1} \right)y} \right]{u_y} - 2{a_1}x{u_x} + \left( {{a_1} - b + 1} \right)z{u_z} - {a_1}bu = 0 \hfill \\
\end{gathered}  \\
  \begin{gathered}
  z\left( {1 + z} \right){u_{zz}} - 2x{u_{xz}} - y\left( {1 - z} \right){u_{yz}}
   + \left[ {1 - b + \left( {{a_1} + {a_2} + 1} \right)z} \right]{u_z} + {a_2}y{u_y} + {a_1}{a_2}u = 0, \hfill \\
\end{gathered}
\end{array}} \right.
$

where $u\equiv \,\,   {F_{37b}}\left( {{a_1},{a_2},b;c;x,y,z} \right)$.

\bigskip

\begin{equation}
{F_{37c}}\left( {{a_1},{a_2},b;c;x,y,z} \right)  = \sum\limits_{m,n,p = 0}^\infty  {} \frac{{{{\left( {{a_1}} \right)}_{2m + n}}{{\left( {{a_2}} \right)}_{n + p}}{{\left( b \right)}_{p - m - n}}}}{{{{\left( c \right)}_p}}}\frac{x^m}{m!}\frac{y^n}{n!}\frac{z^p}{p!},
\end{equation}

region of convergence:
$$
\begin{gathered}
   \left\{ {r < \frac{1}{4},\,\,\,s < \frac{1}{2} + \frac{1}{2}\sqrt {1 - 4r} ,\,\,\,t < \min \left\{ {1,\frac{{{\xi ^2}{\eta ^2}}}{{{r^2}{s^2}\left( {\xi  + \eta } \right)}}{\Phi _1}\left( { - \frac{{\xi \eta }}{{{{\left( {\xi  + \eta } \right)}^2}}}} \right)} \right\}} \right\}, \hfill \\
  \xi  = s\sqrt {1 - 4r} ,\,\,\eta  = r\left( {1 - \frac{{2s}}{{1 + \sqrt {1 - 4r} }}} \right). \hfill \\
\end{gathered}
$$

System of partial differential equations:

$
\left\{ {\begin{array}{*{20}{l}}
  \begin{gathered}
  x\left( {1 + 4x} \right){u_{xx}} - z{u_{xz}} + \left( {1 + 4x} \right)y{u_{xy}} + {y^2}{u_{yy}} \hfill \\
 \,\,\,\,\,\,\,\,\,  + \left[ {1 - b + 2\left( {2{a_1} + 3} \right)x} \right]{u_x} + \left( {{a_1} + {a_1} + 2} \right)y{u_y} + {a_1}\left( {1 + {a_1}} \right)u = 0, \hfill \\
\end{gathered}  \\
  \begin{gathered}
  y\left( {1 + y} \right){u_{yy}} + x\left( {1 + 2y} \right){u_{xy}} + 2xz{u_{xz}} - z\left( {1 - y} \right){u_{yz}} \hfill \\
  \,\,\,\,\,\,\,\,\, + \left[ {1 - b + \left( {{a_1} + {a_2} + 1} \right)y} \right]{u_y} + 2{a_2}x{u_x} + {a_1}z{u_z} + {a_1}{a_2}u = 0, \hfill \\
\end{gathered}  \\
  \begin{gathered}
  z\left( {1 - z} \right){u_{zz}} + xy{u_{xy}} + xz{u_{xz}} + {y^2}{u_{yy}}\hfill\\\,\,\,\,\,\,\,\,\,
   + \left[ {c - \left( {{a_2} + b + 1} \right)z} \right]{u_z} + {a_2}x{u_x} + \left( {{a_2} - b + 1} \right)y{u_y} - {a_2}bu = 0, \hfill \\
\end{gathered}
\end{array}} \right.
$

where $u\equiv \,\,   {F_{37c}}\left( {{a_1},{a_2},b;c;x,y,z} \right)$.

Particular solutions:

$
{u_1} = {F_{37c}}\left( {{a_1},{a_2},b;c;x,y,z} \right),
$

$
{u_2} = {z^{1 - c}}{F_{37c}}\left( {{a_1},1 - c + {a_2},1 - c + b;2 - c;x,y,z} \right).
$

\bigskip

\begin{equation}
{F_{37d}}\left( {a,{b_1},{b_2};x,y,z} \right)  = \sum\limits_{m,n,p = 0}^\infty  {} {{{{\left( a \right)}_{n + p}}{{\left( {{b_1}} \right)}_{2m + n - p}}{{\left( {{b_2}} \right)}_{p - m - n}}}}\frac{x^m}{m!}\frac{y^n}{n!}\frac{z^p}{p!},
\end{equation}

region of convergence:
$$
\begin{gathered}
  \left\{ {r < \frac{1}{4},\,\,\,s < \frac{1}{2} + \frac{1}{2}\sqrt {1 - 4r} ,\,\,\,t < \min \left[ {{w_1}\left( {1 - s{w_1}} \right),{w_2}\left( {1 - s{w_2}} \right)} \right]} \right\}, \hfill \\
  \left. {\begin{array}{*{20}{c}}
  {{w_1}} \\
  {{w_2}}
\end{array}} \right\} = \frac{2}{{1 + \sqrt {1 \pm 4r} }}. \hfill \\
\end{gathered}
$$

System of partial differential equations:

$
\left\{ {\begin{array}{*{20}{l}}
  \begin{gathered}
  x\left( {1 + 4x} \right){u_{xx}} + \left( {1 + 4x} \right)y{u_{xy}} - \left( {1 + 4x} \right)z{u_{xz}} - 2yz{u_{yz}}
   + {y^2}{u_{yy}} + {z^2}{u_{zz}}  \hfill \\
  \,\,\,\,\,\,\,\,\, + \left[ {1 - {b_2} + 2\left( {2{b_1} + 3} \right)x} \right]{u_x}
    + \left( {{b_1} + {b_1} + 2} \right)y{u_y} - 2{b_1}z{u_z} + {b_1}\left( {1 + {b_1}} \right)u = 0, \hfill \\
\end{gathered}  \\
  \begin{gathered}
  y\left( {1 + y} \right){u_{yy}} + x\left( {1 + 2y} \right){u_{xy}} + 2xz{u_{xz}} - z{u_{yz}} - {z^2}{u_{zz}} \hfill \\
  \,\,\,\,\,\,\,\,\, + \left[ {1 - {b_2} + \left( {a + {b_1} + 1} \right)y} \right]{u_y} + 2ax{u_x}
    - \left( {a - {b_1} + 1} \right)z{u_z} + a{b_1}u = 0, \hfill \\
\end{gathered}  \\
  \begin{gathered}
  z\left( {1 + z} \right){u_{zz}} - xy{u_{xy}} - x\left( {2 + z} \right){u_{xz}} - y{u_{yz}} - {y^2}{u_{yy}} \hfill \\
 \,\,\,\,\,\,\,\,\, + \left[ {1 - {b_1} + \left( {a + {b_2} + 1} \right)z} \right]{u_z} - ax{u_x}
  - \left( {a - {b_2} + 1} \right)y{u_y} + a{b_2}u = 0, \hfill \\
\end{gathered}
\end{array}} \right.
$

where $u\equiv \,\,   {F_{37d}}\left( {a,{b_1},{b_2};x,y,z} \right)$.

\bigskip

\begin{equation}
{F_{38a}}\left( {{a_1},{a_2},{a_3};{c_1},{c_2};x,y,z} \right) = \sum\limits_{m,n,p = 0}^\infty  {} \frac{{{{\left( {{a_1}} \right)}_{2m + n}}{{\left( {{a_2}} \right)}_{n + p}}{{\left( {{a_3}} \right)}_p}}}{{{{\left( {{c_1}} \right)}_{m + p}}{{\left( {{c_2}} \right)}_n}}}\frac{x^m}{m!}\frac{y^n}{n!}\frac{z^p}{p!},
\end{equation}

first appearance of this function in the literature, and old notation:
[12], $X_{16}$,

region of convergence:
$$
\begin{gathered}
   \left\{ {s < 1,\,\,\,t  \leq T \left( s \right),\,\,\,r < \frac{1}{4}\left( {1 - {s^2}} \right)} \right\}
  \cup \left\{ {s < 1,\,\,\,T\left( s \right) < t < 1 - s,\,\,\,r \leq \frac{t}{{1 - t}}{\Theta _2}\left( {\frac{s}{{1 - t}}} \right)} \right\}, \hfill \\
  T\left( s \right)=\left\{ {\begin{array}{*{20}{c}}
  {{{\left( {\dfrac{{1 - s}}{{1 - s/2}}} \right)}^2},\,\,\,0 < s < \dfrac{2}{3},} \\
  {\dfrac{{1 - s}}{{2s}},\,\,\,\dfrac{2}{3}\leq s < 1.}
\end{array}} \right. \hfill \\
\end{gathered}
$$

System of partial differential equations:

$
\left\{ {\begin{array}{*{20}{l}}
  \begin{gathered}
  x\left( {1 - 4x} \right){u_{xx}} - 4xy{u_{xy}} + z{u_{xz}} - {y^2}{u_{yy}}\hfill\\\,\,\,\,\,\,\,\,\,
   + \left[ {{c_1} - 2\left( {2{a_1} + 3} \right)x} \right]{u_x} - 2\left( {{a_1} + 1} \right)y{u_y} - {a_1}\left( {1 + {a_1}} \right)u = 0, \hfill \\
\end{gathered}  \\
  \begin{gathered}
  y\left( {1 - y} \right){u_{yy}} - 2xy{u_{xy}} - 2xz{u_{xz}} - yz{u_{yz}}\hfill\\\,\,\,\,\,\,\,\,\,
  + \left[ {{c_2} - \left( {{a_1} + {a_2} + 1} \right)y} \right]{u_y} - 2{a_2}x{u_x} - {a_1}z{u_z} - {a_1}{a_2}u = 0, \hfill \\
\end{gathered}  \\
  \begin{gathered}
  z\left( {1 - z} \right){u_{zz}} + x{u_{xz}} - yz{u_{yz}}
   + \left[ {{c_1} - \left( {{a_2} + {a_3} + 1} \right)z} \right]{u_z} - {a_3}y{u_y} - {a_2}{a_3}u = 0, \hfill \\
\end{gathered}
\end{array}} \right.
$

where $u\equiv \,\,  {F_{38a}}\left( {{a_1},{a_2},{a_3};{c_1},{c_2};x,y,z} \right)$.

Particular solutions:

$
{u_1} = {F_{38a}}\left( {{a_1},{a_2},{a_3};{c_1},{c_2};x,y,z} \right),
$

$
{u_2} = {y^{1 - {c_2}}}{F_{38a}}\left( {1 - {c_2} + {a_1},1 - {c_2} + {a_2},{a_3};{c_1},2 - {c_2};x,y,z} \right).
$

\bigskip

\begin{equation}
{F_{38b}}\left( {{a_1},{a_2},{b};c;x,y,z} \right)  = \sum\limits_{m,n,p = 0}^\infty  {} \frac{{{{\left( {{a_1}} \right)}_{2m + n}}{{\left( {{a_2}} \right)}_{n + p}}{{\left( {{b}} \right)}_{p - n}}}}{{{{\left( c \right)}_{m + p}}}}\frac{x^m}{m!}\frac{y^n}{n!}\frac{z^p}{p!},
\end{equation}

region of convergence:
$$
\begin{gathered}
  \left\{ {s + 2\sqrt r  < 1,\,\,\,t < \min \left\{ {1,\frac{{r\left( {w_1^2 - {s^2}} \right)}}{{2sw_1^2\left( {1 - {w_1}} \right)}}} \right\}} \right\}, \hfill \\
  P\left( w \right) = {w^3} - 3s{w^2} - \left( {1 - 4s - 4r} \right)w - s + 4rs, \hfill \\
{w_1}:\,\,{\rm{the\,\,root\,\,in}}\,\,\left( {s,1} \right)\,\,{\rm{of}}\,\,P\left( w \right) = 0. \hfill \\
\end{gathered}
$$

System of partial differential equations:

$
\left\{ {\begin{array}{*{20}{l}}
  \begin{gathered}
  x\left( {1 - 4x} \right){u_{xx}} - 4xy{u_{xy}} + z{u_{xz}} - {y^2}{u_{yy}}\hfill\\\,\,\,\,\,\,\,\,\,
   + \left[ {c - 2\left( {2{a_1} + 3} \right)x} \right]{u_x} - 2\left( {{a_1} + 1} \right)y{u_y} - {a_1}\left( {1 + {a_1}} \right)u = 0, \hfill \\
\end{gathered}  \\
  \begin{gathered}
  y\left( {1 + y} \right){u_{yy}} + 2xy{u_{xy}} + 2xz{u_{xz}} - \left( {1 - y} \right)z{u_{yz}} \hfill \\
 \,\,\,\,\,\,\,\,\,  + \left[ {1 - b + \left( {{a_1} + {a_2} + 1} \right)y} \right]{u_y} + 2{a_2}x{u_x} + {a_1}z{u_z} + {a_1}{a_2}u = 0, \hfill \\
\end{gathered}  \\
  \begin{gathered}
  z\left( {1 - z} \right){u_{zz}} + x{u_{xz}} + {y^2}{u_{yy}}
   + \left[ {c - \left( {{a_2} + b + 1} \right)z} \right]{u_z} + \left( {{a_2} - b + 1} \right)y{u_y} - {a_2}bu = 0, \hfill \\
\end{gathered}
\end{array}} \right.
$

where $u\equiv \,\,   {F_{38b}}\left( {{a_1},{a_2},{b};c;x,y,z} \right)$.

\bigskip

\begin{equation}
{F_{39a}}\left( {{a_1},{a_2},{a_3};{c_1},{c_2};x,y,z} \right)  = \sum\limits_{m,n,p = 0}^\infty  {} \frac{{{{\left( {{a_1}} \right)}_{2m + n}}{{\left( {{a_2}} \right)}_{n + p}}{{\left( {{a_3}} \right)}_p}}}{{{{\left( {{c_1}} \right)}_{n + p}}{{\left( {{c_2}} \right)}_m}}}\frac{x^m}{m!}\frac{y^n}{n!}\frac{z^p}{p!},
\end{equation}

first appearance of this function in the literature, and old notation:
[12], $X_{15}$,

region of convergence:
$$
 \left\{ {s + 2\sqrt r  < 1,\,\,\,t < 1} \right\}.
$$

System of partial differential equations:

$
\left\{ {\begin{array}{*{20}{l}}
  \begin{gathered}
  x\left( {1 - 4x} \right){u_{xx}} - 4xy{u_{xy}} - {y^2}{u_{yy}}\hfill\\\,\,\,\,\,\,\,\,\,
   + \left[ {{c_2} - 2\left( {2{a_1} + 3} \right)x} \right]{u_x} - 2\left( {{a_1} + 1} \right)y{u_y} - {a_1}\left( {1 + {a_1}} \right)u = 0, \hfill \\
\end{gathered}  \\
  \begin{gathered}
  y\left( {1 - y} \right){u_{yy}} - 2xy{u_{xy}} - 2xz{u_{xz}} + \left( {1 - y} \right)z{u_{yz}} \hfill \\
  \,\,\,\,\,\,\,\,\,\,
   + \left[ {{c_1} - \left( {{a_1} + {a_2} + 1} \right)y} \right]{u_y} - 2{a_2}x{u_x} - {a_1}z{u_z} - {a_1}{a_2}u = 0, \hfill \\
\end{gathered}  \\
  \begin{gathered}
  z\left( {1 - z} \right){u_{zz}} + y\left( {1 - z} \right){u_{yz}}
   + \left[ {{c_1} - \left( {{a_2} + {a_3} + 1} \right)z} \right]{u_z} - {a_3}y{u_y} - {a_2}{a_3}u = 0, \hfill \\
\end{gathered}
\end{array}} \right.
$

where $u\equiv \,\,   {F_{39a}}\left( {{a_1},{a_2},{a_3};{c_1},{c_2};x,y,z} \right)$.

Particular solutions:

$
{u_1} = {F_{39a}}\left( {{a_1},{a_2},{a_3};{c_1},{c_2};x,y,z} \right),
$

$
{u_2} = {x^{1 - {c_2}}}{F_{39a}}\left( {2 - 2{c_2} + {a_1},{a_2},{a_3};{c_1},2 - {c_2};x,y,z} \right).
$

\bigskip

\begin{equation}
{F_{39b}}\left( {{a_1},{a_2},b;c;x,y,z} \right)  = \sum\limits_{m,n,p = 0}^\infty  {} \frac{{{{\left( {{a_1}} \right)}_{2m + n}}{{\left( {{a_2}} \right)}_p}{{\left( b \right)}_{n + p - m}}}}{{{{\left( c \right)}_{n + p}}}}\frac{x^m}{m!}\frac{y^n}{n!}\frac{z^p}{p!},
\end{equation}

region of convergence:
$$
  \left\{ {t < 1,\,\,\,r \leq \frac{{1 + t}}{{{{\left( {2 + 3t} \right)}^2}}},\,\,\,s < \min \left\{ {{\Psi _1}\left( r \right),{\Psi _2}\left( r \right)} \right\}} \right\}
$$
$$
\cup  \left\{ t < 1,\,\,\, \frac{{1 + t}}{{{{\left( {2 + 3t} \right)}^2}}}<r<\frac{1}{4(1+t)},\,\,\,s < t \left(1-2\sqrt{r+rt}\right) \right\}.$$

System of partial differential equations:

$
\left\{ {\begin{array}{*{20}{l}}
  \begin{gathered}
  x\left( {1 + 4x} \right){u_{xx}} - \left( {1 - 4x} \right)y{u_{xy}} - z{u_{xz}} + {y^2}{u_{yy}} \hfill \\
\,\,\,\,\,\,\,\,\,   + \left[ {1 - b + 2\left( {2{a_1} + 3} \right)x} \right]{u_x} + 2\left( {{a_1} + 1} \right)y{u_y} + {a_1}\left( {1 + {a_1}} \right)u = 0, \hfill \\
\end{gathered}  \\
  \begin{gathered}
  y\left( {1 - y} \right){u_{yy}} - xy{u_{xy}} - 2xz{u_{xz}} + \left( {1 - y} \right)z{u_{yz}} + 2{x^2}{u_{xx}} \hfill \\
\,\,\,\,\,\,\,\,\,   + \left[ {c - \left( {{a_1} + b + 1} \right)y} \right]{u_y} + \left( {{a_1} - 2b + 2} \right)x{u_x} - {a_1}z{u_z} - {a_1}bu = 0, \hfill \\
\end{gathered}  \\
  \begin{gathered}
  z\left( {1 - z} \right){u_{zz}} + xz{u_{xz}} + y\left( {1 - z} \right){u_{yz}}
   + \left[ {c - \left( {{a_2} + b + 1} \right)z} \right]{u_z} + {a_2}x{u_x} - {a_2}y{u_y} - {a_2}bu = 0, \hfill \\
\end{gathered}
\end{array}} \right.
$

where $u\equiv \,\,   {F_{39b}}\left( {{a_1},{a_2},b;c;x,y,z} \right)$.

\bigskip

\begin{equation}
{F_{39c}}\left( {{a_1},{a_2},b;c;x,y,z} \right)  = \sum\limits_{m,n,p = 0}^\infty  {} \frac{{{{\left( {{a_1}} \right)}_{2m + n}}{{\left( {{a_2}} \right)}_{n + p}}{{\left( b \right)}_{p - m}}}}{{{{\left( c \right)}_{n + p}}}}\frac{x^m}{m!}\frac{y^n}{n!}\frac{z^p}{p!},
\end{equation}

region of convergence:
$$
\left\{ {s < 1,\,\,\,t < 1,\,\,\,r < \frac{{{{\left( {1 - s} \right)}^2}}}{{4\left( {1 + t} \right)}}} \right\}.
$$

System of partial differential equations:

$
\left\{ {\begin{array}{*{20}{l}}
  \begin{gathered}
  x\left( {1 + 4x} \right){u_{xx}} + 4xy{u_{xy}} - z{u_{xz}} + {y^2}{u_{yy}}\hfill\\\,\,\,\,\,\,\,\,\,
   + \left[ {1 - b + 2\left( {2{a_1} + 3} \right)x} \right]{u_x} + 2\left( {{a_1} + 1} \right)y{u_y} + {a_1}\left( {1 + {a_1}} \right)u = 0, \hfill \\
\end{gathered}  \\
  \begin{gathered}
  y\left( {1 - y} \right){u_{yy}} - 2xy{u_{xy}} - 2xz{u_{xz}} + \left( {1 - y} \right)z{u_{yz}}
   \hfill \\ \,\,\,\,\,\,\,\,\,+ \left[ {c - \left( {{a_1} + {a_2} + 1} \right)y} \right]{u_y} - 2{a_2}x{u_x} - {a_1}z{u_z} - {a_1}{a_2}u = 0, \hfill \\
\end{gathered}  \\
  \begin{gathered}
  z\left( {1 - z} \right){u_{zz}} + xy{u_{xy}} + xz{u_{xz}} + y\left( {1 - z} \right){u_{yz}}\hfill\\\,\,\,\,\,\,\,\,\,
   + \left[ {c - \left( {{a_2} + b + 1} \right)z} \right]{u_z} + {a_2}x{u_x} - by{u_y} - {a_2}bu = 0, \hfill \\
\end{gathered}
\end{array}} \right.
$

where $u\equiv \,\,   {F_{39c}}\left( {{a_1},{a_2},b;c;x,y,z} \right)$.

\bigskip

\begin{equation}
{F_{40a}}\left( {{a_1},{a_2},{a_3};c;x,y,z} \right)  = \sum\limits_{m,n,p = 0}^\infty  {} \frac{{{{\left( {{a_1}} \right)}_{2m + n}}{{\left( {{a_2}} \right)}_{n + p}}{{\left( {{a_3}} \right)}_p}}}{{{{\left( c \right)}_{m + n + p}}}}\frac{x^m}{m!}\frac{y^n}{n!}\frac{z^p}{p!},
\end{equation}

first appearance of this function in the literature, and old notation:
[12], $X_{13}$,

region of convergence:
$$
\left\{ {r < \frac{1}{4},\,\,\,t < 1,\,\,\,s < \frac{1}{2} + \frac{1}{2}\sqrt {1 - 4r} } \right\}.
$$

System of partial differential equations:

$
\left\{ {\begin{array}{*{20}{l}}
  \begin{gathered}
  x\left( {1 - 4x} \right){u_{xx}} + \left( {1 - 4x} \right)y{u_{xy}} + z{u_{xz}} - {y^2}{u_{yy}} \hfill \\
\,\,\,\,\,\,\,\,\,   + \left[ {c - 2\left( {2{a_1} + 3} \right)x} \right]{u_x} - 2\left( {{a_1} + 1} \right)y{u_y} - {a_1}\left( {1 + {a_1}} \right)u = 0, \hfill \\
\end{gathered}  \\
  \begin{gathered}
  y\left( {1 - y} \right){u_{yy}} + x\left( {1 - 2y} \right){u_{xy}} - 2xz{u_{xz}} + z\left( {1 - y} \right){u_{yz}} \hfill \\
 \,\,\,\,\,\,\,\,\,  + \left[ {c - \left( {{a_1} + {a_2} + 1} \right)y} \right]{u_y} - 2{a_2}x{u_x} - {a_1}z{u_z} - {a_1}{a_2}u = 0, \hfill \\
\end{gathered}  \\
  \begin{gathered}
  z\left( {1 - z} \right){u_{zz}} + x{u_{xz}} + y\left( {1 - z} \right){u_{yz}}
   + \left[ {c - \left( {{a_2} + {a_3} + 1} \right)z} \right]{u_z} - {a_3}y{u_y} - {a_2}{a_3}u = 0, \hfill \\
\end{gathered}
\end{array}} \right.
$

where $u\equiv \,\,   {F_{40a}}\left( {{a_1},{a_2},{a_3};c;x,y,z} \right)$.

\bigskip

\begin{equation}
{F_{41a}}\left( {{a_1},{a_2},{a_3};{c_1},{c_2},{c_3};x,y,z} \right)  = \sum\limits_{m,n,p = 0}^\infty  {} \frac{{{{\left( {{a_1}} \right)}_{2m + n + p}}{{\left( {{a_2}} \right)}_n}{{\left( {{a_3}} \right)}_p}}}{{{{\left( {{c_1}} \right)}_m}{{\left( {{c_2}} \right)}_n}{{\left( {{c_3}} \right)}_p}}}\frac{x^m}{m!}\frac{y^n}{n!}\frac{z^p}{p!},
\end{equation}

first appearance of this function in the literature, and old notation:
[11], ${}^{(1)}H_{4}^{(3)}$; [12], $X_{8}$,

region of convergence:
$$
\left\{ {2\sqrt{r} + s + t < 1} \right\}.
$$

System of partial differential equations:

$
\left\{ {\begin{array}{*{20}{l}}
  \begin{gathered}
  x\left( {1 - 4x} \right){u_{xx}} - 4xy{u_{xy}} - 4xz{u_{xz}} - 2yz{u_{yz}} - {y^2}{u_{yy}} - {z^2}{u_{zz}} \hfill \\
 \,\,\,\,\,\,\,\,\,  + \left[ {{c_1} - 2\left( {2{a_1} + 3} \right)x} \right]{u_x}  - 2\left( {{a_1} + 1} \right)y{u_y}
  - 2\left( {{a_1} + 1} \right)z{u_z} - {a_1}\left( {1 + {a_1}} \right)u = 0, \hfill \\
\end{gathered}  \\
  \begin{gathered}
  y\left( {1 - y} \right){u_{yy}} - 2xy{u_{xy}} - yz{u_{yz}}
   + \left[ {{c_2} - \left( {{a_1} + {a_2} + 1} \right)y} \right]{u_y} - 2{a_2}x{u_x} - {a_2}z{u_z} - {a_1}{a_2}u = 0, \hfill \\
\end{gathered}  \\
  \begin{gathered}
  z\left( {1 - z} \right){u_{zz}} - 2xz{u_{xz}} - yz{u_{yz}}
   + \left[ {{c_3} - \left( {{a_1} + {a_3} + 1} \right)z} \right]{u_z} - 2{a_3}x{u_x} - {a_3}y{u_y} - {a_1}{a_3}u = 0, \hfill \\
\end{gathered}
\end{array}} \right.
$

where $u\equiv \,\,   {F_{41a}}\left( {{a_1},{a_2},{a_3};{c_1},{c_2},{c_3};x,y,z} \right)$.

Particular solutions:

$
{u_1} = {F_{41a}}\left( {{a_1},{a_2},{a_3};{c_1},{c_2},{c_3};x,y,z} \right),
$

$
{u_2} = {x^{1 - {c_1}}}{F_{41a}}\left( {2 - 2{c_1} + {a_1},{a_2},{a_3};2 - {c_1},{c_2},{c_3};x,y,z} \right),
$

$
{u_3} = {y^{1 - {c_2}}}{F_{41a}}\left( {1 - {c_2} + {a_1},1 - {c_2} + {a_2},{a_3};{c_1},2 - {c_2},{c_3};x,y,z} \right),
$

$
{u_4} = {z^{1 - {c_3}}}{F_{41a}}\left( {1 - {c_3} + {a_1},{a_2},1 - {c_3} + {a_3};{c_1},{c_2},2 - {c_3};x,y,z} \right),
$

$
{u_5} = {x^{1 - {c_1}}}{y^{1 - {c_2}}}{F_{41a}}\left( {2 - {c_1} - {c_2} + {a_1},1 - {c_2} + {a_2},{a_3};2 - {c_1},2 - {c_2},{c_3};x,y,z} \right),
$

$
{u_6} = {y^{1 - {c_2}}}{z^{1 - {c_3}}}{F_{41a}}\left( {2 - {c_2} - {c_3} + {a_1},1 - {c_2} + {a_2},1 - {c_3} + {a_3};{c_1},2 - {c_2},2 - {c_3};x,y,z} \right),
$

$
{u_7} = {x^{1 - {c_1}}}{z^{1 - {c_3}}}{F_{41a}}\left( {3 - 2{c_1} - {c_3} + {a_1},{a_2},1 - {c_3} + {a_3};2 - {c_1},{c_2},2 - {c_3};x,y,z} \right),
$

$
  {u_8} = {x^{1 - {c_1}}}{y^{1 - {c_2}}}{z^{1 - {c_3}}} \times$
   
   $\,\,\,\,\,\,\,\,\,\times{F_{41a}}\left( {4 - 2{c_1} - {c_2} - {c_3} + {a_1},1 - {c_2} + {a_2},1 - {c_3} + {a_3};2 - {c_1},2 - {c_2},2 - {c_3};x,y,z} \right).
  $

\bigskip

\begin{equation}
{F_{41b}}\left( {{a_1},{a_2},b;{c_1},{c_2};x,y,z} \right) = \sum\limits_{m,n,p = 0}^\infty  {} \frac{{{{\left( {{a_1}} \right)}_{2m + n + p}}{{\left( {{a_2}} \right)}_p}{{\left( b \right)}_{n - m}}}}{{{{\left( {{c_1}} \right)}_n}{{\left( {{c_2}} \right)}_p}}}\frac{x^m}{m!}\frac{y^n}{n!}\frac{z^p}{p!},
\end{equation}

region of convergence:
$$
\left\{ {t + 2\sqrt r  < 1,\,\,\,s < \left( {1 - t} \right)\min \left\{ {{\Psi _1}\left( {\frac{r}{{{{\left( {1 - t} \right)}^2}}}} \right),{\Psi _2}\left( {\frac{r}{{{{\left( {1 - t} \right)}^2}}}} \right)} \right\}} \right\}.
$$

System of partial differential equations:

$
\left\{ {\begin{array}{*{20}{l}}
  \begin{gathered}
  x\left( {1 + 4x} \right){u_{xx}} - \left( {1 - 4x} \right)y{u_{xy}} + 4xz{u_{xz}} + 2yz{u_{yz}}
  + {y^2}{u_{yy}} + {z^2}{u_{zz}} \hfill \\  \,\,\,\,\,\,\,\,\,  + \left[ {1 - b + 2\left( {2{a_1} + 3} \right)x} \right]{u_x} + 2\left( {{a_1} + 1} \right)y{u_y}
   + 2\left( {{a_1} + 1} \right)z{u_z} + {a_1}\left( {1 + {a_1}} \right)u = 0 \hfill \\
\end{gathered}  \\
  \begin{gathered}
  y\left( {1 - y} \right){u_{yy}} - xy{u_{xy}} + xz{u_{xz}} - yz{u_{yz}} + 2{x^2}{u_{xx}} \hfill \\
  \,\,\,\,\,\,\,\,\, + \left[ {{c_1} - \left( {{a_1} + b + 1} \right)y} \right]{u_y} + \left( {{a_1} - 2b + 2} \right)x{u_x}
   - bz{u_z} - {a_1}bu = 0, \hfill \\
\end{gathered}  \\
  \begin{gathered}
  z\left( {1 - z} \right){u_{zz}} - 2xz{u_{xz}} - yz{u_{yz}}
   + \left[ {{c_2} - \left( {{a_1} + {a_2} + 1} \right)z} \right]{u_z} - 2{a_2}x{u_x} - {a_2}y{u_y} - {a_1}{a_2}u = 0, \hfill \\
\end{gathered}
\end{array}} \right.
$

where $u\equiv \,\,   {F_{41b}}\left( {{a_1},{a_2},b;{c_1},{c_2};x,y,z} \right)$.

Particular solutions:

$
{u_1} = {F_{41b}}\left( {{a_1},{a_2},b;{c_1},{c_2};x,y,z} \right),
$

$
{u_2} = {y^{1 - {c_1}}}{F_{41b}}\left( {1 - {c_1} + {a_1},{a_2},1 - {c_1} + b;2 - {c_1},{c_2};x,y,z} \right),
$

$
{u_3} = {z^{1 - {c_2}}}{F_{41b}}\left( {1 - {c_2} + {a_1},1 - {c_2} + {a_2},b;{c_1},2 - {c_2};x,y,z} \right),
$

$
{u_4} = {y^{1 - {c_1}}}{z^{1 - {c_2}}}{F_{41b}}\left( {2 - {c_1} - {c_2} + {a_1},1 - {c_2} + {a_2},1 - {c_1} + b;2 - {c_1},2 - {c_2};x,y,z} \right).
$

\bigskip

\begin{equation}
{F_{41c}}\left( {{a_1},{a_2},b;{c_1},{c_2};x,y,z} \right) = \sum\limits_{m,n,p = 0}^\infty  {} \frac{{{{\left( {{a_1}} \right)}_{2m + n + p}}{{\left( {{a_2}} \right)}_p}{{\left( b \right)}_{n - p}}}}{{{{\left( {{c_1}} \right)}_m}{{\left( {{c_2}} \right)}_n}}}\frac{x^m}{m!}\frac{y^n}{n!}\frac{z^p}{p!} ,
\end{equation}

region of convergence:
$$
\left\{ {s + 2\sqrt r  < 1,\,\,\,\sqrt t  < \sqrt {1 + s - 2\sqrt r }  - \sqrt s } \right\}.
$$

System of partial differential equations:

$
\left\{ {\begin{array}{*{20}{l}}
  \begin{gathered}
  x\left( {1 - 4x} \right){u_{xx}} - 4xy{u_{xy}} - 4xz{u_{xz}} - 2yz{u_{yz}} - {y^2}{u_{yy}} - {z^2}{u_{zz}} \hfill \\
\,\,\,\,\,\,\,\,\,   + \left[ {{c_1} - 2\left( {2{a_1} + 3} \right)x} \right]{u_x} - 2\left( {{a_1} + 1} \right)y{u_y}
    - 2\left( {{a_1} + 1} \right)z{u_z} - {a_1}\left( {1 + {a_1}} \right)u = 0, \hfill \\
\end{gathered}  \\
  \begin{gathered}
  y\left( {1 - y} \right){u_{yy}} - 2xy{u_{xy}} + 2xz{u_{xz}} + {z^2}{u_{zz}} \hfill \\
  \,\,\,\,\,\,\,\,\, + \left[ {{c_2} - \left( {{a_1} + b + 1} \right)y} \right]{u_y} - 2bx{u_x}
    + \left( {{a_1} - b + 1} \right)z{u_z} - {a_1}bu = 0, \hfill \\
\end{gathered}  \\
  \begin{gathered}
  z\left( {1 + z} \right){u_{zz}} + 2xz{u_{xz}} - y\left( {1 - z} \right){u_{yz}}\hfill\\\,\,\,\,\,\,\,\,\,
   + \left[ {1 - b + \left( {{a_1} + {a_2} + 1} \right)z} \right]{u_z} + 2{a_2}x{u_x} + {a_2}y{u_y} + {a_1}{a_2}u = 0, \hfill \\
\end{gathered}
\end{array}} \right.
$

where $u\equiv \,\,   {F_{41c}}\left( {{a_1},{a_2},b;{c_1},{c_2};x,y,z} \right)$.

Particular solutions:

$
{u_1} = {F_{41c}}\left( {{a_1},{a_2},b;{c_1},{c_2};x,y,z} \right),
$

$
{u_2} = {x^{1 - {c_1}}}{F_{41c}}\left( {2 - 2{c_1} + {a_1},{a_2},b;2 - {c_1},{c_2};x,y,z} \right),
$

$
{u_3} = {y^{1 - {c_2}}}{F_{41c}}\left( {1 - {c_2} + {a_1},{a_2},1 - {c_2} + b;{c_1},2 - {c_2};x,y,z} \right),
$

$
{u_4} = {x^{1 - {c_1}}}{y^{1 - {c_2}}}{F_{41c}}\left( {2 - 2{c_1} - {c_2} + {a_1},{a_2},1 - {c_2} + b;2 - {c_1},2 - {c_2};x,y,z} \right).
$

\bigskip

\begin{equation}
{F_{41d}}\left( {a,{b_1},{b_2};c;x,y,z} \right) = \sum\limits_{m,n,p = 0}^\infty  {} \frac{{{{\left( a \right)}_{2m + n + p}}{{\left( {{b_1}} \right)}_{n - m}}{{\left( {{b_2}} \right)}_{p - n}}}}{{{{\left( c \right)}_p}}}\frac{x^m}{m!}\frac{y^n}{n!}\frac{z^p}{p!},
\end{equation}

region of convergence:
$$
\begin{gathered}
  \left\{ {r < \frac{1}{4},\,\,\,s < \min \left\{ {{\Psi _1}\left( r \right),{\Psi _2}\left( r \right)} \right\},\,\,\,\sqrt t  < \min \left\{ {{U^ + }\left( {{w_1}} \right),{U^ + }\left( {{w_2}} \right),{U^ - }\left( {{w_3}} \right)} \right\}} \right\}, \hfill \\
  {P_{rs}}\left( w \right) = 4{w^4} - {w^3} - 6r{w^2} + \left( {1 - s} \right)rw + 2{r^2}, \hfill \\
{w_1}:\,\,{\rm{the\,\,root\,\,in}}\,\,\left( {\sqrt r ,\frac{1}{6}\left[ {1 + a\left( r \right)} \right]} \right)\,\,{\rm{of}}\,\,{P_{rs}}\left( w \right) = 0, \hfill \\
{w_2}:\,\,{\rm{the\,\,smaller\,\,root\,\,in}}\,\,\left( {\sqrt r ,\frac{1}{6}\left[ {1 + a\left( r \right)} \right]} \right)\,\,{\rm{of}}\,\,{P_{r, - s}}\left( w \right) = 0, \hfill \\
{w_3}:\,\,{\rm{the\,\,smaller\,\,root\,\,in}}\,\,\left( {0,\frac{1}{2}} \right)\,\,{\rm{of}}\,\,{P_{ - r, - s}}\left( w \right) = 0, \hfill \\
  {U^ \pm }\left( w \right) = \frac{{\left( { - 3{w^2} + w \pm r} \right)\sqrt {{w^2} \mp r} }}{{w\sqrt {rs} }}. \hfill \\
\end{gathered}
$$

System of partial differential equations:

$
\left\{ {\begin{array}{*{20}{l}}
  \begin{gathered}
  x\left( {1 + 4x} \right){u_{xx}} - \left( {1 - 4x} \right)y{u_{xy}} + 4xz{u_{xz}} + 2yz{u_{yz}}
   + {y^2}{u_{yy}} + {z^2}{u_{zz}} \hfill \\ \,\,\,\,\,\,\,\,\,+ \left[ {1 - {b_1} + 2\left( {2a + 3} \right)x} \right]{u_x} + 2\left( {a + 1} \right)y{u_y}
    + 2\left( {a + 1} \right)z{u_z} + a\left( {1 + a} \right)u = 0, \hfill \\
\end{gathered}  \\
  \begin{gathered}
  y\left( {1 + y} \right){u_{yy}} + xy{u_{xy}} - xz{u_{xz}} - \left( {1 - y} \right)z{u_{yz}} - 2{x^2}{u_{xx}} \hfill \\
 \,\,\,\,\,\,\,\,\,  + \left[ {1 - {b_2} + \left( {a + {b_1} + 1} \right)y} \right]{u_y} - \left( {a - 2{b_1} + 2} \right)x{u_x}
   + {b_1}z{u_z} + a{b_1}u = 0, \hfill \\
\end{gathered}  \\
  \begin{gathered}
  z\left( {1 - z} \right){u_{zz}} + 2xy{u_{xy}} - 2xz{u_{xz}} + {y^2}{u_{yy}} \hfill \\ \,\,\,\,\,\,\,\,\,
   + \left[ {c - \left( {a + {b_2} + 1} \right)z} \right]{u_z} - 2{b_2}x{u_x}
   + \left( {a - {b_2} + 1} \right)y{u_y} - a{b_2}u = 0, \hfill \\
\end{gathered}
\end{array}} \right.
$

where $u\equiv \,\,   {F_{41d}}\left( {a,{b_1},{b_2};c;x,y,z} \right)$.

Particular solutions:

$
{u_1} = {F_{41d}}\left( {a,{b_1},{b_2};c;x,y,z} \right),
$

$
{u_2} = {z^{1 - c}}{F_{41d}}\left( {1 - c + a,{b_1},1 - c + {b_2};c;x,y,z} \right).
$

\bigskip

\begin{equation}
{F_{41e}}\left( {a,{b_1},{b_2};c;x,y,z} \right) = \sum\limits_{m,n,p = 0}^\infty  {} \frac{{{{\left( a \right)}_{2m + n + p}}{{\left( {{b_1}} \right)}_{n - p}}{{\left( {{b_2}} \right)}_{p - n}}}}{{{{\left( c \right)}_m}}}\frac{x^m}{m!}\frac{y^n}{n!}\frac{z^p}{p!},
\end{equation}

region of convergence:
$$
\left\{ {2\sqrt r  + s + t < 1} \right\}.
$$

System of partial differential equations:

$
\left\{ {\begin{array}{*{20}{l}}
  \begin{gathered}
  x\left( {1 - 4x} \right){u_{xx}} - 4xy{u_{xy}} - 4xz{u_{xz}} - 2yz{u_{yz}} - {y^2}{u_{yy}}\hfill \\
  \,\,\,\,\,\,\,\,\, - {z^2}{u_{zz}} + \left[ {c - 2\left( {2a + 3} \right)x} \right]{u_x} - 2\left( {a + 1} \right)y{u_y}
   - 2\left( {a + 1} \right)z{u_z} - a\left( {1 + a} \right)u = 0, \hfill \\
\end{gathered}  \\
  \begin{gathered}
  y\left( {1 + y} \right){u_{yy}} + 2xy{u_{xy}} - 2xz{u_{xz}} - z{u_{yz}} - {z^2}{u_{zz}} \hfill \\
 \,\,\,\,\,\,\,\,\, + \left[ {1 - {b_2} + \left( {a + {b_1} + 1} \right)y} \right]{u_y} + 2{b_1}x{u_x}
   - \left( {a - {b_1} + 1} \right)z{u_z} + a{b_1}u = 0, \hfill \\
\end{gathered}  \\
  \begin{gathered}
  z\left( {1 + z} \right){u_{zz}} - 2xy{u_{xy}} + 2xz{u_{xz}} - y{u_{yz}} - {y^2}{u_{yy}} \hfill \\
 \,\,\,\,\,\,\,\,\, + \left[ {1 - {b_1} + \left( {a + {b_2} + 1} \right)z} \right]{u_z} + 2{b_2}x{u_x}
   - \left( {a - {b_2} + 1} \right)y{u_y} + a{b_2}u = 0, \hfill \\
\end{gathered}
\end{array}} \right.
$

where $u\equiv \,\,   {F_{41e}}\left( {a,{b_1},{b_2};c;x,y,z} \right)$.

Particular solutions:

$
{u_1} = {F_{41e}}\left( {a,{b_1},{b_2};c;x,y,z} \right),
$

$
{u_2} = {x^{1 - c}}{F_{41e}}\left( {2 - 2c + a,{b_1},{b_2};2 - c;x,y,z} \right).
$

\bigskip

\begin{equation}
{F_{42a}}\left( {{a_1},{a_2},{a_3};{c_1},{c_2};x,y,z} \right) = \sum\limits_{m,n,p = 0}^\infty  {} \frac{{{{\left( {{a_1}} \right)}_{2m + n + p}}{{\left( {{a_2}} \right)}_n}{{\left( {{a_3}} \right)}_p}}}{{{{\left( {{c_1}} \right)}_{m + n}}{{\left( {{c_2}} \right)}_p}}}\frac{x^m}{m!}\frac{y^n}{n!}\frac{z^p}{p!},
\end{equation}

first appearance of this function in the literature, and old notation:
[12], $X_{6}$,

region of convergence:
$$
\left\{ {t + 2\sqrt r  < 1,\,\,\,s < \frac{1}{2}\left( {1 - t} \right) + \frac{1}{2}\sqrt {{{\left( {1 - t} \right)}^2} - 4r} } \right\}.
$$

System of partial differential equations:

$
\left\{ {\begin{array}{*{20}{l}}
  \begin{gathered}
  x\left( {1 - 4x} \right){u_{xx}} + \left( {1 - 4x} \right)y{u_{xy}} - 4xz{u_{xz}} - 2yz{u_{yz}} - {y^2}{u_{yy}}
   - {z^2}{u_{zz}} \hfill \\ \,\,\,\,\,\,\,\,\, + \left[ {{c_1} - 2\left( {2{a_1} + 3} \right)x} \right]{u_x}
   - 2\left( {{a_1} + 1} \right)y{u_y} - 2\left( {{a_1} + 1} \right)z{u_z} - {a_1}\left( {1 + {a_1}} \right)u = 0, \hfill \\
\end{gathered}  \\
  \begin{gathered}
  y\left( {1 - y} \right){u_{yy}} + x\left( {1 - 2y} \right){u_{xy}} - yz{u_{yz}}\hfill\\\,\,\,\,\,\,\,\,\,
   + \left[ {{c_1} - \left( {{a_1} + {a_2} + 1} \right)y} \right]{u_y} - 2{a_2}x{u_x} - {a_2}z{u_z} - {a_1}{a_2}u = 0, \hfill \\
\end{gathered}  \\
  \begin{gathered}
  z\left( {1 - z} \right){u_{zz}} - 2xz{u_{xz}} - yz{u_{yz}}\hfill\\\,\,\,\,\,\,\,\,\,
   + \left[ {{c_2} - \left( {{a_1} + {a_3} + 1} \right)z} \right]{u_z} - 2{a_3}x{u_x} - {a_3}y{u_y} - {a_1}{a_3}u = 0, \hfill \\
\end{gathered}
\end{array}} \right.
$

where $u\equiv \,\,   {F_{42a}}\left( {{a_1},{a_2},{a_3};{c_1},{c_2};x,y,z} \right)$.

Particular solutions:

$
{u_1} = {F_{42a}}\left( {{a_1},{a_2},{a_3};{c_1},{c_2};x,y,z} \right),
$

$
{u_2} = {z^{1 - {c_2}}}{F_{42a}}\left( {1 - {c_2} + {a_1},{a_2},1 - {c_2} + {a_3};{c_1},2 - {c_2};x,y,z} \right).
$

\bigskip

\begin{equation}
{F_{42b}}\left( {{a_1},{a_2},b;c;x,y,z} \right) = \sum\limits_{m,n,p = 0}^\infty  {} \frac{{{{\left( {{a_1}} \right)}_{2m + n + p}}{{\left( {{a_2}} \right)}_p}{{\left( b \right)}_{n - p}}}}{{{{\left( c \right)}_{m + n}}}}\frac{x^m}{m!}\frac{y^n}{n!}\frac{z^p}{p!},
\end{equation}

region of convergence:
$$
 \left\{ {t + 2\sqrt r  < 1,\,\,\,s < \min \left\{ {\frac{1}{2} + \frac{1}{2}\sqrt {1 - 4r} ,\,\,\,\frac{{{{\left( {1 - t} \right)}^2} - 4r}}{{4t}}} \right\}} \right\}.
$$

System of partial differential equations:

$
\left\{ {\begin{array}{*{20}{l}}
  \begin{gathered}
  x\left( {1 - 4x} \right){u_{xx}} + \left( {1 - 4x} \right)y{u_{xy}} - 4xz{u_{xz}} - 2yz{u_{yz}} - {y^2}{u_{yy}}
   - {z^2}{u_{zz}}  \hfill \\ \,\,\,\,\,\,\,\,\, + \left[ {c - 2\left( {2{a_1} + 3} \right)x} \right]{u_x}
   - 2\left( {{a_1} + 1} \right)y{u_y} - 2\left( {{a_1} + 1} \right)z{u_z} - {a_1}\left( {1 + {a_1}} \right)u = 0, \hfill \\
\end{gathered}  \\
  \begin{gathered}
  y\left( {1 - y} \right){u_{yy}} + x\left( {1 - 2y} \right){u_{xy}} + 2xz{u_{xz}} + {z^2}{u_{zz}} \hfill \\
 \,\,\,\,\,\,\,\,\,  + \left[ {c - \left( {{a_1} + b + 1} \right)y} \right]{u_y} - 2bx{u_x}
    + \left( {{a_1} - b + 1} \right)z{u_z} - {a_1}bu = 0, \hfill \\
\end{gathered}  \\
  \begin{gathered}
  z\left( {1 + z} \right){u_{zz}} + 2xz{u_{xz}} - y\left( {1 - z} \right){u_{yz}}\hfill\\\,\,\,\,\,\,\,\,\,
   + \left[ {1 - b + \left( {{a_2} + {a_1} + 1} \right)z} \right]{u_z} + 2{a_2}x{u_x} + {a_2}y{u_y} + {a_1}{a_2}u = 0, \hfill \\
\end{gathered}
\end{array}} \right.
$

where $u\equiv \,\,   {F_{42b}}\left( {{a_1},{a_2},b;c;x,y,z} \right)$.

\bigskip

\begin{equation}
{F_{42c}}\left( {{a_1},{a_2},b;c;x,y,z} \right) = \sum\limits_{m,n,p = 0}^\infty  {} \frac{{{{\left( {{a_1}} \right)}_{2m + n + p}}{{\left( {{a_2}} \right)}_n}{{\left( b \right)}_{p - m - n}}}}{{{{\left( c \right)}_p}}}\frac{x^m}{m!}\frac{y^n}{n!}\frac{z^p}{p!},
\end{equation}

region of convergence:
$$
\begin{gathered}
    \left\{ {t < 1,\,\,\,\sqrt s  \leq \sqrt {9t + 8}  - \frac{3}{4}\sqrt t ,\,\,\,r < \min \left\{ {{\Theta _1}\left( t \right),{\Theta _2}\left( t \right)} \right\}} \right\} \hfill \\
   \cup \left\{ {t < 1,\,\,\,\frac{1}{4}\sqrt {9t + 8}  - \frac{3}{4}\sqrt t  < \sqrt s  < \sqrt {1 + t}  - \sqrt t ,} \,\,\, {r < \min \left\{ {{\Theta _2}\left( t \right),s\left( {1 - s - 2\sqrt {st} } \right)} \right\}} \right\}. \hfill \\
\end{gathered}
$$

System of partial differential equations:

$
\left\{ {\begin{array}{*{20}{l}}
  \begin{gathered}
  x\left( {1 + 4x} \right){u_{xx}} + \left( {1 + 4x} \right)y{u_{xy}} - \left( {1 - 4x} \right)z{u_{xz}} + 2yz{u_{yz}}
  + {y^2}{u_{yy}} + {z^2}{u_{zz}} \hfill \\ \,\,\,\,\,\,\,\,\, + \left[ {1 - b + 2\left( {2{a_1} + 3} \right)x} \right]{u_x} + 2\left( {{a_1} + 1} \right)y{u_y}
  + 2\left( {{a_1} + 1} \right)z{u_z} + {a_1}\left( {1 + {a_1}} \right)u = 0, \hfill \\
\end{gathered}  \\
  \begin{gathered}
  y\left( {1 + y} \right){u_{yy}} + x\left( {1 + 2y} \right){u_{xy}} - \left( {1 - y} \right)z{u_{yz}} \hfill \\
 \,\,\,\,\,\,\,\,\, + \left[ {1 - b + \left( {{a_1} + {a_2} + 1} \right)y} \right]{u_y} + 2{a_2}x{u_x} + {a_2}z{u_z} + {a_1}{a_2}u = 0, \hfill \\
\end{gathered}  \\
  \begin{gathered}
  z\left( {1 - z} \right){u_{zz}} + 3xy{u_{xy}} - xz{u_{xz}} + {y^2}{u_{yy}} + 2{x^2}{u_{xx}} \hfill \\
 \,\,\,\,\,\,\,\,\,  + \left[ {c - \left( {{a_1} + b + 1} \right)z} \right]{u_z} + \left( {{a_1} - 2b + 2} \right)x{u_x}
   + \left( {{a_1} - b + 1} \right)y{u_y} - {a_1}bu = 0, \hfill \\
\end{gathered}
\end{array}} \right.
$

where $u\equiv \,\,   {F_{42c}}\left( {{a_1},{a_2},b;c;x,y,z} \right)$.

Particular solutions:

$
{u_1} = {F_{42c}}\left( {{a_1},{a_2},b;c;x,y,z} \right),
$

$
{u_2} = {z^{1 - c}}{F_{42c}}\left( {1 - c + {a_1},{a_2},1 - c + b;2 - c;x,y,z} \right).
$

\bigskip

\begin{equation}
{F_{42d}}\left( {a,{b_1},{b_2};x,y,z} \right) = \sum\limits_{m,n,p = 0}^\infty  {} {{{{\left( a \right)}_{2m + n + p}}{{\left( {{b_1}} \right)}_{p - m - n}}{{\left( {{b_2}} \right)}_{n - p}}}}\frac{x^m}{m!}\frac{y^n}{n!}\frac{z^p}{p!},
\end{equation}

region of convergence:
$$
\begin{gathered}
  \left\{ {r < \frac{1}{4},\,\,\,s < \frac{1}{2} + \frac{1}{2}\sqrt {1 - 4r} ,\,\,\,t < \min \left[ {{U^ + }\left( {{w_1}} \right),{U^ - }\left( {{w_2}} \right), - {U^ + }\left( {{w_3}} \right)} \right]} \right\}, \hfill \\
  {P_{rs}}\left( w \right) = 2s{w^3} + 3r{w^2} - \left( {1 + 2s} \right)rw - {r^2}, \hfill \\
{w_1}:\,\,{\rm{the\,\,root\,\,in}}\,\,\left( {\sqrt r ,\frac{1}{6}\left[ {1 + a\left( r \right)} \right]} \right)\,\,{\rm{of}}\,\,{P_{rs}}\left( w \right) = 0, \hfill \\
{w_2}:\,\,{\rm{the\,\,smaller\,\,root\,\,in}}\,\,\left( {0,\frac{1}{2}} \right)\,\,{\rm{of}}\,\,{P_{ - r, - s}}\left( w \right) = 0, \hfill \\
{w_3}:\,\,{\rm{the\,\,smaller\,\,root\,\,in}}\,\,\left( {\frac{1}{2} - \frac{1}{2}\sqrt {1 - 4r} ,\sqrt r } \right)\,\,{\rm{of}}\,\,{P_{r, - s}}\left( w \right) = 0, \hfill \\
  {U^ \pm }\left( w \right) = \frac{{\left( {{w^2} \mp r} \right)\left( { - {w^2} + w \mp r} \right)}}{{2rw}}. \hfill \\
\end{gathered}
$$

System of partial differential equations:

$
\left\{ {\begin{array}{*{20}{l}}
  \begin{gathered}
  x\left( {1 + 4x} \right){u_{xx}} + \left( {1 + 4x} \right)y{u_{xy}} - \left( {1 - 4x} \right)z{u_{xz}} + 2yz{u_{yz}}
  + {y^2}{u_{yy}} + {z^2}{u_{zz}}  \hfill \\\,\,\,\,\,\,\,\,\,
+ \left[ {1 - {b_1} + 2\left( {2a + 3} \right)x} \right]{u_x} + 2\left( {a + 1} \right)y{u_y}
  + 2\left( {a + 1} \right)z{u_z} + a\left( {1 + a} \right)u = 0, \hfill \\
\end{gathered}  \\
  \begin{gathered}
  y\left( {1 + y} \right){u_{yy}} + x\left( {1 + 2y} \right){u_{xy}} - 2xz{u_{xz}} - z{u_{yz}} - {z^2}{u_{zz}} \hfill \\
  \,\,\,\,\,\,\,\,\, + \left[ {1 - {b_1} + \left( {a + {b_2} + 1} \right)y} \right]{u_y} + 2{b_2}x{u_x}
   - \left( {a - {b_2} + 1} \right)z{u_z} + a{b_2}u = 0, \hfill \\
\end{gathered}  \\
  \begin{gathered}
  z\left( {1 + z} \right){u_{zz}} - 3xy{u_{xy}} + xz{u_{xz}} - y{u_{yz}} - 2{x^2}{u_{xx}} - {y^2}{u_{yy}} \hfill \\
 \,\,\,\,\,\,\,\,\,  + \left[ {1 - {b_2} + \left( {a + {b_1} + 1} \right)z} \right]{u_z}
    - \left( {a - 2{b_1} + 2} \right)x{u_x} - \left( {a - {b_1} + 1} \right)y{u_y} + a{b_1}u = 0, \hfill \\
\end{gathered}
\end{array}} \right.
$

where $u\equiv \,\,   {F_{42d}}\left( {a,{b_1},{b_2};x,y,z} \right)$.

\bigskip

\begin{equation}
{F_{43a}}\left( {{a_1},{a_2},{a_3};{c_1},{c_2};x,y,z} \right) = \sum\limits_{m,n,p = 0}^\infty  {} \frac{{{{\left( {{a_1}} \right)}_{2m + n + p}}{{\left( {{a_2}} \right)}_n}{{\left( {{a_3}} \right)}_p}}}{{{{\left( {{c_1}} \right)}_{n + p}}{{\left( {{c_2}} \right)}_m}}}\frac{x^m}{m!}\frac{y^n}{n!}\frac{z^p}{p!},
\end{equation}

first appearance of this function in the literature, and old notation:
[12], $X_{7}$,

region of convergence:
$$
 \left\{ {s < 1,\,\,\,t < 1,\,\,\,r < \frac{1}{4}\min \left\{ {{{\left( {1 - s} \right)}^2},{{\left( {1 - t} \right)}^2}} \right\}} \right\}.
$$

System of partial differential equations:

$
\left\{ {\begin{array}{*{20}{l}}
  \begin{gathered}
  x\left( {1 - 4x} \right){u_{xx}} - 4xy{u_{xy}} - 4xz{u_{xz}} - 2yz{u_{yz}}
   - {y^2}{u_{yy}} - {z^2}{u_{zz}}  \hfill \\ \,\,\,\,\,\,\,\,\,+ \left[ {{c_2} - 2\left( {2{a_1} + 3} \right)x} \right]{u_x}
    - 2\left( {{a_1} + 1} \right)y{u_y} - 2\left( {{a_1} + 1} \right)z{u_z} - {a_1}\left( {1 + {a_1}} \right)u = 0, \hfill \\
\end{gathered}  \\
  \begin{gathered}
  y\left( {1 - y} \right){u_{yy}} - 2xy{u_{xy}} + \left( {1 - y} \right)z{u_{yz}}\hfill\\\,\,\,\,\,\,\,\,\,
   + \left[ {{c_1} - \left( {{a_1} + {a_2} + 1} \right)y} \right]{u_y} - 2{a_2}x{u_x} - {a_2}z{u_z} - {a_1}{a_2}u = 0, \hfill \\
\end{gathered}  \\
  \begin{gathered}
  z\left( {1 - z} \right){u_{zz}} - 2xz{u_{xz}} + y\left( {1 - z} \right){u_{yz}}\hfill\\\,\,\,\,\,\,\,\,\,
   + \left[ {{c_1} - \left( {{a_1} + {a_3} + 1} \right)z} \right]{u_z} - 2{a_3}x{u_x} - {a_3}y{u_y} - {a_1}{a_3}u = 0, \hfill \\
\end{gathered}
\end{array}} \right.
$

where $u\equiv \,\,   {F_{43a}}\left( {{a_1},{a_2},{a_3};{c_1},{c_2};x,y,z} \right)$.

Particular solutions:

$
{u_1} = {F_{43a}}\left( {{a_1},{a_2},{a_3};{c_1},{c_2};x,y,z} \right),
$

$
{u_2} = {x^{1 - {c_2}}}{F_{43a}}\left( {2 - 2{c_2} + {a_1},{a_2},{a_3};{c_1},2 - {c_2};x,y,z} \right).
$

\bigskip

\begin{equation}
{F_{43b}}\left( {{a_1},{a_2},b;c;x,y,z} \right) = \sum\limits_{m,n,p = 0}^\infty  {} \frac{{{{\left( {{a_1}} \right)}_{2m + n + p}}{{\left( {{a_2}} \right)}_p}{{\left( b \right)}_{n - m}}}}{{{{\left( c \right)}_{n + p}}}}\frac{x^m}{m!}\frac{y^n}{n!}\frac{z^p}{p!},
\end{equation}

region of convergence:
$$
  \left\{ {r < \frac{1}{4},\,\,\,t \leq \frac{1}{3}\left[ {2 - a\left( r \right)} \right],\,\,\,s < \min \left\{ {{\Psi _1}\left( r \right),{\Psi _2}\left( r \right)} \right\}} \right\}
$$
$$
\cup \left\{ t < \frac{1}{4},\,\,\,\frac{1}{3}[2-a(r)]< t  < 1-2\sqrt{r}, \,\,\,  s < \min \left\{ {\Psi_2}(r), \,\,\, \frac{t}{4r}\left[(1-t)^2-4r\right]\right\}\right\}.
$$

System of partial differential equations:

$
\left\{ {\begin{array}{*{20}{l}}
  \begin{gathered}
  x\left( {1 + 4x} \right){u_{xx}} - \left( {1 - 4x} \right)y{u_{xy}} + 4xz{u_{xz}} + 2yz{u_{yz}}
   + {y^2}{u_{yy}} + {z^2}{u_{zz}} \hfill \\  \,\,\,\,\,\,\,\,\, + \left[ {1 - b + 2\left( {2{a_1} + 3} \right)x} \right]{u_x} + 2\left( {{a_1} + 1} \right)y{u_y}
    + 2\left( {{a_1} + 1} \right)z{u_z} + {a_1}\left( {1 + {a_1}} \right)u = 0, \hfill \\
\end{gathered}  \\
  \begin{gathered}
  y\left( {1 - y} \right){u_{yy}} - xy{u_{xy}} + xz{u_{xz}} + \left( {1 - y} \right)z{u_{yz}} + 2{x^2}{u_{xx}} \hfill \\
 \,\,\,\,\,\,\,\,\,  + \left[ {c - \left( {{a_1} + b + 1} \right)y} \right]{u_y} + \left( {{a_1} - 2b + 2} \right)x{u_x}
    - bz{u_z} - {a_1}bu = 0, \hfill \\
\end{gathered}  \\
  \begin{gathered}
  z\left( {1 - z} \right){u_{zz}} - 2xz{u_{xz}} + y\left( {1 - z} \right){u_{yz}}\hfill\\\,\,\,\,\,\,\,\,\,
   + \left[ {c - \left( {{a_1} + {a_2} + 1} \right)z} \right]{u_z} - 2{a_2}x{u_x} - {a_2}y{u_y} - {a_1}{a_2}u = 0, \hfill \\
\end{gathered}
\end{array}} \right.
$

where $u\equiv \,\,   {F_{43b}}\left( {{a_1},{a_2},b;c;x,y,z} \right)$.

\bigskip

\begin{equation}
{F_{44a}}\left( {{a_1},{a_2},{a_3};c;x,y,z} \right) = \sum\limits_{m,n,p = 0}^\infty  {} \frac{{{{\left( {{a_1}} \right)}_{2m + n + p}}{{\left( {{a_2}} \right)}_n}{{\left( {{a_3}} \right)}_p}}}{{{{\left( c \right)}_{m + n + p}}}}\frac{x^m}{m!}\frac{y^n}{n!}\frac{z^p}{p!},
\end{equation}

first appearance of this function in the literature, and old notation:
[11], ${}^{(1)}H_{3}^{(3)}$; [12], $X_{5}$,

region of convergence:
$$
\left\{ {r < \frac{1}{4},\,\,\,\max \left( {s,t} \right) < \frac{1}{2} + \frac{1}{2}\sqrt {1 - 4r} } \right\}.
$$

System of partial differential equations:

$
\left\{ {\begin{array}{*{20}{l}}
  \begin{gathered}
  x\left( {1 - 4x} \right){u_{xx}} + \left( {1 - 4x} \right)y{u_{xy}} + z\left( {1 - 4x} \right){u_{xz}} - 2yz{u_{yz}}
  - {y^2}{u_{yy}} - {z^2}{u_{zz}}
  \hfill \\ \,\,\,\,\,\,\,\,\,
  + \left[ {c - 2\left( {2{a_1} + 3} \right)x} \right]{u_x}
  - 2\left( {{a_1} + 1} \right)y{u_y} - 2\left( {{a_1} + 1} \right)z{u_z} - {a_1}\left( {1 + {a_1}} \right)u = 0, \hfill \\
\end{gathered}  \\
  \begin{gathered}
  y\left( {1 - y} \right){u_{yy}} + x\left( {1 - 2y} \right){u_{xy}} + \left( {1 - y} \right)z{u_{yz}}  \hfill\\ \,\,\,\,\,\,\,\,\,+ \left[ {c - \left( {{a_1} + {a_2} + 1} \right)y} \right]{u_y} - 2{a_2}x{u_x} - {a_2}z{u_z} - {a_1}{a_2}u = 0, \hfill \\
\end{gathered}  \\
  \begin{gathered}
  z\left( {1 - z} \right){u_{zz}} + x\left( {1 - 2z} \right){u_{xz}} + y\left( {1 - z} \right){u_{yz}} \hfill \\ \,\,\,\,\,\,\,\,\, + \left[ {c - \left( {{a_1} + {a_3} + 1} \right)z} \right]{u_z} - 2{a_3}x{u_x} - {a_3}y{u_y} - {a_1}{a_3}u = 0, \hfill \\
\end{gathered}
\end{array}} \right.
$

where $u\equiv \,\,   {F_{44a}}\left( {{a_1},{a_2},{a_3};c;x,y,z} \right)$.

\bigskip

\begin{equation}
{F_{45a}}\left( {{a_1},{a_2};{c_1},{c_2},{c_3};x,y,z} \right) = \sum\limits_{m,n,p = 0}^\infty  {} \frac{{{{\left( {{a_1}} \right)}_{2m + n + p}}{{\left( {{a_2}} \right)}_{n + p}}}}{{{{\left( {{c_1}} \right)}_m}{{\left( {{c_2}} \right)}_n}{{\left( {{c_3}} \right)}_p}}}\frac{x^m}{m!}\frac{y^n}{n!}\frac{z^p}{p!},
\end{equation}

first appearance of this function in the literature, and old notation:
[12], $X_{4}$,

region of convergence:
$$
\left\{ {2\sqrt r  + {{\left( {\sqrt s  + \sqrt t } \right)}^2} < 1} \right\}.
$$

System of partial differential equations:

$
\left\{ {\begin{array}{*{20}{l}}
  \begin{gathered}
  x\left( {1 - 4x} \right){u_{xx}} - 4xy{u_{xy}} - 4xz{u_{xz}} - 2yz{u_{yz}}
   - {y^2}{u_{yy}} - {z^2}{u_{zz}} \hfill \\ \,\,\,\,\,\,\,\,\,+ \left[ {{c_1} - 2\left( {2{a_1} + 3} \right)x} \right]{u_x}
   - 2\left( {{a_1} + 1} \right)y{u_y} - 2\left( {{a_1} + 1} \right)z{u_z} - {a_1}\left( {1 + {a_1}} \right)u = 0, \hfill \\
\end{gathered}  \\
  \begin{gathered}
  y\left( {1 - y} \right){u_{yy}} - 2xy{u_{xy}} - 2xz{u_{xz}} - 2yz{u_{yz}} - {z^2}{u_{zz}} \hfill \\\,\,\,\,\,\,\,\,\,
   + \left[ {{c_2} - \left( {{a_1} + {a_2} + 1} \right)y} \right]{u_y} - 2{a_2}x{u_x}
   - \left( {{a_1} + {a_2} + 1} \right)z{u_z} - {a_1}{a_2}u = 0, \hfill \\
\end{gathered}  \\
  \begin{gathered}
  z\left( {1 - z} \right){u_{zz}} - 2xy{u_{xy}} - 2xz{u_{xz}} - 2yz{u_{yz}} - {y^2}{u_{yy}} \hfill \\\,\,\,\,\,\,\,\,\,
   + \left[ {{c_3} - \left( {{a_1} + {a_2} + 1} \right)z} \right]{u_z} - 2{a_2}x{u_x}
   - \left( {{a_1} + {a_2} + 1} \right)y{u_y} - {a_1}{a_2}u = 0, \hfill \\
\end{gathered}
\end{array}} \right.
$

where $u\equiv \,\,   {F_{45a}}\left( {{a_1},{a_2};{c_1},{c_2},{c_3};x,y,z} \right)$.

Particular solutions:

$
{u_1} = {F_{45a}}\left( {{a_1},{a_2};{c_1},{c_2},{c_3};x,y,z} \right),
$

$
{u_2} = {x^{1 - {c_1}}}{F_{45a}}\left( {2 - 2{c_1} + {a_1},{a_2};2 - {c_1},{c_2},{c_3};x,y,z} \right),
$

$
{u_3} = {y^{1 - {c_2}}}{F_{45a}}\left( {1 - {c_2} + {a_1},1 - {c_2} + {a_2};{c_1},2 - {c_2},{c_3};x,y,z} \right),
$

$
{u_4} = {z^{1 - {c_3}}}{F_{45a}}\left( {1 - {c_3} + {a_1},1 - {c_3} + {a_2};{c_1},{c_2},2 - {c_3};x,y,z} \right),
$

$
{u_5} = {x^{1 - {c_1}}}{y^{1 - {c_2}}}{F_{45a}}\left( {3 - 2{c_1} - {c_2} + {a_1},1 - {c_2} + {a_2};2 - {c_1},2 - {c_2},{c_3};x,y,z} \right),
$

$
{u_6} = {y^{1 - {c_2}}}{z^{1 - {c_3}}}{F_{45a}}\left( {2 - {c_2} - {c_3} + {a_1},2 - {c_2} - {c_3} + {a_2};{c_1},2 - {c_2},2 - {c_3};x,y,z} \right),
$

$
{u_7} = {x^{1 - {c_1}}}{z^{1 - {c_3}}}{F_{45a}}\left( {3 - 2{c_1} - {c_3} + {a_1},1 - {c_3} + {a_2};2 - {c_1},{c_2},2 - {c_3};x,y,z} \right),
$

$
  {u_8} = {x^{1 - {c_1}}}{y^{1 - {c_2}}}{z^{1 - {c_3}}} {F_{45a}}\left( {4 - 2{c_1} - {c_2} - {c_3} + {a_1},2 - {c_2} - {c_3} + {a_2};2 - {c_1},2 - {c_2},2 - {c_3};x,y,z} \right).
  $

\bigskip

\begin{equation}
{F_{45b}}\left( {a,b;{c_1},{c_2};x,y,z} \right) = \sum\limits_{m,n,p = 0}^\infty  {} \frac{{{{\left( a \right)}_{2m + n + p}}{{\left( b \right)}_{n + p - m}}}}{{{{\left( {{c_1}} \right)}_n}{{\left( {{c_2}} \right)}_p}}}\frac{x^m}{m!}\frac{y^n}{n!}\frac{z^p}{p!},
\end{equation}

region of convergence:
$$
\left\{ {r < \frac{1}{4},\,\,\,\sqrt s  + \sqrt t  < \min \left\{ {\sqrt {{\Psi _1}\left( r \right)} ,\sqrt {{\Psi _2}\left( r \right)} } \right\}} \right\}.
$$

System of partial differential equations:

$
\left\{ {\begin{array}{*{20}{l}}
  \begin{gathered}
  x\left( {1 + 4x} \right){u_{xx}} - \left( {1 - 4x} \right)y{u_{xy}} - \left( {1 - 4x} \right)z{u_{xz}}
   + 2yz{u_{yz}} + {y^2}{u_{yy}} + {z^2}{u_{zz}} \hfill \\ \,\,\,\,\,\,\,\,\, + \left[ {1 - b + 2\left( {2a + 3} \right)x} \right]{u_x}
   + 2\left( {a + 1} \right)\left( {y{u_y} + z{u_z}} \right) + a\left( {1 + a} \right)u = 0, \hfill \\
\end{gathered}  \\
  \begin{gathered}
  y\left( {1 - y} \right){u_{yy}} - xy{u_{xy}} - xz{u_{xz}} - 2yz{u_{yz}} + 2{x^2}{u_{xx}} - {z^2}{u_{zz}} \hfill \\\,\,\,\,\,\,\,\,\,
   + \left[ {{c_1} - \left( {a + b + 1} \right)y} \right]{u_y} + \left( {a - 2b + 2} \right)x{u_x}
   - \left( {a + b + 1} \right)z{u_z} - abu = 0, \hfill \\
\end{gathered}  \\
  \begin{gathered}
  z\left( {1 - z} \right){u_{zz}} - xy{u_{xy}} - xz{u_{xz}} - 2yz{u_{yz}} + 2{x^2}{u_{xx}} - {y^2}{u_{yy}} \hfill \\
 \,\,\,\,\,\,\,\,\,  + \left[ {{c_2} - \left( {a + b + 1} \right)z} \right]{u_z} + \left( {a - 2b + 2} \right)x{u_x}
   - \left( {a + b + 1} \right)y{u_y} - abu = 0, \hfill \\
\end{gathered}
\end{array}} \right.
$

where $u\equiv \,\,   {F_{45b}}\left( {a,b;{c_1},{c_2};x,y,z} \right)$.

Particular solutions:

$
{u_1} = {F_{45b}}\left( {a,b;{c_1},{c_2};x,y,z} \right),
$

$
{u_2} = {y^{1 - {c_1}}}{F_{45b}}\left( {1 - {c_1} + a,1 - {c_1} + b;2 - {c_1},{c_2};x,y,z} \right),
$

$
{u_3} = {z^{1 - {c_2}}}{F_{45b}}\left( {1 - {c_2} + a,1 - {c_2} + b;{c_1},2 - {c_2};x,y,z} \right),
$

$
{u_4} = {y^{1 - {c_1}}}{z^{1 - {c_2}}}{F_{45b}}\left( {1 - {c_1} + a,1 - {c_1} + b;2 - {c_1},2 - {c_2};x,y,z} \right).
$

\bigskip

\begin{equation}
{F_{46a}}\left( {{a_1},{a_2};{c_1},{c_2};x,y,z} \right) = \sum\limits_{m,n,p = 0}^\infty  {} \frac{{{{\left( {{a_1}} \right)}_{2m + n + p}}{{\left( {{a_2}} \right)}_{n + p}}}}{{{{\left( {{c_1}} \right)}_{m + n}}{{\left( {{c_2}} \right)}_p}}}\frac{x^m}{m!}\frac{y^n}{n!}\frac{z^p}{p!},
\end{equation}

first appearance of this function in the literature, and old notation:
[11], ${}^{(1)}H_{4}^{(3)}$; [12], $X_{3}$,

region of convergence:
$$
\begin{gathered}
   \left\{ {r < \frac{1}{4},\,\,\,s < \frac{1}{2} + \frac{1}{2}\sqrt {1 - 4r} ,\,\,\,t < \frac{{\left( {{w_1} - 2r} \right)\left( {{w_1} - s} \right)^2}}{{w_1^3}}} \right\}, \hfill \\
  P\left( w \right) = {w^3} - \left( {r + s} \right)w + 2rs, \hfill \\
{w_1}:\,\,{\rm{the\,\,root\,\,in}}\,\,\left( {2r,\sqrt {r + s} } \right)\,\,{\rm{of}}\,\,P\left( w \right) = 0. \hfill \\
\end{gathered}
$$

System of partial differential equations:

$
\left\{ {\begin{array}{*{20}{l}}
  \begin{gathered}
  x\left( {1 - 4x} \right){u_{xx}} + \left( {1 - 4x} \right)y{u_{xy}} - 4xz{u_{xz}}
  - 2yz{u_{yz}} - {y^2}{u_{yy}} - {z^2}{u_{zz}} \hfill \\ \,\,\,\,\,\,\,\,\, + \left[ {{c_1} - 2\left( {2{a_1} + 3} \right)x} \right]{u_x}
   - 2\left( {{a_1} + 1} \right)\left( {y{u_y} + z{u_z}} \right) - {a_1}\left( {1 + {a_1}} \right)u = 0, \hfill \\
\end{gathered}  \\
  \begin{gathered}
  y\left( {1 - y} \right){u_{yy}} + x\left( {1 - 2y} \right){u_{xy}} - 2xz{u_{xz}} - 2yz{u_{yz}} - {z^2}{u_{zz}} \hfill \\
 \,\,\,\,\,\,\,\,\,  + \left[ {{c_1} - \left( {{a_1} + {a_2} + 1} \right)y} \right]{u_y} - 2{a_2}x{u_x}
   - \left( {{a_1} + {a_2} + 1} \right)z{u_z} - {a_1}{a_2}u = 0, \hfill \\
\end{gathered}  \\
  \begin{gathered}
  z\left( {1 - z} \right){u_{zz}} - 2xy{u_{xy}} - 2xz{u_{xz}} - 2yz{u_{yz}} - {y^2}{u_{yy}} \hfill \\
 \,\,\,\,\,\,\,\,\,  + \left[ {{c_2} - \left( {{a_1} + {a_2} + 1} \right)z} \right]{u_z} - 2{a_2}x{u_x}
   - \left( {{a_1} + {a_2} + 1} \right)y{u_y} - {a_1}{a_2}u = 0, \hfill \\
\end{gathered}
\end{array}} \right.
$

where $u\equiv \,\,   {F_{46a}}\left( {{a_1},{a_2};{c_1},{c_2};x,y,z} \right)$.

Particular solutions:

$
{u_1} = {F_{46a}}\left( {{a_1},{a_2};{c_1},{c_2};x,y,z} \right),
$

$
{u_2} = {z^{1 - {c_3}}}{F_{46a}}\left( {1 - {c_3} + {a_1},1 - {c_3} + {a_2};{c_1},2 - {c_2};x,y,z} \right).
$

\bigskip

\begin{equation}
{F_{47b}}\left( {{a_1},{a_2},{b_1},{b_2};c;x,y,z} \right) = \sum\limits_{m,n,p = 0}^\infty  {} \frac{{{{\left( {{a_1}} \right)}_m}{{\left( {{a_2}} \right)}_m}{{\left( {{b_1}} \right)}_{2n - m}}{{\left( {{b_2}} \right)}_{2p - n}}}}{{{{\left( c \right)}_p}}}\frac{x^m}{m!}\frac{y^n}{n!}\frac{z^p}{p!},
\end{equation}

region of convergence:
$$
\left\{ {t < \frac{1}{4},\,\,\,s < \frac{1}{{4\left( {1 + 2\sqrt t } \right)}},\,\,\,r < \frac{1}{{1 + 2\sqrt {s\left( {1 + 2\sqrt t } \right)} }}} \right\}.
$$

System of partial differential equations:

$
\left\{ {\begin{array}{*{20}{l}}
  \begin{gathered}
  x\left( {1 + x} \right){u_{xx}} - 2y{u_{xy}}
   + \left[ {1 - {b_1} + \left( {{a_1} + {a_2} + 1} \right)x} \right]{u_x} + {a_1}{a_2}u = 0, \hfill \\
\end{gathered}  \\
  \begin{gathered}
  y\left( {1 + 4y} \right){u_{yy}} - 4xy{u_{xy}} - 2z{u_{yz}} + {x^2}{u_{xx}}\hfill\\\,\,\,\,\,\,\,\,\,
   + \left[ {1 - {b_2} + 2\left( {2{b_1} + 3} \right)y} \right]{u_y} - 2{b_1}x{u_x} + {b_1}\left( {1 + {b_1}} \right)u = 0, \hfill \\
\end{gathered}  \\
  \begin{gathered}
  z\left( {1 - 4z} \right){u_{zz}} + 4yz{u_{yz}} - {y^2}{u_{yy}}
   + \left[ {c - 2\left( {2{b_2} + 3} \right)z} \right]{u_z} + 2{b_2}y{u_y} - {b_2}\left( {1 + {b_2}} \right)u = 0, \hfill \\
\end{gathered}
\end{array}} \right.
$

where $u\equiv \,\,   {F_{47b}}\left( {{a_1},{a_2},{b_1},{b_2};c;x,y,z} \right)$.

Particular solutions:

$
{u_1} = {F_{47b}}\left( {{a_1},{a_2},{b_1},{b_2};c;x,y,z} \right),
$

$
{u_2} = {z^{1 - c}}{F_{47b}}\left( {{a_1},{a_2},{b_1},2 - 2c + {b_2};2 - c;x,y,z} \right).
$

\bigskip

\begin{equation}
{F_{47c}}\left( {a,{b_1},{b_2},{b_3};x,y,z} \right) = \sum\limits_{m,n,p = 0}^\infty  {} {{{{\left( a \right)}_m}{{\left( {{b_1}} \right)}_{2n - m}}{{\left( {{b_2}} \right)}_{2p - n}}{{\left( {{b_3}} \right)}_{m - p}}}}\frac{x^m}{m!}\frac{y^n}{n!}\frac{z^p}{p!},
\end{equation}

region of convergence:
$$
\begin{gathered}
    \left\{ {s < \frac{1}{4},\,\,\,r < \frac{1}{{1 + 2\sqrt s }},\,} \,\, {t < \min \left\{ {\frac{{{{\left( {1 - 4s} \right)}^2}}}{{64{s^2}}},\frac{1}{{1 + r}}{\Phi _1}\left( {\frac{{{r^2}s}}{{{{\left( {1 + r} \right)}^2}}}} \right),\frac{1}{{1 - r}}{\Phi _2}\left( {\frac{{{r^2}s}}{{{{\left( {1 - r} \right)}^2}}}} \right)} \right\}} \right\}. \hfill \\
\end{gathered}
$$

System of partial differential equations:

$
\left\{ {\begin{array}{*{20}{l}}
  \begin{gathered}
  x\left( {1 + x} \right){u_{xx}} - 2y{u_{xy}} - xz{u_{xz}}
   + \left[ {1 - {b_1} + \left( {a + {b_3} + 1} \right)x} \right]{u_x} - az{u_z} + a{b_3}u = 0, \hfill \\
\end{gathered}  \\
  \begin{gathered}
  y\left( {1 + 4y} \right){u_{yy}} - 4xy{u_{xy}} - 2z{u_{yz}} + {x^2}{u_{xx}}\hfill\\\,\,\,\,\,\,\,\,\,
   + \left[ {1 - {b_2} + 2\left( {2{b_1} + 3} \right)y} \right]{u_y} - 2{b_1}x{u_x} + {b_1}\left( {1 + {b_1}} \right)u = 0, \hfill \\
\end{gathered}  \\
  \begin{gathered}
  z\left( {1 + 4z} \right){u_{zz}} - x{u_{xz}} - 4yz{u_{yz}} + {y^2}{u_{yy}}\hfill \\\,\,\,\,\,\,\,\,\,
   + \left[ {1 - {b_3} + 2\left( {2{b_2} + 3} \right)z} \right]{u_z} - 2{b_2}y{u_y} + {b_2}\left( {1 + {b_2}} \right)u = 0, \hfill \\
\end{gathered}
\end{array}} \right.
$

where $u\equiv \,\,   {F_{47c}}\left( {a,{b_1},{b_2},{b_3};x,y,z} \right)$.

\bigskip

\begin{equation}
{F_{48b}}\left( {{a_1},{a_2},{b_1},{b_2};x,y,z} \right) = \sum\limits_{m,n,p = 0}^\infty  {} {{{{\left( {{a_1}} \right)}_m}{{\left( {{a_2}} \right)}_m}{{\left( {{b_1}} \right)}_{2n - p}}{{\left( {{b_2}} \right)}_{2p - m - n}}}}\frac{x^m}{m!}\frac{y^n}{n!}\frac{z^p}{p!},
\end{equation}

region of convergence:
$$
\left\{ {t < \frac{1}{4},\,\,\,r < \frac{1}{{1 + 2\sqrt t }},\,\,\,s < \min \left\{ {{\Phi _1}\left( t \right),{\Phi _2}\left( t \right),\frac{r}{{1 - r}}{\Phi _2}\left( {\frac{{{r^2}t}}{{{{\left( {1 - r} \right)}^2}}}} \right)} \right\}} \right\}.
$$

System of partial differential equations:

$
\left\{ {\begin{array}{*{20}{l}}
  \begin{gathered}
  x\left( {1 + x} \right){u_{xx}} + y{u_{xy}} - 2z{u_{xz}}
   + \left[ {1 - {b_2} + \left( {{a_1} + {a_2} + 1} \right)x} \right]{u_x} + {a_1}{a_2}u = 0, \hfill \\
\end{gathered}  \\
  \begin{gathered}
  y\left( {1 + 4y} \right){u_{yy}} + x{u_{xy}} - 2\left( {1 + 2y} \right)z{u_{yz}} + {z^2}{u_{zz}}
  \hfill \\ \,\,\,\,\,\,\,\,\,
   + \left[ {1 - {b_2} + 2\left( {2{b_1} + 3} \right)y} \right]{u_y} - 2{b_1}z{u_z} + {b_1}\left( {1 + {b_1}} \right)u = 0, \hfill \\
\end{gathered}  \\
  \begin{gathered}
  z\left( {1 + 4z} \right){u_{zz}} + 2xy{u_{xy}} - 4xz{u_{xz}} - 2y\left( {1 + 2z} \right){u_{yz}}
   + {x^2}{u_{xx}} + {y^2}{u_{yy}} \hfill \\ \,\,\,\,\,\,\,\,\,+ \left[ {1 - {b_1} + 2\left( {2{b_2} + 3} \right)z} \right]{u_z} - 2{b_2}x{u_x}
    - 2{b_2}y{u_y} + {b_2}\left( {1 + {b_2}} \right)u = 0, \hfill \\
\end{gathered}
\end{array}} \right.
$

where $u\equiv \,\,   {F_{48b}}\left( {{a_1},{a_2},{b_1},{b_2};x,y,z} \right)$.

\bigskip

\begin{equation}
{F_{49b}}\left( {{a_1},{a_2},b;{c_1},{c_2};x,y,z} \right) = \sum\limits_{m,n,p = 0}^\infty  {} \frac{{{{\left( {{a_1}} \right)}_{2m + n}}{{\left( {{a_2}} \right)}_n}{{\left( b \right)}_{2p - m}}}}{{{{\left( {{c_1}} \right)}_n}{{\left( {{c_2}} \right)}_p}}}\frac{x^m}{m!}\frac{y^n}{n!}\frac{z^p}{p!},
\end{equation}

region of convergence:
$$
 \left\{ {s < 1,\,\,\,t < \frac{1}{4},\,\,\,r < \frac{{{{\left( {1 - s} \right)}^2}}}{{4\left( {1 + 2\sqrt t } \right)}}} \right\}.
$$

System of partial differential equations:

$
\left\{ {\begin{array}{*{20}{l}}
  \begin{gathered}
  x\left( {1 + 4x} \right){u_{xx}} + 4xy{u_{xy}} - 2z{u_{xz}} + {y^2}{u_{yy}} \hfill \\
 \,\,\,\,\,\,\,\,\,  + \left[ {1 - b + 2\left( {2{a_1} + 3} \right)x} \right]{u_x} + 2\left( {{a_1} + 1} \right)y{u_y}
    + {a_1}\left( {1 + {a_1}} \right)u = 0, \hfill \\
\end{gathered}  \\
  y\left( {1 - y} \right){u_{yy}} - 2xy{u_{xy}} + \left[ {{c_1} - \left( {{a_1} + {a_2} + 1} \right)y} \right]{u_y}
   - 2{a_2}x{u_x} - {a_1}{a_2}u = 0, \\
  \begin{gathered}
  z\left( {1 - 4z} \right){u_{zz}} + 4xz{u_{xz}} - {x^2}{u_{xx}}
   + \left[ {{c_2} - 2\left( {2b + 3} \right)z} \right]{u_z} + 2bx{u_x} - b\left( {1 + b} \right)u = 0, \hfill \\
\end{gathered}
\end{array}} \right.
$

where $u\equiv \,\,   {F_{49b}}\left( {{a_1},{a_2},b;{c_1},{c_2};x,y,z} \right)$.

Particular solutions:

$
{u_1} = {F_{49b}}\left( {{a_1},{a_2},b;{c_1},{c_2};x,y,z} \right),
$

$
{u_2} = {y^{1 - {c_1}}}{F_{49b}}\left( {1 - {c_1} + {a_1},1 - {c_1} + {a_2},b;2 - {c_1},{c_2};x,y,z} \right),
$

$
{u_3} = {z^{1 - {c_2}}}{F_{49b}}\left( {{a_1},{a_2},2 - 2{c_2} + b;{c_1},2 - {c_2};x,y,z} \right),
$

$
{u_4} = {y^{1 - {c_1}}}{z^{1 - {c_2}}}{F_{49b}}\left( {1 - {c_1} + {a_1},1 - {c_1} + {a_2},2 - 2{c_2} + b;2 - {c_1},2 - {c_2};x,y,z} \right).
$

\bigskip

\begin{equation}
{F_{49c}}\left( {{a_1},{a_2},b;{c_1},{c_2};x,y,z} \right) = \sum\limits_{m,n,p = 0}^\infty  {} \frac{{{{\left( {{a_1}} \right)}_{2m + n}}{{\left( {{a_2}} \right)}_n}{{\left( b \right)}_{2p - n}}}}{{{{\left( {{c_1}} \right)}_m}{{\left( {{c_2}} \right)}_p}}}\frac{x^m}{m!}\frac{y^n}{n!}\frac{z^p}{p!},
\end{equation}

region of convergence:
$$
 \left\{ {r < \frac{1}{4},\,\,\,t < \frac{1}{4},\,\,\,s < \frac{{1 - 2\sqrt r }}{{1 + 2\sqrt t }}} \right\}.
$$

System of partial differential equations:

$
\left\{ {\begin{array}{*{20}{l}}
  \begin{gathered}
  x\left( {1 - 4x} \right){u_{xx}} - 4xy{u_{xy}} - {y^2}{u_{yy}}
   + \left[ {{c_1} - 2\left( {2{a_1} + 3} \right)x} \right]{u_x} - 2\left( {{a_1} + 1} \right)y{u_y} - {a_1}\left( {1 + {a_1}} \right)u = 0, \hfill \\
\end{gathered}  \\
  \begin{gathered}
  y\left( {1 + y} \right){u_{yy}} + 2xy{u_{xy}} - 2z{u_{yz}}
   + \left[ {1 - b + \left( {{a_1} + {a_2} + 1} \right)y} \right]{u_y} + 2{a_2}x{u_x} + {a_1}{a_2}u = 0, \hfill \\
\end{gathered}  \\
  \begin{gathered}
  z\left( {1 - 4z} \right){u_{zz}} + 4yz{u_{yz}} - {y^2}{u_{yy}}
  + \left[ {{c_2} - 2\left( {2b + 3} \right)z} \right]{u_z} + 2by{u_y} - b\left( {1 + b} \right)u = 0, \hfill \\
\end{gathered}
\end{array}} \right.
$

where $u\equiv \,\,   {F_{49c}}\left( {{a_1},{a_2},b;{c_1},{c_2};x,y,z} \right)$.

Particular solutions:

$
{u_1} = {F_{49c}}\left( {{a_1},{a_2},b;{c_1},{c_2};x,y,z} \right),
$

$
{u_2} = {x^{1 - {c_1}}}{F_{49c}}\left( {{a_1},{a_2},b;2 - {c_1},{c_2};x,y,z} \right),
$

$
{u_3} = {z^{1 - {c_2}}}{F_{49c}}\left( {{a_1},{a_2},2 - 2{c_2} + b;{c_1},2 - {c_2};x,y,z} \right),
$

$
{u_4} = {x^{1 - {c_1}}}{z^{1 - {c_2}}}{F_{49c}}\left( {2 - 2{c_1} + {a_1},{a_2},2 - 2{c_2} + b;2 - {c_1},2 - {c_2};x,y,z} \right).
$

\bigskip

\begin{equation}
{F_{49d}}\left( {a,{b_1},{b_2};c;x,y,z} \right) = \sum\limits_{m,n,p = 0}^\infty  {} \frac{{{{\left( a \right)}_n}{{\left( {{b_1}} \right)}_{2m + n - p}}{{\left( {{b_2}} \right)}_{2p - m}}}}{{{{\left( c \right)}_n}}}\frac{x^m}{m!}\frac{y^n}{n!}\frac{z^p}{p!},
\end{equation}

region of convergence:
$$
\left\{ {s + 2\sqrt r  < 1,\,\,\,t < \min \left\{ {\frac{1}{{1 + s}}{\Phi _1}\left( {\frac{r}{{{{\left( {1 + s} \right)}^2}}}} \right),\frac{1}{{1 - s}}{\Phi _2}\left( {\frac{r}{{{{\left( {1 - s} \right)}^2}}}} \right)} \right\}} \right\}.
$$

System of partial differential equations:

$
\left\{ {\begin{array}{*{20}{l}}
  \begin{gathered}
  x\left( {1 + 4x} \right){u_{xx}} + 4xy{u_{xy}} - 2\left( {1 + 2x} \right)z{u_{xz}} - 2yz{u_{yz}}
  + {y^2}{u_{yy}} + {z^2}{u_{zz}}  \hfill \\
  \,\,\,\,\,\,\,\,\,+ \left[ {1 - {b_2} + 2\left( {2{b_1} + 3} \right)x} \right]{u_x} + 2\left( {{b_1} + 1} \right)y{u_y}
   - 2{b_1}z{u_z} + {b_1}\left( {1 + {b_1}} \right)u = 0, \hfill \\
\end{gathered}  \\
  \begin{gathered}
  y\left( {1 - y} \right){u_{yy}} - 2xy{u_{xy}} + yz{u_{yz}}
   + \left[ {c - \left( {a + {b_1} + 1} \right)y} \right]{u_y} - 2ax{u_x} + az{u_z} - a{b_1}u = 0, \hfill \\
\end{gathered}  \\
  \begin{gathered}
  z\left( {1 + 4z} \right){u_{zz}} - 2x\left( {1 + 2z} \right){u_{xz}} - y{u_{yz}} + {x^2}{u_{xx}} \hfill \\\,\,\,\,\,\,\,\,\,
   + \left[ {1 - {b_1} + 2\left( {2{b_2} + 3} \right)z} \right]{u_z} - 2{b_2}x{u_x} + {b_2}\left( {1 + {b_2}} \right)u = 0, \hfill \\
\end{gathered}
\end{array}} \right.
$

where $u\equiv \,\,   {F_{49d}}\left( {a,{b_1},{b_2};c;x,y,z} \right)$.

Particular solutions:

$
{u_1} = {F_{49d}}\left( {a,{b_1},{b_2};c;x,y,z} \right),
$

$
{u_2} = {z^{1 - c}}{F_{49d}}\left( {a,{b_1},2 - 2c + {b_2};2 - c;x,y,z} \right).
$

\bigskip

\begin{equation}
{F_{49e}}\left( {a,{b_1},{b_2};c;x,y,z} \right) = \sum\limits_{m,n,p = 0}^\infty  {} \frac{{{{\left( a \right)}_{2m + n}}{{\left( {{b_1}} \right)}_{n - p}}{{\left( {{b_2}} \right)}_{2p - m}}}}{{{{\left( c \right)}_n}}}\frac{x^m}{m!}\frac{y^n}{n!}\frac{z^p}{p!},
\end{equation}

region of convergence:
$$
\begin{gathered}
  \left\{ {s + 2\sqrt r  < 1,\,\,\,t < \min \left\{ {{U^ + }\left( {{w_1}} \right),{U^ - }\left( {{w_2}} \right),{U^ + }\left( {{w_3}} \right)} \right\}} \right\}, \hfill \\
  {P_{rs}}\left( w \right) = 8{w^3} - \left( {8 + 5s} \right){w^2} + 2\left( {1 + s} \right)w + rs, \hfill \\
{w_1}:\,\,{\rm{the\,\,root\,\,in}}\,\,\left( {\sqrt r ,\frac{1}{2}} \right)\,\,{\rm{of}}\,\,{P_{rs}}\left( w \right) = 0, \hfill \\
{w_2}:\,\,{\rm{the\,\,root\,\,in}}\,\,\left( {0,\frac{1}{2}\sqrt r } \right)\,\,{\rm{of}}\,\,{P_{ - r,s}}\left( w \right) = 0, \hfill \\
{w_3}:\,\,{\rm{the\,\,smaller\,\,root\,\,in}}\,\,\left( {\sqrt r ,\frac{1}{2}} \right)\,\,{\rm{of}}\,\,{P_{r, - s}}\left( w \right) = 0, \hfill \\
  {U^ \pm }\left( w \right) = \frac{{w{{\left( {1 - 2w} \right)}^2}\left( {{w^2} \mp r} \right)}}{{2{r^2}s}}. \hfill \\
\end{gathered}
$$

System of partial differential equations:

$
\left\{ {\begin{array}{*{20}{l}}
  \begin{gathered}
  x\left( {1 + 4x} \right){u_{xx}} + 4xy{u_{xy}} - 2z{u_{xz}} + {y^2}{u_{yy}}  \hfill\\ \,\,\,\,\,\,\,\,\, + \left[ {1 - {b_2} + 2\left( {2a + 3} \right)x} \right]{u_x} + 2\left( {a + 1} \right)y{u_y} + a\left( {1 + a} \right)u = 0, \hfill \\
\end{gathered}  \\
  \begin{gathered}
  y\left( {1 - y} \right){u_{yy}} - 2xy{u_{xy}} + 2xz{u_{xz}} + yz{u_{yz}} \hfill\\\,\,\,\,\,\,\,\,\,
   + \left[ {c - \left( {a + {b_1} + 1} \right)y} \right]{u_y} - 2{b_1}x{u_x} + az{u_z} - a{b_1}u = 0, \hfill \\
\end{gathered}  \\
  \begin{gathered}
  z\left( {1 + 4z} \right){u_{zz}} - 4xz{u_{xz}} - y{u_{yz}} + {x^2}{u_{xx}}\hfill\\\,\,\,\,\,\,\,\,\,
   + \left[ {1 - {b_1} + 2\left( {2{b_2} + 3} \right)z} \right]{u_z} - 2{b_2}x{u_x} + {b_2}\left( {1 + {b_2}} \right)u = 0, \hfill \\
\end{gathered}
\end{array}} \right.
$

where $u\equiv \,\,   {F_{49e}}\left( {a,{b_1},{b_2};c;x,y,z} \right)$.

Particular solutions:

$
{u_1} = {F_{49e}}\left( {a,{b_1},{b_2};c;x,y,z} \right),
$

${u_2} = {y^{1 - c}}{F_{49e}}\left( {1 - c + a,1 - c + {b_1},{b_2};2 - c;x,y,z} \right).$

\bigskip

\begin{equation}
{F_{49f}}\left( {a,{b_1},{b_2};c;x,y,z} \right) = \sum\limits_{m,n,p = 0}^\infty  {} \frac{{{{\left( a \right)}_n}{{\left( {{b_1}} \right)}_{2m + n - p}}{{\left( {{b_2}} \right)}_{2p - n}}}}{{{{\left( c \right)}_m}}}\frac{x^m}{m!}\frac{y^n}{n!}\frac{z^p}{p!},
\end{equation}

region of convergence:
$$
 \left\{ {s + 2\sqrt r  < 1,\,\,\,t < \min \left\{ {\frac{1}{{4\left( {1 + 2\sqrt r } \right)}},\frac{{1 - s - 2\sqrt r }}{{{s^2}}}} \right\}} \right\}.
$$

System of partial differential equations:

$
\left\{ {\begin{array}{*{20}{l}}
  \begin{gathered}
  x\left( {1 - 4x} \right){u_{xx}} - 4xy{u_{xy}} + 4xz{u_{xz}} + 2yz{u_{yz}}
  - {y^2}{u_{yy}} - {z^2}{u_{zz}} \hfill \\ \,\,\,\,\,\,\,\,\,+ \left[ {c - 2\left( {2{b_1} + 3} \right)x} \right]{u_x} - 2\left( {{b_1} + 1} \right)y{u_y}
  + 2{b_1}z{u_z} - {b_1}\left( {1 + {b_1}} \right)u = 0, \hfill \\
\end{gathered}  \\
  \begin{gathered}
  y\left( {1 + y} \right){u_{yy}} + 2xy{u_{xy}} - \left( {2 + y} \right)z{u_{yz}}\hfill\\\,\,\,\,\,\,\,\,\,
   + \left[ {1 - {b_2} + \left( {a + {b_1} + 1} \right)y} \right]{u_y} + 2ax{u_x} - az{u_z} + a{b_1}u = 0, \hfill \\
\end{gathered}  \\
  \begin{gathered}
  z\left( {1 + 4z} \right){u_{zz}} - 2x{u_{xz}} - y\left( {1 + 4z} \right){u_{yz}} + {y^2}{u_{yy}}\hfill\\\,\,\,\,\,\,\,\,\,
  + \left[ {1 - {b_1} + 2\left( {2{b_2} + 3} \right)z} \right]{u_z} - 2{b_2}y{u_y} + {b_2}\left( {1 + {b_2}} \right)u = 0, \hfill \\
\end{gathered}
\end{array}} \right.
$

where $u\equiv \,\,   {F_{49f}}\left( {a,{b_1},{b_2};c;x,y,z} \right)$.

Particular solutions:

$
{u_1} = {F_{49f}}\left( {a,{b_1},{b_2};c;x,y,z} \right),
$

$
{u_2} = {x^{1 - c}}{F_{49f}}\left( {a,2 - 2c + {b_1},{b_2};2 - c;x,y,z} \right).
$

\bigskip

\begin{equation}
{F_{49g}}\left( {a,{b_1},{b_2};c;x,y,z} \right) = \sum\limits_{m,n,p = 0}^\infty  {} \frac{{{{\left( a \right)}_{2m + n}}{{\left( {{b_1}} \right)}_{2p - n}}{{\left( {{b_2}} \right)}_{n - m}}}}{{{{\left( c \right)}_p}}}\frac{x^m}{m!}\frac{y^n}{n!}\frac{z^p}{p!},
\end{equation}

region of convergence:
$$
\left\{ {r < \frac{1}{4},\,\,\,t < \frac{1}{4},\,\,\,s < \frac{1}{{1 + 2\sqrt t }}\min \left\{ {{\Phi _1}\left( r \right),{\Phi _2}\left( r \right)} \right\}} \right\}.
$$

System of partial differential equations:

$
\left\{ {\begin{array}{*{20}{l}}
  \begin{gathered}
  x\left( {1 + 4x} \right){u_{xx}} - \left( {1 - 4x} \right)y{u_{xy}} + {y^2}{u_{yy}}\hfill\\\,\,\,\,\,\,\,\,\,
   + \left[ {1 - {b_2} + 2\left( {2a + 3} \right)x} \right]{u_x} + 2\left( {a + 1} \right)y{u_y} + a\left( {1 + a} \right)u = 0, \hfill \\
\end{gathered}  \\
  \begin{gathered}
  y\left( {1 + y} \right){u_{yy}} + xy{u_{xy}} - 2z{u_{yz}} - 2{x^2}{u_{xx}}\hfill\\\,\,\,\,\,\,\,\,\,
   + \left[ {1 - {b_1} + \left( {a + {b_2} + 1} \right)y} \right]{u_y} - \left( {a - 2{b_2} + 2} \right)x{u_x} + a{b_2}u = 0, \hfill \\
\end{gathered}  \\
  \begin{gathered}
  z\left( {1 - 4z} \right){u_{zz}} + 4yz{u_{yz}} - {y^2}{u_{yy}}
   + \left[ {c - 2\left( {2{b_1} + 3} \right)z} \right]{u_z} + 2{b_1}y{u_y} - {b_1}\left( {1 + {b_1}} \right)u = 0, \hfill \\
\end{gathered}
\end{array}} \right.
$

where $u\equiv \,\,   {F_{49g}}\left( {a,{b_1},{b_2};c;x,y,z} \right)$.

Particular solutions:

$
{u_1} = {F_{49g}}\left( {a,{b_1},{b_2};c;x,y,z} \right),
$

$
{u_2} = {z^{1 - c}}{F_{49g}}\left( {a,2 - 2c + {b_1},{b_2};2 - c;x,y,z} \right).
$

\bigskip

\begin{equation}
{F_{49h}}\left( {a,{b_1},{b_2};c;x,y,z} \right) = \sum\limits_{m,n,p = 0}^\infty  {} \frac{{{{\left( a \right)}_{2m + n}}{{\left( {{b_1}} \right)}_{2p - n}}{{\left( {{b_2}} \right)}_{n - p}}}}{{{{\left( c \right)}_m}}}\frac{x^m}{m!}\frac{y^n}{n!}\frac{z^p}{p!},
\end{equation}

region of convergence:
$$
\left\{ {r < \frac{1}{4},\,\,\,t < \frac{1}{4},\,\,\,s < \frac{{\left( {1 - 2\sqrt r } \right)\left( {\sqrt {1 + 4t}  - 1} \right)}}{{2t}}} \right\}.
$$

System of partial differential equations:

$
\left\{ {\begin{array}{*{20}{l}}
  \begin{gathered}
  x\left( {1 - 4x} \right){u_{xx}} - 4xy{u_{xy}} - {y^2}{u_{yy}}
   + \left[ {c - 2\left( {2a + 3} \right)x} \right]{u_x} - 2\left( {a + 1} \right)y{u_y} - a\left( {1 + a} \right) = 0, \hfill \\
\end{gathered}  \\
  \begin{gathered}
  y\left( {1 + y} \right){u_{yy}} + 2xy{u_{xy}} - 2xz{u_{xz}}
   - z\left( {y + 2} \right){u_{yz}}
   \hfill \\ \,\,\,\,\,\,\,\,\,+\left[ {1 - {b_1} + \left( {a + {b_2} + 1} \right)y} \right]{u_y} + 2{b_2}x{u_x}
   - az{u_z} + a{b_2}u = 0, \hfill \\
\end{gathered}  \\
  \begin{gathered}
  z\left( {1 + 4z} \right){u_{zz}} - y\left( {1 + 4z} \right){u_{yz}} + {y^2}{u_{yy}}\hfill\\\,\,\,\,\,\,\,\,\,
   + \left[ {1 - {b_2} + 2\left( {2{b_1} + 3} \right)z} \right]{u_z} - 2{b_1}y{u_y} + {b_1}\left( {1 + {b_1}} \right)u = 0, \hfill \\
\end{gathered}
\end{array}} \right.
$

where $u\equiv \,\,   {F_{49h}}\left( {a,{b_1},{b_2};c;x,y,z} \right)$.

Particular solutions:

$
{u_1} = {F_{49h}}\left( {a,{b_1},{b_2};c;x,y,z} \right),
$

$
{u_2} = {x^{1 - c}}{F_{49h}}\left( {2 - 2c +  + a,{b_1},{b_2};2 - c;x,y,z} \right).
$

\bigskip

\begin{equation}
{F_{49i}}\left( {{b_1},{b_2},{b_3};x,y,z} \right) = \sum\limits_{m,n,p = 0}^\infty  {} {{{{\left( {{b_1}} \right)}_{2m + n - p}}{{\left( {{b_2}} \right)}_{2p - n}}{{\left( {{b_2}} \right)}_{n - m}}}}\frac{x^m}{m!}\frac{y^n}{n!}\frac{z^p}{p!},
\end{equation}

region of convergence:
$$
\begin{gathered}
  \left\{ {r < \frac{1}{4},\,\,\,s < \min \left\{ {{\Psi _1}\left( r \right),{\Psi _2}\left( r \right)} \right\},\,\,\,t < \min \left\{ {{U^ + }\left( {{w_1}} \right), - {U^ + }\left( { - {w_2}} \right),{U^ - }\left( {{w_3}} \right)} \right\}} \right\}, \hfill \\
  {P_{rs}}\left( w \right) = 5{w^4} + 2{w^3} - 6r{w^2} - \left( {2 + s} \right)rw + {r^2}, \hfill \\
{w_1}:\,\,{\rm{the\,\,root\,\,in}}\,\,\left( {\sqrt r, \infty} \right)\,\,{\rm{of}}\,\,{P_{rs}}\left( w \right) = 0, \hfill \\
{w_2}:\,\,{\rm{the\,\,\,root\,\,in}}\,\,\left( {\frac{1}{6}\left[ {1 + a\left( r \right)} \right],\frac{1}{2}} \right)\,\,{\rm{of}}\,\,{P_{ r, s}}\left( -w \right) = 0, \hfill \\
{w_3}:\,\,{\rm{the\,\,smaller\,\,root\,\,in}}\,\,\left( 0, \frac{2}{5}\right)\,\,{\rm{of}}\,\,{P_{-r, - s}}\left( -w \right) = 0, \hfill \\
  {U^ \pm }\left( w \right) = \frac{{\left[ {{3w^2} \pm (w-r)} \right]\left( {{w^2} \mp r} \right)^2}}{{r^2s^2w}}. \hfill \\
\end{gathered}
$$

System of partial differential equations:

$
\left\{ {\begin{array}{*{20}{l}}
  \begin{gathered}
  x\left( {1 + 4x} \right){u_{xx}} - \left( {1 - 4x} \right)y{u_{xy}} - 4xz{u_{xz}} - 2yz{u_{yz}}
   + {y^2}{u_{yy}} + {z^2}{u_{zz}}  \hfill \\ \,\,\,\,\,\,\,\,\,+ \left[ {1 - {b_3} + 2\left( {2{b_1} + 3} \right)x} \right]{u_x} + 2\left( {{b_1} + 1} \right)y{u_y}
    - 2{b_1}z{u_z} + {b_1}\left( {1 + {b_1}} \right)u = 0, \hfill \\
\end{gathered}  \\
  \begin{gathered}
  y\left( {1 + y} \right){u_{yy}} + xy{u_{xy}} + xz{u_{xz}} - \left( {2 + y} \right)z{u_{yz}} - 2{x^2}{u_{xx}} \hfill \\\,\,\,\,\,\,\,\,\,
   + \left[ {1 - {b_2} + \left( {{b_1} + {b_3} + 1} \right)y} \right]{u_y} - \left( {{b_1} - 2{b_3} + 2} \right)x{u_x}
   - {b_3}z{u_z} + {b_1}{b_3}u = 0, \hfill \\
\end{gathered}  \\
  \begin{gathered}
  z\left( {1 + 4z} \right){u_{zz}} - 2x{u_{xz}} - y\left( {1 + 4z} \right){u_{yz}} + {y^2}{u_{yy}}\hfill\\\,\,\,\,\,\,\,\,\,
   + \left[ {1 - {b_1} + 2\left( {2{b_2} + 3} \right)z} \right]{u_z} - 2{b_2}y{u_y} + {b_2}\left( {1 + {b_2}} \right)u = 0, \hfill \\
\end{gathered}
\end{array}} \right.
$

where $u\equiv \,\,   {F_{49i}}\left( {{b_1},{b_2},{b_3};x,y,z} \right)$.

\bigskip

\begin{equation}
{F_{50b}}\left( {{a_1},{a_2},b;c;x,y,z} \right) = \sum\limits_{m,n,p = 0}^\infty  {} \frac{{{{\left( {{a_1}} \right)}_{2m + n}}{{\left( {{a_2}} \right)}_n}{{\left( b \right)}_{2p - m - n}}}}{{{{\left( c \right)}_p}}}\frac{x^m}{m!}\frac{y^n}{n!}\frac{z^p}{p!},
\end{equation}

region of convergence:
$$
\left\{ {t < \frac{1}{4},\,\,\,r < \frac{1}{{4\left( {1 + 2\sqrt t } \right)}},\,\,\,s < \frac{{1 + \sqrt {1 - 4r\left( {1 + 2\sqrt t } \right)} }}{{2\left( {1 + 2\sqrt t } \right)}}} \right\}.
$$

System of partial differential equations:

$
\left\{ {\begin{array}{*{20}{l}}
  \begin{gathered}
  x\left( {1 + 4x} \right){u_{xx}} + \left( {1 + 4x} \right)y{u_{xy}} - 2z{u_{xz}} + {y^2}{u_{yy}} \hfill \\\,\,\,\,\,\,\,\,\,
   + \left[ {1 - b + 2\left( {2{a_1} + 3} \right)x} \right]{u_x} + 2\left( {{a_1} + 1} \right)y{u_y} + {a_1}\left( {1 + {a_1}} \right)u = 0, \hfill \\
\end{gathered}  \\
  \begin{gathered}
  y\left( {1 + y} \right){u_{yy}} + x\left( {1 + 2y} \right){u_{xy}} - 2z{u_{yz}}
   + \left[ {1 - b + \left( {{a_1} + {a_2} + 1} \right)y} \right]{u_y} + 2{a_2}x{u_x} + {a_1}{a_2}u = 0, \hfill \\
\end{gathered}  \\
  \begin{gathered}
  z\left( {1 - 4z} \right){u_{zz}} - 2xy{u_{xy}} + 4yz{u_{yz}} + 4xz{u_{xz}} - {x^2}{u_{xx}} - {y^2}{u_{yy}} \hfill \\\,\,\,\,\,\,\,\,\,
   + \left[ {c - 2\left( {2b + 3} \right)z} \right]{u_z} + 2bx{u_x} + 2by{u_y} - b\left( {1 + b} \right)u = 0, \hfill \\
\end{gathered}
\end{array}} \right.
$

where $u\equiv \,\,   {F_{50b}}\left( {{a_1},{a_2},b;c;x,y,z} \right)$.

Particular solutions:

$
{u_1} = {F_{50b}}\left( {{a_1},{a_2},b;c;x,y,z} \right),
$

$
{u_2} = {z^{1 - c}}{F_{50b}}\left( {{a_1},{a_2},2 - 2c + b;2 - c;x,y,z} \right).
$

\bigskip

\begin{equation}
{F_{50c}}\left( {a,{b_1},{b_2};x,y,z} \right) = \sum\limits_{m,n,p = 0}^\infty  {} {{{{\left( a \right)}_n}{{\left( {{b_1}} \right)}_{2m + n - p}}{{\left( {{b_2}} \right)}_{2p - m - n}}}}\frac{x^m}{m!}\frac{y^n}{n!}\frac{z^p}{p!},
\end{equation}

region of convergence:
$$
\begin{gathered}
  \left\{ {t < \frac{1}{4},\,\,\,s \leq \frac{1}{{1 + A(t)}},\,\,\,r < \min \left\{ {{\Phi _1}\left( t \right),{\Phi _2}\left( t \right)} \right\}} \right\} \hfill \\
   \cup \left\{ {t < \frac{1}{4},\,\,\,\frac{1}{{1 + A\left( t \right)}} < s < \frac{2}{{1 + \sqrt {1 + 4t} }},\,\,\,r < \min \left\{ {{\Phi _2}\left( t \right),s\left( {1 - s - {s^2}t} \right)} \right\}} \right\}. \hfill \\
\end{gathered}
$$

System of partial differential equations:

$
\left\{ {\begin{array}{*{20}{l}}
  \begin{gathered}
  x\left( {1 + 4x} \right){u_{xx}} + \left( {1 + 4x} \right)y{u_{xy}} - 2\left( {1 + 2x} \right)z{u_{xz}} - 2yz{u_{yz}}
  + {z^2}{u_{zz}} + {y^2}{u_{yy}} \hfill \\  \,\,\,\,\,\,\,\,\,+ \left[ {1 - {b_2} + 2\left( {2{b_1} + 3} \right)x} \right]{u_x}
  + 2\left( {{b_1} + 1} \right)y{u_y} - 2{b_1}z{u_z} + {b_1}\left( {1 + {b_1}} \right)u = 0, \hfill \\
\end{gathered}  \\
  \begin{gathered}
  y\left( {1 + y} \right){u_{yy}} + x\left( {1 + 2y} \right){u_{xy}} - \left( {2 + y} \right)z{u_{yz}}\hfill\\\,\,\,\,\,\,\,\,\,
   + \left[ {1 - {b_2} + \left( {a + {b_1} + 1} \right)y} \right]{u_y} + 2ax{u_x} - az{u_z} + a{b_1}u = 0, \hfill \\
\end{gathered}  \\
  \begin{gathered}
  z\left( {1 + 4z} \right){u_{zz}} + 2xy{u_{xy}} - 2x\left( {1 + 2z} \right){u_{xz}} - y\left( {1 + 4z} \right){u_{yz}}
   + {x^2}{u_{xx}} + {y^2}{u_{yy}}  \hfill \\ \,\,\,\,\,\,\,\,\,+ \left[ {1 - {b_1} + 2\left( {2{b_2} + 3} \right)z} \right]{u_z}
   - 2{b_2}x{u_x} - 2{b_2}y{u_y} + {b_2}\left( {1 + {b_2}} \right)u = 0, \hfill \\
\end{gathered}
\end{array}} \right.
$

where $u\equiv \,\,   {F_{50c}}\left( {a,{b_1},{b_2};x,y,z} \right)$.

\bigskip

\begin{equation}
{F_{50d}}\left( {a,{b_1},{b_2};x,y,z} \right) = \sum\limits_{m,n,p = 0}^\infty  {} {{{{\left( a \right)}_{2m + n}}{{\left( {{b_1}} \right)}_{2p - m - n}}{{\left( {{b_2}} \right)}_{n - p}}}}\frac{x^m}{m!}\frac{y^n}{n!}\frac{z^p}{p!},
\end{equation}

region of convergence:
$$
\begin{gathered}
  \left\{ {r < \frac{1}{4},\,\,\,s < \frac{1}{2} + \frac{1}{2}\sqrt {1 - 4r} ,}\,\,\,\, t < \min \left\{1,\frac{{\left( {2r + s} \right)\left( {r + s + {s^2}} \right)}}{{{s^4}}}\right.\times\right. \hfill \\
  \left. \times\left.{\Psi _1}\left( \frac{{r\left( {r + s + {s^2}} \right)}}{{{{\left( {2r + s} \right)}^2}}}  \right),\frac{{\left( {s - 2r} \right)\left( {s  - {s^2}-r} \right)}}{{{s^4}}}{\Psi _2}\left( \frac{{r\left( {s - {s^2}-r} \right)}}{{{{\left( {2r - s} \right)}^2}}}  \right) \right\} \right\}. \hfill \\
 \end{gathered}
$$

System of partial differential equations:

$
\left\{ {\begin{array}{*{20}{l}}
  \begin{gathered}
  x\left( {1 + 4x} \right){u_{xx}} + \left( {1 + 4x} \right)y{u_{xy}} - 2z{u_{xz}} + {y^2}{u_{yy}} \hfill \\\,\,\,\,\,\,\,\,\,
   + \left[ {1 - {b_1} + 2\left( {2a + 3} \right)x} \right]{u_x} + 2\left( {a + 1} \right)y{u_y} + a\left( {1 + a} \right)u = 0, \hfill \\
\end{gathered}  \\
  \begin{gathered}
  y\left( {1 + y} \right){u_{yy}} + x\left( {1 + 2y} \right){u_{xy}} - 2xz{u_{xz}} - z\left( {2 + y} \right){u_{yz}} \hfill \\\,\,\,\,\,\,\,\,\,
   + \left[ {1 - {b_1} + \left( {a + {b_2} + 1} \right)y} \right]{u_y} + 2{b_2}x{u_x} - az{u_z} + a{b_2}u = 0, \hfill \\
\end{gathered}  \\
  \begin{gathered}
  z\left( {1 + 4z} \right){u_{zz}} + 2xy{u_{xy}} - 4xz{u_{xz}} - y\left( {1 + 4z} \right){u_{yz}} + {x^2}{u_{xx}}\hfill \\\,\,\,\,\,\,\,\,\,
   + {y^2}{u_{yy}}  + \left[ {1 - {b_2} + 2\left( {2{b_1} + 3} \right)z} \right]{u_z} - 2{b_1}x{u_x}
   - 2{b_1}y{u_y} + {b_1}\left( {1 + {b_1}} \right)u = 0, \hfill \\
\end{gathered}
\end{array}} \right.
$

where $u\equiv \,\,   {F_{50d}}\left( {a,{b_1},{b_2};x,y,z} \right)$.

\bigskip

\begin{equation}
{F_{51b}}\left( {{a_1},{a_2},b;c;x,y,z} \right) = \sum\limits_{m,n,p = 0}^\infty  {} \frac{{{{\left( {{a_1}} \right)}_{2m + n}}{{\left( {{a_2}} \right)}_n}{{\left( b \right)}_{2p - n}}}}{{{{\left( c \right)}_{m + p}}}}\frac{x^m}{m!}\frac{y^n}{n!}\frac{z^p}{p!},
\end{equation}

region of convergence:
$$
 \left\{ {s + 2\sqrt r  < 1,\,\,\,t < \min \left\{ {\frac{1}{4},\frac{{{{\left( {1 - s} \right)}^2} - 4r}}{{4{s^2}}}} \right\}} \right\}.
$$

System of partial differential equations:

$
\left\{ {\begin{array}{*{20}{l}}
  \begin{gathered}
  x\left( {1 - 4x} \right){u_{xx}} - 4xy{u_{xy}} + z{u_{xz}} - {y^2}{u_{yy}}\hfill\\\,\,\,\,\,\,\,\,\,
   + \left[ {c - 2\left( {2{a_1} + 3} \right)x} \right]{u_x} - 2\left( {{a_1} + 1} \right)y{u_y} - {a_1}\left( {1 + {a_1}} \right)u = 0, \hfill \\
\end{gathered}  \\
  \begin{gathered}
  y\left( {1 + y} \right){u_{yy}} + 2xy{u_{xy}} - 2z{u_{yz}}
   + \left[ {1 - b + \left( {{a_1} + {a_2} + 1} \right)y} \right]{u_y} + 2{a_2}x{u_x} + {a_1}{a_2}u = 0, \hfill \\
\end{gathered}  \\
  \begin{gathered}
  z\left( {1 - 4z} \right){u_{zz}} + x{u_{xz}} + 4yz{u_{yz}} - {y^2}{u_{yy}}
   + \left[ {c - 2\left( {2b + 3} \right)z} \right]{u_z} + 2by{u_y} - b\left( {1 + b} \right)u = 0, \hfill \\
\end{gathered}
\end{array}} \right.
$

where $u\equiv \,\,   {F_{51b}}\left( {{a_1},{a_2},b;c;x,y,z} \right)$.

\bigskip

\begin{equation}
{F_{51c}}\left( {a,{b_1},{b_2};x,y,z} \right) = \sum\limits_{m,n,p = 0}^\infty  {} {{{{\left( a \right)}_{2m + n}}{{\left( {{b_1}} \right)}_{2p - n}}{{\left( {{b_2}} \right)}_{n - m - p}}}}\frac{x^m}{m!}\frac{y^n}{n!}\frac{z^p}{p!},
\end{equation}

region of convergence:
$$
\begin{gathered}
  \left\{ {r < \frac{1}{4},\,\,\,s < \min \left\{ {{\Psi _1}\left( r \right),{\Psi _2}\left( r \right)} \right\},\,\,\,t < \min \left\{ {\frac{1}{4},{U^ + }\left( {{w_1}} \right),{U^ - }\left( {{w_2}} \right)} \right\}} \right\}, \hfill \\
  {P^ \pm }\left( w \right) = 8{w^3} - 6{w^2} + \left( {1 \mp 4r} \right)w \pm 2r + rs, \hfill \\
{w_1}:\,\,{\rm{the\,\,root\,\,in}}\,\,\left( {0,\frac{1}{6}\left[ {1 + a\left( r \right)} \right]} \right)\,\,{\rm{of}}\,\,{P^ + }\left( w \right) = 0, \hfill \\
{w_2}:\,\,{\rm{the\,\,smaller\,\,root\,\,in}}\,\,\left( {\frac{1}{6}\left[ {1 - b\left( r \right)} \right],\frac{1}{6}\left[ {1 + b\left( r \right)} \right]} \right)\,\,{\rm{of}}\,\,{P^ - }\left( w \right) = 0, \hfill \\
  {U^ \pm }\left( w \right) = \frac{{\left( { - 3{w^2} + w \pm r} \right){{\left( {1 - 2w} \right)}^2}}}{{r{s^2}}}. \hfill \\
\end{gathered}
$$

System of partial differential equations:

$
\left\{ {\begin{array}{*{20}{l}}
  \begin{gathered}
  x\left( {1 + 4x} \right){u_{xx}} - \left( {1 - 4x} \right)y{u_{xy}} + z{u_{xz}} + {y^2}{u_{yy}} \hfill \\ \,\,\,\,\,\,\,\,\,
   + \left[ {1 - {b_2} + 2\left( {2a + 3} \right)x} \right]{u_x} + 2\left( {a + 1} \right)y{u_y} + a\left( {1 + a} \right)u = 0, \hfill \\
\end{gathered}  \\
  \begin{gathered}
  y\left( {1 + y} \right){u_{yy}} + xy{u_{xy}} - 2xz{u_{xz}} - \left( {2 + y} \right)z{u_{yz}} - 2{x^2}{u_{xx}} \hfill \\ \,\,\,\,\,\,\,\,\,
   + \left[ {1 - {b_1} + \left( {a + {b_2} + 1} \right)y} \right]{u_y} - \left( {a - 2{b_2} + 2} \right)x{u_x} - az{u_z} + a{b_2}u = 0, \hfill \\
\end{gathered}  \\
  \begin{gathered}
  z\left( {1 + 4z} \right){u_{zz}} + x{u_{xz}} - y\left( {1 + 4z} \right){u_{yz}} + {y^2}{u_{yy}}\hfill \\ \,\,\,\,\,\,\,\,\,
   + \left[ {1 - {b_2} + 2\left( {2{b_1} + 3} \right)z} \right]{u_z} - 2{b_1}y{u_y} + {b_1}\left( {1 + {b_1}} \right)u = 0, \hfill \\
\end{gathered}
\end{array}} \right.
$

where $u\equiv \,\,   {F_{51c}}\left( {a,{b_1},{b_2};x,y,z} \right)$.

\bigskip

\begin{equation}
{F_{52b}}\left( {{a_1},{a_2},b;c;x,y,z} \right) = \sum\limits_{m,n,p = 0}^\infty  {} \frac{{{{\left( {{a_1}} \right)}_{2m + n}}{{\left( {{a_2}} \right)}_n}{{\left( b \right)}_{2p - m}}}}{{{{\left( c \right)}_{n + p}}}}\frac{x^m}{m!}\frac{y^n}{n!}\frac{z^p}{p!},
\end{equation}

region of convergence:
$$
\begin{gathered}
  \left\{ {r < \frac{1}{4},\,\,\,s \leq \frac{1}{4} - r,\,\,\,t < \min \left\{ {\frac{1}{4},\frac{{{{\left( {1 - 4r} \right)}^2}}}{{64{r^2}}}} \right\}} \right\} \hfill \\
   \cup \left\{ {r < \frac{1}{4},\,\,\,\frac{1}{4} - r < s < 1 - 2\sqrt r ,\,\,\,t < \min \left\{ {\frac{1}{4},\frac{s}{{1 - s}}{\Phi _2}\left( {\frac{r}{{{{\left( {1 - s} \right)}^2}}}} \right)} \right\}} \right\}. \hfill \\
\end{gathered}
$$

System of partial differential equations:

$
\left\{ {\begin{array}{*{20}{l}}
  \begin{gathered}
  x\left( {1 + 4x} \right){u_{xx}} + 4xy{u_{xy}} - 2z{u_{xz}} + {y^2}{u_{yy}} \hfill \\\,\,\,\,\,\,\,\,\,
   + \left[ {1 - b + 2\left( {2{a_1} + 3} \right)x} \right]{u_x} + 2\left( {{a_1} + 1} \right)y{u_y} + {a_1}\left( {1 + {a_1}} \right)u = 0, \hfill \\
\end{gathered}  \\
  y\left( {1 - y} \right){u_{yy}} - 2xy{u_{xy}} + z{u_{yz}} + \left[ {c - \left( {{a_1} + {a_2} + 1} \right)y} \right]{u_y}
  - 2{a_2}x{u_x} - {a_1}{a_2}u = 0, \\
  \begin{gathered}
  z\left( {1 - 4z} \right){u_{zz}} + 4xz{u_{xz}} + y{u_{yz}} - {x^2}{u_{xx}}
   + \left[ {c - 2\left( {2b + 3} \right)z} \right]{u_z} + 2bx{u_x} - b\left( {1 + b} \right)u = 0, \hfill \\
\end{gathered}
\end{array}} \right.
$

where $u\equiv \,\,   {F_{52b}}\left( {{a_1},{a_2},b;c;x,y,z} \right)$.

\bigskip

\begin{equation}
{F_{53b}}\left( {{a_1},{a_2},b;{c_1},{c_2};x,y,z} \right) = \sum\limits_{m,n,p = 0}^\infty  {} \frac{{{{\left( {{a_1}} \right)}_p}{{\left( {{a_2}} \right)}_p}{{\left( b \right)}_{2m + 2n - p}}}}{{{{\left( {{c_1}} \right)}_m}{{\left( {{c_2}} \right)}_n}}}\frac{x^m}{m!}\frac{y^n}{n!}\frac{z^p}{p!},
\end{equation}

region of convergence:
$$
\left\{ {\sqrt r  + \sqrt s  < \frac{1}{2},\,\,\,t < \frac{1}{{1 + 2\left( {\sqrt r  + \sqrt s } \right)}}} \right\}.
$$

System of partial differential equations:

$
\left\{ {\begin{array}{*{20}{l}}
  \begin{gathered}
  x\left( {1 - 4x} \right){u_{xx}} - 8xy{u_{xy}} + 4xz{u_{xz}} + 4yz{u_{yz}} - 4{y^2}{u_{yy}} - {z^2}{u_{zz}} \hfill \\
 \,\,\,\,\,\,\,\,\,  + \left[ {{c_1} - 2\left( {2b + 3} \right)x} \right]{u_x} - 2\left( {2b + 3} \right)y{u_y} + 2bz{u_z} - b\left( {1 + b} \right)u = 0, \hfill \\
\end{gathered}  \\
  \begin{gathered}
  y\left( {1 - 4y} \right){u_{yy}} - 8xy{u_{xy}} + 4xz{u_{xz}} + 4yz{u_{yz}} - 4{x^2}{u_{xx}} - {z^2}{u_{zz}} \hfill \\
 \,\,\,\,\,\,\,\,\,  + \left[ {{c_2} - 2\left( {2b + 3} \right)y} \right]{u_y} - 2\left( {2b + 3} \right)x{u_x} + 2bz{u_z} - b\left( {1 + b} \right)u = 0, \hfill \\
\end{gathered}  \\
  z\left( {1 + z} \right){u_{zz}} - 2x{u_{xz}} - 2y{u_{yz}}
  + \left[ {1 - b + \left( {{a_1} + {a_2} + 1} \right)z} \right]{u_z} + {a_1}{a_2}u = 0,
\end{array}} \right.
$

where $u\equiv \,\,   {F_{53b}}\left( {{a_1},{a_2},b;{c_1},{c_2};x,y,z} \right)$.

Particular solutions:

$
{u_1} = {F_{53b}}\left( {{a_1},{a_2},b;{c_1},{c_2};x,y,z} \right),
$

$
{u_2} = {x^{1 - {c_1}}}{F_{53b}}\left( {{a_1},{a_2},2 - 2{c_1} + b;2 - {c_1},{c_2};x,y,z} \right),
$

$
{u_3} = {y^{1 - {c_2}}}{F_{53b}}\left( {{a_1},{a_2},2 - 2{c_2} + b;{c_1},2 - {c_2};x,y,z} \right),
$

$
{u_4} = {x^{1 - {c_1}}}{y^{1 - {c_2}}}{F_{53b}}\left( {{a_1},{a_2},4 - 2{c_1} - 2{c_2} + b;2 - {c_1},2 - {c_2};x,y,z} \right).
$

\bigskip

\begin{equation}
{F_{53c}}\left( {a,{b_1},{b_2};c;x,y,z} \right) = \sum\limits_{m,n,p = 0}^\infty  {} \frac{{{{\left( a \right)}_p}{{\left( {{b_1}} \right)}_{2m + 2n - p}}{{\left( {{b_2}} \right)}_{p - m}}}}{{{{\left( c \right)}_n}}}\frac{x^m}{m!}\frac{y^n}{n!}\frac{z^p}{p!},
\end{equation}

region of convergence:
$$
\left\{ {s < \frac{1}{4},\,\,\,t < \frac{1}{{1 + 2\sqrt s }},\,\,\,r < \min \left[ {\frac{{{{\left( {1 - 2\sqrt s } \right)}^2}}}{4},\frac{{1 - t - 2t\sqrt s }}{{{t^2}}}} \right]} \right\}.
$$

System of partial differential equations:

$
\left\{ {\begin{array}{*{20}{l}}
  \begin{gathered}
  x\left( {1 + 4x} \right){u_{xx}} + 8xy{u_{xy}} - \left( {1 + 4x} \right)z{u_{xz}} - 4yz{u_{yz}} + 4{y^2}{u_{yy}}\hfill \\
\,\,\,\,\,\,\,\,\,   + {z^2}{u_{zz}} + \left[ {1 - {b_2} + 2\left( {2{b_1} + 3} \right)x} \right]{u_x} + 2\left( {2{b_1} + 3} \right)y{u_y}
   - 2{b_1}z{u_z} + {b_1}\left( {1 + {b_1}} \right)u = 0, \hfill \\
\end{gathered}  \\
  \begin{gathered}
  y\left( {1 - 4y} \right){u_{yy}} - 8xy{u_{xy}} + 4xz{u_{xz}} + 4yz{u_{yz}} - 4{x^2}{u_{xx}} \hfill \\
\,\,\,\,\,\,\,\,\,  - {z^2}{u_{zz}} +  \left[ {c - 2\left( {2{b_1} + 3} \right)y} \right]{u_y} - 4{b_1}x{u_x} - 6x{u_x}
   + 2{b_1}z{u_z} - {b_1}\left( {1 + {b_1}} \right)u = 0, \hfill \\
\end{gathered}  \\
  \begin{gathered}
  z\left( {1 + z} \right){u_{zz}} - x\left( {2 + z} \right){u_{xz}} - 2y{u_{yz}}
  + \left[ {1 - {b_1} + \left( {a + {b_2} + 1} \right)z} \right]{u_z} - ax{u_x} + a{b_2}u = 0, \hfill \\
\end{gathered}
\end{array}} \right.
$

where $u\equiv \,\,   {F_{53c}}\left( {a,{b_1},{b_2};c;x,y,z} \right)$.

Particular solutions:

$
{u_1} = {F_{53c}}\left( {a,{b_1},{b_2};c;x,y,z} \right),
$

$
{u_2} = {y^{1 - c}}{F_{53c}}\left( {a,2 - 2c + {b_1},{b_2};2 - c;x,y,z} \right).
$

\bigskip

\begin{equation}
{F_{53d}}\left( {{b_1},{b_2},{b_3};x,y,z} \right) = \sum\limits_{m,n,p = 0}^\infty  {} {{{{\left( {{b_1}} \right)}_{2m + 2n - p}}{{\left( {{b_2}} \right)}_{p - m}}{{\left( {{b_3}} \right)}_{p - n}}}}\frac{x^m}{m!}\frac{y^n}{n!}\frac{z^p}{p!},
\end{equation}

region of convergence:
$$
\begin{gathered}
  \left\{ {\sqrt r  + \sqrt s  < \frac{1}{2},\,\,\,t < \min \left\{ {{U^ + }\left( {{w_1}} \right),{U^ - }\left( {{w_2}} \right)} \right\}} \right\}, \hfill \\
  {P^ \pm }\left( w \right) = w \pm \sqrt {{w^2} \pm 4r}  \pm \sqrt {{w^2} \pm 4s}  \mp 1, \hfill \\
{w_1}:\,\,{\rm{the\,\,positive\,\,root\,\,of}}\,\,{P^ + }\left( w \right) = 0, \hfill \\
{w_2}:\,\,{\rm{the\,\,positive\,\,root\,\,of}}\,\,{P^ - }\left( w \right) = 0, \hfill \\
  {U^ \pm }\left( w \right) = \frac{w}{{4rs}}\left( {w \pm \sqrt {{w^2} \pm 4r} } \right)\left( {w \pm \sqrt {{w^2} \pm 4s} } \right). \hfill \\
\end{gathered}
$$

System of partial differential equations:

$
\left\{ {\begin{array}{*{20}{l}}
  \begin{gathered}
  x\left( {1 + 4x} \right){u_{xx}} + 8xy{u_{xy}} - \left( {1 + 4x} \right)z{u_{xz}} - 4yz{u_{yz}} + 4{y^2}{u_{yy}} \hfill \\
 \,\,\,\,\,\,\,\,\, + {z^2}{u_{zz}}  + \left[ {1 - {b_2} + 2\left( {2{b_1} + 3} \right)x} \right]{u_x} + 2\left( {2{b_1} + 3} \right)y{u_y}
  - 2{b_1}z{u_z} + {b_1}\left( {1 + {b_1}} \right)u = 0, \hfill \\
\end{gathered}  \\
  \begin{gathered}
  y\left( {1 + 4y} \right){u_{yy}} + 8xy{u_{xy}} - 4xz{u_{xz}} - \left( {1 + 4y} \right)z{u_{yz}} + 4{x^2}{u_{xx}}\hfill \\
\,\,\,\,\,\,\,\,\,   + {z^2}{u_{zz}}  + \left[ {1 - {b_3} + 2\left( {2{b_1} + 3} \right)y} \right]{u_y} + 2\left( {2{b_1} + 3} \right)x{u_x}
    - 2{b_1}z{u_z} + {b_1}\left( {1 + {b_1}} \right)u = 0, \hfill \\
\end{gathered}  \\
  \begin{gathered}
  z\left( {1 + z} \right){u_{zz}} + xy{u_{xy}} - x\left( {2 + z} \right){u_{xz}} - y\left( {2 + z} \right){u_{yz}} \hfill \\
 \,\,\,\,\,\,\,\,\,  + \left[ {1 - {b_1} + \left( {{b_2} + {b_3} + 1} \right)z} \right]{u_z} - {b_3}x{u_x} - {b_2}y{u_y} + {b_2}{b_3}u = 0, \hfill \\
\end{gathered}
\end{array}} \right.
$

where $u\equiv \,\,   {F_{53d}}\left( {{b_1},{b_2},{b_3};x,y,z} \right)$.

\bigskip

\begin{equation}
{F_{54a}}\left( {{a_1},{a_2};{c_1},{c_2},{c_3};x,y,z} \right) = \sum\limits_{m,n,p = 0}^\infty  {} \frac{{{{\left( {{a_1}} \right)}_{2m + n}}{{\left( {{a_2}} \right)}_{n + 2p}}}}{{{{\left( {{c_1}} \right)}_m}{{\left( {{c_2}} \right)}_n}{{\left( {{c_3}} \right)}_p}}}\frac{x^m}{m!}\frac{y^n}{n!}\frac{z^p}{p!},
\end{equation}

first appearance of this function in the literature, and old notation:
[12], $X_{12}$,

region of convergence:
$$
 \left\{ {r < \frac{1}{4},\,\,\,t < \frac{1}{4},\,\,\,s < \left( {1 - 2\sqrt r } \right)\left( {1 - 2\sqrt t } \right)} \right\}.
$$

System of partial differential equations:

$
\left\{ {\begin{array}{*{20}{l}}
  \begin{gathered}
  x\left( {1 - 4x} \right){u_{xx}} - 4xy{u_{xy}} - {y^2}{u_{yy}}  + \left[ {{c_1} - 2\left( {2{a_1} + 3} \right)x} \right]{u_x} - 2\left( {{a_1} + 1} \right)y{u_y} - {a_1}\left( {1 + {a_1}} \right)u = 0, \hfill \\
\end{gathered}  \\
  \begin{gathered}
  y\left( {1 - y} \right){u_{yy}} - 2xy{u_{xy}} - 4xz{u_{xz}} - 2yz{u_{yz}}\hfill \\ \,\,\,\,\,\,\,\,\,    + \left[ {{c_2} - \left( {{a_1} + {a_2} + 1} \right)y} \right]{u_y} - 2{a_2}x{u_x} - 2{a_1}z{u_z} - {a_1}{a_2}u = 0, \hfill \\
\end{gathered}  \\
  \begin{gathered}
  z\left( {1 - 4z} \right){u_{zz}} - 4yz{u_{yz}} - {y^2}{u_{yy}}   + \left[ {{c_3} - 2\left( {2{a_2} + 3} \right)z} \right]{u_z} - 2\left( {{a_2} + 1} \right)y{u_y} - {a_2}\left( {1 + {a_2}} \right)u = 0, \hfill \\
\end{gathered}
\end{array}} \right.
$

where $u\equiv \,\,   {F_{54a}}\left( {{a_1},{a_2};{c_1},{c_2},{c_3};x,y,z} \right)$.

Particular solutions:

$
{u_1} = {F_{54a}}\left( {{a_1},{a_2};{c_1},{c_2},{c_3};x,y,z} \right),
$

$
{u_2} = {x^{1 - {c_1}}}{F_{54a}}\left( {2 - 2{c_1} + {a_1},{a_2};2 - {c_1},{c_2},{c_3};x,y,z} \right),
$

$
{u_3} = {y^{1 - {c_2}}}{F_{54a}}\left( {1 - {c_2} + {a_1},1 - {c_2} + {a_2};{c_1},2 - {c_2},{c_3};x,y,z} \right),
$

$
{u_4} = {z^{1 - {c_3}}}{F_{54a}}\left( {{a_1},2 - 2{c_3} + {a_2};{c_1},{c_2},2 - {c_3};x,y,z} \right),
$

$
{u_5} = {x^{1 - {c_1}}}{y^{1 - {c_2}}}{F_{54a}}\left( {3 - 2{c_1} - {c_2} + {a_1},1 - {c_2} + {a_2};2 - {c_1},2 - {c_2},{c_3};x,y,z} \right),
$

$
{u_6} = {x^{1 - {c_1}}}{z^{1 - {c_3}}}{F_{54a}}\left( {2 - 2{c_1} + {a_1},2 - 2{c_3} + {a_2};2 - {c_1},{c_2},2 - {c_3};x,y,z} \right),
$

$
{u_7} = {y^{1 - {c_2}}}{z^{1 - {c_3}}}{F_{54a}}\left( {1 - {c_2} + {a_1},2 - {c_2} - {c_3} + {a_2};{c_1},2 - {c_2},2 - {c_3};x,y,z} \right),
$

$
{u_8} = {x^{1 - {c_1}}}{y^{1 - {c_2}}}{z^{1 - {c_3}}}{F_{54a}}\left( {2 - {c_1} - {c_2} + {a_1},2 - {c_2} - {c_3} + {a_2};2 - {c_1},2 - {c_2},2 - {c_3};x,y,z} \right).
$

\bigskip

\begin{equation}
{F_{54b}}\left( {a,b;{c_1},{c_2};x,y,z} \right) = \sum\limits_{m,n,p = 0}^\infty  {} \frac{{{{\left( a \right)}_{n + 2p}}{{\left( b \right)}_{2m + n - p}}}}{{{{\left( {{c_1}} \right)}_m}{{\left( {{c_2}} \right)}_n}}}\frac{x^m}{m!}\frac{y^n}{n!}\frac{z^p}{p!},
\end{equation}

region of convergence:
$$
\left\{ {s + 2\sqrt r  < 1,\,\,\,t < \min \left\{ {\frac{1}{{1 + 2\sqrt r }}{\Theta _1}\left( {\frac{s}{{1 + 2\sqrt r }}} \right),\frac{1}{{1 - 2\sqrt r }}{\Theta _2}\left( {\frac{s}{{1 - 2\sqrt r }}} \right)} \right\}} \right\}.
$$

System of partial differential equations:

$
\left\{ {\begin{array}{*{20}{l}}
  \begin{gathered}
  x\left( {1 - 4x} \right){u_{xx}} - 4xy{u_{xy}} + 4xz{u_{xz}} + 2yz{u_{yz}} - {y^2}{u_{yy}} - {z^2}{u_{zz}} \hfill \\
   \,\,\,\,\,\,\,\,\,+ \left[ {{c_1} - 2\left( {2b + 3} \right)x} \right]{u_x} - 2by{u_y} - 2y{u_y} + 2bz{u_z} - b\left( {1 + b} \right)u = 0, \hfill \\
\end{gathered}  \\
  \begin{gathered}
  y\left( {1 - y} \right){u_{yy}} - 2xy{u_{xy}} - 4xz{u_{xz}} - yz{u_{yz}} + 2{z^2}{u_{zz}} \hfill \\
  \,\,\,\,\,\,\,\,\, + \left[ {{c_2} - \left( {a + b + 1} \right)y} \right]{u_y} - 2ax{u_x} + \left( {a - 2b + 2} \right)z{u_z} - abu = 0, \hfill \\
\end{gathered}  \\
  \begin{gathered}
  z\left( {1 + 4z} \right){u_{zz}} - 2x{u_{xz}} - y\left( {1 - 4z} \right){u_{yz}} + {y^2}{u_{yy}} \hfill \\
 \,\,\,\,\,\,\,\,\, + \left[ {1 - b + 2\left( {2a + 3} \right)z} \right]{u_z} + 2\left( {a + 1} \right)y{u_y} + a\left( {1 + a} \right)u = 0, \hfill \\
\end{gathered}
\end{array}} \right.
$

where $u\equiv \,\,   {F_{54b}}\left( {a,b;{c_1},{c_2};x,y,z} \right) $.

Particular solutions:

$
{u_1} = {F_{54b}}\left( {a,b;{c_1},{c_2};x,y,z} \right),
$

$
{u_2} = {x^{1 - {c_1}}}{F_{54b}}\left( {a,2 - 2{c_1} + b;2 - {c_1},{c_2};x,y,z} \right),
$

$
{u_3} = {y^{1 - {c_2}}}{F_{54b}}\left( {1 - {c_2} + a,1 - {c_2} + b;{c_1},2 - {c_2};x,y,z} \right),
$

$
{u_4} = {x^{1 - {c_1}}}{y^{1 - {c_2}}}{F_{54b}}\left( {1 - {c_2} + a,2 - {c_1} - {c_2} + b;2 - {c_1},2 - {c_2};x,y,z} \right).
$

\bigskip

\begin{equation}
{F_{54c}}\left( {{b_1},{b_2};c;x,y,z} \right) = \sum\limits_{m,n,p = 0}^\infty  {} \frac{{{{\left( {{b_1}} \right)}_{2m + n - p}}{{\left( {{b_2}} \right)}_{n + 2p - m}}}}{{{{\left( c \right)}_n}}}\frac{x^m}{m!}\frac{y^n}{n!}\frac{z^p}{p!},
\end{equation}

region of convergence:
$$
\begin{gathered}
  \left\{ {r < {\Phi _1}\left( t \right),\,\,\,t < {\Phi _1}\left( r \right),\,\,\,\sqrt s  < \min \left\{ {U_r^ + \left( {{w_{rt}}} \right),U_t^ + \left( {{w_{tr}}} \right),U_r^ - \left( {{w_3}} \right)} \right\}} \right\}, \hfill \\
  P_{rt}^ \pm \left( w \right) = r{w^3} \pm \frac{1}{3}{w^2} + \frac{1}{3}w - t, \hfill \\
{w_{rt}}:\,\,{\rm{the\,\,root\,\,in}}\,\,\left( {0,[{{1 + A\left( r \right)}}]^{-1}} \right)\,\,{\rm{of}}\,\,P_{rt}^ + \left( w \right) = 0, \hfill \\
{w_3}:\,\,the\,\,{\rm{middle\,\,root\,\,of}}\,\,P_{rt}^ - \left( w \right) = 0, \hfill \\
  U_q^ \pm \left( w \right) = \frac{{ - 3q{w^2} \mp 2w + 1}}{{3\sqrt w }}. \hfill \\
\end{gathered}
$$

System of partial differential equations:

$
\left\{ {\begin{array}{*{20}{l}}
  \begin{gathered}
  x\left( {1 + 4x} \right){u_{xx}} - \left( {1 - 4x} \right)y{u_{xy}} - 2\left( {1 + 2x} \right)z{u_{xz}} - 2yz{u_{yz}}
  + {y^2}{u_{yy}} + {z^2}{u_{zz}} \hfill \\ \,\,\,\,\,\,\,\,\,+ \left[ {1 - {b_2} + 2\left( {2{b_1} + 3} \right)x} \right]{u_x}
   + 2\left( {{b_1} + 1} \right)y{u_y} - 2{b_1}z{u_z} + {b_1}\left( {1 + {b_1}} \right)u = 0, \hfill \\
\end{gathered}  \\
  \begin{gathered}
  y\left( {1 - y} \right){u_{yy}} - xy{u_{xy}} - 5xz{u_{xz}} - yz{u_{yz}} + 2{x^2}{u_{xx}}
  + 2{z^2}{u_{zz}} \hfill \\  \,\,\,\,\,\,\,\,\,+ \left[ {c - \left( {{b_1} + {b_2} + 1} \right)y} \right]{u_y} + \left( {{b_1} - 2{b_2} + 2} \right)x{u_x}
  + \left( {{b_2} - 2{b_1} + 2} \right)z{u_z} - {b_1}{b_2}u = 0, \hfill \\
\end{gathered}  \\
  \begin{gathered}
  z\left( {1 + 4z} \right){u_{zz}} - 2xy{u_{xy}} - 2x\left( {1 + 2z} \right){u_{xz}} - y\left( {1 - 4z} \right){u_{yz}}
  + {x^2}{u_{xx}} + {y^2}{u_{yy}}  \hfill \\\,\,\,\,\,\,\,\,\, + \left[ {1 - {b_1} + 2\left( {2{b_2} + 3} \right)z} \right]{u_z} - 2{b_2}x{u_x}
  + 2\left( {{b_2} + 1} \right)y{u_y} + {b_2}\left( {1 + {b_2}} \right)u = 0, \hfill \\
\end{gathered}
\end{array}} \right.
$

where $u\equiv \,\,   {F_{54c}}\left( {{b_1},{b_2};c;x,y,z} \right)$.

Particular solutions:

$
{u_1} = {F_{54c}}\left( {{b_1},{b_2};c;x,y,z} \right),
$

$
{u_2} = {y^{1 - c}}{F_{54c}}\left( {1 - c + {b_1},1 - c + {b_2};2 - c;x,y,z} \right).
$

\bigskip

\begin{equation}
{F_{55a}}\left( {{a_1},{a_2};{c_1},{c_2};x,y,z} \right) = \sum\limits_{m,n,p = 0}^\infty  {} \frac{{{{\left( {{a_1}} \right)}_{2m + n}}{{\left( {{a_2}} \right)}_{n + 2p}}}}{{{{\left( {{c_1}} \right)}_{m + n}}{{\left( {{c_2}} \right)}_p}}}\frac{x^m}{m!}\frac{y^n}{n!}\frac{z^p}{p!},
\end{equation}

first appearance of this function in the literature, and old notation:
 [12], $X_{10}$,

region of convergence:
$$
\left\{ {r < \frac{1}{4},\,\,\,t < \frac{1}{4},\,\,\,s < \left( {1 - 2\sqrt t } \right)\left( {\frac{1}{2} + \frac{1}{2}\sqrt {1 - 4r} } \right)} \right\}.
$$

System of partial differential equations:

$
\left\{ {\begin{array}{*{20}{l}}
  \begin{gathered}
  x\left( {1 - 4x} \right){u_{xx}} + \left( {1 - 4x} \right)y{u_{xy}} - {y^2}{u_{yy}}\hfill \\ \,\,\,\,\,\,\,\,\,
   + \left[ {{c_1} - 2\left( {2{a_1} + 3} \right)x} \right]{u_x} - 2\left( {{a_1} + 1} \right)y{u_y} - {a_1}\left( {1 + {a_1}} \right)u = 0, \hfill \\
\end{gathered}  \\
  \begin{gathered}
  y\left( {1 - y} \right){u_{yy}} + x\left( {1 - 2y} \right){u_{xy}} - 4xz{u_{xz}} \hfill \\ \,\,\,\,\,\,\,\,\,
  - 2yz{u_{yz}}
   + \left[ {{c_1} - \left( {{a_1} + {a_2} + 1} \right)y} \right]{u_y} - 2{a_2}x{u_x} - 2{a_1}z{u_z} - {a_1}{a_2}u = 0, \hfill \\
\end{gathered}  \\
  \begin{gathered}
  z\left( {1 - 4z} \right){u_{zz}} - 4yz{u_{yz}} - {y^2}{u_{yy}}
   + \left[ {{c_2} - 2\left( {2{a_2} + 3} \right)z} \right]{u_z} - 2\left( {{a_2} + 1} \right)y{u_y} - {a_2}\left( {1 + {a_2}} \right)u = 0, \hfill \\
\end{gathered}
\end{array}} \right.
$

where $u\equiv \,\,   {F_{55a}}\left( {{a_1},{a_2};{c_1},{c_2};x,y,z} \right)$.

Particular solutions:

$
{u_1} = {F_{55a}}\left( {{a_1},{a_2};{c_1},{c_2};x,y,z} \right),
$

$
{u_2} = {z^{1 - {c_2}}}{F_{55a}}\left( {{a_1},2 - 2{c_2} + {a_2};{c_1},2 - {c_2};x,y,z} \right).
$

\bigskip

\begin{equation}
{F_{55b}}\left( {a,b;c;x,y,z} \right) = \sum\limits_{m,n,p = 0}^\infty  {} \frac{{{{\left( a \right)}_{n + 2p}}{{\left( b \right)}_{2m + n - p}}}}{{{{\left( c \right)}_{m + n}}}}\frac{x^m}{m!}\frac{y^n}{n!}\frac{z^p}{p!},
\end{equation}

region of convergence:
$$
\begin{gathered}
  \left\{ {r < \frac{1}{4},\,\,\,t < \frac{1}{{4\left( {1 + 2\sqrt r } \right)}},\,\,\,s < \min \left\{ {{U^ - }\left( {{w_1}} \right),{U^ + }\left( {{w_2}} \right)} \right\}} \right\}, \hfill \\
  {P^ \pm }\left( w \right) = 3{w^4} - 4{w^3} - 2\left( {1 \pm 4t} \right){w^2} + 4\left( {1 \pm 4t} \right)w - {\left( {1 \pm 4t} \right)^2} + 64r{t^2}, \hfill \\
{w_1}:\,\,{\rm{the\,\,root\,\,in}}\,\,\left( {0,\frac{{2 - a\left( t \right)}}{3}} \right)\,\,{\rm{of}}\,\,{P^ - }\left( w \right) = 0, \hfill \\
{w_2}:\,\,{\rm{the\,\,greater\,\,root\,\,in}}\,\,\left( {0,1} \right)\,\,{\rm{of}}\,\,{P^+}(w) = 0, \hfill \\
  {U^\pm }(w) = \frac{{w\left( {1 \pm 4t - {w^2}} \right)}}{{8t}}. \hfill \\
\end{gathered}
$$

System of partial differential equations:

$
\left\{ {\begin{array}{*{20}{l}}
  \begin{gathered}
  x\left( {1 - 4x} \right){u_{xx}} + \left( {1 - 4x} \right)y{u_{xy}} + 4xz{u_{xz}} + 2yz{u_{yz}} - {y^2}{u_{yy}} - {z^2}{u_{zz}} \hfill \\
 \,\,\,\,\,\,\,\,\,  + \left[ {c - 2\left( {2b + 3} \right)x} \right]{u_x} - 2\left( {b + 1} \right)y{u_y} + 2bz{u_z} - b\left( {1 + b} \right)u = 0, \hfill \\
\end{gathered}  \\
  \begin{gathered}
  y\left( {1 - y} \right){u_{yy}} + x\left( {1 - 2y} \right){u_{xy}} - 4xz{u_{xz}} - yz{u_{yz}} + 2{z^2}{u_{zz}} \hfill \\
\,\,\,\,\,\,\,\,\,   + \left[ {c - \left( {a + b + 1} \right)y} \right]{u_y} - 2ax{u_x} + \left( {a - 2b + 2} \right)z{u_z} - abu = 0, \hfill \\
\end{gathered}  \\
  \begin{gathered}
  z\left( {1 + 4z} \right){u_{zz}} - 2x{u_{xz}} - y\left( {1 - 4z} \right){u_{yz}} + {y^2}{u_{yy}} \hfill \\
 \,\,\,\,\,\,\,\,\,  + \left[ {1 - b + 2\left( {2a + 3} \right)z} \right]{u_z} + 2\left( {a + 1} \right)y{u_y} + a\left( {1 + a} \right)u = 0, \hfill \\
\end{gathered}
\end{array}} \right.
$

where $u\equiv \,\,   {F_{55b}}\left( {a,b;c;x,y,z} \right)$.

\bigskip

\begin{equation}
{F_{56a}}\left( {{a_1},{a_2};{c_1},{c_2};x,y,z} \right) = \sum\limits_{m,n,p = 0}^\infty  {} \frac{{{{\left( {{a_1}} \right)}_{2m + n}}{{\left( {{a_2}} \right)}_{n + 2p}}}}{{{{\left( {{c_1}} \right)}_{m + p}}{{\left( {{c_2}} \right)}_n}}}\frac{x^m}{m!}\frac{y^n}{n!}\frac{z^p}{p!},
\end{equation}

first appearance of this function in the literature, and old notation:
[12], $X_{11}$,

region of convergence:
$$
\begin{gathered}
   \left\{ {r < \frac{1}{4},\,\,\,t < \frac{1}{4},\,\,\,s < \left( {1 - 2\sqrt {r{w_1}} } \right)\left( {1 - 2\sqrt {t - t{w_1}} } \right)} \right\}, \hfill \\
  P\left( w \right) = 2 - 4w + \sqrt {\frac{w}{r}}  - \sqrt {\frac{{1 - w}}{t}} , \hfill \\
{w_1}:\,\,{\rm{the\,\,root\,\,of}}\,\,P\left( w \right) = 0\,\,{\rm{in}}\,\,\left\{ {\begin{array}{*{20}{c}}
  {\left( \displaystyle {0,\frac{r}{{r + t}}} \right),\,\,r \leq t,} \\
  {\left( \displaystyle {\frac{r}{{r + t}},1} \right),\,\,r > t.}
\end{array}} \right. \hfill \\
\end{gathered}
$$

System of partial differential equations:

$
\left\{ {\begin{array}{*{20}{l}}
  \begin{gathered}
  x\left( {1 - 4x} \right){u_{xx}} - 4xy{u_{xy}} + z{u_{xz}} - {y^2}{u_{yy}}\hfill \\ \,\,\,\,\,\,\,\,\,
   + \left[ {{c_1} - 2\left( {2{a_1} + 3} \right)x} \right]{u_x} - 2\left( {{a_1} + 1} \right)y{u_y} - {a_1}\left( {1 + {a_1}} \right)u = 0, \hfill \\
\end{gathered}  \\
  \begin{gathered}
  y\left( {1 - y} \right){u_{yy}} - 2xy{u_{xy}} - 4xz{u_{xz}} - 2yz{u_{yz}}\hfill \\ \,\,\,\,\,\,\,\,\,
   + \left[ {{c_2} - \left( {{a_1} + {a_2} + 1} \right)y} \right]{u_y} - 2{a_2}x{u_x} - 2{a_1}z{u_z} - {a_1}{a_2}u = 0, \hfill \\
\end{gathered}  \\
  \begin{gathered}
  z\left( {1 - 4z} \right){u_{zz}} + x{u_{xz}} - 4yz{u_{yz}} - {y^2}{u_{yy}}\hfill \\ \,\,\,\,\,\,\,\,\,
   + \left[ {{c_1} - 2\left( {2{a_2} + 3} \right)z} \right]{u_z} - 2\left( {{a_2} + 1} \right)y{u_y} - {a_2}\left( {1 + {a_2}} \right)u = 0, \hfill \\
\end{gathered}
\end{array}} \right.
$

where $u\equiv \,\,   {F_{56a}}\left( {{a_1},{a_2};{c_1},{c_2};x,y,z} \right)$.

Particular solutions:

$
{u_1} = {F_{56a}}\left( {{a_1},{a_2};{c_1},{c_2};x,y,z} \right),
$

$
{u_2} = {y^{1 - {c_2}}}{F_{56a}}\left( {1 - {c_2} + {a_1},1 - {c_2} + {a_2};{c_1},2 - {c_2};x,y,z} \right).
$

\bigskip

\begin{equation}
 {F_{57a}}\left( {{a_1},{a_2};c;x,y,z} \right) = \sum\limits_{m,n,p = 0}^\infty  {} \frac{{{{\left( {{a_1}} \right)}_{2m + n}}{{\left( {{a_2}} \right)}_{n + 2p}}}}{{{{\left( c \right)}_{m + n + p}}}}\frac{x^m}{m!}\frac{y^n}{n!}\frac{z^p}{p!},
\end{equation}

first appearance of this function in the literature, and old notation:
[12], $X_{9}$,

region of convergence:
$$
\left\{ {r < \frac{1}{4},\,\,\,t < \frac{1}{4},\,\,\,s < \frac{1}{2} + \frac{1}{2}\sqrt {\left( {1 - 4r} \right)\left( {1 - 4t} \right)} } \right\}.
$$

System of partial differential equations:

$
\left\{ {\begin{array}{*{20}{l}}
  \begin{gathered}
  x\left( {1 - 4x} \right){u_{xx}} + \left( {1 - 4x} \right)y{u_{xy}} + z{u_{xz}} - {y^2}{u_{yy}}
  \hfill \\ \,\,\,\,\,\,\,\,\, + \left[ {c - 2\left( {2{a_1} + 3} \right)x} \right]{u_x} - 2\left( {{a_1} + 1} \right)y{u_y} - {a_1}\left( {1 + {a_1}} \right)u = 0, \hfill \\
\end{gathered}  \\
  \begin{gathered}
  y\left( {1 - y} \right){u_{yy}} + x\left( {1 - 2y} \right){u_{xy}} - 4xz{u_{xz}} + \left( {1 - 2y} \right)z{u_{yz}}
 \hfill \\ \,\,\,\,\,\,\,\,\,  + \left[ {c - \left( {{a_1} + {a_2} + 1} \right)y} \right]{u_y} - 2{a_2}x{u_x} - 2{a_1}z{u_z} - {a_1}{a_2}u = 0, \hfill \\
\end{gathered}  \\
  \begin{gathered}
  z\left( {1 - 4z} \right){u_{zz}} + x{u_{xz}} + y\left( {1 - 4z} \right){u_{yz}} - {y^2}{u_{yy}}
 \hfill \\  \,\,\,\,\,\,\,\,\, + \left[ {c - 2\left( {2{a_2} + 3} \right)z} \right]{u_z} - 2\left( {{a_2} + 1} \right)y{u_y} - {a_2}\left( {1 + {a_2}} \right)u = 0, \hfill \\
\end{gathered}
\end{array}} \right.
$

where $u\equiv \,\,   {F_{57a}}\left( {{a_1},{a_2};c;x,y,z} \right)$.

\bigskip

\begin{equation}
 {F_{58a}}\left( {{a_1},{a_2};{c_1},{c_2},{c_3};x,y,z} \right) = \sum\limits_{m,n,p = 0}^\infty  {} \frac{{{{\left( {{a_1}} \right)}_{2m + 2n + p}}{{\left( {{a_2}} \right)}_p}}}{{{{\left( {{c_1}} \right)}_m}{{\left( {{c_2}} \right)}_n}{{\left( {{c_3}} \right)}_p}}}\frac{x^m}{m!}\frac{y^n}{n!}\frac{z^p}{p!},
\end{equation}

first appearance of this function in the literature, and old notation:
[11], ${}^{(2)}H_{4}^{(3)}$; [12], $X_{2}$,

region of convergence:
$$
\left\{ {2 {\sqrt r  + \sqrt s }  + t < 1} \right\}.
$$

System of partial differential equations:

$
\left\{ {\begin{array}{*{20}{l}}
  \begin{gathered}
  x\left( {1 - 4x} \right){u_{xx}} - 8xy{u_{xy}} - 4xz{u_{xz}} - 4yz{u_{yz}}
  - 4{y^2}{u_{yy}} - {z^2}{u_{zz}}  \hfill \\ \,\,\,\,\,\,\,\,\, + \left[ {{c_1} - 2\left( {2{a_1} + 3} \right)x} \right]{u_x}
  - 2\left( {2{a_1} + 3} \right)y{u_y} - 2\left( {{a_1} + 1} \right)z{u_z} - {a_1}\left( {1 + {a_1}} \right)u = 0 \hfill \\
\end{gathered}  \\
  \begin{gathered}
  y\left( {1 - 4y} \right){u_{yy}} - 8xy{u_{xy}} - 4xz{u_{xz}} - 4yz{u_{yz}}
   - 4{x^2}{u_{xx}} - {z^2}{u_{zz}} \hfill \\ \,\,\,\,\,\,\,\,\, + \left[ {{c_2} - 2\left( {2{a_1} + 3} \right)y} \right]{u_y}
   - 2\left( {2{a_1} + 3} \right)x{u_x} - 2\left( {{a_1} + 1} \right)z{u_z} - {a_1}\left( {1 + {a_1}} \right)u = 0 \hfill \\
\end{gathered}  \\
  \begin{gathered}
  z\left( {1 - z} \right){u_{zz}} - 2xz{u_{xz}} - 2yz{u_{yz}}
   + \left[ {{c_3} - \left( {{a_1} + {a_2} + 1} \right)z} \right]{u_z} - 2{a_2}x{u_x} - 2{a_2}y{u_y} - {a_1}{a_2}u = 0, \hfill \\
\end{gathered}
\end{array}} \right.
$

where $u\equiv \,\,   {F_{58a}}\left( {{a_1},{a_2};{c_1},{c_2},{c_3};x,y,z} \right)$.

Particular solutions:

$
{u_1} = {F_{58a}}\left( {{a_1},{a_2};{c_1},{c_2},{c_3};x,y,z} \right),
$

$
{u_2} = {x^{1 - {c_1}}}{F_{58a}}\left( {2 - 2{c_1} + {a_1},{a_2};2 - {c_1},{c_2},{c_3};x,y,z} \right),
$

$
{u_3} = {y^{1 - {c_2}}}{F_{58a}}\left( {2 - 2{c_2} + {a_1},{a_2};{c_1},2 - {c_2},{c_3};x,y,z} \right),
$

$
{u_4} = {z^{1 - {c_3}}}{F_{58a}}\left( {1 - {c_3} + {a_1},1 - {c_3} + {a_2};{c_1},{c_2},2 - {c_3};x,y,z} \right),
$

$
{u_5} = {x^{1 - {c_1}}}{y^{1 - {c_2}}}{F_{58a}}\left( {4 - 2{c_1} - 2{c_2} + {a_1},{a_2};2 - {c_1},2 - {c_2},{c_3};x,y,z} \right),
$

$
{u_6} = {x^{1 - {c_1}}}{z^{1 - {c_3}}}{F_{58a}}\left( {3 - 2{c_1} - {c_3} + {a_1},1 - {c_3} + {a_2};2 - {c_1},{c_2},2 - {c_3};x,y,z} \right),
$

$
{u_7} = {y^{1 - {c_2}}}{z^{1 - {c_3}}}{F_{58a}}\left( {3 - 2{c_2} - {c_3} + {a_1},1 - {c_3} + {a_2};{c_1},2 - {c_2},2 - {c_3};x,y,z} \right),
$

$
{u_8} = {x^{1 - {c_1}}}{y^{1 - {c_2}}}{z^{1 - {c_3}}}
{F_{58a}}\left( {5 - 2{c_1} - 2{c_2} - {c_3} + {a_1},1 - {c_3} + {a_2};2 - {c_1},2 - {c_2},2 - {c_3};x,y,z} \right). $

\bigskip

\begin{equation}
{F_{58b}}\left( {a,b;{c_1},{c_2};x,y,z} \right) = \sum\limits_{m,n,p = 0}^\infty  {} \frac{{{{\left( a \right)}_{2m + 2n + p}}{{\left( b \right)}_{p - m}}}}{{{{\left( {{c_1}} \right)}_n}{{\left( {{c_2}} \right)}_p}}}\frac{x^m}{m!}\frac{y^n}{n!}\frac{z^p}{p!},
\end{equation}

region of convergence:
$$
\left\{ {\sqrt r  + \sqrt s  < \frac{1}{2},\,\,\,t < \left( {1 - 2\sqrt s } \right)\min \left\{ {{\Psi _1}\left( {\frac{r}{{{{\left( {1 - 2\sqrt s } \right)}^2}}}} \right),{\Psi _2}\left( {\frac{r}{{{{\left( {1 - 2\sqrt s } \right)}^2}}}} \right)} \right\}} \right\}.
$$

System of partial differential equations:

$
\left\{ {\begin{array}{*{20}{l}}
  \begin{gathered}
  x\left( {1 + 4x} \right){u_{xx}} + 8xy{u_{xy}} - \left( {1 - 4x} \right)z{u_{xz}} + 4yz{u_{yz}}
   + 4{y^2}{u_{yy}} + {z^2}{u_{zz}} \hfill \\\,\,\,\,\,\,\,\,\,+ \left[ {1 - b + 2\left( {2a + 3} \right)x} \right]{u_x} + 2\left( {2a + 3} \right)y{u_y}
    + 2\left( {a + 1} \right)z{u_z} + a\left( {1 + a} \right)u = 0, \hfill \\
\end{gathered}  \\
  \begin{gathered}
  y\left( {1 - 4y} \right){u_{yy}} - 8xy{u_{xy}} - 4xz{u_{xz}} - 4yz{u_{yz}} - 4{x^2}{u_{xx}}
  - {z^2}{u_{zz}}   \hfill \\ \,\,\,\,\,\,\,\,\,+  \left[ {{c_1} - 2\left( {2a + 3} \right)y} \right]{u_y} - 4ax{u_x} - 6x{u_x}
  - 2az{u_z} - 2z{u_z} - a\left( {1 + a} \right)u = 0, \hfill \\
\end{gathered}  \\
  \begin{gathered}
  z\left( {1 - z} \right){u_{zz}} + 2xy{u_{xy}} - xz{u_{xz}} - 2yz{u_{yz}} + 2{x^2}{u_{xx}} \hfill \\
 \,\,\,\,\,\,\,\,\,  + \left[ {{c_2} - \left( {a + b + 1} \right)z} \right]{u_z} + \left( {a - 2b + 2} \right)x{u_x}
   - 2by{u_y} - abu = 0, \hfill \\
\end{gathered}
\end{array}} \right.
$

where $u\equiv \,\,   {F_{58b}}\left( {a,b;{c_1},{c_2};x,y,z} \right)$.

Particular solutions:

$
{u_1} = {F_{58b}}\left( {a,b;{c_1},{c_2};x,y,z} \right),\
$

$
{u_2} = {y^{1 - {c_1}}}{F_{58b}}\left( {2 - 2{c_1} + a,b;2 - {c_1},{c_2};x,y,z} \right),
$

$
{u_3} = {z^{1 - {c_2}}}{F_{58b}}\left( {1 - {c_2} + a,1 - {c_2} + b;{c_1},2 - {c_2};x,y,z} \right),
$

$
{u_4} = {y^{1 - {c_1}}}{z^{1 - {c_2}}}{F_{58b}}\left( {3 - 2{c_1} - {c_2} + a,1 - {c_2} + b;2 - {c_1},2 - {c_2};x,y,z} \right).
$

\bigskip

\begin{equation}
{F_{59a}}\left( {{a_1},{a_2};{c_1},{c_2};x,y,z} \right) = \sum\limits_{m,n,p = 0}^\infty  {} \frac{{{{\left( {{a_1}} \right)}_{2m + 2n - p}}{{\left( {{a_2}} \right)}_p}}}{{{{\left( {{c_1}} \right)}_{n + p}}{{\left( {{c_2}} \right)}_m}}}\frac{x^m}{m!}\frac{y^n}{n!}\frac{z^p}{p!},
\end{equation}

first appearance of this function in the literature, and old notation:
[12], $X_{1}$,

region of convergence:
$$
 \left\{ {\sqrt r  + \sqrt s  < \frac{1}{2},\,\,\,t < \frac{1}{2}\left( {1 - 2\sqrt r } \right) + \frac{1}{2}\sqrt {{{\left( {1 - 2\sqrt r } \right)}^2} - 4s} } \right\}.
$$

System of partial differential equations:

$
\left\{ {\begin{array}{*{20}{l}}
  \begin{gathered}
  x\left( {1 - 4x} \right){u_{xx}} - 8xy{u_{xy}} - 4xz{u_{xz}} - 4yz{u_{yz}} - 4{y^2}{u_{yy}}
  - {z^2}{u_{zz}}  \hfill \\ \,\,\,\,\,\,\,\,\,  +   \left[ {{c_2} - 2\left( {2{a_1} + 3} \right)x} \right]{u_x} - 2\left( {2{a_1} + 3} \right)y{u_y}
  - 2\left( {{a_1} + 1} \right)z{u_z} - {a_1}\left( {1 + {a_1}} \right)u = 0, \hfill \\
\end{gathered}  \\
  \begin{gathered}
  y\left( {1 - 4y} \right){u_{yy}} - 8xy{u_{xy}} - 4xz{u_{xz}} + \left( {1 - 4y} \right)z{u_{yz}} - 4{x^2}{u_{xx}}
  - {z^2}{u_{zz}} \hfill \\ \,\,\,\,\,\,\,\,\,
  + \left[ {{c_1} - 2\left( {2{a_1} + 3} \right)y} \right]{u_y} - 2\left( {2{a_1} + 3} \right)x{u_x}
   - 2\left( {{a_1} + 1} \right)z{u_z} - {a_1}\left( {1 + {a_1}} \right)u = 0, \hfill \\
\end{gathered}  \\
  \begin{gathered}
  z\left( {1 - z} \right){u_{zz}} - 2xz{u_{xz}} + y\left( {1 - 2z} \right){u_{yz}}\hfill \\ \,\,\,\,\,\,\,\,\,
   + \left[ {{c_1} - \left( {{a_1} + {a_2} + 1} \right)z} \right]{u_z} - 2{a_2}x{u_x} - 2{a_2}y{u_y} - {a_1}{a_2}u = 0, \hfill \\
\end{gathered}
\end{array}} \right.
$

where $u\equiv \,\,   {F_{59a}}\left( {{a_1},{a_2};{c_1},{c_2};x,y,z} \right)$.

Particular solutions:

$
{u_1} = {F_{59a}}\left( {{a_1},{a_2};{c_1},{c_2};x,y,z} \right),
$

$
{u_2} = {x^{1 - {c_2}}}{F_{59a}}\left( {1 - {c_2} + {a_1},1 - {c_2} + {a_2};{c_1},2 - {c_2};x,y,z} \right).
$

\bigskip

\begin{equation}
 {F_{59b}}\left( {a,b;c;x,y,z} \right) = \sum\limits_{m,n,p = 0}^\infty  {} \frac{{{{\left( a \right)}_{2m + 2n + p}}{{\left( b \right)}_{p - m}}}}{{{{\left( c \right)}_{n + p}}}}\frac{x^m}{m!}\frac{y^n}{n!}\frac{z^p}{p!},
\end{equation}

region of convergence:
$$
\begin{gathered}
   \left\{ {s < \frac{1}{4},\,\,\,t < \frac{1}{2} + \frac{1}{2}\sqrt {1 - 4s} ,\,\,\,r < \min \left\{ {U\left( {{w_1}} \right),U\left( {{w_2}} \right)} \right\}} \right\}, \hfill \\
  {P^ \pm }\left( w \right) = t{w^3} + \left( {t \pm 2s} \right){w^2} - 3stw \mp 2{s^2}, \hfill \\
{w_1}:\,\,{\rm{the\,\,root\,\,in}}\,\,\left( {\frac{1}{2} - \frac{1}{2}\sqrt {1 - 4s} ,\sqrt s } \right)\,\,{\rm{of}}\,\,{P^ + }\left( w \right) = 0, \hfill \\
{w_2}:\,\,{\rm{the\,\,smaller\,\,root\,\,in}}\,\,\left( \displaystyle {\frac{1}{2} - \frac{1}{2}\sqrt {1 - 4s} ,\sqrt s } \right)\,\,{\rm{of}}\,\,{P^ - }\left( w \right) = 0,\,\,t>\frac{1}{2}, \hfill \\
  U\left( w \right) = \frac{{ s\left( { - {w^2} + w - s} \right)\left( {s - {w^2}} \right)}}{{2t{w^3}}}. \hfill \\
\end{gathered}
$$

System of partial differential equations:

$
\left\{ {\begin{array}{*{20}{l}}
  \begin{gathered}
  x\left( {1 + 4x} \right){u_{xx}} + 8xy{u_{xy}} - \left( {1 - 4x} \right)z{u_{xz}} + 4yz{u_{yz}}
   + 4{y^2}{u_{yy}} + {z^2}{u_{zz}} \hfill \\  \,\,\,\,\,\,\,\,\,+ \left[ {1 - b + 2\left( {2a + 3} \right)x} \right]{u_x} + 2\left( {2a + 3} \right)y{u_y}
    + 2\left( {a + 1} \right)z{u_z} + a\left( {1 + a} \right)u = 0, \hfill \\
\end{gathered}  \\
  \begin{gathered}
  y\left( {1 - 4y} \right){u_{yy}} - 8xy{u_{xy}} - 4xz{u_{xz}} + \left( {1 - 4y} \right)z{u_{yz}}
  - 4{x^2}{u_{xx}} - {z^2}{u_{zz}}  \hfill \\ \,\,\,\,\,\,\,\,\,+ \left[ {c - 2\left( {2a + 3} \right)y} \right]{u_y} - 2\left( {2a + 3} \right)x{u_x}
  - 2\left( {a + 1} \right)z{u_z} - a\left( {1 + a} \right)u = 0, \hfill \\
\end{gathered}  \\
  \begin{gathered}
  z\left( {1 - z} \right){u_{zz}} + 2xy{u_{xy}} - xz{u_{xz}} + y\left( {1 - 2z} \right){u_{yz}} + 2{x^2}{u_{xx}} \hfill \\
 \,\,\,\,\,\,\,\,\, + \left[ {c - \left( {a + b + 1} \right)z} \right]{u_z} + \left( {a - 2b + 2} \right)x{u_x} - 2by{u_y} - abu = 0, \hfill \\
\end{gathered}
\end{array}} \right.
$

where $u\equiv \,\,   {F_{59b}}\left( {a,b;c;x,y,z} \right)$.

\bigskip

\begin{equation}
{F_{60b}}\left( {{b_1},{b_2},{b_3};x,y,z} \right) = \sum\limits_{m,n,p = 0}^\infty  {} {{{{\left( {{b_1}} \right)}_{2m - n}}{{\left( {{b_2}} \right)}_{2n - p}}{{\left( {{b_3}} \right)}_{2p - m}}}}\frac{x^m}{m!}\frac{y^n}{n!}\frac{z^p}{p!},
\end{equation}

region of convergence:
$$
\begin{gathered}
     {{E_{rst}} \cap {E_{str}} \cap {E_{trs}}},
  {E_{rst}} = \left\{ {r < \frac{1}{4},\,\, t < U\left( {r,s} \right)} \right\}, \hfill \\
  {P_{rs}}\left( w \right) = 49s{w^4} - \left( {8 + 56s} \right){w^3} + \left( {14rs + 16s + 4} \right){w^2} - 8rsw + {r^2}s,\,\,r < \frac{1}{4}, \hfill \\
{w_{rs}}:\,\,{\rm{the\,\,root\,\,in}}\,\,\left( {\frac{1}{2},\frac{{2 + \sqrt {4 - 7r} }}{7}} \right)\,\,{\rm{of}}\,\,{P_{rs}}\left( w \right) = 0, \hfill \\
  U\left( {r,s} \right) = \frac{{{w_{rs}}\left( {w_{rs}^2 - r} \right)\sqrt {2{w_{rs}} - 1} }}{{2{r^2}\sqrt s }}. \hfill \\
\end{gathered}
$$

System of partial differential equations:

$
\left\{ {\begin{array}{*{20}{l}}
  \begin{gathered}
  x\left( {1 + 4x} \right){u_{xx}} - 4xy{u_{xy}} - 2z{u_{xz}} + {y^2}{u_{yy}}  \hfill \\ \,\,\,\,\,\,\,\,\, + \left[ {1 - {b_3} + 2\left( {2{b_1} + 3} \right)x} \right]{u_x} - 2{b_1}y{u_y} + {b_1}\left( {1 + {b_1}} \right)u = 0, \hfill \\
\end{gathered}  \\
  \begin{gathered}
  y\left( {1 + 4y} \right){u_{yy}} - 2x{u_{xy}} - 4yz{u_{yz}} + {z^2}{u_{zz}} \hfill \\ \,\,\,\,\,\,\,\,\,  + \left[ {1 - {b_1} + 2\left( {2{b_2} + 3} \right)y} \right]{u_y} - 2{b_2}z{u_z} + {b_2}\left( {1 + {b_2}} \right)u = 0, \hfill \\
\end{gathered}  \\
  \begin{gathered}
  z\left( {1 + 4z} \right){u_{zz}} - 4xz{u_{xz}} - 2y{u_{yz}} + {x^2}{u_{xx}}\hfill \\ \,\,\,\,\,\,\,\,\,
   + \left[ {1 - {b_2} + 2\left( {2{b_3} + 3} \right)z} \right]{u_z} - 2{b_3}x{u_x} + {b_3}\left( {1 + {b_3}} \right)u = 0, \hfill \\
\end{gathered}
\end{array}} \right.
$

where $u\equiv \,\,   {F_{60b}}\left( {{b_1},{b_2},{b_3};x,y,z} \right)$.

\bigskip

\begin{equation}
{F_{61b}}\left( {{b_1},{b_2};c;x,y,z} \right) = \sum\limits_{m,n,p = 0}^\infty  {} \frac{{{{\left( {{b_1}} \right)}_{2m + 2n - p}}{{\left( {{b_2}} \right)}_{2p - m}}}}{{{{\left( c \right)}_n}}}\frac{x^m}{m!}\frac{y^n}{n!}\frac{z^p}{p!},
\end{equation}

region of convergence:
$$
\begin{gathered}
    \left\{ {\sqrt r  + \sqrt s  < \frac{1}{2},\,\,\,t < \min \left\{ {\frac{1}{{1 + 2\sqrt s }}{\Phi _1}\left( \frac{r}{{{{\left( {1 + 2\sqrt s } \right)}^2}}} \right),\frac{1}{{1 - 2\sqrt s }}{\Phi _2}\left( \frac{r}{{{{\left( {1 - 2\sqrt s } \right)}^2}}}  \right)} \right\}} \right\}. \hfill \\
  \end{gathered}
$$

System of partial differential equations:

$
\left\{ {\begin{array}{*{20}{l}}
  \begin{gathered}
  x\left( {1 + 4x} \right){u_{xx}} + 8xy{u_{xy}} - 2\left( {1 + 2x} \right)z{u_{xz}} - 4yz{u_{yz}}
  + 4{y^2}{u_{yy}} + {z^2}{u_{zz}}  \hfill \\ \,\,\,\,\,\,\,\,\,
  + \left[ {1 - {b_2} + 2\left( {2{b_1} + 3} \right)x} \right]{u_x}
   + 2\left( {2{b_1} + 3} \right)y{u_y} - 2{b_1}z{u_z} + {b_1}\left( {1 + {b_1}} \right)u = 0, \hfill \\
\end{gathered}  \\
  \begin{gathered}
  y\left( {1 - 4y} \right){u_{yy}} - 8xy{u_{xy}} + 4xz{u_{xz}} + 4yz{u_{yz}}
  - 4{x^2}{u_{xx}} - {z^2}{u_{zz}}  \hfill \\ \,\,\,\,\,\,\,\,\,
  + \left[ {c - 2\left( {2{b_1} + 3} \right)y} \right]{u_y}
  - 2\left( {2{b_1} + 3} \right)x{u_x} + 2{b_1}z{u_z} - {b_1}\left( {{b_1} + 1} \right)u = 0, \hfill \\
\end{gathered}  \\
  \begin{gathered}
  z\left( {1 + 4z} \right){u_{zz}} - 2\left( {1 + 2z} \right)x{u_{xz}} - 2y{u_{yz}} + {x^2}{u_{xx}} \hfill \\
  \,\,\,\,\,\,\,\,\, + \left[ {1 - {b_1} + 2\left( {2{b_2} + 3} \right)z} \right]{u_z} - 2{b_2}x{u_x} + {b_2}\left( {{b_2} + 1} \right)u = 0, \hfill \\
\end{gathered}
\end{array}} \right.
$

where $u\equiv \,\,   {F_{61b}}\left( {{b_1},{b_2};c;x,y,z} \right)$.

Particular solutions:

$
{u_1} = {F_{61b}}\left( {{b_1},{b_2};c;x,y,z} \right),
$

$
{u_2} = {y^{1 - c}}{F_{61b}}\left( {2 - 2c + {b_1},{b_2};2 - c;x,y,z} \right).$

\bigskip

\section{Systems of partial differential equations, associated with the confluent hypergeometric functions in three variables} \label{S8}

Consider  the following hypergeometric function defined in [7]:

\begin{equation} \label{e1function}
{\rm{E}}_{1}\left(a_1,a_2,a_3,a_4,a_5;c;x,y,z\right)=\sum\limits_{m,n,p=0}^\infty\frac{(a_1)_m(a_2)_m(a_3)_n(a_4)_n(a_5)_p}{(c)_{m+n+p}}\frac{x^m}{m!}\frac{y^n}{n!}\frac{z^p}{p!},
\end{equation}

first appearance of this function in the literature: [21],\, ${}_3\Phi_B^{(1)}$, \,\,[9], \,\,$F_{B1}$,

region of convergence:
$$ \left\{ \frac{1}{r}+\frac{1}{s}>1, \,\,\,t<\infty
\right\}.
$$

To give an example, we construct a system of partial differential equations corresponding to the  confluent hypergeometric function ${\rm{E}}_{1}$ defined in (296).

Assuming
\[
A(m,n,p)=\frac{(a_1)_m(a_2)_m(a_3)_n(a_4)_n(a_5)_p}{m!n!p!(c)_{m+n+p}},
\]
we define
\[
A(m+1,n,p)=\frac{(a_1)_{m+1}(a_2)_{m+1}(a_3)_n(a_4)_n(a_5)_p}{(m+1)!n!p!(c)_{m+n+p+1}},
\]
\[
A(m,n+1,p)=\frac{(a_1)_m(a_2)_m(a_3)_{n+1}(a_4)_{n+1}(a_5)_p}{m!(n+1)!p!(c)_{m+n+p+1}},
\]
\[
A(m,n,p+1)=\frac{(a_1)_m(a_2)_m(a_3)_n(a_4)_n(a_5)_{p+1}}{m!n!(p+1)!(c)_{m+n+p}}.
\]
By virtue of a simple property of the Pochhammer symbol  (2)  in the form
\[
(\lambda)_{\nu+1}=(\lambda)_\nu(\lambda+\nu),
\]
it is easy to determine the following three functions from (87)

\begin{equation}
\label{otn77}
f(m,n,p)=\frac{\left(a_1+m\right)\left(a_2+m\right)}{(m+1)\left(c+m+n+p\right)},\,
\end{equation}
\begin{equation}
\label{otn88}
g(m,n,p)=\frac{\left(a_3+n\right)\left(a_4+n\right)}{(n+1)\left(c+m+n+p\right)},\,
\end{equation}
\begin{equation}
\label{otn99}
h(m,n,p)=\frac{a_5+p}{(p+1)\left(c+m+n+p\right)}.
\end{equation}
Now, according to the definition
(98), by virtue of  (297)--(299), we have
\begin{equation}
\label{otn177}
F(m,n,p)=\left(a_1+m\right)\left(a_2+m\right),\,\,\,\,\,\,F'(m,n,p)={(m+1)\left(c+m+n+p\right)},\,
\end{equation}

\begin{equation}
\label{otn188}
G(m,n,p)=\left(a_3+n\right)\left(a_4+n\right),\,\,\,\,\,\,G'(m,n,p)=(n+1)\left(c+m+n+p\right), ,\,
\end{equation}

\begin{equation}
\label{otn199}
H(m,n,p)=a_5+p, \,\,\,\,\,\,\,\,\,\,\,\,\,\,\,H'(m,n,p)=(p+1)\left(c+m+n+p\right).
\end{equation}

Using (300), we will compose the first equation of the system (90):
\[
 {F'\left( {{\delta_x},{\delta_y}, \delta_z} \right) \left(\frac{u}{x}\right)- F\left( {{\delta_x},{\delta_y}, \delta_z} \right)}u = 0.
\]
Taking into account the definition (89) of the operators $\delta_x$, $\delta_y$ and $\delta_z$, one can write
\[
\left(c+x\frac{\partial}{\partial x}+y\frac{\partial}{\partial y}+z \frac{\partial}{\partial z}\right)\left(1+x\frac{\partial}{\partial x}\right)\left(\frac{u}{x}\right)-\left(a_1+x\frac{\partial}{\partial x}\right)\left(a_2+x\frac{\partial}{\partial x}\right)u=0.
\]
Expanding the brackets in the last equation, we get
\begin{equation} \label{first}
x(1-x)u_{xx}+yu_{xy}+zu_{xz}+\left[c-\left(a_1+a_2+1\right)x\right]u_x-a_1a_2u=0.
\end{equation}

Using (301), we will compose the second equation of the system (90):
\[
 {G'\left( {{\delta_x},{\delta_y}, \delta_z} \right) \left(\frac{u}{y}\right)- G\left( {{\delta_x},{\delta_y}, \delta_z} \right)}u = 0,
\]
which can easily be written in the form
\[
\left(c+x\frac{\partial}{\partial x}+y\frac{\partial}{\partial y}+z \frac{\partial}{\partial z}\right)\frac{\partial u}{\partial y} -\left(a_3+y\frac{\partial}{\partial y}\right)\left(a_4+y\frac{\partial}{\partial y}\right)u=0.
\]
Having performed the necessary transformations in the last equation, we obtain the second equation of the system (90):
\begin{equation} \label{second}
y(1-y)u_{yy}+xu_{xy}+zu_{yz}+\left[c-\left(a_3+a_4+1\right)y\right]u_y-a_3a_4u=0.
\end{equation}

Finally, using (302), we will compose the third equation of the system (90):
\[
 {H'\left( {{\delta_x},{\delta_y}, \delta_z} \right) \left(\frac{u}{z}\right)- H\left( {{\delta_x},{\delta_y}, \delta_z} \right)}u = 0,
\]
or
\[
\left(c+x\frac{\partial}{\partial x}+y\frac{\partial}{\partial y}+z \frac{\partial}{\partial z}\right)\frac{\partial u}{\partial z} -\left(a_5+z\frac{\partial}{\partial z}\right)u=0.
\]
Hence
\begin{equation} \label{three}
zu_{zz}+xu_{xz}+yu_{yz} + (c-z)u_z-a_5u=0.
\end{equation}

Substituting equations (303), (304), and (305) into (90), we obtain a system of partial differential equations satisfied by the confluent hypergeometric function ${\rm{E}}_{1}\left(a_1,a_2,a_3,a_4,a_5;c;x,y,z\right)$ defined in
(306):

$
\left\{
\begin{aligned}
&x(1-x)u_{xx}+yu_{xy}+zu_{xz}+\left[c-\left(a_1+a_2+1\right)x\right]u_x-a_1a_2u=0,\\&
y(1-y)u_{yy}+xu_{xy}+zu_{yz}+\left[c-\left(a_3+a_4+1\right)y\right]u_y-a_3a_4u=0,\\&
zu_{zz}+xu_{xz}+yu_{yz} + (c-z)u_z-a_5u=0,
\end{aligned}
\right.
$

where $u\equiv \,\,{\rm{E}}_{1}\left(a_1,a_2,a_3,a_4,a_5;c;x,y,z\right)$.

\bigskip

In a similar way, systems of partial differential equations are constructed that correspond to confluent hypergeometric functions of three variables from the list [7].

\bigskip

\begin{equation}
{\rm{E}}_{2}\left(a_1,a_2,a_3,a_4;c;x,y,z\right)=\sum\limits_{m,n,p=0}^\infty\frac{(a_1)_m(a_2)_m(a_3)_n(a_4)_n}{(c)_{m+n+p}}\frac{x^m}{m!}\frac{y^n}{n!}\frac{z^p}{p!},
\end{equation}

first appearance of this function in the literature: [21],\, ${}_3\Phi_B^{(5)}$, \,\,[9], \,\,$F_{B5}$,

region of convergence:
$$ \left\{ \frac{1}{r}+\frac{1}{s}>1, \,\,\,t<\infty
\right\}.
$$

System of partial differential equations:

$
\left\{
\begin{aligned}
&x(1-x)u_{xx}+yu_{xy}+zu_{xz}+\left[c-\left(a_1+a_2+1\right)x\right]u_x-a_1a_2u=0,\\&
y(1-y)u_{yy}+xu_{xy}+zu_{yz}+\left[c-\left(a_3+a_4+1\right)y\right]u_y-a_3a_4u=0,\\&
zu_{zz}+xu_{xz}+yu_{yz} + cu_z-u=0,
\end{aligned}
\right.
$

where $u\equiv \,\,{\rm{E}}_{2}\left(a_1,a_2,a_3,a_4;c;x,y,z\right)$.

\bigskip

\begin{equation}
{\rm{E}}_{3}\left(a_1,a_2,a_3,a_4;c;x,y,z\right)=\sum\limits_{m,n,p=0}^\infty\frac{(a_1)_m(a_2)_m(a_3)_n(a_4)_p}{(c)_{m+n+p}}\frac{x^m}{m!}\frac{y^n}{n!}\frac{z^p}{p!},
\end{equation}

first appearance of this function in the literature: [21],\, ${}_3\Phi_B^{(2)}$, \,\,[9], \,\,$F_{B2}$,

region of convergence:
$$ \left\{ r<1, \,\,\,s<\infty,\,\,\,t<\infty
\right\}.
$$

System of partial differential equations:

$
\left\{
\begin{aligned}
&x(1-x)u_{xx}+yu_{xy}+zu_{xz}+\left[c-\left(a_1+a_2+1\right)x\right]u_x-a_1a_2u=0,\\&
yu_{yy}+xu_{xy}+zu_{yz}+(c-y)u_y-a_3u=0,\\&
zu_{zz}+xu_{xz}+yu_{yz} + (c-z)u_z-a_4u=0,
\end{aligned}
\right.
$

where $u\equiv \,\,{\rm{E}}_{3}\left(a_1,a_2,a_3,a_4;c;x,y,z\right)$.

\bigskip

\begin{equation}
{\rm{E}}_{4}\left(a_1,a_2,a_3;c;x,y,z\right)=\sum\limits_{m,n,p=0}^\infty\frac{(a_1)_m(a_2)_m(a_3)_n}{(c)_{m+n+p}}\frac{x^m}{m!}\frac{y^n}{n!}\frac{z^p}{p!},
\end{equation}

first appearance of this function in the literature:  \,\,[9], \,\,$F_{B3}$,

region of convergence:
$$ \left\{ r<1, \,\,\,s<\infty,\,\,\,t<\infty
\right\}.
$$

System of partial differential equations:

$
\left\{
\begin{aligned}
&x(1-x)u_{xx}+yu_{xy}+zu_{xz}+\left[c-\left(a_1+a_2+1\right)x\right]u_x-a_1a_2u=0,\\&
yu_{yy}+xu_{xy}+zu_{yz}+(c-y)u_y-a_3u=0,\\&
zu_{zz}+xu_{xz}+yu_{yz} + cu_z-u=0,
\end{aligned}
\right.
$

where $u\equiv \,\,{\rm{E}}_{4}\left(a_1,a_2,a_3;c;x,y,z\right)$.

\bigskip

\begin{equation}
{\rm{E}}_{5}\left(a_1,a_2,a_3;c;x,y,z\right)=\sum\limits_{m,n,p=0}^\infty\frac{(a_1)_m(a_2)_n(a_3)_p}{(c)_{m+n+p}}\frac{x^m}{m!}\frac{y^n}{n!}\frac{z^p}{p!}
,
\end{equation}

first appearance of this function in the literature: [21],\, ${}_3\Phi_B^{(3)}$, \,\,[9], \,\,$F_{B6}$,

region of convergence:
$$ \left\{ r<\infty,\, \,\,\,s<\infty,\,\,\,t<\infty
\right\}.
$$

System of partial differential equations:

$
\left\{
\begin{aligned}
&xu_{xx}+yu_{xy}+zu_{xz}+(c-x)u_x-a_1u=0,\\&
yu_{yy}+xu_{xy}+zu_{yz}+(c-y)u_y-a_2u=0,\\&
zu_{zz}+xu_{xz}+yu_{yz} + (c-z)u_z-a_3u=0;
\end{aligned}
\right.
$

where $u\equiv \,\,{\rm{E}}_{5}\left(a_1,a_2,a_3;c;x,y,z\right)$.

\bigskip

\begin{equation}
{\rm{E}}_{6}\left(a, b;c;x,y,z\right)=\sum\limits_{m,n,p=0}^\infty\frac{(a)_m(b)_n}{(c)_{m+n+p}}\frac{x^m}{m!}\frac{y^n}{n!}\frac{z^p}{p!}
,
\end{equation}

first appearance of this function in the literature: [21],\, ${}_3\Phi_B^{(4)}$, \,\,[9], \,\,$F_{B4}$,

region of convergence:
$$ \left\{ r<\infty,\, \,\,\,s<\infty,\,\,\,t<\infty
\right\}.
$$

System of partial differential equations:

$
\left\{
\begin{aligned}
&xu_{xx}+yu_{xy}+zu_{xz}+(c-x)u_x-au=0,\\&
yu_{yy}+xu_{xy}+zu_{yz}+(c-y)u_y-bu=0,\\&
zu_{zz}+xu_{xz}+yu_{yz} + cu_z-u=0,
\end{aligned}
\right.
$

where $u\equiv \,\,{\rm{E}}_{6}\left(a, b;c;x,y,z\right)$.

\bigskip

\begin{equation}
{\rm{E}}_{7}\left(a_1,a_2,a_3,a_4; c; x,y,z\right)=\sum\limits_{m,n,p=0}^\infty\frac{(a_1)_{m+n}(a_2)_m(a_3)_n(a_4)_{p}}{(c)_{m+n+p}}\frac{x^m}{m!}\frac{y^n}{n!}\frac{z^p}{p!}
,
\end{equation}

first appearance of this function in the literature: [21],\, ${}_3\Phi_S^{(2)}$,

region of convergence:
$$ \left\{ r<1, \,\,\,\, s<1,\,\,\,t<\infty
\right\}.
$$

System of partial differential equations:

$
\left\{
\begin{aligned}
&x(1-x)u_{xx}+(1-x)yu_{xy}+zu_{xz}+\left[c-\left(a_1+a_2+1\right)x\right]u_x-a_2yu_y-a_1a_2u=0,\\&
y(1-y)u_{yy}+x(1-y)u_{xy}+zu_{yz}-a_3xu_x+\left[c-\left(a_1+a_3+1\right)y\right]u_y-a_1a_3u=0,\\
&zu_{zz}+xu_{xz}+yu_{yz} + (c-z)u_z-a_4u=0,
\end{aligned}
\right.
$

where $u\equiv \,\,{\rm{E}}_{7}\left(a_1,a_2,a_3,a_4; c; x,y,z\right)$.

\bigskip

\begin{equation}
{\rm{E}}_{8}\left(a_1,a_2,a_3, a_4; c; x,y,z\right)=\sum\limits_{m,n,p=0}^\infty\frac{(a_1)_{m+n}(a_2)_m(a_3)_{p}(a_4)_{p}}{(c)_{m+n+p}}\frac{x^m}{m!}\frac{y^n}{n!}\frac{z^p}{p!}
,
\end{equation}

first appearance of this function in the literature: [21],\, ${}_3\Phi_S^{(1)}$,

region of convergence:
$$ \left\{ \frac{1}{r}+\frac{1}{t}>1, \,\,\,s<\infty
\right\}.
$$

System of partial differential equations:

$
\left\{
\begin{aligned}
&x(1-x)u_{xx}+(1-x)yu_{xy}+zu_{xz}+\left[c-\left(a_1+a_2+1\right)x\right]u_x-a_2yu_y-a_1a_2u=0,\\
&yu_{yy}+xu_{xy}+zu_{yz}-xu_x+(c-y)u_y-a_1u=0,\\
&z(1-z)u_{zz}+xu_{xz}+yu_{yz}+\left[c-\left(a_3+a_4+1\right)z\right]u_z-a_3a_4u=0,
\end{aligned}
\right.
$

where $u\equiv \,\,{\rm{E}}_{8}\left(a_1,a_2,a_3, a_4; c; x,y,z\right)$.

\bigskip

\begin{equation}
{\rm{E}}_{9}\left(a_1,a_2,a_3; c; x,y,z\right)=\sum\limits_{m,n,p=0}^\infty\frac{(a_1)_{m+n}(a_2)_m(a_3)_n}{(c)_{m+n+p}}\frac{x^m}{m!}\frac{y^n}{n!}\frac{z^p}{p!}
,
\end{equation}

first appearance of this function in the literature: [21],\, ${}_3\Phi_S^{(4)}$,

region of convergence:
$$ \left\{ r<1, \,\,\,\,\, s<1,\,\,\,t<\infty
\right\}.
$$

System of partial differential equations:

$
\left\{
\begin{aligned}
&x(1-x)u_{xx}+(1-x)yu_{xy}+zu_{xz}+\left[c-\left(a_1+a_2+1\right)x\right]u_x-a_2yu_y-a_1a_2u=0,\\
&y(1-y)u_{yy}+x(1-y)u_{xy}+zu_{yz}-a_3xu_x+\left[c-\left(a_1+a_3+1\right)y\right]u_y-a_1a_3u=0,\\
&zu_{zz}+xu_{xz}+yu_{yz} + cu_z-u=0,
\end{aligned}
\right.
$

where $u\equiv \,\,{\rm{E}}_{9}\left(a_1,a_2,a_3; c; x,y,z\right)$.

\bigskip

\begin{equation}
{\rm{E}}_{10}\left(a_1,a_2,a_3; c; x,y,z\right)=\sum\limits_{m,n,p=0}^\infty\frac{(a_1)_{m+n}(a_2)_m(a_3)_{p}}{(c)_{m+n+p}}\frac{x^m}{m!}\frac{y^n}{n!}\frac{z^p}{p!}
,
\end{equation}

first appearance of this function in the literature: [21],\, ${}_3\Phi_T^{(2)}$,

region of convergence:
$$ \left\{ r<1, \,\,\,s<\infty\,\,\,t<\infty
\right\}.
$$

System of partial differential equations:

$
\left\{
\begin{aligned}
&x(1-x)u_{xx}+(1-x)yu_{xy}+zu_{xz}+\left[c-\left(a_1+a_2+1\right)x\right]u_x-a_2yu_y-a_1a_2u=0,\\
&yu_{yy}+xu_{xy}+zu_{yz}-xu_x+(c-y)u_y-a_1u=0,\\
&zu_{zz}+xu_{xz}+yu_{yz} + (c-z)u_z-a_3u=0,
\end{aligned}
\right.
$

where $u\equiv \,\,{\rm{E}}_{10}\left(a_1,a_2,a_3; c; x,y,z\right)$.

\bigskip

\begin{equation}
{\rm{E}}_{11}\left(a, b; c; x,y,z\right)=\sum\limits_{m,n,p=0}^\infty\frac{(a)_{m+n}(b)_m}{(c)_{m+n+p}}\frac{x^m}{m!}\frac{y^n}{n!}\frac{z^p}{p!}
,
\end{equation}

first appearance of this function in the literature: [21],\, ${}_3\Phi_S^{(5)}$,

region of convergence:
$$ \left\{ r<1, \,\,\,s<\infty\,\,\,t<\infty
\right\}.
$$

System of partial differential equations:

$
\left\{
\begin{aligned}
&x(1-x)u_{xx}+(1-x)yu_{xy}+zu_{xz}+\left[c-\left(a+b+1\right)x\right]u_x-byu_y-abu=0,\\
&yu_{yy}+xu_{xy}+zu_{yz}-xu_x+(c-y)u_y-au=0,\\
&zu_{zz}+xu_{xz}+yu_{yz} + cu_z-u=0=0,
\end{aligned}
\right.
$

where $u\equiv \,\,{\rm{E}}_{11}\left(a, b; c; x,y,z\right)$.

\bigskip

\begin{equation}
{\rm{E}}_{12}\left(a_1,a_2,a_3; c; x,y,z\right)=\sum\limits_{m,n,p=0}^\infty\frac{(a_1)_{m+n}(a_2)_{n+p}(a_3)_m}{(c)_{m+n+p}}\frac{x^m}{m!}\frac{y^n}{n!}\frac{z^p}{p!}
,
\end{equation}

first appearance of this function in the literature: [21],\, ${}_3\Phi_T^{(1)}$,

region of convergence:
$$ \left\{ r<1, \,\,\,\,\, s<1,\,\,\,t<\infty
\right\}.
$$

System of partial differential equations:

$
\left\{
\begin{aligned}
&x(1-x)u_{xx}+(1-x)yu_{xy}+zu_{xz}+\left[c-\left(a_1+a_3+1\right)x\right]u_x-a_3yu_y-a_1a_3u=0,\\
&y(1-y)u_{yy}
+x(1-y)u_{xy}-xzu_{xz}+(1-y)zu_{yz}-a_2xu_x\hfill \\ &\,\,\,\,\,\,\,\,\,+\left[c-\left(a_1+a_2+1\right)y\right]u_y-a_1zu_z-a_1a_2u=0,\\
&zu_{zz}+xu_{xz}+yu_{yz}-yu_y+(c-z)u_z-a_2u=0,
\end{aligned}
\right.
$

where $u\equiv \,\,{\rm{E}}_{12}\left(a_1,a_2,a_3; c; x,y,z\right)$.

\bigskip

\begin{equation}
{\rm{E}}_{13}\left(a,b; c; x,y,z\right)=\sum\limits_{m,n,p=0}^\infty\frac{(a)_{m+n}(b)_{n+p}}{(c)_{m+n+p}}\frac{x^m}{m!}\frac{y^n}{n!}\frac{z^p}{p!}
,
\end{equation}

first appearance of this function in the literature: [21],\, ${}_3\Phi_T^{(3)}$,

region of convergence:
$$ \left\{ r<\infty, \,\,\, s<1,\,\,\,t<\infty
\right\}.
$$

System of partial differential equations:

$
\left\{
\begin{aligned}
&xu_{xx}+yu_{xy}+zu_{xz}+(c-x)u_x-yu_y-au=0,\\
&y(1-y)u_{yy}+x(1-y)u_{xy}+(1-y)zu_{yz}-xzu_{xz}-bxu_x  \\ &\,\,\,\,\,\,\,\,\, +\left[c-\left(a+b+1\right)y\right]u_y-azu_z-abu=0,\\
&zu_{zz}+xu_{xz}+yu_{yz}-yu_y+(c-z)u_z-bu=0,
\end{aligned}
\right.
$

where $u\equiv \,\,{\rm{E}}_{13}\left(a,b; c; x,y,z\right)$.

\bigskip

\begin{equation}
{\rm{E}}_{14}\left(a_1,a_2,a_3; c; x,y,z\right)=\sum\limits_{m,n,p=0}^\infty\frac{(a_1)_{m+n+p}(a_2)_{m}(a_3)_n}{(c)_{m+n+p}}\frac{x^m}{m!}\frac{y^n}{n!}\frac{z^p}{p!}
,
\end{equation}

first appearance of this function in the literature: [21],\, ${}_3\Phi_D^{(1)}$,

region of convergence:
$$ \left\{ r<1, \,\,\,\, s<1,\,\,\,t<\infty
\right\}.
$$

System of partial differential equations:

$
\left\{
\begin{aligned}
&x(1-x)u_{xx}+(1-x)yu_{xy}+(1-x)zu_{xz}  \\ &\,\,\,\,\,\,\,\,\,+\left[c-\left(a_1+a_2+1\right)x\right]u_x-a_2yu_y-a_2zu_z-a_1a_2u=0,\\
&y(1-y)u_{yy}+x(1-y)u_{xy}+(1-y)zu_{yz}  \\ &\,\,\,\,\,\,\,\,\, -a_3xu_x+\left[c-\left(a_1+a_3+1\right)y\right]u_y-a_3zu_z-a_1a_3u=0,\\
&zu_{zz}+xu_{xz}+yu_{yz}-xu_x-yu_y+(c-z)u_z-a_1u=0,
\end{aligned}
\right.
$

where $u\equiv \,\,{\rm{E}}_{14}\left(a_1,a_2,a_3; c; x,y,z\right)$.

\bigskip

\begin{equation}
{\rm{E}}_{15}\left(a,b_1, b_2; c; x,y,z\right)=\sum\limits_{m,n,p=0}^\infty\frac{(a)_{2m+n}(b_1)_{n}(b_2)_p}{(c)_{m+n+p}}\frac{x^m}{m!}\frac{y^n}{n!}\frac{z^p}{p!},
\end{equation}

region of convergence:
$$ \left\{ \left\{r<\frac{1}{4}, \,\,\,\, s<\frac{1}{2}+\frac{1}{2}\sqrt{1-4r}\right\}\cup\left\{s \leq \frac{1}{2}\right\},\,\,\,t<\infty
\right\}.
$$

System of partial differential equations:

$
\left\{
\begin{aligned}
&x(1-4x)u_{xx}+(1-4x)yu_{xy}+zu_{xz}-y^2u_{yy} \\ &\,\,\,\,\,\,\,\,\,+\left[c-\left(4a+6\right)x\right]u_x-2\left(a+1\right)yu_y-a\left(a+1\right)u=0,\\
&y(1-y)u_{yy}+x(1-2y)u_{xy}+zu_{yz}-2b_1xu_x+\left[c-\left(a+b_1+1\right)y\right]u_y-ab_1u=0,\\
&zu_{zz}+xu_{xz}+yu_{yz} + (c-z)u_z-b_2u=0,
\end{aligned}
\right.
$

where $u\equiv \,\,{\rm{E}}_{15}\left(a,b_1, b_2; c; x,y,z\right)$.

\bigskip

\begin{equation}
{\rm{E}}_{16}\left(a, b_1,b_2; c; x,y,z\right)=\sum\limits_{m,n,p=0}^\infty\frac{(a)_{2m+n}(b_1)_p(b_2)_p}{(c)_{m+n+p}}\frac{x^m}{m!}\frac{y^n}{n!}\frac{z^p}{p!},
\end{equation}

region of convergence:
$$ \left\{ \frac{1}{4r}+\frac{1}{t}>1, \,\,\,s<\infty
\right\}.
$$

System of partial differential equations:

$
\left\{
\begin{aligned}
&x(1-4x)u_{xx}+(1-4x)yu_{xy}+zu_{xz}-y^2u_{yy} \\ &\,\,\,\,\,\,\,\,\,+\left[c-\left(4a+6\right)x\right]u_x-2\left(a+1\right)yu_y-a\left(a+1\right)u=0,\\
&yu_{yy}+xu_{xy}+zu_{yz}-2xu_x+(c-y)u_y-au=0,\\
&z(1-z)u_{zz}+xu_{xz}+yu_{yz} +\left[c-\left(b_1+b_2+1\right)z\right]u_z-b_1b_2u=0,
\end{aligned}
\right.
$

where $u\equiv \,\,{\rm{E}}_{16}\left(a, b_1,b_2; c; x,y,z\right)$.

\bigskip

\begin{equation}
{\rm{E}}_{17}\left(a,b; c; x,y,z\right)=\sum\limits_{m,n,p=0}^\infty\frac{(a)_{2m+n}(b)_{n}}{(c)_{m+n+p}}\frac{x^m}{m!}\frac{y^n}{n!}\frac{z^p}{p!},
\end{equation}

region of convergence:
$$ \left\{\left\{ r<\frac{1}{4}, \,\,\,\, s<\frac{1}{2}+\frac{1}{2}\sqrt{1-4r}\right\}\cup\left\{s \leq \frac{1}{2}\right\},\,\,\,t<\infty
\right\}.
$$

System of partial differential equations:

$
\left\{
\begin{aligned}
&x(1-4x)u_{xx}+(1-4x)yu_{xy}+zu_{xz}-y^2u_{yy} \\& \,\,\,\,\,\,\,\,\,+\left[c-\left(4a+6\right)x\right]u_x-2\left(a+1\right)yu_y-a\left(a+1\right)u=0,\\
&y(1-y)u_{yy}+x(1-2y)u_{xy}+zu_{yz}-2bxu_x+\left[c-\left(a+b+1\right)y\right]u_y-abu=0,\\
&zu_{zz}+xu_{xz}+yu_{yz} + cu_z-u=0,
\end{aligned}
\right.
$

where $u\equiv \,\,{\rm{E}}_{17}\left(a,b; c; x,y,z\right)$.

\bigskip

\begin{equation}
{\rm{E}}_{18}\left(a,b; c; x,y,z\right)=\sum\limits_{m,n,p=0}^\infty\frac{(a)_{2m+n}(b)_p}{(c)_{m+n+p}}\frac{x^m}{m!}\frac{y^n}{n!}\frac{z^p}{p!},
\end{equation}

region of convergence:
$$ \left\{ r<\frac{1}{4}, \,\,\,s<\infty,\,\,\,t<\infty
\right\}.
$$

System of partial differential equations:

$
\left\{
\begin{aligned}
&x(1-4x)u_{xx}+(1-4x)yu_{xy}+zu_{xz}-y^2u_{yy} \\ &\,\,\,\,\,\,\,\,\,+\left[c-\left(4a+6\right)x\right]u_x-2\left(a+1\right)yu_y-a\left(a+1\right)u=0,\\
&yu_{yy}+xu_{xy}+zu_{yz}-2xu_x+(c-y)u_y-au=0,\\
&zu_{zz}+xu_{xz}+yu_{yz} + (c-z)u_z-bu=0,
\end{aligned}
\right.
$

where $u\equiv \,\,{\rm{E}}_{18}\left(a,b; c; x,y,z\right)$.

\bigskip

\begin{equation}
{\rm{E}}_{19}\left(a; c; x,y,z\right)=\sum\limits_{m,n,p=0}^\infty\frac{(a)_{2m+n}}{(c)_{m+n+p}}\frac{x^m}{m!}\frac{y^n}{n!}\frac{z^p}{p!},
\end{equation}

region of convergence:
$$ \left\{ r<\frac{1}{4}, \,\,\,s<\infty\,\,\,t<\infty
\right\}.
$$

System of partial differential equations:

$
\left\{
\begin{aligned}
&x(1-4x)u_{xx}+(1-4x)yu_{xy}+zu_{xz}-y^2u_{yy} \\& \,\,\,\,\,\,\,\,\,+\left[c-\left(4a+6\right)x\right]u_x-2\left(a+1\right)yu_y-a\left(a+1\right)u=0,\\
&yu_{yy}+xu_{xy}+zu_{yz}-2xu_x+(c-y)u_y-au=0,\\
&zu_{zz}+xu_{xz}+yu_{yz} + cu_z-u=0,
\end{aligned}
\right.
$

where $u\equiv \,\,{\rm{E}}_{19}\left(a; c; x,y,z\right)$.

\bigskip

\begin{equation}
{\rm{E}}_{20}\left(a,b; c; x,y,z\right)=\sum\limits_{m,n,p=0}^\infty\frac{(a)_{2m+n}(b)_{n+p}}{(c)_{m+n+p}}\frac{x^m}{m!}\frac{y^n}{n!}\frac{z^p}{p!},
\end{equation}

region of convergence:
$$ \left\{\left\{ r<\frac{1}{4},\,\,\,\, s<\frac{1}{2}+\frac{1}{2}\sqrt{1-4r}\right\}\cup\left\{s \leq \frac{1}{2}\right\},\,\,\,t<\infty
\right\}.
$$

System of partial differential equations:

$
\left\{
\begin{aligned}
&x(1-4x)u_{xx}+(1-4x)yu_{xy}+zu_{xz}-y^2u_{yy} \\& \,\,\,\,\,\,\,\,\,+\left[c-\left(4a+6\right)x\right]u_x-2\left(a+1\right)yu_y-a\left(a+1\right)u=0,\\&
y(1-y)u_{yy}+x(1-2y)u_{xy}+(1-y)zu_{yz}-2xzu_{xz}-2bxu_x\\& \,\,\,\,\,\,\,\,\,+\left[c-\left(a+b+1\right)y\right]u_y-azu_z-abu=0,\\&
zu_{zz}+xu_{xz}+yu_{yz}-yu_y+(c-z)u_z-bu=0,
\end{aligned}
\right.
$

where $u\equiv \,\,{\rm{E}}_{20}\left(a,b; c; x,y,z\right)$.

\bigskip

\begin{equation}
{\rm{E}}_{21}\left(a,b; c; x,y,z\right)=\sum\limits_{m,n,p=0}^\infty\frac{(a)_{2m+n+p}(b)_{n}}{(c)_{m+n+p}}\frac{x^m}{m!}\frac{y^n}{n!}\frac{z^p}{p!},
\end{equation}

region of convergence:
$$ \left\{\left\{r<\frac{1}{4}, \,\,\,\, s<\frac{1}{2}+\frac{1}{2}\sqrt{1-4r}\right\}\cup\left\{s \leq \frac{1}{2}\right\},\,\,\,t<\infty
\right\}.
$$

System of partial differential equations:

$
\left\{
\begin{aligned}
&x(1-4x)u_{xx}-y^2u_{yy}-z^2u_{zz}+(1-4x)yu_{xy}+(1-4x)zu_{xz}-2yzu_{yz}\\& \,\,\,\,\,\,\,\,\,+\left[c-\left(4a+6\right)x\right]u_x-2\left(a+1\right)yu_y
-2\left(a+1\right)zu_z-a\left(a+1\right)u=0,\\&
y(1-y)u_{yy}+x(1-2y)u_{xy}+(1-y)zu_{yz}-2bxu_x+\left[c-\left(a+b+1\right)y\right]u_y-bzu_z-abu=0,\\&
zu_{zz}+xu_{xz}+yu_{yz}-2xu_x-yu_y + (c-z)u_z-au=0,
\end{aligned}
\right.\,\,\,\,\,\,\,\,\,\,\,\,\,\,\,\,\,\,\,\,
$

where $u\equiv \,\,{\rm{E}}_{21}\left(a,b; c; x,y,z\right)$.

\bigskip

\begin{equation}
{\rm{E}}_{22}\left(a_1,a_2,a_3,a_4;c_1, c_2; x,y,z\right)=\sum\limits_{m,n,p=0}^\infty\frac{(a_1)_{n+p}(a_2)_m(a_3)_n(a_4)_p}{(c_1)_{m+n}(c_2)_{p}}\frac{x^m}{m!}\frac{y^n}{n!}\frac{z^p}{p!},
\end{equation}

first appearance of this function in the literature: [21],\, ${}_3\Phi_N^{(1)}$,

region of convergence:
$$ \left\{ s+t<1,\,\,\,\,r<\infty
\right\}.
$$

System of partial differential equations:

$
\left\{
\begin{aligned}
&xu_{xx}+yu_{xy}+\left(c_1-x\right)u_x-a_2u=0,\\
&y(1-y)u_{yy}+xu_{xy}-yzu_{yz}+\left[c_1-\left(a_1+a_3+1\right)y\right]u_y-a_3zu_z-a_1a_3u=0,\\
&z(1-z)u_{zz}-yzu_{yz}-a_4yu_{y}+\left[c_2-\left(a_1+a_4+1\right)z\right]u_z-a_1a_4u=0,
\end{aligned}
\right.
$

where $u\equiv \,\,{\rm{E}}_{22}\left(a_1,a_2,a_3,a_4;c_1, c_2; x,y,z\right)$.

Particular solutions:

$
{u_1} ={\rm{E}}_{22}\left(a_1,a_2,a_3,a_4;c_1, c_2; x,y,z\right) ,
$

$
{u_2} = {z^{1 - c_2}}{\rm{E}}_{22}\left(1-c_2+a_1,a_2,a_3,1-c_2+a_4;c_1, 2-c_2; x,y,z\right).$

\bigskip

\begin{equation}
{\rm{E}}_{23}\left(a_1,a_2,a_3,a_4;c_1, c_2; x,y,z\right)=\sum\limits_{m,n,p=0}^\infty\frac{(a_1)_{n+p}(a_2)_m(a_3)_m(a_4)_p}{(c_1)_{m+n}(c_2)_{p}}\frac{x^m}{m!}\frac{y^n}{n!}\frac{z^p}{p!},
\end{equation}

first appearance of this function in the literature: [21],\, ${}_3\Phi_N^{(2)}$,

region of convergence:
$$ \left\{ r<1, \,\,\,\, t<1,\,\,\,\,s<\infty
\right\}.
$$

System of partial differential equations:

$
\left\{
\begin{aligned}
&x(1-x)u_{xx}+yu_{xy}+\left[c_1-\left(a_2+a_3+1\right)x\right]u_x-a_2a_3u=0,\\&
yu_{yy}+xu_{xy}+\left(c_1-y\right)u_y-zu_z-a_1u=0,\\&
z(1-z)u_{zz}-yzu_{yz}-a_4yu_{y}+\left[c_2-\left(a_1+a_4+1\right)z\right]u_z-a_1a_4u=0,
\end{aligned}
\right.
$

where $u\equiv \,\,{\rm{E}}_{23}\left(a_1,a_2,a_3,a_4;c_1, c_2; x,y,z\right)$.

Particular solutions:

$
{u_1} ={\rm{E}}_{23}\left(a_1,a_2,a_3,a_4;c_1, c_2; x,y,z\right) ,
$

$
{u_2} = {z^{1 - c_2}}{\rm{E}}_{23}\left(1-c_2+a_1,a_2,a_3,1-c_2+a_4;c_1, 2-c_2; x,y,z\right).
$

\bigskip

\begin{equation}
{\rm{E}}_{24}\left(a_1,a_2,a_3,a_4;c_1, c_2; x,y,z\right)=\sum\limits_{m,n,p=0}^\infty\frac{(a_1)_{n+p}(a_2)_m(a_3)_m(a_4)_n}{(c_1)_{m+n}(c_2)_{p}}\frac{x^m}{m!}\frac{y^n}{n!}\frac{z^p}{p!},
\end{equation}

region of convergence:
$$ \left\{ \frac{1}{r}+\frac{1}{s}>1,\,\,\,\,t<\infty
\right\}.
$$

System of partial differential equations:

$
\left\{
\begin{aligned}
&x(1-x)u_{xx}+yu_{xy}+\left[c_1-\left(a_2+a_3+1\right)x\right]u_x-a_2a_3u=0,\\&
y(1-y)u_{yy}+xu_{xy}-yzu_{yz}+\left[c_1-\left(a_1+a_4+1\right)y\right]u_y-a_4zu_z-a_1a_4u=0,\\&
zu_{zz}-yu_y+\left(c_2-z\right)u_z-a_1u=0,
\end{aligned}
\right.
$

where $u\equiv \,\,{\rm{E}}_{24}\left(a_1,a_2,a_3,a_4;c_1, c_2; x,y,z\right)$.

Particular solutions:

$
{u_1} ={\rm{E}}_{24}\left(a_1,a_2,a_3,a_4;c_1, c_2; x,y,z\right) ,
$

$
{u_2} = {z^{1 - c_2}}{\rm{E}}_{24}\left(1-c_2+a_1,a_2,a_3,a_4;c_1, 2-c_2; x,y,z\right).
$

\bigskip

\begin{equation}
{\rm{E}}_{25}\left(a_1,a_2,a_3;c_1, c_2; x,y,z\right)=\sum\limits_{m,n,p=0}^\infty\frac{(a_1)_{n+p}(a_2)_m(a_3)_m}{(c_1)_{m+n}(c_2)_{p}}\frac{x^m}{m!}\frac{y^n}{n!}\frac{z^p}{p!},
\end{equation}

first appearance of this function in the literature: [21],\, ${}_3\Phi_N^{(3)}$,

region of convergence:
$$ \left\{ r<1 ,\,\,\,\,s<\infty,\,\,\,\,t<\infty
\right\}.
$$

System of partial differential equations:

$
\left\{
\begin{aligned}
&x(1-x)u_{xx}+yu_{xy}+\left[c_1-\left(a_2+a_3+1\right)x\right]u_x-a_2a_3u=0,\\&
yu_{yy}+xu_{xy}+\left(c_1-y\right)u_y-zu_z-a_1u=0,\\&
zu_{zz}-yu_y+\left(c_2-z\right)u_z-a_1u=0,
\end{aligned}
\right.
$

where $u\equiv \,\,{\rm{E}}_{25}\left(a_1,a_2,a_3;c_1, c_2; x,y,z\right)$.

Particular solutions:

$
{u_1} ={\rm{E}}_{25}\left(a_1,a_2,a_3;c_1, c_2; x,y,z\right) ,
$

$
{u_2} = {z^{1 - c_2}}{\rm{E}}_{25}\left(1-c_2+a_1,a_2,a_3;c_1, 2-c_2; x,y,z\right).
$

\bigskip

\begin{equation}
{\rm{E}}_{26}\left(a_1,a_2,a_3;c_1, c_2; x,y,z\right)=\sum\limits_{m,n,p=0}^\infty\frac{(a_1)_{n+p}(a_2)_m(a_3)_n}{(c_1)_{m+n}(c_2)_{p}}\frac{x^m}{m!}\frac{y^n}{n!}\frac{z^p}{p!},
\end{equation}

first appearance of this function in the literature: [21],\, ${}_3\Phi_P^{(2)}$,

region of convergence:
$$ \left\{ r<\infty ,\,\,\,\,s<1,\,\,\,\,t<\infty
\right\}.
$$

System of partial differential equations:

$
\left\{
\begin{aligned}
&xu_{xx}+yu_{xy}+\left(c_1-x\right)u_x-a_2u=0,\\
&y(1-y)u_{yy}+xu_{xy}-yzu_{yz}+\left[c_1-\left(a_1+a_3+1\right)y\right]u_y-a_3zu_z-a_1a_3u=0,\\
&zu_{zz}-yu_y+\left(c_2-z\right)u_z-a_1u=0,
\end{aligned}
\right.
$

where $u\equiv \,\,{\rm{E}}_{26}\left(a_1,a_2,a_3;c_1, c_2; x,y,z\right)$.

Particular solutions:

$
{u_1} ={\rm{E}}_{26}\left(a_1,a_2,a_3;c_1, c_2; x,y,z\right) ,
$

$
{u_2} = {z^{1 - c_2}}{\rm{E}}_{26}\left(1-c_2+a_1,a_2,a_3;c_1, 2-c_2; x,y,z\right).
$

\bigskip

\begin{equation}
{\rm{E}}_{27}\left(a_1,a_2,a_3;c_1, c_2; x,y,z\right)=\sum\limits_{m,n,p=0}^\infty\frac{(a_1)_{n+p}(a_2)_m(a_3)_p}{(c_1)_{m+n}(c_2)_{p}}\frac{x^m}{m!}\frac{y^n}{n!}\frac{z^p}{p!},
\end{equation}

first appearance of this function in the literature: [21],\, ${}_3\Phi_M^{(2)}$,

region of convergence:
$$ \left\{ r<\infty,\,\,\,\,s<\infty,\,\,\,\,t<1
\right\}.
$$

System of partial differential equations:

$
\left\{
\begin{aligned}
&xu_{xx}+yu_{xy}+\left(c_1-x\right)u_x-a_2u=0,\\&
yu_{yy}+xu_{xy}+\left(c_1-y\right)u_y-zu_z-a_1u=0,\\&
z(1-z)u_{zz}-yzu_{yz}-a_3yu_{y}+\left[c_2-\left(a_1+a_3+1\right)z\right]u_z-a_1a_3u=0,
\end{aligned}
\right.
$

where $u\equiv \,\,{\rm{E}}_{27}\left(a_1,a_2,a_3;c_1, c_2; x,y,z\right)$.

Particular solutions:

$
{u_1} ={\rm{E}}_{27}\left(a_1,a_2,a_3;c_1, c_2; x,y,z\right) ,
$

$
{u_2} = {z^{1 - c_2}}{\rm{E}}_{27}\left(1-c_2+a_1,a_2,1-c_1+a_3;c_1, 2-c_2; x,y,z\right).
$

\bigskip

\begin{equation}
{\rm{E}}_{28}\left(a_1,a_2,a_3;c_1, c_2; x,y,z\right)=\sum\limits_{m,n,p=0}^\infty\frac{(a_1)_{n+p}(a_2)_n(a_3)_p}{(c_1)_{m+n}(c_2)_{p}}\frac{x^m}{m!}\frac{y^n}{n!}\frac{z^p}{p!},
\end{equation}

region of convergence:
$$ \left\{ r<\infty,\,\,\,\,s+t<1
\right\}.
$$

System of partial differential equations:

$
\left\{
\begin{aligned}
&xu_{xx}+yu_{xy}+c_1u_x-u=0,\\&
y(1-y)u_{yy}+xu_{xy}-yzu_{yz}+\left[c_1-\left(a_1+a_2+1\right)y\right]u_y-a_2zu_z-a_1a_2u=0,\\&
z(1-z)u_{zz}-yzu_{yz}-a_3yu_{y}+\left[c_2-\left(a_1+a_3+1\right)z\right]u_z-a_1a_3u=0,
\end{aligned}
\right.
$

where $u\equiv \,\,{\rm{E}}_{28}\left(a_1,a_2,a_3;c_1, c_2; x,y,z\right)$.

Particular solutions:

$
{u_1} ={\rm{E}}_{28}\left(a_1,a_2,a_3;c_1, c_2; x,y,z\right) ,
$

$
{u_2} = {z^{1 - c_2}}{\rm{E}}_{28}\left(1-c_2+a_1,a_2,1-c_1+a_3;c_1, 2-c_2; x,y,z\right).
$

\bigskip

\begin{equation}
{\rm{E}}_{29}\left(a,b;c_1, c_2; x,y,z\right)=\sum\limits_{m,n,p=0}^\infty\frac{(a)_{n+p}(b)_m}{(c_1)_{m+n}(c_2)_{p}}\frac{x^m}{m!}\frac{y^n}{n!}\frac{z^p}{p!}
,
\end{equation}

first appearance of this function in the literature: [21],\, ${}_3\Phi_F^{(3)}$,

region of convergence:
$$ \left\{ r<\infty,\,\,\,\,s<1,\,\,\,\,\,t<\infty
\right\}.
$$

System of partial differential equations:

$
\left\{
\begin{aligned}
&xu_{xx}+yu_{xy}+\left(c_1-x\right)u_x-bu=0,\\&
yu_{yy}+xu_{xy}+\left(c_1-y\right)u_y-zu_z-au=0,\\&
zu_{zz}-yu_y+\left(c_2-z\right)u_z-au=0,
\end{aligned}
\right.
$

where $u\equiv \,\,{\rm{E}}_{29}\left(a,b;c_1, c_2; x,y,z\right)$.

Particular solutions:

$
{u_1} ={\rm{E}}_{29}\left(a,b;c_1, c_2; x,y,z\right) ,
$

$
{u_2} = {z^{1 - c_2}}{\rm{E}}_{29}\left(1-c_2+a,b;c_1, 2-c_2; x,y,z\right).
$

\bigskip

\begin{equation}
{\rm{E}}_{30}\left(a,b;c_1, c_2; x,y,z\right)=\sum\limits_{m,n,p=0}^\infty\frac{(a)_{n+p}(b)_n}{(c_1)_{m+n}(c_2)_{p}}\frac{x^m}{m!}\frac{y^n}{n!}\frac{z^p}{p!},
\end{equation}

region of convergence:
$$ \left\{ r<\infty,\,\,\,\,s<1,\,\,\,\,\,t<\infty
\right\}.
$$

System of partial differential equations:

$
\left\{
\begin{aligned}
&xu_{xx}+yu_{xy}+c_1u_x-u=0,\\&
y(1-y)u_{yy}+xu_{xy}-yzu_{yz}+\left[c_1-\left(a+b+1\right)y\right]u_y-bzu_z-abu=0,\\&
zu_{zz}-yu_y+\left(c_2-z\right)u_z-au=0,
\end{aligned}
\right.
$

where $u\equiv \,\,{\rm{E}}_{30}\left(a,b;c_1, c_2; x,y,z\right)$.

Particular solutions:

$
{u_1} ={\rm{E}}_{30}\left(a,b;c_1, c_2; x,y,z\right) ,
$

$
{u_2} = {z^{1 - c_2}}{\rm{E}}_{30}\left(1-c_2+a,b;c_1, 2-c_2; x,y,z\right).
$

\bigskip

\begin{equation}
{\rm{E}}_{31}\left(a,b;c_1, c_2; x,y,z\right)=\sum\limits_{m,n,p=0}^\infty\frac{(a)_{n+p}(b)_p}{(c_1)_{m+n}(c_2)_{p}}\frac{x^m}{m!}\frac{y^n}{n!}\frac{z^p}{p!}
,
\end{equation}

first appearance of this function in the literature: [21],\, ${}_3\Phi_M^{(4)}$,

region of convergence:
$$ \left\{ r<\infty,\,\,\,\,s<\infty,\,\,\,\,\,t<1
\right\}.
$$

System of partial differential equations:

$
\left\{
\begin{aligned}
&xu_{xx}+yu_{xy}+c_1u_x-u=0,\\
&yu_{yy}+xu_{xy}+\left(c_1-y\right)u_y-zu_z-au=0,\\
&z(1-z)u_{zz}-yzu_{yz}-byu_y+\left[c_2-\left(a+b+1\right)z\right]u_z-abu=0,
\end{aligned}
\right.
$

where $u\equiv \,\,{\rm{E}}_{31}\left(a,b;c_1, c_2; x,y,z\right)$.

Particular solutions:

$
{u_1} ={\rm{E}}_{31}\left(a,b;c_1, c_2; x,y,z\right) ,
$

$
{u_2} = {z^{1 - c_2}}{\rm{E}}_{31}\left(1-c_2+a,1-c_2+b;c_1, 2-c_2; x,y,z\right).
$

\bigskip

\begin{equation}
{\rm{E}}_{32}\left(a;c_1, c_2; x,y,z\right)=\sum\limits_{m,n,p=0}^\infty\frac{(a)_{n+p}}{(c_1)_{m+n}(c_2)_{p}}\frac{x^m}{m!}\frac{y^n}{n!}\frac{z^p}{p!}
,
\end{equation}

first appearance of this function in the literature: [21],\, ${}_3\Phi_F^{(4)}$,

region of convergence:
$$ \left\{ r<\infty,\,\,\,\,s<\infty,\,\,\,\,\,t<\infty
\right\}.
$$

System of partial differential equations:

$
\left\{
\begin{aligned}
&xu_{xx}+yu_{xy}+c_1u_x-u=0,\\&
yu_{yy}+xu_{xy}+\left(c_1-y\right)u_y-zu_z-au=0,\\&
zu_{zz}-yu_y+\left(c_2-z\right)u_z-au=0,
\end{aligned}
\right.
$

where $u\equiv \,\,{\rm{E}}_{32}\left(a;c_1, c_2; x,y,z\right)$.

Particular solutions:

$
{u_1} ={\rm{E}}_{32}\left(a;c_1, c_2; x,y,z\right) ,
$

$
{u_2} = {z^{1 - c_2}}{\rm{E}}_{32}\left(1-c_2+a;c_1, 2-c_2; x,y,z\right).
$

\bigskip

\begin{equation}
{\rm{E}}_{33}\left(a_1, a_2, a_3; c_1, c_2; x,y,z\right)=\sum\limits_{m,n,p=0}^\infty\frac{(a_1)_{n+p}(a_2)_{n+p}(a_3)_{m}}{(c_1)_{m+n}(c_2)_{p}}\frac{x^m}{m!}\frac{y^n}{n!}\frac{z^p}{p!},
\end{equation}

first appearance of this function in the literature: [21],\, ${}_3\Phi_R^{(1)}$,

region of convergence:
$$ \left\{ r<\infty,\,\,\,\,\sqrt{s}+\sqrt{t}<1
\right\}.
$$

System of partial differential equations:

$
\left\{
\begin{aligned}
&xu_{xx}+yu_{xy}+\left(c_1-x\right)u_x-a_3u=0,\\
&y(1-y)u_{yy}-z^2u_{zz}+xu_{xy}-2yzu_{yz}\\& \,\,\,\,\,\,\,\,\,+\left[c_1-\left(a_1+a_2+1\right)y\right]u_y-\left(a_1+a_2+1\right)zu_z-a_1a_2u=0,\\
&z(1-z)u_{zz}-y^2u_{yy}-2yzu_{yz}-\left(a_1+a_2+1\right)yu_y+\left[c_2-\left(a_1+a_2+1\right)z\right]u_z-a_1a_2u=0,
\end{aligned}
\right.
$

where $u\equiv \,\,{\rm{E}}_{33}\left(a_1, a_2, a_3; c_1, c_2; x,y,z\right)$.

Particular solutions:

$
{u_1} ={\rm{E}}_{33}\left(a_1, a_2, a_3; c_1, c_2; x,y,z\right) ,
$

$
{u_2} = {z^{1 - c_2}}{\rm{E}}_{33}\left(1-c_2+a_1, 1-c_2+a_2, a_3; c_1, 2-c_2; x,y,z\right).
$

\bigskip

\begin{equation}
{\rm{E}}_{34}\left(a, b; c_1, c_2; x,y,z\right)=\sum\limits_{m,n,p=0}^\infty\frac{(a)_{n+p}(b)_{n+p}}{(c_1)_{m+n}(c_2)_{p}}\frac{x^m}{m!}\frac{y^n}{n!}\frac{z^p}{p!}
,
\end{equation}

first appearance of this function in the literature: [21],\, ${}_3\Phi_R^{(2)}$,

region of convergence:
$$ \left\{ r<\infty,\,\,\,\,\sqrt{s}+\sqrt{t}<1
\right\}.
$$

System of partial differential equations:

$
\left\{
\begin{aligned}
&xu_{xx}+yu_{xy}+c_1u_x-u=0,\\
&y(1-y)u_{yy}-z^2u_{zz}+xu_{xy}-2yzu_{yz}+\left[c_1-\left(a+b+1\right)y\right]u_y-\left(a+b+1\right)zu_z-abu=0,\\
&z(1-z)u_{zz}-y^2u_{yy}-2yzu_{yz}-\left(a+b+1\right)yu_y+\left[c_2-\left(a+b+1\right)z\right]u_z-abu=0,
\end{aligned}
\right.
$

where $u\equiv \,\,{\rm{E}}_{34}\left(a, b; c_1, c_2; x,y,z\right)$.

Particular solutions:

$
{u_1} ={\rm{E}}_{34}\left(a, b; c_1, c_2; x,y,z\right) ,
$

$
{u_2} = {z^{1 - c_2}}{\rm{E}}_{34}\left(1-c_2+a, 1-c_2+b; c_1, 2-c_2; x,y,z\right).
$

\bigskip

\begin{equation}
{\rm{E}}_{35}\left(a_1,a_2,a_3; c_1, c_2; x,y,z\right)=\sum\limits_{m,n,p=0}^\infty\frac{(a_1)_{m+n}(a_2)_{n+p}(a_3)_m}{(c_1)_{m+n}(c_2)_{p}}\frac{x^m}{m!}\frac{y^n}{n!}\frac{z^p}{p!},
\end{equation}

region of convergence:
$$ \left\{ r<1,\,\,\,\,\, s<1,\,\,\,\,t<\infty
\right\}.
$$

System of partial differential equations:

$
\left\{
\begin{aligned}
&x(1-x)u_{xx}+(1-x)yu_{xy}+\left[c_1-\left(a_1+a_3+1\right)x\right]u_x-a_3yu_y-a_1a_3u=0,\\&
y(1-y)u_{yy}+x(1-y)u_{xy}-xzu_{xz}-yzu_{yz}-a_2xu_x\\& \,\,\,\,\,\,\,\,\,+\left[c_1-\left(a_1+a_2+1\right)y\right]u_y-a_1zu_z-a_1a_2u=0,\\&
zu_{zz}-yu_y+\left(c_2-z\right)u_z-a_2u=0,
\end{aligned}
\right.
$

where $u\equiv \,\,{\rm{E}}_{35}\left(a_1,a_2,a_3; c_1, c_2; x,y,z\right)$.

Particular solutions:

$
{u_1} ={\rm{E}}_{35}\left(a_1,a_2,a_3; c_1, c_2; x,y,z\right) ,
$

$
{u_2} = {z^{1 - c_2}}{\rm{E}}_{35}\left(a_1,1-c_2+a_2,a_3; c_1, 2-c_2; x,y,z\right).
$

\bigskip

\begin{equation}
{\rm{E}}_{36}\left(a_1,a_2, a_3; c_1, c_2; x,y,z\right)=\sum\limits_{m,n,p=0}^\infty\frac{(a_1)_{m+n}(a_2)_{n+p}(a_3)_{p}}{(c_1)_{m+n}(c_2)_{p}}\frac{x^m}{m!}\frac{y^n}{n!}\frac{z^p}{p!},
\end{equation}

first appearance of this function in the literature: [21],\, ${}_3\Phi_M^{(1)}$,

region of convergence:
$$ \left\{ r<\infty,\,\,\,\,s+t<1
\right\}.
$$

System of partial differential equations:

$
\left\{
\begin{aligned}
&xu_{xx}+yu_{xy}+\left(c_1-x\right)u_x-yu_y-a_1u=0,\\&
y(1-y)u_{yy}+x(1-y)u_{xy}-xzu_{xz}-yzu_{yz}-a_2xu_x\\& \,\,\,\,\,\,\,\,\,+\left[c_1-\left(a_1+a_2+1\right)y\right]u_y-a_1zu_z-a_1a_2u=0,\\&
z(1-z)u_{zz}-yzu_{yz}-a_3yu_y+\left[c_2-\left(a_2+a_3+1\right)z\right]u_z-a_2a_3u=0,
\end{aligned}
\right.
$

where $u\equiv \,\,{\rm{E}}_{36}\left(a_1,a_2, a_3; c_1, c_2; x,y,z\right)$.

Particular solutions:

$
{u_1} ={\rm{E}}_{36}\left(a_1,a_2, a_3; c_1, c_2; x,y,z\right) ,
$

$
{u_2} = {z^{1 - c_2}}{\rm{E}}_{36}\left(a_1,1-c_2+a_2, 1-c_2+a_3; c_1, 2-c_2; x,y,z\right).
$

\bigskip

\begin{equation}
{\rm{E}}_{37}\left(a,b; c_1, c_2; x,y,z\right)=\sum\limits_{m,n,p=0}^\infty\frac{(a)_{m+n}(b)_{n+p}}{(c_1)_{m+n}(c_2)_{p}}\frac{x^m}{m!}\frac{y^n}{n!}\frac{z^p}{p!}
 ,
\end{equation}

first appearance of this function in the literature: [37],\, ${}_3\Phi_M^{(3)}$,

region of convergence:
$$ \left\{ r<\infty,\,\,\,\,s<\infty,\,\,\,\,t<\infty
\right\}.
$$

System of partial differential equations:

$
\left\{
\begin{aligned}
&xu_{xx}+yu_{xy}+\left(c_1-x\right)u_x-yu_y-au=0,\\&
y(1-y)u_{yy}+x(1-y)u_{xy}-xzu_{xz}-yzu_{yz}-bxu_x+\left[c_1-\left(a+b+1\right)y\right]u_y-azu_z-abu=0,\\&
zu_{zz}-yu_y+\left(c_2-z\right)u_z-bu=0,
\end{aligned}
\right.
$

where $u\equiv \,\,{\rm{E}}_{37}\left(a,b; c_1, c_2; x,y,z\right)$.

Particular solutions:

$
{u_1} ={\rm{E}}_{37}\left(a,b; c_1, c_2; x,y,z\right) ,
$

$
{u_2} = {z^{1 - c_2}}{\rm{E}}_{37}\left(1-c_2+a,1-c_2+b; c_1, 2-c_2; x,y,z\right).
$

\bigskip

\begin{equation}
{\rm{E}}_{38}\left(a_1,a_2,a_3; c_1, c_2; x,y,z\right)=\sum\limits_{m,n,p=0}^\infty\frac{(a_1)_{m+p}(a_2)_{n+p}(a_3)_m}{(c_1)_{m+n}(c_2)_{p}}\frac{x^m}{m!}\frac{y^n}{n!}\frac{z^p}{p!},
\end{equation}

first appearance of this function in the literature: [21],\, ${}_3\Phi_P^{(1)}$,

region of convergence:
$$ \left\{ r+t<1,\,\,\,\,s<\infty
\right\}.
$$

System of partial differential equations:

$
\left\{
\begin{aligned}
&x(1-x)u_{xx}+yu_{xy}-xzu_{xz}+\left[c_1-\left(a_1+a_3+1\right)x\right]u_x-a_3zu_z-a_1a_3u=0,\\&
yu_{yy}+xu_{xy}+\left(c_1-y\right)u_y-zu_z-a_2u=0,\\&
z(1-z)u_{zz}-xyu_{xy}-xzu_{xz}-yzu_{yz}-a_2xu_x-a_1yu_y+\left[c_2-\left(a_1+a_2+1\right)z\right]u_z-a_1a_2u=0,
\end{aligned}
\right.
$

where $u\equiv \,\,{\rm{E}}_{38}\left(a_1,a_2,a_3; c_1, c_2; x,y,z\right)$.

Particular solutions:

$
{u_1} ={\rm{E}}_{38}\left(a_1,a_2,a_3; c_1, c_2; x,y,z\right) ,
$

$
{u_2} = {z^{1 - c_2}}{\rm{E}}_{38}\left(1-c_2+a_1,1-c_2+a_2,a_3; c_1, 2-c_2; x,y,z\right).
$

\bigskip

\begin{equation}
{\rm{E}}_{39}\left(a,b; c_1, c_2; x,y,z\right)=\sum\limits_{m,n,p=0}^\infty\frac{(a)_{m+p}(b)_{n+p}}{(c_1)_{m+n}(c_2)_{p}}\frac{x^m}{m!}\frac{y^n}{n!}\frac{z^p}{p!}
,
\end{equation}

first appearance of this function in the literature: [21],\, ${}_3\Phi_P^{(3)}$,

region of convergence:
$$ \left\{ r<\infty,\,\,\,\,s<\infty,\,\,\,\,t<1
\right\}.
$$

System of partial differential equations:

$
\left\{
\begin{aligned}
&xu_{xx}+yu_{xy}+\left(c_1-x\right)u_x-zu_z-au=0,\\&
yu_{yy}+xu_{xy}+\left(c_1-y\right)u_y-zu_z-bu=0,\\&
z(1-z)u_{zz}-xyu_{xy}-xzu_{xz}-yzu_{yz}-bxu_x-ayu_y+\left[c_2-\left(a+b+1\right)z\right]u_z-abu=0,
\end{aligned}
\right.
$

where $u\equiv \,\,{\rm{E}}_{39}\left(a,b; c_1, c_2; x,y,z\right)$.

Particular solutions:

$
{u_1} ={\rm{E}}_{39}\left(a,b; c_1, c_2; x,y,z\right) ,
$

$
{u_2} = {z^{1 - c_2}}{\rm{E}}_{39}\left(1-c_2+a,1-c_2+b; c_1, 2-c_2; x,y,z\right).
$

\bigskip

\begin{equation}
{\rm{E}}_{40}\left(a_1,a_2,a_3; c_1, c_2; x,y,z\right)=\sum\limits_{m,n,p=0}^\infty\frac{(a_1)_{m+n+p}(a_2)_{m}(a_3)_n}{(c_1)_{m+n}(c_2)_{p}}\frac{x^m}{m!}\frac{y^n}{n!}\frac{z^p}{p!},
\end{equation}

region of convergence:
$$ \left\{ r<1, \,\,\,\, s<1,\,\,\,\,t<\infty
\right\}.
$$

System of partial differential equations:

$
\left\{
\begin{aligned}
&x(1-x)u_{xx}+(1-x)yu_{xy}-xzu_{xz}+\left[c_1-\left(a_1+a_2+1\right)x\right]u_x-a_2yu_y-a_2zu_z-a_1a_2u=0,\\&
y(1-y)u_{yy}+x(1-y)u_{xy}-yzu_{yz}-a_3xu_x+\left[c_1-\left(a_1+a_3+1\right)y\right]u_y-a_3zu_z-a_1a_3u=0,\\&
zu_{zz}-xu_x-yu_y+\left(c_2-z\right)u_z-a_1u=0,
\end{aligned}
\right.
$

where $u\equiv \,\,{\rm{E}}_{40}\left(a_1,a_2,a_3; c_1, c_2; x,y,z\right)$,

Particular solutions:

$
{u_1} ={\rm{E}}_{40}\left(a_1,a_2,a_3; c_1, c_2; x,y,z\right) ,
$

$
{u_2} = {z^{1 - c_2}}{\rm{E}}_{40}\left(1-c_2+a_1,a_2,a_3; c_1, 2-c_2; x,y,z\right).
$

\bigskip

\begin{equation}
{\rm{E}}_{41}\left(a_1,a_2, a_3; c_1, c_2; x,y,z\right)=\sum\limits_{m,n,p=0}^\infty\frac{(a_1)_{m+n+p}(a_2)_{m}(a_3)_{p}}{(c_1)_{m+n}(c_2)_{p}}\frac{x^m}{m!}\frac{y^n}{n!}\frac{z^p}{p!},
\end{equation}

first appearance of this function in the literature: [21],\, ${}_3\Phi_G^{(1)}$,

region of convergence:
$$ \left\{ r+t<1,\,\,\,\,s<\infty
\right\}.
$$

System of partial differential equations:

$
\left\{
\begin{aligned}
&x(1-x)u_{xx}+(1-x)yu_{xy}-xzu_{xz}+\left[c_1-\left(a_1+a_2+1\right)x\right]u_x-a_2yu_y-a_2zu_z-a_1a_2u=0,\\&
yu_{yy}+xu_{xy}-xu_x+\left(c_1-y\right)u_y-zu_z-a_1u=0,\\&
z(1-z)u_{zz}-xzu_{xz}-yzu_{yz}-a_3xu_x-a_3yu_y+\left[c_2-\left(a_1+a_3+1\right)z\right]u_z-a_1a_3u=0,
\end{aligned}
\right.
$

where $u\equiv \,\,{\rm{E}}_{41}\left(a_1,a_2, a_3; c_1, c_2; x,y,z\right)$.

Particular solutions:

$
{u_1} ={\rm{E}}_{41}\left(a_1,a_2, a_3; c_1, c_2; x,y,z\right) ,
$

$
{u_2} = {z^{1 - c_2}}{\rm{E}}_{41}\left(1-c_2+a_1,a_2, 1-c_2+a_3; c_1, 2-c_2; x,y,z\right).
$

\bigskip

\begin{equation}
{\rm{E}}_{42}\left(a, b; c_1, c_2; x,y,z\right)=\sum\limits_{m,n,p=0}^\infty\frac{(a)_{m+n+p}(b)_{m}}{(c_1)_{m+n}(c_2)_{p}}\frac{x^m}{m!}\frac{y^n}{n!}\frac{z^p}{p!}
,
\end{equation}

first appearance of this function in the literature: [21],\, ${}_3\Phi_F^{(2)}$,

region of convergence:
$$ \left\{ r<1,\,\,\,\,s<\infty,\,\,\,\,t<\infty
\right\}.
$$

System of partial differential equations:

$
\left\{
\begin{aligned}
&x(1-x)u_{xx}+(1-x)yu_{xy}-xzu_{xz}+\left[c_1-\left(a+b+1\right)x\right]u_x-byu_y-bzu_z-abu=0,\\&
yu_{yy}+xu_{xy}-xu_x+\left(c_1-y\right)u_y-zu_z-au=0,\\&
zu_{zz}-xu_x-yu_y+\left(c_2-z\right)u_z-au=0,
\end{aligned}
\right.
$

where $u\equiv \,\,{\rm{E}}_{42}\left(a, b; c_1, c_2; x,y,z\right)$.

Particular solutions:

$
{u_1} ={\rm{E}}_{42}\left(a, b; c_1, c_2; x,y,z\right) ,
$

$
{u_2} = {z^{1 - c_2}}{\rm{E}}_{42}\left(1-c_2+a, b; c_1, 2-c_2; x,y,z\right).
$

\bigskip

\begin{equation}
{\rm{E}}_{43}\left(a, b; c_1, c_2; x,y,z\right)=\sum\limits_{m,n,p=0}^\infty\frac{(a)_{m+n+p}(b)_{n+p}}{(c_1)_{m+n}(c_2)_{p}}\frac{x^m}{m!}\frac{y^n}{n!}\frac{z^p}{p!}
,
\end{equation}

first appearance of this function in the literature: [21],\, ${}_3\Phi_F^{(1)}$,

region of convergence:
$$ \left\{r<\infty,\,\,\,\,\sqrt{s}+\sqrt{t}<1
\right\}.
$$

System of partial differential equations:

$
\left\{
\begin{aligned}
&xu_{xx}+yu_{xy}+\left(c_1-x\right)u_x-yu_y-zu_z-a_1u=0,\\&
y(1-y)u_{yy}-z^2u_{zz}+x(1-y)u_{xy}\\&\,\,\,\,\,\,\,\,\,-2yzu_{yz}-xzu_{xz}-bxu_x+\left[c_1-\left(a+b+1\right)y\right]u_y-\left(a+b+1\right)zu_z-abu=0,\\&
z(1-z)u_{zz}-y^2u_{yy}-xyu_{xy}\\&\,\,\,\,\,\,\,\,\,-xzu_{xz}-2yzu_{yz}-a_2xu_x-\left(a+b+1\right)yu_y+\left[c_2-\left(a+b+1\right)z\right]u_z-abu=0,
\end{aligned}
\right.
$

where $u\equiv \,\,{\rm{E}}_{43}\left(a, b; c_1, c_2; x,y,z\right)$.

Particular solutions:

$
{u_1} ={\rm{E}}_{43}\left(a, b; c_1, c_2; x,y,z\right) ,
$

$
{u_2} = {z^{1 - c_2}}{\rm{E}}_{43}\left(1-c_2+a, 1-c_2+b; c_1, 2-c_2; x,y,z\right).
$

\bigskip

\begin{equation}
{\rm{E}}_{44}\left(a_1, a_2,a_3; c_1, c_2; x,y,z\right)=\sum\limits_{m,n,p=0}^\infty\frac{(a_1)_{2m+p}(a_2)_n(a_3)_n}{(c_1)_{m+n}(c_2)_{p}}\frac{x^m}{m!}\frac{y^n}{n!}\frac{z^p}{p!},
\end{equation}

region of convergence:
$$ \left\{ \frac{1}{4r}+\frac{1}{s}>1,\,\,\,\,t<\infty
\right\}.
$$

System of partial differential equations:

$
\left\{
\begin{aligned}
&x(1-4x)u_{xx}+yu_{xy}-4xzu_{xz}-z^2u_{zz}\\& \,\,\,\,\,\,\,\,\,+\left[c_1-\left(4a_1+6\right)x\right]u_x-2\left(a_1+1\right)zu_z-a_1\left(a_1+1\right)u=0,\\
&y(1-y)u_{yy}+xu_{xy}+\left[c_1-\left(a_2+a_3+1\right)y\right]u_y-a_2a_3u=0,\\&
zu_{zz}-2xu_x+\left(c_2-z\right)u_z-a_1u=0,
\end{aligned}
\right.
$

where $u\equiv \,\,{\rm{E}}_{44}\left(a_1, a_2,a_3; c_1, c_2; x,y,z\right)$.

Particular solutions:

$
{u_1} ={\rm{E}}_{44}\left(a_1, a_2,a_3; c_1, c_2; x,y,z\right) ,
$

$
{u_2} = {z^{1 - c_2}}{\rm{E}}_{44}\left(1-c_2+a_1, a_2,a_3; c_1, 2-c_2; x,y,z\right).
$

\bigskip

\begin{equation}
{\rm{E}}_{45}\left(a_1,a_2, a_3; c_1, c_2; x,y,z\right)=\sum\limits_{m,n,p=0}^\infty\frac{(a_1)_{2m+p}(a_2)_{n}(a_3)_p}{(c_1)_{m+n}(c_2)_{p}}\frac{x^m}{m!}\frac{y^n}{n!}\frac{z^p}{p!},
\end{equation}

region of convergence:
$$ \left\{ 2\sqrt{r}+t<1,\,\,\,\,s<\infty
\right\}.
$$

System of partial differential equations:

$
\left\{
\begin{aligned}
&x(1-4x)u_{xx}+yu_{xy}-4xzu_{xz}-z^2u_{zz}\\& \,\,\,\,\,\,\,\,\,+\left[c_1-\left(4a_1+6\right)x\right]u_x-2\left(a_1+1\right)zu_z-a_1\left(a_1+1\right)u=0,\\&
yu_{yy}+xu_{xy}+\left(c_1-y\right)u_y-a_2u=0,\\&
z(1-z)u_{zz}-2xzu_{xz}-2a_3xu_x+\left[c_2-\left(a_1+a_3+1\right)z\right]u_z-a_1a_3u=0,
\end{aligned}
\right.
$

where $u\equiv \,\,{\rm{E}}_{45}\left(a_1,a_2, a_3; c_1, c_2; x,y,z\right)$.

Particular solutions:

$
{u_1} ={\rm{E}}_{45}\left(a_1,a_2, a_3; c_1, c_2; x,y,z\right) ,
$

$
{u_2} = {z^{1 - c_2}}{\rm{E}}_{45}\left(1-c_2+a_1,a_2, 1-c_1+a_3; c_1, 2-c_2; x,y,z\right).
$

\bigskip

\begin{equation}
{\rm{E}}_{46}\left(a, b; c_1, c_2; x,y,z\right)=\sum\limits_{m,n,p=0}^\infty\frac{(a)_{2m+p}(b)_{n}}{(c_1)_{m+n}(c_2)_{p}}\frac{x^m}{m!}\frac{y^n}{n!}\frac{z^p}{p!},
\end{equation}

region of convergence:
$$ \left\{ r<\frac{1}{4},\,\,\,\,s<\infty,\,\,\,\,t<\infty
\right\}.
$$

System of partial differential equations:

$
\left\{
\begin{aligned}
&x(1-4x)u_{xx}+yu_{xy}-4xzu_{xz}-z^2u_{zz}+\left[c_1-\left(4a+6\right)x\right]u_x-2\left(a+1\right)zu_z-a\left(a+1\right)u=0,\\&
yu_{yy}+xu_{xy}+\left(c_1-y\right)u_y-bu=0,\\&
zu_{zz}-2xu_x+\left(c_2-z\right)u_z-au=0,
\end{aligned}
\right.
$

where $u\equiv \,\,{\rm{E}}_{46}\left(a, b; c_1, c_2; x,y,z\right)$.

Particular solutions:

$
{u_1} ={\rm{E}}_{46}\left(a, b; c_1, c_2; x,y,z\right) ,
$

$
{u_2} = {z^{1 - c_2}}{\rm{E}}_{46}\left(1-c_2+a, b; c_1, 2-c_2; x,y,z\right).
$

\bigskip

\begin{equation}
{\rm{E}}_{47}\left(a, b; c_1, c_2; x,y,z\right)=\sum\limits_{m,n,p=0}^\infty\frac{(a)_{2m+p}(b)_p}{(c_1)_{m+n}(c_2)_{p}}\frac{x^m}{m!}\frac{y^n}{n!}\frac{z^p}{p!},
\end{equation}

region of convergence:
$$ \left\{ 2\sqrt{r}+t<1,\,\,\,\,s<\infty
\right\}.
$$

System of partial differential equations:

$
\left\{
\begin{aligned}
&x(1-4x)u_{xx}+yu_{xy}-4xzu_{xz}-z^2u_{zz}+\left[c_1-\left(4a+6\right)x\right]u_x-2\left(a+1\right)zu_z-a\left(a+1\right)u=0,\\&
yu_{yy}+xu_{xy}+c_1u_y-u=0,\\&
z(1-z)u_{zz}-2xzu_{xz}-2bxu_x+\left[c_2-\left(a+b+1\right)z\right]u_z-abu=0,
\end{aligned}
\right.
$

where $u\equiv \,\,{\rm{E}}_{47}\left(a, b; c_1, c_2; x,y,z\right)$.

Particular solutions:

$
{u_1} ={\rm{E}}_{47}\left(a, b; c_1, c_2; x,y,z\right) ,
$

$
{u_2} = {z^{1 - c_2}}{\rm{E}}_{47}\left(1-c_2+a, 1-c_2+b; c_1, 2-c_2; x,y,z\right).
$

\bigskip

\begin{equation}
{\rm{E}}_{48}\left(a; c_1, c_2; x,y,z\right)=\sum\limits_{m,n,p=0}^\infty\frac{(a)_{2m+p}}{(c_1)_{m+n}(c_2)_{p}}\frac{x^m}{m!}\frac{y^n}{n!}\frac{z^p}{p!},
\end{equation}

region of convergence:
$$ \left\{ r<\frac{1}{4},\,\,\,\,s<\infty,\,\,\,\,t<\infty
\right\}.
$$

System of partial differential equations:

$
\left\{
\begin{aligned}
&x(1-4x)u_{xx}+yu_{xy}-4xzu_{xz}-z^2u_{zz}+\left[c_1-\left(4a+6\right)x\right]u_x-2\left(a+1\right)zu_z-a\left(a+1\right)u=0,\\&
yu_{yy}+xu_{xy}+c_1u_y-u=0,\\&
zu_{zz}-2xu_x+\left(c_2-z\right)u_z-au=0,
\end{aligned}
\right.
$

where $u\equiv \,\,{\rm{E}}_{48}\left(a; c_1, c_2; x,y,z\right)$.

Particular solutions:

$
{u_1} ={\rm{E}}_{48}\left(a; c_1, c_2; x,y,z\right) ,
$

$
{u_2} = {z^{1 - c_2}}{\rm{E}}_{48}\left(1-c_2+a; c_1, 2-c_2; x,y,z\right).
$

\bigskip

\begin{equation}
{\rm{E}}_{49}\left(a,b; c_1, c_2; x,y,z\right)=\sum\limits_{m,n,p=0}^\infty\frac{(a)_{2m+n}(b)_{n+p}}{(c_1)_{m+n}(c_2)_{p}}\frac{x^m}{m!}\frac{y^n}{n!}\frac{z^p}{p!},
\end{equation}

region of convergence:
$$ \left\{ \left\{r<\frac{1}{4}, \,\,\,\, s<\frac{1}{2}+\frac{1}{2}\sqrt{1-4r}\right\}\cup\left\{s \leq \frac{1}{2}\right\},\,\,\,\,t<\infty
\right\}.
$$

System of partial differential equations:

$
\left\{
\begin{aligned}
&x(1-4x)u_{xx}-y^2u_{yy}+(1-4x)yu_{xy}+\left[c_1-\left(4a+6\right)x\right]u_x-2\left(a+1\right)yu_y-a\left(a+1\right)u=0,\\&
y(1-y)u_{yy}+x(1-2y)u_{xy}-2xzu_{xz}-yzu_{yz}-2bxu_x\\& \,\,\,\,\,\,\,\,\,+\left[c_1-\left(a+b+1\right)y\right]u_y-azu_z-abu=0,\\&
zu_{zz}-yu_y+\left(c_2-z\right)u_z-bu=0,
\end{aligned}
\right.
$

where $u\equiv \,\,{\rm{E}}_{49}\left(a,b; c_1, c_2; x,y,z\right)$.

Particular solutions:

$
{u_1} ={\rm{E}}_{49}\left(a,b; c_1, c_2; x,y,z\right) ,
$

$
{u_2} = {z^{1 - c_2}}{\rm{E}}_{49}\left(a,1-c_2+b; c_1, 2-c_2; x,y,z\right).
$

\bigskip

\begin{equation}
{\rm{E}}_{50}\left(a, b; c_1, c_2; x,y,z\right)=\sum\limits_{m,n,p=0}^\infty\frac{(a)_{2m+p}(b)_{n+p}}{(c_1)_{m+n}(c_2)_{p}}\frac{x^m}{m!}\frac{y^n}{n!}\frac{z^p}{p!},
\end{equation}

region of convergence:
$$ \left\{ 2\sqrt{r}+t<1,\,\,\,\,s<\infty
\right\}.
$$

System of partial differential equations:

$
\left\{
\begin{aligned}
&x(1-4x)u_{xx}+yu_{xy}-4xzu_{xz}-z^2u_{zz}+\left[c_1-\left(4a+6\right)x\right]u_x-2\left(a+1\right)zu_z-a\left(a+1\right)u=0,\\&
yu_{yy}+xu_{xy}+\left(c_1-y\right)u_y-zu_z-bu=0,\\&
z(1-z)u_{zz}-2xyu_{xy}-2xzu_{xz}-yzu_{yz}-2bxu_x-ayu_y+\left[c_2-\left(a+b+1\right)z\right]u_z-abu=0,
\end{aligned}
\right.
$

where ${\rm{E}}_{50}\left(a, b; c_1, c_2; x,y,z\right)$.

Particular solutions:

$
{u_1} ={\rm{E}}_{50}\left(a, b; c_1, c_2; x,y,z\right) ,
$

$
{u_2} = {z^{1 - c_2}}{\rm{E}}_{50}\left(1-c_2+a, 1-c_2+b; c_1, 2-c_2; x,y,z\right).
$

\bigskip

\begin{equation}
{\rm{E}}_{51}\left(a, b; c_1, c_2; x,y,z\right)=\sum\limits_{m,n,p=0}^\infty\frac{(a)_{2m+n+p}(b)_{n}}{(c_1)_{m+n}(c_2)_{p}}\frac{x^m}{m!}\frac{y^n}{n!}\frac{z^p}{p!},
\end{equation}

region of convergence:
$$ \left\{\left\{ r<\frac{1}{4}, \,\,\,\, s<\frac{1}{2}+\frac{1}{2}\sqrt{1-4r}\right\}\cup\left\{s \leq \frac{1}{2}\right\},\,\,\,\,t<\infty
\right\}.
$$

System of partial differential equations:

$
\left\{
\begin{aligned}
&x(1-4x)u_{xx}+(1-4x)yu_{xy}-4xzu_{xz}-2yzu_{yz}\\&\,\,\,\,\,\,\,\,\,-y^2u_{yy}-z^2u_{zz}+\left[c_1-\left(4a+6\right)x\right]u_x
-2\left(a+1\right)yu_y-2\left(a+1\right)zu_z-a\left(a+1\right)u=0,\\&
y(1-y)u_{yy}+x(1-2y)u_{xy}-yzu_{yz}-2bxu_x+\left[c_1-\left(a+b+1\right)y\right]u_y-bzu_z-abu=0,\\&
zu_{zz}-2xu_x-yu_y+\left(c_2-z\right)u_z-au=0,
\end{aligned}
\right.
$

where ${\rm{E}}_{51}\left(a, b; c_1, c_2; x,y,z\right)$,

Particular solutions:

$
{u_1} ={\rm{E}}_{51}\left(a, b; c_1, c_2; x,y,z\right) ,
$

$
{u_2} = {z^{1 - c_2}}{\rm{E}}_{51}\left(1-c_2+a, b; c_1, 2-c_2; x,y,z\right).
$

\bigskip

\begin{equation}
{\rm{E}}_{52}\left(a, b; c_1, c_2; x,y,z\right)=\sum\limits_{m,n,p=0}^\infty\frac{(a)_{2m+n+p}(b)_p}{(c_1)_{m+n}(c_2)_{p}}\frac{x^m}{m!}\frac{y^n}{n!}\frac{z^p}{p!},
\end{equation}

region of convergence:
$$ \left\{ 2\sqrt{r}+t<1,\,\,\,\,s<\infty
\right\}.
$$

System of partial differential equations:

$
\left\{
\begin{aligned}
&x(1-4x)u_{xx}+y(1-4x)u_{xy}-4xzu_{xz}-2yzu_{yz}\\&\,\,\,\,\,\,\,\,\,-y^2u_{yy}-z^2u_{zz}+\left[c_1-\left(4a+6\right)x\right]u_x
-2\left(a+1\right)yu_y-2\left(a+1\right)zu_z-a\left(a+1\right)u=0,\\&
yu_{yy}+xu_{xy}-2xu_x+\left(c_1-y\right)u_y-zu_z-au=0,\\&
z(1-z)u_{zz}-2xzu_{xz}-yzu_{yz}-2bxu_x-byu_y+\left[c_2-\left(a+b+1\right)z\right]u_z-abu=0,
\end{aligned}
\right.
$

where $u\equiv \,\,{\rm{E}}_{52}\left(a, b; c_1, c_2; x,y,z\right)$.

Particular solutions:

$
{u_1} ={\rm{E}}_{52}\left(a, b; c_1, c_2; x,y,z\right) ,
$

$
{u_2} = {z^{1 - c_2}}{\rm{E}}_{52}\left(1-c_2+a, 1-c_2+b; c_1, 2-c_2; x,y,z\right).
$

\bigskip

\begin{equation}
{\rm{E}}_{53}\left(a; c_1, c_2; x,y,z\right)=\sum\limits_{m,n,p=0}^\infty\frac{(a)_{2m+n+p}}{(c_1)_{m+n}(c_2)_{p}}\frac{x^m}{m!}\frac{y^n}{n!}\frac{z^p}{p!},
\end{equation}

region of convergence:
$$ \left\{ r<\frac{1}{4},\,\,\,\,s<\infty,\,\,\,\,t<\infty
\right\}.
$$

System of partial differential equations:

$
\left\{
\begin{aligned}
&x(1-4x)u_{xx}+y(1-4x)u_{xy}-4xzu_{xz}-2yzu_{yz}\\&\,\,\,\,\,\,\,\,\,-y^2u_{yy}-z^2u_{zz}+\left[c_1-\left(4a+6\right)x\right]u_x
-2\left(a+1\right)yu_y-2\left(a+1\right)zu_z-a\left(a+1\right)u=0,\\&
yu_{yy}+xu_{xy}-2xu_x+\left(c_1-y\right)u_y-zu_z-au=0,\\&
zu_{zz}-2xu_x-yu_y+\left(c_2-z\right)u_z-au=0,
\end{aligned}
\right.
$

where $u\equiv \,\,{\rm{E}}_{53}\left(a; c_1, c_2; x,y,z\right)$.

Particular solutions:

$
{u_1} ={\rm{E}}_{53}\left(a; c_1, c_2; x,y,z\right) ,
$

$
{u_2} = {z^{1 - c_2}}{\rm{E}}_{53}\left(1-c_2+a; c_1, 2-c_2; x,y,z\right).
$

\bigskip

\begin{equation}
{\rm{E}}_{54}\left(a,b; c_1, c_2; x,y,z\right)=\sum\limits_{m,n,p=0}^\infty\frac{(a)_{m+n}(b)_{n+2p}}{(c_1)_{m+n}(c_2)_{p}}\frac{x^m}{m!}\frac{y^n}{n!}\frac{z^p}{p!} ,
\end{equation}

region of convergence:
$$ \left\{ r<\infty,\,\,\,\,2\sqrt{t}+s<1
\right\}.
$$

System of partial differential equations:

$
\left\{
\begin{aligned}
&xu_{xx}+yu_{xy}+\left(c_1-x\right)u_x-yu_y-au=0,\\&
y(1-y)u_{yy}+x(1-y)u_{xy}-2xzu_{xz}-2yzu_{yz}-bxu_x\\& \,\,\,\,\,\,\,\,\,+\left[c_1-\left(a+b+1\right)y\right]u_y-2azu_z-abu=0,\\&
z(1-4z)u_{zz}-y^2u_{yy}-4yzu_{yz}-2\left(b+1\right)yu_y+\left[c_2-\left(4b+6\right)z\right]u_z-b\left(b+1\right)u=0,
\end{aligned}
\right.
$

where $u\equiv \,\,{\rm{E}}_{54}\left(a,b; c_1, c_2; x,y,z\right)$.

Particular solutions:

$
{u_1} ={\rm{E}}_{54}\left(a,b; c_1, c_2; x,y,z\right) ,
$

$
{u_2} = {z^{1 - c_2}}{\rm{E}}_{54}\left(a,2-2c_2+b; c_1, 2-c_2; x,y,z\right).
$

\bigskip

\begin{equation}
{\rm{E}}_{55}\left(a_1, a_2,a_3; c_1, c_2; x,y,z\right)=\sum\limits_{m,n,p=0}^\infty\frac{(a_1)_{n+2p}(a_2)_m(a_3)_m}{(c_1)_{m+n}(c_2)_{p}}\frac{x^m}{m!}\frac{y^n}{n!}\frac{z^p}{p!},
\end{equation}

region of convergence:
$$ \left\{r<1, \,\,\,s<\infty,\,\,\,t<\frac{1}{4}
\right\}.
$$

System of partial differential equations:

$
\left\{
\begin{aligned}
&x(1-x)u_{xx}+yu_{xy}+\left[c_1-\left(a_2+a_3+1\right)x\right]u_x-a_2a_3u=0,\\&
yu_{yy}+xu_{xy}+\left(c_1-y\right)u_y-2zu_z-a_1u=0,\\&
z(1-4z)u_{zz}-y^2u_{yy}-4yzu_{yz}-2\left(a_1+1\right)yu_y+\left[c_2-\left(4a_1+6\right)z\right]u_z-a_1\left(a_1+1\right)u=0,
\end{aligned}
\right.
$

where $u\equiv \,\,{\rm{E}}_{55}\left(a_1, a_2,a_3; c_1, c_2; x,y,z\right)$.

Particular solutions:

$
{u_1} ={\rm{E}}_{55}\left(a_1, a_2,a_3; c_1, c_2; x,y,z\right) ,
$

$
{u_2} = {z^{1 - c_2}}{\rm{E}}_{55}\left(2-2c_2+a_1, a_2,a_3; c_1, 2-c_2; x,y,z\right).
$

\bigskip

\begin{equation}
{\rm{E}}_{56}\left(a_1,a_2, a_3; c_1, c_2; x,y,z\right)=\sum\limits_{m,n,p=0}^\infty\frac{(a_1)_{n+2p}(a_2)_{m}(a_3)_n}{(c_1)_{m+n}(c_2)_{p}}\frac{x^m}{m!}\frac{y^n}{n!}\frac{z^p}{p!},
\end{equation}

region of convergence:
$$ \left\{ r<\infty,\,\,\,\,2\sqrt{t}+s<1
\right\}.
$$

System of partial differential equations:

$
\left\{
\begin{aligned}
&xu_{xx}+yu_{xy}+\left(c_1-x\right)u_x-a_2u=0,\\&
y(1-y)u_{yy}+xu_{xy}-2yzu_{yz}+\left[c_1-\left(a_1+a_3+1\right)y\right]u_y-2a_3zu_z-a_1a_3u=0,\\&
z(1-4z)u_{zz}-y^2u_{yy}-4yzu_{yz}-2\left(a_1+1\right)yu_y+\left[c_2-\left(4a_1+6\right)z\right]u_z-a_1\left(a_1+1\right)u=0,
\end{aligned}
\right.
$

where $u\equiv \,\,{\rm{E}}_{56}\left(a_1,a_2, a_3; c_1, c_2; x,y,z\right)$.

Particular solutions:

$
{u_1} ={\rm{E}}_{56}\left(a_1,a_2, a_3; c_1, c_2; x,y,z\right) ,
$

$
{u_2} = {z^{1 - c_2}}{\rm{E}}_{56}\left(2-2c_2+a_1,a_2, a_3; c_1, 2-c_2; x,y,z\right).
$

\bigskip

\begin{equation}
{\rm{E}}_{57}\left(a, b; c_1, c_2; x,y,z\right)=\sum\limits_{m,n,p=0}^\infty\frac{(a)_{n+2p}(b)_m}{(c_1)_{m+n}(c_2)_{p}}\frac{x^m}{m!}\frac{y^n}{n!}\frac{z^p}{p!},
\end{equation}

region of convergence:
$$ \left\{ r<\infty,\,\,\,s<\infty,\,\,\,\,t<\frac{1}{4}
\right\}.
$$

System of partial differential equations:

$
\left\{
\begin{aligned}
&xu_{xx}+yu_{xy}+\left(c_1-x\right)u_x-bu=0,\\&
yu_{yy}+xu_{xy}+\left(c_1-y\right)u_y-2zu_z-au=0,\\&
z(1-4z)u_{zz}-y^2u_{yy}-4yzu_{yz}-2\left(a+1\right)yu_y+\left[c_2-\left(4a+6\right)z\right]u_z-a\left(a+1\right)u=0,
\end{aligned}
\right.
$

where $u\equiv \,\,{\rm{E}}_{57}\left(a, b; c_1, c_2; x,y,z\right)$.

Particular solutions:

$
{u_1} ={\rm{E}}_{57}\left(a, b; c_1, c_2; x,y,z\right) ,
$

$
{u_2} = {z^{1 - c_2}}{\rm{E}}_{57}\left(2-2c_2+a, b; c_1, 2-c_2; x,y,z\right).
$

\bigskip

\begin{equation}
{\rm{E}}_{58}\left(a, b; c_1, c_2; x,y,z\right)=\sum\limits_{m,n,p=0}^\infty\frac{(a)_{n+2p}(b)_{n}}{(c_1)_{m+n}(c_2)_{p}}\frac{x^m}{m!}\frac{y^n}{n!}\frac{z^p}{p!},
\end{equation}

region of convergence:
$$ \left\{ r<\infty,\,\,\,2\sqrt{t}+s<1
\right\}.
$$

System of partial differential equations:

$
\left\{
\begin{aligned}
&xu_{xx}+yu_{xy}+c_1u_x-u=0,\\&
y(1-y)u_{yy}+xu_{xy}-2yzu_{yz}+\left[c_1-\left(a+b+1\right)y\right]u_y-2bzu_z-abu=0,\\&
z(1-4z)u_{zz}-y^2u_{yy}-4yzu_{yz}-2\left(a+1\right)yu_y+\left[c_2-\left(4a+6\right)z\right]u_z-a\left(a+1\right)u=0,
\end{aligned}
\right.
$

where $u\equiv \,\,{\rm{E}}_{58}\left(a, b; c_1, c_2; x,y,z\right)$.

Particular solutions:

$
{u_1} ={\rm{E}}_{58}\left(a, b; c_1, c_2; x,y,z\right) ,
$

$
{u_2} = {z^{1 - c_2}}{\rm{E}}_{58}\left(2-2c_2+a, b; c_1, 2-c_2; x,y,z\right).
$

\bigskip

\begin{equation}
{\rm{E}}_{59}\left(a; c_1, c_2; x,y,z\right)=\sum\limits_{m,n,p=0}^\infty\frac{(a)_{n+2p}}{(c_1)_{m+n}(c_2)_{p}}\frac{x^m}{m!}\frac{y^n}{n!}\frac{z^p}{p!},
\end{equation}

region of convergence:
$$ \left\{ r<\infty,\,\,\,\,s<\infty,\,\,\,t<\frac{1}{4}
\right\}.
$$

System of partial differential equations:

$
\left\{
\begin{aligned}
&xu_{xx}+yu_{xy}+c_1u_x-u=0,\\&
yu_{yy}+xu_{xy}+\left(c_1-y\right)u_y-2zu_z-au=0,\\&
z(1-4z)u_{zz}-y^2u_{yy}-4yzu_{yz}-2\left(a+1\right)yu_y+\left[c_2-\left(4a+6\right)z\right]u_z-a\left(a+1\right)u=0,
\end{aligned}
\right.
$

where $u\equiv \,\,{\rm{E}}_{59}\left(a; c_1, c_2; x,y,z\right)$.

Particular solutions:

$
{u_1} ={\rm{E}}_{59}\left(a; c_1, c_2; x,y,z\right) ,
$

$
{u_2} = {z^{1 - c_2}}{\rm{E}}_{59}\left(2-2c_2+a; c_1, 2-c_2; x,y,z\right).
$

\bigskip

\begin{equation}
{\rm{E}}_{60}\left(a, b; c_1, c_2; x,y,z\right)=\sum\limits_{m,n,p=0}^\infty\frac{(a)_{m+n+2p}(b)_{m}}{(c_1)_{m+n}(c_2)_{p}}\frac{x^m}{m!}\frac{y^n}{n!}\frac{z^p}{p!},
\end{equation}

region of convergence:
$$ \left\{ r+2\sqrt{t}<1,\,\,\,s<\infty
\right\}.
$$

System of partial differential equations:

$
\left\{
\begin{aligned}
&x(1-x)u_{xx}+(1-x)yu_{xy}-2xzu_{xz}+\left[c_1-\left(a+b+1\right)x\right]u_x-byu_y-2bzu_z-abu=0,\\&
yu_{yy}+xu_{xy}-xu_x+\left(c_1-y\right)u_y-2zu_z-au=0,\\&
z(1-4z)u_{zz}-x^2u_{xx}-y^2u_{yy}-2xyu_{xy}-4xzu_{xz}\\&\,\,\,\,\,\,\,\,\,-4yzu_{yz}-2\left(a+1\right)xu_x-2\left(a+1\right)yu_y+\left[c_2-\left(4a+6\right)z\right]u_z-a\left(a+1\right)u=0,
\end{aligned}
\right.
$

where $u\equiv \,\,{\rm{E}}_{60}\left(a, b; c_1, c_2; x,y,z\right)$.

Particular solutions:

$
{u_1} ={\rm{E}}_{60}\left(a, b; c_1, c_2; x,y,z\right) ,
$

$
{u_2} = {z^{1 - c_2}}{\rm{E}}_{60}\left(2-2c_2+a, b; c_1, 2-c_2; x,y,z\right).
$

\bigskip

\begin{equation}
{\rm{E}}_{61}\left(a; c_1, c_2; x,y,z\right)=\sum\limits_{m,n,p=0}^\infty\frac{(a)_{m+2n+2p}}{(c_1)_{m+n}(c_2)_{p}}\frac{x^m}{m!}\frac{y^n}{n!}\frac{z^p}{p!},
\end{equation}

region of convergence:
$$ \left\{ r<\infty,\,\,\,\sqrt{s}+\sqrt{t}<\frac{1}{2}
\right\}.
$$

System of partial differential equations:

$
\left\{
\begin{aligned}
&xu_{xx}+yu_{xy}+\left(c_1-x\right)u_x-2yu_{y}-2zu_z-au=0,\\&
y(1-4y)u_{yy}-x^2u_{xx}-4z^2u_{zz}+x(1-4y)u_{xy}\\&\,\,\,\,\,-4xzu_{xz}-8yzu_{yz}-2(a+1)xu_x+\left[c_1-(4a+6)y\right]u_y-(4a+6)zu_z-a(a+1)u=0,\\&
z(1-4z)u_{zz}-x^2u_{xx}-4y^2u_{yy}-4xyu_{xy}\\&\,\,\,\,\,-4xzu_{xz}-8yzu_{yz}-2\left(a+1\right)xu_x-\left(4a+6\right)yu_y+\left[c_2-\left(4a+6\right)z\right]u_z-a\left(a+1\right)u=0,
\end{aligned}
\right.
$

where $u\equiv \,\,{\rm{E}}_{61}\left(a; c_1, c_2; x,y,z\right)$.

Particular solutions:

$
{u_1} ={\rm{E}}_{61}\left(a; c_1, c_2; x,y,z\right) ,
$

$
{u_2} = {z^{1 - c_2}}{\rm{E}}_{61}\left(2-2c_2+a; c_1, 2-c_2; x,y,z\right).
$

\bigskip

\begin{equation}
{\rm{E}}_{62}\left(a_1,a_2,a_3; c_1, c_2, c_3; x,y,z\right)=\sum\limits_{m,n,p=0}^\infty\frac{(a_1)_{m+n}(a_2)_{n+p}(a_3)_m}{(c_1)_{m}(c_2)_{n}(c_3)_p}\frac{x^m}{m!}\frac{y^n}{n!}\frac{z^p}{p!},
\end{equation}

first appearance of this function in the literature: [21],\, ${}_3\Phi_K^{(1)}$,

region of convergence:
$$ \left\{ r+s<1,\,\,\,t<\infty
\right\}.
$$

System of partial differential equations:

$
\left\{
\begin{aligned}
&x(1-x)u_{xx}-xyu_{xy}+\left[c_1-\left(a_1+a_3+1\right)x\right]u_x-a_3yu_{y}-a_1a_3u=0,\\&
y(1-y)u_{yy}-xyu_{xy}-xzu_{xz}-yzu_{yz}-a_2xu_x+\left[c_2-\left(a_1+a_2+1\right)y\right]u_y-a_1zu_z-a_1a_2u=0,\\&
zu_{zz}-yu_y+\left(c_3-z\right)u_z-a_2u=0,
\end{aligned}
\right.
$

where $u\equiv \,\,{\rm{E}}_{62}\left(a_1,a_2,a_3; c_1, c_2, c_3; x,y,z\right)$.

Particular solutions:

$
{u_1} ={\rm{E}}_{62}\left(a_1,a_2,a_3; c_1, c_2, c_3; x,y,z\right) ,
$

$
{u_2} = {x^{1 - c_1}}{\rm{E}}_{62}\left(1-c_1+a_1,a_2,1-c_1+a_3; 2-c_1, c_2, c_3; x,y,z\right),
$

$
{u_3} = {y^{1 - c_2}}{\rm{E}}_{62}\left(1-c_2+a_1,1-c_2+a_2,a_3; c_1, 2-c_2, c_3; x,y,z\right),
$

$
{u_4} = {z^{1 - c_3}}{\rm{E}}_{62}\left(a_1,1-c_3+a_2,a_3; c_1, c_2, 2-c_3; x,y,z\right),
$

$
{u_5} = {x^{1 - c_1}}{y^{1 - c_2}}{\rm{E}}_{62}\left(2-c_1-c_2+a_1,1-c_2+a_2,a_3; 2-c_1, 2-c_2, c_3; x,y,z\right),
$

$
{u_6} = {y^{1 - c_2}}{z^{1 - c_3}}{\rm{E}}_{62}\left(1-c_2+a_1,1-c_2+a_2,a_3; c_1, 2-c_2, 2-c_3; x,y,z\right),
$

$
{u_7} = {x^{1 - c_1}}{z^{1 - c_3}}{\rm{E}}_{62}\left(1-c_1+a_1,1-c_3+a_2,1-c_1+a_3; 2-c_1, c_2, 2-c_3; x,y,z\right),
$

$
{u_8} = {x^{1 - c_1}}{y^{1 - c_2}}{z^{1 - c_3}}\times$

$\,\,\,\,\,\,\,\,\,\times{\rm{E}}_{62}\left(2-c_1-c_2+a_1,2-c_2-c_3+a_2,1-c_1+a_3; 2-c_1, 2-c_2, 2-c_3; x,y,z\right).
$

\bigskip

\begin{equation}
{\rm{E}}_{63}\left(a, b; c_1, c_2, c_3; x,y,z\right)=\sum\limits_{m,n,p=0}^\infty\frac{(a)_{m+n}(b)_{n+p}}{(c_1)_{m}(c_2)_{n}(c_3)_p}\frac{x^m}{m!}\frac{y^n}{n!}\frac{z^p}{p!},
\end{equation}

first appearance of this function in the literature: [21],\, ${}_3\Phi_K^{(2)}$,

region of convergence:
$$ \left\{ r<\infty,\,\,\,s<1,\,\,\,t<\infty
\right\}.
$$

System of partial differential equations:

$
\left\{
\begin{aligned}
&xu_{xx}+\left(c_1-x\right)u_x-yu_{y}-au=0,\\&
y(1-y)u_{yy}-xyxu_{xy}-xzu_{xz}-yzu_{yz}-bxu_x+\left[c_2-\left(a+b+1\right)y\right]u_y-azu_z-abu=0,\\&
zu_{zz}-yu_y+\left(c_3-z\right)u_z-bu=0,
\end{aligned}
\right.
$

where $u\equiv \,\,{\rm{E}}_{63}\left(a, b; c_1, c_2, c_3; x,y,z\right)$.

Particular solutions:

$
{u_1} ={\rm{E}}_{63}\left(a, b; c_1, c_2, c_3; x,y,z\right) ,
$

$
{u_2} = {x^{1 - c_1}}{\rm{E}}_{63}\left(1-c_1+a, b; 2-c_1, c_2, c_3; x,y,z\right),
$

$
{u_3} = {y^{1 - c_2}}{\rm{E}}_{63}\left(1-c_2+a, 1-c_2+b; c_1, 2-c_2, c_3; x,y,z\right),
$

$
{u_4} = {z^{1 - c_3}}{\rm{E}}_{63}\left(a, 1-c_3+b; c_1, c_2, 2-c_3; x,y,z\right),
$

$
{u_5} = {x^{1 - c_1}}{y^{1 - c_2}}{\rm{E}}_{63}\left(2-c_1-c_2+a, 1-c_2+b; 2-c_1, 2-c_2, c_3; x,y,z\right),
$

$
{u_6} = {y^{1 - c_2}}{z^{1 - c_3}}{\rm{E}}_{63}\left(1-c_2+a, 1-c_3+b; c_1, 2-c_2, 2-c_3; x,y,z\right),
$

$
{u_7} = {x^{1 - c_1}}{z^{1 - c_3}}{\rm{E}}_{63}\left(1-c_1+a, 1-c_3+b; 2-c_1, c_2, 2-c_3; x,y,z\right),
$

$
{u_8} = {x^{1 - c_1}}{y^{1 - c_2}}{z^{1 - c_3}}{\rm{E}}_{63}\left(2-c_1-c_2+a, 2-c_2-c_3+b; 2-c_1, 2-c_2, 2-c_3; x,y,z\right).
$

\bigskip

\begin{equation}
{\rm{E}}_{64}\left(a_1,a_2,a_3; c_1, c_2, c_3; x,y,z\right)=\sum\limits_{m,n,p=0}^\infty\frac{(a_1)_{m+n+p}(a_2)_{m}(a_3)_n}{(c_1)_{m}(c_2)_{n}(c_3)_p}\frac{x^m}{m!}\frac{y^n}{n!}\frac{z^p}{p!},
\end{equation}

first appearance of this function in the literature: [21],\, ${}_3\Phi_A^{(1)}$, \,\,[9], \,\,$F_{A1}$,

region of convergence:
$$ \left\{ r+s<1,\,\,\,t<\infty
\right\}.
$$

System of partial differential equations:

$
\left\{
\begin{aligned}
&x(1-x)u_{xx}-xyu_{xy}-xzu_{xz}+\left[c_1-\left(a_1+a_2+1\right)x\right]u_x-a_2yu_{y}-a_2zu_{z}-a_1a_2u=0,\\&
y(1-y)u_{yy}-xyu_{xy}-yzu_{yz}-a_3xu_x+\left[c_2-\left(a_1+a_3+1\right)y\right]u_y-a_3zu_z-a_1a_3u=0,\\&
zu_{zz}-xu_{x}-yu_y+\left(c_3-z\right)u_z-a_1u=0,
\end{aligned}
\right.
$

where $u\equiv \,\,{\rm{E}}_{64}\left(a_1,a_2,a_3; c_1, c_2, c_3; x,y,z\right)$.

Particular solutions:

$
{u_1} ={\rm{E}}_{64}\left(a_1,a_2,a_3; c_1, c_2, c_3; x,y,z\right) ,
$

$
{u_2} = {x^{1 - c_1}}{\rm{E}}_{64}\left(1-c_1+a_1,1-c_1+a_2,a_3; 2-c_1, c_2, c_3; x,y,z\right),
$

$
{u_3} = {y^{1 - c_2}}{\rm{E}}_{64}\left(1-c_2+a_1,a_2,1-c_2+a_3; c_1, 2-c_2, c_3; x,y,z\right),
$

$
{u_4} = {z^{1 - c_3}}{\rm{E}}_{64}\left(1-c_3+a_1,a_2,a_3; c_1, c_2, 2-c_3; x,y,z\right),
$

$
{u_5} = {x^{1 - c_1}}{y^{1 - c_2}}{\rm{E}}_{64}\left(2-c_1-c_2+a_1,1-c_1+a_2,1-c_2+a_3; 2-c_1, 2-c_2, c_3; x,y,z\right),
$

$
{u_6} = {y^{1 - c_2}}{z^{1 - c_3}}{\rm{E}}_{64}\left(2-c_2-c_3+a_1,a_2,1-c_2+a_3; c_1, 2-c_2, 2-c_3; x,y,z\right),
$

$
{u_7} = {x^{1 - c_1}}{z^{1 - c_3}}{\rm{E}}_{64}\left(2-c_1-c_3+a_1,1-c_1+a_2,a_3; 2-c_1, c_2, 2-c_3; x,y,z\right),
$

$
{u_8} = {x^{1 - c_1}}{y^{1 - c_2}}{z^{1 - c_3}}\times$

$
\,\,\,\,\,\,\,\,\,\,\times{\rm{E}}_{64}\left(3-c_1-c_2-c_3+a_1,1-c_1+a_2,1-c_2+a_3; 2-c_1, 2-c_2, 2-c_3; x,y,z\right).
$

\bigskip

\begin{equation}
{\rm{E}}_{65}\left(a, b; c_1, c_2, c_3; x,y,z\right)=\sum\limits_{m,n,p=0}^\infty\frac{(a)_{m+n+p}(b)_{m}}{(c_1)_{m}(c_2)_{n}(c_3)_p}\frac{x^m}{m!}\frac{y^n}{n!}\frac{z^p}{p!},
\end{equation}

first appearance of this function in the literature: [21],\, ${}_3\Phi_A^{(2)}$, \,\,[9], \,\,$F_{A2}$,

region of convergence:
$$ \left\{ r<1,\,\,\,s<\infty,\,\,\,t<\infty
\right\}.
$$

System of partial differential equations:

$
\left\{
\begin{aligned}
&x(1-x)u_{xx}-xyu_{xy}-xzu_{xz}+\left[c_1-\left(a+b+1\right)x\right]u_x-byu_{y}-bzu_{z}-abu=0,\\&
yu_{yy}-xu_x+\left(c_2-y\right)u_y-zu_z-au=0,\\&
zu_{zz}-xu_{x}-yu_y+\left(c_3-z\right)u_z-au=0,
\end{aligned}
\right.
$

where $u\equiv \,\,{\rm{E}}_{65}\left(a, b; c_1, c_2, c_3; x,y,z\right)$.

Particular solutions:

$
{u_1} ={\rm{E}}_{65}\left(a, b; c_1, c_2, c_3; x,y,z\right) ,
$

$
{u_2} = {x^{1 - c_1}}{\rm{E}}_{65}\left(1-c_1+a, 1-c_1+b; 2-c_1, c_2, c_3; x,y,z\right),
$

$
{u_3} = {y^{1 - c_2}}{\rm{E}}_{65}\left(1-c_2+a, b; c_1, 2-c_2, c_3; x,y,z\right),
$

$
{u_4} = {z^{1 - c_3}}{\rm{E}}_{65}\left(1-c_3+a, b; c_1, c_2, 2-c_3; x,y,z\right),
$

$
{u_5} = {x^{1 - c_1}}{y^{1 - c_2}}{\rm{E}}_{65}\left(2-c_1-c_2+a, 1-c_1+b; 2-c_1, 2-c_2, c_3; x,y,z\right),
$

$
{u_6} = {y^{1 - c_2}}{z^{1 - c_3}}{\rm{E}}_{65}\left(2-c_2-c_3+a, b; c_1, 2-c_2, 2-c_3; x,y,z\right),
$

$
{u_7} = {x^{1 - c_1}}{z^{1 - c_3}}{\rm{E}}_{65}\left(2-c_1-c_3+a, 1-c_1+b; 2-c_1, c_2, 2-c_3; x,y,z\right),
$

$
{u_8} = {x^{1 - c_1}}{y^{1 - c_2}}{z^{1 - c_3}}{\rm{E}}_{65}\left(3-c_1-c_2-c_3+a, 1-c_1+b; 2-c_1, 2-c_2, 2-c_3; x,y,z\right).
$

\bigskip

\begin{equation}
{\rm{E}}_{66}\left(a; c_1, c_2, c_3; x,y,z\right)=\sum\limits_{m,n,p=0}^\infty\frac{(a)_{m+n+p}}{(c_1)_{m}(c_2)_{n}(c_3)_p}\frac{x^m}{m!}\frac{y^n}{n!}\frac{z^p}{p!}
,
\end{equation}

first appearance of this function in the literature: [21],\, ${}_3\Phi_A^{(3)}$, \,\,[9], \,\,$F_{A3}$,

region of convergence:
$$ \left\{ r<\infty,\,\,\,s<\infty,\,\,\,t<\infty
\right\}.
$$

System of partial differential equations:

$
\left\{
\begin{aligned}
&xu_{xx}+\left(c_1-x\right)u_x-yu_{y}-zu_{z}-au=0,\\&
yu_{yy}-xu_x+\left(c_2-y\right)u_y-zu_z-au=0,\\&
zu_{zz}-xu_{x}-yu_y+\left(c_3-z\right)u_z-au=0,
\end{aligned}
\right.
$

where $u\equiv \,\,{\rm{E}}_{66}\left(a; c_1, c_2, c_3; x,y,z\right)$.

Particular solutions:

$
{u_1} ={\rm{E}}_{66}\left(a; c_1, c_2, c_3; x,y,z\right),
$

$
{u_2} = {x^{1 - c_1}}{\rm{E}}_{66}\left(1-c_1+a; 2-c_1, c_2, c_3; x,y,z\right),
$

$
{u_3} = {y^{1 - c_2}}{\rm{E}}_{66}\left(1-c_2+a; c_1, 2-c_2, c_3; x,y,z\right),
$

$
{u_4} = {z^{1 - c_3}}{\rm{E}}_{66}\left(1-c_3+a; c_1, c_2, 2-c_3; x,y,z\right),
$

$
{u_5} = {x^{1 - c_1}}{y^{1 - c_2}}{\rm{E}}_{66}\left(2-c_1-c_2+a; 2-c_1, 2-c_2, c_3; x,y,z\right),
$

$
{u_6} = {y^{1 - c_2}}{z^{1 - c_3}}{\rm{E}}_{66}\left(2-c_2-c_3+a; c_1, 2-c_2, 2-c_3; x,y,z\right),
$

$
{u_7} = {x^{1 - c_1}}{z^{1 - c_3}}{\rm{E}}_{66}\left(2-c_1-c_3+a; 2-c_1, c_2, 2-c_3; x,y,z\right),
$

$
{u_8} = {x^{1 - c_1}}{y^{1 - c_2}}{z^{1 - c_3}}{\rm{E}}_{66}\left(3-c_1-c_2-c_3+a; 2-c_1, 2-c_2, 2-c_3; x,y,z\right).
$

\bigskip

\begin{equation}
{\rm{E}}_{67}\left(a,b; c_1, c_2,c_3; x,y,z\right)=\sum\limits_{m,n,p=0}^\infty\frac{(a)_{m+n+p}(b)_{m+n}}{(c_1)_{m}(c_2)_{n}(c_3)_p}\frac{x^m}{m!}\frac{y^n}{n!}\frac{z^p}{p!},
\end{equation}

first appearance of this function in the literature: [21],\, ${}_3\Phi_E^{(1)}$,

region of convergence:
$$ \left\{ \sqrt{r}+\sqrt{s}<1,\,\,\,t<\infty
\right\}.
$$

System of partial differential equations:

$
\left\{
\begin{aligned}
&x(1-x)u_{xx}-2xyu_{xy}\\&\,\,\,\,-xzu_{xz}-yzu_{yz}-y^2u_{yy} +\left[c_1-\left(a+b+1\right)x\right]u_x-\left(a+b+1\right)yu_y-bzu_{z}-abu=0,\\&
y(1-y)u_{yy}-2xyu_{xy}\\&\,\,\,\, -xzu_{xz}-yzu_{yz}-x^2u_{xx}- \left(a+b+1\right)xu_x
+\left[c_2-\left(a+b+1\right)y\right]u_y-bzu_z-abu=0,\\&
zu_{zz}-xu_{x}-yu_y+\left(c_3-z\right)u_z-au=0,
\end{aligned}
\right.
$

where $u\equiv \,\,{\rm{E}}_{67}\left(a,b; c_1, c_2,c_3; x,y,z\right)$.

Particular solutions:

$
{u_1} ={\rm{E}}_{67}\left(a,b; c_1, c_2,c_3; x,y,z\right) ,
$

$
{u_2} = {x^{1 - c_1}}{\rm{E}}_{67}\left(1-c_1+a,1-c_1+b; 2-c_1, c_2,c_3; x,y,z\right),
$

$
{u_3} = {y^{1 - c_2}}{\rm{E}}_{67}\left(1-c_2+a,1-c_2+b; c_1, 2-c_2,c_3; x,y,z\right),
$

$
{u_4} = {z^{1 - c_3}}{\rm{E}}_{67}\left(1-c_3+a,b; c_1, c_2,2-c_3; x,y,z\right),
$

$
{u_5} = {x^{1 - c_1}}{y^{1 - c_2}}{\rm{E}}_{67}\left(2-c_1-c_2+a,2-c_1-c_2+b; 2-c_1, 2-c_2,c_3; x,y,z\right),
$

$
{u_6} = {y^{1 - c_2}}{z^{1 - c_3}}{\rm{E}}_{67}\left(2-c_2-c_3+a,1-c_2+b; c_1, 2-c_2,2-c_3; x,y,z\right),
$

$
{u_7} = {x^{1 - c_1}}{z^{1 - c_3}}{\rm{E}}_{67}\left(2-c_1-c_3+a,1-c_1+b; 2-c_1, c_2,2-c_3; x,y,z\right),
$

$
{u_8} = {x^{1 - c_1}}{y^{1 - c_2}}{z^{1 - c_3}}{\rm{E}}_{67}\left(3-c_1-c_2-c_3+a,2-c_1-c_2+b; 2-c_1, 2-c_2,2-c_3; x,y,z\right).
$

\bigskip

\begin{equation}
{\rm{E}}_{68}\left(a,b; c_1, c_2, c_3; x,y,z\right)=\sum\limits_{m,n,p=0}^\infty\frac{(a)_{2m+n}(b)_{n+p}}{(c_1)_{m}(c_2)_{n}(c_3)_{p}}\frac{x^m}{m!}\frac{y^n}{n!}\frac{z^p}{p!},
\end{equation}

region of convergence:
$$ \left\{ \sqrt{r}+\sqrt{s}<1,\,\,\,t<\infty
\right\}.
$$

System of partial differential equations:

$
\left\{
\begin{aligned}
&x(1-4x)u_{xx}-4xyu_{xy}-y^2u_{yy}+ \left[c_1-\left(4a+6\right)x\right]u_x-2\left(a+1\right)yu_y-a\left(a+1\right)u=0,\\&
y(1-y)u_{yy}-2xyu_{xy}-2xzu_{xz}-yzu_{yz} -2bxu_x
+\left[c_2-\left(a+b+1\right)y\right]u_y-azu_z-abu=0,\\&
zu_{zz}-yu_y+\left(c_3-z\right)u_z-bu=0,
\end{aligned}
\right.
$

where $u\equiv \,\,{\rm{E}}_{68}\left(a,b; c_1, c_2, c_3; x,y,z\right)$.

Particular solutions:

$
{u_1} ={\rm{E}}_{68}\left(a,b; c_1, c_2, c_3; x,y,z\right) ,
$

$
{u_2} = {x^{1 - c_1}}{\rm{E}}_{68}\left(2-2c_1+a,b; 2-c_1, c_2, c_3; x,y,z\right),
$

$
{u_3} = {y^{1 - c_2}}{\rm{E}}_{68}\left(1-c_2+a,1-c_2+b; c_1, 2-c_2, c_3; x,y,z\right),
$

$
{u_4} = {z^{1 - c_3}}{\rm{E}}_{68}\left(a,1-c_3+b; c_1, c_2, 2-c_3; x,y,z\right),
$

$
{u_5} = {x^{1 - c_1}}{y^{1 - c_2}}{\rm{E}}_{68}\left(3-2c_1-c_2+a,1-c_2+b; 2-c_1, 2-c_2, c_3; x,y,z\right),
$

$
{u_6} = {y^{1 - c_2}}{z^{1 - c_3}}{\rm{E}}_{68}\left(1-c_2+a,2-c_2-c_3+b; c_1, 2-c_2, 2-c_3; x,y,z\right),
$

$
{u_7} = {x^{1 - c_1}}{z^{1 - c_3}}{\rm{E}}_{68}\left(2-2c_1+a,1-c_3+b; 2-c_1, c_2, 2-c_3; x,y,z\right),
$

$
{u_8} = {x^{1 - c_1}}{y^{1 - c_2}}{z^{1 - c_3}}{\rm{E}}_{68}\left(3-2c_1-c_2+a,2-c_2-c_3+b; 2-c_1, 2-c_2, 2-c_3; x,y,z\right).
$

\bigskip

\begin{equation}
{\rm{E}}_{69}\left(a, b; c_1, c_2,c_3; x,y,z\right)=\sum\limits_{m,n,p=0}^\infty\frac{(a)_{2m+n+p}(b)_{n}}{(c_1)_{m}(c_2)_{n}(c_3)_p}\frac{x^m}{m!}\frac{y^n}{n!}\frac{z^p}{p!},
\end{equation}

region of convergence:
$$ \left\{ 2\sqrt{r}+s<1,\,\,\,t<\infty
\right\}.
$$

System of partial differential equations:

$
\left\{
\begin{aligned}
&x(1-4x)u_{xx}-4xyu_{xy}-4xzu_{xz}-2yzu_{yz}\\&\,\,\,\,\,\,\,\,\,-y^2u_{yy}-z^2u_{zz}+ \left[c_1-\left(4a+6\right)x\right]u_x-2\left(a+1\right)yu_y-2\left(a+1\right)zu_z-a\left(a+1\right)u=0,\\&
y(1-y)u_{yy}-2xyu_{xy}-yzu_{yz} -2bxu_x
+\left[c_2-\left(a+b+1\right)y\right]u_y-bzu_z-abu=0,\\&
zu_{zz}-2xu_x-yu_y+\left(c_3-z\right)u_z-au=0,
\end{aligned}
\right.
$

where $u\equiv \,\,{\rm{E}}_{69}\left(a, b; c_1, c_2,c_3; x,y,z\right)$.

Particular solutions:

$
{u_1} ={\rm{E}}_{69}\left(a, b; c_1, c_2,c_3; x,y,z\right) ,
$

$
{u_2} = {x^{1 - c_1}}{\rm{E}}_{69}\left(2-2c_1+a, b; 2-c_1, c_2,c_3; x,y,z\right);
$

$
{u_3} = {y^{1 - c_2}}{\rm{E}}_{69}\left(1-c_2+a, 1-c_2+b; c_1, 2-c_2,c_3; x,y,z\right);
$

$
{u_4} = {z^{1 - c_3}}{\rm{E}}_{69}\left(1-c_3+a, b; c_1, c_2,2-c_3; x,y,z\right);
$

$
{u_5} = {x^{1 - c_1}}{y^{1 - c_2}}{\rm{E}}_{69}\left(3-2c_1-c_2+a, 1-c_2+b; 2-c_1, 2-c_2,c_3; x,y,z\right);
$

$
{u_6} = {y^{1 - c_2}}{z^{1 - c_3}}{\rm{E}}_{69}\left(2-c_2-c_3+a, 1-c_2+b; c_1, 2-c_2,2-c_3; x,y,z\right);
$

$
{u_7} = {x^{1 - c_1}}{z^{1 - c_3}}{\rm{E}}_{69}\left(3-2c_1-c_3+a, b; 2-c_1, c_2,2-c_3; x,y,z\right);
$

$
{u_8} = {x^{1 - c_1}}{y^{1 - c_2}}{z^{1 - c_3}}{\rm{E}}_{69}\left(4-2c_1-c_2-c_3+a, 1-c_2+b; 2-c_1, 2-c_2,2-c_3; x,y,z\right).
$

\bigskip

\begin{equation}
{\rm{E}}_{70}\left(a; c_1, c_2,c_3; x,y,z\right)=\sum\limits_{m,n,p=0}^\infty\frac{(a)_{2m+n+p}}{(c_1)_{m}(c_2)_{n}(c_3)_p}\frac{x^m}{m!}\frac{y^n}{n!}\frac{z^p}{p!},
\end{equation}

region of convergence:
$$ \left\{ r<\frac{1}{4},\,\,\,s<\infty,\,\,\,t<\infty
\right\}.
$$

System of partial differential equations:

$
\left\{
\begin{aligned}
&x(1-4x)u_{xx}-4xyu_{xy}-4xzu_{xz}-2yzu_{yz}-y^2u_{yy}\\ &\,\,\,\,-z^2u_{zz}+ \left[c_1-\left(4a+6\right)x\right]u_x-2\left(a+1\right)yu_y-2\left(a+1\right)zu_z-a\left(a+1\right)u=0,\\&
yu_{yy}-2xu_x
+\left(c_2-y\right)u_y-zu_z-au=0,\\&
zu_{zz}-2xu_x-yu_y+\left(c_3-z\right)u_z-au=0,
\end{aligned}
\right.
$

where $u\equiv \,\,{\rm{E}}_{70}\left(a; c_1, c_2,c_3; x,y,z\right)$.

Particular solutions:

$
{u_1} ={\rm{E}}_{70}\left(a; c_1, c_2,c_3; x,y,z\right) ,
$

$
{u_2} = {x^{1 - c_1}}{\rm{E}}_{70}\left(2-2c_1+a; 2-c_1, c_2,c_3; x,y,z\right),
$

$
{u_3} = {y^{1 - c_2}}{\rm{E}}_{70}\left(1-c_2+a; c_1, 2-c_2,c_3; x,y,z\right),
$

$
{u_4} = {z^{1 - c_3}}{\rm{E}}_{70}\left(1-c_3+a; c_1, c_2,2-c_3; x,y,z\right),
$

$
{u_5} = {x^{1 - c_1}}{y^{1 - c_2}}{\rm{E}}_{70}\left(3-2c_1-c_2+a; 2-c_1, 2-c_2,c_3; x,y,z\right),
$

$
{u_6} = {y^{1 - c_2}}{z^{1 - c_3}}{\rm{E}}_{70}\left(2-c_2-c_3+a; c_1, 2-c_2,2-c_3; x,y,z\right),
$

$
{u_7} = {x^{1 - c_1}}{z^{1 - c_3}}{\rm{E}}_{70}\left(3-2c_1-c_3+a; 2-c_1, c_2,2-c_3; x,y,z\right),
$

$
{u_8} = {x^{1 - c_1}}{y^{1 - c_2}}{z^{1 - c_3}}{\rm{E}}_{70}\left(4-2c_1-c_2-c_3+a; 2-c_1, 2-c_2,2-c_3; x,y,z\right).
$

\bigskip

\begin{equation}
{\rm{E}}_{71}\left(a; c_1, c_2, c_3; x,y,z\right)=\sum\limits_{m,n,p=0}^\infty\frac{(a)_{2m+2n+p}}{(c_1)_{m}(c_2)_{n}(c_3)_{p}}\frac{x^m}{m!}\frac{y^n}{n!}\frac{z^p}{p!},
\end{equation}

region of convergence:
$$ \left\{ \sqrt{r}+\sqrt{s}<\frac{1}{2},\,\,\,t<\infty
\right\}.
$$

System of partial differential equations:

$
\left\{
\begin{aligned}
&x(1-4x)u_{xx}-8xyu_{xy}-4xzu_{xz}-4yzu_{yz}\\ &\,\,\,\,-4y^2u_{yy}-z^2u_{zz}+ \left[c_1-\left(4a+6\right)x\right]u_x-2\left(2a+3\right)yu_y-2\left(a+1\right)zu_z-a\left(a+1\right)u=0,\\&
y(1-4y)u_{yy}-4x^2u_{xx}-z^2u_{zz}-8xyu_{xy}\\ &\,\,\,\,-4xzu_{xz}-4yzu_{yz}-2\left(2a+3\right)xu_x+ \left[c_2-\left(4a+6\right)y\right]u_y-2\left(a+1\right)zu_z-a\left(a+1\right)u=0,\\&
zu_{zz}-2xu_x-2yu_y+\left(c_3-z\right)u_z-au=0,
\end{aligned}
\right.
$

where $u\equiv \,\,{\rm{E}}_{71}\left(a; c_1, c_2, c_3; x,y,z\right)$.

Particular solutions:

$
{u_1} ={\rm{E}}_{71}\left(a; c_1, c_2, c_3; x,y,z\right),
$

$
{u_2} = {x^{1 - c_1}}{\rm{E}}_{71}\left(2-2c_1+a; 2-c_1, c_2, c_3; x,y,z\right),
$

$
{u_3} = {y^{1 - c_2}}{\rm{E}}_{71}\left(2-2c_2+a; c_1, 2-c_2, c_3; x,y,z\right),
$

$
{u_4} = {z^{1 - c_3}}{\rm{E}}_{71}\left(1-c_3+a; c_1, c_2, 2-c_3; x,y,z\right),
$

$
{u_5} = {x^{1 - c_1}}{y^{1 - c_2}}{\rm{E}}_{71}\left(4-2c_1-2c_2+a; 2-c_1, 2-c_2, c_3; x,y,z\right),
$

$
{u_6} = {y^{1 - c_2}}{z^{1 - c_3}}{\rm{E}}_{71}\left(3-2c_2-c_3+a; c_1, 2-c_2, 2-c_3; x,y,z\right),
$

$
{u_7} = {x^{1 - c_1}}{z^{1 - c_3}}{\rm{E}}_{71}\left(3-2c_1-c_3+a; 2-c_1, c_2, 2-c_3; x,y,z\right),
$

$
{u_8} = {x^{1 - c_1}}{y^{1 - c_2}}{z^{1 - c_3}}{\rm{E}}_{71}\left(5-2c_1-2c_2-c_3+a; 2-c_1, 2-c_2, 2-c_3; x,y,z\right).
$

\bigskip

\begin{equation}
{\rm{E}}_{72}\left(a_1,a_2,a_3,a_4,b,c;x,y,z\right)=\sum\limits_{m,n,p=0}^\infty\frac{(a_1)_m(a_2)_n(a_3)_n(a_4)_p(b)_{m-p}}{(c)_{m+n}}\frac{x^m}{m!}\frac{y^n}{n!}\frac{z^p}{p!},
\end{equation}

region of convergence:
$$ \left\{ \frac{1}{r}+\frac{1}{s}>1,\,\,\,t<\infty
\right\}.
$$

System of partial differential equations:

$
\left\{
\begin{aligned}
&x(1-x)u_{xx}+yu_{xy}+xzu_{xz}+ \left[c-\left(a_1+b+1\right)x\right]u_x+a_1zu_z-a_1bu=0,\\&
y(1-y)u_{yy}+xu_{xy}+ \left[c-\left(a_2+a_3+1\right)y\right]u_y-a_2a_3u=0,\\&
zu_{zz}-xu_{xz}+\left(1-b+z\right)u_z+a_4u=0,
\end{aligned}
\right.
$

where $u\equiv \,\,{\rm{E}}_{72}\left(a_1,a_2,a_3,a_4,b,c;x,y,z\right)$.

\bigskip

\begin{equation}
{\rm{E}}_{73}\left(a_1,a_2,a_3,a_4,b,c;x,y,z\right)=\sum\limits_{m,n,p=0}^\infty\frac{(a_1)_m(a_2)_n(a_3)_p(a_4)_p(b)_{m-p}}{(c)_{m+n}}\frac{x^m}{m!}\frac{y^n}{n!}\frac{z^p}{p!},
\end{equation}

region of convergence:
$$ \left\{ (1+r)t<1,\,\,\,s<\infty
\right\}.
$$

System of partial differential equations:

$
\left\{
\begin{aligned}
&x(1-x)u_{xx}+yu_{xy}+xzu_{xz}+ \left[c-\left(a_1+b+1\right)x\right]u_x+a_1zu_z-a_1bu=0,\\&
yu_{yy}+xu_{xy}+ \left(c-y\right)u_y-a_2u=0,\\&
z(1+z)u_{zz}-xu_{xz}+\left[1-b+\left(a_3+a_4+1\right)z\right]u_z+a_3a_4u=0,
\end{aligned}
\right.
$

where $u\equiv \,\,{\rm{E}}_{73}\left(a_1,a_2,a_3,a_4,b,c;x,y,z\right)$.

\bigskip

\begin{equation}
{\rm{E}}_{74}\left(a_1,a_2,a_3,a_4,b,c;x,y,z\right)=\sum\limits_{m,n,p=0}^\infty\frac{(a_1)_n(a_2)_n(a_3)_p(a_4)_p(b)_{m-p}}{(c)_{m+n}}\frac{x^m}{m!}\frac{y^n}{n!}\frac{z^p}{p!},
\end{equation}

region of convergence:
$$ \left\{ r<\infty,\,\,\,t<1, \,\,\,\, s<1
\right\}.
$$

System of partial differential equations:

$
\left\{
\begin{aligned}
&xu_{xx}+yu_{xy}+ \left(c-x\right)u_x+zu_z-bu=0,\\&
y(1-y)u_{yy}+xu_{xy}+ \left[c-\left(a_1+a_2+1\right)y\right]u_y-a_1a_2u=0,\\&
z(1+z)u_{zz}-xu_{xz}+\left[1-b+\left(a_3+a_4+1\right)z\right]u_z+a_3a_4u=0,
\end{aligned}
\right.
$

where $u\equiv \,\,{\rm{E}}_{74}\left(a_1,a_2,a_3,a_4,b,c;x,y,z\right)$.

\bigskip

\begin{equation}
{\rm{E}}_{75}\left(a_1,a_2,a_3,b,c;x,y,z\right)=\sum\limits_{m,n,p=0}^\infty\frac{(a_1)_m(a_2)_n(a_3)_n(b)_{m-p}}{(c)_{m+n}}\frac{x^m}{m!}\frac{y^n}{n!}\frac{z^p}{p!},
\end{equation}

region of convergence:
$$ \left\{ \frac{1}{r}+\frac{1}{s}>1,\,\,\,t<\infty
\right\}.
$$

System of partial differential equations:

$
\left\{
\begin{aligned}
&x(1-x)u_{xx}+yu_{xy}+xzu_{xz}+ \left[c-\left(a_1+b+1\right)x\right]u_x+a_1zu_z-a_1bu=0,\\&
y(1-y)u_{yy}+xu_{xy}+ \left[c-\left(a_2+a_3+1\right)y\right]u_y-a_2a_3u=0,\\&
zu_{zz}-xu_{xz}+\left(1-b\right)u_z+u=0,
\end{aligned}
\right.
$

where $u\equiv \,\,{\rm{E}}_{75}\left(a_1,a_2,a_3,b,c;x,y,z\right)$.

\bigskip

\begin{equation}
{\rm{E}}_{76}\left(a_1,a_2,a_3,b,c;x,y,z\right)=\sum\limits_{m,n,p=0}^\infty\frac{(a_1)_m(a_2)_n(a_3)_p(b)_{m-p}}{(c)_{m+n}}\frac{x^m}{m!}\frac{y^n}{n!}\frac{z^p}{p!},
\end{equation}

region of convergence:
$$ \left\{
r<1,\,\,\,s<\infty,\,\,\,t<\infty\right\}.
$$

System of partial differential equations:

$
\left\{
\begin{aligned}
&x(1-x)u_{xx}+yu_{xy}+xzu_{xz}+ \left[c-\left(a_1+b+1\right)x\right]u_x+a_1zu_z-a_1bu=0,\\&
yu_{yy}+xu_{xy}+ \left(c-y\right)u_y-a_2u=0,\\&
zu_{zz}-xu_{xz}+\left(1-b+z\right)u_z+a_3u=0,
\end{aligned}
\right.
$

where $u\equiv \,\,{\rm{E}}_{76}\left(a_1,a_2,a_3,b,c;x,y,z\right)$.

\bigskip

\begin{equation}
{\rm{E}}_{77}\left(a_1,a_2,a_3,b,c;x,y,z\right)=\sum\limits_{m,n,p=0}^\infty\frac{(a_1)_n(a_2)_n(a_3)_p(b)_{m-p}}{(c)_{m+n}}\frac{x^m}{m!}\frac{y^n}{n!}\frac{z^p}{p!},
\end{equation}

region of convergence:
$$ \left\{ r<\infty,\,\,\,s<1,\,\,\,t<\infty
\right\}.
$$

System of partial differential equations:

$
\left\{
\begin{aligned}
&xu_{xx}+yu_{xy}+ \left(c-x\right)u_x+zu_z-bu=0,\\&
y(1-y)u_{yy}+xu_{xy}+ \left[c-\left(a_1+a_2+1\right)y\right]u_y-a_1a_2u=0,\\&
zu_{zz}-xu_{xz}+\left(1-b+z\right)u_z+a_3u=0,
\end{aligned}
\right.
$

where $u\equiv \,\,{\rm{E}}_{77}\left(a_1,a_2,a_3,b,c;x,y,z\right)$.

\bigskip

\begin{equation}
{\rm{E}}_{78}\left(a_1,a_2,a_3,b,c;x,y,z\right)=\sum\limits_{m,n,p=0}^\infty\frac{(a_1)_n(a_2)_p(a_3)_p(b)_{m-p}}{(c)_{m+n}}\frac{x^m}{m!}\frac{y^n}{n!}\frac{z^p}{p!},
\end{equation}

region of convergence:
$$ \left\{ r<\infty,\,\,\,s<\infty,\,\,\,t<1
\right\}.
$$

System of partial differential equations:

$
\left\{
\begin{aligned}
&xu_{xx}+yu_{xy}+ \left(c-x\right)u_x+zu_z-bu=0,\\&
yu_{yy}+xu_{xy}+ \left(c-y\right)u_y-a_1u=0,\\&
z(1+z)u_{zz}-xu_{xz}+\left[1-b+\left(a_2+a_3+1\right)z\right]u_z+a_2a_3u=0,
\end{aligned}
\right.
$

where $u\equiv \,\,{\rm{E}}_{78}\left(a_1,a_2,a_3,b,c;x,y,z\right)$.

\bigskip

\begin{equation}
{\rm{E}}_{79}\left(a_1,a_2,a_3,b,c;x,y,z\right)=\sum\limits_{m,n,p=0}^\infty\frac{(a_1)_m(a_2)_p(a_3)_p(b)_{m-p}}{(c)_{m+n}}\frac{x^m}{m!}\frac{y^n}{n!}\frac{z^p}{p!},
\end{equation}

region of convergence:
$$ \left\{ (1+r)t<1,\,\,\,s<\infty
\right\}.
$$

System of partial differential equations:

$
\left\{
\begin{aligned}
&x(1-x)u_{xx}+yu_{xy}+xzu_{xz}+ \left[c-\left(a_1+b+1\right)x\right]u_x+a_1zu_z-a_1bu=0,\\&
yu_{yy}+xu_{xy}+ cu_y-u=0,\\&
z(1+z)u_{zz}-xu_{xz}+\left[1-b+\left(a_2+a_3+1\right)z\right]u_z+a_2a_3u=0,
\end{aligned}
\right.
$

where $u\equiv \,\,{\rm{E}}_{79}\left(a_1,a_2,a_3,b,c;x,y,z\right)$.

\bigskip

\begin{equation}
{\rm{E}}_{80}\left(a_1,a_2,b,c;x,y,z\right)=\sum\limits_{m,n,p=0}^\infty\frac{(a_1)_m(a_2)_n(b)_{m-p}}{(c)_{m+n}}\frac{x^m}{m!}\frac{y^n}{n!}\frac{z^p}{p!},
\end{equation}

region of convergence:
$$ \left\{ r<1,\,\,\,s<\infty,\,\,\,t<\infty
\right\}.
$$

System of partial differential equations:

$
\left\{
\begin{aligned}
&x(1-x)u_{xx}+yu_{xy}+xzu_{xz}+ \left[c-\left(a_1+b+1\right)x\right]u_x+a_1zu_z-a_1bu=0,\\&
yu_{yy}+xu_{xy}+ \left(c-y\right)u_y-a_2u=0,\\&
zu_{zz}-xu_{xz}+\left(1-b\right)u_z+u=0,
\end{aligned}
\right.
$

where $u\equiv \,\,{\rm{E}}_{80}\left(a_1,a_2,b,c;x,y,z\right)$.

\bigskip

\begin{equation}
{\rm{E}}_{81}\left(a_1,a_2,b,c;x,y,z\right)=\sum\limits_{m,n,p=0}^\infty\frac{(a_1)_m(a_2)_p(b)_{m-p}}{(c)_{m+n}}\frac{x^m}{m!}\frac{y^n}{n!}\frac{z^p}{p!},
\end{equation}

region of convergence:
$$ \left\{ r<1,\,\,\,s<\infty,\,\,\,t<\infty
\right\}.
$$

System of partial differential equations:

$
\left\{
\begin{aligned}
&x(1-x)u_{xx}+yu_{xy}+xzu_{xz}+ \left[c-\left(a_1+b+1\right)x\right]u_x+a_1zu_z-a_1bu=0,\\&
yu_{yy}+xu_{xy}+ cu_y-u=0,\\&
zu_{zz}-xu_{xz}+\left(1-b+z\right)u_z+a_2u=0,
\end{aligned}
\right.
$

where $u\equiv \,\,{\rm{E}}_{81}\left(a_1,a_2,b,c;x,y,z\right)$.

\bigskip

\begin{equation}
{\rm{E}}_{82}\left(a_1,a_2,b,c;x,y,z\right)=\sum\limits_{m,n,p=0}^\infty\frac{(a_1)_n(a_2)_n(b)_{m-p}}{(c)_{m+n}}\frac{x^m}{m!}\frac{y^n}{n!}\frac{z^p}{p!},
\end{equation}

region of convergence:
$$ \left\{ r<\infty,\,\,\,s<1,\,\,\,t<\infty
\right\}.
$$

System of partial differential equations:

$
\left\{
\begin{aligned}
&xu_{xx}+yu_{xy}+ \left(c-x\right)u_x+zu_z-bu=0,\\&
y(1-y)u_{yy}+xu_{xy}+ \left[c-\left(a_1+a_2+1\right)y\right]u_y-a_1a_2u=0,\\&
zu_{zz}-xu_{xz}+\left(1-b\right)u_z+u=0,
\end{aligned}
\right.
$

where $u\equiv \,\,{\rm{E}}_{82}\left(a_1,a_2,b,c;x,y,z\right)$.

\bigskip

\begin{equation}
{\rm{E}}_{83}\left(a_1,a_2,b,c;x,y,z\right)=\sum\limits_{m,n,p=0}^\infty\frac{(a_1)_n(a_2)_p(b)_{m-p}}{(c)_{m+n}}\frac{x^m}{m!}\frac{y^n}{n!}\frac{z^p}{p!},
\end{equation}

region of convergence:
$$ \left\{ r<\infty,\,\,\,s<\infty,\,\,\,t<\infty
\right\}.
$$

System of partial differential equations:

$
\left\{
\begin{aligned}
&xu_{xx}+yu_{xy}+ \left(c-x\right)u_x+zu_z-bu=0,\\&
yu_{yy}+xu_{xy}+ \left(c-y\right)u_y-a_1u=0,\\&
zu_{zz}-xu_{xz}+\left(1-b+z\right)u_z+a_2u=0,
\end{aligned}
\right.
$

where $u\equiv \,\,{\rm{E}}_{83}\left(a_1,a_2,b,c;x,y,z\right)$.

\bigskip

\begin{equation}
{\rm{E}}_{84}\left(a_1,a_2,b,c;x,y,z\right)=\sum\limits_{m,n,p=0}^\infty\frac{(a_1)_p(a_2)_p(b)_{m-p}}{(c)_{m+n}}\frac{x^m}{m!}\frac{y^n}{n!}\frac{z^p}{p!},
\end{equation}

region of convergence:
$$ \left\{(1+r)t<1,\,\,\,s<\infty
\right\}.
$$

System of partial differential equations:

$
\left\{
\begin{aligned}
&xu_{xx}+yu_{xy}+ \left(c-x\right)u_x+zu_z-bu=0,\\&
yu_{yy}+xu_{xy}+ cu_y-u=0,\\&
z(1+z)u_{zz}-xu_{xz}+\left[1-b+\left(a_1+a_2+1\right)z\right]u_z+a_1a_2u=0,
\end{aligned}
\right.
$

where $u\equiv \,\,{\rm{E}}_{84}\left(a_1,a_2,b,c;x,y,z\right)$.

\bigskip

\begin{equation}
{\rm{E}}_{85}\left(a,b,c;x,y,z\right)=\sum\limits_{m,n,p=0}^\infty\frac{(a)_m(b)_{m-p}}{(c)_{m+n}}\frac{x^m}{m!}\frac{y^n}{n!}\frac{z^p}{p!},
\end{equation}

region of convergence:
$$ \left\{ r<1,\,\,\,s<\infty,\,\,\,t<\infty
\right\}.
$$

System of partial differential equations:

$
\left\{
\begin{aligned}
&x(1-x)u_{xx}+yu_{xy}+xzu_{xz}+ \left[c-\left(a+b+1\right)x\right]u_x+azu_z-abu=0,\\&
yu_{yy}+xu_{xy}+ cu_y-u=0,\\&
zu_{zz}-xu_{xz}+\left(1-b\right)u_z+u=0,
\end{aligned}
\right.
$

where $u\equiv \,\,{\rm{E}}_{85}\left(a,b,c;x,y,z\right)$.

\bigskip

\begin{equation}
{\rm{E}}_{86}\left(a,b,c;x,y,z\right)=\sum\limits_{m,n,p=0}^\infty\frac{(a)_n(b)_{m-p}}{(c)_{m+n}}\frac{x^m}{m!}\frac{y^n}{n!}\frac{z^p}{p!},
\end{equation}

region of convergence:
$$ \left\{ r<\infty,\,\,\,s<\infty,\,\,\,t<\infty
\right\}.
$$

System of partial differential equations:

$
\left\{
\begin{aligned}
&xu_{xx}+yu_{xy}+ \left(c-x\right)u_x+zu_z-bu=0,\\&
yu_{yy}+xu_{xy}+ \left(c-y\right)u_y-au=0,\\&
zu_{zz}-xu_{xz}+\left(1-b\right)u_z+u=0,
\end{aligned}
\right.
$

where $u\equiv \,\,{\rm{E}}_{86}\left(a,b,c;x,y,z\right)$.

\bigskip

\begin{equation}
{\rm{E}}_{87}\left(a,b,c;x,y,z\right)=\sum\limits_{m,n,p=0}^\infty\frac{(a)_p(b)_{m-p}}{(c)_{m+n}}\frac{x^m}{m!}\frac{y^n}{n!}\frac{z^p}{p!},
\end{equation}

region of convergence:
$$ \left\{ r<\infty,\,\,\,s<\infty,\,\,\,t<\infty
\right\}.
$$

System of partial differential equations:

$
\left\{
\begin{aligned}
&xu_{xx}+yu_{xy}+ \left(c-x\right)u_x+zu_z-bu=0,\\&
yu_{yy}+xu_{xy}+ cu_y-u=0,\\&
zu_{zz}-xu_{xz}+\left(1-b+z\right)u_z+au=0,
\end{aligned}
\right.
$

where $u\equiv \,\,{\rm{E}}_{87}\left(a,b,c;x,y,z\right)$.

\bigskip

\begin{equation}
{\rm{E}}_{88}\left(b,c;x,y,z\right)=\sum\limits_{m,n,p=0}^\infty\frac{(b)_{m-p}}{(c)_{m+n}}\frac{x^m}{m!}\frac{y^n}{n!}\frac{z^p}{p!},
\end{equation}

region of convergence:
$$ \left\{ r<\infty,\,\,\,s<\infty,\,\,\,t<\infty
\right\}.
$$

System of partial differential equations:

$
\left\{
\begin{aligned}
&xu_{xx}+yu_{xy}+ \left(c-x\right)u_x+zu_z-bu=0,\\&
yu_{yy}+xu_{xy}+ cu_y-u=0,\\&
zu_{zz}-xu_{xz}+\left(1-b\right)u_z+u=0,
\end{aligned}
\right.
$

where $u\equiv \,\,{\rm{E}}_{88}\left(b,c;x,y,z\right)$.

\bigskip

\begin{equation}
{\rm{E}}_{89}\left(a_1,a_2,a_3,b; c;x,y,z\right)=
\sum\limits_{m,n,p=0}^\infty\frac{(a_1)_m(a_2)_n(a_3)_p(b)_{m+n-p}}{(c)_{m+n}}\frac{x^m}{m!}\frac{y^n}{n!}\frac{z^p}{p!},
\end{equation}

region of convergence:
$$ \left\{ r<1, \,\,\,\,\, s<1,\,\,\,t<\infty
\right\}.
$$

System of partial differential equations:

$
\left\{
\begin{aligned}
&x(1-x)u_{xx}+(1-x)yu_{xy}+xzu_{xz}+ \left[c-\left(a_1+b+1\right)x\right]u_x-a_1yu_y+a_1zu_z-a_1bu=0,\\&
y(1-y)u_{yy}+x(1-y)u_{xy}+yzu_{yz}-a_2xu_x+\left[c-\left(a_2+b+1\right)y\right]u_y+a_2zu_z-a_2bu=0,\\&
zu_{zz}-xu_{xz}-yu_{yz}+\left(1-b+z\right)u_z+a_3u=0,
\end{aligned}
\right.
$

where $u\equiv \,\,{\rm{E}}_{89}\left(a_1,a_2,a_3,b; c;x,y,z\right)$.

\bigskip

\begin{equation}
{\rm{E}}_{90}\left(a_1,a_2,a_3,b; c;x,y,z\right)=
\sum\limits_{m,n,p=0}^\infty\frac{(a_1)_m(a_2)_p(a_3)_p(b)_{m+n-p}}{(c)_{m+n}}\frac{x^m}{m!}\frac{y^n}{n!}\frac{z^p}{p!},
\end{equation}

region of convergence:
$$ \left\{ (1+r)t<1,\,\,\,s<\infty
\right\}.
$$

System of partial differential equations:

$
\left\{
\begin{aligned}
&x(1-x)u_{xx}+(1-x)yu_{xy}+xzu_{xz}+ \left[c-\left(a_1+b+1\right)x\right]u_x-a_1yu_y+a_1zu_z-a_1bu=0,\\&
yu_{yy}+xu_{xy}-xu_x+\left(c-y\right)u_y+zu_z-bu=0,\\&
z(1+z)u_{zz}-xu_{xz}-yu_{yz}+\left[1-b+\left(a_2+a_3+1\right)z\right]u_z-a_2a_3u=0,
\end{aligned}
\right.
$

where $u\equiv \,\,{\rm{E}}_{90}\left(a_1,a_2,a_3,b; c;x,y,z\right)$.

\bigskip

\begin{equation}
{\rm{E}}_{91}\left(a_1,a_2,b; c;x,y,z\right)=
\sum\limits_{m,n,p=0}^\infty\frac{(a_1)_m(a_2)_n(b)_{m+n-p}}{(c)_{m+n}}\frac{x^m}{m!}\frac{y^n}{n!}\frac{z^p}{p!},
\end{equation}

region of convergence:
$$ \left\{ r<1,\,\,\,\, s<1,\,\,\,t<\infty
\right\}.
$$

System of partial differential equations:

$
\left\{
\begin{aligned}
&x(1-x)u_{xx}+(1-x)yu_{xy}+xzu_{xz}+ \left[c-\left(a_1+b+1\right)x\right]u_x-a_1yu_y+a_1zu_z-a_1bu=0,\\&
y(1-y)u_{yy}+x(1-y)u_{xy}+yzu_{yz}-a_2xu_x+\left[c-\left(a_2+b+1\right)y\right]u_y+a_2zu_z-a_2bu=0,\\&
zu_{zz}-xu_{xz}-yu_{yz}+(1-b)u_z+u=0,
\end{aligned}
\right.
$

where $u\equiv \,\,{\rm{E}}_{91}\left(a_1,a_2,b; c;x,y,z\right)$.

\bigskip

\begin{equation}
{\rm{E}}_{92}\left(a_1,a_2,b; c;x,y,z\right)=
\sum\limits_{m,n,p=0}^\infty\frac{(a_1)_m(a_2)_p(b)_{m+n-p}}{(c)_{m+n}}\frac{x^m}{m!}\frac{y^n}{n!}\frac{z^p}{p!},
\end{equation}

region of convergence:
$$ \left\{ r<1,\,\,\,s<\infty,\,\,\,t<\infty
\right\}.
$$

System of partial differential equations:

$
\left\{
\begin{aligned}
&x(1-x)u_{xx}+(1-x)yu_{xy}+xzu_{xz}+ \left[c-\left(a_1+b+1\right)x\right]u_x-a_1yu_y+a_1zu_z-a_1bu=0,\\&
yu_{yy}+xu_{xy}-xu_x+\left(c-y\right)u_y+zu_z-bu=0,\\&
zu_{zz}-xu_{xz}-yu_{yz}+\left(1-b+z\right)u_z+a_2u=0,
\end{aligned}
\right.
$

where $u\equiv \,\,{\rm{E}}_{92}\left(a_1,a_2,b; c;x,y,z\right)$.

\bigskip

\begin{equation}
{\rm{E}}_{93}\left(a,b; c;x,y,z\right)=
\sum\limits_{m,n,p=0}^\infty\frac{(a)_m(b)_{m+n-p}}{(c)_{m+n}}\frac{x^m}{m!}\frac{y^n}{n!}\frac{z^p}{p!},
\end{equation}

region of convergence:
$$ \left\{ r<1,\,\,\,s<\infty,\,\,\,t<\infty
\right\}.
$$

System of partial differential equations:

$
\left\{
\begin{aligned}
&x(1-x)u_{xx}+(1-x)yu_{xy}+xzu_{xz}+ \left[c-\left(a+b+1\right)x\right]u_x-ayu_y+azu_z-abu=0,\\&
yu_{yy}+xu_{xy}-xu_x+\left(c-y\right)u_y+zu_z-bu=0,\\&
zu_{zz}-xu_{xz}-yu_{yz}+(1-b)u_z+u=0,
\end{aligned}
\right.
$

where $u\equiv \,\,{\rm{E}}_{93}\left(a,b; c;x,y,z\right)$.

\bigskip

\begin{equation}
{\rm{E}}_{94}\left(a_1,a_2,a_3,b; c;x,y,z\right)=\sum\limits_{m,n,p=0}^\infty\frac{(a_1)_{m+n}(a_2)_m(a_3)_p(b)_{n-p}}{(c)_{m+n}}\frac{x^m}{m!}\frac{y^n}{n!}\frac{z^p}{p!},
\end{equation}

region of convergence:
$$ \left\{ r<1, \,\,\,\,\, s<1,\,\,\,t<\infty
\right\}.
$$

System of partial differential equations:

$
\left\{
\begin{aligned}
&x(1-x)u_{xx}+(1-x)yu_{xy}+ \left[c-\left(a_1+a_2+1\right)x\right]u_x-a_2yu_y-a_1a_2u=0,\\&
y(1-y)u_{yy}+x(1-y)u_{xy}\\& \,\,\,\,\,\,\,\,\,+xzu_{xz}+yzu_{yz}-bxu_x+\left[c-\left(a_1+b+1\right)y\right]u_y+a_1zu_z-a_1bu=0,\\&
zu_{zz}-yu_{yz}+\left(1-b+z\right)u_z+a_3u=0,
\end{aligned}
\right.
$

where $u\equiv \,\,{\rm{E}}_{94}\left(a_1,a_2,a_3,b; c;x,y,z\right)$.

\bigskip

\begin{equation}
{\rm{E}}_{95}\left(a_1,a_2,a_3,b; c;x,y,z\right)=\sum\limits_{m,n,p=0}^\infty\frac{(a_1)_{m+n}(a_2)_p(a_3)_p(b)_{n-p}}{(c)_{m+n}}\frac{x^m}{m!}\frac{y^n}{n!}\frac{z^p}{p!},
\end{equation}

region of convergence:
$$ \left\{ r<\infty,\,\,\,(1+s)t<1
\right\}.
$$

System of partial differential equations:

$
\left\{
\begin{aligned}
&xu_{xx}+yu_{xy}+ (c-x)u_x-yu_y-a_1u=0,\\&
y(1-y)u_{yy}+x(1-y)u_{xy}\\& \,\,\,\,\,\,\,\,\,+xzu_{xz}+yzu_{yz}-bxu_x+\left[c-\left(a_1+b+1\right)y\right]u_y+a_1zu_z-a_1bu=0,\\&
z(1+z)u_{zz}-yu_{yz}+\left[1-b+\left(a_2+a_3+1\right)z\right]u_z+a_2a_3u=0,
\end{aligned}
\right.
$

where $u\equiv \,\,{\rm{E}}_{95}\left(a_1,a_2,a_3,b; c;x,y,z\right)$.

\bigskip

\begin{equation}
{\rm{E}}_{96}\left(a_1,a_2,b; c;x,y,z\right)=\sum\limits_{m,n,p=0}^\infty\frac{(a_1)_{m+n}(a_2)_m(b)_{n-p}}{(c)_{m+n}}\frac{x^m}{m!}\frac{y^n}{n!}\frac{z^p}{p!},
\end{equation}

region of convergence:
$$ \left\{ r<1, \,\,\,\, s<1,\,\,\,t<\infty
\right\}.
$$

System of partial differential equations:

$
\left\{
\begin{aligned}
&x(1-x)u_{xx}+(1-x)yu_{xy}+ \left[c-\left(a_1+a_2+1\right)x\right]u_x-a_2yu_y-a_1a_2u=0,\\&
y(1-y)u_{yy}+x(1-y)u_{xy}+xzu_{xz}+yzu_{yz}\\& \,\,\,\,\,\,\,\,\,-bxu_x+\left[c-\left(a_1+b+1\right)y\right]u_y+a_1zu_z-a_1bu=0,\\&
zu_{zz}-yu_{yz}+\left(1-b\right)u_z+u=0,
\end{aligned}
\right.
$

where $u\equiv \,\,{\rm{E}}_{96}\left(a_1,a_2,b; c;x,y,z\right)$.

\bigskip

\begin{equation}
{\rm{E}}_{97}\left(a_1,a_2,b; c;x,y,z\right)=\sum\limits_{m,n,p=0}^\infty\frac{(a_1)_{m+n}(a_2)_p(b)_{n-p}}{(c)_{m+n}}\frac{x^m}{m!}\frac{y^n}{n!}\frac{z^p}{p!},
\end{equation}

region of convergence:
$$ \left\{ r<\infty,\,\,\, s<1,\,\,\,t<\infty
\right\}.
$$

System of partial differential equations:

$
\left\{
\begin{aligned}
&xu_{xx}+yu_{xy}+ (c-x)u_x-yu_y-a_1u=0,\\&
y(1-y)u_{yy}+x(1-y)u_{xy}+xzu_{xz}+yzu_{yz}\\& \,\,\,\,\,\,\,\,\,-bxu_x+\left[c-\left(a_1+b+1\right)y\right]u_y+a_1zu_z-a_1bu=0,\\&
zu_{zz}-yu_{yz}+\left(1-b+z\right)u_z+a_2u=0,
\end{aligned}
\right.
$

where $u\equiv \,\,{\rm{E}}_{97}\left(a_1,a_2,b; c;x,y,z\right)$.

\bigskip

\begin{equation}
{\rm{E}}_{98}\left(a,b; c;x,y,z\right)=\sum\limits_{m,n,p=0}^\infty\frac{(a)_{m+n}(b)_{n-p}}{(c)_{m+n}}\frac{x^m}{m!}\frac{y^n}{n!}\frac{z^p}{p!},
\end{equation}

region of convergence:
$$ \left\{ r<\infty,\,\,\, s<1,\,\,\,t<\infty
\right\}.
$$

System of partial differential equations:

$
\left\{
\begin{aligned}
&xu_{xx}+yu_{xy}+ (c-x)u_x-yu_y-au=0,\\&
y(1-y)u_{yy}+x(1-y)u_{xy}+xzu_{xz}+yzu_{yz}-bxu_x+\left[c-\left(a+b+1\right)y\right]u_y+azu_z-abu=0,\\&
zu_{zz}-yu_{yz}+\left(1-b\right)u_z+u=0,
\end{aligned}
\right.
$

where $u\equiv \,\,{\rm{E}}_{98}\left(a,b; c;x,y,z\right)$.

\bigskip

\begin{equation}
{\rm{E}}_{99}\left(a_1,a_2,a_3,b;c; x,y,z\right)=\sum\limits_{m,n,p=0}^\infty\frac{(a_1)_{n+p}(a_2)_m(a_3)_m(b)_{n-p}}{(c)_{m+n}}\frac{x^m}{m!}\frac{y^n}{n!}\frac{z^p}{p!},
\end{equation}

region of convergence:
$$ \left\{ r<\infty,\,\,\,\frac{1}{r}+\frac{1}{s}>1
\right\}.
$$

System of partial differential equations:

$
\left\{
\begin{aligned}
&x(1-x)u_{xx}+yu_{xy}+\left[c_1-\left(a_2+a_3+1\right)x\right]u_x-a_2a_3u=0,\\&
y(1-y)u_{yy}+z^2u_{zz}+xu_{xy}+\left[c-\left(a_1+b+1\right)y\right]u_y+\left(a_1-b+1\right)zu_z-a_1bu=0,\\&
zu_{zz}-yu_{yz}+yu_y+(1-b+z)u_z+a_1u=0,
\end{aligned}
\right.
$

where $u\equiv \,\,{\rm{E}}_{99}\left(a_1,a_2,a_3,b;c; x,y,z\right)$.

\bigskip

\begin{equation}
{\rm{E}}_{100}\left(a_1,a_2,a_3,b;c; x,y,z\right)=\sum\limits_{m,n,p=0}^\infty\frac{(a_1)_{n+p}(a_2)_m(a_3)_p(b)_{n-p}}{(c)_{m+n}}\frac{x^m}{m!}\frac{y^n}{n!}\frac{z^p}{p!},
\end{equation}

region of convergence:
$$ \left\{r<\infty,\,\,\,\left\{t+2\sqrt{st}<1\right\}=\left\{\sqrt{t}<\sqrt{1+s}-\sqrt{s}\right\}
\right\}.
$$

System of partial differential equations:

$
\left\{
\begin{aligned}
&xu_{xx}+yu_{xy}+\left(c-x\right)u_x-a_2u=0,\\&
y(1-y)u_{yy}+z^2u_{zz}+xu_{xy}+\left[c-\left(a_1+b+1\right)y\right]u_y+\left(a_1-b+1\right)zu_z-a_1bu=0,\\&
z(1+z)u_{zz}-y(1-z)u_{yz}+a_3yu_y+\left[1-b+\left(a_1+a_3+1\right)z\right]u_z+a_1a_3u=0,
\end{aligned}
\right.
$

where $u\equiv \,\,{\rm{E}}_{100}\left(a_1,a_2,a_3,b;c; x,y,z\right)$.

\bigskip

\begin{equation}
{\rm{E}}_{101}\left(a_1,a_2,b;c; x,y,z\right)=\sum\limits_{m,n,p=0}^\infty\frac{(a_1)_{n+p}(a_2)_m(b)_{n-p}}{(c)_{m+n}}\frac{x^m}{m!}\frac{y^n}{n!}\frac{z^p}{p!},
\end{equation}

region of convergence:
$$ \left\{ r<\infty,\,\,\,s<1,\,\,\,t<\infty
\right\}.
$$

System of partial differential equations:

$
\left\{
\begin{aligned}
&xu_{xx}+yu_{xy}+\left(c-x\right)u_x-a_2u=0,\\&
y(1-y)u_{yy}+z^2u_{zz}+xu_{xy}+\left[c-\left(a_1+b+1\right)y\right]u_y+\left(a_1-b+1\right)zu_z-a_1bu=0,\\&
zu_{zz}-yu_{yz}+yu_y+(1-b+z)u_z+a_1u=0,
\end{aligned}
\right.
$

where $u\equiv \,\,{\rm{E}}_{101}\left(a_1,a_2,b;c; x,y,z\right)$.

\bigskip

\begin{equation}
{\rm{E}}_{102}\left(a_1,a_2,b;c; x,y,z\right)=\sum\limits_{m,n,p=0}^\infty\frac{(a_1)_{n+p}(a_2)_p(b)_{n-p}}{(c)_{m+n}}\frac{x^m}{m!}\frac{y^n}{n!}\frac{z^p}{p!},
\end{equation}

region of convergence:
$$ \left\{ r<\infty,\,\,\,s<\infty,\,\,\,t<\infty
\right\}.
$$

System of partial differential equations:

$
\left\{
\begin{aligned}
&xu_{xx}+yu_{xy}+cu_x-u=0,\\&
y(1-y)u_{yy}+z^2u_{zz}+xu_{xy}+\left[c-\left(a_1+b+1\right)y\right]u_y+\left(a_1-b+1\right)zu_z-a_1bu=0,\\&
z(1+z)u_{zz}-y(1-z)u_{yz}+a_2yu_y+\left[1-b+\left(a_1+a_2+1\right)z\right]u_z+a_1a_2u=0,
\end{aligned}
\right.
$

where $u\equiv \,\,{\rm{E}}_{102}\left(a_1,a_2,b;c; x,y,z\right)$.

\bigskip

\begin{equation}
{\rm{E}}_{103}\left(a,b;c; x,y,z\right)=\sum\limits_{m,n,p=0}^\infty\frac{(a)_{n+p}(b)_{n-p}}{(c)_{m+n}}\frac{x^m}{m!}\frac{y^n}{n!}\frac{z^p}{p!},
\end{equation}

region of convergence:
$$ \left\{ r<\infty,\,\,\,s<1,\,\,\,t<\infty
\right\}.
$$

System of partial differential equations:

$
\left\{
\begin{aligned}
&xu_{xx}+yu_{xy}+cu_x-u=0,\\&
y(1-y)u_{yy}+z^2u_{zz}+xu_{xy}+\left[c-\left(a+b+1\right)y\right]u_y+\left(a-b+1\right)zu_z-abu=0,\\&
zu_{zz}-yu_{yz}+yu_y+(1-b+z)u_z+au=0,
\end{aligned}
\right.
$

where $u\equiv \,\,{\rm{E}}_{103}\left(a,b;c; x,y,z\right)$.

\bigskip

\begin{equation}
{\rm{E}}_{104}\left(a_1,a_2,a_3,b;c;x,y,z\right)=\sum\limits_{m,n,p=0}^\infty\frac{(a_1)_{n+p}(a_2)_m(a_3)_n(b)_{m-p}}{(c)_{m+n}}\frac{x^m}{m!}\frac{y^n}{n!}\frac{z^p}{p!},
\end{equation}

region of convergence:
$$ \left\{ \frac{1}{r}+\frac{1}{s}>1,\,\,\,t<\infty
\right\}.
$$

System of partial differential equations:

$
\left\{
\begin{aligned}
&x(1-x)u_{xx}+yu_{xy}+xzu_{xz}+\left[c-\left(a_2+b+1\right)x\right]u_x+a_2zu_z-a_2bu=0,\\&
y(1-y)u_{yy}+xu_{xy}-yzu_{yz}+\left[c-\left(a_1+a_3+1\right)y\right]u_y-a_3zu_z-a_1a_3u=0,\\&
zu_{zz}-xu_{xz}+yu_y+(1-b+z)u_z+a_1u=0,
\end{aligned}
\right.
$

where $u\equiv \,\,{\rm{E}}_{104}\left(a_1,a_2,a_3,b;c;x,y,z\right)$.

\bigskip

\begin{equation}
{\rm{E}}_{105}\left(a_1,a_2,a_3,b;c;x,y,z\right)=\sum\limits_{m,n,p=0}^\infty\frac{(a_1)_{n+p}(a_2)_m(a_3)_p(b)_{m-p}}{(c)_{m+n}}\frac{x^m}{m!}\frac{y^n}{n!}\frac{z^p}{p!},
\end{equation}

region of convergence:
$$ \left\{ (1+r)t<1,\,\,\,s<\infty
\right\}.
$$

System of partial differential equations:

$
\left\{
\begin{aligned}
&x(1-x)u_{xx}+yu_{xy}+xzu_{xz}+\left[c-\left(a_2+b+1\right)x\right]u_x+a_2zu_z-a_2bu=0,\\&
yu_{yy}+xu_{xy}+(c-y)u_y-zu_z-a_1u=0,\\&
z(1+z)u_{zz}-xu_{xz}+yzu_{yz}+a_3yu_y+\left[1-b+\left(a_1+a_3+1\right)z\right]u_z+a_1a_3u=0,
\end{aligned}
\right.
$

where $u\equiv \,\,{\rm{E}}_{105}\left(a_1,a_2,a_3,b;c;x,y,z\right)$.

\bigskip

\begin{equation}
{\rm{E}}_{106}\left(a_1,a_2,a_3,b;c; x,y,z\right)=\sum\limits_{m,n,p=0}^\infty\frac{(a_1)_{n+p}(a_2)_n(a_3)_p(b)_{m-p}}{(c)_{m+n}}\frac{x^m}{m!}\frac{y^n}{n!}\frac{z^p}{p!},
\end{equation}

region of convergence:
$$ \left\{ r<\infty,\,\,\,s+t<1
\right\}.
$$

System of partial differential equations:

$
\left\{
\begin{aligned}
&xu_{xx}+yu_{xy}+(c-x)u_x+zu_z-bu=0,\\&
y(1-y)u_{yy}+xu_{xy}-yzu_{yz}+\left[c-\left(a_1+a_2+1\right)y\right]u_y-a_2zu_z-a_1a_2u=0,\\&
z(1+z)u_{zz}-xu_{xz}+yzu_{yz}+a_3yu_y+\left[1-b+\left(a_1+a_3+1\right)z\right]u_z+a_1a_3u=0,
\end{aligned}
\right.
$

where $u\equiv \,\,{\rm{E}}_{106}\left(a_1,a_2,a_3,b;c; x,y,z\right)$.

\bigskip

\begin{equation}
{\rm{E}}_{107}\left(a_1,a_2,b;c; x,y,z\right)=\sum\limits_{m,n,p=0}^\infty\frac{(a_1)_{n+p}(a_2)_m(b)_{m-p}}{(c)_{m+n}}\frac{x^m}{m!}\frac{y^n}{n!}\frac{z^p}{p!},
\end{equation}

region of convergence:
$$ \left\{ r<1,\,\,\,s<\infty,\,\,\,t<\infty
\right\}.
$$

System of partial differential equations:

$
\left\{
\begin{aligned}
&x(1-x)u_{xx}+yu_{xy}+xzu_{xz}+\left[c-\left(a_2+b+1\right)x\right]u_x+a_2zu_z-a_2bu=0,\\&
yu_{yy}+xu_{xy}+(c-y)u_y-zu_z-a_1u=0,\\&
zu_{zz}-xu_{xz}+yu_y+(1-b+z)u_z+a_1u=0,
\end{aligned}
\right.
$

where $u\equiv \,\,{\rm{E}}_{107}\left(a_1,a_2,b;c; x,y,z\right)$.

\bigskip

\begin{equation}
{\rm{E}}_{108}\left(a_1,a_2, b;c; x,y,z\right)=\sum\limits_{m,n,p=0}^\infty\frac{(a_1)_{n+p}(a_2)_n(b)_{m-p}}{(c)_{m+n}}\frac{x^m}{m!}\frac{y^n}{n!}\frac{z^p}{p!},
\end{equation}

region of convergence:
$$ \left\{ r<\infty,\,\,\,s<1,\,\,\,t<\infty
\right\}.
$$

System of partial differential equations:

$
\left\{
\begin{aligned}
&xu_{xx}+yu_{xy}+(c-x)u_x+zu_z-bu=0,\\&
y(1-y)u_{yy}+xu_{xy}-yzu_{yz}+\left[c-\left(a_1+a_2+1\right)y\right]u_y-a_2zu_z-a_1a_2u=0,\\&
zu_{zz}-xu_{xz}+yu_y+(1-b+z)u_z+a_1u=0,
\end{aligned}
\right.
$

where $u\equiv \,\,{\rm{E}}_{108}\left(a_1,a_2, b;c; x,y,z\right)$.

\bigskip

\begin{equation}
{\rm{E}}_{109}\left(a_1,a_2, b;c; x,y,z\right)=\sum\limits_{m,n,p=0}^\infty\frac{(a_1)_{n+p}(a_2)_p(b)_{m-p}}{(c)_{m+n}}\frac{x^m}{m!}\frac{y^n}{n!}\frac{z^p}{p!},
\end{equation}

region of convergence:
$$ \left\{ r<\infty,\,\,\,s<\infty,\,\,\,t<1
\right\}.
$$

System of partial differential equations:

$
\left\{
\begin{aligned}
&xu_{xx}+yu_{xy}+(c-x)u_x+zu_z-bu=0,\\&
yu_{yy}+xu_{xy}+(c-y)u_y-zu_z-a_1u=0,\\&
z(1+z)u_{zz}-xu_{xz}+yzu_{yz}+a_2yu_y+\left[1-b+\left(a_1+a_2+1\right)z\right]u_z+a_1a_2u=0,
\end{aligned}
\right.
$

where $u\equiv \,\,{\rm{E}}_{109}\left(a_1,a_2, b;c; x,y,z\right)$.

\bigskip

\begin{equation}
{\rm{E}}_{110}\left(a,b;c; x,y,z\right)=\sum\limits_{m,n,p=0}^\infty\frac{(a)_{n+p}(b)_{m-p}}{(c)_{m+n}}\frac{x^m}{m!}\frac{y^n}{n!}\frac{z^p}{p!},
\end{equation}

region of convergence:
$$ \left\{ r<\infty,\,\,\,s<\infty,\,\,\,t<\infty
\right\}.
$$

System of partial differential equations:

$
\left\{
\begin{aligned}
&xu_{xx}+yu_{xy}+(c-x)u_x+zu_z-bu=0,\\&
yu_{yy}+xu_{xy}+(c-y)u_y-zu_z-au=0,\\&
zu_{zz}-xu_{xz}+yu_y+(1-b+z)u_z+au=0,
\end{aligned}
\right.
$

where $u\equiv \,\,{\rm{E}}_{110}\left(a,b;c; x,y,z\right)$.

\bigskip

\begin{equation}
{\rm{E}}_{111}\left(a_1, a_2, b; c; x,y,z\right)=\sum\limits_{m,n,p=0}^\infty\frac{(a_1)_{n+p}(a_2)_{n+p}(b)_{m-p}}{(c)_{m+n}}\frac{x^m}{m!}\frac{y^n}{n!}\frac{z^p}{p!},
\end{equation}

region of convergence:
$$ \left\{ r<\infty,\,\,\,\sqrt{s}+\sqrt{t}<1
\right\}.
$$

System of partial differential equations:

$
\left\{
\begin{aligned}
&xu_{xx}+yu_{xy}+(c-x)u_x+zu_z-bu=0,\\&
y(1-y)u_{yy}-z^2u_{zz}+xu_{xy}-2yzu_{yz}\\& \,\,\,\,\,\,\,\,\,+\left[c-\left(a_1+a_2+1\right)y\right]u_y-\left(a_1+a_2+1\right)zu_z-a_1a_2u=0,\\&
z(1+z)u_{zz}+y^2u_{zz}-xu_{xz}+2yzu_{yz}\\& \,\,\,\,\,\,\,\,\,+\left(a_1+a_2+1\right)yu_y+\left[1-b+\left(a_1+a_2+1\right)z\right]u_z+a_1a_2u=0,
\end{aligned}
\right.
$

where $u\equiv \,\,{\rm{E}}_{111}\left(a_1, a_2, b; c; x,y,z\right)$.

\bigskip

\begin{equation}
{\rm{E}}_{112}\left(a_1,a_2, b; c; x,y,z\right)=\sum\limits_{m,n,p=0}^\infty\frac{(a_1)_{n+p}(a_2)_{m}(b)_{m+n-p}}{(c)_{m+n}}\frac{x^m}{m!}\frac{y^n}{n!}\frac{z^p}{p!},
\end{equation}

region of convergence:
$$ \left\{ r<1, \,\,\,\, s<1,\,\,\,t<\infty
\right\}.
$$

System of partial differential equations:

$
\left\{
\begin{aligned}
&x(1-x)u_{xx}+(1-x)yu_{xy}+xzu_{xz}+\left[c-\left(a_2+b+1\right)x\right]u_x-a_2yu_y+a_2zu_z-a_2bu=0,\\&
y(1-y)u_{yy}+z^2u_{zz}+x(1-y)u_{xy}\\&\,\,\,\,\,\,\,\,\,-xzu_{xz}-a_1xu_x+\left[c-\left(a_1+b+1\right)y\right]u_y+\left(a_1-b+1\right)zu_z-a_1bu=0,\\&
zu_{zz}-xu_{xz}-yu_{yz}+yu_y+\left(1-b+z\right)u_z+a_1u=0,
\end{aligned}
\right.
$

where $u\equiv \,\,{\rm{E}}_{112}\left(a_1,a_2, b; c; x,y,z\right)$.

\bigskip

\begin{equation}
{\rm{E}}_{113}\left(a_1, a_2, b; c; x,y,z\right)=\sum\limits_{m,n,p=0}^\infty\frac{(a_1)_{n+p}(a_2)_p(b)_{m+n-p}}{(c)_{m+n}}\frac{x^m}{m!}\frac{y^n}{n!}\frac{z^p}{p!},
\end{equation}

region of convergence:
$$ \left\{ r<\infty,\,\,\,\left\{t+2\sqrt{st}<1\right\}=\left\{\sqrt{t}<\sqrt{1+s}-\sqrt{s}\right\}
\right\}.
$$

System of partial differential equations:

$
\left\{
\begin{aligned}
&xu_{xx}+yu_{xy}+\left(c-x\right)u_x-yu_y+zu_z-bu=0,\\&
y(1-y)u_{yy}+z^2u_{zz}+x(1-y)u_{xy}\\&\,\,\,\,\,\,\,\,\,-xzu_{xz}-a_1xu_x+\left[c-\left(a_1+b+1\right)y\right]u_y+\left(a_1-b+1\right)zu_z-a_1bu=0,\\&
z(1+z)u_{zz}-xu_{xz}-y(1-z)u_{yz}+a_2yu_y+\left[1-b+\left(a_1+a_2+1\right)z\right]u_z+a_1a_2u=0,
\end{aligned}
\right.
$

where $u\equiv \,\,{\rm{E}}_{113}\left(a_1, a_2, b; c; x,y,z\right)$.

\bigskip

\begin{equation}
{\rm{E}}_{114}\left(a, b; c; x,y,z\right)=\sum\limits_{m,n,p=0}^\infty\frac{(a)_{n+p}(b)_{m+n-p}}{(c)_{m+n}}\frac{x^m}{m!}\frac{y^n}{n!}\frac{z^p}{p!},
\end{equation}

region of convergence:
$$ \left\{ r<\infty,\,\,\,s<1,\,\,\,r<\infty
\right\}.
$$

System of partial differential equations:

$
\left\{
\begin{aligned}
&xu_{xx}+yu_{xy}+\left(c-x\right)u_x-yu_y+zu_z-bu=0,\\&
y(1-y)u_{yy}+z^2u_{zz}+x(1-y)u_{xy}\\&\,\,\,\,\,\,\,\,\,-xzu_{xz}-axu_x+\left[c-\left(a+b+1\right)y\right]u_y+\left(a-b+1\right)zu_z-abu=0,\\&
zu_{zz}-xu_{xz}-yu_{yz}+yu_y+\left(1-b+z\right)u_z+au=0,
\end{aligned}
\right.
$

where $u\equiv \,\,{\rm{E}}_{114}\left(a, b; c; x,y,z\right)$.

\bigskip

\begin{equation}
{\rm{E}}_{115}\left(a_1,a_2, b; c; x,y,z\right)=\sum\limits_{m,n,p=0}^\infty\frac{(a_1)_{m+n}(a_2)_{n+p}(b)_{m-p}}{(c)_{m+n}}\frac{x^m}{m!}\frac{y^n}{n!}\frac{z^p}{p!},
\end{equation}

region of convergence:
$$ \left\{ r<1, \,\,\,\, s<1,\,\,\,t<\infty
\right\}.
$$

System of partial differential equations:

$
\left\{
\begin{aligned}
&x(1-x)u_{xx}+(1-x)yu_{xy}+xzu_{xz}+yzu_{yz}\\& \,\,\,\,\,\,\,\,\,+\left[c-\left(a_1+b+1\right)x\right]u_x-byu_y+a_1zu_z-a_1bu=0,\\&
y(1-y)u_{yy}+x(1-y)u_{xy}-xzu_{xz}-yzu_{yz}\\& \,\,\,\,\,\,\,\,\,-a_2xu_x+\left[c-\left(a_1+a_2+1\right)y\right]u_y-a_1zu_z-a_1a_2u=0,\\&
zu_{zz}-xu_{xz}+yu_y+\left(1-b+z\right)u_z+a_2u=0,
\end{aligned}
\right.
$

where $u\equiv \,\,{\rm{E}}_{115}\left(a_1,a_2, b; c; x,y,z\right)$.

\bigskip

\begin{equation}
{\rm{E}}_{116}\left(a_1,a_2, b; c; x,y,z\right)=\sum\limits_{m,n,p=0}^\infty\frac{(a_1)_{m+p}(a_2)_{n+p}(b)_{m-p}}{(c)_{m+n}}\frac{x^m}{m!}\frac{y^n}{n!}\frac{z^p}{p!},
\end{equation}

region of convergence:
$$ \left\{ \left\{t+2\sqrt{rt}<1\right\}=\left\{\sqrt{t}<\sqrt{1+r}-\sqrt{r}\right\},\,\,\,s<\infty
\right\}.
$$

System of partial differential equations:

$
\left\{
\begin{aligned}
&x(1-x)u_{xx}+yu_{xy}+z^2u_{zz}+\left[c-\left(a_1+b+1\right)x\right]u_x+\left(a_1-b+1\right)zu_z-a_1bu=0,\\&
yu_{yy}+xu_{xy}+(c-y)u_y-zu_z-a_2u=0,\\&
z(1+z)u_{zz}+xyu_{xy}-x(1-z)u_{xz}+yzu_{yz}\\& \,\,\,\,\,\,\,\,\,+a_2xu_x+a_1yu_y+\left[1-b+\left(a_1+a_2+1\right)z\right]u_z+a_1a_2u=0,
\end{aligned}
\right.
$

where $u\equiv \,\,{\rm{E}}_{116}\left(a_1,a_2, b; c; x,y,z\right)$.

\bigskip

\begin{equation}
{\rm{E}}_{117}\left(a_1,a_2, b; c; x,y,z\right)=\sum\limits_{m,n,p=0}^\infty\frac{(a_1)_{m+n+p}(a_2)_{n}(b)_{m-p}}{(c)_{m+n}}\frac{x^m}{m!}\frac{y^n}{n!}\frac{z^p}{p!},
\end{equation}

region of convergence:
$$ \left\{ r<1, \,\,\,\,s<1,\,\,\,t<\infty
\right\}.
$$

System of partial differential equations:

$
\left\{
\begin{aligned}
&x(1-x)u_{xx}+(1-x)yu_{xy}+z^2u_{zz}+yzu_{yz}\\&\,\,\,\,\,\,\,\,\,+\left[c-\left(a_1+b+1\right)x\right]u_x-byu_y+\left(a_1-b+1\right)zu_z-a_1bu=0,\\&
y(1-y)u_{yy}+x(1-y)u_{xy}-yzu_{yz}-a_2xu_x+\left[c-\left(a_1+a_2+1\right)y\right]u_y-a_2zu_z-a_1a_2u=0,\\&
zu_{zz}-xu_{xz}+xu_x+yu_y+\left(1-b+z\right)u_z+a_1u=0,
\end{aligned}
\right.
$

where $u\equiv \,\,{\rm{E}}_{117}\left(a_1,a_2, b; c; x,y,z\right)$.

\bigskip

\begin{equation}
{\rm{E}}_{118}\left(a_1,a_2, b; c; x,y,z\right)=\sum\limits_{m,n,p=0}^\infty\frac{(a_1)_{m+n+p}(a_2)_p(b)_{m-p}}{(c)_{m+n}}\frac{x^m}{m!}\frac{y^n}{n!}\frac{z^p}{p!},
\end{equation}

region of convergence:
$$ \left\{ \left\{t+2\sqrt{rt}<1\right\}=\left\{\sqrt{t}<\sqrt{1+r}-\sqrt{r}\right\},\,\,\,s<\infty
\right\}.
$$

System of partial differential equations:

$
\left\{
\begin{aligned}
&x(1-x)u_{xx}+(1-x)yu_{xy}+z^2u_{zz}+yzu_{yz}\\& \,\,\,\,\,\,\,\,+\left[c-\left(a_1+b+1\right)x\right]u_x-byu_y+\left(a_1-b+1\right)zu_z-a_1bu=0,\\&
yu_{yy}+xu_{xy}-xu_x+(c-y)u_y-zu_z-a_1u=0,\\&
z(1+z)u_{zz}-x(1-z)u_{xz}+yzu_{yz}\\& \,\,\,\,\,\,\,\,\,+a_2xu_x+a_2yu_y+\left[1-b+\left(a_1+a_2+1\right)z\right]u_z+a_1a_2u=0,
\end{aligned}
\right.
$

where $u\equiv \,\,{\rm{E}}_{118}\left(a_1,a_2, b; c; x,y,z\right)$.

\bigskip

\begin{equation}
{\rm{E}}_{119}\left(a, b; c; x,y,z\right)=\sum\limits_{m,n,p=0}^\infty\frac{(a)_{m+n+p}(b)_{m-p}}{(c)_{m+n}}\frac{x^m}{m!}\frac{y^n}{n!}\frac{z^p}{p!},
\end{equation}

region of convergence:
$$ \left\{ r<1,\,\,\,s<\infty,\,\,\,t<\infty
\right\}.
$$

System of partial differential equations:

$
\left\{
\begin{aligned}
&x(1-x)u_{xx}+(1-x)yu_{xy}+z^2u_{zz}+yzu_{yz}\\& \,\,\,\,\,\,\,\,\,+\left[c-\left(a+b+1\right)x\right]u_x-byu_y+\left(a-b+1\right)zu_z-abu=0,\\&
yu_{yy}+xu_{xy}-xu_x+(c-y)u_y-zu_z-au=0,\\&
zu_{zz}-xu_{xz}+xu_x+yu_y+\left(1-b+z\right)u_z+au=0,
\end{aligned}
\right.
$

where $u\equiv \,\,{\rm{E}}_{119}\left(a, b; c; x,y,z\right)$.

\bigskip

\begin{equation}
{\rm{E}}_{120}\left(a_1,a_2,b; c; x,y,z\right)=\sum\limits_{m,n,p=0}^\infty\frac{(a_1)_{n+p}(a_2)_{n}(b)_{2m-p}}{(c)_{m+n}}\frac{x^m}{m!}\frac{y^n}{n!}\frac{z^p}{p!},
\end{equation}

region of convergence:
$$ \left\{ \frac{1}{4r}+\frac{1}{s}>1,\,\,\,t<\infty
\right\}.
$$

System of partial differential equations:

$
\left\{
\begin{aligned}
&x(1-4x)u_{xx}+yu_{xy}-z^2u_{zz}+4xzu_{xz}+\left[c-(4b+6)x\right]u_x+2bzu_z-b(b+1)u=0,\\&
y(1-y)u_{yy}+xu_{xy}-yzu_{yz}+\left[c-\left(a_1+a_2+1\right)y\right]u_y-a_2zu_z-a_1a_2u=0,\\&
zu_{zz}-2xu_{xz}+yu_y+\left(1-b+z\right)u_z+a_1u=0,
\end{aligned}
\right.
$

where $u\equiv \,\,{\rm{E}}_{120}\left(a_1,a_2,b; c; x,y,z\right)$.

\bigskip

\begin{equation}
{\rm{E}}_{121}\left(a_1, a_2,b; c; x,y,z\right)=\sum\limits_{m,n,p=0}^\infty\frac{(a_1)_{n+p}(a_2)_p(b)_{2m-p}}{(c)_{m+n}}\frac{x^m}{m!}\frac{y^n}{n!}\frac{z^p}{p!},
\end{equation}

region of convergence:
$$ \left\{ \left(1+2\sqrt{r}\right)t<1,\,\,\,s<\infty
\right\}.
$$

System of partial differential equations:

$
\left\{
\begin{aligned}
&x(1-4x)u_{xx}+yu_{xy}-z^2u_{zz}+4xzu_{xz}+\left[c-(4b+6)x\right]u_x+2bzu_z-b(b+1)u=0,\\&
yu_{yy}+xu_{xy}+(c-y)u_y-zu_z-a_1u=0,\\&
z(1+z)u_{zz}-2xu_{xz}+yzu_{yz}+a_2yu_y+\left[1-b+\left(a_1+a_2+1\right)z\right]u_z+a_1a_2u=0,
\end{aligned}
\right.
$

where $u\equiv \,\,{\rm{E}}_{121}\left(a_1, a_2,b; c; x,y,z\right)$.

\bigskip

\begin{equation}
{\rm{E}}_{122}\left(a_1,a_2,b; c; x,y,z\right)=\sum\limits_{m,n,p=0}^\infty\frac{(a_1)_{n}(a_2)_p(b)_{2m-p}}{(c)_{m+n}}\frac{x^m}{m!}\frac{y^n}{n!}\frac{z^p}{p!},
\end{equation}

region of convergence:
$$ \left\{ r<\frac{1}{4},\,\,\,s<\infty,\,\,\,t<\infty
\right\}.
$$

System of partial differential equations:

$
\left\{
\begin{aligned}
&x(1-4x)u_{xx}+yu_{xy}-z^2u_{zz}+4xzu_{xz}+\left[c-(4b+6)x\right]u_x+2bzu_z-b(b+1)u=0,\\&
yu_{yy}+xu_{xy}+(c-y)u_y-a_1u=0,\\&
zu_{zz}-2xu_{xz}+(1-b+z)u_z+a_2u=0,
\end{aligned}
\right.
$

where $u\equiv \,\,{\rm{E}}_{122}\left(a_1,a_2,b; c; x,y,z\right)$.

\bigskip

\begin{equation}
{\rm{E}}_{123}\left(a,b; c; x,y,z\right)=\sum\limits_{m,n,p=0}^\infty\frac{(a)_{n}(b)_{2m-p}}{(c)_{m+n}}\frac{x^m}{m!}\frac{y^n}{n!}\frac{z^p}{p!},
\end{equation}

region of convergence:
$$ \left\{ r<\frac{1}{4},\,\,\,s<\infty,\,\,\,t<\infty
\right\}.
$$

System of partial differential equations:

$
\left\{
\begin{aligned}
&x(1-4x)u_{xx}+yu_{xy}-z^2u_{zz}+4xzu_{xz}+\left[c-(4b+6)x\right]u_x+2bzu_z-b(b+1)u=0,\\&
yu_{yy}+xu_{xy}+(c-y)u_y-au=0,\\&
zu_{zz}-2xu_{xz}+(1-b)u_z+u=0,
\end{aligned}
\right.
$

where $u\equiv \,\,{\rm{E}}_{123}\left(a,b; c; x,y,z\right)$.

\bigskip

\begin{equation}
{\rm{E}}_{124}\left(a,b; c; x,y,z\right)=\sum\limits_{m,n,p=0}^\infty\frac{(a)_p(b)_{2m-p}}{(c)_{m+n}}\frac{x^m}{m!}\frac{y^n}{n!}\frac{z^p}{p!},
\end{equation}

region of convergence:
$$ \left\{ r<\frac{1}{4},\,\,\,s<\infty,\,\,\,t<\infty
\right\}.
$$

System of partial differential equations:

$
\left\{
\begin{aligned}
&x(1-4x)u_{xx}+yu_{xy}-z^2u_{zz}+4xzu_{xz}+\left[c-(4b+6)x\right]u_x+2bzu_z-b(b+1)u=0,\\&
yu_{yy}+xu_{xy}+cu_y-u=0,\\&
zu_{zz}-2xu_{xz}+(1-b+z)u_z+au=0,
\end{aligned}
\right.
$

where $u\equiv \,\,{\rm{E}}_{124}\left(a,b; c; x,y,z\right)$.

\bigskip

\begin{equation}
{\rm{E}}_{125}\left(b; c; x,y,z\right)=\sum\limits_{m,n,p=0}^\infty\frac{(b)_{2m-p}}{(c)_{m+n}}\frac{x^m}{m!}\frac{y^n}{n!}\frac{z^p}{p!},
\end{equation}

region of convergence:
$$ \left\{ r<\frac{1}{4},\,\,\,s<\infty,\,\,\,t<\infty
\right\}.
$$

System of partial differential equations:

$
\left\{
\begin{aligned}
&x(1-4x)u_{xx}+yu_{xy}-z^2u_{zz}+4xzu_{xz}+\left[c-(4b+6)x\right]u_x+2bzu_z-b(b+1)u=0,\\&
yu_{yy}+xu_{xy}+cu_y-u=0,\\&
zu_{zz}-2xu_{xz}+(1-b)u_z+u=0,
\end{aligned}
\right.
$

where $u\equiv \,\,{\rm{E}}_{125}\left(b; c; x,y,z\right)$.

\bigskip

\begin{equation}
{\rm{E}}_{126}\left(a_1,a_2,b; c; x,y,z\right)=\sum\limits_{m,n,p=0}^\infty\frac{(a_1)_{n}(a_2)_{p}(b)_{2m+n-p}}{(c)_{m+n}}\frac{x^m}{m!}\frac{y^n}{n!}\frac{z^p}{p!},
\end{equation}

region of convergence:
$$ \left\{ \left\{r<\frac{1}{4}, \,\,\,s<\frac{1}{2}+\frac{1}{2}\sqrt{1-4r}\right\}\cup\left\{s \leq \frac{1}{2}\right\},\,\,\,t<\infty
\right\}.
$$

System of partial differential equations:

$
\left\{
\begin{aligned}
&x(1-4x)u_{xx}+(1-4x)yu_{xy}+4xzu_{xz}\\
&\,\,\,\,\,\,\,\,+2yzu_{yz}-y^2u_{yy}-z^2u_{zz}+\left[c-(4b+6)x\right]u_x-2(b+1)yu_y+2bzu_z-b(b+1)u=0,\\
&y(1-y)u_{yy}+x(1-2y)u_{xy}+yzu_{yz}-2a_1xu_x+\left[c-\left(a_1+b+1\right)y\right]u_y+a_1zu_z-a_1bu=0,\\
&zu_{zz}-2xu_{xz}-yu_{yz}+(1-b+z)u_z+a_2u=0,
\end{aligned}
\right.
$

where $u\equiv \,\,{\rm{E}}_{126}\left(a_1,a_2,b; c; x,y,z\right)$.

\bigskip

\begin{equation}
{\rm{E}}_{127}\left(a_1,a_2,b; c; x,y,z\right)=\sum\limits_{m,n,p=0}^\infty\frac{(a_1)_{p}(a_2)_p(b)_{2m+n-p}}{(c)_{m+n}}\frac{x^m}{m!}\frac{y^n}{n!}\frac{z^p}{p!},
\end{equation}

region of convergence:
$$ \left\{ \left(1+2\sqrt{r}\right)t<1,\,\,\,s<\infty
\right\}.
$$

System of partial differential equations:

$
\left\{
\begin{aligned}
&x(1-4x)u_{xx}+(1-4x)yu_{xy}+4xzu_{xz}\\&\,\,\,\,\,\,\,\,\,+2yzu_{yz}-y^2u_{yy}-z^2u_{zz}+\left[c-(4b+6)x\right]u_x-2(b+1)yu_y+2bzu_z-b(b+1)u=0,\\&
yu_{yy}+xu_{xy}-2xu_x+(c-y)u_y+zu_z-bu=0,\\&
z(1+z)u_{zz}-2xu_{xz}-yu_{yz}+\left[1-b+\left(a_1+a_2+1\right)z\right]u_z+a_1a_2u=0,
\end{aligned}
\right.
$

where $u\equiv \,\,{\rm{E}}_{127}\left(a_1,a_2,b; c; x,y,z\right)$.

\bigskip

\begin{equation}
{\rm{E}}_{128}\left(a,b; c; x,y,z\right)=\sum\limits_{m,n,p=0}^\infty\frac{(a)_{n}(b)_{2m+n-p}}{(c)_{m+n}}\frac{x^m}{m!}\frac{y^n}{n!}\frac{z^p}{p!},
\end{equation}

region of convergence:
$$ \left\{\left\{ r<\frac{1}{4}, \,\,\,s<\frac{1}{2}+\frac{1}{2}\sqrt{1-4r}\right\}\cup\left\{s \leq \frac{1}{2}\right\},\,\,\,t<\infty
\right\}.
$$

System of partial differential equations:

$
\left\{
\begin{aligned}
&x(1-4x)u_{xx}+(1-4x)yu_{xy}+4xzu_{xz}\\
&\,\,\,\,\,\,\,\,\,+2yzu_{yz}-y^2u_{yy}-z^2u_{zz}+\left[c-(4b+6)x\right]u_x-2(b+1)yu_y+2bzu_z-b(b+1)u=0,\\
&y(1-y)u_{yy}+x(1-2y)u_{xy}+yzu_{yz}-2axu_x+\left[c-\left(a+b+1\right)y\right]u_y+azu_z-abu=0,\\
&zu_{zz}-2xu_{xz}-yu_{yz}+(1-b)u_z+u=0,
\end{aligned}
\right.
$

where $u\equiv \,\,{\rm{E}}_{128}\left(a,b; c; x,y,z\right)$.

\bigskip

\begin{equation}
{\rm{E}}_{129}\left(a,b; c; x,y,z\right)=\sum\limits_{m,n,p=0}^\infty\frac{(a)_{p}(b)_{2m+n-p}}{(c)_{m+n}}\frac{x^m}{m!}\frac{y^n}{n!}\frac{z^p}{p!},
\end{equation}

region of convergence:
$$ \left\{ r<\frac{1}{4},\,\,\,s<\infty,\,\,\,t<\infty
\right\}.
$$

System of partial differential equations:

$
\left\{
\begin{aligned}
&x(1-4x)u_{xx}+(1-4x)yu_{xy}+4xzu_{xz}\\
&\,\,\,\,\,\,\,\,\,+2yzu_{yz}-y^2u_{yy}-z^2u_{zz}+\left[c-(4b+6)x\right]u_x-2(b+1)yu_y+2bzu_z-b(b+1)u=0,\\
&yu_{yy}+xu_{xy}-2xu_x+(c-y)u_y+zu_z-bu=0,\\
&zu_{zz}-2xu_{xz}-yu_{yz}+(1-b+z)u_z+au=0,
\end{aligned}
\right.
$

where $u\equiv \,\,{\rm{E}}_{129}\left(a,b; c; x,y,z\right)$.

\bigskip

\begin{equation}
{\rm{E}}_{130}\left(b; c; x,y,z\right)=\sum\limits_{m,n,p=0}^\infty\frac{(b)_{2m+n-p}}{(c)_{m+n}}\frac{x^m}{m!}\frac{y^n}{n!}\frac{z^p}{p!},
\end{equation}

region of convergence:
$$ \left\{ r<\frac{1}{4},\,\,\,s<\infty,\,\,\,t<\infty
\right\}.
$$

System of partial differential equations:

$
\left\{
\begin{aligned}
&x(1-4x)u_{xx}+(1-4x)yu_{xy}+4xzu_{xz}\\&\,\,\,\,\,\,\,\,\,+2yzu_{yz}-y^2u_{yy}-z^2u_{zz}+\left[c-(4b+6)x\right]u_x-2(b+1)yu_y+2bzu_z-b(b+1)u=0,\\&
yu_{yy}+xu_{xy}-2xu_x+(c-y)u_y+zu_z-bu=0,\\&
zu_{zz}-2xu_{xz}-yu_{yz}+(1-b)u_z+u=0,
\end{aligned}
\right.
$

where $u\equiv \,\,{\rm{E}}_{130}\left(b; c; x,y,z\right)$.

\bigskip

\begin{equation}
{\rm{E}}_{131}\left(a_1,a_2,b; c; x,y,z\right)=\sum\limits_{m,n,p=0}^\infty\frac{(a_1)_{2m+n}(a_2)_{p}(b)_{n-p}}{(c)_{m+n}}\frac{x^m}{m!}\frac{y^n}{n!}\frac{z^p}{p!},
\end{equation}

region of convergence:
$$ \left\{ \left\{r<\frac{1}{4}, \,\,\,\, s<\frac{1}{2}+\frac{1}{2}\sqrt{1-4r}\right\}\cup\left\{s \leq \frac{1}{2}\right\},\,\,\,t<\infty
\right\}.
$$

System of partial differential equations:

$
\left\{
\begin{aligned}
&x(1-4x)u_{xx}+(1-4x)yu_{xy}-y^2u_{yy}\\& \,\,\,\,\,\,\,\,\,+\left[c-\left(4a_1+6\right)x\right]u_x-2\left(a_1+1\right)yu_y-a_1\left(a_1+1\right)u=0,\\&
y(1-y)u_{yy}+x(1-2y)u_{xy}+2xzu_{xz}+yzu_{yz}\\& \,\,\,\,\,\,\,\,\,-2bxu_x+\left[c-\left(a_1+b+1\right)y\right]u_y+a_1zu_z-a_1bu=0,\\&
zu_{zz}-yu_{yz}+(1-b+z)u_z+a_2u=0,
\end{aligned}
\right.
$

where $u\equiv \,\,{\rm{E}}_{131}\left(a_1,a_2,b; c; x,y,z\right)$.

\bigskip

\begin{equation}
{\rm{E}}_{132}\left(a,b; c; x,y,z\right)=\sum\limits_{m,n,p=0}^\infty\frac{(a)_{2m+n}(b)_{n-p}}{(c)_{m+n}}\frac{x^m}{m!}\frac{y^n}{n!}\frac{z^p}{p!},
\end{equation}

region of convergence:
$$ \left\{ \left\{r<\frac{1}{4}, \,\,\,\, s<\frac{1}{2}+\frac{1}{2}\sqrt{1-4r}\right\}\cup\left\{s \leq \frac{1}{2}\right\},\,\,\,t<\infty
\right\}.
$$

System of partial differential equations:

$
\left\{
\begin{aligned}
&x(1-4x)u_{xx}+(1-4x)yu_{xy}-y^2u_{yy}+\left[c-(4a+6)x\right]u_x-2(a+1)yu_y-a(a+1)u=0,\\&
y(1-y)u_{yy}+x(1-2y)u_{xy}+2xzu_{xz}+yzu_{yz}\\& \,\,\,\,\,\,\,\,\,-2bxu_x+\left[c-\left(a+b+1\right)y\right]u_y+azu_z-abu=0,\\&
zu_{zz}-yu_{yz}+(1-b)u_z+u=0,
\end{aligned}
\right.
$

where $u\equiv \,\,{\rm{E}}_{132}\left(a,b; c; x,y,z\right)$.

\bigskip

\begin{equation}
{\rm{E}}_{133}\left(a_1,a_2, b; c; x,y,z\right)=\sum\limits_{m,n,p=0}^\infty\frac{(a_1)_{2m+p}(a_2)_{n} (b)_{n-p}}{(c)_{m+n}}\frac{x^m}{m!}\frac{y^n}{n!}\frac{z^p}{p!},
\end{equation}

region of convergence:
$$ \left\{ \frac{1}{4r}+\frac{1}{s}>1,\,\,\,\,t<\infty
\right\}.
$$

System of partial differential equations:

$
\left\{
\begin{aligned}
&x(1-4x)u_{xx}+yu_{xy}-4xzu_{xz}-z^2u_{zz}\\& \,\,\,\,\,\,\,\,\,+\left[c-(4a_1+6)x\right]u_x-2\left(a_1+1\right)zu_z-a_1\left(a_1+1\right)u=0,\\&
y(1-y)u_{yy}+xu_{xy}+yzu_{yz}+\left[c-\left(a_2+b+1\right)y\right]u_y+a_2zu_z-a_2bu=0,\\&
zu_{zz}-yu_{yz}+2xu_x+(1-b+z)u_z+a_1u=0,
\end{aligned}
\right.
$

where $u\equiv \,\,{\rm{E}}_{133}\left(a_1,a_2, b; c; x,y,z\right)$.

\bigskip

\begin{equation}
{\rm{E}}_{134}\left(a_1, a_2, b; c; x,y,z\right)=\sum\limits_{m,n,p=0}^\infty\frac{(a_1)_{2m+p}(a_2)_p (b)_{n-p}}{(c)_{m+n}}\frac{x^m}{m!}\frac{y^n}{n!}\frac{z^p}{p!},
\end{equation}

region of convergence:
$$ \left\{ 2\sqrt{r}+t<1,\,\,\,s<\infty
\right\}.
$$

System of partial differential equations:

$
\left\{
\begin{aligned}
&x(1-4x)u_{xx}+yu_{xy}-4xzu_{xz}-z^2u_{zz}\\& \,\,\,\,\,\,\,\,\,+\left[c-(4a_1+6)x\right]u_x-2\left(a_1+1\right)zu_z-a_1\left(a_1+1\right)u=0,\\&
yu_{yy}+xu_{xy}+(c-y)u_y+zu_z-bu=0,\\&
z(1+z)u_{zz}+2xzu_{xz}-yu_{yz}+2a_2u_x+\left[1-b+\left(a_1+a_2+1\right)z\right]u_z+a_1a_2u=0,
\end{aligned}
\right.
$

where $u\equiv \,\,{\rm{E}}_{134}\left(a_1, a_2, b; c; x,y,z\right)$.

\bigskip

\begin{equation}
{\rm{E}}_{135}\left(a, b; c; x,y,z\right)=\sum\limits_{m,n,p=0}^\infty\frac{(a)_{2m+p} (b)_{n-p}}{(c)_{m+n}}\frac{x^m}{m!}\frac{y^n}{n!}\frac{z^p}{p!},
\end{equation}

region of convergence:
$$ \left\{ r<\frac{1}{4},\,\,\,s<\infty,\,\,\,t<\infty
\right\}.
$$

System of partial differential equations:

$
\left\{
\begin{aligned}
&x(1-4x)u_{xx}+yu_{xy}-4xzu_{xz}-z^2u_{zz}+\left[c-(4a+6)x\right]u_x-2\left(a+1\right)zu_z-a\left(a+1\right)u=0,\\&
yu_{yy}+xu_{xy}+(c-y)u_y+zu_z-bu=0,\\&
zu_{zz}-yu_{yz}+2xu_x+(1-b+z)u_z+au=0,
\end{aligned}
\right.
$

where $u\equiv \,\,{\rm{E}}_{135}\left(a, b; c; x,y,z\right)$.

\bigskip

\begin{equation}
{\rm{E}}_{136}\left(a_1,a_2, b; c; x,y,z\right)=\sum\limits_{m,n,p=0}^\infty\frac{(a_1)_{n+2p}(a_2)_{m}(b)_{n-p}}{(c)_{m+n}}\frac{x^m}{m!}\frac{y^n}{n!}\frac{z^p}{p!},
\end{equation}

region of convergence:
$$ \left\{ r<\infty,\,\,\,\left\{t<\frac{1}{4},\,\,\,\, s<\min\left\{\Psi_1(t), \Psi_2(t)\right\}\right\}=\left\{s<1,\,\,\,\, t<\min\left\{\Theta_1(s), \Theta_2(s)\right\}\right\}
\right\}.
$$

System of partial differential equations:

$
\left\{
\begin{aligned}
&xu_{xx}+yu_{xy}+(c-x)u_x-a_2u=0,\\&
y(1-y)u_{yy}+2z^2u_{zz}+xu_{xy}-yzu_{yz}+\left[c-\left(a_1+b+1\right)y\right]u_y+\left(a_1-2b+2\right)zu_z-a_1bu=0,\\&
z(1+4z)u_{zz}+y^2u_{yy}-y(1-4z)u_{yz}\\& \,\,\,\,\,\,\,\,\,+2\left(a_1+1\right)yu_y+\left[1-b+\left(4a_1+6\right)z\right]u_z+a_1\left(a_1+1\right)u=0,
\end{aligned}
\right.
$

where $u\equiv \,\,{\rm{E}}_{136}\left(a_1,a_2, b; c; x,y,z\right)$.

\bigskip

\begin{equation}
{\rm{E}}_{137}\left(a, b; c; x,y,z\right)=\sum\limits_{m,n,p=0}^\infty\frac{(a)_{n+2p}(b)_{n-p}}{(c)_{m+n}}\frac{x^m}{m!}\frac{y^n}{n!}\frac{z^p}{p!},
\end{equation}

region of convergence:
$$ \left\{ r<\infty,\,\,\,\left\{t<\frac{1}{4}, \,\,\,\, s<\min\left\{\Psi_1(t), \Psi_2(t)\right\}\right\}=\left\{s<1, \,\,\,\, t<\min\left\{\Theta_1(s), \Theta_2(s)\right\}\right\}
\right\}.
$$

System of partial differential equations:

$
\left\{
\begin{aligned}
&{xu_{xx} +yu_{xy} +cu_{x} -u=0,} \\& {y(1-y)u_{yy} +2z^{2} u_{zz} +xu_{xy} -yzu_{yz}+ \left[c-\left(a +b+1\right)y\right]u_y} {+ \left(a -2b+2\right)zu_z-a bu}=0, \\
&{z(1+4z)u_{zz}+y^2u_{yy} -y(1-4z)u_{yz} +2\left(a+1\right)yu_y} {+\left[1-b+\left(4a+6\right)z\right]u_z+a \left(a+1\right)u=0},
\end{aligned}
\right.
$

where $u\equiv \,\,{\rm{E}}_{137}\left(a, b; c; x,y,z\right)$.

\bigskip

\begin{equation}
{\rm{E}}_{138}\left(a_1,a_2, b; c; x,y,z\right)=\sum\limits_{m,n,p=0}^\infty\frac{(a_1)_{n+2p}(a_2)_{m}(b)_{m-p}}{(c)_{m+n}}\frac{x^m}{m!}\frac{y^n}{n!}\frac{z^p}{p!},
\end{equation}

region of convergence:
$$ \left\{ (1+r)t<1,\,\,\,s<\infty
\right\}.
$$

System of partial differential equations:

$
\left\{
\begin{aligned}
&x(1-x)u_{xx}  +yu_{xy}+xzu_{xz} +\left[c-\left(a_{2} +b+1\right)x\right]u_{x}+a_{2} zu_{z} -a_{2} bu=0, \\& 
{yu_{yy} +xu_{xy} + (c-y)u_y-2zu_{z} -a_{1} u=0,} \\
&z(1+4z)u_{zz} -xu_{xz}+4yzu_{yz} +y^{2} u_{yy} \\& \,\,\,\,\,\,\,\,\,+2 \left(a_{1} +1\right)yu_y+ [1-b+\left(4a_{1} +6\right)z]u_z+a_{1}\left(a_{1} +1\right)u=0,
\end{aligned}
\right.
$

where $u\equiv \,\,{\rm{E}}_{138}\left(a_1,a_2, b; c; x,y,z\right)$.

\bigskip

\begin{equation}
{\rm{E}}_{139}\left(a_1, a_2, b; c; x,y,z\right)=\sum\limits_{m,n,p=0}^\infty\frac{(a_1)_{n+2p}(a_2)_n(b)_{m-p}}{(c)_{m+n}}\frac{x^m}{m!}\frac{y^n}{n!}\frac{z^p}{p!},
\end{equation}

region of convergence:
$$ \left\{r<\infty,\,\,\,2\sqrt{t}+s<1
\right\}.
$$

System of partial differential equations:

$
\left\{
\begin{aligned}
&{xu_{xx} +yu_{xy} +(c-x)u_{x} +zu_{z} -bu=0,} \\ &{y(1-y)u_{yy} +xu_{xy} -2yzu_{yz} + \left[c-\left(a_{1} +a_{2} +1\right)y\right]u_y}{-2a_{2} u_{z} -a_{1} a_{2} u=0,} \\
&z(1+4z)u_{zz} -xu_{xz}+4yzu_{yz} +y^{2} u_{yy}\\& \,\,\,\,\,\,\,\,\, +2 \left(a_{1} +1\right)yu_y+ [1-b+\left(4a_{1} +6\right)z]u_z+a_{1}\left(a_{1} +1\right)u=0,
\end{aligned}
\right.
$

where $u\equiv \,\,{\rm{E}}_{139}\left(a_1, a_2, b; c; x,y,z\right)$.

\bigskip

\begin{equation}
{\rm{E}}_{140}\left(a, b; c; x,y,z\right)=\sum\limits_{m,n,p=0}^\infty\frac{(a)_{n+2p}(b)_{m-p}}{(c)_{m+n}}\frac{x^m}{m!}\frac{y^n}{n!}\frac{z^p}{p!},
\end{equation}

region of convergence:
$$ \left\{ r<\infty,\,\,\,s<\infty,\,\,\,t<\frac{1}{4}
\right\}.
$$

System of partial differential equations:

$
\left\{
\begin{aligned}
&{xu_{xx} +yu_{xy} +(c-x)u_{x} +zu_{z} -bu=0,} \\
&{yu_{yy} +xu_{xy} + (c-y)u_y-2zu_{z} -au=0,} \\
&z(1+4z)u_{zz} -xu_{xz}+4yzu_{yz} +y^{2} u_{yy}\\& \,\,\,\,\,\,\,\,\, +2 \left(a_{1} +1\right)yu_y+ [1-b+\left(4a_{1} +6\right)z]u_z+a_{1}\left(a_{1} +1\right)u=0,
\end{aligned}
\right.
$

where $u\equiv \,\,{\rm{E}}_{140}\left(a, b; c; x,y,z\right)$.

\bigskip

\begin{equation}
{\rm{R}}_{141}\left(a, b; c; x,y,z\right)=\sum\limits_{m,n,p=0}^\infty\frac{(a)_{n+p}(b)_{2m-p}}{(c)_{m+n}}\frac{x^m}{m!}\frac{y^n}{n!}\frac{z^p}{p!},
\end{equation}

region of convergence:
$$ \left\{ r<\infty,\,\,\,s<\infty,\,\,\,t<\frac{1}{4}
\right\}.
$$

System of partial differential equations:

$
\left\{
\begin{aligned}
&{x(1-4x)u_{xx} -z^{2} u_{zz} +yu_{xy} +4xzu_{xz}} {+\left[c-(4b+6)x\right]u_{x} +2bzu_{z} -b(b+1)u=0,} \\ &{yu_{yy} +xu_{xy} + (c-y)u_y-zu_{z} -au=0,} \\& {zu_{zz} -2xu_{xz} +yu_{y} +(1-b+z)u_{z} +au=0,}
\end{aligned}
\right.
$

where $u\equiv \,\,{\rm{R}}_{141}\left(a, b; c; x,y,z\right)$.

\bigskip

\begin{equation}
{\rm{E}}_{142}\left(a, b; c; x,y,z\right)=\sum\limits_{m,n,p=0}^\infty\frac{(a)_{n+p}(b)_{2m+n-p}}{(c)_{m+n}}\frac{x^m}{m!}\frac{y^n}{n!}\frac{z^p}{p!},
\end{equation}

region of convergence:
$$ \left\{ r<\infty,\,\,\,s<\infty,\,\,\,t<\frac{1}{4}
\right\}.
$$

System of partial differential equations:

$
\left\{
\begin{aligned}
&{x(1-4x)u_{xx} -y^{2} u_{yy} -z^{2} u_{zz}}+(1-4x)yu_{xy} \\ &\,\,\,\,\,\,\,\,\,\,\,\, +4xzu_{xz} +2yzu_{yz}{+\left[c-(4b+6)x\right]u_{x} -2(b+1)yu_{y} +2bzu_{z} -b(b+1)u=0,} \\ & {y(1-y)u_{yy} +z^{2} u_{zz} +x(1-2y)u_{xy}  } \\&\,\,\,\,\,\,\,\,\,\,\,\, -2xzu_{xz}{-2axu_{x}+\left[c-(a+b+1)y\right]u_{y}  +(a-b+1)zu_{z} -abu=0,} \\ &{zu_{zz} -2xu_{xz} -yu_{yz} +yu_{y} +(1-b+z)u_{z} +au=0,}
\end{aligned}
\right.
$

where $u\equiv \,\,{\rm{E}}_{142}\left(a, b; c; x,y,z\right)$.

\bigskip

\begin{equation}
{\rm{E}}_{143}\left(a, b; c; x,y,z\right)=\sum\limits_{m,n,p=0}^\infty\frac{(a)_{n+2p}(b)_{m+n-p}}{(c)_{m+n}}\frac{x^m}{m!}\frac{y^n}{n!}\frac{z^p}{p!},
\end{equation}

region of convergence:
$$ \left\{ r<\infty,\,\,\,\left\{t<\frac{1}{4}, \,\,\, s<\min\left\{\Psi_1(t), \Psi_2(t)\right\}\right\}=\left\{s<1, \,\,\, t<\min\left\{\Theta_1(s), \Theta_2(s)\right\}\right\}
\right\}.
$$

System of partial differential equations:

$
\left\{
\begin{aligned}
&{xu_{xx} +yu_{xy} +(c-x)u_{x} -yu_{y} +zu_{z} -bu=0,} \\& {y(1-y)u_{yy} -2z^{2} u_{zz} }+x(1-y)u_{xy} \\
&\,\,\,\,\,\,\,\,\,\,\,\,\, -2xzu_{xz} -yzu_{yz} {-axu_{x} +\left[c-(a+b+1)y\right]u_{y} +(a-2b+2)zu_{z} -abu=0,} \\
&  z(1+4z)u_{zz}+{y^{2} u_{yy} -xu_{xz} }-y(1-4z)u_{yz} \\& \,\,\,\,\,\,\,\,\,+2(a+1)yu_{y}{ +\left[1-b+(4a+6)z\right]u_{z} +a(a+1)u=0,}
\end{aligned}
\right.
$

where $u\equiv \,\,{\rm{E}}_{143}\left(a, b; c; x,y,z\right)$.

\bigskip

\begin{equation}
{\rm{E}}_{144}\left(a, b; c; x,y,z\right)=\sum\limits_{m,n,p=0}^\infty\frac{(a)_{2m+n+p}(b)_{n-p}}{(c)_{m+n}}\frac{x^m}{m!}\frac{y^n}{n!}\frac{z^p}{p!},
\end{equation}

region of convergence:
$$ \left\{\left\{ r<\frac{1}{4},\,\,\, s<\frac{1}{2}+\frac{1}{2}\sqrt{1-4r}\right\}\cup\left\{s \leq \frac{1}{2}\right\},\,\,\,t<\infty
\right\}.
$$

System of partial differential equations:

$
\left\{
\begin{aligned}
&x(1-4x)u_{xx} -y^{2} u_{yy} -z^{2} u_{zz} +(1-4x)yu_{xy} -4xzu_{xz}  -2yzu_{yz} \\& \,\,\,\,\,\,\,\,\,+\left[c-(4a+6)x\right]u_{x} -2(a+1)yu_{y} -2(a+1)zu_{z} -a(a+1)u=0, \\ &y(1-y)u_{yy} +z^{2} u_{zz} +x(1-2y)u_{xy}\\&\,\,\,\,\,\,\,\,\,\,\,\,  +2xzu_{xz}-2bxu_{x}+\left[c-(a+b+1)y\right]u_{y} {  +(a-b+1)zu_{z} -abu=0,} \\& {zu_{zz} -yu_{yz} +2xu_{x} +yu_{y} +(1-b+z)u_{z} +au=0,}
\end{aligned}
\right.
$

where $u\equiv \,\,{\rm{E}}_{144}\left(a, b; c; x,y,z\right)$.

\bigskip

\begin{equation}
{\rm{E}}_{145}\left(a, b; c; x,y,z\right)=\sum\limits_{m,n,p=0}^\infty\frac{(a)_{m+n+2p}(b)_{n-p}}{(c)_{m+n}}\frac{x^m}{m!}\frac{y^n}{n!}\frac{z^p}{p!},
\end{equation}

region of convergence:
$$ \left\{ r<\infty,\,\,\,\left\{t<\frac{1}{4}, \,\,\,\, s<\min\left\{\Psi_1(t), \Psi_2(t)\right\}\right\}=\left\{s<1,\,\,\,\, t<\min\left\{\Theta_1(s), \Theta_2(s)\right\}\right\}
\right\}.
$$

System of partial differential equations:

$
\left\{
\begin{aligned}
&{xu_{xx} +yu_{xy} +(c-x)u_{x} -yu_{y} -2zu_{z} -au=0,} \\
&{y(1-y)u_{yy} +2z^{2} u_{zz} -x(1-y)u_{xy} } \\
&\,\,\,\,\,\,\,\,\,\,\,\,\, +xzu_{xz} -yzu_{yz} {-bxu_{x}+\left[c-(a+b+1)y\right]u_{y}  +(a-2b+2)zu_{z} -abu=0,} \\
&  z(1+4z)u_{zz}+{x^{2} u_{xx} +y^{2} u_{yy} +2xyu_{xy} +4xzu_{xz} -y(1-4z)u_{yz} } \\
& \,\,\,\,\,\,\,\,\,\,\,\,\,{+2(a+1)xu_{x}+2(a+1)yu_{y} +\left[1-b+(4a+6)z\right]u_{z} +a(a+1)u=0,}
\end{aligned}
\right.
$

where $u\equiv \,\,{\rm{E}}_{145}\left(a, b; c; x,y,z\right)$.

\bigskip

\begin{equation}
{\rm{E}}_{146}\left(a,b; c; x,y,z\right)=\sum\limits_{m,n,p=0}^\infty\frac{(a)_{2m+p}(b)_{2n-p}}{(c)_{m+n}}\frac{x^m}{m!}\frac{y^n}{n!}\frac{z^p}{p!},
\end{equation}

region of convergence:
$$ \left\{\frac{1}{r}+\frac{1}{s}>4,\,\,\,t<\infty
\right\}.
$$

System of partial differential equations:

$
\left\{
\begin{aligned}
&{x(1-4x)u_{xx} -z^{2} u_{zz} +yu_{xy} -4xzu_{xz}}{ +\left[c-(4a +6)x\right]u_{x} -2(a+1)zu_{z} -a(a+1)u=0,} \\&
{{y(1-4y)u_{yy} -z^{2} u_{zz} +xu_{xy} +4yzu_{yz}}} { +\left[c-(4b +6)y\right]u_{y} +2bzu_{z} -b(b+1)u=0,}
\\& {z u_{zz} -2yu_{yz} +2xu_{x}+(1-b+z)u_{z} +au=0,}
\end{aligned}
\right.
$

where $u\equiv \,\,{\rm{E}}_{146}\left(a,b; c; x,y,z\right)$.

\bigskip

\begin{equation}
{\rm{E}}_{147}\left(a,b; c; x,y,z\right)=\sum\limits_{m,n,p=0}^\infty\frac{(a)_{n+2p}(b)_{2m-p}}{(c)_{m+n}}\frac{x^m}{m!}\frac{y^n}{n!}\frac{z^p}{p!},
\end{equation}

region of convergence:
$$ \left\{ \frac{1}{r}+\frac{1}{s}>4,\,\,\,t<\infty
\right\}.
$$

System of partial differential equations:

$
\left\{
\begin{aligned}
&{x(1-4x)u_{xx} -z^{2} u_{zz} +yu_{xy} +4xzu_{xz}}+\left[c-(4b+6)x\right]u_{x} +2bzu_{z} -b(b+1)u=0, \\
&{yu_{yy}+x u_{xy}}{ +\left(c- y\right)u_{y} -2zu_z-a u=0,}\\
& { z(1+4z)u_{zz} +y^{2} u_{yy}-2xu_{xz} +4yzu_{yz}}{+2(a+1)yu_{y} +\left[1-b+(4a+6)z\right]u_{z} +a(a+1)u=0,}
\end{aligned}
\right.
$

where $u\equiv \,\,{\rm{E}}_{147}\left(a,b; c; x,y,z\right)$.

\bigskip

\begin{equation}
{\rm{E}}_{148}\left(a_1,a_2, a_3, b; c; x,y,z\right)=\sum\limits_{m,n,p=0}^\infty\frac{(a_1)_{n}(a_2)_n(a_3)_p(b)_{2m-p}  }{(c)_{m+n}}\frac{x^m}{m!}\frac{y^n}{n!}\frac{z^p}{p!},
\end{equation}

region of convergence:
$$ \left\{ \frac{1}{4r}+\frac{1}{s}>4,\,\,\,t<\infty
\right\}.
$$

System of partial differential equations:

$
\left\{
\begin{aligned}
&{x(1-4x)u_{xx} -z^{2} u_{zz} +yu_{xy} +4xzu_{xz}}{ +\left[c-(4b +6)x\right]u_{x} +2bzu_{z} -b(b+1)u=0,} \\&
{y(1-y)u_{yy}+x u_{xy}}{ +\left[c-(a_1+a_2 +1)y\right]u_{y} -a_1a_2 u=0,}
\\& {z u_{zz} -2xu_{xz} +(1-b+z)u_{z} +a_3u=0,}
\end{aligned}
\right.
$

where $u\equiv \,\,{\rm{E}}_{148}\left(a_1,a_2, a_3, b; c; x,y,z\right)$.

\bigskip

\begin{equation}
{\rm{E}}_{149}\left(a_1,a_2,a_3, b; c; x,y,z\right)=\sum\limits_{m,n,p=0}^\infty\frac{(a_1)_{n}(a_2)_{p}(a_3)_p(b)_{2m-p}  }{(c)_{m+n}}\frac{x^m}{m!}\frac{y^n}{n!}\frac{z^p}{p!},
\end{equation}

region of convergence:
$$ \left\{ \left(1+2\sqrt{r}\right)t<1,\,\,\,s<\infty
\right\}.
$$

System of partial differential equations:

$
\left\{
\begin{aligned}
&{x(1-4x)u_{xx} -z^{2} u_{zz} +yu_{xy} +4xzu_{xz}}{ +\left[c-(4b +6)x\right]u_{x} +2bzu_{z} -b(b+1)u=0,} \\
&{yu_{yy}+x u_{xy}}{ +\left(c- y\right)u_{y} -a_1 u=0,}\\
& {z(1+z) u_{zz} -2xu_{xz} +\left[1-b+\left(a_2+a_3+1\right)z\right]u_{z} +a_2a_3u=0,}
\end{aligned}
\right.
$

where $u\equiv \,\,{\rm{E}}_{149}\left(a_1,a_2,a_3, b; c; x,y,z\right)$.

\bigskip

\begin{equation}
{\rm{E}}_{150}\left(a_1,a_2, b; c; x,y,z\right)=\sum\limits_{m,n,p=0}^\infty\frac{(a_1)_{p}(a_2)_{p}(b)_{2m-p}  }{(c)_{m+n}}\frac{x^m}{m!}\frac{y^n}{n!}\frac{z^p}{p!},
\end{equation}

region of convergence:
$$ \left\{ \left(1+2\sqrt{r}\right)t<1,\,\,\,s<\infty
\right\}.
$$

System of partial differential equations:

$
\left\{
\begin{aligned}
&{x(1-4x)u_{xx} -z^{2} u_{zz} +yu_{xy} +4xzu_{xz}}{ +\left[c-(4b +6)x\right]u_{x} +2bzu_{z} -b(b+1)u=0,} \\&
{yu_{yy}+x u_{xy}}{ +cu_{y} - u=0,}
\\& {z(1+z) u_{zz} -2xu_{xz} +\left[1-b+\left(a_1+a_2+1\right)z\right]u_{z} +a_1a_2u=0,}
\end{aligned}
\right.
$

where $u\equiv \,\,{\rm{E}}_{150}\left(a_1,a_2, b; c; x,y,z\right)$.

\bigskip

\begin{equation}
{\rm{E}}_{151}\left(a_1, a_2, b; c; x,y,z\right)=\sum\limits_{m,n,p=0}^\infty\frac{(a_1)_n(a_2)_n(b)_{2m-p}  }{(c)_{m+n}}\frac{x^m}{m!}\frac{y^n}{n!}\frac{z^p}{p!},
\end{equation}

region of convergence:
$$ \left\{ \frac{1}{4r}+\frac{1}{s}>4,\,\,\,t<\infty
\right\}.
$$

System of partial differential equations:

$
\left\{
\begin{aligned}
&{x(1-4x)u_{xx} -z^{2} u_{zz} +yu_{xy} +4xzu_{xz}}{ +\left[c-(4b +6)x\right]u_{x} +2bzu_{z} -b(b+1)u=0,} \\&
{y(1-y)u_{yy}+x u_{xy}}{ +\left[c-\left(a_1+a_2+1\right)y\right]u_{y} - a_1a_2 u=0,}
\\& {z u_{zz} -2xu_{xz} +(1-b)u_{z} +u=0,}
\end{aligned}
\right.
$

where $u\equiv \,\,{\rm{E}}_{151}\left(a_1, a_2, b; c; x,y,z\right)$.

\bigskip

\begin{equation}
{\rm{E}}_{152}\left(a, b; c; x,y,z\right)=\sum\limits_{m,n,p=0}^\infty\frac{(a)_{n+p}(b)_{2m-p}}{(c)_{m+n}}\frac{x^m}{m!}\frac{y^n}{n!}\frac{z^p}{p!},
\end{equation}

region of convergence:
$$ \left\{ r<\frac{1}{4},\,\,\,s<\infty,\,\,\,t<\infty
\right\}.
$$

System of partial differential equations:

$
\left\{
\begin{aligned}
&x(1-4x)u_{xx} -z^{2} u_{zz} +yu_{xy} +4xzu_{xz} +\left[c-(4b+6)x\right]u_{x} +2bzu_{z} -b(b+1)u=0, \\&
{{yu_{yy} +xu_{xy} +(c-y)u_{y}-zu_{z} -au}=0,} \\ & {zu_{zz} -2xu_{xz} +yu_{y} +(1-b+z)u_{z} +au=0,}
\end{aligned}
\right.
$

where $u\equiv \,\,{\rm{E}}_{152}\left(a, b; c; x,y,z\right)$.

\bigskip

\begin{equation}
{\rm{E}}_{153}\left(a_1,a_2,a_3,b; c_1,c_2;;x,y,z\right)=
\sum\limits_{m,n,p=0}^\infty\frac{(a_1)_m(a_2)_n(a_3)_p(b)_{m+n-p}}{(c_1)_m(c_2)_n}\frac{x^m}{m!}\frac{y^n}{n!}\frac{z^p}{p!},,
\end{equation}

first appearance of this function in the literature: [15],\, $A_1$,

region of convergence:
$$ \left\{ r+s<1
,\,\,\,t<\infty
\right\}.
$$

System of partial differential equations:

$
\left\{
\begin{aligned}
&{x(1-x)u_{xx} -xyu_{xy} +xzu_{xz} +\left[c_{1} -\left(a_{1} +b+1\right)x\right]u_{x}}{ -a_{1} yu_{y} +a_{1} zu_{z} -a_{1} bu=0,} \\ &{y(1-y)u_{yy} -xyu_{xy} +yzu_{yz} -a_{2} xu_{x}}{+\left[c_2-\left(a_{2} +b+1\right)y\right]u_{y} +a_{2} zu_{z} -a_{2} bu=0,} \\& {zu_{zz} -xu_{xz} -yu_{yz} +(1-b+z)u_{z} +a_{3} u=0,}
\end{aligned}
\right.
$

where $u\equiv \,\,{\rm{E}}_{153}\left(a_1,a_2,a_3,b; c_1,c_2;;x,y,z\right)$.

Particular solutions:

$
{u_1} ={\rm{E}}_{153}\left(a_1,a_2,a_3,b; c_1,c_2;;x,y,z\right),
$

$
{u_2} = {x^{1 - c_1}}{\rm{E}}_{153}\left(1-c_1+a_1,a_2,a_3,1-c_1+b; 2-c_1,c_2;;x,y,z\right),
$

$
{u_3} = {y^{1 - c_2}}{\rm{E}}_{153}\left(a_1,1-c_2+a_2,a_3,1-c_2+b; c_1,2-c_2;;x,y,z\right),
$

$
{u_4} = {x^{1 - c_1}}{y^{1 - c_2}}{\rm{E}}_{153}\left(1-c_1+a_1,1-c_2+a_2,a_3,2-c_1-c_2+b; 2-c_1,2-c_2;;x,y,z\right).
$

\bigskip

\begin{equation}
{\rm{E}}_{154}\left(a_1,a_2,a_3,b; c_1,c_2;;x,y,z\right)=
\sum\limits_{m,n,p=0}^\infty\frac{(a_1)_m(a_2)_p(a_3)_p(b)_{m+n-p}}{(c_1)_m(c_2)_n}\frac{x^m}{m!}\frac{y^n}{n!}\frac{z^p}{p!},
\end{equation}

region of convergence:
$$ \left\{ (1+r)t<1,\,\,\,s<\infty
\right\}.
$$

System of partial differential equations:

$
\left\{
\begin{aligned}
&{x(1-x)u_{xx} -xyu_{xy} +xzu_{xz} +\left[c_{1} -\left(a_{1} +b+1\right)x\right]u_{x}}{ -a_{1} yu_{y} +a_{1} zu_{z} -a_{1} bu=0,} \\& {yu_{yy} -xu_{x} +\left(c_{2} -y\right)u_{y} +zu_{z} -bu=0,} \\
& z(1+z)u_{zz} -xu_{xz} -yu_{yz}+\left[1-b+\left(a_{2} +a_{3} +1\right)z\right]u_{z} +a_{2} a_{3} u=0,
\end{aligned}
\right.
$

where $u\equiv \,\,{\rm{E}}_{154}\left(a_1,a_2,a_3,b; c_1,c_2;;x,y,z\right)$.

Particular solutions:

$
{u_1} ={\rm{E}}_{154}\left(a_1,a_2,a_3,b; c_1,c_2;;x,y,z\right) ,
$

$
{u_2} = {x^{1 - c_1}}{\rm{E}}_{154}\left(1-c_1+a_1,a_2,a_3,1-c_1+b; 2-c_1,c_2;;x,y,z\right),
$

$
{u_3} = {y^{1 - c_2}}{\rm{E}}_{154}\left(a_1,a_2,a_3,1-c_2+b; c_1,2-c_2;;x,y,z\right),
$

$
{u_4} = {x^{1 - c_1}}{y^{1 - c_2}}{\rm{E}}_{154}\left(1-c_1+a_1,a_2,a_3,2-c_1-c_2+b; 2-c_1,2-c_2;;x,y,z\right).
$

\bigskip

\begin{equation}
{\rm{E}}_{155}\left(a_1,a_2,b; c_1,c_2;x,y,z\right)=
\sum\limits_{m,n,p=0}^\infty\frac{(a_1)_m(a_2)_n(b)_{m+n-p}}{(c_1)_m(c_2)_n}\frac{x^m}{m!}\frac{y^n}{n!}\frac{z^p}{p!},
\end{equation}

first appearance of this function in the literature: [15],\, $A_2$,

region of convergence:
$$ \left\{ r+s<1
\,\,
\,\,t<\infty
\right\}.
$$

System of partial differential equations:

$
\left\{
\begin{aligned}
&{x(1-x)u_{xx} -xyu_{xy} +xzu_{xz} +\left[c_{1} -\left(a_{1} +b+1\right)x\right]u_{x}}{ -a_{1} yu_{y} +a_{1} zu_{z} -a_{1} bu=0,} \\ &{y(1-y)u_{yy} -xyu_{xy} +yzu_{yz} -a_{2} xu_{x}}{+\left[c_2-\left(a_{2} +b+1\right)y\right]u_{y} +a_{2} zu_{z} -a_{2} bu=0,} \\& {zu_{zz} -xu_{xz} -yu_{yz} +(1-b)u_{z} +u=0,}
\end{aligned}
\right.
$

where $u\equiv \,\,{\rm{E}}_{155}\left(a_1,a_2,b; c_1,c_2;x,y,z\right)$.

Particular solutions:

$
{u_1} ={\rm{E}}_{155}\left(a_1,a_2,b; c_1,c_2;x,y,z\right),
$

$
{u_2} = {x^{1 - c_1}}{\rm{E}}_{155}\left(1-c_1+a_1,a_2,1-c_1+b; 2-c_1,c_2;x,y,z\right),
$

$
{u_3} = {y^{1 - c_2}}{\rm{E}}_{155}\left(a_1,1-c_2+a_2,1-c_2+b; c_1,2-c_2;x,y,z\right),
$

$
{u_4} = {x^{1 - c_1}}{y^{1 - c_2}}{\rm{E}}_{155}\left(1-c_1+a_1,1-c_2+a_2,2-c_1-c_2+b; 2-c_1,2-c_2;x,y,z\right).
$

\bigskip

\begin{equation}
{\rm{E}}_{156}\left(a_1,a_2,b; c_1,c_2;x,y,z\right)=
\sum\limits_{m,n,p=0}^\infty\frac{(a_1)_m(a_2)_p(b)_{m+n-p}}{(c_1)_m(c_2)_n}\frac{x^m}{m!}\frac{y^n}{n!}\frac{z^p}{p!},
\end{equation}

region of convergence:
$$ \left\{ r<1,\,\,\,s<\infty,\,\,\,t<\infty
\right\}.
$$

System of partial differential equations:

$
\left\{
\begin{aligned}
&{x(1-x)u_{xx} -xyu_{xy} +xzu_{xz} +\left[c_{1} -\left(a_{1} +b+1\right)x\right]u_{x}}{ -a_{1} yu_{y} +a_{1} zu_{z} -a_{1} bu=0,} \\& {yu_{yy} -xu_{x} +\left(c_{2} -y\right)u_{y} +zu_{z} -bu=0,} \\& {zu_{zz} -xu_{xz} -yu_{yz} +(1-b+z)u_{z} +a_{2} u=0,}
\end{aligned}
\right.
$

where $u\equiv \,\,{\rm{E}}_{156}\left(a_1,a_2,b; c_1,c_2;x,y,z\right)$.

Particular solutions:

$
{u_1} ={\rm{E}}_{156}\left(a_1,a_2,b; c_1,c_2;x,y,z\right),
$

$
{u_2} = {x^{1 - c_1}}{\rm{E}}_{156}\left(1-c_1+a_1,a_2,1-c_1+b; 2-c_1,c_2;x,y,z\right),
$

$
{u_3} = {y^{1 - c_2}}{\rm{E}}_{156}\left(a_1,a_2,1-c_2+b; c_1,2-c_2;x,y,z\right),
$

$
{u_4} = {x^{1 - c_1}}{y^{1 - c_2}}{\rm{E}}_{156}\left(1-c_1+a_1,1-c_2+a_2,2-c_1-c_2+b; 2-c_1,2-c_2;x,y,z\right).
$

\bigskip

\begin{equation}
{\rm{E}}_{157}\left(a,b; c_1,c_2;x,y,z\right)=
\sum\limits_{m,n,p=0}^\infty\frac{(a)_m(b)_{m+n-p}}{(c_1)_m(c_2)_n}\frac{x^m}{m!}\frac{y^n}{n!}\frac{z^p}{p!},
\end{equation}

region of convergence:
$$ \left\{ r<1,\,\,\,s<\infty,\,\,\,t<\infty
\right\}.
$$

System of partial differential equations:

$
\left\{
\begin{aligned}
&{x(1-x)u_{xx} -xyu_{xy} +xzu_{xz} +\left[c_{1} -\left(a +b+1\right)x\right]u_{x}}{ -a yu_{y} +azu_{z} -abu=0,} \\& {yu_{yy} -xu_{x} +\left(c_{2} -y\right)u_{y} +zu_{z} -bu=0,} \\& {zu_{zz} -xu_{xz} -yu_{yz} +(1-b)u_{z} +u=0,}
\end{aligned}
\right.
$

where $u\equiv \,\,{\rm{E}}_{157}\left(a,b; c_1,c_2;x,y,z\right)$.

Particular solutions:

$
{u_1} ={\rm{E}}_{157}\left(a,b; c_1,c_2;x,y,z\right),
$

$
{u_2} = {x^{1 - c_1}}{\rm{E}}_{157}\left(1-c_1+a,1-c_1+b; 2-c_1,c_2;x,y,z\right),
$

$
{u_3} = {y^{1 - c_2}}{\rm{E}}_{157}\left(a,1-c_2+b; c_1,2-c_2;x,y,z\right),
$

$
{u_4} = {x^{1 - c_1}}{y^{1 - c_2}}{\rm{E}}_{157}\left(1-c_1+a,2-c_1-c_2+b; 2-c_1,2-c_2;x,y,z\right).
$

\bigskip

\begin{equation}
{\rm{E}}_{158}\left(a_1,a_2,b; c_1,c_2;x,y,z\right)=
\sum\limits_{m,n,p=0}^\infty\frac{(a_1)_p(a_2)_p(b)_{m+n-p}}{(c_1)_m(c_2)_n}\frac{x^m}{m!}\frac{y^n}{n!}\frac{z^p}{p!},
\end{equation}

region of convergence:
$$ \left\{ r<\infty,\,\,\,s<\infty,\,\,\,t<1
\right\}.
$$

System of partial differential equations:

$
\left\{
\begin{aligned}
&{xu_{xx} +(c_{1} -x)u_{x} -yu_{y} +zu_{z} -bu=0,} \\& {yu_{yy} -xu_{x} +\left(c_{2} -y\right)u_{y} +zu_{z} -bu=0,} \\& {z(1+z)u_{zz} -xu_{xz} -yu_{yz}}{ +\left[1-b+\left(a_{1} +a_{2} +1\right)z\right]u_{z} +a_{1} a_{2} u=0,}
\end{aligned}
\right.
$

where $u\equiv \,\,{\rm{E}}_{158}\left(a_1,a_2,b; c_1,c_2;x,y,z\right)$.

Particular solutions:

$
{u_1} ={\rm{E}}_{158}\left(a_1,a_2,b; c_1,c_2;x,y,z\right),
$

$
{u_2} = {x^{1 - c_1}}{\rm{E}}_{158}\left(a_1,a_2,1-c_1+b; 2-c_1,c_2;x,y,z\right),
$

$
{u_3} = {y^{1 - c_2}}{\rm{E}}_{158}\left(a_1,a_2,1-c_2+b; c_1,2-c_2;x,y,z\right),
$

$
{u_4} = {x^{1 - c_1}}{y^{1 - c_2}}{\rm{E}}_{158}\left(a_1,a_2,2-c_1-c_2+b; 2-c_1,2-c_2;x,y,z\right).
$

\bigskip

\begin{equation}
{\rm{E}}_{159}\left(a,b; c_1,c_2;x,y,z\right)=
\sum\limits_{m,n,p=0}^\infty\frac{(a)_p(b)_{m+n-p}}{(c_1)_m(c_2)_n}\frac{x^m}{m!}\frac{y^n}{n!}\frac{z^p}{p!},
\end{equation}

region of convergence:
$$ \left\{ r<\infty,\,\,\,s<\infty,\,\,\,t<\infty
\right\}.
$$

System of partial differential equations:

$
\left\{
\begin{aligned}
&{xu_{xx} +(c_{1} -x)u_{x} -yu_{y} +zu_{z} -bu=0,} \\& {yu_{yy} -xu_{x} +\left(c_{2} -y\right)u_{y} +zu_{z} -bu=0,} \\& {zu_{zz} -xu_{xz} -yu_{yz} +(1-b+z)u_{z} +a u=0,}
\end{aligned}
\right.
$

where $u\equiv \,\,{\rm{E}}_{159}\left(a,b; c_1,c_2;x,y,z\right)$.

Particular solutions:

$
{u_1} ={\rm{E}}_{159}\left(a,b; c_1,c_2;x,y,z\right),
$

$
{u_2} = {x^{1 - c_1}}{\rm{E}}_{159}\left(a,1-c_1+b; 2-c_1,c_2;x,y,z\right),
$

$
{u_3} = {y^{1 - c_2}}{\rm{E}}_{159}\left(a,1-c_2+b; c_1,2-c_2;x,y,z\right),
$

$
{u_4} = {x^{1 - c_1}}{y^{1 - c_2}}{\rm{E}}_{159}\left(a,2-c_1-c_2+b; 2-c_1,2-c_2;x,y,z\right).
$

\bigskip

\begin{equation}
{\rm{E}}_{160}\left(b; c_1,c_2;x,y,z\right)=
\sum\limits_{m,n,p=0}^\infty\frac{(b)_{m+n-p}}{(c_1)_m(c_2)_n}\frac{x^m}{m!}\frac{y^n}{n!}\frac{z^p}{p!},
\end{equation}

region of convergence:
$$ \left\{ r<\infty,\,\,\,s<\infty,\,\,\,t<\infty
\right\}.
$$

System of partial differential equations:

$
\left\{
\begin{aligned}
&{xu_{xx} +(c_{1} -x)u_{x} -yu_{y} +zu_{z} -bu=0,} \\& {yu_{yy} -xu_{x} +\left(c_{2} -y\right)u_{y} +zu_{z} -bu=0,} \\& {zu_{zz} -xu_{xz} -yu_{yz} +(1-b)u_{z} +u=0,}
\end{aligned}
\right.
$

where $u\equiv \,\,{\rm{E}}_{160}\left(b; c_1,c_2;x,y,z\right)$.

Particular solutions:

$
{u_1} ={\rm{E}}_{160}\left(b; c_1,c_2;x,y,z\right),
$

$
{u_2} = {x^{1 - c_1}}{\rm{E}}_{160}\left(1-c_1+b; 2-c_1,c_2;x,y,z\right),
$

$
{u_3} = {y^{1 - c_2}}{\rm{E}}_{160}\left(1-c_2+b; c_1,2-c_2;x,y,z\right),
$

$
{u_4} = {x^{1 - c_1}}{y^{1 - c_2}}{\rm{E}}_{160}\left(2-c_1-c_2+b; 2-c_1,2-c_2;x,y,z\right).
$

\bigskip

\begin{equation}
{\rm{E}}_{161}\left(a_1,a_2,a_3,b; c_1,c_2;x,y,z\right)=\sum\limits_{m,n,p=0}^\infty\frac{(a_1)_{m+n}(a_2)_n(a_3)_p(b)_{m-p}}{(c_1)_m(c_2)_n}\frac{x^m}{m!}\frac{y^n}{n!}\frac{z^p}{p!},
\end{equation}

region of convergence:
$$ \left\{ r+s<1,\,\,\,t<\infty
\right\}.
$$

System of partial differential equations:

$
\left\{
\begin{aligned}
&{x(1-x)u_{xx} -xyu_{xy} +xzu_{xz} +yzu_{yz}} {+\left[c_{1} -\left(a_{1} +b+1\right)x\right]u_{x} -byu_{y} +a_{1} zu_{z} -a_{1} bu=0,} \\& {y(1-y)u_{yy} -xyu_{xy} -a_{2} xu_{x} +\left[c_{2} -\left(a_{1} +a_{2} +1\right)y\right]u_{y} -a_{1} a_{2} u=0,} \\& {zu_{zz} -xu_{xz} +(1-b+z)u_{z} +a_{3} u=0,}
\end{aligned}
\right.
$

where $u\equiv \,\,{\rm{E}}_{161}\left(a_1,a_2,a_3,b; c_1,c_2;x,y,z\right)$.

Particular solutions:

$
{u_1} ={\rm{E}}_{161}\left(a_1,a_2,a_3,b; c_1,c_2;x,y,z\right),
$

$
{u_2} = {x^{1 - c_1}}{\rm{E}}_{161}\left(1-c_1+a_1,a_2,a_3,1-c_1+b; 2-c_1,c_2;x,y,z\right),
$

$
{u_3} = {y^{1 - c_2}}{\rm{E}}_{161}\left(1-c_2+a_1,1-c_2+a_2,a_3,b; c_1,2-c_2;x,y,z\right),
$

$
{u_4} = {x^{1 - c_1}}{y^{1 - c_2}}{\rm{E}}_{161}\left(2-c_1-c_2+a_1,1-c_2+a_2,a_3,1-c_1+b; 2-c_1,2-c_2;x,y,z\right).
$

\bigskip

\begin{equation}
{\rm{E}}_{162}\left(a_1,a_2,a_3,b; c_1,c_2;x,y,z\right)=\sum\limits_{m,n,p=0}^\infty\frac{(a_1)_{m+n}(a_2)_p(a_3)_p(b)_{m-p}}{(c_1)_m(c_2)_n}\frac{x^m}{m!}\frac{y^n}{n!}\frac{z^p}{p!},
\end{equation}

region of convergence:
$$ \left\{ (1+r)t<1,\,\,\,s<\infty
\right\}.
$$

System of partial differential equations:

$
\left\{
\begin{aligned}
&{x(1-x)u_{xx} -xyu_{xy} +xzu_{xz} +yzu_{yz}} {+\left[c_{1} -\left(a_{1} +b+1\right)x\right]u_{x} -byu_{y} +a_{1} zu_{z} -a_{1} bu=0,} \\& {yu_{yy} -xu_{x} +\left(c_{2} -y\right)u_{y} -a_{1} u=0,} \\& {z(1+z)u_{zz} -xu_{xz} +\left[1-b+\left(a_{2} +a_{3} +1\right)z\right]u_{z} +a_{2} a_{3} u=0,}
\end{aligned}
\right.
$

where $u\equiv \,\,{\rm{E}}_{162}\left(a_1,a_2,a_3,b; c_1,c_2;x,y,z\right)$.

Particular solutions:

$
{u_1} ={\rm{E}}_{162}\left(a_1,a_2,a_3,b; c_1,c_2;x,y,z\right),
$

$
{u_2} = {x^{1 - c_1}}{\rm{E}}_{162}\left(1-c_1+a_1,a_2,a_3,1-c_1+b; 2-c_1,c_2;x,y,z\right),
$

$
{u_3} = {y^{1 - c_2}}{\rm{E}}_{162}\left(1-c_2+a_1,a_2,a_3,b; c_1,2-c_2;x,y,z\right),
$

$
{u_4} = {x^{1 - c_1}}{y^{1 - c_2}}{\rm{E}}_{162}\left(2-c_1-c_2+a_1,a_2,a_3,1-c_1+b; 2-c_1,2-c_2;x,y,z\right).
$

\bigskip

\begin{equation}
{\rm{E}}_{163}\left(a_1,a_2,b; c_1,c_2;x,y,z\right)=\sum\limits_{m,n,p=0}^\infty\frac{(a_1)_{m+n}(a_2)_n(b)_{m-p}}{(c_1)_m(c_2)_n}\frac{x^m}{m!}\frac{y^n}{n!}\frac{z^p}{p!},
\end{equation}

region of convergence:
$$ \left\{ r+s<1,\,\,\,t<\infty
\right\}.
$$

System of partial differential equations:

$
\left\{
\begin{aligned}
&{x(1-x)u_{xx} -xyu_{xy} +xzu_{xz} +yzu_{yz}}{+\left[c_{1} -\left(a_{1} +b+1\right)x\right]u_{x} -byu_{y} +a_{1} zu_{z} -a_{1} bu=0,} \\& {y(1-y)u_{yy} -xyu_{xy} -a_{2} xu_{x} +\left[c_{2} -\left(a_{1} +a_{2} +1\right)y\right]u_{y} -a_{1} a_{2} u=0,}  \\& {zu_{zz} -xu_{xz} +(1-b)u_{z} + u=0,}
\end{aligned}
\right.
$

where $u\equiv \,\,{\rm{E}}_{163}\left(a_1,a_2,b; c_1,c_2;x,y,z\right)$.

Particular solutions:

$
{u_1} ={\rm{E}}_{163}\left(a_1,a_2,b; c_1,c_2;x,y,z\right),
$

$
{u_2} = {x^{1 - c_1}}{\rm{E}}_{163}\left(1-c_1+a_1,a_2,1-c_1+b; 2-c_1,c_2;x,y,z\right),
$

$
{u_3} = {y^{1 - c_2}}{\rm{E}}_{163}\left(1-c_2+a_1,1-c_2+a_2,b; c_1,2-c_2;x,y,z\right),
$

$
{u_4} = {x^{1 - c_1}}{y^{1 - c_2}}{\rm{E}}_{163}\left(2-c_1-c_2+a_1,1-c_2+a_2,1-c_1+b; 2-c_1,2-c_2;x,y,z\right).
$

\bigskip

\begin{equation}
{\rm{E}}_{164}\left(a_1,a_2,b; c_1,c_2;x,y,z\right)=\sum\limits_{m,n,p=0}^\infty\frac{(a_1)_{m+n}(a_2)_p(b)_{m-p}}{(c_1)_m(c_2)_n}\frac{x^m}{m!}\frac{y^n}{n!}\frac{z^p}{p!},
\end{equation}

region of convergence:
$$ \left\{ r<1,\,\,\,s<\infty,\,\,\,t<\infty
\right\}.
$$

System of partial differential equations:

$
\left\{
\begin{aligned}
&{x(1-x)u_{xx} -xyu_{xy} +xzu_{xz} +yzu_{yz}}{+\left[c_{1} -\left(a_{1} +b+1\right)x\right]u_{x} -byu_{y} +a_{1} zu_{z} -a_{1} bu=0,} \\& {yu_{yy} -xu_{x} +\left(c_{2} -y\right)u_{y} -a_{1} u=0,}  \\& {zu_{zz} -xu_{xz} +(1-b+z)u_{z} +a_{2} u=0,}
\end{aligned}
\right.
$

where $u\equiv \,\,{\rm{E}}_{164}\left(a_1,a_2,b; c_1,c_2;x,y,z\right)$.

Particular solutions:

$
{u_1} ={\rm{E}}_{164}\left(a_1,a_2,b; c_1,c_2;x,y,z\right),
$

$
{u_2} = {x^{1 - c_1}}{\rm{E}}_{164}\left(1-c_1+a_1,a_2,1-c_1+b; 2-c_1,c_2;x,y,z\right),
$

$
{u_3} = {y^{1 - c_2}}{\rm{E}}_{164}\left(1-c_2+a_1,a_2,b; c_1,2-c_2;x,y,z\right),
$

$
{u_4} = {x^{1 - c_1}}{y^{1 - c_2}}{\rm{E}}_{164}\left(2-c_1-c_2+a_1,a_2,1-c_1+b; 2-c_1,2-c_2;x,y,z\right).
$

\bigskip

\begin{equation}
{\rm{E}}_{165}\left(a,b; c_1,c_2;x,y,z\right)=\sum\limits_{m,n,p=0}^\infty\frac{(a)_{m+n}(b)_{m-p}}{(c_1)_m(c_2)_n}\frac{x^m}{m!}\frac{y^n}{n!}\frac{z^p}{p!},
\end{equation}

region of convergence:
$$ \left\{ r<1,\,\,\,s<\infty,\,\,\,t<\infty
\right\}.
$$

System of partial differential equations:

$
\left\{
\begin{aligned}
&{x(1-x)u_{xx} -xyu_{xy} +xzu_{xz} +yzu_{yz}}{+\left[c_{1} -\left(a+b+1\right)x\right]u_{x} -byu_{y} +a zu_{z} -a bu=0,} \\& {yu_{yy} -xu_{x} +\left(c_{2} -y\right)u_{y} -a u=0,}  \\& {zu_{zz} -xu_{xz} +(1-b)u_{z} + u=0,}
\end{aligned}
\right.
$

where $u\equiv \,\,{\rm{E}}_{165}\left(a,b; c_1,c_2;x,y,z\right)$.

Particular solutions:

$
{u_1} ={\rm{E}}_{165}\left(a,b; c_1,c_2;x,y,z\right),
$

$
{u_2} = {x^{1 - c_1}}{\rm{E}}_{165}\left(1-c_1+a,1-c_1+b; 2-c_1,c_2;x,y,z\right),
$

$
{u_3} = {y^{1 - c_2}}{\rm{E}}_{165}\left(1-c_2+a,b; c_1,2-c_2;x,y,z\right),
$

$
{u_4} = {x^{1 - c_1}}{y^{1 - c_2}}{\rm{E}}_{165}\left(2-c_1-c_2+a,1-c_1+b; 2-c_1,2-c_2;x,y,z\right).
$

\bigskip

\begin{equation}
{\rm{E}}_{166}\left(a_1,a_2,a_3,b; c_1,c_2;x,y,z\right)=\sum\limits_{m,n,p=0}^\infty\frac{(a_1)_{n+p}(a_2)_m(a_3)_n(b)_{m-p}}{(c_1)_m(c_2)_n}\frac{x^m}{m!}\frac{y^n}{n!}\frac{z^p}{p!},
\end{equation}

region of convergence:
$$ \left\{ r<1,\,\,\,\, s<1,\,\,\,t<\infty
\right\}.
$$

System of partial differential equations:

$
\left\{
\begin{aligned}
&{x(1-x)u_{xx} +xzu_{xz} +\left[c_{1} -\left(a_{2} +b+1\right)x\right]u_{x} +a_{2} zu_{z} -a_{2} bu=0,} \\& {y(1-y)u_{yy} -yzu_{yz} +\left[c_{2} -\left(a_{1} +a_{3} +1\right)y\right]u_{y} -a_{3} zu_{z} -a_{1} a_{3} u=0,} \\& {zu_{zz} -xu_{xz}+yu_{y} +(1-b+z)u_{z}  +a_{1} u=0,}
\end{aligned}
\right.
$

where $u\equiv \,\,{\rm{E}}_{166}\left(a_1,a_2,a_3,b; c_1,c_2;x,y,z\right)$.

Particular solutions:

$
{u_1} ={\rm{E}}_{166}\left(a_1,a_2,a_3,b; c_1,c_2;x,y,z\right),
$

$
{u_2} = {x^{1 - c_1}}{\rm{E}}_{166}\left(a_1,1-c_1+a_2,a_3,1-c_1+b; 2-c_1,c_2;x,y,z\right),
$

$
{u_3} = {y^{1 - c_2}}{\rm{E}}_{166}\left(1-c_2+a_1,a_2,1-c_2+a_3,b; c_1,2-c_2;x,y,z\right),
$

$
{u_4} = {x^{1 - c_1}}{y^{1 - c_2}}{\rm{E}}_{166}\left(1-c_2+a_1,1-c_1+a_2,1-c_2+a_3,1-c_1+b; 2-c_1,2-c_2;x,y,z\right).
$

\bigskip

\begin{equation}
{\rm{E}}_{167}\left(a_1,a_2,a_3,b; c_1,c_2;x,y,z\right)=\sum\limits_{m,n,p=0}^\infty\frac{(a_1)_{n+p}(a_2)_m(a_3)_p(b)_{m-p}}{(c_1)_m(c_2)_n}\frac{x^m}{m!}\frac{y^n}{n!}\frac{z^p}{p!},
\end{equation}

region of convergence:
$$ \left\{ (1+r)t<1,\,\,\,s<\infty
\right\}.
$$

System of partial differential equations:

$
\left\{
\begin{aligned}
&{x(1-x)u_{xx} +xzu_{xz} +\left[c_{1} -\left(a_{2} +b+1\right)x\right]u_{x} +a_{2} zu_{z} -a_{2} bu=0,} \\ &{yu_{yy} +(c_{2} -y)u_{y} -zu_{z} -a_{1} u=0,} \\ & {z(1+z)u_{zz} -xu_{xz} +yzu_{yz} +a_{3} yu_{y}}{+\left[1-b+\left(a_{1} +a_{3} +1\right)z\right]u_{z} +a_{1} a_{3} u=0,}
\end{aligned}
\right.
$

where $u\equiv \,\,{\rm{E}}_{167}\left(a_1,a_2,a_3,b; c_1,c_2;x,y,z\right)$.

Particular solutions:

$
{u_1} ={\rm{E}}_{167}\left(a_1,a_2,a_3,b; c_1,c_2;x,y,z\right),
$

$
{u_2} = {x^{1 - c_1}}{\rm{E}}_{167}\left(a_1,1-c_1+a_2,a_3,1-c_1+b; 2-c_1,c_2;x,y,z\right),
$

$
{u_3} = {y^{1 - c_2}}{\rm{E}}_{167}\left(1-c_2+a_1,a_2,a_3,b; c_1,2-c_2;x,y,z\right),
$

$
{u_4} = {x^{1 - c_1}}{y^{1 - c_2}}{\rm{E}}_{167}\left(1-c_2+a_1,1-c_1+a_2,a_3,1-c_1+b; 2-c_1,2-c_2;x,y,z\right).
$

\bigskip

\begin{equation}
{\rm{E}}_{168}\left(a_1,a_2,a_3,b; c_1,c_2;x,y,z\right)=\sum\limits_{m,n,p=0}^\infty\frac{(a_1)_{n+p}(a_2)_n(a_3)_p(b)_{m-p}}{(c_1)_m(c_2)_n}\frac{x^m}{m!}\frac{y^n}{n!}\frac{z^p}{p!},
\end{equation}

region of convergence:
$$ \left\{ r<\infty,\,\,\,s+t<1
\right\}.
$$

System of partial differential equations:

$
\left\{
\begin{aligned}
&{xu_{xx} +\left(c_{1} -x\right)u_{x} +zu_{z} -bu=0,} \\
&{y(1-y)u_{yy} -yzu_{yz} +\left[c_{2} -\left(a_{1} +a_{2} +1\right)y\right]u_{y} }{-a_{2} zu_{z} -a_{1} a_{2} u=0,} \\ & {z(1+z)u_{zz} -xu_{xz} +yzu_{yz} +a_{3} yu_{y}}{+\left[1-b+\left(a_{1} +a_{3} +1\right)z\right]u_{z} +a_{1} a_{3} u=0,}
\end{aligned}
\right.
$

where $u\equiv \,\,{\rm{E}}_{168}\left(a_1,a_2,a_3,b; c_1,c_2;x,y,z\right)$.

Particular solutions:

$
{u_1} ={\rm{E}}_{168}\left(a_1,a_2,a_3,b; c_1,c_2;x,y,z\right),
$

$
{u_2} = {x^{1 - c_1}}{\rm{E}}_{168}\left(a_1,a_2,a_3,1-c_1+b; 2-c_1,c_2;x,y,z\right),
$

$
{u_3} = {y^{1 - c_2}}{\rm{E}}_{168}\left(1-c_2+a_1,1-c_2+a_2,a_3,b; c_1,2-c_2;x,y,z\right),
$

$
{u_4} = {x^{1 - c_1}}{y^{1 - c_2}}{\rm{E}}_{168}\left(1-c_2+a_1,1-c_2+a_2,a_3,1-c_1+b; 2-c_1,2-c_2;x,y,z\right).
$

\bigskip

\begin{equation}
{\rm{E}}_{169}\left(a_1,a_2,b; c_1,c_2;x,y,z\right)=\sum\limits_{m,n,p=0}^\infty\frac{(a_1)_{n+p}(a_2)_m(b)_{m-p}}{(c_1)_m(c_2)_n}\frac{x^m}{m!}\frac{y^n}{n!}\frac{z^p}{p!}
\end{equation}

region of convergence:
$$ \left\{ r<1,\,\,\,s<\infty,\,\,\,t<\infty
\right\}.
$$

System of partial differential equations:

$
\left\{
\begin{aligned}
&{x(1-x)u_{xx} +xzu_{xz} +\left[c_{1} -\left(a_{2} +b+1\right)x\right]u_{x} +a_{2} zu_{z} -a_{2} bu=0,} \\
&{yu_{yy} +(c_{2} -y)u_{y} -zu_{z} -a_{1} u=0,} \\
& {zu_{zz} -xu_{xz}+yu_{y} +(1-b+z)u_{z}  +a_{1} u=0,}
\end{aligned}
\right.
$

where $u\equiv \,\,{\rm{E}}_{169}\left(a_1,a_2,b; c_1,c_2;x,y,z\right)$.

Particular solutions:

$
{u_1} ={\rm{E}}_{169}\left(a_1,a_2,b; c_1,c_2;x,y,z\right),
$

$
{u_2} = {x^{1 - c_1}}{\rm{E}}_{169}\left(a_1,1-c_1+a_2,1-c_1+b; 2-c_1,c_2;x,y,z\right),
$

$
{u_3} = {y^{1 - c_2}}{\rm{E}}_{169}\left(1-c_2+a_1,a_2,b; c_1,2-c_2;x,y,z\right),
$

$
{u_4} = {x^{1 - c_1}}{y^{1 - c_2}}{\rm{E}}_{169}\left(1-c_2+a_1,1-c_1+a_2,1-c_1+b; 2-c_1,2-c_2;x,y,z\right).
$

\bigskip

\begin{equation}
{\rm{E}}_{170}\left(a_1,a_2,b; c_1,c_2;x,y,z\right)=\sum\limits_{m,n,p=0}^\infty\frac{(a_1)_{n+p}(a_2)_n(b)_{m-p}}{(c_1)_m(c_2)_n}\frac{x^m}{m!}\frac{y^n}{n!}\frac{z^p}{p!},
\end{equation}

region of convergence:
$$ \left\{ r<\infty,\,\,\,s<1,\,\,\,t<\infty
\right\}.
$$

System of partial differential equations:

$
\left\{
\begin{aligned}
&{xu_{xx} +\left(c_{1} -x\right)u_{x} +zu_{z} -bu=0,} \\
&{y(1-y)u_{yy} -yzu_{yz} +\left[c_{2} -\left(a_{1} +a_{2} +1\right)y\right]u_{y} }{-a_{2} zu_{z} -a_{1} a_{2} u=0,} \\
& {zu_{zz} -xu_{xz}+yu_{y} +(1-b+z)u_{z}  +a_{1} u=0,}
\end{aligned}
\right.
$

where $u\equiv \,\,{\rm{E}}_{170}\left(a_1,a_2,b; c_1,c_2;x,y,z\right)$.

Particular solutions:

$
{u_1} ={\rm{E}}_{170}\left(a_1,a_2,b; c_1,c_2;x,y,z\right),
$

$
{u_2} = {x^{1 - c_1}}{\rm{E}}_{170}\left(a_1,a_2,1-c_1+b; 2-c_1,c_2;x,y,z\right),
$

$
{u_3} = {y^{1 - c_2}}{\rm{E}}_{170}\left(1-c_2+a_1,1-c_2+a_2,b; c_1,2-c_2;x,y,z\right),
$

$
{u_4} = {x^{1 - c_1}}{y^{1 - c_2}}{\rm{E}}_{170}\left(1-c_2+a_1,1-c_2+a_2,1-c_1+b; 2-c_1,2-c_2;x,y,z\right).
$

\bigskip

\begin{equation}
{\rm{E}}_{171}\left(a_1,a_2,b; c_1,c_2;x,y,z\right)=\sum\limits_{m,n,p=0}^\infty\frac{(a_1)_{n+p}(a_2)_p(b)_{m-p}}{(c_1)_m(c_2)_n}\frac{x^m}{m!}\frac{y^n}{n!}\frac{z^p}{p!},
\end{equation}

region of convergence:
$$ \left\{ r<\infty,\,\,\,s<\infty,\,\,\,t<1
\right\}.
$$

System of partial differential equations:

$
\left\{
\begin{aligned}
&{xu_{xx} +\left(c_{1} -x\right)u_{x} +zu_{z} -bu=0,} \\ &{yu_{yy} +(c_{2} -y)u_{y} -zu_{z} -a_{1} u=0,} \\ & {z(1+z)u_{zz} -xu_{xz} +yzu_{yz} +a_{2} yu_{y}}{+\left[1-b+\left(a_{1} +a_{2} +1\right)z\right]u_{z} +a_{1} a_{2} u=0,}
\end{aligned}
\right.
$

where $u\equiv \,\,{\rm{E}}_{171}\left(a_1,a_2,b; c_1,c_2;x,y,z\right)$.

Particular solutions:

$
{u_1} ={\rm{E}}_{171}\left(a_1,a_2,b; c_1,c_2;x,y,z\right),
$

$
{u_2} = {x^{1 - c_1}}{\rm{E}}_{171}\left(a_1,a_2,1-c_1+b; 2-c_1,c_2;x,y,z\right),
$

$
{u_3} = {y^{1 - c_2}}{\rm{E}}_{171}\left(1-c_2+a_1,a_2,b; c_1,2-c_2;x,y,z\right),
$

$
{u_4} = {x^{1 - c_1}}{y^{1 - c_2}}{\rm{E}}_{171}\left(1-c_2+a_1,a_2,1-c_1+b; 2-c_1,2-c_2;x,y,z\right).
$

\bigskip

\begin{equation}
{\rm{E}}_{172}\left(a,b; c_1,c_2;x,y,z\right)=\sum\limits_{m,n,p=0}^\infty\frac{(a)_{n+p}(b)_{m-p}}{(c_1)_m(c_2)_n}\frac{x^m}{m!}\frac{y^n}{n!}\frac{z^p}{p!},
\end{equation}

region of convergence:
$$ \left\{ r<\infty,\,\,\,s<\infty,\,\,\,t<\infty
\right\}.
$$

System of partial differential equations:

$
\left\{
\begin{aligned}
&{xu_{xx} +\left(c_{1} -x\right)u_{x} +zu_{z} -bu=0,} \\ &{yu_{yy} +(c_{2} -y)u_{y} -zu_{z} -a u=0,} \\ & {zu_{zz} -xu_{xz} +(1-b+z)u_{z} +yu_{y} +au=0,}
\end{aligned}
\right.
$

where $u\equiv \,\,{\rm{E}}_{172}\left(a,b; c_1,c_2;x,y,z\right)$.

Particular solutions:

$
{u_1} ={\rm{E}}_{172}\left(a,b; c_1,c_2;x,y,z\right),
$

$
{u_2} = {x^{1 - c_1}}{\rm{E}}_{172}\left(a,1-c_1+b; 2-c_1,c_2;x,y,z\right),
$

$
{u_3} = {y^{1 - c_2}}{\rm{E}}_{172}\left(1-c_2+a,b; c_1,2-c_2;x,y,z\right),
$

$
{u_4} = {x^{1 - c_1}}{y^{1 - c_2}}{\rm{E}}_{172}\left(1-c_2+a,1-c_1+b; 2-c_1,2-c_2;x,y,z\right).
$

\bigskip

\begin{equation}
{\rm{E}}_{173}\left(a, b; c_1, c_2; x,y,z\right)=\sum\limits_{m,n,p=0}^\infty\frac{(a)_{m+n}(b)_{m+n-p}}{(c_1)_{m}(c_2)_{n}}\frac{x^m}{m!}\frac{y^n}{n!}\frac{z^p}{p!},
\end{equation}

region of convergence:
$$ \left\{ \sqrt{r}+\sqrt{s}<1,\,\,\,t<\infty
\right\}.
$$

System of partial differential equations:

$
\left\{
\begin{aligned}
&{x(1-x)u_{xx} -y^{2} u_{yy} -2xyu_{xy} }\\&\,\,\,\,\,\,\,\,\,\,\,\,\,+xzu_{xz}+yzu_{yz}{ +\left[c_1-(a +b+1)x\right]u_{x} -(a+b+1)yu_y+azu_{z} -abu=0,} \\&
{y(1-y)u_{yy} -x^{2} u_{xx} -2xyu_{xy} } \\&\,\,\,\,\,\,\,\,\,\,\,\,\,+xzu_{xz}+yzu_{yz} -(a+b+1)xu_x{ +\left[c_2-(a+b+1)y\right]u_{y} +azu_{z} -abu=0,}
\\& {z u_{zz} -xu_{xz}-yu_{yz}+(1-b)u_{z} +u=0,}
\end{aligned}
\right.
$

where $u\equiv \,\,{\rm{E}}_{173}\left(a, b; c_1, c_2; x,y,z\right)$.

Particular solutions:

$
{u_1} ={\rm{E}}_{173}\left(a, b; c_1, c_2; x,y,z\right),
$

$
{u_2} = {x^{1 - c_1}}{\rm{E}}_{173}\left(1-c_1+a, 1-c_1+b; 2-c_1, c_2; x,y,z\right),
$

$
{u_3} = {y^{1 - c_2}}{\rm{E}}_{173}\left(1-c_2+a, 1-c_2+b; c_1, 2-c_2; x,y,z\right),
$

$
{u_4} = {x^{1 - c_1}}{y^{1 - c_2}}{\rm{E}}_{173}\left(2-c_1-c_2+a, 2-c_1-c_2+b; 2-c_1, 2-c_2; x,y,z\right).
$

\bigskip

\begin{equation}
{\rm{E}}_{174}\left(a, b; c_1, c_2; x,y,z\right)=\sum\limits_{m,n,p=0}^\infty\frac{(a)_{n+p}(b)_{m+n-p}}{(c_1)_{m}(c_2)_{n}}\frac{x^m}{m!}\frac{y^n}{n!}\frac{z^p}{p!},
\end{equation}

region of convergence:
$$ \left\{ r<\infty,\,\,\,s<1,\,\,\,t<\infty
\right\}.
$$

System of partial differential equations:

$
\left\{
\begin{aligned}
&{xu_{xx}}{ +\left(c_1-x\right)u_{x} -yu_y+zu_{z} -bu=0,} \\&
{{y(1-y)u_{yy} +z^{2} u_{zz} -xyu_{xy} -xzu_{xz}} -axu_x}\\& \,\,\,\,\,\,\,\,\, { +\left[c_2-(a+b+1)y\right]u_{y} +(a-b+1)zu_{z} -abu=0,}
\\& {z u_{zz} -xu_{xz}-yu_{yz}+yu_y+(1-b+z)u_{z} +au=0,}
\end{aligned}
\right.
$

where $u\equiv \,\,{\rm{E}}_{174}\left(a, b; c_1, c_2; x,y,z\right)$.

Particular solutions:

$
{u_1} ={\rm{E}}_{174}\left(a, b; c_1, c_2; x,y,z\right),
$

$
{u_2} = {x^{1 - c_1}}{\rm{E}}_{174}\left(a, 1-c_1+b; 2-c_1, c_2; x,y,z\right),
$

$
{u_3} = {y^{1 - c_2}}{\rm{E}}_{174}\left(1-c_2+a, 1-c_2+b; c_1, 2-c_2; x,y,z\right),
$

$
{u_4} = {x^{1 - c_1}}{y^{1 - c_2}}{\rm{E}}_{174}\left(1-c_2+a, 2-c_1-c_2+b; 2-c_1, 2-c_2; x,y,z\right).
$

\bigskip

\begin{equation}
{\rm{E}}_{175}\left(a_1,a_2,b; c_1, c_2; x,y,z\right)=\sum\limits_{m,n,p=0}^\infty\frac{(a_1)_{n+p}(a_2)_{n+p}(b)_{m-p}}{(c_1)_{m}(c_2)_{n}}\frac{x^m}{m!}\frac{y^n}{n!}\frac{z^p}{p!},
\end{equation}

region of convergence:
$$ \left\{ r<\infty,\,\,\,\sqrt{s}+\sqrt{t}<1
\right\}.
$$

System of partial differential equations:

$
\left\{
\begin{aligned}
&{xu_{xx} +\left(c_{1} -x\right)u_{x} +zu_{z} -bu=0,} \\
&  {y(1-y)u_{yy} -z^{2} u_{zz} -2yzu_{yz} +\left[c_{2} -\left(a_{1} +a_{2} +1\right)y\right]u_{y} } {-\left(a_{1} +a_{2} +1\right)zu_{z} -a_{1} a_{2} u=0,} \\
& {z(1+z)u_{zz} +y^{2} u_{yy}-xu_{xz}  +2yzu_{yz} +\left(a_{1} +a_{2} +1\right)yu_{y}}\\& \,\,\,\,\,\,\,\,\,{+\left[1-b+\left(a_{1} +a_{2} +1\right)z\right]u_{z}  +a_{1} a_{2} u=0,}
\end{aligned}
\right.
$

where $u\equiv \,\,{\rm{E}}_{175}\left(a_1,a_2,b; c_1, c_2; x,y,z\right)$,

Particular solutions:

$
{u_1} ={\rm{E}}_{175}\left(a_1,a_2,b; c_1, c_2; x,y,z\right),
$

$
{u_2} = {x^{1 - c_1}}{\rm{E}}_{175}\left(a_1,a_2,1-c_1+b; 2-c_1, c_2; x,y,z\right),
$

$
{u_3} = {y^{1 - c_2}}{\rm{E}}_{175}\left(1-c_2+a_1,1-c_2+a_2,b; c_1, 2-c_2; x,y,z\right),
$

$
{u_4} = {x^{1 - c_1}}{y^{1 - c_2}}{\rm{E}}_{175}\left(1-c_2+a_1,1-c_2+a_2,1-c_1+b; 2-c_1, 2-c_2; x,y,z\right).
$

\bigskip

\begin{equation}
{\rm{E}}_{176}\left(a_1,a_2, b; c_1, c_2; x,y,z\right)=\sum\limits_{m,n,p=0}^\infty\frac{(a_1)_{n+p}(a_2)_{m}(b)_{m+n-p}}{(c_1)_{m}(c_2)_{n}}\frac{x^m}{m!}\frac{y^n}{n!}\frac{z^p}{p!},
\end{equation}

region of convergence:
$$ \left\{r+s<1,\,\,\,t<\infty
\right\}.
$$

System of partial differential equations:

$
\left\{
\begin{aligned}
&{x(1-x)u_{xx} -xyu_{xy} +xzu_{xz} +\left[c_{1} -\left(a_{2} +b+1\right)x\right]u_{x}}{-a_{2} yu_{y} +a_{2} zu_{z} -a_{2} bu=0,} \\
&{y(1-y)u_{yy} +z^{2} u_{zz}  }  -xyu_{xy} -xzu_{xz} -a_{1} xu_{x}\\& \,\,\,\,\,\,\,\,\,{ +\left[c_{2} -\left(a_{1} +b +1\right)y\right]u_{y}+(a_{1} -b+1)zu_{z} -a_{1} bu=0,} \\
& {zu_{zz} -xu_{xz} -yu_{yz}+yu_{y} +(1-b+z)u_{z}  +a_{1} u=0,} \end{aligned}\right.
$

where $u\equiv \,\,{\rm{E}}_{176}\left(a_1,a_2, b; c_1, c_2; x,y,z\right)$.

Particular solutions:

$
{u_1} ={\rm{E}}_{176}\left(a_1,a_2, b; c_1, c_2; x,y,z\right),
$

$
{u_2} = {x^{1 - c_1}}{\rm{E}}_{176}\left(a_1,1-c_1+a_2, 1-c_1+b; 2-c_1, c_2; x,y,z\right),
$

$
{u_3} = {y^{1 - c_2}}{\rm{E}}_{176}\left(1-c_2+a_1,a_2, 1-c_2+b; c_1, 2-c_2; x,y,z\right),
$

$
{u_4} = {x^{1 - c_1}}{y^{1 - c_2}}{\rm{E}}_{176}\left(1-c_2+a_1,1-c_1+a_2, 2-c_1-c_2+b; 2-c_1, 2-c_2; x,y,z\right).
$

\bigskip

\begin{equation}
{\rm{E}}_{177}\left(a_1,a_2, b; c_1, c_2; x,y,z\right)=\sum\limits_{m,n,p=0}^\infty\frac{(a_1)_{n+p}(a_2)_p(b)_{m+n-p}}{(c_1)_{m}(c_2)_{n}}\frac{x^m}{m!}\frac{y^n}{n!}\frac{z^p}{p!},
\end{equation}

region of convergence:
$$ \left\{ r<\infty,\,\,\,\left\{t+2\sqrt{st}<1\right\}=\left\{\sqrt{t}<\sqrt{1+s}-\sqrt{s}\right\}
\right\}.
$$

System of partial differential equations:

$
\left\{
\begin{aligned}
&{xu_{xx} +(c_{1} -x)u_{x} -yu_{y} +zu_{z} -bu=0,} \\&
{y(1-y)u_{yy} +z^{2} u_{zz} -xyu_{xy} -xzu_{xz} -a_{1} xu_{x} } \\& \,\,\,\,\,\,\,\,\, { +\left[c_{2} -\left(a_{1} +b +1\right)y\right]u_{y}+(a_{1} -b+1)zu_{z} -a_{1} bu=0,}  \\& {z(1+z)u_{zz} -xu_{xz} -y(1-z)u_{yz}+a_{2} yu_{y}}{+\left[1-b+\left(a_{1} +a_{2} +1\right)z\right]u_{z}  +a_{1} a_{2} u=0,}
\end{aligned}
\right.
$

where $u\equiv \,\,{\rm{E}}_{177}\left(a_1,a_2, b; c_1, c_2; x,y,z\right)$.

Particular solutions:

$
{u_1} ={\rm{E}}_{177}\left(a_1,a_2, b; c_1, c_2; x,y,z\right),
$

$
{u_2} = {x^{1 - c_1}}{\rm{E}}_{177}\left(a_1,a_2, 1-c_1+b; 2-c_1, c_2; x,y,z\right),
$

$
{u_3} = {y^{1 - c_2}}{\rm{E}}_{177}\left(1-c_2+a_1,a_2, 1-c_2+b; c_1, 2-c_2; x,y,z\right),
$

$
{u_4} = {x^{1 - c_1}}{y^{1 - c_2}}{\rm{E}}_{177}\left(1-c_2+a_1,a_2, 2-c_1-c_2+b; 2-c_1, 2-c_2; x,y,z\right).
$

\bigskip

\begin{equation}
{\rm{E}}_{178}\left(a_1,a_2, b; c_1, c_2; x,y,z\right)=\sum\limits_{m,n,p=0}^\infty\frac{(a_1)_{m+n}(a_2)_p(b)_{m+n-p}}{(c_1)_{m}(c_2)_{n}}\frac{x^m}{m!}\frac{y^n}{n!}\frac{z^p}{p!},
\end{equation}

region of convergence:
$$ \left\{ \sqrt{r}+\sqrt{s}<1,\,\,\,t<\infty
\right\}.
$$

System of partial differential equations:

$
\left\{
\begin{aligned}
&x(1-x)u_{xx} -y^{2} u_{yy} -2xyu_{xy} +xzu_{xz}\\&\,\,\,\,\,\,\,\,\,\,\,\, +yzu_{yz}  +\left[c_{1} -\left(a_{1} +b+1\right)x\right]u_{x}  {-\left(a_{1} +b+1\right)yu_{y} +a_{1} zu_{z} -a_{1} bu=0,}  \\& y(1-y)u_{yy} -x^{2} u_{xx}-2xyu_{xy} +xzu_{xz}   \\& \,\,\,\,\,\,\,\,\,\,\,\,\,+yzu_{yz}-\left(a_{1} +b+1\right)xu_{x}+\left[c_{2} -\left(a_{1} +b+1\right)y\right]u_{y}  +a_{1} zu_{z} -a_{1} bu=0, \\& {zu_{zz} -xu_{xz} -yu_{yz} +(1-b+z)u_{z} +a_{2} u=0,}
\end{aligned}
\right.
$

where $u\equiv \,\,{\rm{E}}_{178}\left(a_1,a_2, b; c_1, c_2; x,y,z\right)$.

Particular solutions:

$
{u_1} ={\rm{E}}_{178}\left(a_1,a_2, b; c_1, c_2; x,y,z\right),
$

$
{u_2} = {x^{1 - c_1}}{\rm{E}}_{178}\left(1-c_1+a_1,a_2, 1-c_1+b; 2-c_1, c_2; x,y,z\right),
$

$
{u_3} = {y^{1 - c_2}}{\rm{E}}_{178}\left(1-c_2+a_1,a_2, 1-c_2+b; c_1, 2-c_2; x,y,z\right),
$

$
{u_4} = {x^{1 - c_1}}{y^{1 - c_2}}{\rm{E}}_{178}\left(2-c_1-c_2+a_1,a_2, 2-c_1-c_2+b; 2-c_1, 2-c_2; x,y,z\right).
$

\bigskip

\begin{equation}
{\rm{E}}_{179}\left(a_1,a_2, b; c_1, c_2; x,y,z\right)=\sum\limits_{m,n,p=0}^\infty\frac{(a_1)_{m+n}(a_2)_{n+p}(b)_{m-p}}{(c_1)_{m}(c_2)_{n}}\frac{x^m}{m!}\frac{y^n}{n!}\frac{z^p}{p!},
\end{equation}

region of convergence:
$$ \left\{ r+s<1,\,\,\,t<\infty
\right\},
$$

System of partial differential equations:

$
\left\{
\begin{aligned}
&{x(1-x)u_{xx} -xyu_{xy} +xzu_{xz} +yzu_{yz}  } {+\left[c_{1} -\left(a_{1} +b+1\right)x\right]u_{x}-byu_{y} +a_{1} zu_{z} -a_{1} bu=0,}  \\
& {y(1-y)u_{yy} -xyu_{xy} -xzu_{xz} -yzu_{yz} -a_{2} xu_{x}}{+\left[c_{2} -\left(a_{1} +a_{2} +1\right)y\right]u_{y} -a_{1} zu_{z} -a_{1} a_{2} u=0,}  \\
& {zu_{zz} -xu_{xz}+yu_{y} +(1-b+z)u_{z}  +a_{2} u=0,}
\end{aligned}
\right.
$

where $u\equiv \,\,{\rm{E}}_{179}\left(a_1,a_2, b; c_1, c_2; x,y,z\right)$.

Particular solutions:

$
{u_1} ={\rm{E}}_{179}\left(a_1,a_2, b; c_1, c_2; x,y,z\right),
$

$
{u_2} = {x^{1 - c_1}}{\rm{E}}_{179}\left(1-c_1+a_1,a_2, 1-c_1+b; 2-c_1, c_2; x,y,z\right),
$

$
{u_3} = {y^{1 - c_2}}{\rm{E}}_{179}\left(1-c_2+a_1,1-c_2+a_2, b; c_1, 2-c_2; x,y,z\right),
$

$
{u_4} = {x^{1 - c_1}}{y^{1 - c_2}}{\rm{E}}_{179}\left(2-c_1-c_2+a_1,1-c_2+a_2, 1-c_1+b; 2-c_1, 2-c_2; x,y,z\right).
$

\bigskip

\begin{equation}
{\rm{E}}_{180}\left(a_1,a_2, b; c_1, c_2; x,y,z\right)=\sum\limits_{m,n,p=0}^\infty\frac{(a_1)_{m+p}(a_2)_{n+p}(b)_{m-p}}{(c_1)_{m}(c_2)_{n}}\frac{x^m}{m!}\frac{y^n}{n!}\frac{z^p}{p!},
\end{equation}

region of convergence:
$$ \left\{ \left\{t+2\sqrt{rt}<1\right\}=\left\{\sqrt{t}<\sqrt{1+r}-\sqrt{r}\right\},\,\,\,s<\infty
\right\}.
$$

System of partial differential equations:

$
\left\{
\begin{aligned}
&{x(1-x)u_{xx} +z^{2} u_{zz} +\left[c_{1} -\left(a_{1} +b+1\right)x\right]u_{x}}{+\left(a_{1} -b+1\right)zu_{z} -a_{1} bu=0,} \\
& {yu_{yy} +\left(c_{2} -y\right)u_{y} -zu_{z} -a_{2} u=0,} \\
& {z(1+z)u_{zz} +xyu_{xy}
}-x(1-z)u_{xz}\\& \,\,\,\,\,\,\,\,\,  +yzu_{yz}+a_{2} xu_{x} +a_{1} yu_{y}{+\left[1-b+\left(a_{1} +a_{2} +1\right)z\right]u_{z} + a_{1} a_{2} u=0,}
\end{aligned}
\right.
$

where $u\equiv \,\,{\rm{E}}_{180}\left(a_1,a_2, b; c_1, c_2; x,y,z\right)$.

Particular solutions:

$
{u_1} ={\rm{E}}_{180}\left(a_1,a_2, b; c_1, c_2; x,y,z\right),
$

$
{u_2} = {x^{1 - c_1}}{\rm{E}}_{180}\left(1-c_1+a_1,a_2, 1-c_1+b; 2-c_1, c_2; x,y,z\right),
$

$
{u_3} = {y^{1 - c_2}}{\rm{E}}_{180}\left(a_1,1-c_2+a_2, b; c_1, 2-c_2; x,y,z\right),
$

$
{u_4} = {x^{1 - c_1}}{y^{1 - c_2}}{\rm{E}}_{180}\left(1-c_1+a_1,1-c_2+a_2, 1-c_1+b; 2-c_1, 2-c_2; x,y,z\right).
$

\bigskip

\begin{equation}
{\rm{E}}_{181}\left(a_1,a_2,b; c_1, c_2; x,y,z\right)=\sum\limits_{m,n,p=0}^\infty\frac{(a_1)_{m+n+p}(a_2)_{n}(b)_{m-p}}{(c_1)_{m}(c_2)_{n}}\frac{x^m}{m!}\frac{y^n}{n!}\frac{z^p}{p!},
\end{equation}

region of convergence:
$$ \left\{ r+s<1,\,\,\,t<\infty
\right\},
$$

System of partial differential equations:

$
\left\{
\begin{aligned}
&x(1-x)u_{xx} +z^{2} u_{zz} -xyu_{xy} +yzu_{yz} \\& \,\,\,\,\,\,\,\,\,+\left[c_{1} -\left(a_{1} +b+1\right)x\right]u_{x} {-byu_{y} +\left(a_{1} -b+1\right)zu_{z} -a_{1} bu=0,} \\
& {y(1-y)u_{yy} -xyu_{xy} -yzu_{yz} -a_{2} xu_{x} } {+\left[c_{2} -\left(a_{1} +a_{2} +1\right)y\right]u_{y} -a_{2} zu_{z} -a_{1} a_{2} u=0,} \\
& {zu_{zz} -xu_{xz} +xu_{x} +yu_{y} +(1-b+z)u_{z} +a_{1} u=0,}
\end{aligned}
\right.
$

where $u\equiv \,\,{\rm{E}}_{181}\left(a_1,a_2,b; c_1, c_2; x,y,z\right)$.

Particular solutions:

$
{u_1} ={\rm{E}}_{181}\left(a_1,a_2,b; c_1, c_2; x,y,z\right),
$

$
{u_2} = {x^{1 - c_1}}{\rm{E}}_{181}\left(1-c_1+a_1,a_2,1-c_1+b; 2-c_1, c_2; x,y,z\right),
$

$
{u_3} = {y^{1 - c_2}}{\rm{E}}_{181}\left(1-c_2+a_1,1-c_2+a_2,b; c_1, 2-c_2; x,y,z\right),
$

$
{u_4} = {x^{1 - c_1}}{y^{1 - c_2}}{\rm{E}}_{181}\left(2-c_1-c_2+a_1,1-c_2+a_2,1-c_1+b; 2-c_1, 2-c_2; x,y,z\right).
$

\bigskip

\begin{equation}
{\rm{E}}_{182}\left(a_1,a_2,b; c_1, c_2; x,y,z\right)=\sum\limits_{m,n,p=0}^\infty\frac{(a_1)_{m+n+p}(a_2)_p(b)_{m-p}}{(c_1)_{m}(c_2)_{n}}\frac{x^m}{m!}\frac{y^n}{n!}\frac{z^p}{p!},
\end{equation}

region of convergence:
$$ \left\{ \left\{t+2\sqrt{rt}<1\right\}=\left\{\sqrt{t}<\sqrt{1+r}-\sqrt{r}\right\},\,\,\,s<\infty
\right\}.
$$

System of partial differential equations:

$
\left\{
\begin{aligned}
&x(1-x)u_{xx} +z^{2} u_{zz} -xyu_{xy} +yzu_{yz} \\& \,\,\,\,\,\,\,\,\,+\left[c_{1} -\left(a_{1} +b+1\right)x\right]u_{x} {-byu_{y} +\left(a_{1} -b+1\right)zu_{z} -a_{1} bu=0,}  \\& {yu_{yy} -xu_{x} +\left(c_{2} -y\right)u_{y} -zu_{z} -a_{1} u=0,} \\& z(1+z)u_{zz} -x(1-z)u_{xz} +yzu_{yz}\\& \,\,\,\,\,\,\,\,\, +a_{2} xu_{x} +a_{2} yu_{y}   {+\left[1-b+\left(a_{1} +a_{2} +1\right)z\right]u_{z}  +a_{1} a_{2} u=0,}
\end{aligned}
\right.
$

where $u\equiv \,\,{\rm{E}}_{182}\left(a_1,a_2,b; c_1, c_2; x,y,z\right)$.

Particular solutions:

$
{u_1} ={\rm{E}}_{182}\left(a_1,a_2,b; c_1, c_2; x,y,z\right),
$

$
{u_2} = {x^{1 - c_1}}{\rm{E}}_{182}\left(1-c_1+a_1,a_2,1-c_1+b; 2-c_1, c_2; x,y,z\right),
$

$
{u_3} = {y^{1 - c_2}}{\rm{E}}_{182}\left(1-c_2+a_1,a_2,b; c_1, 2-c_2; x,y,z\right),
$

$
{u_4} = {x^{1 - c_1}}{y^{1 - c_2}}{\rm{E}}_{182}\left(2-c_1-c_2+a_1,a_2,1-c_1+b; 2-c_1, 2-c_2; x,y,z\right).
$

\bigskip

\begin{equation}
{\rm{E}}_{183}\left(a,b; c_1, c_2; x,y,z\right)=\sum\limits_{m,n,p=0}^\infty\frac{(a)_{m+n+p}(b)_{m-p}}{(c_1)_{m}(c_2)_{n}}\frac{x^m}{m!}\frac{y^n}{n!}\frac{z^p}{p!},
\end{equation}

region of convergence:
$$ \left\{ r<1,\,\,\,s<\infty,\,\,\,t<\infty
\right\}.
$$

System of partial differential equations:

$
\left\{
\begin{aligned}
&x(1-x)u_{xx} +z^{2} u_{zz} -xyu_{xy} +yzu_{yz}\\& \,\,\,\,\,\,\,\,\, +\left[c_{1} -\left(a +b+1\right)x\right]u_{x} {-byu_{y} +\left(a -b+1\right)zu_{z} -a bu=0,}  \\& {yu_{yy} -xu_{x} +\left(c_{2} -y\right)u_{y} -zu_{z} -a u=0,} \\& {zu_{zz} -xu_{xz} +xu_{x} +yu_{y} +(1-b+z)u_{z} +a u=0,}
\end{aligned}
\right.
$

where $u\equiv \,\,{\rm{E}}_{183}\left(a,b; c_1, c_2; x,y,z\right)$.

Particular solutions:

$
{u_1} ={\rm{E}}_{183}\left(a,b; c_1, c_2; x,y,z\right),
$

$
{u_2} = {x^{1 - c_1}}{\rm{E}}_{183}\left(1-c_1+a,1-c_1+b; 2-c_1, c_2; x,y,z\right),
$

$
{u_3} = {y^{1 - c_2}}{\rm{E}}_{183}\left(1-c_2+a,b; c_1, 2-c_2; x,y,z\right),
$

$
{u_4} = {x^{1 - c_1}}{y^{1 - c_2}}{\rm{E}}_{183}\left(2-c_1-c_2+a,1-c_1+b; 2-c_1, 2-c_2; x,y,z\right).
$

\bigskip

\begin{equation}
{\rm{E}}_{184}\left(a,b; c_1, c_2; x,y,z\right)=\sum\limits_{m,n,p=0}^\infty\frac{(a)_{m+n+p}(b)_{m+n-p}}{(c_1)_{m}(c_2)_{n}}\frac{x^m}{m!}\frac{y^n}{n!}\frac{z^p}{p!},
\end{equation}

region of convergence:
$$ \left\{ \sqrt{r}+\sqrt{s}<1,\,\,\,t<\infty
\right\}.
$$

System of partial differential equations:

$
\left\{
\begin{aligned}
&{x(1-x)u_{xx} -y^{2} u_{yy}+z^{2} u_{zz} -2xyu_{xy}  }\\
&\,\,\,\,\,\,\,\,\,\,\,\,+\left[c_{1} -\left(a+b+1\right)x\right]u_{x}{ -(a+b+1)yu_{y} +(a-b+1)zu_{z} -abu=0,} \\
& {y(1-y)u_{yy} -x^{2} u_{xx} +z^{2} u_{zz} -2xyu_{xy}}\\
&\,\,\,\,\,\,\,\,\,\,\,\, -(a+b+1)xu_{x}{ +\left[c_{2} -(a+b+1)y\right]u_{y} +(a-b+1)zu_{z} -abu=0,} \\ &{zu_{zz} -xu_{xz} -yu_{yz} +xu_{x} +yu_{y} +(1-b+z)u_{z} +au=0,}
\end{aligned}
\right.
$

where $u\equiv \,\,{\rm{E}}_{184}\left(a,b; c_1, c_2; x,y,z\right)$.

Particular solutions:

$
{u_1} ={\rm{E}}_{184}\left(a,b; c_1, c_2; x,y,z\right),
$

$
{u_2} = {x^{1 - c_1}}{\rm{E}}_{184}\left(1-c_1+a,1-c_1+b; 2-c_1, c_2; x,y,z\right),
$

$
{u_3} = {y^{1 - c_2}}{\rm{E}}_{184}\left(1-c_2+a,1-c_2+b; c_1, 2-c_2; x,y,z\right),
$

$
{u_4} = {x^{1 - c_1}}{y^{1 - c_2}}{\rm{E}}_{184}\left(2-c_1-c_2+a,2-c_1-c_2+b; 2-c_1, 2-c_2; x,y,z\right).
$

\bigskip

\begin{equation}
{\rm{E}}_{185}\left(a_1,a_2,b; c_1, c_2; x,y,z\right)=\sum\limits_{m,n,p=0}^\infty\frac{(a_1)_{n+p}(a_2)_{n}(b)_{2m-p}}{(c_1)_{m}(c_2)_{n} }\frac{x^m}{m!}\frac{y^n}{n!}\frac{z^p}{p!},
\end{equation}

region of convergence:
$$ \left\{ r<\frac{1}{4}, \,\,\, s<1,\,\,\,t<\infty
\right\}.
$$

System of partial differential equations:

$
\left\{
\begin{aligned}
&{x(1-4x)u_{xx} -z^{2} u_{zz} +4xzu_{xz}} {+\left[c_{1} -(4b+6)x\right]u_{x} +2bzu_{z} -b(b+1)u=0,} \\
& {y(1-y)u_{yy} -yzu_{yz} +\left[c_{2} -\left(a_{1} +a_{2} +1\right)y\right]u_{y}-a_{2} zu_{z} -a_{1} a_{2} u=0,} \\
& {zu_{zz} -2xu_{xz} +yu_{y}  +(1-b+z)u_{z} +a_{1} u=0,}
\end{aligned}
\right.
$

where $u\equiv \,\,{\rm{E}}_{185}\left(a_1,a_2,b; c_1, c_2; x,y,z\right)$.

Particular solutions:

$
{u_1} ={\rm{E}}_{185}\left(a_1,a_2,b; c_1, c_2; x,y,z\right),
$

$
{u_2} = {x^{1 - c_1}}{\rm{E}}_{185}\left(a_1,a_2,2-2c_1+b; 2-c_1, c_2; x,y,z\right),
$

$
{u_3} = {y^{1 - c_2}}{\rm{E}}_{185}\left(1-c_2+a_1,1-c_2+a_2,b; c_1, 2-c_2; x,y,z\right),
$

$
{u_4} = {x^{1 - c_1}}{y^{1 - c_2}}{\rm{E}}_{185}\left(1-c_2+a_1,1-c_2+a_2,2-2c_1+b; 2-c_1, 2-c_2; x,y,z\right).
$

\bigskip

\begin{equation}
{\rm{E}}_{186}\left(a_1,a_2,b; c_1, c_2; x,y,z\right)=\sum\limits_{m,n,p=0}^\infty\frac{(a_1)_{n+p}(a_2)_p(b)_{2m-p}}{(c_1)_{m}(c_2)_{n} }\frac{x^m}{m!}\frac{y^n}{n!}\frac{z^p}{p!},
\end{equation}

region of convergence:
$$ \left\{ \left(1+2\sqrt{r}\right)t<1,\,\,\,s<\infty
\right\}.
$$

System of partial differential equations:

$
\left\{
\begin{aligned}
&{x(1-4x)u_{xx} -z^{2} u_{zz} +4xzu_{xz}}{+\left[c_{1} -(4b+6)x\right]u_{x} +2bzu_{z} -b(b+1)u=0,} \\
& {yu_{yy} +\left(c_{2} -y\right)u_{y} -zu_{z} -a_{1} u=0,} \\
& {z(1+z)u_{zz} -2xu_{xz} +yzu_{yz} +a_{2} yu_{y}}{+\left[1-b+\left(a_{1} +a_{2} +1\right)z\right]u_{z} +a_{1} a_{2} u=0,}
\end{aligned}
\right.
$

where $u\equiv \,\,{\rm{E}}_{186}\left(a_1,a_2,b; c_1, c_2; x,y,z\right)$.

Particular solutions:

$
{u_1} ={\rm{E}}_{186}\left(a_1,a_2,b; c_1, c_2; x,y,z\right),
$

$
{u_2} = {x^{1 - c_1}}{\rm{E}}_{186}\left(a_1,a_2,2-2c_1+b; 2-c_1, c_2; x,y,z\right),
$

$
{u_3} = {y^{1 - c_2}}{\rm{E}}_{186}\left(1-c_2+a_1,1-c_2+a_2,b; c_1, 2-c_2; x,y,z\right),
$

$
{u_4} = {x^{1 - c_1}}{y^{1 - c_2}}{\rm{E}}_{186}\left(1-c_2+a_1,1-c_2+a_2,2-2c_1+b; 2-c_1, 2-c_2; x,y,z\right).
$

\bigskip

\begin{equation}
{\rm{E}}_{187}\left(a_1,a_2,b; c_1, c_2; x,y,z\right)=\sum\limits_{m,n,p=0}^\infty\frac{(a_1)_{n}(a_2)_{p}(b)_{2m+n-p}}{(c_1)_{m}(c_2)_{n} }\frac{x^m}{m!}\frac{y^n}{n!}\frac{z^p}{p!},
\end{equation}

region of convergence:
$$ \left\{ 2\sqrt{r}+s<1,\,\,\,t<\infty
\right\}.
$$

System of partial differential equations:

$
\left\{
\begin{aligned}
&{x(1-4x)u_{xx} -y^{2} u_{yy}-z^{2} u_{zz} -4xyu_{xy}  }\\&\,\,\,\,\,\,\,\,\,\,\,\,\,+4xzu_{xz} +2yzu_{yz} {+\left[c_{1} -(4b+6)x\right]u_{x} } {-2(b+1)yu_{y} +2bzu_{z} -b(b+1)u=0,} \\ &{y(1-y)u_{yy} -2xyu_{xy} +yzu_{yz} -2a_{1} xu_{x}} {+\left[c_{2} -\left(a_{1} +b+1\right)y\right]u_{y} +a_{1} zu_{z} -a_{1} bu=0,} \\& {zu_{zz} -2xu_{xz} -yu_{yz} +(1-b+z)u_{z} +a_{2} u=0,}
\end{aligned}
\right.
$

where $u\equiv \,\,{\rm{E}}_{187}\left(a_1,a_2,b; c_1, c_2; x,y,z\right)$.

Particular solutions:

$
{u_1} ={\rm{E}}_{187}\left(a_1,a_2,b; c_1, c_2; x,y,z\right),
$

$
{u_2} = {x^{1 - c_1}}{\rm{E}}_{187}\left(a_1,a_2,2-2c_1+b; 2-c_1, c_2; x,y,z\right),
$

$
{u_3} = {y^{1 - c_2}}{\rm{E}}_{187}\left(1-c_2+a_1,a_2,1-c_2+b; c_1, 2-c_2; x,y,z\right),
$

$
{u_4} = {x^{1 - c_1}}{y^{1 - c_2}}{\rm{E}}_{187}\left(1-c_2+a_1,a_2,3-2c_1-c_2+b; 2-c_1, 2-c_2; x,y,z\right).
$

\bigskip

\begin{equation}
{\rm{E}}_{188}\left(a_1,a_2,b; c_1, c_2; x,y,z\right)=\sum\limits_{m,n,p=0}^\infty\frac{(a_1)_{p}(a_2)_p(b)_{2m+n-p}}{(c_1)_{m}(c_2)_{n}}\frac{x^m}{m!}\frac{y^n}{n!}\frac{z^p}{p!},
\end{equation}

region of convergence:
$$ \left\{ \left(1+2\sqrt{r}\right)t<1,\,\,\,s<\infty
\right\}.
$$

System of partial differential equations:

$
\left\{
\begin{aligned}
&{x(1-4x)u_{xx} -y^{2} u_{yy} -z^{2} u_{zz} -4xyu_{xy} }\\& \,\,\,\,\,\,\,\,\,\,\,\,\,+4xzu_{xz} +2yzu_{yz}{+\left[c_{1} -(4b+6)x\right]u_{x} }  {-2(b+1)yu_{y} +2bzu_{z} -b(b+1)u=0,} \\ & {yu_{yy} -2xu_{x} +\left(c_{2} -y\right)u_{y} +zu_{z} -bu=0,} \\ &{z(1+z)u_{zz} -2xu_{xz} -yu_{yz} +\left[1-b+\left(a_{1} +a_{2} +1\right)z\right]u_{z} +a_{1} a_{2} u=0,}
\end{aligned}
\right.
$

where $u\equiv \,\,{\rm{E}}_{188}\left(a_1,a_2,b; c_1, c_2; x,y,z\right)$.

Particular solutions:

$
{u_1} ={\rm{E}}_{188}\left(a_1,a_2,b; c_1, c_2; x,y,z\right),
$

$
{u_2} = {x^{1 - c_1}}{\rm{E}}_{188}\left(a_1,a_2,2-2c_1+b; 2-c_1, c_2; x,y,z\right),
$

$
{u_3} = {y^{1 - c_2}}{\rm{E}}_{188}\left(a_1,a_2,1-c_2+b; c_1, 2-c_2; x,y,z\right),
$

$
{u_4} = {x^{1 - c_1}}{y^{1 - c_2}}{\rm{E}}_{188}\left(a_1,a_2,3-2c_1-c_2+b; 2-c_1, 2-c_2; x,y,z\right).
$

\bigskip

\begin{equation}
{\rm{E}}_{189}\left(a,b; c_1, c_2; x,y,z\right)=\sum\limits_{m,n,p=0}^\infty\frac{(a)_{n}(b)_{2m+n-p}}{(c_1)_{m}(c_2)_{n} }\frac{x^m}{m!}\frac{y^n}{n!}\frac{z^p}{p!},
\end{equation}

region of convergence:
$$ \left\{ 2\sqrt{r}+s<1,\,\,\,t<\infty
\right\}.
$$

System of partial differential equations:

$
\left\{
\begin{aligned}
& x(1-4x)u_{xx} -y^{2} u_{yy} -z^{2} u_{zz} -4xyu_{xy} \\
&\,\,\,\,\,\,\,\,\,\,\,\,\,+4xzu_{xz} +2yzu_{yz} +\left[c_{1} -(4b+6)x\right]u_{x} {-2(b+1)yu_{y} +2bzu_{z} -b(b+1)u=0,} \\
& y(1-y)u_{yy}-2xyu_{xy}+yzu_{yz}-2axu_x+\left[c_2-(a+b+1)y\right]u_y+azu_z-abu=0, \\
&{zu_{zz} -2xu_{xz} -yu_{yz} +(1-b)u_z+u=0,}
\end{aligned}
\right.
$

where $u\equiv \,\,{\rm{E}}_{189}\left(a,b; c_1, c_2; x,y,z\right)$.

Particular solutions:

$
{u_1} ={\rm{E}}_{189}\left(a,b; c_1, c_2; x,y,z\right),
$

$
{u_2} = {x^{1 - c_1}}{\rm{E}}_{189}\left(a,2-2c_1+b; 2-c_1, c_2; x,y,z\right),
$

$
{u_3} = {y^{1 - c_2}}{\rm{E}}_{189}\left(1-c_2+a,1-c_2+b; c_1, 2-c_2; x,y,z\right),
$

$
{u_4} = {x^{1 - c_1}}{y^{1 - c_2}}{\rm{E}}_{189}\left(1-c_2+a,3-2c_1-c_2+b; 2-c_1, 2-c_2; x,y,z\right).
$

\bigskip

\begin{equation}
{\rm{E}}_{190}\left(a,b; c_1, c_2; x,y,z\right)=\sum\limits_{m,n,p=0}^\infty\frac{(a)_{p}(b)_{2m+n-p}}{(c_1)_{m}(c_2)_{n} }\frac{x^m}{m!}\frac{y^n}{n!}\frac{z^p}{p!},
\end{equation}

region of convergence:
$$ \left\{ r<\frac{1}{4},\,\,\,s<\infty,\,\,\,t<\infty
\right\}.
$$

System of partial differential equations:

$
\left\{
\begin{aligned}
&x(1-4x)u_{xx} -y^{2} u_{yy}-z^{2} u_{zz}-4xyu_{xy} +4xzu_{xz}\\&\,\,\,\,\,\,\,\,\,\,\,\,\,   +2yzu_{yz} +\left[c_{1} -(4b+6)x\right]u_{x} {-2(b+1)yu_{y} +2bzu_{z} -b(b+1)u=0,} \\ & {yu_{yy} -2xu_{x} +\left(c_{2} -y\right)u_{y} +zu_{z} -bu=0,} \\& {zu_{zz} -2xu_{xz} -yu_{yz} +(1-b+z)u_{z} +au=0,}
\end{aligned}
\right.
$

where $u\equiv \,\,{\rm{E}}_{190}\left(a,b; c_1, c_2; x,y,z\right)$.

Particular solutions:

$
{u_1} ={\rm{E}}_{190}\left(a,b; c_1, c_2; x,y,z\right),
$

$
{u_2} = {x^{1 - c_1}}{\rm{E}}_{190}\left(a,2-2c_1+b; 2-c_1, c_2; x,y,z\right),
$

$
{u_3} = {y^{1 - c_2}}{\rm{E}}_{190}\left(a,1-c_2+b; c_1, 2-c_2; x,y,z\right),
$

$
{u_4} = {x^{1 - c_1}}{y^{1 - c_2}}{\rm{E}}_{190}\left(a,3-2c_1-c_2+b; 2-c_1, 2-c_2; x,y,z\right).
$

\bigskip

\begin{equation}
{\rm{E}}_{191}\left(b; c_1, c_2; x,y,z\right)=\sum\limits_{m,n,p=0}^\infty\frac{(b)_{2m+n-p}}{(c_1)_{m}(c_2)_{n} }\frac{x^m}{m!}\frac{y^n}{n!}\frac{z^p}{p!},
\end{equation}

region of convergence:
$$ \left\{ r<\frac{1}{4},\,\,\,s<\infty,\,\,\,t<\infty
\right\}.
$$

System of partial differential equations:

$
\left\{
\begin{aligned}
&x(1-4x)u_{xx} -y^{2} u_{yy}-z^{2} u_{zz}-4xyu_{xy} +4xzu_{xz} \\
&\,\,\,\,\,\,\,\,\,\,\,\,\,  +2yzu_{yz} +\left[c_{1} -(4b+6)x\right]u_{x} {-2(b+1)yu_{y} +2bzu_{z} -b(b+1)u=0,} \\
& {yu_{yy} -2xu_{x} +\left(c_{2} -y\right)u_{y} +zu_{z} -bu=0,} \\
& {zu_{zz} -2xu_{xz} -yu_{yz} +(1-b)u_z+u=0,}
\end{aligned}
\right.
$

where $u\equiv \,\,{\rm{E}}_{191}\left(b; c_1, c_2; x,y,z\right)$.

Particular solutions:

$
{u_1} ={\rm{E}}_{191}\left(b; c_1, c_2; x,y,z\right),
$

$
{u_2} = {x^{1 - c_1}}{\rm{E}}_{191}\left(2-2c_1+b; 2-c_1, c_2; x,y,z\right),
$

$
{u_3} = {y^{1 - c_2}}{\rm{E}}_{191}\left(1-c_2+b; c_1, 2-c_2; x,y,z\right),
$

$
{u_4} = {x^{1 - c_1}}{y^{1 - c_2}}{\rm{E}}_{191}\left(3-2c_1-c_2+b; 2-c_1, 2-c_2; x,y,z\right).
$

\bigskip

\begin{equation}
{\rm{E}}_{192}\left(a_1,a_2,b; c_1, c_2; x,y,z\right)=\sum\limits_{m,n,p=0}^\infty\frac{(a_1)_{2m+n}(a_2)_{p}(b)_{n-p}}{(c_1)_{m}(c_2)_{n} }\frac{x^m}{m!}\frac{y^n}{n!}\frac{z^p}{p!},
\end{equation}

region of convergence:
$$ \left\{ 2\sqrt{r}+s<1,\,\,\,t<\infty
\right\}.
$$

System of partial differential equations:

$
\left\{
\begin{aligned}
&{x(1-4x)u_{xx} -y^{2} u_{yy} -4xyu_{xy} +\left[c_{1} -\left(4a_{1} +6\right)x\right]u_{x}}{-2\left(a_{1} +1\right)yu_{y} -a_{1} \left(a_{1} +1\right)u=0,} \\
& y(1-y)u_{yy}-2xyu_{xy} +2xzu_{xz} +yzu_{yz} \\& \,\,\,\,\,\,\,\,\, -2bxu_{x}+\left[c_{2} -\left(a_{1} +b+1\right)y\right]u_{y} +a_{1} zu_{z} -a_{1} bu=0, \\
&{zu_{zz} -yu_{yz} +(1-b+z)u_{z} +a_{2} u=0,}
\end{aligned}
\right.
$

where $u\equiv \,\,{\rm{E}}_{192}\left(a_1,a_2,b; c_1, c_2; x,y,z\right)$.

Particular solutions:

$
{u_1} ={\rm{E}}_{192}\left(a_1,a_2,b; c_1, c_2; x,y,z\right),
$

$
{u_2} = {x^{1 - c_1}}{\rm{E}}_{192}\left(2-2c_1+a_1,a_2,b; 2-c_1, c_2; x,y,z\right),
$

$
{u_3} = {y^{1 - c_2}}{\rm{E}}_{192}\left(1-c_2+a_1,a_2,1-c_2+b; c_1, 2-c_2; x,y,z\right),
$

$
{u_4} = {x^{1 - c_1}}{y^{1 - c_2}}{\rm{E}}_{192}\left(3-2c_1-c_2+a_1,a_2,1-c_2+b; 2-c_1, 2-c_2; x,y,z\right).
$

\bigskip

\begin{equation}
{\rm{E}}_{193}\left(a,b; c_1, c_2; x,y,z\right)=\sum\limits_{m,n,p=0}^\infty\frac{(a)_{2m+n}(b)_{n-p}}{(c_1)_{m}(c_2)_{n} }\frac{x^m}{m!}\frac{y^n}{n!}\frac{z^p}{p!},
\end{equation}

region of convergence:
$$ \left\{ 2\sqrt{r}+s<1,\,\,\,t<\infty
\right\}.
$$

System of partial differential equations:

$
\left\{
\begin{aligned}
&{x(1-4x)u_{xx} -y^{2} u_{yy} -4xyu_{xy} +\left[c_{1} -(4a+6)x\right]u_{x} } {-2(a+1)yu_{y} -a(a+1)u=0,} \\
& y(1-y)u_{yy}-2xyu_{xy} +2xzu_{xz} +yzu_{yz}  -2bxu_{x}+\left[c_{2} -(a+b+1)y\right]u_{y} +azu_{z} -abu=0, \\
& {zu_{zz} -yu_{yz} +(1-b)u_{z} +u=0,}
\end{aligned}
\right.
$

where $u\equiv \,\,{\rm{E}}_{193}\left(a,b; c_1, c_2; x,y,z\right)$.

Particular solutions:

$
{u_1} ={\rm{E}}_{193}\left(a,b; c_1, c_2; x,y,z\right),
$

$
{u_2} = {x^{1 - c_1}}{\rm{E}}_{193}\left(2-2c_1+a,b; 2-c_1, c_2; x,y,z\right),
$

$
{u_3} = {y^{1 - c_2}}{\rm{E}}_{193}\left(1-c_2+a,1-c_2+b; c_1, 2-c_2; x,y,z\right),
$

$
{u_4} = {x^{1 - c_1}}{y^{1 - c_2}}{\rm{E}}_{193}\left(3-2c_1-c_2+a,1-c_2+b; 2-c_1, 2-c_2; x,y,z\right).
$

\bigskip

\begin{equation}
{\rm{E}}_{194}\left(a_1,a_2,b; c_1, c_2; x,y,z\right)=\sum\limits_{m,n,p=0}^\infty\frac{(a_1)_{2m+p}(a_2)_{n}(b)_{n-p}}{(c_1)_{m}(c_2)_{n} }\frac{x^m}{m!}\frac{y^n}{n!}\frac{z^p}{p!},
\end{equation}

region of convergence:
$$ \left\{ r<\frac{1}{4},\,s<1,\,\,\,t<\infty
\right\}.
$$

System of partial differential equations:

$
\left\{
\begin{aligned}
&{x(1-4x)u_{xx} -z^{2} u_{zz} -4xzu_{xz} +\left[c_{1} -\left(4a_{1} +6\right)x\right]u_{x}}{-2\left(a_{1} +1\right)zu_{z} -a_{1} (a_{1} +1)u=0,} \\& {y(1-y)u_{yy} +yzu_{yz} +\left[c_{2} -\left(a_{2} +b+1\right)y\right]u_{y} +a_{2} zu_{z} -a_{2} bu=0,} \\& {zu_{zz} -yu_{yz}+2xu_{x} +(1-b+z)u_{z} +a_{1} u=0,}
\end{aligned}
\right.
$

where $u\equiv \,\,{\rm{E}}_{194}\left(a_1,a_2,b; c_1, c_2; x,y,z\right)$.

Particular solutions:

$
{u_1} ={\rm{E}}_{194}\left(a_1,a_2,b; c_1, c_2; x,y,z\right),
$

$
{u_2} = {x^{1 - c_1}}{\rm{E}}_{194}\left(2-2c_1+a_1,a_2,b; 2-c_1, c_2; x,y,z\right),
$

$
{u_3} = {y^{1 - c_2}}{\rm{E}}_{194}\left(a_1,1-c_2+a_2,1-c_2+b; c_1, 2-c_2; x,y,z\right),
$

$
{u_4} = {x^{1 - c_1}}{y^{1 - c_2}}{\rm{E}}_{194}\left(2-2c_1+a_1,1-c_2+a_2,1-c_2+b; 2-c_1, 2-c_2; x,y,z\right).
$

\bigskip

\begin{equation}
{\rm{E}}_{195}\left(a_1,a_2,b; c_1, c_2; x,y,z\right)=\sum\limits_{m,n,p=0}^\infty\frac{(a_1)_{2m+p}(a_2)_p(b)_{n-p}}{(c_1)_{m}(c_2)_{n}}\frac{x^m}{m!}\frac{y^n}{n!}\frac{z^p}{p!},
\end{equation}

region of convergence:
$$ \left\{ 2\sqrt{r}+s<1,\,\,\,t<\infty
\right\}.
$$

System of partial differential equations:

$
\left\{
\begin{aligned}
&x(1-4x)u_{xx} -z^{2} u_{zz} -4xzu_{xz} +\left[c_{1} -(4a_1+6)x\right]u_{x}  {-2(a_1+1)zu_{z} -a_1(a_1+1)u=0,} \\
& {yu_{yy} +\left(c_{2} - y\right)u_{y} +zu_{z} -bu=0,} \\
& z(1+z)u_{zz} +2xzu_{xz}-yu_{yz}  +2a_{2} xu_{x}+\left[1-b+\left(a_{1} +a_{2} +1\right)z\right]u_{z} +a_{1} a_{2} u=0,
\end{aligned}
\right.
$

where $u\equiv \,\,{\rm{E}}_{195}\left(a_1,a_2,b; c_1, c_2; x,y,z\right)$.

Particular solutions:

$
{u_1} ={\rm{E}}_{195}\left(a_1,a_2,b; c_1, c_2; x,y,z\right),
$

$
{u_2} = {x^{1 - c_1}}{\rm{E}}_{195}\left(2-2c_1+a_1,a_2,b; 2-c_1, c_2; x,y,z\right),
$

$
{u_3} = {y^{1 - c_2}}{\rm{E}}_{195}\left(a_1,a_2,1-c_2+b; c_1, 2-c_2; x,y,z\right),
$

$
{u_4} = {x^{1 - c_1}}{y^{1 - c_2}}{\rm{E}}_{195}\left(2-2c_1+a_1,a_2,1-c_2+b; 2-c_1, 2-c_2; x,y,z\right).
$

\bigskip

\begin{equation}
{\rm{E}}_{196}\left(a,b; c_1, c_2; x,y,z\right)=\sum\limits_{m,n,p=0}^\infty\frac{(a)_{2m+p}(b)_{n-p}}{(c_1)_{m}(c_2)_{n}}\frac{x^m}{m!}\frac{y^n}{n!}\frac{z^p}{p!},
\end{equation}

region of convergence:
$$ \left\{ r<\frac{1}{4},\,\,\,s<\infty,\,\,\,t<\infty
\right\}.
$$

System of partial differential equations:

$
\left\{
\begin{aligned}
&x(1-4x)u_{xx} -z^{2} u_{zz} -4xzu_{xz} +\left[c_{1} -\left(4a +6\right)x\right]u_{x}{-2\left(a+1\right)zu_{z} -a(a +1)u=0,} \\
& {yu_{yy} +\left(c_{2}- y\right)u_{y} +zu_{z} -bu=0,} \\
&{zu_{zz} -yu_{yz}+2xu_{x} +(1-b+z)u_{z} +au=0,}
\end{aligned}
\right.
$

where $u\equiv \,\,{\rm{E}}_{196}\left(a,b; c_1, c_2; x,y,z\right)$.

Particular solutions:

$
{u_1} ={\rm{E}}_{196}\left(a,b; c_1, c_2; x,y,z\right),
$

$
{u_2} = {x^{1 - c_1}}{\rm{E}}_{196}\left(2-2c_1+a,b; 2-c_1, c_2; x,y,z\right),
$

$
{u_3} = {y^{1 - c_2}}{\rm{E}}_{196}\left(a,1-c_2+b; c_1, 2-c_2; x,y,z\right),
$

$
{u_4} = {x^{1 - c_1}}{y^{1 - c_2}}{\rm{E}}_{196}\left(2-2c_1+a,1-c_2+b; 2-c_1, 2-c_2; x,y,z\right).
$

\bigskip

\begin{equation}
{\rm{E}}_{197}\left(a_1,a_2,b; c_1, c_2; x,y,z\right)=\sum\limits_{m,n,p=0}^\infty\frac{(a_1)_{n+2p}(a_2)_{m}(b)_{m-p}}{(c_1)_{m}(c_2)_{n}}\frac{x^m}{m!}\frac{y^n}{n!}\frac{z^p}{p!},
\end{equation}

region of convergence:
$$ \left\{ r<1,\,\,\,s<\infty,\,\,\,t<\frac{1}{4}
\right\}.
$$

System of partial differential equations:

$
\left\{
\begin{aligned}
&{x(1-x)u_{xx} +xzu_{xz} +\left[c_{1} -\left(a_{2} +b+1\right)x\right]u_{x} +a_{2} zu_{z} -a_{2} bu=0,} \\& {yu_{yy} +\left(c_{2} -y\right)u_{y} -2zu_{z} -a_{1} u=0,} \\& z(1+4z)u_{zz} +y^{2} u_{yy} -xu_{xz} +4yzu_{yz} \\& \,\,\,\,\,\,\,\,\,+2(a_1+1)yu_{y} +\left[1-b+(4a_1+6)z\right]u_{z} +a_1\left(a_1+1\right)u=0,
\end{aligned}
\right.
$

where $u\equiv \,\,{\rm{E}}_{197}\left(a_1,a_2,b; c_1, c_2; x,y,z\right)$.

Particular solutions:

$
{u_1} ={\rm{E}}_{197}\left(a_1,a_2,b; c_1, c_2; x,y,z\right),
$

$
{u_2} = {x^{1 - c_1}}{\rm{E}}_{197}\left(a_1,1-c_1+a_2,1-c_1+b; 2-c_1, c_2; x,y,z\right),
$

$
{u_3} = {y^{1 - c_2}}{\rm{E}}_{197}\left(1-c_2+a_1,a_2,b; c_1, 2-c_2; x,y,z\right),
$

$
{u_4} = {x^{1 - c_1}}{y^{1 - c_2}}{\rm{E}}_{197}\left(1-c_2+a_1,1-c_1+a_2,1-c_1+b; 2-c_1, 2-c_2; x,y,z\right).
$

\bigskip

\begin{equation}
{\rm{E}}_{198}\left(a_1,a_2,b; c_1, c_2; x,y,z\right)=\sum\limits_{m,n,p=0}^\infty\frac{(a_1)_{n+2p}(a_2)_n(b)_{m-p}}{(c_1)_{m}(c_2)_{n}}\frac{x^m}{m!}\frac{y^n}{n!}\frac{z^p}{p!},
\end{equation}

region of convergence:
$$ \left\{ 2\sqrt{r}+s<1,\,\,\,t<\infty
\right\}.
$$

System of partial differential equations:

$
\left\{
\begin{aligned}
&{xu_{xx} +(c_{1} -x)u_{x} +zu_{z} -bu=0,} \\&
{y(1-y)u_{yy} -2yzu_{yz} +\left[c_{2} -\left(a_{1} +a_{2} +1\right)y\right]u_{y} -2a_{2} zu_{z} -a_{1} a_{2} u=0,} \\ &z(1+4z)u_{zz} +y^{2} u_{yy}-xu_{xz} +4yzu_{yz}\\& \,\,\,\,\,\,\,\,\, +2(a_1+1)yu_{y} +\left[1-b+\left(4a_1+6\right)z\right]u_{z} +a_1\left(a_1+1\right)u=0,
\end{aligned}
\right.
$

where $u\equiv \,\,{\rm{E}}_{198}\left(a_1,a_2,b; c_1, c_2; x,y,z\right)$.

Particular solutions:

$
{u_1} ={\rm{E}}_{198}\left(a_1,a_2,b; c_1, c_2; x,y,z\right),
$

$
{u_2} = {x^{1 - c_1}}{\rm{E}}_{198}\left(a_1,a_2,1-c_1+b; 2-c_1, c_2; x,y,z\right),
$

$
{u_3} = {y^{1 - c_2}}{\rm{E}}_{198}\left(1-c_2+a_1,1-c_2+a_2,b; c_1, 2-c_2; x,y,z\right),
$

$
{u_4} = {x^{1 - c_1}}{y^{1 - c_2}}{\rm{E}}_{198}\left(1-c_2+a_1,1-c_2+a_2,1-c_1+b; 2-c_1, 2-c_2; x,y,z\right).
$

\bigskip

\begin{equation}
{\rm{E}}_{199}\left(a,b; c_1, c_2; x,y,z\right)=\sum\limits_{m,n,p=0}^\infty\frac{(a)_{n+2p}(b)_{m-p}}{(c_1)_{m}(c_2)_{n}}\frac{x^m}{m!}\frac{y^n}{n!}\frac{z^p}{p!},
\end{equation}

region of convergence:
$$ \left\{ r<\infty,\,\,\,s<\infty,\,\,\,t<\frac{1}{4}
\right\}.
$$

System of partial differential equations:

$
\left\{
\begin{aligned}
&{xu_{xx} +\left(c_{1} -x\right)u_{x} +zu_{z} -bu=0,} \\ &{yu_{yy} +(c_{2} -y)u_{y} -2zu_{z} -a u=0,} \\& z(1+4z)u_{zz} +y^{2} u_{yy} -xu_{xz} +4yzu_{yz} +2(a+1)yu_{y} +\left[1-b+(4a+6)z\right]u_{z} +a(a+1)u=0,
\end{aligned}
\right.
$

where
$U \equiv \,\,{\rm{E}}_{199}\left(a,b; c_1, c_2; x,y,z\right)$.

Particular solutions:

$
{u_1} ={\rm{E}}_{199}\left(a,b; c_1, c_2; x,y,z\right),
$

$
{u_2} = {x^{1 - c_1}}{\rm{E}}_{199}\left(a,1-c_1+b; 2-c_1, c_2; x,y,z\right),
$

$
{u_3} = {y^{1 - c_2}}{\rm{E}}_{199}\left(1-c_2+a,b; c_1, 2-c_2; x,y,z\right),
$

$
{u_4} = {x^{1 - c_1}}{y^{1 - c_2}}{\rm{E}}_{199}\left(1-c_2+a,1-c_1+b; 2-c_1, 2-c_2; x,y,z\right).
$

\bigskip

\begin{equation}
{\rm{E}}_{200}\left(a,b; c_1, c_2; x,y,z\right)=\sum\limits_{m,n,p=0}^\infty\frac{(a)_{n+p}(b)_{2m-p}}{(c_1)_{m}(c_2)_{n}}\frac{x^m}{m!}\frac{y^n}{n!}\frac{z^p}{p!},
\end{equation}

region of convergence:
$$ \left\{ r<\frac{1}{4},\,\,\,s<\infty,\,\,\,t<\infty
\right\}.
$$

System of partial differential equations:

$
\left\{
\begin{aligned}
&x(1-4x)u_{xx} -z^{2} u_{zz} +4xzu_{xz}+\left[c_{1} -(4b+6)x\right]u_{x} +2bzu_{z} -b(b+1)u=0, \\& {yu_{yy} +(c_{2} -y)u_{y} -zu_{z} -au=0,} \\& {zu_{zz} -2xu_{xz} +yu_{y} +(1-b+z)u_{z} +au=0,}
\end{aligned}
\right.
$

where $u\equiv \,\,{\rm{E}}_{200}\left(a,b; c_1, c_2; x,y,z\right).$

Particular solutions:

$
{u_1} ={\rm{E}}_{200}\left(a,b; c_1, c_2; x,y,z\right),
$

$
{u_2} = {x^{1 - c_1}}{\rm{E}}_{200}\left(a,2-2c_1+b; 2-c_1, c_2; x,y,z\right),
$

$
{u_3} = {y^{1 - c_2}}{\rm{E}}_{200}\left(1-c_2+a,b; c_1, 2-c_2; x,y,z\right),
$

$
{u_4} = {x^{1 - c_1}}{y^{1 - c_2}}{\rm{E}}_{200}\left(1-c_2+a,2-2c_1+b; 2-c_1, 2-c_2; x,y,z\right).
$

\bigskip

\begin{equation}
{\rm{E}}_{201}\left(a, b; c_1, c_2; x,y,z\right)=\sum\limits_{m,n,p=0}^\infty\frac{(a)_{n+p}(b)_{2m+n-p}}{(c_1)_{m}(c_2)_{n}}\frac{x^m}{m!}\frac{y^n}{n!}\frac{z^p}{p!},
\end{equation}

region of convergence:
$$ \left\{ 2\sqrt{r}+s<1,\,\,\,t<\infty
\right\}.
$$

System of partial differential equations:

$
\left\{
\begin{aligned}
&{x(1-4x)u_{xx} -y^{2} u_{yy} }-z^{2} u_{zz} -4xyu_{xy} +4xzu_{xz}\\&\,\,\,\,\,\,\,\,\,\,\,\,\,\,  +2yzu_{yz} {+\left[c_{1} -(4b+6)x\right]u_{x} -2(b+1)yu_{y}} {+2bzu_{z} -b(b+1)u=0,} \\& {y(1-y)u_{yy} +z^{2} u_{zz} -2xyu_{xy} -2xzu_{xz} } \\& \,\,\,\,\,\,\,\,\,\,\,\,\,\,  -2axu_{x}{+\left[c_{2} -(a+b+1)y\right]u_{y} +(a-b+1)zu_{z} -abu=0,} \\& {z u_{zz} -2xu_{xz} -yu_{yz} +yu_{y} +(1-b+z)u_{z} +au=0,}
\end{aligned}
\right.
$

where $u\equiv \,\,{\rm{E}}_{201}\left(a, b; c_1, c_2; x,y,z\right)$.

Particular solutions:

$
{u_1} ={\rm{E}}_{201}\left(a,b; c_1, c_2; x,y,z\right),
$

$
{u_2} = {x^{1 - c_1}}{\rm{E}}_{201}\left(a,2-2c_1+b; 2-c_1, c_2; x,y,z\right),
$

$
{u_3} = {y^{1 - c_2}}{\rm{E}}_{201}\left(1-c_2+a,1-c_2+b; c_1, 2-c_2; x,y,z\right),
$

$
{u_4} = {x^{1 - c_1}}{y^{1 - c_2}}{\rm{E}}_{201}\left(1-c_2+a,3-2c_1-c_2+b; 2-c_1, 2-c_2; x,y,z\right).
$

\bigskip

\begin{equation}
{\rm{E}}_{202}\left(a, b; c_1, c_2; x,y,z\right)=\sum\limits_{m,n,p=0}^\infty\frac{(a)_{n+2p}(b)_{m+n-p}}{(c_1)_{m}(c_2)_{n}}\frac{x^m}{m!}\frac{y^n}{n!}\frac{z^p}{p!},
\end{equation}

region of convergence:
$$ \left\{ r<\infty,\,\,\,\left\{t<\frac{1}{4},\,\,\, s<\min\left\{\Psi_1(t), \Psi_2(t)\right\}\right\}=\left\{s<1, \,\,\,\, t<\min\left\{\Theta_1(s), \Theta_2(s)\right\}\right\}
\right\}.
$$

System of partial differential equations:

$
\left\{
\begin{aligned}
&{xu_{xx} +(c_{1} -x)u_{x} -yu_{y} +zu_{z} -bu=0,} \\
& y(1-y)u_{yy} +2z^{2} u_{zz} -xyu_{xy}-2xzu_{xz}\\&\,\,\,\,\,\,\,\,\,\,\,\,\,\,\,-yzu_{yz}-axu_x+\left[c_{2} -(a+b+1)y\right]u_{y} +(a-2b+2)zu_{z} -abu=0, \\
& z(1+4z)u_{zz} +y^{2} u_{yy} -xu_{xz}-y(1-4z)u_{yz}\\& \,\,\,\,\,\,\,\,\, +2(a+1)yu_{y} {+\left[1-b+(4a+6)z\right]u_{z} +a(a+1)u=0,}
\end{aligned}
\right.
$

where $u\equiv \,\,{\rm{E}}_{202}\left(a, b; c_1, c_2; x,y,z\right)$.

Particular solutions:

$
{u_1} ={\rm{E}}_{202}\left(a,b; c_1, c_2; x,y,z\right),
$

$
{u_2} = {x^{1 - c_1}}{\rm{E}}_{202}\left(a,1-c_1+b; 2-c_1, c_2; x,y,z\right),
$

$
{u_3} = {y^{1 - c_2}}{\rm{E}}_{202}\left(1-c_2+a,1-c_2+b; c_1, 2-c_2; x,y,z\right),
$

$
{u_4} = {x^{1 - c_1}}{y^{1 - c_2}}{\rm{E}}_{202}\left(1-c_2+a,2-c_1-c_2+b; 2-c_1, 2-c_2; x,y,z\right).
$

\bigskip

\begin{equation}
{\rm{E}}_{203}\left(a, b; c_1, c_2; x,y,z\right)=\sum\limits_{m,n,p=0}^\infty\frac{(a)_{m+n+2p}(b)_{n-p}}{(c_1)_{m}(c_2)_{n}}\frac{x^m}{m!}\frac{y^n}{n!}\frac{z^p}{p!},
\end{equation}

region of convergence:
$$ \left\{ r<\infty,\,\,\,\left\{t<\frac{1}{4}, \,\,\,\, s<\min\left\{\Psi_1(t), \Psi_2(t)\right\}\right\}=\left\{s<1,\,\,\,\, t<\min\left\{\Theta_1(s), \Theta_2(s)\right\}\right\}
\right\}.
$$

System of partial differential equations:

$
\left\{
\begin{aligned}
&{xu_{xx} +(c_{1} -x)u_{x} -yu_{y} -2zu_{z} -au=0,} \\
& {y(1-y)u_{yy} +2z^{2} u_{zz} -xyu_{xy} -yzu_{yz} }\\
&\,\,\,\,\,\,\,\,\,\,\,\,\,\, +xzu_{xz} -bxu_{x}{+\left[c_{2} -(a+b+1)y\right]u_{y}}{ +(a-2b+2)zu_{z} -abu=0,} \\
& z(1+4z)u_{zz} +x^{2} u_{xx} +y^{2} u_{yy} +2xyu_{xy} +4xzu_{xz} -y(1-4z)u_{yz}\\
&\,\,\,\,\,\,\,\,\,\,\,\,\,\,\, +2(a+1)xu_{x}  {+2(a+1)yu_{y} +\left[1-b+(4a+6)z\right]u_{z} +a(a+1)u=0,}
\end{aligned}
\right.
$

where $u\equiv \,\,{\rm{E}}_{203}\left(a, b; c_1, c_2; x,y,z\right)$.

Particular solutions:

$
{u_1} ={\rm{E}}_{203}\left(a,b; c_1, c_2; x,y,z\right),
$

$
{u_2} = {x^{1 - c_1}}{\rm{E}}_{203}\left(1-c_1+a,b; 2-c_1, c_2; x,y,z\right),
$

$
{u_3} = {y^{1 - c_2}}{\rm{E}}_{203}\left(1-c_2+a,1-c_2+b; c_1, 2-c_2; x,y,z\right),
$

$
{u_4} = {x^{1 - c_1}}{y^{1 - c_2}}{\rm{E}}_{203}\left(2-c_1-c_2+a,1-c_2+b; 2-c_1, 2-c_2; x,y,z\right).
$

\bigskip

\begin{equation}
{\rm{E}}_{204}\left(a, b; c_1, c_2; x,y,z\right)=\sum\limits_{m,n,p=0}^\infty\frac{(a)_{2m+n+p}(b)_{n-p}}{(c_1)_{m}(c_2)_{n}}\frac{x^m}{m!}\frac{y^n}{n!}\frac{z^p}{p!},
\end{equation}

region of convergence:
$$ \left\{ 2\sqrt{r}+s<1,\,\,\,t<\infty
\right\}.
$$

System of partial differential equations:

$
\left\{
\begin{aligned}
&{x(1-4x)u_{xx} -y^{2} u_{yy} -z^{2} u_{zz} -4xyu_{xy}-4xzu_{xz}  }\\& \,\,\,\,\,\,\,\,\,\,\,\,\,-2yzu_{yz} {+\left[c_{1} -(4a+6)x\right]u_{x}}{ -2(a+1)yu_{y} -2(a+1)zu_{z} -a(a+1)u=0,} \\ & {y(1-y)u_{yy} +z^{2} u_{zz}  -2xyu_{xy} +2xzu_{xz}}\\&\,\,\,\,\,\,\,\,\,\,\,\,\, -2bxu_{x}{+\left[c_{2} -(a+b+1)y\right]u_{y} +(a-b+1)zu_{z} -abu=0,} \\& {z u_{zz}  -yu_{yz} +2xu_x+yu_{y} +(1-b+z)u_{z} +au=0,}
\end{aligned}
\right.
$

where $u\equiv \,\,{\rm{E}}_{204}\left(a, b; c_1, c_2; x,y,z\right)$.

Particular solutions:

$
{u_1} ={\rm{E}}_{204}\left(a,b; c_1, c_2; x,y,z\right),
$

$
{u_2} = {x^{1 - c_1}}{\rm{E}}_{204}\left(2-2c_1+a,b; 2-c_1, c_2; x,y,z\right),
$

$
{u_3} = {y^{1 - c_2}}{\rm{E}}_{204}\left(1-c_2+a,1-c_2+b; c_1, 2-c_2; x,y,z\right),
$

$
{u_4} = {x^{1 - c_1}}{y^{1 - c_2}}{\rm{E}}_{204}\left(3-2c_1-c_2+a,1-c_2+b; 2-c_1, 2-c_2; x,y,z\right).
$

\bigskip

\begin{equation}
{\rm{E}}_{205}\left(a,b; c_1, c_2; x,y,z\right)=\sum\limits_{m,n,p=0}^\infty\frac{(a)_{n+2p}(b)_{2m-p}}{(c_1)_{m}(c_2)_{n}}\frac{x^m}{m!}\frac{y^n}{n!}\frac{z^p}{p!},
\end{equation}

region of convergence:
$$ \left\{ \left(1+2\sqrt{r}\right)t<\frac{1}{4},\,\,\,s<\infty
\right\},
$$

System of partial differential equations:

$
\left\{
\begin{aligned}
&{x(1-4x)u_{xx} -z^{2} u_{zz} +4xzu_{xz}} {+\left[c_{1} -(4b+6)x\right]u_{x}} {+2bzu_{z} -b(b+1)u=0,} \\
&{yu_{yy} +\left(c_{2} -y\right)u_{y} -2zu_{z} -a u=0,} \\
& z(1+4z)u_{zz} +y^{2} u_{yy} -2xu_{xz} +4yzu_{yz}\\& \,\,\,\,\,\,\,\,\, +2(a+1)yu_{y}  {+\left[1-b+(4a+6)z\right]u_{z} +a(a+1)u=0,}
\end{aligned}
\right.
$

where $u\equiv \,\,{\rm{E}}_{205}\left(a,b; c_1, c_2; x,y,z\right)$.

Particular solutions:

$
{u_1} ={\rm{E}}_{205}\left(a,b; c_1, c_2; x,y,z\right),
$

$
{u_2} = {x^{1 - c_1}}{\rm{E}}_{205}\left(a,2-2c_1+b; 2-c_1, c_2; x,y,z\right),
$

$
{u_3} = {y^{1 - c_2}}{\rm{E}}_{205}\left(1-c_2+a,b; c_1, 2-c_2; x,y,z\right),
$

$
{u_4} = {x^{1 - c_1}}{y^{1 - c_2}}{\rm{E}}_{205}\left(1-c_2+a,2-2c_1+b; 2-c_1, 2-c_2; x,y,z\right).
$

\bigskip

\begin{equation}
{\rm{E}}_{206}\left(a, b; c_1, c_2; x,y,z\right)=\sum\limits_{m,n,p=0}^\infty\frac{(a)_{2m+p}(b)_{2n-p}}{(c_1)_{m}(c_2)_{n}}\frac{x^m}{m!}\frac{y^n}{n!}\frac{z^p}{p!},
\end{equation}

region of convergence:
$$ \left\{ r<\frac{1}{4},\,\,\,\, s<\frac{1}{4},\,\,\,t<\infty
\right\}.
$$

System of partial differential equations:

$
\left\{
\begin{aligned}
&x(1-4x)u_{xx} -z^{2} u_{zz} -4xzu_{xz} +\left[c_{1} -(4a +6)x\right]u_{x}  {-2\left(a +1\right)zu_{z} -a \left(a +1\right)u=0,} \\
& y(1-4y)u_{yy} +z^{2} u_{zz}+4yzu_{yz} +\left[c_{2} -(4b+6)y\right]u_{y} +2bzu_{z} -b(b+1)u=0, \\
& {zu_{zz} -2yu_{yz} +2xu_{x} +(1-b+z)u_{z} +au=0,}
\end{aligned}
\right.
$

where $u\equiv \,\,{\rm{E}}_{206}\left(a, b; c_1, c_2; x,y,z\right)$.

Particular solutions:

$
{u_1} ={\rm{E}}_{206}\left(a,b; c_1, c_2; x,y,z\right),
$

$
{u_2} = {x^{1 - c_1}}{\rm{E}}_{206}\left(2-2c_1+a,b; 2-c_1, c_2; x,y,z\right),
$

$
{u_3} = {y^{1 - c_2}}{\rm{E}}_{206}\left(a,2-2c_2+b; c_1, 2-c_2; x,y,z\right),
$

$
{u_4} = {x^{1 - c_1}}{y^{1 - c_2}}{\rm{E}}_{206}\left(2-2c_1+a,2-2c_2+b; 2-c_1, 2-c_2; x,y,z\right).
$

\bigskip

\begin{equation}
{\rm{E}}_{207}\left(a, b; c_1, c_2; x,y,z\right)=\sum\limits_{m,n,p=0}^\infty\frac{(a)_{p}(b)_{2m+2n-p}}{(c_1)_{m}(c_2)_{n}}\frac{x^m}{m!}\frac{y^n}{n!}\frac{z^p}{p!},
\end{equation}

region of convergence:
$$ \left\{ \sqrt{r}+\sqrt{s}<\frac{1}{2},\,\,\,t<\infty
\right\}.
$$

System of partial differential equations:

$
\left\{
\begin{aligned}
&{x(1-4x)u_{xx}  -4y^{2} u_{yy} -z^{2} u_{zz}}-8xyu_{xy} +4xzu_{xz}\\
&\,\,\,\,\,\,\,\,\,\,\,\,\,   +4yzu_{yz} {+\left[c_{1} -(4b +6)x\right]u_{x} -(4b+6)yu_{y} +2bzu_{z} -b(b+1)u=0,} \\
& {y(1-4y)u_{yy} -4x^{2} u_{xx}}-z^{2} u_{zz}-8xyu_{xy} +4xzu_{xz}\\
& \,\,\,\,\,\,\,\,\,\,\,\,   +4yzu_{yz} {-(4b+6)xu_{x} +\left[c_{2} -(4b+6)y\right]u_{y} +2bzu_{z} -b(b+1)u=0,} \\
& {zu_{zz} -2xu_{xz} -2yu_{yz} +(1-b+z)u_{z} +au=0,}
\end{aligned}
\right.
$

where $u\equiv \,\,{\rm{E}}_{207}\left(a, b; c_1, c_2; x,y,z\right)$.

Particular solutions:

$
{u_1} ={\rm{E}}_{207}\left(a,b; c_1, c_2; x,y,z\right),
$

$
{u_2} = {x^{1 - c_1}}{\rm{E}}_{207}\left(a,2-2c_1+b; 2-c_1, c_2; x,y,z\right),
$

$
{u_3} = {y^{1 - c_2}}{\rm{E}}_{207}\left(a,2-2c_2+b; c_1, 2-c_2; x,y,z\right),
$

$
{u_4} = {x^{1 - c_1}}{y^{1 - c_2}}{\rm{E}}_{207}\left(a,4-2c_1-2c_2+b; 2-c_1, 2-c_2; x,y,z\right).
$

\bigskip

\begin{equation}
{\rm{E}}_{208}\left(b; c_1, c_2; x,y,z\right)=\sum\limits_{m,n,p=0}^\infty\frac{(b)_{2m+2n-p}}{(c_1)_{m}(c_2)_{n}}\frac{x^m}{m!}\frac{y^n}{n!}\frac{z^p}{p!},
\end{equation}

region of convergence:
$$ \left\{ \sqrt{r}+\sqrt{s}<\frac{1}{2},\,\,\,t<\infty
\right\}.
$$

System of partial differential equations:

$
\left\{
\begin{aligned}
&{x(1-4x)u_{xx}-4y^{2} u_{yy} -z^{2} u_{zz}  }  -8xyu_{xy} +4xzu_{xz} +4yzu_{yz}\\& \,\,\,\,\,\,\,\,\,{+\left[c_{1} -(4b +6)x\right]u_{x} -(4b+6)yu_{y} +2bzu_{z} -b(b+1)u=0,}\\
&{y(1-4y)u_{yy} -4x^{2} u_{xx}}  -z^{2} u_{zz} -8xyu_{xy} +4xzu_{xz} +4yzu_{yz}\\& \,\,\,\,\,\,\,\,\,{-(4b+6)xu_{x} +\left[c_{2} -(4b+6)y\right]u_{y} +2bzu_{z} -b(b+1)u=0,} \\
& {zu_{zz} -2xu_{xz} -2yu_{yz} +(1-b)u_{z} +u=0,}
\end{aligned}
\right.
$

where $u\equiv \,\,{\rm{E}}_{208}\left(b; c_1, c_2; x,y,z\right)$.

Particular solutions:

$
{u_1} ={\rm{E}}_{208}\left(b; c_1, c_2; x,y,z\right),
$

$
{u_2} = {x^{1 - c_1}}{\rm{E}}_{208}\left(2-2c_1+b; 2-c_1, c_2; x,y,z\right),
$

$
{u_3} = {y^{1 - c_2}}{\rm{E}}_{208}\left(2-2c_2+b; c_1, 2-c_2; x,y,z\right),
$

$
{u_4} = {x^{1 - c_1}}{y^{1 - c_2}}{\rm{E}}_{208}\left(4-2c_1-2c_2+b; 2-c_1, 2-c_2; x,y,z\right).
$

\bigskip

\begin{equation}
{\rm{E}}_{209}\left(a_1,a_2,a_3,b_1,b_2;c;x,y,z\right)=\sum\limits_{m,n,p=0}^\infty\frac{(a_1)_n(a_2)_n(a_3)_p(b_1)_{m-n}(b_2)_{m-p}}{(c)_m}\frac{x^m}{m!}\frac{y^n}{n!}\frac{z^p}{p!},
\end{equation}

region of convergence:
$$ \left\{ s(1+r)<1,\,\,\,t<\infty
\right\}.
$$

System of partial differential equations:

$
\left\{
\begin{aligned}
&{x(1-x)u_{xx} +xyu_{xy} +xzu_{xz} -yzu_{yz}} {+\left[c-\left(b_{1} +b_{2} +1\right)x\right]u_{x} +b_{2} yu_{y} +b_{1} zu_{z} -b_{1} b_{2} u=0,} \\& {y(1+y)u_{yy} -xu_{xy} +\left[1-b_{1} +\left(a_{1} +a_{2} +1\right)y\right]u_{y} +a_{1} a_{2} u=0,} \\& {zu_{zz} -xu_{xz} +\left(1-b_{2}+z\right)u_{z} +a_{3} u=0,}
\end{aligned}
\right.
$

where $u\equiv \,\,{\rm{E}}_{209}\left(a_1,a_2,a_3,b_1,b_2;c;x,y,z\right)$.

Particular solutions:

$
{u_1} ={\rm{E}}_{209}\left(a_1,a_2,a_3,b_1,b_2;c;x,y,z\right),
$

$
{u_2} = {x^{1 - c}}{\rm{E}}_{209}\left(a_1,a_2,a_3,1-c+b_1,1-c+b_2;2-c;x,y,z\right).
$

\bigskip

\begin{equation}
{\rm{E}}_{210}\left(a_1,a_2,b_1,b_2;c;x,y,z\right)=\sum\limits_{m,n,p=0}^\infty\frac{(a_1)_n(a_2)_n(b_1)_{m-n}(b_2)_{m-p}}{(c)_m}\frac{x^m}{m!}\frac{y^n}{n!}\frac{z^p}{p!},\,\,s(1+r)<1,
\end{equation}

region of convergence:
$$ \left\{ s(1+r)<1,\,\,\,t<\infty
\right\}.
$$

System of partial differential equations:

$
\left\{
\begin{aligned}
&{x(1-x)u_{xx} +xyu_{xy} +xzu_{xz} -yzu_{yz}} {+\left[c-\left(b_{1} +b_{2} +1\right)x\right]u_{x} +b_{2} yu_{y} +b_{1} zu_{z} -b_{1} b_{2} u=0,} \\& {y(1+y)u_{yy} -xu_{xy} +\left[1-b_{1} +\left(a_{1} +a_{2} +1\right)y\right]u_{y} +a_{1} a_{2} u=0,} \\& {zu_{zz} -xu_{xz} +\left(1-b_{2} \right)u_{z} +u=0,}
\end{aligned}
\right.
$

where $u\equiv \,\,{\rm{E}}_{210}\left(a_1,a_2,b_1,b_2;c;x,y,z\right)$.

Particular solutions:

$
{u_1} ={\rm{E}}_{210}\left(a_1,a_2,b_1,b_2;c;x,y,z\right),
$

$
{u_2} = {x^{1 - c}}{\rm{E}}_{210}\left(a_1,a_2,1-c+b_1,1-c+b_2;2-c;x,y,z\right).
$

\bigskip

\begin{equation}
{\rm{E}}_{211}\left(a_1,a_2,b_1,b_2;c;x,y,z\right)=\sum\limits_{m,n,p=0}^\infty\frac{(a_1)_n(a_2)_p(b_1)_{m-n}(b_2)_{m-p}}{(c)_m}\frac{x^m}{m!}\frac{y^n}{n!}\frac{z^p}{p!},
\end{equation}

region of convergence:
$$ \left\{r<1,\,\,\,s<\infty,\,\,\,t<\infty
\right\}.
$$

System of partial differential equations:

$
\left\{
\begin{aligned}
&{x(1-x)u_{xx} +xyu_{xy} +xzu_{xz} -yzu_{yz}}{+\left[c-\left(b_{1} +b_{2} +1\right)x\right]u_{x} +b_{2} yu_{y} +b_{1} zu_{z} -b_{1} b_{2} u=0,} \\& {yu_{yy} -xu_{xy} +(1-b_{1} +y)u_{y} +a_{1} u=0,} \\& {zu_{zz} -xu_{xz} +\left(1-b_{2}+z \right)u_{z} +a_2u=0,}
\end{aligned}
\right.
$

where $u\equiv \,\,{\rm{E}}_{211}\left(a_1,a_2,b_1,b_2;c;x,y,z\right)$.

Particular solutions:

$
{u_1} ={\rm{E}}_{211}\left(a_1,a_2,b_1,b_2;c;x,y,z\right),
$

$
{u_2} = {x^{1 - c}}{\rm{E}}_{211}\left(a_1,a_2,1-c+b_1,1-c+b_2;2-c;x,y,z\right).
$

\bigskip

\begin{equation}
{\rm{E}}_{212}\left(a,b_1,b_2;c;x,y,z\right)=\sum\limits_{m,n,p=0}^\infty\frac{(a)_n(b_1)_{m-n}(b_2)_{m-p}}{(c)_m}\frac{x^m}{m!}\frac{y^n}{n!}\frac{z^p}{p!},
\end{equation}

region of convergence:
$$ \left\{ r<1,\,\,\,s<\infty,\,\,\,t<\infty
\right\}.
$$

System of partial differential equations:

$
\left\{
\begin{aligned}
&{x(1-x)u_{xx} +xyu_{xy} +xzu_{xz} -yzu_{yz}} {+\left[c-\left(b_{1} +b_{2} +1\right)x\right]u_{x} +b_{2} yu_{y} +b_{1} zu_{z} -b_{1} b_{2} u=0,} \\& {yu_{yy} -xu_{xy} +(1-b_{1} +y)u_{y} +a u=0,} \\& {zu_{zz} -xu_{xz} +\left(1-b_{2} \right)u_{z} +u=0,}
\end{aligned}
\right.
$

where $u\equiv \,\,{\rm{E}}_{212}\left(a,b_1,b_2;c;x,y,z\right)$.

Particular solutions:

$
{u_1} ={\rm{E}}_{212}\left(a,b_1,b_2;c;x,y,z\right),
$

$
{u_2} = {x^{1 - c}}{212}\left(a,1-c+b_1,1-c+b_2;2-c;x,y,z\right).
$

\bigskip

\begin{equation}
{\rm{E}}_{213}\left(a,b;c;x,y,z\right)=\sum\limits_{m,n,p=0}^\infty\frac{(a)_{m-n}(b)_{m-p}}{(c)_m}\frac{x^m}{m!}\frac{y^n}{n!}\frac{z^p}{p!},
\end{equation}

region of convergence:
$$ \left\{ r<1,\,\,\,s<\infty,\,\,\,t<\infty
\right\}.
$$

System of partial differential equations:

$
\left\{
\begin{aligned}
&{x(1-x)u_{xx} +xyu_{xy} +xzu_{xz} -yzu_{yz}}{+\left[c-\left(a +b +1\right)x\right]u_{x} +b yu_{y} +azu_{z} -ab u=0,} \\& {yu_{yy} -xu_{xy} +(1-a )u_{y} +u=0,} \\& {zu_{zz}-xu_{xz}+\left(1-b \right)u_{z} +u=0,}
\end{aligned}
\right.
$

where $u\equiv \,\,{\rm{E}}_{213}\left(a,b;c;x,y,z\right)$.

Particular solutions:

$
{u_1} ={\rm{E}}_{213}\left(a,b;c;x,y,z\right),
$

$
{u_2} = {x^{1 - c}}{213}\left(1-c+a,1-c+b;2-c;x,y,z\right).
$

\bigskip

\begin{equation}
{\rm{E}}_{214}\left(a_1,a_2,a_3,b_1,b_2;c;x,y,z\right)=\sum\limits_{m,n,p=0}^\infty\frac{(a_1)_m(a_2)_n(a_3)_p(b_1)_{m-n}(b_2)_{n-p}}{(c)_m}\frac{x^m}{m!}\frac{y^n}{n!}\frac{z^p}{p!},
\end{equation}

region of convergence:
$$ \left\{ s(1+r)<1,\,\,\,t<\infty
\right\}.
$$

System of partial differential equations:

$
\left\{
\begin{aligned}
&{x(1-x)u_{xx} +xyu_{xy} +\left[c-\left(a_{1} +b_{1} +1\right)x\right]u_{x} +a_{1} yu_{y} -a_{1} b_{1} u=0,} \\& y(1+y)u_{yy} -xu_{xy} -yzu_{yz} +\left[1-b_{1} +\left(a_{2} +b_{2} +1\right)y\right]u_{y} -a_{2} zu_{z} +a_{2} b_{2} u=0, \\& {zu_{zz} -yu_{yz} +\left(1-b_{2}+z \right)u_{z} +a_{3} u=0,}
\end{aligned}
\right.
$

where $u\equiv \,\,{\rm{E}}_{214}\left(a_1,a_2,a_3,b_1,b_2;c;x,y,z\right)$.

Particular solutions:

$
{u_1} ={\rm{E}}_{214}\left(a_1,a_2,a_3,b_1,b_2;c;x,y,z\right),
$

$
{u_2} = {x^{1 - c}}{\rm{E}}_{214}\left(1-c+a_1,a_2,a_3,1-c+b_1,b_2;2-c;x,y,z\right).
$

\bigskip

\begin{equation}
{\rm{E}}_{215}\left(a_1,a_2,a_3,b_1,b_2;c;x,y,z\right)=\sum\limits_{m,n,p=0}^\infty\frac{(a_1)_m(a_2)_p(a_3)_p(b_1)_{m-n}(b_2)_{n-p}}{(c)_m}\frac{x^m}{m!}\frac{y^n}{n!}\frac{z^p}{p!},
\end{equation}

region of convergence:
$$ \left\{ r<1, \,\,\,
, t<1,\,\,\,s<\infty
\right\}.
$$

System of partial differential equations:

$
\left\{
\begin{aligned}
&{x(1-x)u_{xx} +xyu_{xy} +\left[c-\left(a_{1} +b_{1} +1\right)x\right]u_{x} +a_{1} yu_{y} -a_{1} b_{1} u=0,} \\& {yu_{yy} -xu_{xy} +\left(1-b_{1} +y\right)u_{y} -zu_{z} +b_{2} u=0,} \\& {z(1+z)u_{zz} -yu_{yz} +\left[1-b_{2} +\left(a_{2} +a_{3} +1\right)z\right]u_{z} +a_{2} a_{3} u=0,}
\end{aligned}
\right.
$

where $u\equiv \,\,{\rm{E}}_{215}\left(a_1,a_2,a_3,b_1,b_2;c;x,y,z\right)$.

Particular solutions:

$
{u_1} ={\rm{E}}_{215}\left(a_1,a_2,a_3,b_1,b_2;c;x,y,z\right),
$

$
{u_2} = {x^{1 - c}}{\rm{E}}_{215}\left(1-c+a_1,a_2,a_3,1-c+b_1,b_2;2-c;x,y,z\right).
$

\bigskip

\begin{equation}
{\rm{E}}_{216}\left(a_1,a_2,a_3,b_1,b_2;c;x,y,z\right)=\sum\limits_{m,n,p=0}^\infty\frac{(a_1)_n(a_2)_p(a_3)_p(b_1)_{m-n}(b_2)_{n-p}}{(c)_m}\frac{x^m}{m!}\frac{y^n}{n!}\frac{z^p}{p!},
\end{equation}

region of convergence:
$$ \left\{ r<\infty,\,\,\,s<1 \vee t<1
\right\}.
$$

System of partial differential equations:

$
\left\{
\begin{aligned}
&{xu_{xx} +(c-x)u_{x} +yu_{y} -b_{1} u=0,} \\& y(1+y)u_{yy} -xu_{xy} -yzu_{yz} +\left[1-b_{1} +\left(a_{1} +b_{2} +1\right)y\right]u_{y} -a_{1} zu_{z} +a_{1} b_{2} u=0, \\& {z(1+z)u_{zz} -yu_{yz} +\left[1-b_{2} +\left(a_{2} +a_{3} +1\right)z\right]u_{z} +a_{2} a_{3} u=0,}
\end{aligned}
\right.
$

where $u\equiv \,\,{\rm{E}}_{216}\left(a_1,a_2,a_3,b_1,b_2;c;x,y,z\right)$.

Particular solutions:

$
{u_1} ={\rm{E}}_{216}\left(a_1,a_2,a_3,b_1,b_2;c;x,y,z\right),
$

$
{u_2} = {x^{1 - c}}{\rm{E}}_{216}\left(a_1,a_2,a_3,1-c+b_1,b_2;2-c;x,y,z\right).
$

\bigskip

\begin{equation}
{\rm{E}}_{217}\left(a_1,a_2,b_1,b_2;c;x,y,z\right)=\sum\limits_{m,n,p=0}^\infty\frac{(a_1)_m(a_2)_n(b_1)_{m-n}(b_2)_{n-p}}{(c)_m}\frac{x^m}{m!}\frac{y^n}{n!}\frac{z^p}{p!},
\end{equation}

region of convergence:
$$ \left\{ s(1+r)<1,\,\,\,t<\infty
\right\}.
$$

System of partial differential equations:

$
\left\{
\begin{aligned}
&{x(1-x)u_{xx} +xyu_{xy} +\left[c-\left(a_{1} +b_{1} +1\right)x\right]u_{x} +a_{1} yu_{y} -a_{1} b_{1} u=0,} \\
& y(1+y)u_{yy} -xu_{xy} -yzu_{yz}+\left[1-b_{1} +\left(a_{2} +b_{2} +1\right)y\right]u_{y} -a_{2} zu_{z} +a_{2} b_{2} u=0, \\
& {zu_{zz} - yu_{yz} +\left(1-b_{2} \right)u_{z} +u=0,}
\end{aligned}
\right.
$

where $u\equiv \,\,{\rm{E}}_{217}\left(a_1,a_2,b_1,b_2;c;x,y,z\right)$.

Particular solutions:

$
{u_1} ={\rm{E}}_{217}\left(a_1,a_2,b_1,b_2;c;x,y,z\right),
$

$
{u_2} = {x^{1 - c}}{\rm{E}}_{217}\left(1-c+a_1,a_2,1-c+b_1,b_2;2-c;x,y,z\right).
$

\bigskip

\begin{equation}
{\rm{E}}_{218}\left(a_1,a_2,b_1,b_2;c;x,y,z\right)=\sum\limits_{m,n,p=0}^\infty\frac{(a_1)_m(a_2)_p(b_1)_{m-n}(b_2)_{n-p}}{(c)_m}\frac{x^m}{m!}\frac{y^n}{n!}\frac{z^p}{p!},
\end{equation}

region of convergence:
$$ \left\{ r<1,\,\,\,s<\infty,\,\,\,\,t<\infty
\right\}.
$$

System of partial differential equations:

$
\left\{
\begin{aligned}
&{x(1-x)u_{xx} +xyu_{xy} +\left[c-\left(a_{1} +b_{1} +1\right)x\right]u_{x} +a_{1} yu_{y} -a_{1} b_{1} u=0,} \\& {yu_{yy} -xu_{xy} +\left(1-b_{1} +y\right)u_{y} -zu_{z} +b_{2} u=0,} \\& {zu_{zz} -yu_{yz} +\left(1-b_{2}+z\right)u_{z} +a_{2} u=0,}
\end{aligned}
\right.
$

where $u\equiv \,\,{\rm{E}}_{218}\left(a_1,a_2,b_1,b_2;c;x,y,z\right)$.

Particular solutions:

$
{u_1} ={\rm{E}}_{218}\left(a_1,a_2,b_1,b_2;c;x,y,z\right) ,
$

$
{u_2} = {x^{1 - c}}{\rm{E}}_{218}\left(1-c+a_1,a_2,1-c+b_1,b_2;2-c;x,y,z\right).
$

\bigskip

\begin{equation}
{\rm{E}}_{219}\left(a_1,a_2,b_1,b_2;c;x,y,z\right)=\sum\limits_{m,n,p=0}^\infty\frac{(a_1)_n(a_2)_p(b_1)_{m-n}(b_2)_{n-p}}{(c)_m}\frac{x^m}{m!}\frac{y^n}{n!}\frac{z^p}{p!},
\end{equation}

region of convergence:
$$ \left\{ r<\infty,\,\,\,s<1,\,\,\,\,t<\infty
\right\}.
$$

System of partial differential equations:

$
\left\{
\begin{aligned}
&{xu_{xx} +(c-x)u_{x} +yu_{y} -b_{1} u=0,} \\& y(1+y)u_{yy} -xu_{xy} -yzu_{yz} +\left[1-b_{1} +\left(a_{1} +b_{2} +1\right)y\right]u_{y} -a_{1} zu_{z} +a_{1} b_{2} u=0, \\& {zu_{zz} -yu_{yz} +\left(1-b_{2}+z\right)u_{z} +a_{2} u=0,}
\end{aligned}
\right.
$

where $u\equiv \,\,{\rm{E}}_{219}\left(a_1,a_2,b_1,b_2;c;x,y,z\right)$.

Particular solutions:

$
{u_1} ={\rm{E}}_{219}\left(a_1,a_2,b_1,b_2;c;x,y,z\right),
$

$
{u_2} = {x^{1 - c}}{\rm{E}}_{219}\left(a_1,a_2,1-c+b_1,b_2;2-c;x,y,z\right).
$

\bigskip

\begin{equation}
{\rm{E}}_{220}\left(a_1,a_2,b_1,b_2;c;x,y,z\right)=\sum\limits_{m,n,p=0}^\infty\frac{(a_1)_p(a_2)_p(b_1)_{m-n}(b_2)_{n-p}}{(c)_m}\frac{x^m}{m!}\frac{y^n}{n!}\frac{z^p}{p!},
\end{equation}

region of convergence:
$$ \left\{ r<\infty,\,\,\,s<\infty,\,\,\,\,t<1
\right\}.
$$

System of partial differential equations:

$
\left\{
\begin{aligned}
&{xu_{xx} +(c-x)u_{x} +yu_{y} -b_{1} u=0,} \\& {yu_{yy} -xu_{xy} +\left(1-b_{1} +y\right)u_{y} -zu_{z} +b_{2} u=0,} \\& {z(1+z)u_{zz} -yu_{yz} +\left[1-b_{2} +\left(a_{1} +a_{2} +1\right)z\right]u_{z} +a_{1} a_{2} u=0,}
\end{aligned}
\right.
$

where $u\equiv \,\,{\rm{E}}_{220}\left(a_1,a_2,b_1,b_2;c;x,y,z\right)$.

Particular solutions:

$
{u_1} ={\rm{E}}_{220}\left(a_1,a_2,b_1,b_2;c;x,y,z\right) ,
$

$
{u_2} = {x^{1 - c}}{\rm{E}}_{220}\left(a_1,a_2,1-c+b_1,b_2;2-c;x,y,z\right).
$

\bigskip

\begin{equation}
{\rm{E}}_{221}\left(a,b_1,b_2;c;x,y,z\right)=\sum\limits_{m,n,p=0}^\infty\frac{(a)_m(b_1)_{m-n}(b_2)_{n-p}}{(c)_m}\frac{x^m}{m!}\frac{y^n}{n!}\frac{z^p}{p!},
\end{equation}

region of convergence:
$$ \left\{ r<1,\,\,\,s<\infty,\,\,\,\,t<\infty
\right\}.
$$

System of partial differential equations:

$
\left\{
\begin{aligned}
&{x(1-x)u_{xx} +xyu_{xy} +\left[c-\left(a+b_{1} +1\right)x\right]u_{x} +ayu_{y} -ab_{1} u=0,} \\
& {yu_{yy} -xu_{xy} +\left(1-b_{1} +y\right)u_{y} -zu_{z} +b_{2} u=0,} \\
& {zu_{zz} - yu_{yz} +\left(1-b_{2}\right)u_{z} +u=0,}
\end{aligned}
\right.
$

where $u\equiv \,\,{\rm{E}}_{221}\left(a,b_1,b_2;c;x,y,z\right)$.

Particular solutions:

$
{u_1} ={\rm{E}}_{221}\left(a,b_1,b_2;c;x,y,z\right) ,
$

$
{u_2} = {x^{1 - c}}{\rm{E}}_{221}\left(1-c+a,1-c+b_1,b_2;2-c;x,y,z\right).
$

\bigskip

\begin{equation}
{\rm{E}}_{222}\left(a,b_1,b_2;c;x,y,z\right)=\sum\limits_{m,n,p=0}^\infty\frac{(a)_n(b_1)_{m-n}(b_2)_{n-p}}{(c)_m}\frac{x^m}{m!}\frac{y^n}{n!}\frac{z^p}{p!},
\end{equation}

region of convergence:
$$ \left\{ r<\infty,\,\,\,s<1,\,\,\,\,t<\infty
\right\}.
$$

System of partial differential equations:

$
\left\{
\begin{aligned}
&{xu_{xx} +(c-x)u_{x} +yu_{y} -b_{1} u=0,} \\
& y(1+y)u_{yy} -xu_{xy} -yzu_{yz}  +\left[1-b_{1}+\left(a+b_{2} +1\right)y\right]u_{y} -azu_{z} +ab_{2} u=0, \\
& {zu_{zz} - yu_{yz} +\left(1-b_{2}\right)u_{z} +u=0,}
\end{aligned}
\right.
$

where $u\equiv \,\,{\rm{E}}_{222}\left(a,b_1,b_2;c;x,y,z\right)$.

Particular solutions:

$
{u_1} ={\rm{E}}_{222}\left(a,b_1,b_2;c;x,y,z\right),
$

$
{u_2} = {x^{1 - c}}{\rm{E}}_{222}\left(a,1-c+b_1,b_2;2-c;x,y,z\right).
$

\bigskip

\begin{equation}
{\rm{E}}_{223}\left(a,b_1,b_2;c;x,y,z\right)=\sum\limits_{m,n,p=0}^\infty\frac{(a)_p(b_1)_{m-n}(b_2)_{n-p}}{(c)_m}\frac{x^m}{m!}\frac{y^n}{n!}\frac{z^p}{p!},
\end{equation}

region of convergence:
$$ \left\{ r<\infty,\,\,\,s<\infty,\,\,\,\,t<\infty
\right\}.
$$

System of partial differential equations:

$
\left\{
\begin{aligned}
&{xu_{xx} +(c-x)u_{x} +yu_{y} -b_{1} u=0,} \\& {yu_{yy} -xu_{xy} +\left(1-b_{1} +y\right)u_{y} -zu_{z} +b_{2} u=0,} \\& {zu_{zz} -yu_{yz} +\left(1-b_{2}+z\right)u_{z} +au=0,}
\end{aligned}
\right.
$

where $u\equiv \,\,{\rm{E}}_{223}\left(a,b_1,b_2;c;x,y,z\right)$.

Particular solutions:

$
{u_1} ={\rm{E}}_{223}\left(a,b_1,b_2;c;x,y,z\right) ,
$

$
{u_2} = {x^{1 - c}}{\rm{E}}_{223}\left(a,1-c+b_1,b_2;2-c;x,y,z\right).
$

\bigskip

\begin{equation}
{\rm{E}}_{224}\left(a,b;c;x,y,z\right)=\sum\limits_{m,n,p=0}^\infty\frac{(a)_{m-n}(b)_{n-p}}{(c)_m}\frac{x^m}{m!}\frac{y^n}{n!}\frac{z^p}{p!},
\end{equation}

region of convergence:
$$ \left\{ r<\infty,\,\,\,s<\infty,\,\,\,\,t<\infty
\right\}.
$$

System of partial differential equations:

$
\left\{
\begin{aligned}
&{xu_{xx} +(c-x)u_{x} +yu_{y} -a u=0,} \\
& {yu_{yy} -xu_{xy} +\left(1-a +y\right)u_{y} -zu_{z} +b u=0,} \\
& {zu_{zz} - yu_{yz} +\left(1-b\right)u_{z} +u=0,}
\end{aligned}
\right.
$

where $u\equiv \,\,{\rm{E}}_{224}\left(a,b;c;x,y,z\right)$.

Particular solutions:

$
{u_1} ={\rm{E}}_{224}\left(a,b;c;x,y,z\right),
$

$
{u_2} = {x^{1 - c}}{\rm{E}}_{224}\left(1-c+a,b;2-c;x,y,z\right).
$

\bigskip

\begin{equation}
{\rm{E}}_{225}\left(a_1,a_2,a_3,a_4,b,c;x,y,z\right)=\sum\limits_{m,n,p=0}^\infty\frac{(a_1)_m(a_2)_n(a_3)_p(a_4)_p(b)_{m-n-p}}{(c)_{m}}\frac{x^m}{m!}\frac{y^n}{n!}\frac{z^p}{p!},
\end{equation}

region of convergence:
$$ \left\{ t(1+r)<1,\,\,\,s<\infty
\right\}.
$$

System of partial differential equations:

$
\left\{
\begin{aligned}
&x(1-x)u_{xx} +xyu_{xy} +xzu_{xz}+\left[c-\left(a_{1} +b+1\right)x\right]u_{x} +a_{1} yu_{y} +a_{1} zu_{z} -a_{1} bu=0, \\& {yu_{yy} -xu_{xy} +zu_{yz} +\left(1-b +y\right)u_{y} +a_{2} u=0,} \\& {z(1+z)u_{zz} -xu_{xz} +yu_{yz} +\left[1-b+\left(a_{3} +a_{4} +1\right)z\right]u_{z} +a_{3} a_{4} u=0,}
\end{aligned}
\right.
$

where $u\equiv \,\,{\rm{E}}_{225}\left(a_1,a_2,a_3,a_4,b,c;x,y,z\right)$.

Particular solutions:

$
{u_1} ={\rm{E}}_{225}\left(a_1,a_2,a_3,a_4,b,c;x,y,z\right),
$

$
{u_2} = {x^{1 - c}}{\rm{E}}_{225}\left(1-c+a_1,a_2,a_3,a_4,1-c+b,2-c;x,y,z\right).
$

\bigskip

\begin{equation}
{\rm{E}}_{226}\left(a_1,a_2,a_3,a_4,b,c;x,y,z\right)=\sum\limits_{m,n,p=0}^\infty\frac{(a_1)_n(a_2)_n(a_3)_p(a_4)_p(b)_{m-n-p}}{(c)_{m}}\frac{x^m}{m!}\frac{y^n}{n!}\frac{z^p}{p!},
\end{equation}

region of convergence:
$$ \left\{ r<\infty,\,\,\,\frac{1}{s}+\frac{1}{t}>1
\right\}.
$$

System of partial differential equations:

$
\left\{
\begin{aligned}
&{xu_{xx} +(c-x)u_{x} +yu_{y} +zu_{z} -bu=0,} \\& {y(1+y)u_{yy} -xu_{xy} +zu_{yz} +\left[1-b+\left(a_{1} +a_{2} +1\right)y\right]u_{y} +a_{1} a_{2} u=0,} \\& {z(1+z)u_{zz} -xu_{xz}+yu_{yz} +\left[1-b +\left(a_{3} +a_{4} +1\right)z\right]u_{z} +a_{3} a_{4} u=0,}
\end{aligned}
\right.
$

where $u\equiv \,\,{\rm{E}}_{226}\left(a_1,a_2,a_3,a_4,b,c;x,y,z\right)$.

Particular solutions:

$
{u_1} ={\rm{E}}_{226}\left(a_1,a_2,a_3,a_4,b,c;x,y,z\right) ,
$

$
{u_2} = {x^{1 - c}}{\rm{E}}_{226}\left(a_1,a_2,a_3,a_4,1-c+b,2-c;x,y,z\right).
$

\bigskip

\begin{equation}
{\rm{E}}_{227}\left(a_1,a_2,a_3,b,c;x,y,z\right)=\sum\limits_{m,n,p=0}^\infty\frac{(a_1)_m(a_2)_n(a_3)_p(b)_{m-n-p}}{(c)_{m}}\frac{x^m}{m!}\frac{y^n}{n!}\frac{z^p}{p!},
\end{equation}

region of convergence:
$$ \left\{ r<1,\,\,\,s<\infty,\,\,\,t<\infty
\right\}.
$$

System of partial differential equations:

$
\left\{
\begin{aligned}
&x(1-x)u_{xx} +xyu_{xy} +xzu_{xz} +\left[c-\left(a_{1} +b+1\right)x\right]u_{x}+a_{1} yu_{y} +a_{1} zu_{z} -a_{1} bu=0, \\ &{yu_{yy} -xu_{xy} +zu_{yz} +\left(1-b +y\right)u_{y} +a_{2} u=0,} \\ &{zu_{zz} -xu_{xz} +yu_{yz} +(1-b+z)u_{z} +a_{3} u=0,}
\end{aligned}
\right.
$

where $u\equiv \,\,{\rm{E}}_{227}\left(a_1,a_2,a_3,b,c;x,y,z\right)$.

Particular solutions:

$
{u_1} ={\rm{E}}_{227}\left(a_1,a_2,a_3,b,c;x,y,z\right) ,
$

$
{u_2} = {x^{1 - c}}{\rm{E}}_{227}\left(1-c+a_1,a_2,a_3,1-c+b,2-c;x,y,z\right).
$

\bigskip

\begin{equation}
{\rm{E}}_{228}\left(a_1,a_2,a_3,b,c;x,y,z\right)=\sum\limits_{m,n,p=0}^\infty\frac{(a_1)_m(a_2)_p(a_3)_p(b)_{m-n-p}}{(c)_{m}}\frac{x^m}{m!}\frac{y^n}{n!}\frac{z^p}{p!},
\end{equation}

region of convergence:
$$ \left\{ r<\infty,\,\,\,\frac{1}{s}+\frac{1}{t}>1
\right\}.
$$

System of partial differential equations:

$
\left\{
\begin{aligned}
&x(1-x)u_{xx} +xyu_{xy} +xzu_{xz}+\left[c-\left(a_{1} +b+1\right)x\right]u_{x} +a_{1} yu_{y} +a_{1} zu_{z} -a_{1} bu=0, \\& {yu_{yy} -xu_{xy} +zu_{yz} +(1-b)u_{y} +u=0,} \\& {z(1+z)u_{zz} -xu_{xz} +yu_{yz} +\left[1-b+\left(a_{2} +a_{3} +1\right)z\right]u_{z} +a_{2} a_{3} u=0,}
\end{aligned}
\right.
$

where $u\equiv \,\,{\rm{E}}_{228}\left(a_1,a_2,a_3,b,c;x,y,z\right)$.

Particular solutions:

$
{u_1} ={\rm{E}}_{228}\left(a_1,a_2,a_3,b,c;x,y,z\right),
$

$
{u_2} = {x^{1 - c}}{\rm{E}}_{228}\left(1-c+a_1,a_2,a_3,1-c+b,2-c;x,y,z\right).
$

\bigskip

\begin{equation}
{\rm{E}}_{229}\left(a_1,a_2,a_3,b,c;x,y,z\right)=\sum\limits_{m,n,p=0}^\infty\frac{(a_1)_n(a_2)_p(a_3)_p(b)_{m-n-p}}{(c)_{m}}\frac{x^m}{m!}\frac{y^n}{n!}\frac{z^p}{p!},
\end{equation}

region of convergence:
$$ \left\{ r<\infty,\,\,\,s<\infty,\,\,\,t<1
\right\}.
$$

System of partial differential equations:

$
\left\{
\begin{aligned}
&{xu_{xx} +(c-x)u_{x} +yu_{y} +zu_{z} -bu=0,} \\& {yu_{yy} -xu_{xy} +zu_{yz} +\left(1-b +y\right)u_{y} +a_{1} u=0,} \\& {z(1+z)u_{zz} -xu_{xz} +yu_{yz} +\left[1-b+\left(a_{2} +a_{3} +1\right)z\right]u_{z} +a_{2} a_{3} u=0,}
\end{aligned}
\right.
$

where $u\equiv \,\,{\rm{E}}_{229}\left(a_1,a_2,a_3,b,c;x,y,z\right)$.

Particular solutions:

$
{u_1} ={\rm{E}}_{229}\left(a_1,a_2,a_3,b,c;x,y,z\right),
$

$
{u_2} = {x^{1 - c}}{\rm{E}}_{229}\left(a_1,a_2,a_3,1-c+b,2-c;x,y,z\right).
$

\bigskip

\begin{equation}
{\rm{E}}_{230}\left(a_1,a_2,b,c;x,y,z\right)=\sum\limits_{m,n,p=0}^\infty\frac{(a_1)_m(a_2)_p(b)_{m-n-p}}{(c)_{m}}\frac{x^m}{m!}\frac{y^n}{n!}\frac{z^p}{p!},
\end{equation}

region of convergence:
$$ \left\{ r<1,\,\,\,s<\infty,\,\,\,t<\infty
\right\}.
$$

System of partial differential equations:

$
\left\{
\begin{aligned}
&x(1-x)u_{xx} +xyu_{xy} +xzu_{xz} +\left[c-\left(a_{1} +b+1\right)x\right]u_{x}  +a_{1} yu_{y} +a_{1} zu_{z} -a_{1} bu=0, \\ &{yu_{yy} -xu_{xy} +zu_{yz} +(1-b)u_{y} +u=0,} \\& {zu_{zz} -xu_{xz} +yu_{yz} +(1-b+z)u_{z} +a_{2} u=0,}
\end{aligned}
\right.
$

where $u\equiv \,\,{\rm{E}}_{230}\left(a_1,a_2,b,c;x,y,z\right)$.

Particular solutions:

$
{u_1} ={\rm{E}}_{230}\left(a_1,a_2,b,c;x,y,z\right) ,
$

$
{u_2} = {x^{1 - c}}{\rm{E}}_{230}\left(1-c+a_1,a_2,1-c+b,2-c;x,y,z\right).
$

\bigskip

\begin{equation}
{\rm{E}}_{231}\left(a_1,a_2,b,c;x,y,z\right)=\sum\limits_{m,n,p=0}^\infty\frac{(a_1)_n(a_2)_p(b)_{m-n-p}}{(c)_{m}}\frac{x^m}{m!}\frac{y^n}{n!}\frac{z^p}{p!},
\end{equation}

region of convergence:
$$ \left\{ r<\infty,\,\,\,s<\infty,\,\,\,t<\infty
\right\}.
$$

System of partial differential equations:

$
\left\{
\begin{aligned}
&{xu_{xx} +(c-x)u_{x} +yu_{y} +zu_{z} -bu=0,} \\& {yu_{yy} -xu_{xy} +zu_{yz} +\left(1-b +y\right)u_{y} +a_{1} u=0,} \\& {zu_{zz} -xu_{xz} +yu_{yz} +(1-b+z)u_{z} +a_{2} u=0,}
\end{aligned}
\right.
$

where $u\equiv \,\,{\rm{E}}_{231}\left(a_1,a_2,b,c;x,y,z\right)$.

Particular solutions:

$
{u_1} ={\rm{E}}_{231}\left(a_1,a_2,b,c;x,y,z\right) ,
$

$
{u_2} = {x^{1 - c}}{\rm{E}}_{231}\left(a_1,a_2,1-c+b,2-c;x,y,z\right).
$

\bigskip

\begin{equation}
{\rm{E}}_{232}\left(a_1,a_2,b,c;x,y,z\right)=\sum\limits_{m,n,p=0}^\infty\frac{(a_1)_p(a_2)_p(b)_{m-n-p}}{(c)_{m}}\frac{x^m}{m!}\frac{y^n}{n!}\frac{z^p}{p!},
\end{equation}

region of convergence:
$$ \left\{ r<\infty,\,\,\,s<\infty,\,\,\,t<1
\right\}.
$$

System of partial differential equations:

$
\left\{
\begin{aligned}
&{xu_{xx} +(c-x)u_{x} +yu_{y} +zu_{z}-bu=0,} \\& {yu_{yy} -xu_{xy} +zu_{yz} +(1-b)u_{y} +u=0,} \\& {z(1+z)u_{zz} -xu_{xz} +yu_{yz} +\left[1-b+\left(a_{1} +a_{2} +1\right)z\right]u_{z} +a_{1} a_{2} u=0,}
\end{aligned}
\right.
$

where $u\equiv \,\,{\rm{E}}_{232}\left(a_1,a_2,b,c;x,y,z\right)$.

Particular solutions:

$
{u_1} ={\rm{E}}_{232}\left(a_1,a_2,b,c;x,y,z\right) ,
$

$
{u_2} = {x^{1 - c}}{\rm{E}}_{232}\left(a_1,a_2,1-c+b,2-c;x,y,z\right).
$

\bigskip

\begin{equation}
{\rm{E}}_{233}\left(a,b,c;x,y,z\right)=\sum\limits_{m,n,p=0}^\infty\frac{(a)_p(b)_{m-n-p}}{(c)_{m}}\frac{x^m}{m!}\frac{y^n}{n!}\frac{z^p}{p!},
\end{equation}

region of convergence:
$$ \left\{ r<\infty,\,\,\,s<\infty,\,\,\,t<\infty
\right\}.
$$

System of partial differential equations:

$
\left\{
\begin{aligned}
&{xu_{xx} +(c-x)u_{x} +yu_{y} +zu_{z} -bu=0,} \\& {yu_{yy} -xu_{xy} +zu_{yz} +(1-b)u_{y} +u=0,} \\& {zu_{zz} -xu_{xz} +yu_{yz} +(1-b+z)u_{z} +au=0,}
\end{aligned}
\right.
$

where $u\equiv \,\,{\rm{E}}_{233}\left(a,b,c;x,y,z\right)$.

Particular solutions:

$
{u_1} ={\rm{E}}_{233}\left(a,b,c;x,y,z\right) ,
$

$
{u_2} = {x^{1 - c}}{\rm{E}}_{233}\left(a,1-c+b,2-c;x,y,z\right).
$

\bigskip

\begin{equation}
{\rm{E}}_{234}\left(a_1,a_2,b_1,b_2; c; x,y,z\right)=\sum\limits_{m,n,p=0}^\infty\frac{(a_1)_n(a_2)_p(b_1)_{m+n-p}(b_2)_{m-n}}{(c)_m}\frac{x^m}{m!}\frac{y^n}{n!}\frac{z^p}{p!},
\end{equation}

region of convergence:
$$ \left\{ \left\{s+2\sqrt{rs}<1\right\}=\left\{\sqrt{s}<\sqrt{1+r}-\sqrt{r}\right\},\,\,\,t<\infty
\right\}.
$$

System of partial differential equations:

$
\left\{
\begin{aligned}
&x(1-x)u_{xx} +y^{2} u_{yy}-yzu_{yz} +xzu_{xz}  \\
&\,\,\,\,\,\,\,\,\,\,\,\,\, +\left[c-\left(b_{1} +b_{2} +1\right)x\right]u_{x}{ +\left(b_{1} -b_{2} +1\right)yu_y+b_{2} zu_{z} -b_{1} b_{2} u=0,} \\
& y(1+y)u_{yy} -x(1-y)u_{xy} -yzu_{yz} +a_{1} xu_{x} \\
&\,\,\,\,\,\,\,\,\,\,\,\,\, +\left[1-b_{2} +\left(a_{1} +b_{1} +1\right)y\right]u_{y}-a_1zu_z+a_{1} b_{1} u=0, \\
& {zu_{zz} -xu_{xz} -yu_{yz}  +\left(1-b_{1}+z\right)u_{z} +a_{2} u=0,}
\end{aligned}
\right.
$

where $u\equiv \,\,{\rm{E}}_{234}\left(a_1,a_2,b_1,b_2; c; x,y,z\right)$.

Particular solutions:

$
{u_1} ={\rm{E}}_{234}\left(a_1,a_2,b_1,b_2; c; x,y,z\right) ,
$

$
{u_2} = {x^{1 - c}}{\rm{E}}_{234}\left(a_1,a_2,1-c+b_1,1-c+b_2; 2-c; x,y,z\right).
$

\bigskip

\begin{equation}
{\rm{E}}_{235}\left(a_1,a_2,b_1,b_2; c; x,y,z\right)=\sum\limits_{m,n,p=0}^\infty\frac{(a_1)_p(a_2)_p(b_1)_{m+n-p}(b_2)_{m-n}}{(c)_m}\frac{x^m}{m!}\frac{y^n}{n!}\frac{z^p}{p!},
\end{equation}

region of convergence:
$$ \left\{ t(1+r)<1,\,\,\,s<\infty
\right\}.
$$

System of partial differential equations:

$
\left\{
\begin{aligned}
&{x(1-x)u_{xx} +y^{2} u_{yy} -yzu_{yz} +xzu_{xz} }\\& \,\,\,\,\,\,\,\,\,+\left[c-\left(b_{1} +b_{2} +1\right)x\right]u_{x} +(b_{1} -b_{2} +1)yu_y+b_{2} zu_{z} -b_{1} b_{2} u=0, \\
& {yu_{yy} -xu_{xy} +xu_{x} +(1-b_{2} +y)u_{y} -zu_{z} +b_{1} u=0,} \\
& {z(1+z)u_{zz} -xu_{xz} -yu_{yz}}{+\left[1-b_{1}+\left(a_1+a_2+1 \right)z\right]u_{z} +a_1a_{2} u=0,}
\end{aligned}
\right.
$

where $u\equiv \,\,{\rm{E}}_{235}\left(a_1,a_2,b_1,b_2; c; x,y,z\right)$.

Particular solutions:

$
{u_1} ={\rm{E}}_{235}\left(a_1,a_2,b_1,b_2; c; x,y,z\right) ,
$

$
{u_2} = {x^{1 - c}}{\rm{E}}_{235}\left(a_1,a_2,1-c+b_1,1-c+b_2; 2-c; x,y,z\right).
$

\bigskip

\begin{equation}
{\rm{E}}_{236}\left(a,b_1,b_2; c; x,y,z\right)=\sum\limits_{m,n,p=0}^\infty\frac{(a)_n(b_1)_{m+n-p}(b_2)_{m-n}}{(c)_m}\frac{x^m}{m!}\frac{y^n}{n!}\frac{z^p}{p!},
\end{equation}

region of convergence:
$$ \left\{ \left\{s+2\sqrt{rs}<1\right\}=\left\{\sqrt{s}<\sqrt{1+r}-\sqrt{r}\right\},\,\,\,t<\infty
\right\}.
$$

System of partial differential equations:

$
\left\{
\begin{aligned}
&{x(1-x)u_{xx} +y^{2} u_{yy} -yzu_{yz} +xzu_{xz} }\\& \,\,\,\,\,\,\,\,\,+\left[c-\left(b_{1} +b_{2} +1\right)x\right]u_{x} +(b_{1} -b_{2} +1)yu_y+b_{2} zu_{z} -b_{1} b_{2} u=0, \\
& y(1+y)u_{yy} -x(1-y)u_{xy} -yzu_{yz} +a xu_{x} +\left[1-b_{2} +(a +b_{1} +1)y\right]u_{y}-azu_z +a b_{1} u=0, \\
& {zu_{zz} -xu_{xz} -yu_{yz} +(1-b_{1} )u_{z} +u=0,}
\end{aligned}
\right.
$

where $u\equiv \,\,{\rm{E}}_{236}\left(a,b_1,b_2; c; x,y,z\right)$.

Particular solutions:

$
{u_1} ={\rm{E}}_{236}\left(a,b_1,b_2; c; x,y,z\right),
$

$
{u_2} = {x^{1 - c}}{\rm{E}}_{236}\left(a,1-c+b_1,1-c+b_2; 2-c; x,y,z\right).
$

\bigskip

\begin{equation}
{\rm{E}}_{237}\left(a,b_1,b_2; c; x,y,z\right)=\sum\limits_{m,n,p=0}^\infty\frac{(a)_p(b_1)_{m+n-p}(b_2)_{m-n}}{(c)_m}\frac{x^m}{m!}\frac{y^n}{n!}\frac{z^p}{p!},
\end{equation}

region of convergence:
$$ \left\{ r<1,\,\,\,s<\infty,\,\,\,t<\infty
\right\}.
$$

System of partial differential equations:

$
\left\{
\begin{aligned}
&x(1-x)u_{xx} +y^{2} u_{yy} +xzu_{xz} -yzu_{yz} \\& \,\,\,\,\,\,\,\,\,+\left[c-\left(b_{1} +b_{2} +1\right)x\right]u_{x} +(b_{1} -b_{2} +1)yu_y+b_{2} zu_{z} -b_{1} b_{2} u=0, \\&
 {yu_{yy} -xu_{xy} +xu_{x} } {+\left(1-b_{2} +y\right)u_{y}-zu_z + b_{1} u=0,} \\&
  {zu_{zz} -xu_{xz} -yu_{yz} +\left(1-b_{1}+z \right)u_{z} +au=0,}
\end{aligned}
\right.
$

where $u\equiv \,\,{\rm{E}}_{237}\left(a,b_1,b_2; c; x,y,z\right)$.

Particular solutions:

$
{u_1} ={\rm{E}}_{237}\left(a,b_1,b_2; c; x,y,z\right) ,
$

$
{u_2} = {x^{1 - c}}{\rm{E}}_{237}\left(a,1-c+b_1,1-c+b_2; 2-c; x,y,z\right).
$

\bigskip

\begin{equation}
{\rm{E}}_{238}\left(a,b; c; x,y,z\right)=\sum\limits_{m,n,p=0}^\infty\frac{(a)_{m+n-p}(b)_{m-n}}{(c)_m}\frac{x^m}{m!}\frac{y^n}{n!}\frac{z^p}{p!},
\end{equation}

region of convergence:
$$ \left\{ r<1,\,\,\,s<\infty,\,\,\,t<\infty
\right\}.
$$

System of partial differential equations:

$
\left\{
\begin{aligned}
&x(1-x)u_{xx} +y^{2} u_{yy} +xzu_{xz} -yzu_{yz} \\& \,\,\,\,\,\,\,\,\, +\left[c-\left(a +b +1\right)x\right]u_{x} +(a -b+1)yu_y+b zu_{z} -a b u=0, \\&
 {yu_{yy} -xu_{xy} +xu_{x} } {+\left(1-b +y\right)u_{y}-zu_z + a u=0,} \\&
 {zu_{zz} -xu_{xz} -yu_{yz} +\left(1-a \right)u_{z} +u=0,}
\end{aligned}
\right.
$

where $u\equiv \,\,{\rm{E}}_{238}\left(a,b; c; x,y,z\right)$.

Particular solutions:

$
{u_1} ={\rm{E}}_{238}\left(a,b; c; x,y,z\right) ,
$

$
{u_2} = {x^{1 - c}}{\rm{E}}_{238}\left(1-c+a,1-c+b; 2-c; x,y,z\right).
$

\bigskip

\begin{equation}
{\rm{E}}_{239}\left(a_1,a_2,b_1,b_2; c; x,y,z\right)=\sum\limits_{m,n,p=0}^\infty\frac{(a_1)_m(a_2)_n(b_1)_{m+n-p}(b_2)_{p-n}}{(c)_m}\frac{x^m}{m!}\frac{y^n}{n!}\frac{z^p}{p!},
\end{equation}

region of convergence:
$$ \left\{ r+s<1,\,\,\,\,t<\infty
\right\}.
$$

System of partial differential equations:

$
\left\{
\begin{aligned}
&x(1-x)u_{xx} -xyu_{xy} +xzu_{xz}+\left[c-\left(a_{1} +b_{1} +1\right)x\right]u_{x} -a_{1} yu_{y} +a_{1} zu_{z} -a_{1} b_{1} u=0, \\
& y(1+y)u_{yy} +xyu_{xy} -(y+1)zu_{yz} \\& \,\,\,\,\,\,\,\,\,+a_{2} xu_{x}+\left[1-b_{2} +\left(a_{2} +b_{1} +1\right)y\right]u_{y}-a_2zu_z +a_{2} b_{1} u=0, \\
& {zu_{zz} -xu_{xz} -yu_{yz} -yu_{y} +\left(1-b_{1} +z\right)u_{z} +b_{2} u=0,}
\end{aligned}
\right.
$

where $u\equiv \,\,{\rm{E}}_{239}\left(a_1,a_2,b_1,b_2; c; x,y,z\right)$.

Particular solutions:

$
{u_1} ={\rm{E}}_{239}\left(a_1,a_2,b_1,b_2; c; x,y,z\right) ,
$

$
{u_2} = {x^{1 - c}}{\rm{E}}_{239}\left(1-c+a_1,a_2,1-c+b_1,b_2; 2-c; x,y,z\right).
$

\bigskip

\begin{equation}
{\rm{E}}_{240}\left(a_1,a_2,b_1,b_2; c; x,y,z\right)=\sum\limits_{m,n,p=0}^\infty\frac{(a_1)_m(a_2)_p(b_1)_{m+n-p}(b_2)_{p-n}}{(c)_m}\frac{x^m}{m!}\frac{y^n}{n!}\frac{z^p}{p!},
\end{equation}

region of convergence:
$$ \left\{ t(1+r)<1,\,\,\,s<\infty
\right\}.
$$

System of partial differential equations:

$
\left\{
\begin{aligned}
&x(1-x)u_{xx} -xyu_{xy} +xzu_{xz} +\left[c-\left(a_{1} +b_{1} +1\right)x\right]u_{x} -a_{1} yu_y+a_{1} zu_{z} -a_{1} b_{1} u=0, \\& {yu_{yy} -zu_{yz} +xu_{x} +\left(1-b_{2} +y\right)u_{y} -zu_{z} +b_{1} u=0,} \\& z(1+z)u_{zz} -xu_{xz} -y(1+z)u_{yz} -a_{2} yu_{y} +\left[1-b_{1} +\left(a_{2} +b_{2} +1\right)z\right]u_{z} +a_{2} b_{2} u=0,
\end{aligned}
\right.
$

where $u\equiv \,\,{\rm{E}}_{240}\left(a_1,a_2,b_1,b_2; c; x,y,z\right)$.

Particular solutions:

$
{u_1} ={\rm{E}}_{240}\left(a_1,a_2,b_1,b_2; c; x,y,z\right) ,
$

$
{u_2} = {x^{1 - c}}{\rm{E}}_{240}\left(1-c+a_1,a_2,1-c+b_1,b_2; 2-c; x,y,z\right).
$

\bigskip

\begin{equation}
{\rm{E}}_{241}\left(a_1,a_2,b_1,b_2; c; x,y,z\right)=\sum\limits_{m,n,p=0}^\infty\frac{(a_1)_n(a_2)_p(b_1)_{m+n-p}(b_2)_{p-n}}{(c)_m}\frac{x^m}{m!}\frac{y^n}{n!}\frac{z^p}{p!},
\end{equation}

region of convergence:
$$ \left\{ r<\infty,\,\,\,s<1 \vee t<1
\right\}.
$$

System of partial differential equations:

$
\left\{
\begin{aligned}
&{xu_{xx} +(c-x)u_{x} -yu_{y} +zu_{z} -b_{1} u=0,} \\
& y(1+y)u_{yy} +xyu_{xy} -(y+1)zu_{yz}\\& \,\,\,\,\,\,\,\,\, +a_{1} xu_{x}  {+\left[1-b_{2} +\left(a_{1} +b_{1} +1\right)y\right]u_{y}-a_1zu_z +a_{1} b_{1} u=0,} \\
& z(1+z)u_{zz} -xu_{xz} -y(1+z)u_{yz} -a_{2} yu_{y}+\left[1-b_{1} +\left(a_{2} +b_{2} +1\right)z\right]u_{z} +a_{2} b_{2} u=0,
\end{aligned}
\right.
$

where $u\equiv \,\,{\rm{E}}_{241}\left(a_1,a_2,b_1,b_2; c; x,y,z\right)$.

Particular solutions:

$
{u_1} ={\rm{E}}_{241}\left(a_1,a_2,b_1,b_2; c; x,y,z\right) ,
$

$
{u_2} = {x^{1 - c}}{\rm{E}}_{241}\left(a_1,a_2,1-c+b_1,b_2; 2-c; x,y,z\right).
$

\bigskip

\begin{equation}
{\rm{E}}_{242}\left(a,b_1,b_2; c; x,y,z\right)=\sum\limits_{m,n,p=0}^\infty\frac{(a)_m(b_1)_{m+n-p}(b_2)_{p-n}}{(c)_m}\frac{x^m}{m!}\frac{y^n}{n!}\frac{z^p}{p!},
\end{equation}

region of convergence:
$$ \left\{ r<1,\,\,\,s<\infty,\,\,\,t<\infty
\right\}.
$$

System of partial differential equations:

$
\left\{
\begin{aligned}
&{x(1-x)u_{xx} -xyu_{xy} +xzu_{xz} +\left[c-\left(a +b_{1} +1\right)x\right]u_{x} }  {-a yu_{y} +a zu_{z} -a b_{1} u=0,} \\& {yu_{yy} -zu_{yz} +xu_{x} +\left(1-b_{2} +y\right)u_{y} -zu_{z} +b_{1} u=0,} \\& {zu_{zz} -xu_{xz} -yu_{yz} -yu_{y} +\left(1-b_{1} +z\right)u_{z} +b_{2} u=0,}
\end{aligned}
\right.
$

where $u\equiv \,\,{\rm{E}}_{242}\left(a,b_1,b_2; c; x,y,z\right)$.

Particular solutions:

$
{u_1} ={\rm{E}}_{242}\left(a,b_1,b_2; c; x,y,z\right) ,
$

$
{u_2} = {x^{1 - c}}{\rm{E}}_{242}\left(1-c+a,1-c+b_1,b_2; 2-c; x,y,z\right).
$

\bigskip

\begin{equation}
{\rm{E}}_{243}\left(a,b_1,b_2; c; x,y,z\right)=\sum\limits_{m,n,p=0}^\infty\frac{(a)_n(b_1)_{m+n-p}(b_2)_{p-n}}{(c)_m}\frac{x^m}{m!}\frac{y^n}{n!}\frac{z^p}{p!},
\end{equation}

region of convergence:
$$ \left\{ r<\infty,\,\,\,s<1,\,\,\,t<\infty
\right\}.
$$

System of partial differential equations:

$
\left\{
\begin{aligned}
&{xu_{xx} +(c-x)u_{x} -yu_{y} +zu_{z} -b_{1} u=0,} \\
& {y(1+y)u_{yy} +xyu_{xy} -(y+1)zu_{yz} +axu_{x}} {+\left[1-b_{2} +\left(a+b_{1} +1\right)y\right]u_{y}-azu_z +ab_{1} u=0,} \\
& {zu_{zz} -xu_{xz} -yu_{yz} -yu_{y} +\left(1-b_{1} +z\right)u_{z} +b_{2} u=0,}
\end{aligned}
\right.
$

where $u\equiv \,\,{\rm{E}}_{243}\left(a,b_1,b_2; c; x,y,z\right)$.

Particular solutions:

$
{u_1} ={\rm{E}}_{243}\left(a,b_1,b_2; c; x,y,z\right) ,
$

$
{u_2} = {x^{1 - c}}{\rm{E}}_{243}\left(a,1-c+b_1,b_2; 2-c; x,y,z\right).
$

\bigskip

\begin{equation}
{\rm{E}}_{244}\left(a,b_1,b_2; c; x,y,z\right)=\sum\limits_{m,n,p=0}^\infty\frac{(a)_p(b_1)_{m+n-p}(b_2)_{p-n}}{(c)_m}\frac{x^m}{m!}\frac{y^n}{n!}\frac{z^p}{p!},
\end{equation}

region of convergence:
$$ \left\{ r<\infty,\,\,\,s<\infty,\,\,\,t<1
\right\}.
$$

System of partial differential equations:

$
\left\{
\begin{aligned}
&{xu_{xx} +(c-x)u_{x} -yu_{y} +zu_{z} -b_{1} u=0,} \\& {yu_{yy} -zu_{yz} +xu_{x} +\left(1-b_{2} +y\right)u_{y} -zu_{z} +b_{1} u=0,} \\& {z(1+z)u_{zz} -xu_{xz} -y(1+z)u_{yz} -ayu_{y}} {+\left[1-b_{1} +\left(a+b_{2} +1\right)z\right]u_{z} +ab_{2} u=0,}
\end{aligned}
\right.
$

where $u\equiv \,\,{\rm{E}}_{244}\left(a,b_1,b_2; c; x,y,z\right)$.

Particular solutions:

$
{u_1} ={\rm{E}}_{244}\left(a,b_1,b_2; c; x,y,z\right) ,
$

$
{u_2} = {x^{1 - c}}{\rm{E}}_{244}\left(a,1-c+b_1,b_2; 2-c; x,y,z\right).
$

\bigskip

\begin{equation}
{\rm{E}}_{245}\left(a,b; c; x,y,z\right)=\sum\limits_{m,n,p=0}^\infty\frac{(a)_{m+n-p}(b)_{p-n}}{(c)_m}\frac{x^m}{m!}\frac{y^n}{n!}\frac{z^p}{p!},
\end{equation}

region of convergence:
$$ \left\{ r<\infty,\,\,\,s<\infty,\,\,\,t<\infty
\right\}.
$$

System of partial differential equations:

$
\left\{
\begin{aligned}
&{xu_{xx} +(c-x)u_{x} -yu_{y} +zu_{z} - a u=0,} \\& {yu_{yy} -zu_{yz} +xu_{x} +\left(1-b +y\right)u_{y} -zu_{z} +a u=0,} \\& {zu_{zz} -xu_{xz} -yu_{yz} -yu_{y} +\left(1- a +z\right)u_{z} +b u=0,}
\end{aligned}
\right.
$

where $u\equiv \,\,{\rm{E}}_{245}\left(a,b; c; x,y,z\right)$.

Particular solutions:

$
{u_1} ={\rm{E}}_{245}\left(a,b; c; x,y,z\right),
$

$
{u_2} = {x^{1 - c}}{\rm{E}}_{245}\left(1-c+a,b; 2-c; x,y,z\right).
$

\bigskip

\begin{equation}
{\rm{E}}_{246}\left(a_1,a_2,b_1,b_2; c; x,y,z\right)=\sum\limits_{m,n,p=0}^\infty\frac{(a_1)_{m+n}(a_2)_p(b_1)_{n-p}(b_2)_{m-n}}{(c)_m}\frac{x^m}{m!}\frac{y^n}{n!}\frac{z^p}{p!},
\end{equation}

region of convergence:
$$ \left\{ r+s<1,\,\,\,t<\infty
\right\}.
$$

System of partial differential equations:

$
\left\{
\begin{aligned}
&x(1-x)u_{xx} +y^{2} u_{yy} +\left[c-\left(a_{1} +b_{2} +1\right)x\right]u_{x}+(a_{1} -b_{2} +1)yu_{y} -a_{1} b_{2} u=0, \\& {y(1+y)u_{yy} -x(1-y)u_{xy} -xzu_{xz} -yzu_{yz} +b_{1} xu_{x} } \\&\,\,\,\,\,\,\,\,\,\,\,\,\,\,{+\left[1-b_{2} +\left(a_{1} +b_{1} +1\right)y\right]u_{y} -a_{1}zu_z+a_{1} b_{1} u=0,} \\& {zu_{zz} -yu_{yz} +(1+z-b_{1} )u_{z} +a_{2} u=0,}
\end{aligned}
\right.
$

where $u\equiv \,\,{\rm{E}}_{246}\left(a_1,a_2,b_1,b_2; c; x,y,z\right)$.

Particular solutions:

$
{u_1} ={\rm{E}}_{246}\left(a_1,a_2,b_1,b_2; c; x,y,z\right),
$

$
{u_2} = {x^{1 - c}}{\rm{E}}_{246}\left(1-c+a_1,a_2,b_1,1-c+b_2; 2-c; x,y,z\right).
$

\bigskip

\begin{equation}
{\rm{E}}_{247}\left(a,b_1,b_2; c; x,y,z\right)=\sum\limits_{m,n,p=0}^\infty\frac{(a)_{m+n}(b_1)_{n-p}(b_2)_{m-n}}{(c)_m}\frac{x^m}{m!}\frac{y^n}{n!}\frac{z^p}{p!},
\end{equation}

region of convergence:
$$ \left\{ \left\{s+2\sqrt{rs}<1\right\}=\left\{\sqrt{s}<\sqrt{1+r}-\sqrt{r}\right\},\,\,\,t<\infty
\right\}.
$$

System of partial differential equations:

$
\left\{
\begin{aligned}
&x(1-x)u_{xx} +y^{2} u_{yy} +\left[c-\left(a +b_{2} +1\right)x\right]u_{x}+\left(a -b_{2} +1\right)yu_{y} -a b_{2} u=0, \\& y(1+y)u_{yy} -x(1-y)u_{xy} -xzu_{xz} -yzu_{yz}\\& \,\,\,\,\,\,\,\,\, +b_{1} xu_{x} +\left[1-b_{2} +\left(a+b_{1} +1\right)y\right]u_{y} -azu_z+ab_{1} u=0, \\& {zu_{zz} -yu_{yz} +\left(1-b_{1} \right)u_{z} +u=0;}
\end{aligned}
\right.
$

where $u\equiv \,\,{\rm{E}}_{247}\left(a,b_1,b_2; c; x,y,z\right)$.

Particular solutions:

$
{u_1} ={\rm{E}}_{247}\left(a,b_1,b_2; c; x,y,z\right),
$

$
{u_2} = {x^{1 - c}}{\rm{E}}_{247}\left(1-c+a,b_1,1-c+b_2; 2-c; x,y,z\right).
$

\bigskip

\begin{equation}
{\rm{E}}_{248}\left(a_1,a_2,b_1,b_2; c; x,y,z\right)=\sum\limits_{m,n,p=0}^\infty\frac{(a_1)_{m+n}(a_2)_m(b_1)_{n-p}(b_2)_{p-n}}{(c)_m}\frac{x^m}{m!}\frac{y^n}{n!}\frac{z^p}{p!},
\end{equation}

region of convergence:
$$ \left\{ r+s<1,\,\,\,t<\infty
\right\}.
$$

System of partial differential equations:

$
\left\{
\begin{aligned}
&{x(1-x)u_{xx} -xyu_{xy} +\left[c-\left(a_{1} +a_{2}+1\right)x\right]u_{x} -a_{2} yu_{y} -a_{1} a_{2} u=0,} \\
& y(1+y)u_{yy} +xyu_{xy} -xzu_{xz} -(1+y)zu_{yz}\\& \,\,\,\,\,\,\,\,\, +b_{1} xu_{x} {+\left[1-b_{2} +\left(a_{1} +b_{1} +1\right)y\right]u_{y} -a_{1}zu_z+a_{1} b_{1} u=0,} \\
& {zu_{zz} -yu_{yz} -yu_{y} +\left(1-b_{1} +z\right)u_{z}+b_{2} u=0,}
\end{aligned}
\right.
$

where $u\equiv \,\,{\rm{E}}_{248}\left(a_1,a_2,b_1,b_2; c; x,y,z\right)$.

Particular solutions:

$
{u_1} ={\rm{E}}_{248}\left(a_1,a_2,b_1,b_2; c; x,y,z\right) ,
$

$
{u_2} = {x^{1 - c}}{\rm{E}}_{248}\left(1-c+a_1,1-c+a_2,b_1,b_2; 2-c; x,y,z\right).
$

\bigskip

\begin{equation}
{\rm{E}}_{249}\left(a_1,a_2,b_1,b_2; c; x,y,z\right)=\sum\limits_{m,n,p=0}^\infty\frac{(a_1)_{m+n}(a_2)_p(b_1)_{n-p}(b_2)_{p-n}}{(c)_m}\frac{x^m}{m!}\frac{y^n}{n!}\frac{z^p}{p!},
\end{equation}

region of convergence:
$$ \left\{ r<\infty,\,\,\,s<1 \vee t<1
\right\}.
$$

System of partial differential equations:

$
\left\{
\begin{aligned}
&{xu_{xx} +(c-x)u_{x} -yu_{y} -a_{1} u=0,} \\& y(1+y)u_{yy} +xyu_{xy} -xzu_{xz} -(1+y)zu_{yz}\\& \,\,\,\,\,\,\,\,\, +b_{1} xu_{x} +\left[1-b_{2} +\left(a_{1} +b_{1} +1\right)y\right]u_{y} -a_{1}zu_z+a_{1} b_1 u=0,\\& {z(1+z)u_{zz} -y(1+z)u_{yz} -a_{2} yu_{y}} {+\left[1-b_{1} +\left(a_{2} +b_{2} +1\right)z\right]u_{z} +a_{2} b_{2} u=0,}
\end{aligned}
\right.
$

where $u\equiv \,\,{\rm{E}}_{249}\left(a_1,a_2,b_1,b_2; c; x,y,z\right)$.

Particular solutions:

$
{u_1} ={\rm{E}}_{249}\left(a_1,a_2,b_1,b_2; c; x,y,z\right),
$

$
{u_2} = {x^{1 - c}}{\rm{E}}_{249}\left(1-c+a_1,a_2,b_1,b_2; 2-c; x,y,z\right).
$

\bigskip

\begin{equation}
{\rm{E}}_{250}\left(a,b_1,b_2; c; x,y,z\right)=\sum\limits_{m,n,p=0}^\infty\frac{(a)_{m+n}(b_1)_{n-p}(b_2)_{p-n}}{(c)_m}\frac{x^m}{m!}\frac{y^n}{n!}\frac{z^p}{p!},
\end{equation}

region of convergence:
$$ \left\{ r<\infty,\,\,\,s<1,\,\,\,t<\infty
\right\}.
$$

System of partial differential equations:

$
\left\{
\begin{aligned}
&{xu_{xx} +(c-x)u_{x} -yu_{y} -au=0,} \\& y(1+y)u_{yy} +xyu_{xy} -xzu_{xz} -(1+y)zu_{yz}\\& \,\,\,\,\,\,\,\,\, +b_{1} xu_{x} {+\left[1-b_{2} +\left(a +b_{1} +1\right)y\right]u_{y} -a zu_{z} +ab_{1}u=0,} \\& {zu_{zz} -yu_{yz} -yu_{y}+\left(1-b_1+z\right)u_{z}-b_2u=0, }
\end{aligned}
\right.
$

where $u\equiv \,\,{\rm{E}}_{250}\left(a,b_1,b_2; c; x,y,z\right)$.

Particular solutions:

$
{u_1} ={\rm{E}}_{250}\left(a,b_1,b_2; c; x,y,z\right),
$

$
{u_2} = {x^{1 - c}}{\rm{E}}_{250}\left(1-c+a,b_1,b_2; 2-c; x,y,z\right).
$

\bigskip

\begin{equation}
{\rm{E}}_{251}\left(a_1,a_2,b_1,b_2; c; x,y,z\right)=\sum\limits_{m,n,p=0}^\infty\frac{(a_1)_{n+p}(a_2)_m(b_1)_{m-p}(b_2)_{p-n}}{(c)_m}\frac{x^m}{m!}\frac{y^n}{n!}\frac{z^p}{p!},
\end{equation}

region of convergence:
$$ \left\{ t(1+r)<1,\,\,\,s<\infty
\right\}.
$$

System of partial differential equations:

$
\left\{
\begin{aligned}
&{x(1-x)u_{xx} +xzu_{xz} +\left[c-\left(a_{2} +b_{1}+1\right)x\right]u_{x} } {+a_{2} zu_{z} -a_{2} b_{1} u=0,} \\& {yu_{yy} -zu_{yz} +\left(1-b_{2} +y\right)u_{y} +zu_{z} +a_{1} u=0,} \\& {z(1+z)u_{zz} -y^{2} u_{yy} -xu_{xz} -\left(a_{1} -b_{2} +1\right)yu_{y} } {+\left[1-b_{1} +\left(a_{1} +b_{2} +1\right)z\right]u_{z} +a_{1} b_{2} u=0;}
\end{aligned}
\right.
$

where $u\equiv \,\,{\rm{E}}_{251}\left(a_1,a_2,b_1,b_2; c; x,y,z\right)$.

Particular solutions:

$
{u_1} ={\rm{E}}_{251}\left(a_1,a_2,b_1,b_2; c; x,y,z\right),
$

$
{u_2} = {x^{1 - c}}{\rm{E}}_{251}\left(a_1,1-c+a_2,1-c+b_1,b_2; 2-c; x,y,z\right).
$

\bigskip

\begin{equation}
{\rm{E}}_{252}\left(a_1,a_2,b_1,b_2; c; x,y,z\right)=\sum\limits_{m,n,p=0}^\infty\frac{(a_1)_{n+p}(a_2)_n(b_1)_{m-p}(b_2)_{p-n}}{(c)_m}\frac{x^m}{m!}\frac{y^n}{n!}\frac{z^p}{p!},
\end{equation}

region of convergence:
$$ \left\{ r<\infty,\,\,\,\left\{s+2\sqrt{st}<1\right\}=\left\{\sqrt{s}<\sqrt{1+t}-\sqrt{t}\right\}
\right\}.
$$

System of partial differential equations:

$
\left\{
\begin{aligned}
&{xu_{xx} +(c-x)u_{x} +zu_{z} -b_{1} u=0,} \\& y(1+y)u_{yy} -(1-y)zu_{yz} +\left[1-b_{2} +(a_{1} +a_{2}+1 )y\right]u_{y}+a_{2} zu_{z} +a_{1} a_{2} u=0, \\ & {z(1+z)u_{zz} -y^{2} u_{yy} -xu_{xz} -\left(a_{1} -b_{2} +1\right)yu_{y}} {+\left[1-b_{1} +\left(a_{2} +b_{2} +1\right)z\right]u_{z} +a_{1} b_{2} u=0;}
\end{aligned}
\right.
$

where $u\equiv \,\,{\rm{E}}_{252}\left(a_1,a_2,b_1,b_2; c; x,y,z\right)$.

Particular solutions:

$
{u_1} ={\rm{E}}_{252}\left(a_1,a_2,b_1,b_2; c; x,y,z\right) ,
$

$
{u_2} = {x^{1 - c}}{\rm{E}}_{252}\left(a_1,a_2,1-c+b_1,b_2; 2-c; x,y,z\right).
$

\bigskip

\begin{equation}
{\rm{E}}_{253}\left(a,b_1,b_2; c; x,y,z\right)=\sum\limits_{m,n,p=0}^\infty\frac{(a)_{n+p}(b_1)_{m-p}(b_2)_{p-n}}{(c)_m}\frac{x^m}{m!}\frac{y^n}{n!}\frac{z^p}{p!},
\end{equation}

region of convergence:
$$ \left\{ r<\infty,\,\,\,s<\infty,\,\,\,t<1
\right\}.
$$

System of partial differential equations:

$
\left\{
\begin{aligned}
&{xu_{xx} +(c-x)u_{x} +zu_{z} -b_{1} u=0,} \\& {yu_{yy} -zu_{yz} +\left(1-b_{2} +y\right)u_{y} +zu_{z} +au=0,} \\& {z(1+z)u_{zz} -y^{2} u_{yy} -xu_{xz} -(a-b_{2} +1)yu_{y} } {+\left[1-b_{1} +\left(a+b_{2} +1\right)z\right]u_{z} +ab_{2} u=0,}
\end{aligned}
\right.
$

where $u\equiv \,\,{\rm{E}}_{253}\left(a,b_1,b_2; c; x,y,z\right)$.

Particular solutions:

$
{u_1} ={\rm{E}}_{253}\left(a,b_1,b_2; c; x,y,z\right),
$

$
{u_2} = {x^{1 - c}}{\rm{E}}_{253}\left(a,1-c+b_1,b_2; 2-c; x,y,z\right).
$

\bigskip

\begin{equation}
{\rm{E}}_{254}\left(a_1,a_2,b_1,b_2; c; x,y,z\right)=\sum\limits_{m,n,p=0}^\infty\frac{(a_1)_{n+p}(a_2)_p(b_1)_{m-p}(b_2)_{m-n}}{(c)_m}\frac{x^m}{m!}\frac{y^n}{n!}\frac{z^p}{p!},
\end{equation}

region of convergence:
$$ \left\{ t(1+r)<1,\,\,\,s<\infty
\right\}.
$$

System of partial differential equations:

$
\left\{
\begin{aligned}
&x(1-x)u_{xx} +xyu_{xy} +xzu_{xz} -yzu_{yz} +\left[c-\left(b_{1} +b_{2}+1 \right)x\right]u_{x}  {+b_{1} yu_{y} +b_{2} zu_{z} -b_{1} b_{2} u=0,} \\& {yu_{yy} -xu_{xy} +\left(1-b_{2} +y\right)u_{y} +zu_{z} +a_{1} u=0,} \\& {z(1+z)u_{zz} -xu_{xz} +yzu_{yz} +a_{2} yu_{y}} {+\left[1-b_{1} +\left(a_{1} +a_{2} +1\right)z\right]u_{z} +a_{1} a_{2} u=0,}
\end{aligned}
\right.
$

where $u\equiv \,\,{\rm{E}}_{254}\left(a_1,a_2,b_1,b_2; c; x,y,z\right)$.

Particular solutions:

$
{u_1} ={\rm{E}}_{254}\left(a_1,a_2,b_1,b_2; c; x,y,z\right),
$

$
{u_2} = {x^{1 - c}}{\rm{E}}_{254}\left(a_1,a_2,1-c+b_1,1-c+b_2; 2-c; x,y,z\right).
$

\bigskip

\begin{equation}
{\rm{E}}_{255}\left(a,b_1,b_2; c; x,y,z\right)=\sum\limits_{m,n,p=0}^\infty\frac{(a)_{n+p}(b_1)_{m-p}(b_2)_{m-n}}{(c)_m}\frac{x^m}{m!}\frac{y^n}{n!}\frac{z^p}{p!},
\end{equation}

region of convergence:
$$ \left\{ r<1,\,\,\,s<\infty,\,\,\,t<\infty
\right\}.
$$

System of partial differential equations:

$
\left\{
\begin{aligned}
&{x(1-x)u_{xx} +xyu_{xy} +xzu_{xz} -yzu_{yz} } {+\left[c-\left(b_{1} +b_{2}+1\right)x\right]u_{x} +b_{1} yu_{y} +b_{2} zu_{z} -b_{1} b_{2} u=0,} \\& {yu_{yy} -xu_{xy} +\left(1-b_{2} +y\right)u_{y} +zu_{z} +au=0,} \\& {zu_{zz} -xu_{xz} +yu_{y} +\left(1-b_{1} +z\right)u_{z} +au=0,}
\end{aligned}
\right.
$

where $u\equiv \,\,{\rm{E}}_{255}\left(a,b_1,b_2; c; x,y,z\right)$.

Particular solutions:

$
{u_1} ={\rm{E}}_{255}\left(a,b_1,b_2; c; x,y,z\right) ,
$

$
{u_2} = {x^{1 - c}}{\rm{E}}_{255}\left(a,1-c+b_1,1-c+b_2; 2-c; x,y,z\right).
$

\bigskip

\begin{equation}
{\rm{E}}_{256}\left(a_1,a_2,b_1,b_2; c; x,y,z\right)=\sum\limits_{m,n,p=0}^\infty\frac{(a_1)_{m+n}(a_2)_n(b_1)_{m-p}(b_2)_{p-n}}{(c)_m}\frac{x^m}{m!}\frac{y^n}{n!}\frac{z^p}{p!},
\end{equation}

region of convergence:
$$ \left\{ r+s<1,\,\,\,t<\infty
\right\}.
$$

System of partial differential equations:

$
\left\{
\begin{aligned}
&{x(1-x)u_{xx} -xyu_{xy} +xzu_{xz} +yzu_{yz}} {+\left[c-\left(a_{1} +b_{1}+1 \right)x\right]u_{x} -b_{1} yu_{y}+a_1zu_z -a_{1} b_{1} u=0,} \\& {y(1+y)u_{yy} +xyu_{xy} -zu_{yz} +a_{2} xu_{x}}{+\left[1-b_{2} +\left(a_{1} +a_{2}+1\right)y\right]u_{y} +a_{1} a_{2} u=0,} \\& {zu_{zz} -xu_{xz} -yu_{y} +\left(1-b_{1} +z\right)u_{z} +b_{2} u=0,}
\end{aligned}
\right.
$

where $u\equiv \,\,{\rm{E}}_{256}\left(a_1,a_2,b_1,b_2; c; x,y,z\right)$.

Particular solutions:

$
{u_1} ={\rm{E}}_{256}\left(a_1,a_2,b_1,b_2; c; x,y,z\right),
$

$
{u_2} = {x^{1 - c}}{\rm{E}}_{256}\left(1-c+a_1,a_2,1-c+b_1,b_2; 2-c; x,y,z\right).
$

\bigskip

\begin{equation}
{\rm{E}}_{257}\left(a_1,a_2,b_1,b_2; c; x,y,z\right)=\sum\limits_{m,n,p=0}^\infty\frac{(a_1)_{m+n}(a_2)_p(b_1)_{m-p}(b_2)_{p-n}}{(c)_m}\frac{x^m}{m!}\frac{y^n}{n!}\frac{z^p}{p!},
\end{equation}

region of convergence:
$$ \left\{ t(1+r)<1,\,\,\,s<\infty
\right\}.
$$

System of partial differential equations:

$
\left\{
\begin{aligned}
&{x(1-x)u_{xx} -xyu_{xy} +xzu_{xz} +yzu_{yz}} {+\left[c-\left(a_{1} +b_1+1 \right)x\right]u_{x} -b_{1} yu_{y}+a_1zu_z -a_{1} b_1 u=0,} \\& {yu_{yy} -zu_{yz} +xu_{x} +\left(1-b_{2} +y\right)u_{y} +a_{1} u=0,} \\& {z(1+z)u_{zz} -xu_{xz} -yzu_{yz} -a_{2} yu_{y} } {+\left[1-b_{1} +\left(a_{2} +b_{2} +1\right)z\right]u_{z} +a_{2} b_{2} u=0,}
\end{aligned}
\right.
$

where $u\equiv \,\,{\rm{E}}_{257}\left(a_1,a_2,b_1,b_2; c; x,y,z\right)$.

Particular solutions:

$
{u_1} ={\rm{E}}_{257}\left(a_1,a_2,b_1,b_2; c; x,y,z\right) ,
$

$
{u_2} = {x^{1 - c}}{\rm{E}}_{257}\left(1-c+a_1,a_2,1-c+b_1,b_2; 2-c; x,y,z\right).
$

\bigskip

\begin{equation}
{\rm{E}}_{258}\left(a,b_1,b_2; c; x,y,z\right)=\sum\limits_{m,n,p=0}^\infty\frac{(a)_{m+n}(b_1)_{m-p}(b_2)_{p-n}}{(c)_m}\frac{x^m}{m!}\frac{y^n}{n!}\frac{z^p}{p!},
\end{equation}

region of convergence:
$$ \left\{ r<1,\,\,\,s<\infty,\,\,\,t<\infty
\right\}.
$$

System of partial differential equations:

$
\left\{
\begin{aligned}
&x(1-x)u_{xx} -xyu_{xy} +xzu_{xz} +yzu_{yz}  +\left[c-\left(a+b_1+1\right)x\right]u_{x} -b_{1} yu_{y}+azu_z -ab_1 u=0, \\& {yu_{yy} -zu_{yz} +xu_{x} +\left(1-b_{2} +y\right)u_{y} +au=0,} \\&  {zu_{zz} -xu_{xz} -yu_{y} +\left(1-b_{1} +z\right)u_{z} +b_{2} u=0,}
\end{aligned}
\right.
$

where $u\equiv \,\,{\rm{E}}_{258}\left(a,b_1,b_2; c; x,y,z\right)$.

Particular solutions:

$
{u_1} ={\rm{E}}_{258}\left(a,b_1,b_2; c; x,y,z\right) ,
$

$
{u_2} = {x^{1 - c}}{\rm{E}}_{258}\left(1-c+a,1-c+b_1,b_2; 2-c; x,y,z\right).
$

\bigskip

\begin{equation}
{\rm{E}}_{259}\left(a_1,a_2,a_3,b; c;x,y,z\right)=\sum\limits_{m,n,p=0}^\infty\frac{(a_1)_{n+p}(a_2)_m(a_3)_p(b)_{m-n-p}}{(c)_m}\frac{x^m}{m!}\frac{y^n}{n!}\frac{z^p}{p!},
\end{equation}

region of convergence:
$$ \left\{ r<\infty,\,\,\,s<\infty,\,\,\,t(1+t)<1
\right\}.
$$

System of partial differential equations:

$
\left\{
\begin{aligned}
&{x(1-x)u_{xx} +xyu_{xy} +xzu_{xz}}{+\left[c-\left(a_{2} +b+1\right)x\right]u_{x} +a_{2} yu_{y} +a_{2} zu_{z} -a_{2} bu=0,} \\& {yu_{yy} -xu_{xy} +zu_{yz} +(1-b+y)u_{y} +zu_{z} +a_{1} u=0,} \\& {z(1+z)u_{zz} -xu_{xz} +y(1+z)u_{yz} +a_{3} yu_{y}}{+\left[1-b+\left(a_{1} +a_{3} +1\right)z\right]u_{z} +a_{1} a_{3} u=0,}
\end{aligned}
\right.
$

where $u\equiv \,\,{\rm{E}}_{259}\left(a_1,a_2,a_3,b; c;x,y,z\right)$.

Particular solutions:

$
{u_1} ={\rm{E}}_{259}\left(a_1,a_2,a_3,b; c;x,y,z\right) ,
$

$
{u_2} = {x^{1 - c}}{\rm{E}}_{259}\left(a_1,1-c+a_2,a_3,1-c+b; 2-c;x,y,z\right).
$

\bigskip

\begin{equation}
{\rm{E}}_{260}\left(a_1,a_2,a_3,b; c;x,y,z\right)=\sum\limits_{m,n,p=0}^\infty\frac{(a_1)_{n+p}(a_2)_n(a_3)_p(b)_{m-n-p}}{(c)_m}\frac{x^m}{m!}\frac{y^n}{n!}\frac{z^p}{p!},
\end{equation}

region of convergence:
$$ \left\{ r<\infty,\,\,\,s<1, \,\,\,t<1
\right\}.
$$

System of partial differential equations:

$
\left\{
\begin{aligned}
&{xu_{xx} +(c-x)u_{x} +yu_{y} +zu_{z} -bu=0,} \\& {y(1+y)u_{yy} -xu_{xy} +(1+y)zu_{yz}} {+\left[1-b+\left(a_{1} +a_{2}+1\right)y\right]u_{y} +a_2zu_z+a_{1} a_{2} u=0,} \\& {z(1+z)u_{zz} -xu_{xz} +y(1+z)u_{yz} +a_{3} yu_{y}}{+\left[1-b+\left(a_{1} +a_{3} +1\right)z\right]u_{z} +a_{1} a_{3} u=0,}
\end{aligned}
\right.
$

where $u\equiv \,\,{\rm{E}}_{260}\left(a_1,a_2,a_3,b; c;x,y,z\right)$.

Particular solutions:

$
{u_1} ={\rm{E}}_{260}\left(a_1,a_2,a_3,b; c;x,y,z\right) ,
$

$
{u_2} = {x^{1 - c}}{\rm{E}}_{260}\left(a_1,a_2,a_3,1-c+b; 2-c;x,y,z\right).
$

\bigskip

\begin{equation}
{\rm{E}}_{261}\left(a_1,a_2,b; c;x,y,z\right)=\sum\limits_{m,n,p=0}^\infty\frac{(a_1)_{n+p}(a_2)_p(b)_{m-n-p}}{(c)_m}\frac{x^m}{m!}\frac{y^n}{n!}\frac{z^p}{p!},
\end{equation}

region of convergence:
$$ \left\{ r<\infty,\,\,\,s<\infty,\,\,\,t<1
\right\}.
$$

System of partial differential equations:

$
\left\{
\begin{aligned}
&{xu_{xx} +(c-x)u_{x} +yu_{y} +zu_{z} -bu=0,} \\& {yu_{yy} -xu_{xy} +zu_{yz} +(1-b+y)u_{y} + zu_{z} +a_{1} u=0,} \\& {z(1+z)u_{zz} -xu_{xz} +y(1+z)u_{yz} +a_{2} yu_{y}}  {+\left[(1-b+\left(a_{1} +a_{2} +1\right)z\right])u_{z} +a_{1} a_{2} u=0,}
\end{aligned}
\right.
$

where $u\equiv \,\,{\rm{E}}_{261}\left(a_1,a_2,b; c;x,y,z\right)$.

Particular solutions:

$
{u_1} ={\rm{E}}_{261}\left(a_1,a_2,b; c;x,y,z\right) ,
$

$
{u_2} = {x^{1 - c}}{\rm{E}}_{261}\left(a_1,a_2,1-c+b; 2-c;x,y,z\right).
$

\bigskip

\begin{equation}
{\bf{E}}_{262}\left(a_1,a_2,a_3,b;c; x,y,z\right)=\sum\limits_{m,n,p=0}^\infty\frac{(a_1)_{m+n}(a_2)_n(a_3)_p(b)_{m-n-p}}{(c)_m}\frac{x^m}{m!}\frac{y^n}{n!}\frac{z^p}{p!},
\end{equation}

$$ \left\{ \left\{s+2\sqrt{rs}<1\right\}=\left\{\sqrt{s}<\sqrt{1+r}-\sqrt{r}\right\},\,\,\,t<\infty
\right\}.
$$

System of partial differential equations:

$
\left\{
\begin{aligned}
&{x(1-x)u_{xx} +y^{2} u_{yy} +xzu_{xz} +yzu_{yz} } \\&\,\,\,\,\,\,\,\,\, +\left[c-\left(1+a_{1} +b\right)x\right]u_{x}+\left(1+a_{1} -b\right)yu_{y} +a_{1} zu_{z} -a_{1} bu=0, \\& {y(1+y)u_{yy} -x(1-y)u_{xy} +zu_{yz} +a_{2} xu_{x} } {+\left[1-b+\left(a_{1} +a_{2}+1 \right)y\right]u_{y} +a_{1} a_{2} u=0,} \\& {zu_{zz} -xu_{xz} +yu_{yz} +(1-b+z)u_{z} +a_{3} u=0,}
\end{aligned}
\right.
$

where $u\equiv \,\,{\bf{E}}_{262}\left(a_1,a_2,a_3,b;c; x,y,z\right)$.

Particular solutions:

$
{u_1} ={\rm{E}}_{262}\left(a_1,a_2,a_3,b; c;x,y,z\right) ,
$

$
{u_2} = {x^{1 - c}}{\rm{E}}_{262}\left(1-c+a_1,a_2,a_3,1-c+b; 2-c;x,y,z\right).
$

\bigskip

\begin{equation}
{\rm{E}}_{263}\left(a_1,a_2,a_3,b;c; x,y,z\right)=\sum\limits_{m,n,p=0}^\infty\frac{(a_1)_{m+n}(a_2)_p(a_3)_p(b)_{m-n-p}}{(c)_m}\frac{x^m}{m!}\frac{y^n}{n!}\frac{z^p}{p!},
\end{equation}

$$ \left\{ t(1+r)<1,\,\,\,s<\infty
\right\}.
$$

System of partial differential equations:

$
\left\{
\begin{aligned}
&{x(1-x)u_{xx} +y^{2} u_{yy} +xzu_{xz} +yzu_{yz} }\\& \,\,\,\,\,\,\,\,\,{+\left[c-(a_{1} +b+1)x\right]u_{x} +\left(1-b+a_{1}\right)yu_{y} +a_{1} zu_{z} -a_{1} bu=0,} \\
& {yu_{yy} -xu_{xy} +zu_{yz} +xu_{x} +(1-b+y)u_{y} +a_{1} u=0,} \\
& {z(1+z)u_{zz} -xu_{xz} +yu_{yz} +\left[1-b+\left(1+a_{2} +a_{3}\right)z\right]u_{z} +a_{2} a_{3} u=0,}
\end{aligned}
\right.
$

where $u\equiv \,\,{\rm{E}}_{263}\left(a_1,a_2,a_3,b;c; x,y,z\right)$.

Particular solutions:

$
{u_1} ={\rm{E}}_{263}\left(a_1,a_2,a_3,b; c;x,y,z\right) ,
$

$
{u_2} = {x^{1 - c}}{\rm{E}}_{263}\left(1-c+a_1,a_2,a_3,1-c+b; 2-c;x,y,z\right).
$

\bigskip

\begin{equation}
{\rm{E}}_{264}\left(a_1,a_2,b;c; x,y,z\right)=\sum\limits_{m,n,p=0}^\infty\frac{(a_1)_{m+n}(a_2)_n(b)_{m-n-p}}{(c)_m}\frac{x^m}{m!}\frac{y^n}{n!}\frac{z^p}{p!},
\end{equation}

region of convergence:
$$ \left\{ \left\{s+2\sqrt{rs}<1\right\}=\left\{\sqrt{s}<\sqrt{1+r}-\sqrt{r}\right\},\,\,\,t<\infty
\right\}.
$$

System of partial differential equations:

$
\left\{
\begin{aligned}
&{x(1-x)u_{xx} +y^{2} u_{yy} +xzu_{xz} +yzu_{yz} }\\& \,\,\,\,\,\,\,\,\, {+\left[c-\left(a_{1} +b+1\right)x\right]u_{x} +\left(1+a_{1} -b\right)yu_{y} +a_{1} zu_{z} -a_{1} bu=0,} \\& {y(1+y)u_{yy} -x(1-y)u_{xy} +zu_{yz} +a_{2} xu_{x} } {+\left[1-b+\left(a_{1} +a_{2}+1\right)y\right]u_{y} +a_{1} a_{2} u=0,} \\& {zu_{zz} -xu_{xz} +yu_{yz} +(1-b)u_{z} +u=0,}
\end{aligned}
\right.
$

where $u\equiv \,\,{\rm{E}}_{264}\left(a_1,a_2,b;c; x,y,z\right)$.

Particular solutions:

$
{u_1} ={\rm{E}}_{264}\left(a_1,a_2,b; c;x,y,z\right) ,
$

$
{u_2} = {x^{1 - c}}{\rm{E}}_{264}\left(1-c+a_1,a_2,1-c+b; 2-c;x,y,z\right).
$

\bigskip

\begin{equation}
{\rm{E}}_{265}\left(a_1,a_2,b;c; x,y,z\right)=\sum\limits_{m,n,p=0}^\infty\frac{(a_1)_{m+n}(a_2)_p(b)_{m-n-p}}{(c)_m }\frac{x^m}{m!}\frac{y^n}{n!}\frac{z^p}{p!},
\end{equation}

region of convergence:
$$ \left\{ r<1,\,\,\,s<\infty,\,\,\,t<\infty
\right\}.
$$

System of partial differential equations:

$
\left\{
\begin{aligned}
&{x(1-x)u_{xx} +y^{2} u_{yy} +xzu_{xz} +yzu_{yz} } \\& \,\,\,\,\,\,\,\,\,{+\left[c-\left(a_{1} +b+1\right)x\right]u_{x} +\left(a_{1} -b+1\right)yu_{y} +a_{1} zu_{z} -a_{1} bu=0,} \\& {yu_{yy} -xu_{xy} +zu_{yz} +xu_{x} +(1-b+y)u_{y} +a_{1} u=0,} \\& {zu_{zz} -xu_{xz} +yu_{yz} +(1-b+z)u_{z} +a_{2} u=0,}
\end{aligned}
\right.
$

where $u\equiv \,\,{\rm{E}}_{265}\left(a_1,a_2,b;c; x,y,z\right)$.

Particular solutions:

$
{u_1} ={\rm{E}}_{265}\left(a_1,a_2,b; c;x,y,z\right) ,
$

$
{u_2} = {x^{1 - c}}{\rm{E}}_{265}\left(1-c+a_1,a_2,1-c+b; 2-c;x,y,z\right).
$

\bigskip

\begin{equation}
{\rm{E}}_{266}\left(a,b;c; x,y,z\right)=\sum\limits_{m,n,p=0}^\infty\frac{(a)_{m+n}(b)_{m-n-p}}{(c)_m}\frac{x^m}{m!}\frac{y^n}{n!}\frac{z^p}{p!},
\end{equation}

region of convergence:
$$ \left\{ r<1,\,\,\,s<\infty,\,\,\,t<\infty
\right\}.
$$

System of partial differential equations:

$
\left\{
\begin{aligned}
&{x(1-x)u_{xx} +y^{2} u_{yy} +xzu_{xz} +yzu_{yz} } \\& \,\,\,\,\,\,\,\,\,{+\left[c-\left(a_{1} +b+1\right)x\right]u_{x} +\left(a_{1} -b+1\right)yu_{y} +a_{1} zu_{z} -a_{1} bu=0,} \\& {yu_{yy} -xu_{xy} +zu_{yz} +xu_{x} +(1-b+y)u_{y} +au=0,} \\& {zu_{zz} -xu_{xz} +yu_{yz} +(1-b)u_{z} +u=0,}
\end{aligned}
\right.
$

where $u\equiv \,\,{\rm{E}}_{266}\left(a,b;c; x,y,z\right)$.

Particular solutions:

$
{u_1} ={\rm{E}}_{266}\left(a,b;c; x,y,z\right) ,
$

$
{u_2} = {x^{1 - c}}{\rm{E}}_{266}\left(1-c+a,1-c+b;2-c; x,y,z\right).
$

\bigskip

\begin{equation}
{\rm{E}}_{267}\left(a, b_1, b_2; c; x,y,z\right)=\sum\limits_{m,n,p=0}^\infty\frac{(a)_{m+n}(b_1)_{m+n-p}(b_2)_{p-n}}{(c)_{m}}\frac{x^m}{m!}\frac{y^n}{n!}\frac{z^p}{p!},
\end{equation}

region of convergence:
$$ \left\{ \sqrt{r}+\sqrt{s}<1,\,\,\,t<\infty
\right\}.
$$

System of partial differential equations:

$
\left\{
\begin{aligned}
&{x(1-x)u_{xx} -y^{2} u_{yy} -2xyu_{xy} +xzu_{xz} +yzu_{yz} } \\& \,\,\,\,\,\,\,\,\,{+\left[c-\left(a+b_{1}+1 \right)x\right]u_{x} -\left(a+b_{1}+1\right)yu_{y} +azu_{z} -ab_{1} u=0,} \\& {y(1+y)u_{yy} +x^{2} u_{xx} +2xyu_{xy} -xzu_{xz} -(1+y)zu_{yz}} \\&\,\,\,\,\,\,\,\,\, {+\left(1+a+b_{1}\right)xu_{x} +\left[1-b+\left(1+a+b_{1}\right)y\right]u_{y} -azu_{z} +ab_{1} u=0,} \\& {zu_{zz} -xu_{xz} -yu_{yz} -yu_{y} +\left(1-b_{1} +z\right)u_{z} +b_{2} u=0,}
\end{aligned}
\right.
$

where $u\equiv \,\,{\rm{E}}_{267}\left(a, b_1, b_2; c; x,y,z\right)$.

Particular solutions:

$
{u_1} ={\rm{E}}_{267}\left(a, b_1, b_2; c; x,y,z\right) ,
$

$
{u_2} = {x^{1 - c}}{\rm{E}}_{267}\left(1-c+a, 1-c+b_1, b_2; 2-c; x,y,z\right).
$

\bigskip

\begin{equation}
{\rm{E}}_{268}\left(a,b_1, b_2; c; x,y,z\right)=\sum\limits_{m,n,p=0}^\infty\frac{(a)_{n}(b_1)_{m+n-p}(b_2)_{m+p-n}}{(c)_{m}}\frac{x^m}{m!}\frac{y^n}{n!}\frac{z^p}{p!},
\end{equation}

region of convergence:
$$ \left\{ \left\{s+2\sqrt{rs}<1\right\}=\left\{\sqrt{s}<\sqrt{1+r}-\sqrt{r}\right\},\,\,\,t<\infty
\right\}.
$$

System of partial differential equations:

$
\left\{
\begin{aligned}
&{x(1-x)u_{xx} +y^{2} u_{yy} +z^{2} u_{zz} -2yzu_{yz} } \\
& \,\,\,\,\,\,\,\,\,+\left[c-\left(b_{1} +b_{2}+1\right)x\right]u_{x}+\left(b_{1} -b_{2}+1 \right)yu_{y} +\left(b_{2} -b_{1}+1 \right)zu_{z} -b_{1} b_{2} u=0, \\
& {y(1+y)u_{yy} -x(1-y)u_{xy} -(1+y)zu_{yz} +axu_{x} } \\
&\,\,\,\,\,\,\,\,\, {+\left[1-b_{2} +\left(a+b_{1}+1\right)y\right]u_{y} -azu_{z} +ab_{1} u=0,} \\
& {zu_{zz} -xu_{xz} -yu_{yz} +xu_{x} -yu_{y} +\left(1-b_{1} +z\right)u_{z} +b_{2} u=0,}
\end{aligned}
\right.
$

where $u\equiv \,\,{\rm{E}}_{268}\left(a,b_1, b_2; c; x,y,z\right)$.

Particular solutions:

$
{u_1} ={\rm{E}}_{268}\left(a, b_1, b_2; c; x,y,z\right) ,
$

$
{u_2} = {x^{1 - c}}{\rm{E}}_{268}\left(a, 1-c+b_1, 1-c+b_2; 2-c; x,y,z\right).
$

\bigskip

\begin{equation}
{\rm{E}}_{269}\left(a, b; c; x,y,z\right)=\sum\limits_{m,n,p=0}^\infty\frac{(a)_{m+n-p}(b)_{m+p-n}}{(c)_{m}}\frac{x^m}{m!}\frac{y^n}{n!}\frac{z^p}{p!},
\end{equation}

region of convergence:
$$ \left\{ r<1,\,\,\,s<\infty,\,\,\,t<\infty
\right\}.
$$

System of partial differential equations:

$
\left\{
\begin{aligned}
&{x(1-x)u_{xx} +y^{2} u_{yy} +z^{2} u_{zz} -2yzu_{yz} } \\&\,\,\,\,\,\,\,\,\, +\left[c-\left(a +b+1\right)x\right]u_{x}+\left(a -b+1 \right)yu_{y} +\left(1-a -b \right)zu_{z} -a b u=0, \\& {yu_{yy}-xu_{xy}-zu_{yz}+xu_x+\left(1-b+y\right)u_y-zu_z+ au=0,} \\& {zu_{zz} -xu_{xz} -yu_{yz} +xu_{x} -yu_{y} +\left(1- a +z\right)u_{z} +b u=0,}
\end{aligned}
\right.
$

where $u\equiv \,\,{\rm{E}}_{269}\left(a, b; c; x,y,z\right)$.

Particular solutions:

$
{u_1} ={\rm{E}}_{269}\left(a, b; c; x,y,z\right) ,
$

$
{u_2} = {x^{1 - c}}{\rm{E}}_{269}\left(1-c+a, 1-c+b; 2-c; x,y,z\right).
$

\bigskip

\begin{equation}
{\rm{E}}_{270}\left(a,b_1, b_2; c; x,y,z\right)=\sum\limits_{m,n,p=0}^\infty\frac{(a)_{m+p}(b_1)_{m+n-p}(b_2)_{p-n}}{(c)_{m}}\frac{x^m}{m!}\frac{y^n}{n!}\frac{z^p}{p!},\,\,
\end{equation}

region of convergence:
$$ \left\{ \left\{t+2\sqrt{rt}<1\right\}=\left\{\sqrt{t}<\sqrt{1+r}-\sqrt{r}\right\},\,\,\,s<\infty
\right\}.
$$

System of partial differential equations:

$
\left\{
\begin{aligned}
& x(1-x)u_{xx} +z^{2} u_{zz} -xyu_{xy} -yzu_{yz}\\& \,\,\,\,\,\,\,\,\, +\left[c-\left(a+b_{1} +1\right)x\right]u_{x}-{ayu_{y} +\left(a-b_{1} +1\right)zu_{z} -ab_{1} u=0,} \\
& {yu_{yy} -zu_{yz} +xu_{x} +\left(1-b_{2} +y\right)u_{y} -zu_{z} +b_{1} u=0,} \\
&z(1+z)u_{zz}-xyu_{xy} -x(1-z)u_{xz} -y(1+z)u_{yz}  \\
& {\,\,\,\,\,\,\,\,\,-b_{2} xu_{x} -ayu_{y}+\left[1-b_{1} +\left(a+b_{2} +1\right)z\right]u_{z} +ab_{2} u=0,}
\end{aligned}
\right.
$

where $u\equiv \,\,{\rm{E}}_{270}\left(a,b_1, b_2; c; x,y,z\right)$.

Particular solutions:

$
{u_1} ={\rm{E}}_{270}\left(a,b_1, b_2; c; x,y,z\right) ,
$

$
{u_2} = {x^{1 - c}}{\rm{E}}_{270}\left(1-c+a,1-c+b_1, b_2; 2-c; x,y,z\right).
$

\bigskip

\begin{equation}
{\bf{E}}_{271}\left(a, b_1, b_2; c; x,y,z\right)=\sum\limits_{m,n,p=0}^\infty\frac{(a)_{n+p}(b_1)_{m+n-p}(b_2)_{m-n}}{(c)_{m}}\frac{x^m}{m!}\frac{y^n}{n!}\frac{z^p}{p!},
\end{equation}

region of convergence:
$$ \left\{ \left\{s+2\sqrt{rs}<1\right\}=\left\{\sqrt{s}<\sqrt{1+r}-\sqrt{r}\right\},\,\,\,t<\infty
\right\}.
$$

System of partial differential equations:

$
\left\{
\begin{aligned}
& x(1-x)u_{xx} +y^{2} u_{yy}+xzu_{xz} -yzu_{yz}  \\& \,\,\,\,\,\,\,\,\, +\left[c-\left(b_{1}+b_2 +1\right)x\right]u_{x}+ {\left(b_1-b_2+1\right)yu_{y} +b_2zu_{z} -b_{1}b_2 u=0,} \\& y(1+y)u_{yy}- z^2u_{zz}-x(1-y)u_{xy} +xzu_{xz}\\&\,\,\,\,\,\,\,\,\,  +axu_x +\left[1-b_{2}+\left(a+b_1+1\right)y\right]u_{y}-\left(a-b_1+1\right)zu_{z}+ab_1 u=0, \\ &{zu_{zz} -xu_{xz} -yu_{yz} +yu_y+\left(1-b_1+z\right) u_{z} +au=0, }
\end{aligned}
\right.
$

where $u\equiv \,\,{\bf{E}}_{271}\left(a, b_1, b_2; c; x,y,z\right)$.

Particular solutions:

$
{u_1} ={\rm{E}}_{271}\left(a,b_1, b_2; c; x,y,z\right) ,
$

$
{u_2} = {x^{1 - c}}{\rm{E}}_{271}\left(a,1-c+b_1, 1-c+b_2; 2-c; x,y,z\right).
$

\bigskip

\begin{equation}
{\rm{E}}_{272}\left(a,b_1, b_2; c; x,y,z\right)=\sum\limits_{m,n,p=0}^\infty\frac{(a)_{m+p}(b_1)_{n-p}(b_2)_{p-n}}{(c)_{m}}\frac{x^m}{m!}\frac{y^n}{n!}\frac{z^p}{p!},
\end{equation}

region of convergence:
$$ \left\{ r<\infty,\,\,\,s<\infty,\,\,\,t<1
\right\}.
$$

System of partial differential equations:

$
\left\{
\begin{aligned}
& {xu_{xx} +(c-x)u_{x} -zu_{z} -au=0,} \\
& {yu_{yy} -zu_{yz} +\left(1-b_{2} +y\right)u_{y} -zu_{z} +b_{1} u=0,} \\
&z(1+z)u_{zz} -xyu_{xy} +xzu_{xz} -y(1+z)u_{yz}\\& \,\,\,\,\,\,\,\,\,+ b_{2} xu_{x} -ayu_{y} {+\left[1-b_{1} +\left(a+b_{2} +1\right)z\right]u_{z} +ab_{2} u=0,}
\end{aligned}
\right.
$

where $u\equiv \,\,{\rm{E}}_{272}\left(a,b_1, b_2; c; x,y,z\right)$.

Particular solutions:

$
{u_1} ={\rm{E}}_{272}\left(a,b_1, b_2; c; x,y,z\right) ,
$

$
{u_2} = {x^{1 - c}}{\rm{E}}_{272}\left(1-c+a,b_1, b_2; 2-c; x,y,z\right).
$

\bigskip

\begin{equation}
{\rm{E}}_{273}\left(a,b_1, b_2; c; x,y,z\right)=\sum\limits_{m,n,p=0}^\infty\frac{(a)_{n+p}(b_1)_{m+n-p}(b_2)_{p-n}}{(c)_{m}}\frac{x^m}{m!}\frac{y^n}{n!}\frac{z^p}{p!},
\end{equation}

region of convergence:
$$ \left\{ r<\infty,\,\,\,s+t<1
\right\}.
$$

System of partial differential equations:

$
\left\{
\begin{aligned}
& {xu_{xx} +(c-x)u_{x} -yu_{y} +zu_{z} -b_{1} u=0,} \\& {y(1+y)u_{yy} -z^{2} u_{zz} +xyu_{xy}+xzu_{xz} -zu_{yz} }\\&\,\,\,\,\,\,\,\,\, +axu_{x} {+\left[1-b_2+\left(a+b_{1} +1\right)y\right]u_{y} -\left(a-b_{1} +1\right)zu_{z} +ab_{1} u=0,} \\ &z(1+z)u_{zz} -y^{2} u_{yy} -xu_{xz} -yu_{yz} -(a-b_{2} +1)yu_{y}{+\left[1-b_{1} +\left(a+b_{2} +1\right)z\right]u_{z} +ab_{2} u=0,}
\end{aligned}
\right.
$

where $u\equiv \,\,{\rm{E}}_{273}\left(a,b_1, b_2; c; x,y,z\right)$.

Particular solutions:

$
{u_1} ={\rm{E}}_{273}\left(a,b_1, b_2; c; x,y,z\right) ,
$

$
{u_2} = {x^{1 - c}}{\rm{E}}_{273}\left(a,1-c+b_1, b_2; 2-c; x,y,z\right).
$

\bigskip

\begin{equation}
{\rm{E}}_{274}\left(a_1,a_2, b; c; x,y,z\right)=\sum\limits_{m,n,p=0}^\infty\frac{(a_1)_{m+n}(a_2)_{n+p}(b)_{m-n-p}}{(c)_{m}}\frac{x^m}{m!}\frac{y^n}{n!}\frac{z^p}{p!},
\end{equation}

region of convergence:
$$ \left\{ \left\{s+2\sqrt{rs}<1\right\}=\left\{\sqrt{s}<\sqrt{1+r}-\sqrt{r}\right\},\,\,\,t<\infty
\right\}.
$$

System of partial differential equations:

$
\left\{
\begin{aligned}
& x(1-x)u_{xx} +y^{2} u_{yy}  +xzu_{xz} +yzu_{yz}\\&\,\,\,\,\,\,\,\,\,+\left[c-\left(a_{1} +b+1\right)x\right]u_{x} {+\left(a_{1} -b+1\right)yu_{y} +a_{1} zu_{z} -a_{1} bu=0,}
 \\& {y(1+y)u_{yy} -x(1-y)u_{xy} +xzu_{xz} +(1+y)zu_{yz} }\\&\,\,\,\,\,\,\,\,\,  +a_{2} xu_{x}+\left[1-b +\left(a_{1} +a_{2} +1\right)y\right]u_{y} +a_{1} zu_{z} +a_{1} a_{2} u=0, \\ &{zu_{zz} -xu_{xz} +yu_{yz} +yu_{y} +(1-b +z)u_{z} +a_{2} u=0,}
\end{aligned}
\right.
$

where $u\equiv \,\,{\rm{E}}_{274}\left(a_1,a_2, b; c; x,y,z\right)$.

Particular solutions:

$
{u_1} ={\rm{E}}_{274}\left(a_1,a_2, b; c; x,y,z\right) ,
$

$
{u_2} = {x^{1 - c}}{\rm{E}}_{274}\left(1-c+a_1,a_2, 1-c+b; 2-c; x,y,z\right).
$

\bigskip

\begin{equation}
{\rm{E}}_{275}\left(a,b_1, b_2; c; x,y,z\right)=\sum\limits_{m,n,p=0}^\infty\frac{(a)_{m+n+p}(b_1)_{n-p}(b_2)_{p-n}}{(c)_m}\frac{x^m}{m!}\frac{y^n}{n!}\frac{z^p}{p!},
\end{equation}

region of convergence:
$$ \left\{ r<\infty,\,\,\,s+t<1
\right\}.
$$

System of partial differential equations:

$
\left\{
\begin{aligned}
& {xu_{xx} +(c-x)u_{x} -yu_{y} -zu_{z} -au=0,}
 \\& {y(1+y)u_{yy} -z^{2} u_{zz} +xyu_{xy} -xzu_{xz} -zu_{yz} } \\& \,\,\,\,\,\,\,\,\, +b_{1} xu_{x}{+\left[1-b_{2} +\left(a+b_{1} +1\right)y\right]u_{y} -\left(a-b_{1} +1\right)zu_{z} +ab_{1} u=0,} \\ &{z(1+z)u_{zz} -y^{2} u_{yy} -xyu_{xy}+xzu_{xz}-yu_{yz}} \\& \,\,\,\,\,\,\,\,\,+b_{2} xu_{x}{-\left(a-b_{2} +1\right)yu_{y} } {+\left[1-b_{1} +\left(a+b_{2} +1\right)z\right]u_{z} +ab_{2} u=0,}
\end{aligned}
\right.
$

where $u\equiv \,\,{\rm{E}}_{275}\left(a,b_1, b_2; c; x,y,z\right)$.

Particular solutions:

$
{u_1} ={\rm{E}}_{275}\left(a,b_1, b_2; c; x,y,z\right) ,
$

$
{u_2} = {x^{1 - c}}{\rm{E}}_{275}\left(1-c+a,b_1, b_2; 2-c; x,y,z\right).
$

\bigskip

\begin{equation}
{\rm{E}}_{276}\left(a,b_1, b_2; c; x,y,z\right)=\sum\limits_{m,n,p=0}^\infty\frac{(a)_{m+n+p}(b_1)_{m-n}(b_2)_{n-p}}{(c)_m}\frac{x^m}{m!}\frac{y^n}{n!}\frac{z^p}{p!},
\end{equation}

region of convergence:
$$ \left\{ \left\{s+2\sqrt{rs}<1\right\}=\left\{\sqrt{s}<\sqrt{1+r}-\sqrt{r}\right\},\,\,\,t<\infty
\right\}.
$$

System of partial differential equations:

$
\left\{
\begin{aligned}
& x(1-x)u_{xx} +y^{2} u_{yy} -xzu_{xz} +yzu_{yz} \\
& \,\,\,\,\,\,\,\,\, +\left[c-\left(a+b_{1} +1\right)x\right]u_{x} {{+\left(a-b_{1}+1 \right)yu_{y} -b_{1} zu_{z} -ab_{1} u=0,}} \\
& {y(1+y)u_{yy} -z^{2} u_{zz}-x(1-y)u_{xy} - xzu_{xz} }\\&
\,\,\,\,\,\,\,\,\, +b_{2} xu_{x}+\left[1-b_{1} +\left(a+b_{2} +1\right)y\right]u_{y} -\left(a-b_{2} +1\right)zu_{z} +ab_{2} u=0, \\
&{zu_{zz} -yu_{yz} +xu_{x} +yu_{y} +(1-b_{2} +z)u_{z} +au=0, }
\end{aligned}
\right.
$

where $u\equiv \,\,{\rm{E}}_{276}\left(a,b_1, b_2; c; x,y,z\right)$.

Particular solutions:

$
{u_1} ={\rm{E}}_{276}\left(a,b_1, b_2; c; x,y,z\right) ,
$

$
{u_2} = {x^{1 - c}}{\rm{E}}_{276}\left(1-c+a,1-c+b_1, b_2; 2-c; x,y,z\right).
$

\bigskip

\begin{equation}
{\rm{E}}_{277}\left(a_1,a_2, b; c; x,y,z\right)=\sum\limits_{m,n,p=0}^\infty\frac{(a_1)_{m+n+p}(a_2)_{p}(b)_{m-n-p}}{(c)_{m}}\frac{x^m}{m!}\frac{y^n}{n!}\frac{z^p}{p!},
\end{equation}

region of convergence:
$$ \left\{ \left\{t+2\sqrt{rt}<1\right\}=\left\{\sqrt{t}<\sqrt{1+r}-\sqrt{r}\right\},\,\,\,s<\infty
\right\}.
$$

System of partial differential equations:

$
\left\{
\begin{aligned}
& x(1-x)u_{xx} +y^{2} u_{yy} +z^2u_{zz}+2yzu_{yz}\\&\,\,\,\,\,\,\,\,\, +\left[c-\left(a_1+b +1\right)x\right]u_{x} {{+\left(a_1-b+1 \right)yu_{y} +\left(a_{1}-b+1\right) zu_{z} -a_1b u=0,}} \\
& {yu_{yy}  -xu_{xy} +zu_{yz} + xu_{x}}+\left(1-b+y\right)u_{y} +zu_{z} +a_{1} u=0, \\
&z(1+z)u_{zz} -x(1-z)u_{xz} +y(1+z)u_{yz}\\& \,\,\,\,\,\,\,\,\,+a_2xu_{x} +a_2yu_{y} +\left[1-b +\left(a_1+a_2+1\right)z\right]u_{z} +a_1a_2u=0,
\end{aligned}
\right.
$

where $u\equiv \,\,{\rm{E}}_{277}\left(a_1,a_2, b; c; x,y,z\right)$.

Particular solutions:

$
{u_1} ={\rm{E}}_{277}\left(a_1,a_2, b; c; x,y,z\right) ,
$

$
{u_2} = {x^{1 - c}}{\rm{E}}_{277}\left(1-c+a_1,a_2, 1-c+b; 2-c; x,y,z\right).
$

\bigskip

\begin{equation}
{\rm{E}}_{278}\left(a_1,a_2, b_1, b_2; c; x,y,z\right)=\sum\limits_{m,n,p=0}^\infty\frac{(a_1)_{n}(a_2)_{p}(b_1)_{p-n}(b_2)_{2m-p}  }{(c)_m}\frac{x^m}{m!}\frac{y^n}{n!}\frac{z^p}{p!},
\end{equation}

region of convergence:
$$ \left\{ t(1+2\sqrt{r})<1,\,\,\,s<\infty
\right\}.
$$

System of partial differential equations:

$
\left\{
\begin{aligned}
& {x(1-4x)u_{xx} -z^{2} u_{zz} +4xzu_{xz}} {+\left[c-\left(4b_{2} +6\right)x\right]u_{x}} { +2b_{2} zu_{z} -b_{2} \left(b_{2} +1\right)u=0,} \\
& {yu_{yy} -zu_{yz} +\left(1-b_{1} +y\right)u_{y} +a_{1} u=0,} \\
&z(1+z)u_{zz} -2xu_{xz} -yzu_{yz} -a_{2} yu_{y}{+\left[1-b_{2} +\left(a_{2} +b_{1} +1\right)z\right]u_{z} +a_{2} b_{1} u=0,}
\end{aligned}
\right.
$

where $u\equiv \,\,{\rm{E}}_{278}\left(a_1,a_2, b_1, b_2; c; x,y,z\right)$.

Particular solutions:

$
{u_1} ={\rm{E}}_{278}\left(a_1,a_2, b_1, b_2; c; x,y,z\right) ,
$

$
{u_2} = {x^{1 - c}}{\rm{E}}_{278}\left(a_1,a_2, b_1, 2-2c+b_2; 2-c; x,y,z\right).
$

\bigskip

\begin{equation}
{\rm{E}}_{279}\left(a_1, a_2, b_1, b_2; c; x,y,z\right)=\sum\limits_{m,n,p=0}^\infty\frac{(a_1)_{n}(a_2)_n(b_1)_{p-n}(b_2)_{2m-p}  }{(c)_m}\frac{x^m}{m!}\frac{y^n}{n!}\frac{z^p}{p!},
\end{equation}

region of convergence:
$$ \left\{ r<\frac{1}{4},\,s<1,\,\,\,t<\infty
\right\}.
$$

System of partial differential equations:

$
\left\{
\begin{aligned}
& {x(1-4x)u_{xx} -z^{2} u_{zz} +4xzu_{xz}+\left[c-\left(4b_{2} +6\right)x\right]u_{x}}{ +2b_{2} zu_{z} -b_{2} \left(b_{2} +1\right)u=0,} \\
& {y(1+y)u_{yy} -zu_{yz} +\left[1-b_{1} +\left(a_{1} +a_{2} +1\right)y\right]u_{y} +a_{1} a_{2} u=0,} \\
&{zu_{zz} -2xu_{xz} -yu_{y} +\left(1-b_{2} +z\right)u_{z} +b_{1} u=0,}
\end{aligned}
\right.
$

where $u\equiv \,\,{\rm{E}}_{279}\left(a_1, a_2, b_1, b_2; c; x,y,z\right)$.

Particular solutions:

$
{u_1} ={\rm{E}}_{279}\left(a_1,a_2, b_1, b_2; c; x,y,z\right) ,
$

$
{u_2} = {x^{1 - c}}{\rm{E}}_{279}\left(a_1,a_2, b_1, 2-2c+b_2; 2-c; x,y,z\right).
$

\bigskip

\begin{equation}
{\rm{E}}_{280}\left(a, b_1, b_2; c; x,y,z\right)=\sum\limits_{m,n,p=0}^\infty\frac{(a)_{n}(b_1)_{p-n}(b_2)_{2m-p}  }{(c)_m}\frac{x^m}{m!}\frac{y^n}{n!}\frac{z^p}{p!},
\end{equation}

region of convergence:
$$ \left\{ r<\frac{1}{4},\,\,\,s<\infty,\,\,\,t<\infty
\right\}.
$$

System of partial differential equations:

$
\left\{
\begin{aligned}
& {x(1-4x)u_{xx} -z^{2} u_{zz} +4xzu_{xz} +\left[c-\left(4b_{2} +6\right)x\right]u_{x}}{ +2b_{2} zu_{z} -b_{2} \left(b_{2} +1\right)u=0,} \\
& {yu_{yy} -zu_{yz} +\left(1-b_{1} +y\right)u_{y} +au=0,} \\
&{zu_{zz} -2xu_{xz} -yu_{y} +\left(1-b_{2} +z\right)u_{z} +b_{1} u=0,}
\end{aligned}
\right.
$

where $u\equiv \,\,{\rm{E}}_{280}\left(a, b_1, b_2; c; x,y,z\right)$.

Particular solutions:

$
{u_1} ={\rm{E}}_{280}\left(a, b_1, b_2; c; x,y,z\right) ,
$

$
{u_2} = {x^{1 - c}}{\rm{E}}_{280}\left(a, b_1, 2-2c+b_2; 2-c; x,y,z\right).
$

\bigskip

\begin{equation}
{\rm{E}}_{281}\left(a, b_1, b_2; c; x,y,z\right)=\sum\limits_{m,n,p=0}^\infty\frac{(a)_{p}(b_1)_{p-n}(b_2)_{2m-p}  }{(c)_m}\frac{x^m}{m!}\frac{y^n}{n!}\frac{z^p}{p!},
\end{equation}

region of convergence:
$$ \left\{ t(1+2\sqrt{r})<1,\,\,\,s<\infty
\right\}.
$$

System of partial differential equations:

$
\left\{
\begin{aligned}
& {x(1-4x)u_{xx} -z^{2} u_{zz} +4xzu_{xz} +\left[c-\left(4b_{2} +6\right)x\right]u_{x}}{ +2b_{2} zu_{z} -b_{2} \left(b_{2} +1\right)u=0,} \\
& {yu_{yy} -zu_{yz} +\left(1-b_{1}\right)u_{y} +u=0,} \\
&{z(1+z)u_{zz} -2xu_{xz} -yzu_{yz} -ayu_{y}} {+\left[1-b_{2} +\left(a+b_{1} +1\right)z\right]u_{z} +ab_{1} u=0,}
\end{aligned}
\right.
$

where $u\equiv \,\,{\rm{E}}_{281}\left(a, b_1, b_2; c; x,y,z\right)$.

Particular solutions:

$
{u_1} ={\rm{E}}_{281}\left(a, b_1, b_2; c; x,y,z\right) ,
$

$
{u_2} = {x^{1 - c}}{\rm{E}}_{281}\left(a, b_1, 2-2c+b_2; 2-c; x,y,z\right).
$

\bigskip

\begin{equation}
{\rm{E}}_{282}\left( a, b; c; x,y,z\right)=\sum\limits_{m,n,p=0}^\infty\frac{(a)_{p-n}(b)_{2m-p}  }{(c)_m}\frac{x^m}{m!}\frac{y^n}{n!}\frac{z^p}{p!},
\end{equation}

region of convergence:
$$ \left\{ r<\frac{1}{4},\,\,\,s<\infty,\,\,\,t<\infty
\right\}.
$$

System of partial differential equations:

$
\left\{
\begin{aligned}
& {x(1-4x)u_{xx} -z^{2} u_{zz} +4xzu_{xz} +\left[c-\left(4b +6\right)x\right]u_{x}}{ +2b zu_{z} -b \left(b+1\right)u=0,} \\
& {yu_{yy} -zu_{yz} +\left(1-a \right)u_{y} +u=0,} \\
&{zu_{zz} -2xu_{xz} -yu_{y} +\left(1-b +z\right)u_{z} +a u=0,}
\end{aligned}
\right.
$

where $u\equiv \,\,{\rm{E}}_{282}\left( a, b; c; x,y,z\right)$.

Particular solutions:

$
{u_1} ={\rm{E}}_{282}\left(a, b; c; x,y,z\right) ,
$

$
{u_2} = {x^{1 - c}}{\rm{E}}_{282}\left(a, 2-2c+b; 2-c; x,y,z\right).
$

\bigskip

\begin{equation}
{\rm{E}}_{283}\left(a_1,a_2, b_1, b_2; c; x,y,z\right)=\sum\limits_{m,n,p=0}^\infty\frac{(a_1)_{m}(a_2)_{n}(b_1)_{m-p}(b_2)_{2p-n}  }{(c)_m}\frac{x^m}{m!}\frac{y^n}{n!}\frac{z^p}{p!},
\end{equation}

region of convergence:
$$ \left\{ 4t(1+r)<1,\,\,\,s<\infty
\right\}.
$$

System of partial differential equations:

$
\left\{
\begin{aligned}
& {x(1-x)u_{xx} +xzu_{xz} +\left[c-\left(a_{1} +b_{1} +1\right)x\right]u_{x} +a_{1} zu_{z} -a_{1} b_{1} u=0,} \\
& {yu_{yy} -2zu_{yz} +\left(1-b_{2} +y\right)u_{y} +a_{2} u=0,} \\
&z(1+4z)u_{zz} +y^{2} u_{yy} -xu_{xz} - 4yzu_{yz} -2b_{2} yu_{y} { +\left[1-b_{1} +\left(4b_{2} +6\right)z\right]u_{z} +b_{2} \left(b_{2} +1\right)u=0,}
\end{aligned}
\right.
$

where $u\equiv \,\,{\rm{E}}_{283}\left(a_1,a_2, b_1, b_2; c; x,y,z\right)$.

Particular solutions:

$
{u_1} ={\rm{E}}_{283}\left(a_1,a_2, b_1, b_2; c; x,y,z\right) ,
$

$
{u_2} = {x^{1 - c}}{\rm{E}}_{283}\left(1-c+a_1,a_2, 1-c+b_1, b_2; 2-c; x,y,z\right).
$

\bigskip

\begin{equation}
{\rm{E}}_{284}\left(a_1, a_2, b_1, b_2; c; x,y,z\right)=\sum\limits_{m,n,p=0}^\infty\frac{(a_1)_{n}(a_2)_n(b_1)_{m-p}(b_2)_{2p-n}  }{(c)_m}\frac{x^m}{m!}\frac{y^n}{n!}\frac{z^p}{p!},
\end{equation}

region of convergence:
$$ \left\{ r<\infty,\,\,\,s(1+2\sqrt{t})<1,
\right\}.
$$

System of partial differential equations:

$
\left\{
\begin{aligned}
& {xu_{xx} +(c-x)u_{x} +zu_{z} -b_{1} u=0,} \\
& {y(1+y)u_{yy} -2zu_{yz} +\left[1-b_{2} +\left(1+a_{1}+a_{2} \right)y\right]u_{y} +a_{1} a_{2} u=0,} \\
&z(1+4z)u_{zz} +y^{2} u_{yy}-xu_{xz} -4yzu_{yz}-2b_{2} yu_{y} {+\left[1-b_{1} +\left(4b_{2} +6\right)z\right]u_{z} +b_{2} \left(b_{2} +1\right)u=0,}
\end{aligned}
\right.
$

where $u\equiv \,\,{\rm{E}}_{284}\left(a_1, a_2, b_1, b_2; c; x,y,z\right)$.

Particular solutions:

$
{u_1} ={\rm{E}}_{284}\left(a_1,a_2, b_1, b_2; c; x,y,z\right) ,
$

$
{u_2} = {x^{1 - c}}{\rm{E}}_{284}\left(a_1,a_2, 1-c+b_1, b_2; 2-c; x,y,z\right).
$

\bigskip

\begin{equation}
{\rm{E}}_{285}\left(a, b_1, b_2; c; x,y,z\right)=\sum\limits_{m,n,p=0}^\infty\frac{(a)_{m}(b_1)_{m-p}(b_2)_{2p-n}  }{(c)_m}\frac{x^m}{m!}\frac{y^n}{n!}\frac{z^p}{p!},
\end{equation}

region of convergence:
$$ \left\{ 4t(1+r)<1,\,\,\,s<\infty
\right\}.
$$

System of partial differential equations:

$
\left\{
\begin{aligned}
& {x(1-x)u_{xx} +xzu_{xz} +\left[c-\left(a+b_{1} +1\right)x\right]u_{x}+azu_z -ab_{1} u=0,} \\
& {yu_{yy} -2zu_{yz} +\left(1-b_{2} \right)u_{y} +u=0,} \\
&{z(1+4z)u_{zz} +y^{2} u_{yy}-xu_{xz} -4yzu_{yz}  } -2b_{2} yu_{y} { +\left[1-b_{1} +\left(4b_{2} +6\right)z\right]u_{z} +b_{2} \left(b_{2} +1\right)u=0,}
\end{aligned}
\right.
$

where $u\equiv \,\,{\rm{E}}_{285}\left(a, b_1, b_2; c; x,y,z\right)$.

Particular solutions:

$
{u_1} ={\rm{E}}_{285}\left(a, b_1, b_2; c; x,y,z\right) ,
$

$
{u_2} = {x^{1 - c}}{\rm{E}}_{285}\left(1-c+a, 1-c+b_1, b_2; 2-c; x,y,z\right).
$

\bigskip

\begin{equation}
{\rm{E}}_{286}\left(a, b_1, b_2; c; x,y,z\right)=\sum\limits_{m,n,p=0}^\infty\frac{(a)_n(b_1)_{m-p}(b_2)_{2p-n}  }{(c)_m}\frac{x^m}{m!}\frac{y^n}{n!}\frac{z^p}{p!},
\end{equation}

region of convergence:
$$ \left\{ t(1+2\sqrt{r})<1,\,\,\,s<\infty
\right\}.
$$

System of partial differential equations:

$
\left\{
\begin{aligned}
& {xu_{xx} +(c-x)u_{x} +zu_{z} -b_{1} u=0,} \\
& {yu_{yy} -2zu_{yz} +\left(1-b_{2} +y\right)u_{y} +au=0,} \\
&{z(1+4z)u_{zz} +y^{2} u_{yy}-xu_{xz} -4yzu_{yz}  }-2b_{2} yu_{y}{+\left[1-b_{1} +\left(4b_{2} +6\right)z\right]u_{z} +b_{2} \left(b_{2} +1\right)u=0,}
\end{aligned}
\right.
$

where $u\equiv \,\,{\rm{E}}_{286}\left(a, b_1, b_2; c; x,y,z\right)$.

Particular solutions:

$
{u_1} ={\rm{E}}_{286}\left(a, b_1, b_2; c; x,y,z\right) ,
$

$
{u_2} = {x^{1 - c}}{\rm{E}}_{286}\left(a, 1-c+b_1, b_2; 2-c; x,y,z\right).
$

\bigskip

\begin{equation}
{\rm{E}}_{287}\left(a, b; c; x,y,z\right)=\sum\limits_{m,n,p=0}^\infty\frac{(a)_{m-p}(b)_{2p-n}  }{(c)_m}\frac{x^m}{m!}\frac{y^n}{n!}\frac{z^p}{p!},
\end{equation}

region of convergence:
$$ \left\{ r<\infty,\,\,\,s<\infty,\,\,\,t<\frac{1}{4}
\right\}.
$$

System of partial differential equations:

$
\left\{
\begin{aligned}
& {xu_{xx} +(c-x)u_{x} +zu_{z} - a u=0,} \\
& {yu_{yy} -2zu_{yz} +\left(1-b \right)u_{y} +u=0,} \\
&{z(1+4z)u_{zz} +y^{2} u_{yy}-xu_{xz} -4yzu_{yz}  }-2b yu_{y}{+\left[1- a+\left(4b+6\right)z\right]u_{z} +b\left(b+1\right)u=0,}
\end{aligned}
\right.
$

where $u\equiv \,\,{\rm{E}}_{287}\left(a, b; c; x,y,z\right)$.

Particular solutions:

$
{u_1} ={\rm{E}}_{287}\left(a, b; c; x,y,z\right) ,
$

$
{u_2} = {x^{1 - c}}{\rm{E}}_{287}\left(1-c+a, b; 2-c; x,y,z\right).
$

\bigskip

\begin{equation}
{\rm{E}}_{288}\left(a_1,a_2,a_3,  b; c; x,y,z\right)=\sum\limits_{m,n,p=0}^\infty\frac{(a_1)_{n}(a_2)_{p}(a_3)_p(b)_{2m-n-p}  }{(c)_m}\frac{x^m}{m!}\frac{y^n}{n!}\frac{z^p}{p!},
\end{equation}

region of convergence:
$$ \left\{ t(1+2\sqrt{r})<1,\,\,\,s<\infty
\right\}.
$$

System of partial differential equations:

$
\left\{
\begin{aligned}
& {x(1-4x)u_{xx} -y^{2} u_{yy}-z^2u_{zz}+4xyu_{xy} }\\&\,\,\,\,\,\,\,\,\, +4xzu_{xz} -2yzu_{yz} {+\left[c-(4b+6)x\right]u_{x} +2byu_{y} +2bzu_{z} -b(b+1)u=0,} \\
& {yu_{yy} -2xu_{xy} +zu_{yz} +(1-b+y)u_{y} +a_{1} u=0,} \\
&{z(1+z)u_{zz} -2xu_{xz} +yu_{yz} +\left[1-b+\left(a_{2} +a_{3} +1\right)z\right]u_{z} +a_{2} a_{3} u=0,}
\end{aligned}
\right.
$

where $u\equiv \,\,{\rm{E}}_{288}\left(a_1,a_2,a_3,  b; c; x,y,z\right)$.

Particular solutions:

$
{u_1} ={\rm{E}}_{288}\left(a_1,a_2,a_3,  b; c; x,y,z\right) ,
$

$
{u_2} = {x^{1 - c}}{\rm{E}}_{288}\left(a_1,a_2,a_3,  2-2c+b; 2-c; x,y,z\right).
$

\bigskip

\begin{equation}
{\rm{E}}_{289}\left(a_1, a_2,  b; c; x,y,z\right)=\sum\limits_{m,n,p=0}^\infty\frac{(a_1)_{n}(a_2)_p(b)_{2m-n-p}  }{(c)_m}\frac{x^m}{m!}\frac{y^n}{n!}\frac{z^p}{p!},
\end{equation}

region of convergence:
$$ \left\{ r<\frac{1}{4},\,\,\,s<\infty,\,\,\,t<\infty
\right\}.
$$

System of partial differential equations:

$
\left\{
\begin{aligned}
& {x(1-4x)u_{xx} -y^{2} u_{yy} -z^2u_{zz}+4xyu_{xy}}\\&\,\,\,\,\,\,\,\,\, +4xzu_{xz}  -2yzu_{yz}{+\left[c-(4b+6)x\right]u_{x} +2byu_{y} +2bzu_{z} -b(b+1)u=0,} \\
& {yu_{yy} -2xu_{xy} +zu_{yz} +(1-b+y)u_{y} +a_{1} u=0,} \\
&{zu_{zz} -2xu_{xz} +yu_{yz} +(1-b+z)u_{z} +a_{2} u=0,}
\end{aligned}
\right.
$

where $u\equiv \,\,{\rm{E}}_{289}\left(a_1, a_2,  b; c; x,y,z\right)$.

Particular solutions:

$
{u_1} ={\rm{E}}_{289}\left(a_1,a_2,  b; c; x,y,z\right) ,
$

$
{u_2} = {x^{1 - c}}{\rm{E}}_{289}\left(a_1,a_2,  2-2c+b; 2-c; x,y,z\right).
$

\bigskip

\begin{equation}
{\rm{E}}_{290}\left(a_1,a_2, b; c; x,y,z\right)=\sum\limits_{m,n,p=0}^\infty\frac{(a_1)_{n}(a_2)_{n}(b)_{2m-n-p}  }{(c)_m}\frac{x^m}{m!}\frac{y^n}{n!}\frac{z^p}{p!},
\end{equation}

region of convergence:
$$ \left\{ s(1+2\sqrt{r})<1,\,\,\,t<\infty
\right\}.
$$

System of partial differential equations:

$
\left\{
\begin{aligned}
& {x(1-4x)u_{xx} -y^{2} u_{yy}-z^2u_{zz}+4xyu_{xy}  } \\&\,\,\,\,\,\,\,\,\,+4xzu_{xz} -2yzu_{yz}{+\left[c-(4b+6)x\right]u_{x} +2byu_{y} +2bzu_{z} -b(b+1)u=0,} \\
& {y(1+y)u_{yy} -2xu_{xy} +zu_{yz} +\left[1-b+\left(a_1+a_2+1\right)y\right]u_{y} +a_1a_2u=0,} \\
&{zu_{zz} -2xu_{xz} +yu_{yz} +(1-b)u_{z}+ u=0,}
\end{aligned}
\right.
$

where $u\equiv \,\,{\rm{E}}_{290}\left(a_1,a_2, b; c; x,y,z\right)$.

Particular solutions:

$
{u_1} ={\rm{E}}_{290}\left(a_1,a_2,  b; c; x,y,z\right) ,
$

$
{u_2} = {x^{1 - c}}{\rm{E}}_{290}\left(a_1,a_2,  2-2c+b; 2-c; x,y,z\right).
$

\bigskip

\begin{equation}
{\rm{E}}_{291}\left(a_1,a_2,b; c; x,y,z\right)=\sum\limits_{m,n,p=0}^\infty\frac{(a_1)_{n+p}(a_2)_{p}(b)_{2m-n-p}}{(c)_{m} }\frac{x^m}{m!}\frac{y^n}{n!}\frac{z^p}{p!},
\end{equation}

region of convergence:
$$ \left\{ t(1+2\sqrt{r})<1,\,\,\,s<\infty
\right\}.
$$

System of partial differential equations:

$
\left\{
\begin{aligned}
& {x(1-4x)u_{xx} -y^{2} u_{yy}}-z^2u_{zz}+4xyu_{xy}\\&\,\,\,\,\,\,\,\,\,  +4xzu_{xz} -2yzu_{yz} {+\left[c-(4b+6)x\right]u_{x} +2byu_{y} +2bzu_{z} -b(b+1)u=0,} \\
& {yu_{yy} -2xu_{xy} +zu_{yz} +(1-b+y)u_{y} +zu_{z} +a_{1} u=0,} \\
&{z(1+z)u_{zz} -2xu_{xz} +y(1+z)u_{yz} }+a_{2} yu_{y} { +\left[1-b+\left(a_{1} +a_{2} +1\right)z\right]u_{z} +a_{1} a_{2} u=0,}
\end{aligned}
\right.
$

where $u\equiv \,\,{\rm{E}}_{291}\left(a_1,a_2,b; c; x,y,z\right)$.

Particular solutions:

$
{u_1} ={\rm{E}}_{291}\left(a_1,a_2,  b; c; x,y,z\right) ,
$

$
{u_2} = {x^{1 - c}}{\rm{E}}_{291}\left(a_1,a_2,  2-2c+b; 2-c; x,y,z\right).
$

\bigskip

\begin{equation}
{\rm{E}}_{292}\left(a,b; c; x,y,z\right)=\sum\limits_{m,n,p=0}^\infty\frac{(a)_n(b)_{2m-n-p}}{(c)_{m}}\frac{x^m}{m!}\frac{y^n}{n!}\frac{z^p}{p!},
\end{equation}

region of convergence:
$$ \left\{ r<\frac{1}{4},\,\,\,s<\infty,\,\,\,t<\infty
\right\}.
$$

System of partial differential equations:

$
\left\{
\begin{aligned}
& {x(1-4x)u_{xx} -y^{2} u_{yy} -z^2u_{zz}+4xyu_{xy}}\\&\,\,\,\,\,\,\,\,\,  +4xzu_{xz} -2yzu_{yz}{+\left[c-(4b+6)x\right]u_{x} +2byu_{y} +2bzu_{z} -b(b+1)u=0,} \\
& {yu_{yy} -2xu_{xy} +zu_{yz} +(1-b+y)u_{y} +au=0,} \\
&{zu_{zz} -2xu_{xz} +yu_{yz} +(1-b)u_{z} +u=0,}
\end{aligned}
\right.
$

where $u\equiv \,\,{\rm{E}}_{292}\left(a,b; c; x,y,z\right)$.

Particular solutions:

$
{u_1} ={\rm{E}}_{292}\left(a,  b; c; x,y,z\right) ,
$

$
{u_2} = {x^{1 - c}}{\rm{E}}_{292}\left(a,  2-2c+b; 2-c; x,y,z\right).
$

\bigskip

\begin{equation}
{\rm{E}}_{293}\left(a, b_1, b_2; c; x,y,z\right)=\sum\limits_{m,n,p=0}^\infty\frac{(a)_{m}(b_1)_{2p-n}(b_2)_{m+n-p}}{(c)_m}\frac{x^m}{m!}\frac{y^n}{n!}\frac{z^p}{p!},
\end{equation}

region of convergence:
$$ \left\{ 4t(1+r)<1,\,\,\,s<\infty
\right\}.
$$

System of partial differential equations:

$
\left\{
\begin{aligned}
& {x(1-x)u_{xx} -xyu_{xy} +xzu_{xz}} {+\left[c-\left(a+b_{2} +1\right)x\right]u_{x} -ayu_{y} +azu_{z} -ab_{2} u=0,} \\
& {yu_{yy} -2zu_{yz} +xu_{x} +\left(1-b_{1} +y\right)u_{y} -zu_{z} +b_{2} u=0,} \\
&{z(1+4z)u_{zz} +y^{2} u_{yy} -xu_{xz} -y(1+4z)u_{yz} } \\& \,\,\,\,\,\,\,\,\,-2b_{1} yu_{y}{ +\left[1-b_{2} +\left(4b_{1} +6\right)z\right]u_{z} +b_{1} \left(b_{1} +1\right)u=0,}
\end{aligned}
\right.
$

where $u\equiv \,\,{\rm{E}}_{293}\left(a, b_1, b_2; c; x,y,z\right)$.

Particular solutions:

$
{u_1} ={\rm{E}}_{293}\left(a, b_1, b_2; c; x,y,z\right) ,
$

$
{u_2} = {x^{1 - c}}{\rm{E}}_{293}\left(1-c+a, b_1, 1-c+b_2; 2-c; x,y,z\right).
$

\bigskip

\begin{equation}
{\rm{E}}_{294}\left(a, b_1, b_2; c; x,y,z\right)=\sum\limits_{m,n,p=0}^\infty\frac{(a)_{n}(b_1)_{2p-n}(b_2)_{m+n-p}}{(c)_m}\frac{x^m}{m!}\frac{y^n}{n!}\frac{z^p}{p!},
\end{equation}

region of convergence:
$$ \left\{ r<\infty,\,\,\,\left\{ts^2+s<1\right\}=\left\{s<\frac{\sqrt{1+4t}-1}{2t}\right\}
\right\}.
$$

System of partial differential equations:

$
\left\{
\begin{aligned}
& {xu_{xx} +(c-x)u_{x} -yu_{y} +zu_{z} -b_{2} u=0,} \\
& {y(1+y)u_{yy} +xyu_{xy} -(y+2)zu_{yz} +axu_{x}} { +\left[1-b_{1} +\left(a+b_{2} +1\right)y\right]u_{y} -azu_{z} +ab_{2} u=0,} \\
&{z(1+4z)u_{zz} +y^{2} u_{yy} -xu_{xz} -y(1+4z)u_{yz}}\\& \,\,\,\,\,\,\,\,\,-2b_{1} yu_{y}{  +\left[1-b_{2} +\left(4b_{1} +6\right)z\right]u_{z} +b_{1} \left(b_{1} +1\right)u=0,}
\end{aligned}
\right.
$

where $u\equiv \,\,{\rm{E}}_{294}\left(a, b_1, b_2; c; x,y,z\right)$.

Particular solutions:

$
{u_1} ={\rm{E}}_{294}\left(a, b_1, b_2; c; x,y,z\right) ,
$

$
{u_2} = {x^{1 - c}}{\rm{E}}_{294}\left(a, b_1, 1-c+b_2; 2-c; x,y,z\right).
$

\bigskip

\begin{equation}
{\rm{E}}_{295}\left(a, b; c; x,y,z\right)=\sum\limits_{m,n,p=0}^\infty\frac{(a)_{2p-n}(b)_{m+n-p}}{(c)_m}\frac{x^m}{m!}\frac{y^n}{n!}\frac{z^p}{p!},
\end{equation}

region of convergence:
$$ \left\{ r<\infty,\,\,\,s<\infty,\,\,\,t<\frac{1}{4}
\right\}.
$$

System of partial differential equations:

$
\left\{
\begin{aligned}
& {xu_{xx} +(c-x)u_{x} -yu_{y} +zu_{z} -b u=0,} \\
& {yu_{yy} -2zu_{yz} +xu_{x} +\left(1- a +y\right)u_{y} -zu_{z} +b u=0,} \\
&{z(1+4z)u_{zz} +y^{2} u_{yy} -xu_{xz} -y(1+4z)u_{yz} }-2 ayu_{y}{ +\left[1-b +\left(4a+6\right)z\right]u_{z} +a \left(a +1\right)u=0,}
\end{aligned}
\right.
$

where $u\equiv \,\,{\rm{E}}_{295}\left(a, b; c; x,y,z\right)$.

Particular solutions:

$
{u_1} ={\rm{E}}_{295}\left(a, b; c; x,y,z\right) ,
$

$
{u_2} = {x^{1 - c}}{\rm{E}}_{295}\left(a, 1-c+b; 2-c; x,y,z\right).
$

\bigskip

\begin{equation}
{\rm{E}}_{296}\left(a, b_1, b_2; c; x,y,z\right)=\sum\limits_{m,n,p=0}^\infty\frac{(a)_{m+n}(b_1)_{2p-n}(b_2)_{n-p}}{(c)_m}\frac{x^m}{m!}\frac{y^n}{n!}\frac{z^p}{p!},
\end{equation}

region of convergence:
$$ \left\{ r<\infty,\,\,\,\left\{ts^2+s<1\right\}=\left\{s<\frac{\sqrt{1+4t}-1}{2t}\right\}
\right\}.
$$

System of partial differential equations:

$
\left\{
\begin{aligned}
& {xu_{xx} +(c-x)u_{x} -yu_{y} -au=0,} \\
& {y(1+y)u_{yy} +xyu_{xy} -xzu_{xz} -(y+2)zu_{yz}}\\& \,\,\,\,\,\,\,\,\, +b_{2} xu_{x}{ +\left[1-b_{1} +\left(a+b_{2} +1\right)y\right]u_{y} -azu_{z} +ab_{2} u=0,} \\
&{z(1+4z)u_{zz} +y^{2} u_{yy} -y(1+4z)u_{yz} } -2b_{1} yu_{y}{ +\left[1-b_{2} +\left(4b_{1} +6\right)z\right]u_{z} +b_{1} \left(b_{1} +1\right)u=0,}
\end{aligned}
\right.
$

where $u\equiv \,\,{\rm{E}}_{296}\left(a, b_1, b_2; c; x,y,z\right)$.

Particular solutions:

$
{u_1} ={\rm{E}}_{296}\left(a, b_1, b_2; c; x,y,z\right) ,
$

$
{u_2} = {x^{1 - c}}{\rm{E}}_{296}\left(1-c+a, b_1, b_2; 2-c; x,y,z\right).
$

\bigskip

\begin{equation}
{\rm{E}}_{297}\left(a, b_1, b_2; c; x,y,z\right)=\sum\limits_{m,n,p=0}^\infty\frac{(a)_{m+n}(b_1)_{2p-n}(b_2)_{m-p}}{(c)_m}\frac{x^m}{m!}\frac{y^n}{n!}\frac{z^p}{p!},
\end{equation}

region of convergence:
$$ \left\{ 4t(1+r)<1,\,\,\,s<\infty
\right\}.
$$

System of partial differential equations:

$
\left\{
\begin{aligned}
& {x(1-x)u_{xx} -xyu_{xy} +xzu_{xz} +yzu_{yz}}{ +\left[c-\left(a+b_{2} +1\right)x\right]u_{x} -b_{2} yu_{y} +azu_{z} -ab_{2} u=0,} \\
& {yu_{yy} -2zu_{yz} +xu_{x} +\left(1-b_{1} +y\right)u_{y} +au=0,} \\
&{z(1+4z)u_{zz} +y^{2} u_{yy}-xzu_{xz}-4yzu_{yz} -2b_{1} yu_{y}} { +\left[1-b_{2} +\left(4b_{1} +6\right)z\right]u_{z} +b_{1} \left(b_{1} +1\right)u=0,}
\end{aligned}
\right.
$

where $u\equiv \,\,{\rm{E}}_{297}\left(a, b_1, b_2; c; x,y,z\right)$.

Particular solutions:

$
{u_1} ={\rm{E}}_{297}\left(a, b_1, b_2; c; x,y,z\right) ,
$

$
{u_2} = {x^{1 - c}}{\rm{E}}_{297}\left(1-c+a, b_1, 1-c+b_2; 2-c; x,y,z\right).
$

\bigskip

\begin{equation}
{\rm{E}}_{298}\left(a, b_1, b_2; c; x,y,z\right)=\sum\limits_{m,n,p=0}^\infty\frac{(a)_{n+p}(b_1)_{2m-p}(b_2)_{p-n}}{(c)_m}\frac{x^m}{m!}\frac{y^n}{n!}\frac{z^p}{p!},
\end{equation}

region of convergence:
$$ \left\{ t(1+2\sqrt{r})<1,\,\,\,s<\infty
\right\}.
$$

System of partial differential equations:

$
\left\{
\begin{aligned}
& {x(1-4x)u_{xx} -z^{2} u_{zz} +4xzu_{xz}}{ +\left[c-\left(4b_{1} +6\right)x\right]u_{x} +2b_{1} zu_{z} -b_{1} \left(b_{1} +1\right)u=0,} \\
& {yu_{yy} -zu_{yz} +\left(1-b_{2} +y\right)u_{y} +zu_{z} +au=0,} \\
&{z(1+z)u_{zz} -y^2u_{yy}-2xu_{xz}  -\left(a-b_2+1\right) yu_{y}}{ +\left[1-b_1+\left(a +b_{2} +1\right)z\right]u_{z} +a b_{2} u=0,}
\end{aligned}
\right.
$

where $u\equiv \,\,{\rm{E}}_{298}\left(a, b_1, b_2; c; x,y,z\right)$.

Particular solutions:

$
{u_1} ={\rm{E}}_{298}\left(a, b_1, b_2; c; x,y,z\right) ,
$

$
{u_2} = {x^{1 - c}}{\rm{E}}_{298}\left(a, 2-2c+b_1, b_2; 2-c; x,y,z\right).
$

\bigskip

\begin{equation}
{\rm{E}}_{299}\left(a, b_1, b_2; c; x,y,z\right)=\sum\limits_{m,n,p=0}^\infty\frac{(a)_{n}(b_1)_{2m+n-p}(b_2)_{p-n}}{(c)_m}\frac{x^m}{m!}\frac{y^n}{n!}\frac{z^p}{p!},
\end{equation}

region of convergence:
$$ \left\{ 2\sqrt{r}+s<1,\,\,\,t<\infty
\right\}.
$$

System of partial differential equations:

$
\left\{
\begin{aligned}
& {x(1-4x)u_{xx} -y^2u_{xx}-z^{2} u_{zz}-4xyu_{xy} +4xzu_{xz} }\\
& \,\,\,\,\,\,\,\,\,+2yzu_{yz}{+\left[c-\left(4b_{1} +6\right)x\right]u_{x}-2\left(b_1+1\right)yu_y+2b_{1} zu_{z} -b_{1} \left(b_{1} +1\right)u=0,} \\
& {y(1+y)u_{yy} +2xyu_{xy}-(1+y)zu_{yz} +2axu_x} {+\left[1-b_{2} +\left(1+a+b_1\right)y\right]u_{y} -azu_{z} +ab_1u=0,} \\
&{zu_{zz} -2xu_{xz} -yu_{yz}-yu_{y} +(1-b_{1} +z)u_{z} +b_{2} u=0,}
\end{aligned}
\right.
$

where $u\equiv \,\,{\rm{E}}_{299}\left(a, b_1, b_2; c; x,y,z\right)$.

Particular solutions:

$
{u_1} ={\rm{E}}_{299}\left(a, b_1, b_2; c; x,y,z\right) ,
$

$
{u_2} = {x^{1 - c}}{\rm{E}}_{299}\left(a, 2-2c+b_1, b_2; 2-c; x,y,z\right).
$

\bigskip

\begin{equation}
{\rm{E}}_{300}\left(a, b_1, b_2; c; x,y,z\right)=\sum\limits_{m,n,p=0}^\infty\frac{(a)_{p}(b_1)_{2m+n-p}(b_2)_{p-n}}{(c)_m}\frac{x^m}{m!}\frac{y^n}{n!}\frac{z^p}{p!},
\end{equation}

region of convergence:
$$ \left\{ t(1+2\sqrt{r})<1,\,\,\,s<\infty
\right\}.
$$

System of partial differential equations:

$
\left\{
\begin{aligned}
& {x(1-4x)u_{xx} -y^2u_{xx}-z^{2} u_{zz} -4xyu_{xy}+4xzu_{xz} }\\& \,\,\,\,\,\,\,\,\,+2yzu_{yz}{+\left[c-\left(4b_{1} +6\right)x\right]u_{x}-2\left(b_1+1\right)yu_y+2b_{1} zu_{z} -b_{1} \left(b_{1} +1\right)u=0,} \\
& {yu_{yy} -zu_{yz}+2xu_x+\left(1-b_2+y\right)u_{y} -zu_{z} +b_1u=0,} \\
&{z(1+z)u_{zz} -2xzu_{xz} -y(1+z)u_{yz}-ayu_{y}}{ +\left[1-b_{1} +\left(1+a+b_2\right)z\right]u_{z} +ab_{2} u=0,}
\end{aligned}
\right.
$

where $u\equiv \,\,{\rm{E}}_{300}\left(a, b_1, b_2; c; x,y,z\right)$.

Particular solutions:

$
{u_1} ={\rm{E}}_{300}\left(a, b_1, b_2; c; x,y,z\right) ,
$

$
{u_2} = {x^{1 - c}}{\rm{E}}_{300}\left(a, 2-2c+b_1, b_2; 2-c; x,y,z\right).
$

\bigskip

\begin{equation}
{\rm{E}}_{301}\left(a, b_1, b_2; c; x,y,z\right)=\sum\limits_{m,n,p=0}^\infty\frac{(a)_{p}(b_1)_{m+2n-p}(b_2)_{m-n}}{(c)_m}\frac{x^m}{m!}\frac{y^n}{n!}\frac{z^p}{p!},
\end{equation}

region of convergence:
$$ \left\{\left\{ s<\frac{1}{4}, \,\,\,\, r<\min\left\{\Psi_1(s), \Psi_2(s)\right\}\right\}= \left\{r<1 \wedge s<\min\left\{\Theta_1(r), \Theta_2(r)\right\}\right\},\,\,\,t<\infty
\right\}.
$$

System of partial differential equations:

$
 \left\{\begin{aligned}
 & {x(1-x)u_{xx} +2y^{2} u_{yy} -xyu_{xy}+xzu_{xz} -yzu_{yz}} \\
 &{\,\,\,\,\,\,\,\,\,+\left[c-\left(b_{1} +b_{2} +1\right)x\right]u_{x} +\left(b_{1} -2b_{2} +2\right)yu_{y} +b_{2} zu_{z} -b_{1} b_{2} u=0,} \\
 & { y(1+4y)u_{yy}+x^{2} u_{xx} +z^{2} u_{zz} -x(1-4y)u_{xy} -2xzu_{xz}  } \\
 & { \,\,\,\,\,\,\,\,\,-4yzu_{yz}+2\left(b_{1} +1\right)xu_{x} +\left[1-b_{2} +\left(4b_{1} +6\right)y\right]u_{y} -2b_{1} zu_{z} +b_{1} \left(b_{1} +1\right)u=0,} \\
 & {zu_{zz} -xu_{xz} -2yu_{yz} +(1-b_{1} +z)u_{z} +au=0,} \end{aligned}\right.
$

where $u\equiv \,\,{\rm{E}}_{301}\left(a, b_1, b_2; c; x,y,z\right)$.

Particular solutions:

$
{u_1} ={\rm{E}}_{301}\left(a, b_1, b_2; c; x,y,z\right) ,
$

$
{u_2} = {x^{1 - c}}{\rm{E}}_{301}\left(a, 1-c+b_1, 1-c+b_2; 2-c; x,y,z\right).
$

\bigskip

\begin{equation}
{\rm{E}}_{302}\left(a, b; c; x,y,z\right)=\sum\limits_{m,n,p=0}^\infty\frac{(a)_{m+2n-p}(b)_{m-n}}{(c)_m}\frac{x^m}{m!}\frac{y^n}{n!}\frac{z^p}{p!},
\end{equation}

region of convergence:
$$ \left\{ \left\{s<\frac{1}{4},\,\,\,\, r<\min\left\{\Psi_1(s), \Psi_2(s)\right\}\right\}= \left\{r<1, \,\,\,\, s<\min\left\{\Theta_1(r), \Theta_2(r)\right\}\right\},\,\,\,t<\infty
\right\}.
$$

System of partial differential equations:

$
\left\{\begin{aligned}
& {x(1-x)u_{xx} +2y^{2} u_{yy}-xyu_{xy} +xzu_{xz}   } \\
&{\,\,\,\,\,\,\,\,\,-yzu_{yz}+\left[c-\left(a +b +1\right)x\right]u_{x} +\left(a -2b +2\right)yu_{y} +b zu_{z} -a b u=0,} \\
& { y(1+4y)u_{yy}+x^{2} u_{xx} +z^{2} u_{zz} -x(1-4y)u_{xy} -2xzu_{xz}} \\
 & { \,\,\,\,\,\,\,\,\,-4yzu_{yz}+2\left(a +1\right)xu_{x} +\left[1-b +\left(4a +6\right)y\right]u_{y} -2a zu_{z} + a\left(a +1\right)u=0,} \\
& {zu_{zz} -xu_{xz} -2yu_{yz}  +\left(1-a\right)u_{z} + u=0,} \end{aligned}\right.
$

where $u\equiv \,\,{\rm{E}}_{302}\left(a, b; c; x,y,z\right)$.

Particular solutions:

$
{u_1} ={\rm{E}}_{302}\left(a, b; c; x,y,z\right) ,
$

$
{u_2} = {x^{1 - c}}{\rm{E}}_{302}\left(1-c+a, 1-c+b; 2-c; x,y,z\right).
$

\bigskip

\begin{equation}
{\rm{E}}_{303}\left(a, b_1, b_2; c; x,y,z\right)=\sum\limits_{m,n,p=0}^\infty\frac{(a)_{m}(b_1)_{m+2n-p}(b_2)_{p-n}}{(c)_m}\frac{x^m}{m!}\frac{y^n}{n!}\frac{z^p}{p!},
\end{equation}

region of convergence:
$$ \left\{ r+2\sqrt{s}<1,\,\,\,t<\infty
\right\}.
$$

System of partial differential equations:

$
\left\{\begin{aligned}
& {x(1-x)u_{xx}  -2xyu_{xy} +xzu_{xz}} {+\left[c-\left(a +b_{1} +1\right)x\right]u_{x}  -2ayu_{y} +azu_{z} -ab_{1} u=0,} \\
& { y(1+4y)u_{yy}+x^{2} u_{xx} +z^{2} u_{zz}+ 4xyu_{xy} } -2xzu_{xz}\\
& \,\,\,\,\,\, -(1+4y)zu_{yz}{+2\left(b_{1} +1\right)xu_{x} +\left[1-b_{2} +\left(4b_{1} +6\right)y\right]u_{y} -2b_{1} zu_{z} +b_{1} \left(b_{1} +1\right)u=0,} \\
&{zu_{zz} -xu_{xz} -2yu_{yz} -yu_y +(1-b_{1} +z)u_{z} + b_2u=0,} \end{aligned}\right.
$

where $u\equiv \,\,{\rm{E}}_{303}\left(a, b_1, b_2; c; x,y,z\right)$.

Particular solutions:

$
{u_1} ={\rm{E}}_{303}\left(a, b_1, b_2; c; x,y,z\right) ,
$

$
{u_2} = {x^{1 - c}}{\rm{E}}_{303}\left(1-c+a, 1-c+b_1, b_2; 2-c; x,y,z\right).
$

\bigskip

\begin{equation}
{\rm{E}}_{304}\left(a, b_1, b_2; c; x,y,z\right)=\sum\limits_{m,n,p=0}^\infty\frac{(a)_{p}(b_1)_{m+2n-p}(b_2)_{p-n}}{(c)_m}\frac{x^m}{m!}\frac{y^n}{n!}\frac{z^p}{p!},
\end{equation}

region of convergence:
$$ \left\{ r<\infty,\,\,\,\left\{st^2+t<1\right\}=\left\{t<\frac{\sqrt{1+4s}-1}{2s}\right\}
\right\}.
$$

System of partial differential equations:

$
\left\{\begin{aligned}
& {xu_{xx} +(c-x)u_{x}  -2yu_{y} +zu_{z} -b_{1} u=0,} \\
& { y(1+4y)u_{yy}+x^{2} u_{xx} +z^{2} u_{zz}+ 4xyu_{xy} }-2xzu_{xz} \\
& \,\,\,\,\,\, -(1+4y)zu_{yz}{+2\left(b_{1} +1\right)xu_{x} +\left[1-b_{2} +\left(4b_{1} +6\right)y\right]u_{y} -2b_{1} zu_{z} +b_{1} \left(b_{1} +1\right)u=0,} \\
& {z(1+z)u_{zz} -xu_{xz} -y(2+z)u_{yz} -ayu_y} +\left[1-b_{1} +\left(a+b_2+1\right)z\right]u_{z} +ab_2u=0, \end{aligned}\right.
$

where $u\equiv \,\,{\rm{E}}_{304}\left(a, b_1, b_2; c; x,y,z\right)$.

Particular solutions:

$
{u_1} ={\rm{E}}_{304}\left(a, b_1, b_2; c; x,y,z\right) ,
$

$
{u_2} = {x^{1 - c}}{\rm{E}}_{304}\left(a, 1-c+b_1, b_2; 2-c; x,y,z\right).
$

\bigskip

\begin{equation}
{\rm{E}}_{305}\left(a, b; c; x,y,z\right)=\sum\limits_{m,n,p=0}^\infty\frac{(a)_{m+2n-p}(b)_{p-n}}{(c)_m}\frac{x^m}{m!}\frac{y^n}{n!}\frac{z^p}{p!},
\end{equation}

region of convergence:
$$ \left\{ r<\infty,\,\,\,s<\frac{1}{4},\,\,\,t<\infty
\right\}.
$$

System of partial differential equations:

$
\left\{\begin{aligned}
& {xu_{xx} +(c-x)u_{x}  -2yu_{y} +zu_{z} - a u=0,} \\
& { y(1+4y)u_{yy}+x^{2} u_{xx} +z^{2} u_{zz}}+ 4xyu_{xy}-2xzu_{xz} \\
& \,\,\,\,\,\,  -(1+4y)zu_{yz}{+2\left(a +1\right)xu_{x} +\left[1-b +\left(4a +6\right)y\right]u_{y} -2a zu_{z} +a\left(a+1\right)u=0,} \\
& {zu_{zz} -xu_{xz} -2yu_{yz} -yu_y +\left(1-a+z\right)u_{z} + bu=0,} \end{aligned}\right.
$

where $u\equiv \,\,{\rm{E}}_{305}\left(a, b; c; x,y,z\right)$.

Particular solutions:

$
{u_1} ={\rm{E}}_{305}\left(a, b; c; x,y,z\right) ,
$

$
{u_2} = {x^{1 - c}}{\rm{E}}_{305}\left(1-c+a, b; 2-c; x,y,z\right).
$

\bigskip

\begin{equation}
{\rm{E}}_{306}\left(a, b_1, b_2; c; x,y,z\right)=\sum\limits_{m,n,p=0}^\infty\frac{(a)_{m+2n}(b_1)_{m-p}(b_2)_{p-n}}{(c)_m}\frac{x^m}{m!}\frac{y^n}{n!}\frac{z^p}{p!},
\end{equation}

region of convergence:
$$ \left\{ r+2\sqrt{s}<1,\,\,\,t<\infty
\right\}.
$$

System of partial differential equations:

$
\left\{\begin{aligned}
& {x(1-x)u_{xx} -2xyu_{xy} +xzu_{xz} +2yzu_{yz} }  {+\left[c-\left(a+b_{1} +1\right)x\right]u_{x} -2b_{1} yu_{y} +azu_{z} -ab_{1} u=0,} \\
& {y(1+4y)u_{yy} +x^{2} u_{xx} +4xyu_{xy} -zu_{yz} } +2(a+1)xu_{x} {+\left[1-b_{2} +(4a+6)y \right]u_{y} +a(a+1)u=0,} \\
& {zu_{zz} -xu_{xz} -yu_{y} +(1-b_{1} +z)u_{z} +b_{2} u=0,} \end{aligned}\right.
$

where $u\equiv \,\,{\rm{E}}_{306}\left(a, b_1, b_2; c; x,y,z\right)$.

Particular solutions:

$
{u_1} ={\rm{E}}_{306}\left(a, b_1, b_2; c; x,y,z\right) ,
$

$
{u_2} = {x^{1 - c}}{\rm{E}}_{306}\left(1-c+a, 1-c+b_1, b_2; 2-c; x,y,z\right).
$

\bigskip

\begin{equation}
{\rm{E}}_{307}\left(a, b_1, b_2; c; x,y,z\right)=\sum\limits_{m,n,p=0}^\infty\frac{(a)_{2m+n}(b_1)_{n-p}(b_2)_{p-n}}{(c)_m}\frac{x^m}{m!}\frac{y^n}{n!}\frac{z^p}{p!},
\end{equation}

region of convergence:
$$ \left\{ s+2\sqrt{r}<1,\,\,\,t<\infty
\right\}.
$$

System of partial differential equations:

$
\left\{\begin{aligned}
& {x(1-4x)u_{xx} -y^{2} u_{yy} -4xyu_{xy} +\left[c-(4a+6)x\right]u_{x}} {-2(a+1)yu_{y} -a(a+1)u=0,} \\
& {y(1+y)u_{yy} +2xyu_{xy} -2xzu_{xz} -(1+y)zu_{yz} } \\& \,\,\,\,\,\,\,\,\, +2b_{1} xu_{x}{+\left[1-b_{2} +\left(a+b_{1} +1\right)y\right]u_{y} -azu_{z} +ab_1u=0,} \\
& {zu_{zz} -yu_{yz} -yu_{y} +(1-b_{1} +z)u_{z} +b_{2} u=0,} \end{aligned}\right.
$

where $u\equiv \,\,{\rm{E}}_{307}\left(a, b_1, b_2; c; x,y,z\right)$.

Particular solutions:

$
{u_1} ={\rm{E}}_{307}\left(a, b_1, b_2; c; x,y,z\right) ,
$

$
{u_2} = {x^{1 - c}}{\rm{E}}_{307}\left(2-2c+a, b_1, b_2; 2-c; x,y,z\right).
$

\bigskip

\begin{equation}
{\rm{E}}_{308}\left(a, b_1, b_2; c; x,y,z\right)=\sum\limits_{m,n,p=0}^\infty\frac{(a)_{n+2p}(b_1)_{m-n}(b_2)_{n-p}}{(c)_m}\frac{x^m}{m!}\frac{y^n}{n!}\frac{z^p}{p!},
\end{equation}

region of convergence:
$$ \left\{ r<\infty,\,\,\,\left\{t<\frac{1}{4}, \,\,\, s<\min\left\{\Psi_1(t), \Psi_2(t)\right\}\right\}=\left\{s<1, \,\,\, t<\min\left\{\Theta_1(s), \Theta_2(s)\right\}\right\}
\right\}.
$$

System of partial differential equations:

$
\left\{\begin{aligned} &{xu_{xx} +(c-x)u_{x} +yu_y-b_1u=0,} \\
& {y(1+y)u_{yy} -2z^2u_{zz} -xu_{xy} +yzu_{yz} } \\& \,\,\,\,\,\,\,\,\,+\left[1-b_1+\left(a+b_2+1\right)y\right]u_{y}{- (a-2b_2+2)zu_{z} +ab_2u=0,}\\
& {z(1+4z)u_{zz} +y^2u_{yy} -y(1-4z)u_{yz} } +2(a+1)yu_{y}+\left[1-b_{2} +(4a+6)z\right]u_{z} +a(a+1)u=0, \end{aligned}\right.
$

where $u\equiv \,\,{\rm{E}}_{308}\left(a, b_1, b_2; c; x,y,z\right)$.

Particular solutions:

$
{u_1} ={\rm{E}}_{308}\left(a, b_1, b_2; c; x,y,z\right) ,
$

$
{u_2} = {x^{1 - c}}{\rm{E}}_{308}\left(a, 1-c+b_1, b_2; 2-c; x,y,z\right).
$

\bigskip

\begin{equation}
{\rm{E}}_{309}\left(a, b_1, b_2; c; x,y,z\right)=\sum\limits_{m,n,p=0}^\infty\frac{(a)_{n+2p}(b_1)_{m-n}(b_2)_{m-p}}{(c)_m}\frac{x^m}{m!}\frac{y^n}{n!}\frac{z^p}{p!},
\end{equation}

region of convergence:
$$ \left\{ t(1+r)<\frac{1}{4},\,\,\,s<\infty
\right\}.
$$

System of partial differential equations:

$
 \left\{\begin{aligned}& {x(1-x)u_{xx} +xyu_{xy} +xzu_{xz} -yzu_{yz} }  +\left[c-\left(b_{1} +b_{2} +1\right)x\right]u_{x}{+b_{2} yu_{y} +b_{1} zu_{z} -b_{1} b_{2} u=0,} \\
& {yu_{yy} -xu_{xy} +(1-b_{1} +y)u_{y} +2zu_{z} +au=0,} \\
& { z(1+4z)u_{zz} +y^{2} u_{yy}- xu_{xz} +4yzu_{yz} } +2(a+1)yu_{y} {+\left[1-b_{2} +(4a+6)z\right]u_{z} +a(a+1)u=0,} \end{aligned}\right.
$

where $u\equiv \,\,{\rm{E}}_{309}\left(a, b_1, b_2; c; x,y,z\right)$.

Particular solutions:

$
{u_1} ={\rm{E}}_{309}\left(a, b_1, b_2; c; x,y,z\right) ,
$

$
{u_2} = {x^{1 - c}}{\rm{E}}_{309}\left(a, 1-c+b_1, 1-c+b_2; 2-c; x,y,z\right).
$

\bigskip

\begin{equation}
{\rm{E}}_{310}\left(a_1,a_2,b; c; x,y,z\right)=\sum\limits_{m,n,p=0}^\infty\frac{(a_1)_{n+2p}(a_2)_{m}(b)_{m-n-p}}{(c)_{m}}\frac{x^m}{m!}\frac{y^n}{n!}\frac{z^p}{p!},
\end{equation}

region of convergence:
$$ \left\{ t(1+r)<\frac{1}{4},\,\,\,s<\infty
\right\}.
$$

System of partial differential equations:

$
\left\{\begin{aligned}
& {x(1-x)u_{xx} +xyu_{xy} +xzu_{xz} +\left[c-\left(a_{2} +b+1\right)x\right]u_{x} } {+a_{2} yu_{y} +a_{2} zu_{z} -a_{2} bu=0,} \\
& {yu_{yy} -xu_{xy} +zu_{yz} +(1-b+y)u_{y} +2zu_{z} +a_{1} u=0,} \\
& { z(1+4z)u_{zz}+y^{2} u_{yy} -xu_{xz} +y(1+4z)u_{yz} }\\& \,\,\,\,\,\,\,\,\,+2(a_{1} +1)yu_{y} {+\left[1-b+(4a_{1} +6)z\right]u_{z} +a_{1} \left(a_{1} +1\right)u=0,} \end{aligned}\right.
$

where $u\equiv \,\,{\rm{E}}_{310}\left(a_1,a_2,b; c; x,y,z\right)$.

Particular solutions:

$
{u_1} ={\rm{E}}_{310}\left(a_1,a_2,b; c; x,y,z\right) ,
$

$
{u_2} = {x^{1 - c}}{\rm{E}}_{310}\left(a_1,1-c+a_2,1-c+b; 2-c; x,y,z\right).
$

\bigskip

\begin{equation}
{\rm{E}}_{311}\left(a_1,a_2,b; c; x,y,z\right)=\sum\limits_{m,n,p=0}^\infty\frac{(a_1)_{n+2p}(a_2)_n(b)_{m-n-p}}{(c)_{m}}\frac{x^m}{m!}\frac{y^n}{n!}\frac{z^p}{p!},
\end{equation}

region of convergence:
$$ \left\{ r<\infty,\,\,\,\left\{t<\frac{1}{4}, \,\,\,\, s<\frac{1}{2}+\frac{1}{2}\sqrt{1-4t}\right\}\cup\left\{s\leq\frac{1}{2}\right\}
\right\}.
$$

System of partial differential equations:

$
\left\{\begin{aligned}& {xu_{xx} +(c-x)u_{x} +yu_{y} +zu_{z} -bu=0,} \\
& {y(1+y)u_{yy} -xu_{xy} +(1+2y)zu_{yz}} {+(1-b+y(a_{1} +a_{2} +1))u_{y} +2a_{2} zu_{z} +a_{1} a_{2} u=0,} \\
& { z(1+4z)u_{zz} +y^{2} u_{yy}-xu_{xz} +y(1+4z)u_{yz} }\\& \,\,\,\,\,\,\,\,\,+2(a_{1} +1)yu_{y} {+\left[1-b+\left(4a_{1} +6\right)z\right]u_{z} +a_1(a_1+1)u=0,} \end{aligned}\right.
$

where $u\equiv \,\,{\rm{E}}_{311}\left(a_1,a_2,b; c; x,y,z\right)$.

Particular solutions:

$
{u_1} ={\rm{E}}_{311}\left(a_1,a_2,b; c; x,y,z\right) ,
$

$
{u_2} = {x^{1 - c}}{\rm{E}}_{311}\left(a_1,a_2,1-c+b; 2-c; x,y,z\right).
$

\bigskip

\begin{equation}
{\rm{E}}_{312}\left(a,b; c; x,y,z\right)=\sum\limits_{m,n,p=0}^\infty\frac{(a)_{n+2p}(b)_{m-n-p}}{(c)_{m} }\frac{x^m}{m!}\frac{y^n}{n!}\frac{z^p}{p!},
\end{equation}

region of convergence:
$$ \left\{ r<\infty,\,\,\,s<\infty,\,\,\,t<\frac{1}{4}
\right\}.
$$

System of partial differential equations:

$
\left\{\begin{aligned}
& {xu_{xx} +(c-x)u_{x} +yu_{y} +zu_{z} -bu=0,} \\
& {yu_{yy} -xu_{xy} +zu_{yz} +(1-b+y)u_{y} +2zu_{z} +a u=0,} \\
& { z(1+4z)u_{zz} +y^{2} u_{yy}-xu_{xz} +y(1+4z)u_{yz} }\\& \,\,\,\,\,\,\,\,\,+2(a +1)yu_{y} {+\left[1-b+\left(4a +6\right)z\right]u_{z} +a(a+1)u=0,}
 \end{aligned}\right.
$

where $u\equiv \,\,{\rm{E}}_{312}\left(a,b; c; x,y,z\right)$.

Particular solutions:

$
{u_1} ={\rm{E}}_{312}\left(a,b; c; x,y,z\right) ,
$

$
{u_2} = {x^{1 - c}}{\rm{E}}_{312}\left(a,1-c+b; 2-c; x,y,z\right).
$

\bigskip

\begin{equation}
{\rm{E}}_{313}\left(a_1,a_2, b; c; x,y,z\right)=\sum\limits_{m,n,p=0}^\infty\frac{(a_1)_{m+2n}(a_2)_{p}(b)_{m-n-p}}{(c)_{m}}\frac{x^m}{m!}\frac{y^n}{n!}\frac{z^p}{p!},
\end{equation}

region of convergence:
$$ \left\{\left\{ s<{1}/{4}, \,\,\,\, r<\min\left\{\Psi_1(s), \Psi_2(s)\right\}\right\}=\left\{r<1,\,\,\,\, s<\min\left\{\Theta_1(r), \Theta_2(r)\right\}\right\},\,\,\,t<\infty
\right\}.
$$

System of partial differential equations:

$
\left\{\begin{aligned}
& {x(1-x)u_{xx} +2y^{2} u_{yy}+2yzu_{yz} -xyu_{xy} +xzu_{xz} }\\
&\,\,\,\,\,\,\,\,\, {+\left[c-(a_{1} +b+1)x\right]u_{x} +(a_{1} -2b+2)yu_{y} +a_{1} zu_{z} -a_{1} bu=0,} \\
&  y(1+4y)u_{yy} +{x^{2} u_{xx}-x(1-4y)u_{xy} +zu_{yz}  } \\
&\,\,\,\,\,\,\,\,\,+2\left(a_{1} +1\right)xu_{x} {+\left[1-b+\left(4a_{1} +6\right)y\right]u_{y} +a_{1} \left(a_{1} +1\right)u=0,} \\
& {zu_{zz} -xu_{xz} +yu_{yz} +(1-b+z)u_{z} +a_{2} u=0,} \end{aligned}\right.
 $

where $u\equiv \,\,{\rm{E}}_{313}\left(a_1,a_2, b; c; x,y,z\right)$.

Particular solutions:

$
{u_1} ={\rm{E}}_{313}\left(a_1,a_2,b; c; x,y,z\right) ,
$

$
{u_2} = {x^{1 - c}}{\rm{E}}_{313}\left(1-c+a_1,a_2,1-c+b; 2-c; x,y,z\right).
$

\bigskip

\begin{equation}
{\rm{E}}_{314}\left(a, b; c; x,y,z\right)=\sum\limits_{m,n,p=0}^\infty\frac{(a)_{m+2n} (b)_{m-n-p}}{(c)_{m}}\frac{x^m}{m!}\frac{y^n}{n!}\frac{z^p}{p!},
\end{equation}

region of convergence:
$$ \left\{\left\{ s<{1}/{4}, \,\,\,\, r<\min\left\{\Psi_1(s), \Psi_2(s)\right\}\right\}=\left\{ r<1, \,\,\,\, s<\min\left\{\Theta_1(r), \Theta_2(r)\right\}\right\},\,\,\,t<\infty
\right\}.
$$

System of partial differential equations:

$
\left\{\begin{aligned}
& {x(1-x)u_{xx} +2y^{2} u_{yy}}-xyu_{xy} +2yzu_{yz}\\
&\,\,\,\,\,\,\,\,\,   +xzu_{xz}{+\left[c-\left(a_{1} +b+1\right)x\right]u_{x} +\left(a_{1} -2b+2\right)yu_{y} +a_{1} zu_{z} -a_{1} bu=0,} \\
&  y(1+4y)u_{yy} +{x^{2} u_{xx}  } -x(1-4y)u_{xy}\\
&\,\,\,\,\,\,\,\,\, +zu_{yz}+2\left(a_{1} +1\right)xu_{x} {+\left[1-b+\left(4a_{1} +6\right)y\right]u_{y} +a_{1} \left(a_{1} +1\right)u=0,} \\
& {zu_{zz} -xu_{xz} +yu_{yz} +(1-b)u_{z} +u=0,} \end{aligned}\right.
$

where $u\equiv \,\,{\rm{E}}_{314}\left(a, b; c; x,y,z\right)$.

Particular solutions:

$
{u_1} ={\rm{E}}_{314}\left(a,b; c; x,y,z\right) ,
$

$
{u_2} = {x^{1 - c}}{\rm{E}}_{314}\left(1-c+a,1-c+b; 2-c; x,y,z\right).
$

\bigskip

\begin{equation}
{\rm{E}}_{315}\left( a, b; c; x,y,z\right)=\sum\limits_{m,n,p=0}^\infty\frac{(a)_{m+2n-p}(b)_{m+p-n}}{(c)_m}\frac{x^m}{m!}\frac{y^n}{n!}\frac{z^p}{p!},
\end{equation}

region of convergence:
$$ \left\{\left\{ s<{1}/{4} ,\,\,\, r<\min\left\{\Psi_1(s), \Psi_2(s)\right\}\right\}= \left\{r<1, \,\,\, s<\min\left\{\Theta_1(r), \Theta_2(r)\right\}\right\},\,\,\,t<\infty
\right\}.
$$

System of partial differential equations:

$ \left\{\begin{aligned}
& {x(1-x)u_{xx} +2y^{2} u_{yy} +z^{2} u_{{zz} } -xyu_{xy} } \\
&\,\,\,\,\,\,\,\,\,-3yzu_{yz} {+\left[c-\left(a +b +1\right)x\right]u_{x} +\left(a -2b +2\right)yu_{y} +\left(b - a+1 \right)zu_{z} - a b u=0,} \\
&  y(1+4y)u_{yy} +{x^{2} u_{xx}+z^{2} u_{{zz} } - x(1-4y)u_{xy}} -(1+4y)zu_{yz}\\& \,\,\,\,\,\,\,\,\,-2xzu_{xz} { +2\left(a +1\right)xu_{x} +\left[1-b +\left(4a +6\right)y\right]u_{y} -2a zu_{z} + a \left(a +1\right)u=0,} \\
& {zu_{zz} -xu_{xz} -2yu_{yz} +xu_{x} -yu_{y} +(1-a +z)u_{z} +bu=0,} \end{aligned}\right.
$

where $u\equiv \,\,{\rm{E}}_{315}\left( a, b; c; x,y,z\right)$.

Particular solutions:

$
{u_1} ={\rm{E}}_{315}\left(a,b; c; x,y,z\right) ,
$

$
{u_2} = {x^{1 - c}}{\rm{E}}_{315}\left(1-c+a,1-c+b; 2-c; x,y,z\right).
$

\bigskip

\begin{equation}
{\rm{E}}_{316}\left(a, b; c; x,y,z\right)=\sum\limits_{m,n,p=0}^\infty\frac{(a)_{m+n+2p}(b)_{m-n-p}}{(c)_{m}}\frac{x^m}{m!}\frac{y^n}{n!}\frac{z^p}{p!},
\end{equation}

region of convergence:
$$ \left\{ \left\{t<{1}/{4}, \,\,\, r<\min\left\{\Psi_1(t), \Psi_2(t)\right\}\right\}=\left\{ r<1, \,\,\, t<\min\left\{\Theta_1(r), \Theta_2(r)\right\}\right\},\,\,\,s<\infty
\right\}.
$$

System of partial differential equations:

$ \left\{\begin{aligned} & {x(1-x)u_{xx} +y^{2} u_{yy}+2z^{2} u_{{zz}} } \\&\,\,\,\,\,\,\,\,\, -xzu_{xz} +3yzu_{yz}+\left[c-\left(a +b +1\right)x\right]u_{x} { +\left(a-b +1\right)yu_{y} +\left(a-2b+2 \right)zu_{z} -ab u=0,} \\
& {yu_{yy} -xu_{xy} +zu_{yz} +xu_{x} +(1-b+y)u_{y} +2zu_{z} +au=0,} \\
& { z(1+4z)u_{zz} +x^{2} u_{xx} +y^{2} u_{yy} +2xyu_{xy}-x(1-4z)u_{xz}   } \\
&\,\,\,\,\,\,\,\,\, +y(1+4z)u_{yz}{+2(a+1)xu_{x} +2(a+1)yu_{y} +\left[1-b+(4a+6)z\right]u_{z} +a(a+1)u=0,} \end{aligned}\right.
$

where $u\equiv \,\,{\rm{E}}_{316}\left(a, b; c; x,y,z\right)$.

Particular solutions:

$
{u_1} ={\rm{E}}_{316}\left(a,b; c; x,y,z\right) ,
$

$
{u_2} = {x^{1 - c}}{\rm{E}}_{316}\left(1-c+a,1-c+b; 2-c; x,y,z\right).
$

\bigskip

\begin{equation}
{\rm{E}}_{317}\left(a, b_1, b_2; c; x,y,z\right)=\sum\limits_{m,n,p=0}^\infty\frac{(a)_{p}(b_1)_{2m-n}(b_2)_{2n-p}}{(c)_{m}}\frac{x^m}{m!}\frac{y^n}{n!}\frac{z^p}{p!},
\end{equation}

region of convergence:
$$ \left\{ 4s(1+2\sqrt{r})<1,\,\,\,t<\infty
\right\}.
$$

System of partial differential equations:

$ \left\{\begin{aligned}
&{x(1-4x)u_{xx} -y^{2} u_{yy} +4xyu_{xy} +\left[c-\left(4b_{1} +6\right)x\right]u_{x}} {+2b_{1}y u_{y} -b_{1} (b_{1} +1)u=0,} \\
& {y(1+4y)u_{yy} +z^{2} u_{_{zz} } -2xu_{xy} -4yzu_{yz}}{+\left[1-b_{1} +\left(4b_{2} +6\right)y\right]u_{y} -2b_{2} zu_{z} +b_{2} \left(b_{2} +1\right)u=0,} \\
& {zu_{zz} -2yu_{yz} +(1-b_{2} +z)u_{z} +au=0,} \end{aligned}\right.
$

where $u\equiv \,\,{\rm{E}}_{317}\left(a, b_1, b_2; c; x,y,z\right)$.

Particular solutions:

$
{u_1} ={\rm{E}}_{317}\left(a, b_1, b_2; c; x,y,z\right) ,
$

$
{u_2} = {x^{1 - c}}{\rm{E}}_{317}\left(a, 2-2c+b_1, b_2; 2-c; x,y,z\right).
$

\bigskip

\begin{equation}
{\rm{E}}_{318}\left( a, b; c; x,y,z\right)=\sum\limits_{m,n,p=0}^\infty\frac{(a)_{2m-n}(b)_{2n-p}}{(c)_{m}}\frac{x^m}{m!}\frac{y^n}{n!}\frac{z^p}{p!},
\end{equation}

region of convergence:
$$ \left\{ 4s(1+2\sqrt{r})<1,\,\,\,t<\infty
\right\}.
$$

System of partial differential equations:

$  \left\{\begin{aligned}
&{x(1-4x)u_{xx} -y^{2} u_{yy} +4xyu_{xy} +\left[c-\left(4a +6\right)x\right]u_{x}} {+2ay u_{y} -a (a +1)u=0,} \\
& {y(1+4y)u_{yy} +z^{2} u_{_{zz} } - 2xu_{xy} -4yzu_{yz}} {+\left[1- a +\left(4b +6\right)y\right]u_{y} -2b zu_{z} +b \left(b +1\right)u=0,} \\
& {zu_{zz} -2yu_{yz} +(1-b)u_{z} +u=0,} \end{aligned}\right.
$

where $u\equiv \,\,{\rm{E}}_{318}\left( a, b; c; x,y,z\right)$.

Particular solutions:

$
{u_1} ={\rm{E}}_{318}\left(a,b; c; x,y,z\right) ,
$

$
{u_2} = {x^{1 - c}}{\rm{E}}_{318}\left(2-2c+a,b; 2-c; x,y,z\right).
$

\bigskip

\begin{equation}
{\rm{E}}_{319}\left(a, b; c; x,y,z\right)=\sum\limits_{m,n,p=0}^\infty\frac{(a)_{m+2n-p}(b)_{2p-n}}{(c)_{m}}\frac{x^m}{m!}\frac{y^n}{n!}\frac{z^p}{p!},
\end{equation}

region of convergence:
$$ \left\{ r<\infty,\,\,\,Z_1\cap Z_2
\right\},
$$

$$
Z_1=\left\{s<\Phi_1(t)\right\}=\left\{s<\frac{1}{4},\,\,\, t<\Phi_2(s)\right\}
$$

$$
Z_2=\left\{t<\Phi_1(s)\right\}=\left\{t<\frac{1}{4}, \,\,\,\, s<\Phi_2(t)\right\}.
$$

System of partial differential equations:

$\left\{\begin{aligned}
 &{xu_{xx} +(c-x)u_{x} -2yu_{y} +zu_{z} - a u=0,} \\
 & {y(1+4y)u_{yy} +x^2u_{xx}+z^{2} u_{_{zz} } +4xyu_{xy} }-2xzu_{xz} \\& \,\,\,\,\,\,\,\,\,-2(1+2y)zu_{yz}{+2\left(a+1\right)xu_x +\left[1-b +\left(4a+6\right)y\right]u_{y} -2a zu_{z} +a (a +1)u=0,} \\
 & { z(1+4z)u_{zz}+y^{2} u_{yy} -xu_{xz} -2y(1+2z)u_{yz}  }-2b yu_{y}{+\left[1-a +\left(4b +6\right)z\right]u_{z} +b(b+1)u=0,} \end{aligned}\right.
 $

where $u\equiv \,\,{\rm{E}}_{319}\left(a, b; c; x,y,z\right)$.

Particular solutions:

$
{u_1} ={\rm{E}}_{319}\left(a,b; c; x,y,z\right) ,
$

$
{u_2} = {x^{1 - c}}{\rm{E}}_{319}\left(1-c+a,b; 2-c; x,y,z\right).
$

\bigskip

\begin{equation}
{\rm{E}}_{320}\left(a, b; c; x,y,z\right)=\sum\limits_{m,n,p=0}^\infty\frac{(a)_{2m+n-p}(b)_{2p-n}}{(c)_{m}}\frac{x^m}{m!}\frac{y^n}{n!}\frac{z^p}{p!},
\end{equation}

region of convergence:
$$ \left\{ 4t(1+2\sqrt{r})<1,\,\,\,s<\infty
\right\}.
$$

System of partial differential equations:

$
  \left\{\begin{aligned}
 &{x(1-4x)u_{xx} -y^{2} u_{yy}-z^{2} u_{zz}  }  -4xyu_{xy} +4xzu_{xz}\\& \,\,\,\,\,\,\,\,\, +2yzu_{yz}{+\left[c-\left(4a +6\right)x\right]u_{x} -2\left(a +1\right)yu_{y} +2azu_{z} -a \left(a +1\right)u=0,} \\
 & {yu_{yy} -2zu_{yz} +2xu_{x} +\left(1-b +y\right)u_{y} - zu_{z} +a u=0,} \\
 & { z(1+4z)u_{zz} +y^{2} u_{yy}-2xu_{xz} -y(1+4z)u_{yz}  }-2b yu_{y} {+\left[1-a +\left(4b+6\right)z\right])u_{z} +b (b +1)u=0,} \end{aligned}\right.
$

where $u\equiv \,\,{\rm{E}}_{320}\left(a, b; c; x,y,z\right)$.

Particular solutions:

$
{u_1} ={\rm{E}}_{320}\left(a,b; c; x,y,z\right) ,
$

$
{u_2} = {x^{1 - c}}{\rm{E}}_{320}\left(2-2c+a,b; 2-c; x,y,z\right).
$

\bigskip

\begin{equation}
{\rm{E}}_{321}\left(a,b; c; x,y,z\right)=\sum\limits_{m,n,p=0}^\infty\frac{(a)_{n+2p}(b)_{2m-n-p}}{(c)_{m}}\frac{x^m}{m!}\frac{y^n}{n!}\frac{z^p}{p!},
\end{equation}

region of convergence:
$$ \left\{ 4t(1+2\sqrt{r})<1,\,\,\,s<\infty
\right\}.
$$

System of partial differential equations:

$
 \left\{\begin{aligned}
&{x(1-4x)u_{xx} -y^{2} u_{yy}   }-z^2u_{zz}+4xyu_{xy}\\
&\,\,\,\,\,\,\,\,\, +4xzu_{xz} -2yzu_{yz}{+\left[c-\left(4b +6\right)x\right]u_{x} +2byu_{y} +2bzu_{z} -b(b+1)u=0,} \\
& {yu_{yy} -2xu_{xy} +zu_{yz} +(1-b+y)u_{y} +2zu_{z} +au=0,} \\
& { z(1+4z)u_{zz} +y^{2} u_{yy}-2xu_{xz} +y(4z+1)u_{yz}  } \\& \,\,\,\,\,\,\,\,\,+2(a+1)yu_{y} {+\left[1-b+(4a+6)z\right]u_{z} +a(a+1)u=0,} \end{aligned}\right.
$

where $u\equiv \,\,{\rm{E}}_{321}\left(a,b; c; x,y,z\right)$.

Particular solutions:

$
{u_1} ={\rm{E}}_{321}\left(a,b; c; x,y,z\right) ,
$

$
{u_2} = {x^{1 - c}}{\rm{E}}_{321}\left(a,2-2c+b; 2-c; x,y,z\right).
$

\bigskip

\begin{equation}
{\rm{E}}_{322}\left(a, b; c; x,y,z\right)=\sum\limits_{m,n,p=0}^\infty\frac{(a)_{2m+2n-p}(b)_{p-n}}{(c)_{m}}\frac{x^m}{m!}\frac{y^n}{n!}\frac{z^p}{p!},
\end{equation}

region of convergence:
$$ \left\{ \sqrt{r}+\sqrt{s}<\frac{1}{2},\,\,\,t<\infty
\right\}.
$$

System of partial differential equations:

$
  \left\{\begin{aligned}
 & {x(1-4x)u_{xx}  -4y^{2} u_{yy}  }-z^2u_{zz}-8xyu_{xy} \\
 &\,\,\,\,\,\,\,\,\, +4xzu_{xz}+4yzu_{yz} {+\left[c-\left(4a +6\right)x\right]u_{x} -(4a+6)yu_{y} +2azu_{z} -a(a+1)u=0,} \\
 &  y(1+4y)u_{yy}+{4x^{2} u_{xx} +z^{2} u_{zz}}+8xyu_{xy} -4xzu_{xz}  \\
 &\,\,\,\,\,\,\,\,\, -(1+4y)zu_{xy}{+(4a+6)xu_{x} +\left[1-b +(4a+6)y\right]u_{y} -2azu_{z} +a(a+1)u=0,} \\
 & {zu_{zz}-2xu_{xz}-2yu_{yz} -yu_y+(1-a +z)u_{z} + b u=0,} \end{aligned}\right.
$

where $u\equiv \,\,{\rm{E}}_{322}\left(a, b; c; x,y,z\right)$.

Particular solutions:

$
{u_1} ={\rm{E}}_{322}\left(a,b; c; x,y,z\right) ,
$

$
{u_2} = {x^{1 - c}}{\rm{E}}_{322}\left(2-2c+a,b; 2-c; x,y,z\right).
$

\bigskip

\begin{equation}
{\rm{E}}_{323}\left(a_1,a_2,b_1,b_2, b_3;x,y,z\right)=\sum\limits_{m,n,p=0}^\infty{(a_1)_m(a_2)_n(b_1)_{m-n}(b_2)_{n-p}(b_3)_{p-m}}\frac{x^m}{m!}\frac{y^n}{n!}\frac{z^p}{p!},
\end{equation}

region of convergence:
$$ \left\{ s(1+r)<1,\,\,\,t<\infty
\right\}.
$$

System of partial differential equations:

$
 \left\{\begin{aligned}
 &{x(1+x)u_{xx} -xyu_{xy} -zu_{xz} } {+\left[1-b_{3} +\left(a_{1} +b_{1} +1\right)x\right]u_{x} -a_{1} yu_{y} +a_{1} b_{1} u=0,} \\
 & {y(1+y)u_{yy} -xu_{xy} -yzu_{yz} }{+\left[1-b_{1} +\left(a_{2} +b_{2} +1\right)y\right]u_{y} -a_{2} zu_{z} +a_{2} b_{2} u=0,} \\
 & {zu_{zz} -yu_{yz} -xu_{x} +(1-b_{2} +z)u_{z} +b_{3} u=0,} \end{aligned}\right.
 $

where $u\equiv \,\,{\rm{E}}_{323}\left(a_1,a_2,b_1,b_2, b_3;x,y,z\right)$.

\bigskip

\begin{equation}
{\rm{E}}_{324}\left(a,b_1,b_2, b_3;x,y,z\right)=\sum\limits_{m,n,p=0}^\infty{(a)_m(b_1)_{m-n}(b_2)_{n-p}(b_3)_{p-m}}\frac{x^m}{m!}\frac{y^n}{n!}\frac{z^p}{p!},
\end{equation}

region of convergence:
$$ \left\{ r<1,\,\,\,s<\infty,\,\,\,t<\infty
\right\}.
$$

System of partial differential equations:

$\left\{\begin{aligned}
&{x(1+x)u_{xx} -xyu_{xy} -zu_{xz} }{+\left[1-b_{3} +\left(a +b_{1} +1\right)x\right]u_{x} -a yu_{y} +a b_{1} u=0,} \\
& {yu_{yy} -xu_{xy} +\left(1-b_{1} +y\right)u_{y} -zu_{z} +b_{2} u=0,} \\
& {zu_{zz} -yu_{yz} -xu_{x} +(1-b_{2} +z)u_{z} +b_{3} u=0,} \end{aligned}\right.
$

where $u\equiv \,\,{\rm{E}}_{324}\left(a,b_1,b_2, b_3;x,y,z\right)$.

\bigskip

\begin{equation}
{\rm{E}}_{325}\left(b_1,b_2, b_3;x,y,z\right)=\sum\limits_{m,n,p=0}^\infty{(b_1)_{m-n}(b_2)_{n-p}(b_3)_{p-m}}\frac{x^m}{m!}\frac{y^n}{n!}\frac{z^p}{p!},
\end{equation}

region of convergence:
$$ \left\{ r<\infty,\,\,\,s<\infty,\,\,\,t<\infty
\right\}.
$$

System of partial differential equations:

$\left\{\begin{aligned}
& {xu_{xx} -zu_{xz} +\left(1-b_{3} +x\right)u_{x} -yu_{y} +b_{1} u=0,} \\
& {yu_{yy} -xu_{xy} +\left(1-b_{1} +y\right)u_{y} -zu_{z} +b_{2} u=0,} \\
& {zu_{zz} -yu_{yz} -xu_{x} +(1-b_{2} +z)u_{z} +b_{3} u=0,} \end{aligned}\right.
$

where $u\equiv \,\,{\rm{E}}_{325}\left(b_1,b_2, b_3;x,y,z\right)$.

\bigskip

\begin{equation}
{\rm{E}}_{326}\left(a_1,a_2,b_1,b_2, b_3;x,y,z\right)=\sum\limits_{m,n,p=0}^\infty{(a_1)_m(a_2)_p(b_1)_{m-n}(b_2)_{n-p}(b_3)_{n-m}}\frac{x^m}{m!}\frac{y^n}{n!}\frac{z^p}{p!},
\end{equation}

region of convergence:
$$ \left\{ r<1,\,\,\, s<1,\,\,\,t<\infty
\right\}.
$$

System of partial differential equations:

$
  \left\{\begin{aligned}
 & {x(1+x)u_{xx} -(x+1)yu_{xy}}{+\left[1-b_{3} +\left(a_{1} +b_{1} +1\right)x\right]u_{x} -a_{1} yu_{y} +a_{1} b_{1} u=0,} \\
 & {y(1+y)u_{yy} -x(y+1)u_{xy} }+xzu_{xz} -yzu_{yz}\\& \,\,\,\,\,\,\,\,\, -b_{2} xu_{x}{+\left[1-b_{1} +\left(b_{2} +b_{3} +1\right)y \right]u_{y} -b_{3} zu_{z} +b_{2} b_{3} u=0,} \\
 & {zu_{zz} -yu_{yz} +(1-b_{2} +z)u_{z} +a_{2} u=0,} \end{aligned}\right.
 $

where $u\equiv \,\,{\rm{E}}_{326}\left(a_1,a_2,b_1,b_2, b_3;x,y,z\right)$.

\bigskip

\begin{equation}
{\rm{E}}_{327}\left(a_1,a_2,b_1,b_2, b_3;x,y,z\right)=\sum\limits_{m,n,p=0}^\infty{(a_1)_p(a_2)_p(b_1)_{m-n}(b_2)_{n-p}(b_3)_{n-m}}\frac{x^m}{m!}\frac{y^n}{n!}\frac{z^p}{p!},
\end{equation}

region of convergence:
$$ \left\{ r<\infty,\,\,\,t(1+s)<1
\right\}.
$$

System of partial differential equations:

$  \left\{\begin{aligned}&{xu_{xx} -yu_{xy} +(1-b_{3} +x)u_{x} -yu_{y} +b_{1} u=0,} \\& y(1+y)u_{yy} -x(1+y)u_{xy} -yzu_{yz}\\& \,\,\,\,\,\,\,\,\, +xzu_{xz}   -b_{2} xu_{x}{+\left[1-b_{1} +\left(b_{2} +b_{3} +1\right)y\right]u_{y} -b_{3} zu_{z} +b_{2} b_{3} u=0,} \\&
 {z(1+z)u_{zz} -yu_{yz} +\left[1-b_{2} +\left(a_{1} +a_{2} +1\right)z\right]u_{z} +a_{1} a_{2} u=0,} \end{aligned}\right.
$

where $u\equiv \,\,{\rm{E}}_{327}\left(a_1,a_2,b_1,b_2, b_3;x,y,z\right)$.

\bigskip

\begin{equation}
{\rm{E}}_{328}\left(a,b_1,b_2, b_3;x,y,z\right)=\sum\limits_{m,n,p=0}^\infty{(a)_m(b_1)_{m-n}(b_2)_{n-p}(b_3)_{n-m}}\frac{x^m}{m!}\frac{y^n}{n!}\frac{z^p}{p!},
\end{equation}

region of convergence:
$$ \left\{ r<1 \vee s<1,\,\,\,t<\infty
\right\}.
$$

System of partial differential equations:

$  \left\{\begin{aligned}
&{x(1+x)u_{xx} -(x+1)yu_{xy} } {+\left[1-b_{3} +\left(a +b_{1} +1\right)x\right]u_{x} -a yu_{y} +a b_{1} u=0,} \\
&{y(1+y)u_{yy} -x(y+1)u_{xy} -yzu_{yz} +xzu_{xz}  }\\& \,\,\,\,\,\,\,\,\, -b_{2} xu_{x}{+\left[1-b_{1} +\left(b_{2} +b_{3} +1\right)y \right]u_{y} -b_{3} zu_{z} +b_{2} b_{3} u=0,} \\
& {zu_{zz} -yu_{yz} +(1-b_{2} )u_{z} +u=0,} \end{aligned}\right.
$

where $u\equiv \,\,{\rm{E}}_{328}\left(a,b_1,b_2, b_3;x,y,z\right)$.

\bigskip

\begin{equation}
{\rm{E}}_{329}\left(a,b_1,b_2, b_3;x,y,z\right)=\sum\limits_{m,n,p=0}^\infty{(a)_p(b_1)_{m-n}(b_2)_{n-p}(b_3)_{n-m}}\frac{x^m}{m!}\frac{y^n}{n!}\frac{z^p}{p!},
\end{equation}

region of convergence:
$$ \left\{ r<\infty,\,\,\,s<1,\,\,\,t<\infty
\right\}.
$$

System of partial differential equations:

$
 \left\{\begin{aligned}&{xu_{xx} -yu_{xy} +(1-b_{3} +x)u_{x} -yu_{y} +b_{1} u=0,} \\& {y(1+y)u_{yy} -x(y+1)u_{xy} -yzu_{yz} +xzu_{xz} }\\& \,\,\,\,\,\,\,\,\,-b_{2} xu_{x} {+\left[1-b_{1} +\left(b_{2} +b_{3} +1\right)y\right]u_{y} -b_{3} zu_{z} +b_{2} b_{3} u=0,} \\& {zu_{zz} -yu_{yz} +\left(1-b_{2} +z\right)u_{z} +a u=0,} \end{aligned}\right.
 $

where $u\equiv \,\,{\rm{E}}_{329}\left(a,b_1,b_2, b_3;x,y,z\right)$.

\bigskip

\begin{equation}
{\rm{E}}_{330}\left(a,b, c;x,y,z\right)=\sum\limits_{m,n,p=0}^\infty{(a)_{m-n}(b)_{n-p}(c)_{n-m}}\frac{x^m}{m!}\frac{y^n}{n!}\frac{z^p}{p!},
\end{equation}

region of convergence:
$$ \left\{ r<\infty,\,\,\,s<1,\,\,\,t<\infty
\right\}.
$$

System of partial differential equations:

$
 \left\{\begin{aligned}& {xu_{xx} -yu_{xy} +(1-c +x)u_{x} -yu_{y} +a u=0,} \\& y(1+y)u_{yy} -x(y+1)u_{xy} -yzu_{yz} \\& \,\,\,\,\,\,\,\,\,+xzu_{xz}  -b xu_{x}{+\left[1-a+\left(b +c+1\right)y\right]u_{y} -c zu_{z} +b c u=0,} \\& {zu_{zz} -yu_{yz} +(1-b)u_{z} +u=0,} \end{aligned}\right.
$

where $u\equiv \,\,{\rm{E}}_{330}\left(a,b, c;x,y,z\right)$.

\bigskip

\begin{equation}
{\rm{E}}_{331}\left(a_1,a_2,a_3,b_1,b_2;x,y,z\right)=\sum\limits_{m,n,p=0}^\infty{(a_1)_m(a_2)_n(a_3)_n(b_1)_{p-m-n}(b_2)_{m-p}}\frac{x^m}{m!}\frac{y^n}{n!}\frac{z^p}{p!},
\end{equation}

region of convergence:
$$ \left\{ \frac{1}{r}+\frac{1}{s}>1,\,\,\,t<\infty
\right\}.
$$

System of partial differential equations:

$
 \left\{\begin{aligned}
& {x(1+x)u_{xx} +yu_{xy} -(1+x)zu_{xz} }{+\left[1-b_{1} +\left(a_{1} +b_{2} +1\right)x\right]u_{x} -a_{1} zu_{z} +a_{1} b_{2} u=0,} \\
& {y(1+y)u_{yy} +xu_{xy} -zu_{yz} +\left[1-b_{1} +\left(a_{2} +a_{3} +1\right)y\right]u_{y} +a_{2} a_{3} u=0,} \\
& {zu_{zz} -xu_{xz} -xu_{x} -yu_{y} +\left(1-b_{2} +z\right)u_{z} +b_{1} u=0,} \end{aligned}\right. $

where $u\equiv \,\,{\rm{E}}_{331}\left(a_1,a_2,a_3,b_1,b_2;x,y,z\right)$.

\bigskip

\begin{equation}
{\rm{E}}_{332}\left(a_1,a_2,a_3,b_1,b_2;x,y,z\right)=\sum\limits_{m,n,p=0}^\infty{(a_1)_n(a_2)_n(a_3)_p(b_1)_{p-m-n}(b_2)_{m-p}}\frac{x^m}{m!}\frac{y^n}{n!}\frac{z^p}{p!},
\end{equation}

region of convergence:
$$ \left\{ t(1+r)<1,\,\,\,s<\infty
\right\},
$$

System of partial differential equations:

$
 \left\{\begin{aligned}
&{xu_{xx} +yu_{xy} -zu_{xz} +(1-b_{1} +x)u_{x} -zu_{z} +b_{2} u=0,} \\
& {y(1+y)u_{yy} +xu_{xy} -zu_{yz} +\left[1-b_{1} + \left(a_{1} +a_{2} +1\right)y \right]u_{y} +a_{1} a_{2} u=0,} \\
& {z(1+z)u_{zz} -x(1+z)u_{xz} -yzu_{yz}-a_{3} xu_{x}  } {-a_{3} yu_{y}+\left[1-b_{2} +\left(a_{3} +b_{1} +1\right)z\right]u_{z}  +a_{3} b_{1} u=0,} \end{aligned}\right.
$

where $u\equiv \,\,{\rm{E}}_{332}\left(a_1,a_2,a_3,b_1,b_2;x,y,z\right)$.

\bigskip

\begin{equation}
{\rm{E}}_{333}\left(a_1,a_2,a_3,b_1,b_2;x,y,z\right)=\sum\limits_{m,n,p=0}^\infty{(a_1)_m(a_2)_n(a_3)_p(b_1)_{p-m-n}(b_2)_{m-p}}\frac{x^m}{m!}\frac{y^n}{n!}\frac{z^p}{p!},
\end{equation}

region of convergence:
$$ \left\{ r<1 \vee t<1,\,\,\,s<\infty
\right\}.
$$

System of partial differential equations:

$
 \left\{\begin{aligned}& {x(1+x)u_{xx} +yu_{xy} -(1+x)zu_{xz} } {+\left[1-b_{1} +\left(a_{1} +b_{2} +1\right)x\right]u_{x} -a_{1} zu_{z} +a_{1} b_{2} u=0,} \\
 & {yu_{yy} +xu_{xy} -zu_{yz} +(1-b_{1} +y)u_{y} +a_{2} u=0,} \\
 & {z(1+z)u_{zz} -x(1+z)u_{xz} -yzu_{yz}-a_{3} xu_{x}  } {-a_{3} yu_{y}+\left[1-b_{2} +\left(a_{3} +b_{1} +1\right)z\right]u_{z}  +a_{3} b_{1} u=0,} \end{aligned}\right.
  $

where $u\equiv \,\,{\rm{E}}_{333}\left(a_1,a_2,a_3,b_1,b_2;x,y,z\right)$.

\bigskip

\begin{equation}
{\rm{E}}_{334}\left(a_1,a_2,b_1,b_2;x,y,z\right)=\sum\limits_{m,n,p=0}^\infty{(a_1)_m(a_2)_n(b_1)_{p-m-n}(b_2)_{m-p}}\frac{x^m}{m!}\frac{y^n}{n!}\frac{z^p}{p!},
\end{equation}

region of convergence:
$$ \left\{ r<1 ,\,\,\,s<\infty,\,\,\,t<\infty
\right\}.
$$

System of partial differential equations:

$
 \left\{\begin{aligned} & {x(1+x)u_{xx} +yu_{xy} -(1+x)zu_{xz} } {+\left[1-b_{1} +\left(a_{1} +b_{2} +1\right)x\right]u_{x} -a_{1} zu_{z} +a_{1} b_{2} u=0,} \\
& {yu_{yy} +xu_{xy} -zu_{yz} +(1-b_{1} +y)u_{y} +a_{2} u=0,} \\
& {zu_{zz} -xu_{xz} -xu_{x} -yu_{y} +\left(1-b_{2} +z\right)u_{z} +b_{1} u=0,} \end{aligned}\right.
$

where $u\equiv \,\,{\rm{E}}_{334}\left(a_1,a_2,b_1,b_2;x,y,z\right)$.

\bigskip

\begin{equation}
{\rm{E}}_{335}\left(a_1,a_2,b_1,b_2;x,y,z\right)=\sum\limits_{m,n,p=0}^\infty{(a_1)_n(a_2)_n(b_1)_{p-m-n}(b_2)_{m-p}}\frac{x^m}{m!}\frac{y^n}{n!}\frac{z^p}{p!},
\end{equation}

region of convergence:
$$ \left\{ r<\infty ,\,\,\,s<1,\,\,\,t<\infty
\right\}.
$$

System of partial differential equations:

$
 \left\{\begin{aligned}
 &{xu_{xx} +yu_{xy} -zu_{xz} +(1-b_{1} +x)u_{x} -zu_{z} +b_{2} u=0,} \\
 & {y(1+y)u_{yy} +xu_{xy} -zu_{yz} +\left[1-b_{1} + \left(a_{1} +a_{2} +1\right)y \right]u_{y}  +a_{1} a_{2} u=0,} \\
 & {zu_{zz} -xu_{xz} -xu_{x} -yu_{y} +\left(1-b_{2} +z\right)u_{z} +b_{1} u=0,} \end{aligned}\right.
 $

where $u\equiv \,\,{\rm{E}}_{335}\left(a_1,a_2,b_1,b_2;x,y,z\right)$.

\bigskip

\begin{equation}
{\rm{E}}_{336}\left(a_1,a_2,b_1,b_2;x,y,z\right)=\sum\limits_{m,n,p=0}^\infty{(a_1)_n(a_2)_p(b_1)_{p-m-n}(b_2)_{m-p}}\frac{x^m}{m!}\frac{y^n}{n!}\frac{z^p}{p!},
\end{equation}

region of convergence:
$$ \left\{ r<\infty ,\,\,\,s<\infty,\,\,\,t<1
\right\}.
$$

System of partial differential equations:

$
\left\{\begin{aligned}
 &{xu_{xx} +yu_{xy} -zu_{xz} +(1-b_{1} +x)u_{x} -zu_{z} +b_{2} u=0,} \\
 & {yu_{yy} +xu_{xy} -zu_{yz} +(1-b_{1} +y)u_{y} +a_{1} u=0,} \\
 & {z(1+z)u_{zz} -x(1+z)u_{xz} -yzu_{yz}-a_{2} xu_{x}  } {-a_{2} yu_{y}+\left[1-b_{2} +\left(a_{2} +b_{1} +1\right)z\right]u_{z}  +a_{2} b_{1} u=0,} \end{aligned}\right.
 $

where $u\equiv \,\,{\rm{E}}_{336}\left(a_1,a_2,b_1,b_2;x,y,z\right)$.

\bigskip

\begin{equation}
{\rm{E}}_{337}\left(a,b, c;x,y,z\right)=\sum\limits_{m,n,p=0}^\infty{(a)_m(b)_{p-m-n}(c)_{m-p}}\frac{x^m}{m!}\frac{y^n}{n!}\frac{z^p}{p!},
\end{equation}

region of convergence:
$$ \left\{ r<1 ,\,\,\,s<\infty,\,\,\,t<\infty
\right\}.
$$

System of partial differential equations:

$ \left\{\begin{aligned}
 & {x(1+x)u_{xx} +yu_{xy} -(1+x)zu_{xz} } {+\left[1-b +\left(a +c +1\right)x\right]u_{x} -azu_{z} +a cu=0,} \\
 & {yu_{yy} +xu_{xy} -zu_{yz} +(1-b )u_{y} +u=0,} \\
 & {zu_{zz} -xu_{xz} -xu_{x} -yu_{y} +\left(1- c +z\right)u_{z} +b u=0,} \end{aligned}\right.
 $

where $u\equiv \,\,{\rm{E}}_{337}\left(a,b, c;x,y,z\right)$.

\bigskip

\begin{equation}
{\rm{E}}_{338}\left(a,b,c;x,y,z\right)=\sum\limits_{m,n,p=0}^\infty\frac{(a)_n(b)_{p-m-n}(c)_{m-p}}{m!n!p!}x^my^nz^p,
\end{equation}

region of convergence:
$$ \left\{ r<\infty,\,\,\,s<\infty,\,\,\,t<\infty
\right\},
$$

System of partial differential equations:

$ \left\{\begin{aligned} &{xu_{xx} +yu_{xy} -zu_{xz} +(1-b +x)u_{x} -zu_{z} +c u=0,} \\& {yu_{yy} +xu_{xy} -zu_{yz} +(1-b +y)u_{y} +au=0,} \\& {zu_{zz} -xu_{xz} -xu_{x} -yu_{y} +(1- c +z)u_{z} +b u=0,} \end{aligned}\right.
 $

where $u\equiv \,\,{\rm{E}}_{338}\left(a,b,c;x,y,z\right)$.

\bigskip

\begin{equation}
{\rm{E}}_{339}\left(a,b,c;x,y,z\right)=\sum\limits_{m,n,p=0}^\infty{(a)_p(b)_{p-m-n}(c)_{m-p}}\frac{x^m}{m!}\frac{y^n}{n!}\frac{z^p}{p!},
\end{equation}

region of convergence:
$$ \left\{ r<\infty,\,\,\,s<\infty,\,\,\,t<1
\right\},
$$

System of partial differential equations:

$ \left\{\begin{aligned}
& {xu_{xx} +yu_{xy} -zu_{xz} +(1-b +x)u_{x} -zu_{z} + c u=0,} \\
& {yu_{yy} +xu_{xy} -zu_{yz} +(1-b )u_{y} +u=0,} \\
 & {z(1+z)u_{zz} -x(1+z)u_{xz} -yzu_{yz}-a xu_{x}  } {-a yu_{y}+\left[1-c +\left(a +b +1\right)z\right]u_{z}  +a bu=0,} \end{aligned}\right.
$

where $u\equiv \,\,{\rm{E}}_{339}\left(a,b,c;x,y,z\right)$.

\bigskip

\begin{equation}
{\rm{E}}_{340}\left(a,b;x,y,z\right)=\sum\limits_{m,n,p=0}^\infty{(a)_{p-m-n}(b)_{m-p}}\frac{x^m}{m!}\frac{y^n}{n!}\frac{z^p}{p!},
\end{equation}

region of convergence:
$$ \left\{ r<\infty,\,\,\,s<\infty,\,\,\,t<\infty
\right\}.
$$

System of partial differential equations:

$
 \left\{\begin{aligned}
&{xu_{xx} +yu_{xy} -zu_{xz} +(1-a +x)u_{x} -zu_{z} +b u=0,} \\
& {yu_{yy} +xu_{xy} -zu_{yz} +(1-a )u_{y} +u=0,} \\
& {zu_{zz} -xu_{xz} -xu_{x} -yu_{y} +(1-b+z)u_{z} +a u=0,} \end{aligned}\right.
$

where $u\equiv \,\,{\rm{E}}_{340}\left(a,b;x,y,z\right)$.

\bigskip

\begin{equation}
{\rm{E}}_{341}\left(a,b_1,b_2,b_3; x,y,z\right)=\sum\limits_{m,n,p=0}^\infty{(a)_{p}(b_1)_{m+n-p}(b_2)_{m-n}(b_3)_{n-m}}\frac{x^m}{m!}\frac{y^n}{n!}\frac{z^p}{p!},
\end{equation}

region of convergence:
$$ \left\{ r+s<1,\,\,\,t<\infty
\right\}.
$$

System of partial differential equations:

$
  \left\{\begin{aligned} &{x(1+x)u_{xx}-y^{2} u_{yy}}-yu_{xy} -xzu_{xz} +yzu_{yz}\\
 &\,\,\,\,\,\,\,\,\,{+\left[1-b_{3} +\left(b_{1} +b_{2} +1\right)x\right]u_{x} +\left(b_{2} -b_{1} -1\right)yu_y-b_{2} zu_{z} +b_{1} b_{2} u=0,} \\
 & {y(1+y)u_{yy} -x^{2} u_{xx}  } -xu_{xy} -yzu_{yz} +xzu_{xz}\\
 &\,\,\,\,\,\,\,\,\,+(b_{3} -b_{1} -1)xu_{x} {+(1-b_{2} +y(b_{1} +b_{3} +1))u_{y} -b_{3} zu_{z} +b_{1} b_{3} u=0,} \\
 & {zu_{zz} -xu_{xz} -yu_{yz} +(1-b_{1} +z)u_{z} +au=0,} \end{aligned}\right.
  $

where $u\equiv \,\,{\rm{E}}_{341}\left(a,b_1,b_2,b_3; x,y,z\right)$.

\bigskip

\begin{equation}
{\rm{E}}_{342}\left(a,b,c; x,y,z\right)=\sum\limits_{m,n,p=0}^\infty{(a)_{m+n-p}(b)_{m-n}(c)_{n-m}}\frac{x^m}{m!}\frac{y^n}{n!}\frac{z^p}{p!},
\end{equation}

region of convergence:
$$ \left\{ r+s<1,\,\,\,t<\infty
\right\}.
$$

System of partial differential equations:

$
\left\{\begin{aligned}
&{x(1+x)u_{xx}-y^{2} u_{yy} }-yu_{xy} -xzu_{xz}\\
 &\,\,\,\,\,\,\,\,\,  +yzu_{yz}{+\left[1-c +\left(a +b +1\right)x\right]u_{x} +\left(b -a -1\right)yu_y-b zu_{z} + a b u=0,} \\
& {y(1+y)u_{yy} -x^{2} u_{xx} }-xu_{xy}+xzu_{xz} \\
&\,\,\,\,\,\,\,\,\,   -yzu_{yz}-(a -c +1)xu_{x} {+\left[1-b +\left(a +c +1\right)y\right]u_{y} - c zu_{z} + a c u=0,} \\& {zu_{zz} -xu_{xz} -yu_{yz} +(1- a )u_{z} +u=0,} \end{aligned}\right.
 $

where $u\equiv \,\,{\rm{E}}_{342}\left(a,b,c; x,y,z\right)$.

\bigskip

\begin{equation}
{\rm{E}}_{343}\left(a,b_1,b_2,b_3; x,y,z\right)=\sum\limits_{m,n,p=0}^\infty{(a)_{n}(b_1)_{m+n-p}(b_2)_{m-n}(b_3)_{p-m}}\frac{x^m}{m!}\frac{y^n}{n!}\frac{z^p}{p!},
\end{equation}

region of convergence:
$$ \left\{ \left\{s+2\sqrt{rs}<1\right\}=\left\{\sqrt{s}<\sqrt{1+r}-\sqrt{r}\right\},\,\,\,t<\infty
\right\}.
$$

System of partial differential equations:

$
\left\{\begin{aligned}
&{x(1+x)u_{xx} -y^{2} u_{yy}}-(x+1)zu_{xz} \\
&\,\,\,\,\,\,\,\,\,   +yzu_{yz} {+\left[1-b_{3} +\left(b_{1} +b_{2} +1\right)x\right]u_{x} +(b_{2} -b_{1} -1)yu_y-b_{2} zu_{z} +b_{1} b_{2} u=0,} \\
& {y(1+y)u_{yy} -x(1-y)u_{xy} -yzu_{yz} +axu_{x} } {+\left[1-b_{2} +\left(a+b_{1} +1\right)y\right]u_{y} -azu_{z} +ab_{1} u=0,} \\
& {zu_{zz} -xu_{xz} -yu_{yz} -xu_{x} +(1-b_{1} +z)u_{z} +b_{3} u=0,} \end{aligned}\right.
$

where $u\equiv \,\,{\rm{E}}_{343}\left(a,b_1,b_2,b_3; x,y,z\right)$.

\bigskip

\begin{equation}
{\rm{E}}_{344}\left(a,b_1,b_2,b_3; x,y,z\right)=\sum\limits_{m,n,p=0}^\infty{(a)_p(b_1)_{m+n-p}(b_2)_{m-n}(b_3)_{p-m}}\frac{x^m}{m!}\frac{y^n}{n!}\frac{z^p}{p!},
\end{equation}

region of convergence:
$$ \left\{ r<1 ,\,\,\, t<1,\,\,\,s<\infty
\right\}.
$$

System of partial differential equations:

$
  \left\{\begin{aligned}
 &{x(1+x)u_{xx}}-y^{2} u_{yy}-(x+1)zu_{xz}  \\
&\,\,\,\,\,\,\,\,\,   +yzu_{yz} {+\left[1-b_{3} +\left(b_{1} +b_{2} +1\right)x\right]u_{x} +(b_{2} -b_{1} -1)yu_y-b_{2} zu_{z} +b_{1} b_{2} u=0,} \\
 & {yu_{yy} -xu_{xy} +xu_{x} +(1-b_{2} +y)u_{y} -zu_{z} +b_{1} u=0,} \\
 & {z(1+z)u_{zz} -x(1+z)u_{xz} -yu_{yz}-a xu_{x}  } {+\left[1-b_{1} +\left(a +b_{3} +1\right)z\right]u_{z}  +a b_{3} u=0,} \end{aligned}\right.
$

where
$u\equiv \,\,{\rm{E}}_{344}\left(a,b_1,b_2,b_3; x,y,z\right)$

\bigskip

\begin{equation}
{\rm{E}}_{345}\left(a,b,c; x,y,z\right)=\sum\limits_{m,n,p=0}^\infty{(a)_{m+n-p}(b)_{m-n}(c)_{p-m}}\frac{x^m}{m!}\frac{y^n}{n!}\frac{z^p}{p!},
\end{equation}

region of convergence:
$$ \left\{ r<1,\,\,\,\,s<\infty,\,\,\,\,t<\infty
\right\}.
$$

System of partial differential equations:

$
 \left\{\begin{aligned}
&{x(1+x)u_{xx} -y^{2} u_{yy} } \\
&\,\,\,\,\,\,\,\,\,  -(x+1)zu_{xz} +yzu_{yz}{+\left[1- c +\left(a +b +1\right)x\right]u_{x} +(b - a -1)yu_y-b zu_{z} + a b u=0,} \\
& {yu_{yy} -xu_{xy} +xu_{x} +(1-b +y)u_{y} -zu_{z} + a u=0,} \\
& {zu_{zz} -xu_{xz} -yu_{yz} -xu_{x} +(1-a +z)u_{z} + c u=0,} \end{aligned}\right.
 $

where $u\equiv \,\,{\rm{E}}_{345}\left(a,b,c; x,y,z\right)$.

\bigskip

\begin{equation}
{\rm{E}}_{346}\left(a,b_1,b_2,b_3; x,y,z\right)=\sum\limits_{m,n,p=0}^\infty{(a)_{m}(b_1)_{m+n-p}(b_2)_{p-m}(b_3)_{p-n}}\frac{x^m}{m!}\frac{y^n}{n!}\frac{z^p}{p!},
\end{equation}

region of convergence:
$$ \left\{ r<1 \vee t<1,\,\,\,s<\infty
\right\}.
$$

System of partial differential equations:

$
 \left\{\begin{aligned}
& {x(1+x)u_{xx} +xyu_{xy} -(1+x)zu_{xz}}  {+\left[1-b_{2} +\left(a+b_{1} +1\right)x\right]u_{x} +ayu_{y} -azu_{z} +ab_{1} u=0;} \\
& {yu_{yy} -zu_{yz} +xu_{x} +(1-b_{3} +y)u_{y} -zu_{z} +b_{1} u=0,} \\
& {z(1+z)u_{zz}}+xyu_{xy} \\
&\,\,\,\,\,\,\,\,\,  -x(1+z)u_{xz} -y(1+z)u_{yz} {-b_{3} xu_{x} -b_{2} yu_{y} +\left[1-b_{1} +\left(b_{2} +b_{3} +1\right)z\right]u_{z} +b_{2} b_{3} u=0,} \end{aligned}\right.
$

where $u\equiv \,\,{\rm{E}}_{346}\left(a,b_1,b_2,b_3; x,y,z\right)$.

\bigskip

\begin{equation}
{\rm{E}}_{347}\left(a,b,c; x,y,z\right)=\sum\limits_{m,n,p=0}^\infty{(a)_{m+n-p}(b)_{p-m}(c)_{p-n}}\frac{x^m}{m!}\frac{y^n}{n!}\frac{z^p}{p!},
\end{equation}

region of convergence:
$$ \left\{ r<\infty,\,\,\,\,s<\infty,\,\,\,\,t<1
\right\}.
$$

System of partial differential equations:

$
  \left\{\begin{aligned}
 &{xu_{xx} -zu_{xz} +(1-b +x)u_{x} +yu_{y} -zu_{z} +a u=0,} \\
 & {yu_{yy} -zu_{yz} +xu_{x} +(1-c +y)u_{y} -zu_{z} +a u=0,} \\
& {z(1+z)u_{zz} +xyu_{xy}-x(1+z)u_{xz} -y(1+z)u_{yz} }\\& \,\,\,\,\,\,\,\,\, {-c xu_{x} -b yu_{y} +\left[1-a +\left(b +c +1\right)z\right]u_{z} +b c u=0,} \end{aligned}\right.
 $

where $u\equiv \,\,{\rm{E}}_{347}\left(a,b,c; x,y,z\right)$.

\bigskip

\begin{equation}
{\rm{E}}_{348}\left(a,b_1,b_2,b_3; x,y,z\right)=\sum\limits_{m,n,p=0}^\infty{(a)_{m+n}(b_1)_{n-p}(b_2)_{p-n}(b_3)_{p-m}}\frac{x^m}{m!}\frac{y^n}{n!}\frac{z^p}{p!},
\end{equation}

region of convergence:
$$ \left\{ r<\infty,\,\,\,s<1 \vee t<1
\right\}.
$$

System of partial differential equations:

$ \left\{\begin{aligned}
& {xu_{xx} -zu_{xz} +(1-b_{3} +x)u_{x} +yu_{y}  +au=0,} \\
& {y(1+y)u_{yy} +xyu_{xy} -xzu_{xz}-(y+1)zu_{yz}}\\& \,\,\,\,\,\,\,\,\, {+b_1xu_x+\left[1-b_{2} +\left(a+b_{1} +1\right)y\right]u_{y} -azu_{z} +ab_{1} u=0,} \\
& z(1+z)u_{zz} +xyu_{xy}-xzu_{xz}-y(1+z)u_{yz}\\& \,\,\,\,\,\,\,\,\, -b_{2} xu_{x}  { -b_{3} yu_{y} +\left[1-b_{1} +\left(b_{2} +b_{3} +1\right)z\right]u_{z} +b_{2} b_{3} u=0,} \end{aligned}\right.
 $

where $u\equiv \,\,{\rm{E}}_{348}\left(a,b_1,b_2,b_3; x,y,z\right)$.

\bigskip

\begin{equation}
{\rm{E}}_{349}\left(a,b_1,b_2,b_3; x,y,z\right)=\sum\limits_{m,n,p=0}^\infty{(a)_{m+n}(b_1)_{m-n}(b_2)_{n-p}(b_3)_{p-m}}\frac{x^m}{m!}\frac{y^n}{n!}\frac{z^p}{p!},
\end{equation}

region of convergence:
$$ \left\{ \left\{s+2\sqrt{rs}<1\right\}=\left\{\sqrt{s}<\sqrt{1+r}-\sqrt{r}\right\},\,\,\,t<\infty
\right\}.
$$

System of partial differential equations:

$
 \left\{\begin{aligned}
& {x(1+x)u_{xx} -y^2u_{yy} -zu_{xz} +\left[1-b_{3} +\left(a+b_{1}+1 \right)x\right]u_{x}} {-\left(a-b_{1}+1\right)yu_{y} +ab_{1} u=0,} \\
& {y(1+y)u_{yy} -x(1-y)u_{xy} -xzu_{xz} -yzu_{yz} }\\& \,\,\,\,\,\,\,\,\, {+b_{2} xu_{x} +\left[1-b_{1} +\left(a+b_{2} +1\right)y\right]u_{y} -azu_{z} +ab_{2} u=0,} \\
& {zu_{zz} -yu_{yz} -xu_{x} +(1-b_{2} +z)u_{z} +b_{3} u=0,} \end{aligned}\right.
$

where $u\equiv \,\,{\rm{E}}_{349}\left(a,b_1,b_2,b_3; x,y,z\right)$.

\bigskip

\begin{equation}
{\rm{E}}_{350}\left(a_1,a_2,b_1,b_2; x,y,z\right)=\sum\limits_{m,n,p=0}^\infty{(a_1)_m(a_2)_n(b_1)_{m+n-p}(b_2)_{p-m-n}}\frac{x^m}{m!}\frac{y^n}{n!}\frac{z^p}{p!},
\end{equation}

region of convergence:
$$ \left\{ r<1, \,\,\, s<1,\,\,\,t<\infty
\right\}.
$$

System of partial differential equations:

$
 \left\{\begin{aligned}
&{x(1+x)u_{xx} +(1+x)yu_{xy} -(1+x)zu_{xz} }\\& \,\,\,\,\,\,\,\,\, {+\left[1-b_{2} +\left(a_{1} +b_{1} +1\right)x\right]u_{x} +a_{1} yu_{y} -a_{1} zu_{z} +a_{1} b_{1} u=0,} \\
& {y(1+y)u_{yy} +x(1+y)u_{xy} -(1+y)zu_{yz} } \\& \,\,\,\,\,\,\,\,\, {+a_{2} xu_{x} +\left[1-b_{2} +\left(a_{2} +b_{1} +1\right)y\right]u_{y} -a_{2} zu_{z} +a_{2} b_{1} u=0,} \\
& {zu_{zz} -xu_{xz} -yu_{yz} -xu_{x} -yu_{y} +(1-b_{1} +z)u_{z} +b_{2} u=0,} \end{aligned}\right.
$

where $u\equiv \,\,{\rm{E}}_{349}\left(a,b_1,b_2,b_3; x,y,z\right)$.

\bigskip

\begin{equation}
{\rm{E}}_{351}\left(a_1,a_2,b_1,b_2; x,y,z\right)=\sum\limits_{m,n,p=0}^\infty{(a_1)_m(a_2)_p(b_1)_{m+n-p}(b_2)_{p-m-n}}\frac{x^m}{m!}\frac{y^n}{n!}\frac{z^p}{p!},
\end{equation}

region of convergence:
$$ \left\{ r<1 ,  \,\,t<1,\,\,\,s<\infty
\right\}.
$$

System of partial differential equations:

$
 \left\{\begin{aligned}
 &{x(1+x)u_{xx} +(1+x)yu_{xy} } -(1+x)zu_{xz} \\& \,\,\,\,\,\,\,\,\, {+\left[1-b_{2} +\left(a_{1} +b_{1} +1\right)x\right]u_{x} +a_{1} yu_{y} -a_{1} zu_{z} +a_{1} b_{1} u=0,} \\
 & {yu_{yy} +xu_{xy} -zu_{yz} +xu_{x} +\left(1-b_{2} +y\right)u_{y} -zu_{z} +b_{1} u=0,} \\
 & {z(1+z)u_{zz} -x(1+z)u_{xz} }   -y(1+z)u_{yz}\\& \,\,\,\,\,\,\,\,\,-a_{2} xu_{x}-a_{2} yu_{y} +\left[1-b_{1} +\left(a_{2} +b_{2} +1\right)z\right]u_{z} +a_{2} b_{2} u=0, \end{aligned}\right.
 $

where $u\equiv \,\,{\rm{E}}_{351}\left(a_1,a_2,b_1,b_2; x,y,z\right)$.

\bigskip

\begin{equation}
{\rm{E}}_{352}\left(a,b,c; x,y,z\right)=\sum\limits_{m,n,p=0}^\infty{(a)_m(b)_{m+n-p}(c)_{p-m-n}}\frac{x^m}{m!}\frac{y^n}{n!}\frac{z^p}{p!},
\end{equation}

region of convergence:
$$ \left\{ r<1,\,\,\,s<\infty,\,\,\,t<\infty
\right\}.
$$

System of partial differential equations:

$
  \left\{\begin{aligned}
 &{x(1+x)u_{xx} }+(1+x)yu_{xy}  -(1+x)zu_{xz} {+\left[1- c +\left(a +b +1\right)x\right]u_{x} +a yu_{y} -azu_{z} +a b u=0,} \\
 & {yu_{yy} +xu_{xy} -zu_{yz} +xu_{x} +\left(1-c +y\right)u_{y} -zu_{z} +b u=0, } \\
 & {zu_{zz} -xu_{xz} -yu_{yz} -xu_{x} -yu_{y} +(1-b +z)u_{z} + c u=0,} \end{aligned}\right.
  $

where $u\equiv \,\,{\rm{E}}_{352}\left(a,b,c; x,y,z\right)$.

\bigskip

\begin{equation}
{\rm{E}}_{353}\left(a_1,a_2,b_1,b_2; x,y,z\right)=\sum\limits_{m,n,p=0}^\infty{(a_1)_{m+n}(a_2)_m(b_1)_{n-p}(b_2)_{p-m-n}}\frac{x^m}{m!}\frac{y^n}{n!}\frac{z^p}{p!},
\end{equation}

region of convergence:
$$ \left\{ r<1, \,\, s<1,\,\,\,t<\infty
\right\},
$$

System of partial differential equations:

$
 \left\{\begin{aligned}
 & {x(1+x)u_{xx} +(1+x)yu_{xy} -zu_{xz} }  {+\left[1-b_{2} +\left(a_{1} +b_{1} +1\right)x\right]u_{x} +a_{2} yu_{y} +a_{1} a_{2} u=0,} \\
 & {y(1+y)u_{yy} +x(1+y)u_{xy} }\\
&\,\,\,\,\,\,\,\,\,-xzu_{xz} -(1+y)zu_{yz} {+b_{1} xu_{x} +\left[1-b_{2} +\left(a_{1} +b_{1} +1\right)y\right]u_{y} -a_1zu_z+a_{1} b_{1} u=0,} \\
 & {zu_{zz} -yu_{yz} -xu_{x} -yu_{y} +(1-b_{1} +z)u_{z} +b_{2} u=0,} \end{aligned}\right.
 $

where $u\equiv \,\,{\rm{E}}_{353}\left(a_1,a_2,b_1,b_2; x,y,z\right)$.

\bigskip

\begin{equation}
{\rm{E}}_{354}\left(a_1,a_2,b_1,b_2; x,y,z\right)=\sum\limits_{m,n,p=0}^\infty{(a_1)_{m+n}(a_2)_p(b_1)_{n-p}(b_2)_{p-m-n}}\frac{x^m}{m!}\frac{y^n}{n!}\frac{z^p}{p!},
\end{equation}

region of convergence:
$$ \left\{
r<\infty,\,\,\,s<1 \,\,\, t<1
\right\},
$$

System of partial differential equations:

$  \left\{\begin{aligned}
& {xu_{xx} +yu_{xy} -zu_{xz} +(1-b_{2} +x)u_{x} +yu_{y} +a_{1} u=0,} \\
& {y(1+y)u_{yy} +x(1+y)u_{xy}  }\\
&\,\,\,\,\,\,\,\,\,-xzu_{xz} -(1+y)zu_{yz}{+b_{1} xu_{x} +\left[1-b_{2} +\left(a_{1} +b_{1} +1\right)y\right]u_{y} -a_1zu_z+a_{1} b_{1} u=0,} \\
& {z(1+z)u_{zz} -xzu_{xz} -y(1+z)u_{yz}}  {-a_{2} xu_{x} -a_{2} yu_{y} +\left[1-b_{1} +\left(a_{2} +b_{2} +1\right)z\right]u_{z} +a_{2} b_{2} u=0,} \end{aligned}\right.
$

where $u\equiv \,\,{\rm{E}}_{354}\left(a_1,a_2,b_1,b_2; x,y,z\right)$.

\bigskip

\begin{equation}
{\rm{E}}_{355}\left(a,b,c; x,y,z\right)=\sum\limits_{m,n,p=0}^\infty{(a)_{m+n}(b)_{n-p}(c)_{p-m-n}}\frac{x^m}{m!}\frac{y^n}{n!}\frac{z^p}{p!},
\end{equation}

region of convergence:
$$ \left\{ r<\infty,\,\,\,s<1,\,\,\,t<\infty
\right\}.
$$

System of partial differential equations:

$
 \left\{\begin{aligned}
 & {xu_{xx} +yu_{xy} -zu_{xz} +(1- c +x)u_{x} +yu_{y} +au=0,} \\
 & {y(1+y)u_{yy}+x(1+y)u_{xy}  }\\& \,\,\,\,\,\,\,\,\,-xzu_{xz} -(1+y)zu_{yz} {+b xu_{x} +\left[1-c +(a +b +1)y\right]u_{y} -azu_z+a b u=0,} \\
 & {zu_{zz} -yu_{yz} -xu_{x} -yu_{y} +(1-b +z)u_{z} + c u=0,} \end{aligned}\right.
 $

where $u\equiv \,\,{\rm{E}}_{355}\left(a,b,c; x,y,z\right)$.

\bigskip

\begin{equation}
{\rm{E}}_{356}\left(a_1,a_2,b_1, b_2; x,y,z\right)=\sum\limits_{m,n,p=0}^\infty{(a_1)_{m+n}(a_2)_p(b_1)_{n-m}(b_2)_{m-n-p}}\frac{x^m}{m!}\frac{y^n}{n!}\frac{z^p}{p!},
\end{equation}

region of convergence:
$$ \left\{ r+s<1,\,\,\,t<\infty
\right\}.
$$

System of partial differential equations:

$
 \left\{\begin{aligned}
& {x(1+x)u_{xx} -y^2u_{yy} },-yu_{xy} -xzu_{xz}-yzu_{yz} \\& \,\,\,\,\,\,\,\,\,{+\left[1-b_{1} +\left(1+a_{1} +b_{2} \right)x\right]u_{x} -\left(a_{1} -b_{2} +1\right)yu_{y}-a_1zu_z +a_{1} b_{2} u=0,} \\
& {y(1+y)u_{yy} -x^{2} u_{xx} -xu_{xy} +zu_{yz} }\\& \,\,\,\,\,\,\,\,\, {-\left(a_{1} -b_{1} +1\right)xu_{x} +\left[1-b_{2} +\left(a_{1} +b_{1} +1\right)y\right]u_{y} +a_{1} b_{1} u=0,} \\
& {zu_{zz}-xu_{xz} +yu_{yz}  +(1-b_{2} +z)u_{z} +a_{2} u=0,} \end{aligned}\right.
 $

where $u\equiv \,\,{\rm{E}}_{356}\left(a_1,a_2,b_1, b_2; x,y,z\right)$.

\bigskip

\begin{equation}
{\rm{E}}_{357}\left(a,b, c; x,y,z\right)=\sum\limits_{m,n,p=0}^\infty{(a)_{m+n}(b)_{n-m}(c)_{m-n-p}}\frac{x^m}{m!}\frac{y^n}{n!}\frac{z^p}{p!},
\end{equation}

region of convergence:
$$ \left\{ r+s<1,\,\,\,t<\infty
\right\}.
$$

System of partial differential equations:

$
 \left\{\begin{aligned}
& {x(1+x)u_{xx} -y^2u_{yy}}-yu_{xy} \\
&\,\,\,\,\,\,\,\,\,   -xzu_{xz}-yzu_{yz}{+\left[1-b +\left(1+a+c \right)x\right]u_{x} -\left(a - c +1\right)yu_{y}-azu_z +a c u=0,} \\
& {y(1+y)u_{yy} -x^{2} u_{xx} -xu_{xy} +zu_{yz} }  {-\left(a-b +1\right)xu_{x} +\left[1- c +\left(a+b +1\right)y\right]u_{y} +ab u=0,} \\
& {zu_{zz} -xu_{xz} +yu_{yz} +(1- c )u_{z} +u=0,} \end{aligned}\right.
$

where $u\equiv \,\,{\rm{E}}_{357}\left(a,b, c; x,y,z\right)$.

\bigskip

\begin{equation}
{\rm{E}}_{358}\left(a_1,a_2,b_1, b_2; x,y,z\right)=\sum\limits_{m,n,p=0}^\infty{(a_1)_{m+n}(a_2)_n(b_1)_{p-m}(b_2)_{m-n-p}}\frac{x^m}{m!}\frac{y^n}{n!}\frac{z^p}{p!},
\end{equation}

region of convergence:
$$ \left\{ \left\{s+2\sqrt{rs}<1\right\}=\left\{\sqrt{s}<\sqrt{1+r}-\sqrt{r}\right\},\,\,\,t<\infty
\right\}.
$$

System of partial differential equations:

$  \left\{\begin{aligned}
&x(1+x)u_{xx} -y^2u_{yy}-(x+1)zu_{xz} \\
&\,\,\,\,\,\,\,\,\, -yzu_{yz} {+\left[1-b_{1} +\left(1+a_1+b_{2} \right)x \right]u_{x} -\left(a_1-b_{2} +1\right)yu_{y} -a_1zu_z+a_1b_{2} u=0,} \\
& {y(1+y)u_{yy} -x(1-y)u_{xy} +zu_{yz} } {+a_2xu_{x} +\left[1-b_{2} +\left(a_{1} +a_{2} +1\right)y\right]u_{y} +a_{1} a_{2} u=0,} \\
& {zu_{zz} -xu_{xz} +yu_{yz} -xu_{x} +(1-b_{2} +z)u_{z} +b_{1} u=0,} \end{aligned}\right.
$

where $u\equiv \,\,{\rm{E}}_{358}\left(a_1,a_2,b_1, b_2; x,y,z\right)$.

\bigskip

\begin{equation}
{\rm{E}}_{359}\left(a_1,a_2,b_1, b_2; x,y,z\right)=\sum\limits_{m,n,p=0}^\infty{(a_1)_{m+n}(a_2)_p(b_1)_{p-m}(b_2)_{m-n-p}}\frac{x^m}{m!}\frac{y^n}{n!}\frac{z^p}{p!},
\end{equation}

region of convergence:
$$ \left\{ r<1 \vee t<1,\,\,\,s<\infty
\right\}.
$$

System of partial differential equations:

$ \left\{\begin{aligned}
&x(1+x)u_{xx} -y^2u_{yy} -(x+1)zu_{xz}\\
&\,\,\,\,\,\,\,\,\,  -yzu_{yz}{+\left[1-b_{1} +\left(1+a_1+b_{2} \right)x \right]u_{x} -\left(a_1-b_{2} +1\right)yu_{y} -a_1zu_z+a_1b_{2} u=0,} \\
& {yu_{yy} -xu_{xy} +zu_{yz} +xu_{x} +(1-b_{2} +y)u_{y} +a_{1} u=0,} \\
& {z(1+z)u_{zz} -x(1+z)u_{xz} +yu_{yz}}  {-a_{2} xu_{x} +\left[1-b_{2} +\left(a_{2} +b_{1} +1\right)z\right]u_{z} +a_{2} b_{1} u=0,} \end{aligned}\right.
$

where $u\equiv \,\,{\rm{E}}_{359}\left(a_1,a_2,b_1, b_2; x,y,z\right)$.

\bigskip

\begin{equation}
{\rm{E}}_{360}\left(a,b, c; x,y,z\right)=\sum\limits_{m,n,p=0}^\infty{(a)_{m+n-p}(b)_{n+p-m}(c)_{m-n}}\frac{x^m}{m!}\frac{y^n}{n!}\frac{z^p}{p!},
\end{equation}

region of convergence:
$$ \left\{ r+s<1,\,\,\,t<\infty
\right\}.
$$

System of partial differential equations:

$
 \left\{\begin{aligned}
 &{x(1+x)u_{xx} -y^{2} u_{yy} -yu_{xy} }-(1+x)zu_{xz} \\
 &\,\,\,\,\,\,\,\,\, +yzu_{yz}{+\left[1-b + \left(1+a +c \right)x \right]u_{x} +\left(c - a -1\right)yu_{y} - c zu_{z} +ac u=0,} \\
 & {y(1+y)u_{yy}-x^2u_{xx} -z^{2} u_{zz} -xu_{xy}  } \\
&\,\,\,\,\,\,\,\,\,+2xzu_{xz} -\left(a-b+1\right)xu_x { +\left[1-c +\left(a+b +1\right)y\right]u_{y} -\left(b - a+1\right)zu_z+ab u=0,} \\
 & {zu_{zz} -xu_{xz} -yu_{yz} -xu_{x} +yu_{y} +(1-a +z)u_{z} +b u=0,} \end{aligned}\right.
  $

where $u\equiv \,\,{\rm{E}}_{360}\left(a,b, c; x,y,z\right)$.

\bigskip

\begin{equation}
{\rm{E}}_{361}\left(a, b, c; x,y,z\right)=\sum\limits_{m,n,p=0}^\infty{(a)_{n+p}(b)_{m+n-p}(c)_{p-m-n}}\frac{x^m}{m!}\frac{y^n}{n!}\frac{z^p}{p!},
\end{equation}

region of convergence:
$$ \left\{ r<\infty,\,\,\,s+t<1
\right\}.
$$

System of partial differential equations:

$ \left\{\begin{aligned}
&{xu_{xx} +yu_{xy} -zu_{xz} +\left(1-c +x\right)u_{x} +yu_{y} -zu_{z} +b u=0,} \\
& {y(1+y)u_{yy} -z^{2} u_{zz} +x(1+y)u_{xy}  } \\
&\,\,\,\,\,\,\,\,\, +xzu_{xz} -zu_{yz}+axu_{x}{ +\left[1-c +\left(a+b +1\right)y\right]u_{y} -(a-b+1)zu_z+ab u=0,} \\
& {z(1+z)u_{zz} -y^{2} u_{yy}-xyu_{xy}  } \\
&\,\,\,\,\,\,\,\,\,-x(1+z)u_{xz} -yu_{yz}{-axu_{x} -\left(a- c +1\right)yu_y+ \left[1-b + \left(a+ c +1\right)z \right]u_{z} +a cu=0,} \end{aligned}\right.
 $

where $u\equiv \,\,{\rm{E}}_{361}\left(a, b, c; x,y,z\right)$.

\bigskip

\begin{equation}
{\rm{E}}_{362}\left(a, b, c; x,y,z\right)=\sum\limits_{m,n,p=0}^\infty{(a)_{n+p}(b)_{p-m-n}(c)_{m-p}}\frac{x^m}{m!}\frac{y^n}{n!}\frac{z^p}{p!},
\end{equation}

region of convergence:
$$ \left\{ r<\infty,\,\,\,s<\infty,\,\,\,t<1
\right\}.
$$

System of partial differential equations:

$ \left\{\begin{aligned}
&{xu_{xx} +yu_{xy} -zu_{xz} +(1-b+x )u_{x} -zu_{z} + cu=0,} \\
& {y u_{yy} +xu_{xy} -zu_{yz} +(1-b +y)u_{y} +zu_{z} +au=0,} \\
& {z(1+z)u_{zz} -y^{2} u_{yy}}-xyu_{xy} -x(1+z)u_{xz} \\
&\,\,\,\,\,\,\,\,\, {-axu_{x} -\left(a-b +1\right)yu_y+ \left[1- c + \left(a+b +1\right)z \right]u_{z} +ab u=0,} \end{aligned}\right.
$

where $u\equiv \,\,{\rm{E}}_{362}\left(a, b, c; x,y,z\right)$.

\bigskip

\begin{equation}
{\rm{E}}_{363}\left(a, b, c; x,y,z\right)=\sum\limits_{m,n,p=0}^\infty{(a)_{m+n+p}(b)_{m-p}(c)_{p-m-n}}\frac{x^m}{m!}\frac{y^n}{n!}\frac{z^p}{p!},
\end{equation}

region of convergence:
$$ \left\{ r+t<1,\,\,\,s<\infty
\right\}.
$$

System of partial differential equations:

$
 \left\{\begin{aligned}& {x(1+x)u_{xx} -z^{2} u_{zz} +(1+x)yu_{xy} } \\
&\,\,\,\,\,\,\,\,\,-yzu_{yz} -zu_{xz} {+\left[1- c +\left(a+b+1 \right)x\right]u_{x} +b yu_{y} -\left(a-b+1 \right)zu_{z} +ab u=0,} \\
& {yu_{yy} +xu_{xy} -zu_{yz} +xu_{x} +(1-c +y)u_{y} +zu_{z} +au=0,} \\
& { z(1+z)u_{zz}- x^{2} u_{xx} -y^{2} u_{yy}} -2xyu_{xy}\\
&\,\,\,\,\,\,\,\,\,  -xu_{xz} -(a-c +1)xu_{x} {-(a-c +1)yu_y+ \left[1-b +(a+c +1)z\right]u_{z} +ac u=0,} \end{aligned}\right.
 $

where $u\equiv \,\,{\rm{E}}_{363}\left(a, b, c; x,y,z\right)$.

\bigskip

\begin{equation}
{\rm{E}}_{364}\left(a, b_1, b_2, b_3; x,y,z\right)=\sum\limits_{m,n,p=0}^\infty{(a)_{n}(b_1)_{m-n}(b_2)_{m-p}(b_3)_{2p-m}  }\frac{x^m}{m!}\frac{y^n}{n!}\frac{z^p}{p!},
\end{equation}

region of convergence:
$$ \left\{ \left\{r^2t+r<1\right\}=\left\{r<\frac{\sqrt{1+4t}-1}{2t}\right\},\,\,\,s<\infty
\right\}.
$$

System of partial differential equations:

$
 \left\{\begin{aligned} &{x(1+x)u_{xx} -xyu_{xy} }+yzu_{yz} -(2+x)zu_{xz}\\& \,\,\,\,\,\,\,\,\,{+\left[1-b_{3} + \left(b_{1} +b_{2}+1 \right)x \right]u_{x} -b_{2} yu_{y} -b_{1} zu_{z} +b_{1} b_{2} u=0,} \\
 & {yu_{yy} -xu_{xy} +(1-b_{1} +y)u_{y} +au=0,} \\
 & {z(1+4z)u_{zz} +x^{2} u_{xx} -x(1+4z)u_{xz} -2b_{3} xu_{x}} {+\left[1-b_{2} +\left(4b_{3} +6\right)z\right]u_{z} -b_{3} (1+b_{3} )u=0,} \end{aligned}\right.
  $

where $u\equiv \,\,{\rm{E}}_{364}\left(a, b_1, b_2, b_3; x,y,z\right)$.

\bigskip

\begin{equation}
{\rm{E}}_{365}\left(a, b, c; x,y,z\right)=\sum\limits_{m,n,p=0}^\infty{(a)_{m-n}(b)_{m-p}(c)_{2p-m}}\frac{x^m}{m!}\frac{y^n}{n!}\frac{z^p}{p!},
\end{equation}

region of convergence:
$$ \left\{ \left\{r^2t+r<1\right\}=\left\{r<\frac{\sqrt{1+4t}-1}{2t}\right\},\,\,\,s<\infty
\right\}.
$$

System of partial differential equations:

$
 \left\{\begin{aligned} &{x(1+x)u_{xx} -xyu_{xy} } -(2+x)zu_{xz} +yzu_{yz}\\& \,\,\,\,\,\,\,\,\,{+\left[1-c +\left(a +b+1 \right)x \right]u_{x} -b yu_{y} - a zu_{z} + a b u=0,} \\
 & {yu_{yy} -xu_{xy} +(1-a )u_{y} +u=0,} \\
 & { z(1+4z)u_{zz}+ x^{2} u_{xx} -x(1+4z)u_{xz} -2 c xu_{x}}{+\left[1-b +\left(4 c+6\right)z\right]u_{z} + c \left(1+c \right)u=0,} \end{aligned}\right.
 $

where $u\equiv \,\,{\rm{E}}_{365}\left(a, b, c; x,y,z\right)$.

\bigskip

\begin{equation}
{\rm{E}}_{366}\left(a, b_1, b_2, b_3; x,y,z\right)=\sum\limits_{m,n,p=0}^\infty{(a)_{m}(b_1)_{m-n}(b_2)_{n-p}(b_3)_{2p-m}  }\frac{x^m}{m!}\frac{y^n}{n!}\frac{z^p}{p!},
\end{equation}

region of convergence:
$$ \left\{ r(1+2\sqrt{t})<1,\,\,\,s<\infty
\right\}.
$$

System of partial differential equations:

$
 \left\{\begin{aligned}
 & {x(1+x)u_{xx} -xyu_{xy} -2zu_{xz}} {+\left[1-b_{3} +\left(a+b_{1}+1 \right)x \right]u_{x} -ayu_{y} +ab_{1} u=0,} \\
 & {yu_{yy} -xu_{xy} +(1-b_{1} +y)u_{y} -zu_{z} +b_{2} u=0,} \\
 & { z(1+4z)u_{zz}+x^{2} u_{xx} -yu_{yz} -4xzu_{xz} -2b_{3} xu_{x}}  {+\left[1-b_{2} +\left(4b_{3} +6\right)z\right]u_{z} +b_{3} \left(1+b_{3} \right)u=0,} \end{aligned}\right.
  $

where $u\equiv \,\,{\rm{E}}_{366}\left(a, b_1, b_2, b_3; x,y,z\right)$.

\bigskip

\begin{equation}
{\rm{E}}_{367}\left(a, b_1, b_2, b_3; x,y,z\right)=\sum\limits_{m,n,p=0}^\infty{(a)_{n}(b_1)_{m-n}(b_2)_{n-p}(b_3)_{2p-m}  }\frac{x^m}{m!}\frac{y^n}{n!}\frac{z^p}{p!},
\end{equation}

region of convergence:
$$ \left\{ r<\infty,\,\,\,4t(1+s)<1
\right\}.
$$

System of partial differential equations:

$
\left\{\begin{aligned}
&{xu_{xx} -2zu_{xz} +\left(1-b_{3} +x\right)u_{x} -yu_{y} +b_{1} u=0,} \\
& y(1+y)u_{yy} -xu_{xy} -yzu_{yz}  +\left[1-b_{1}+\left(a+b_{2} +1\right)y\right]u_{y} {-azu_{z} +ab_{2} u=0,} \\
& { z(1+4z)u_{zz}+x^{2} u_{xx} -yu_{yz} -4xzu_{xz} -2b_{3} xu_{x} } {+\left[1-b_{2} +\left(4b_{3} +6\right)z\right]u_{z} +b_{3} \left(1+b_{3} \right)u=0,} \end{aligned}\right.
$

where $u\equiv \,\,{\rm{E}}_{367}\left(a, b_1, b_2, b_3; x,y,z\right)$.

\bigskip

\begin{equation}
{\rm{E}}_{368}\left( a, b, c; x,y,z\right)=\sum\limits_{m,n,p=0}^\infty{(a)_{m-n}(b)_{n-p}(c)_{2p-m}}\frac{x^m}{m!}\frac{y^n}{n!}\frac{z^p}{p!},
\end{equation}

region of convergence:
$$ \left\{ r<\infty,\,\,\,s<\infty,\,\,\,t<\frac{1}{4}
\right\}.
$$

System of partial differential equations:

$
 \left\{\begin{aligned}
&{xu_{xx} -2zu_{xz} +\left(1-c +x\right)u_{x} -yu_{y} + a u=0,} \\
& {yu_{yy} -xu_{xy} +(1-a +y)u_{y} -zu_{yz} +b u=0,} \\
& { z(1+4z)u_{zz}+x^{2} u_{xx} -yu_{yz} -4xzu_{xz} -2 c xu_{x}}  {+\left[1-b +\left(4c+6\right)z\right]u_{z} + c \left(1+c\right)u=0,} \end{aligned}\right.
$

where $u\equiv \,\,{\rm{E}}_{368}\left( a, b, c; x,y,z\right)$.

\bigskip

\begin{equation}
{\rm{E}}_{369}\left(a_1,a_2, b_1, b_2; x,y,z\right)=\sum\limits_{m,n,p=0}^\infty{(a_1)_{m}(a_2)_{n}(b_1)_{m-p}(b_2)_{2p-m-n}}\frac{x^m}{m!}\frac{y^n}{n!}\frac{z^p}{p!},
\end{equation}

region of convergence:
$$ \left\{ \left\{r^2t+r<1\right\}=\left\{r<\frac{\sqrt{1+4t}-1}{2t}\right\},\,\,\,s<\infty
\right\}.
$$

System of partial differential equations:

$
 \left\{\begin{aligned}
&{x(1+x)u_{xx} +yu_{xy} -(x+2)zu_{xz}}  {+\left[1-b_{2} + \left(a_{1} +b_{1}+1 \right)x\right]u_{x} -a_{1} zu_{z} +a_{1} b_{1} u=0,} \\
& {yu_{yy} +xu_{xy} -2zu_{yz} +(1-b_{2} +y)u_{y} +a_{2} u=0,} \\
& { z(1+4z)u_{zz}+x^{2} u_{xx} +y^{2} u_{yy} -x(1+4z)u_{xz} } \\
&\,\,\,\,\,\,\,\,\,+2xyu_{xy} -4yzu_{yz} {-2b_{2} xu_{x} -2b_{2} yu_{y} +\left[1-b_{1} +\left(4b_{2} +6\right)z\right]u_{z} +b_{2} \left(1+b_{2} \right)u=0,} \end{aligned}\right.
$

where $u\equiv \,\,{\rm{E}}_{369}\left(a_1,a_2, b_1, b_2; x,y,z\right)$.

\bigskip

\begin{equation}
{\rm{E}}_{370}\left(a_1,a_2, b_1, b_2; x,y,z\right)=\sum\limits_{m,n,p=0}^\infty{(a_1)_{n}(a_2)_n(b_1)_{m-p}(b_2)_{2p-m-n}  }\frac{x^m}{m!}\frac{y^n}{n!}\frac{z^p}{p!},
\end{equation}

region of convergence:
$$ \left\{ r<\infty,\,\,\,s(1+2\sqrt{t})<1
\right\}.
$$

System of partial differential equations:

$
\left\{\begin{aligned}
&{xu_{xx} +yu_{xy} -2zu_{xz} +\left(1-b_{2}+x \right)u_{x} -zu_{z} +b_{1} u=0,} \\
& y(1+y)u_{yy} +xu_{xy} -2zu_{yz} +\left[1-b_{2} +(a_{1} +a_{2} +1)y\right]u_{y} +a_{1} a_{2} u=0,\\
& { z(1+4z)u_{zz}+ x^{2} u_{xx} +y^{2} u_{yy}-x(1+4z)u_{xz}  } \\
&\,\,\,\,\,\,\,\,\,+2xyu_{xy} -4yzu_{yz}{-2b_{2} xu_{x} -2b_{2} yu_{y} +\left[1-b_{1} +\left(4b_{2} +6\right)z\right]u_{z} +b_{2} (1+b_{2} )u=0,} \end{aligned}\right.
$

where $u\equiv \,\,{\rm{E}}_{370}\left(a_1,a_2, b_1, b_2; x,y,z\right)$.

\bigskip

\begin{equation}
{\rm{E}}_{371}\left(a, b, c; x,y,z\right)=\sum\limits_{m,n,p=0}^\infty{(a)_{m}(b)_{m-p}(c)_{2p-m-n}  }\frac{x^m}{m!}\frac{y^n}{n!}\frac{z^p}{p!},
\end{equation}

region of convergence:
$$ \left\{ \left\{r^2t+r<1\right\}=\left\{r<\frac{\sqrt{1+4t}-1}{2t}\right\},\,\,\,s<\infty
\right\}.
$$

System of partial differential equations:

$
 \left\{\begin{aligned}
 &{x(1+x)u_{xx} +yu_{xy} -(x+2)zu_{xz}} {+\left[1-c + \left(1+a +b \right)x \right]u_{x} -a zu_{z} +a b u=0,} \\
 & {yu_{yy} +xu_{xy} -2zu_{yz} +(1- c )u_{y} +u=0,} \\
 & { z(1+4z)u_{zz} + x^{2} u_{xx} +y^{2} u_{yy}-x(1+4z)u_{xz} } \\
 &\,\,\,\,\,\,\,\,\,+2xyu_{xy} -4yzu_{yz} {-2c xu_{x} -2c yu_{y} +\left[1-b+\left(4c +6\right)z\right]u_{z} +c(1+c)u=0,} \end{aligned}\right.
 $

where $u\equiv \,\,{\rm{E}}_{371}\left(a, b, c; x,y,z\right)$.

\bigskip

\begin{equation}
{\rm{E}}_{372}\left(a, b, c; x,y,z\right)=\sum\limits_{m,n,p=0}^\infty{(a)_{n}(b)_{m-p}(c)_{2p-m-n}}\frac{x^m}{m!}\frac{y^n}{n!}\frac{z^p}{p!},
\end{equation}

region of convergence:
$$ \left\{ r<\infty,\,\,\,s<\infty,\,\,\,t<\frac{1}{4}
\right\}.
$$

System of partial differential equations:

$
 \left\{\begin{aligned}
&{xu_{xx} +yu_{xy} -2zu_{xz} +\left(1-c+x \right)u_{x} -zu_{z} +b u=0,} \\
& {yu_{yy} +xu_{xy} -2zu_{yz} +(1-c +y)u_{y} +au=0,} \\
& { z(1+4z)u_{zz}+x^{2} u_{xx} +y^{2} u_{yy} -x(1+4z)u_{xz}}\\
&\,\,\,\,\,\,\,\,\, +2xyu_{xy} -4yzu_{yz} {-2c xu_{x} -2c yu_{y} +\left[1-b +\left(4c+6\right)z\right]u_{z} +c (1+c )u=0,} \end{aligned}\right.
$

where $u\equiv \,\,{\rm{E}}_{372}\left(a, b, c; x,y,z\right)$.

\bigskip

\begin{equation}
{\rm{E}}_{373}\left( a, b; x,y,z\right)=\sum\limits_{m,n,p=0}^\infty{(a)_{m-p}(b)_{2p-m-n}}\frac{x^m}{m!}\frac{y^n}{n!}\frac{z^p}{p!},
\end{equation}

region of convergence:
$$ \left\{ r<\infty,\,\,\,s<\infty,\,\,\,t<\frac{1}{4}
\right\}.
$$

System of partial differential equations:

$
 \left\{\begin{aligned}
 &{xu_{xx} +yu_{xy} -2zu_{xz} +\left(1-b+x \right)u_{x} -zu_{z} +a u=0,} \\
 & {yu_{yy} +xu_{xy} -2zu_{yz} +(1-b )u_{y} +u=0,} \\
 & { z(1+4z)u_{zz}+x^{2} u_{xx} +y^{2} u_{yy}  } -x(1+4z)u_{xz}\\
 &\,\,\,\,\,\,\,\,\,+2xyu_{xy} -4yzu_{yz} {-2b xu_{x} -2b yu_{y} +\left[1-a +\left(4b+6\right)z\right]u_{z} +b (1+b )u=0,} \end{aligned}\right.
 $

where $u\equiv \,\,{\rm{E}}_{373}\left( a, b; x,y,z\right)$.

\bigskip

\begin{equation}
{\rm{E}}_{374}\left(a_1,a_2, b_1, b_2; x,y,z\right)=\sum\limits_{m,n,p=0}^\infty{(a_1)_{m}(a_2)_{n}(b_1)_{2p-m}(b_2)_{m-n-p}}\frac{x^m}{m!}\frac{y^n}{n!}\frac{z^p}{p!},
\end{equation}

region of convergence:
$$ \left\{ \left\{r^2t+r<1\right\}=\left\{r<\frac{\sqrt{1+4t}-1}{2t}\right\},\,\,\,s<\infty
\right\}.
$$

System of partial differential equations:

$
 \left\{\begin{aligned}
 &{x(1+x)u_{xx} -xyu_{xy} } -(x+2)zu_{xz}\\& \,\,\,\,\,\,\,\,\, {+\left[1-b_{1} +\left(a_{1} +b_{2}+1 \right)x\right]u_{x} -a_{1} yu_{y} -a_{1} zu_{z} +a_{1} b_{2} u=0,} \\
 & {yu_{yy} -xu_{xy} +zu_{yz} +(1-b_{2} +y)u_{y} +a_{2} u=0,} \\
 & { z(1+4z)u_{zz}+ x^{2} u_{xx} -x(1+4z)u_{xz} +yu_{yz} }\\& \,\,\,\,\,\,\,\,\, {-2b_{1} xu_{x} +\left[1-b_{2} +\left(4b_{1} +6\right)z\right]u_{z} +b_{1} (1+b_{1} )u=0,} \end{aligned}\right.
 $

where $u\equiv \,\,{\rm{E}}_{374}\left(a_1,a_2, b_1, b_2; x,y,z\right)$.

\bigskip

\begin{equation}
{\rm{E}}_{375}\left(a_1,a_2, b_1, b_2; x,y,z\right)=\sum\limits_{m,n,p=0}^\infty{(a_1)_{n}(a_2)_n(b_1)_{2p-m}(b_2)_{m-n-p}  }\frac{x^m}{m!}\frac{y^n}{n!}\frac{z^p}{p!},
\end{equation}

region of convergence:
$$ \left\{ r<\infty,\,\,\,\frac{1}{s}+\frac{1}{4t}>1
\right\}.
$$

System of partial differential equations:

$
  \left\{\begin{aligned}
 & {xu_{xx} -2zu_{xz} +\left(1-b_{1} +x\right)u_{x} -yu_{y} -zu_{z} +b_{2} u=0,} \\
 & {y(1+y)u_{yy} -xu_{xy} +zu_{yz}} {+\left[1-b_{2} +\left(a_{1} +a_{2} +1\right)y\right]u_{y} +a_{1} a_{2} u=0,} \\
 & z(1+4z)u_{zz}+ x^{2} u_{xx} \\& \,\,\,\,\,\,\,\,\,-x(1+4z)u_{xz} +yu_{yz}  {-2b_{1} xu_{x} +\left[1-b_{2} +\left(4b_{1} +6\right)z\right]u_{z} +b_{1} (1+b_{1} )u=0,} \end{aligned}\right.
 $

where $u\equiv \,\,{\rm{E}}_{375}\left(a_1,a_2, b_1, b_2; x,y,z\right)$.

\bigskip

\begin{equation}
{\rm{E}}_{376}\left(a,b, c; x,y,z\right)=\sum\limits_{m,n,p=0}^\infty{(a)_{m}(b)_{2p-m}(c)_{m-n-p}  }\frac{x^m}{m!}\frac{y^n}{n!}\frac{z^p}{p!},
\end{equation}

region of convergence:
$$ \left\{ \left\{r^2t+r<1\right\}=\left\{r<\frac{\sqrt{1+4t}-1}{2t}\right\},\,\,\,s<\infty
\right\}.
$$

System of partial differential equations:

$
 \left\{\begin{aligned}
 &{x(1+x)u_{xx} -xyu_{xy} -(2+x)zu_{xz}} {+\left[1-b+\left(1+a+c\right)x\right]u_{x} -ayu_{y} -azu_{z} +ac u=0,} \\
 & {yu_{yy} -xu_{xy} +zu_{yz} +(1-c)u_{y} +u=0,} \\
 & { z(1+4z)u_{zz}+ x^{2} u_{xx}-x(1+4z)u_{xz} +yu_{yz} } {-2bxu_{x} +\left[1-c+\left(4b+6\right)z\right]u_{z} +b (1+b)u=0,} \end{aligned}\right.
  $

where $u\equiv \,\,{\rm{E}}_{376}\left(a,b, c; x,y,z\right)$.

\bigskip

\begin{equation}
{\rm{E}}_{377}\left(a, b, c; x,y,z\right)=\sum\limits_{m,n,p=0}^\infty{(a)_{n}(b)_{2p-m}(c)_{m-n-p}  }\frac{x^m}{m!}\frac{y^n}{n!}\frac{z^p}{p!},
\end{equation}

region of convergence:
$$ \left\{ r<\infty,\,\,\,s<\infty,\,\,\,t<\frac{1}{4}
\right\}.
$$

System of partial differential equations:

$ \left\{\begin{aligned}
& {xu_{xx} -2zu_{xz} +\left(1-b+x\right)u_{x} -yu_{y} -zu_{z} + c u=0,} \\
& {yu_{yy} -xu_{xy} +zu_{yz} +(1-c +y)u_{y} +au=0,} \\
& { z(1+4z)u_{zz}+x^{2} u_{xx} -x(1+4z)u_{xz} +yu_{yz} } {-2b xu_{x}+\left[1-c +\left(4b +6\right)z\right]u_{z} +b(1+b)u=0,} \end{aligned}\right.
$

where $u\equiv \,\,{\rm{E}}_{377}\left(a, b, c; x,y,z\right)$.

\bigskip

\begin{equation}
{\rm{E}}_{378}\left(a, b; x,y,z\right)=\sum\limits_{m,n,p=0}^\infty{(a)_{2p-m}(b)_{m-n-p}}\frac{x^m}{m!}\frac{y^n}{n!}\frac{z^p}{p!},
\end{equation}

region of convergence:
$$ \left\{ r<\infty,\,\,\,s<\infty,\,\,\,t<\frac{1}{4}
\right\},
$$

System of partial differential equations:

$
 \left\{\begin{aligned}
& {xu_{xx} -2zu_{xz} +\left(1- a +x\right)u_{x} -yu_{y} -zu_{z} +b u=0,} \\
& {yu_{yy} -xu_{xy} +zu_{yz} +(1-b )u_{y} +u=0,} \\
& {z(1+4z)u_{zz}+x^{2} u_{xx}  -x(1+4z)u_{xz} +yu_{yz} } {-2axu_{x} +\left[1-b +\left(4a +6\right)z\right]u_{z} +a (1+a)u=0,} \end{aligned}\right.
$

where $u\equiv \,\,{\rm{E}}_{378}\left(a, b; x,y,z\right)$.

\bigskip

\begin{equation}
{\rm{E}}_{379}\left(a, b, c; x,y,z\right)=\sum\limits_{m,n,p=0}^\infty{(a)_{2p-m}(b)_{m-n}(c)_{m+n-p}}\frac{x^m}{m!}\frac{y^n}{n!}\frac{z^p}{p!},
\end{equation}

region of convergence:
$$ \left\{
\left\{r^2t+r<1\right\}=\left\{r<\frac{\sqrt{1+4t}-1}{2t}\right\},\,\,\,s<\infty \right\}.
$$

System of partial differential equations:

$
  \left\{\begin{aligned}
 & {x(1+x)u_{xx} -y^{2} u_{yy} }  -(2+x)zu_{xz} +yzu_{yz}\\& \,\,\,\,\,\,\,\,\,+\left[1-a + \left(b +c+1 \right)x \right]u_{x} +\left(b - c+1 \right)yu_{y} -b zu_{z} +b c u=0, \\
 & {yu_{yy} -xu_{xy} +xu_{x} +(1-b +y)u_{y} -zu_{z} + c u=0,} \\
 & { z(1+4z)u_{zz}+x^{2} u_{xx} -x(1+4z)u_{xz} -yu_{yz}} {-2a xu_{x} +\left[1- c +\left(4a+6\right)z\right]u_{z} + a (1+ a )u=0,} \end{aligned}\right.
$

where $u\equiv \,\,{\rm{E}}_{379}\left(a, b, c; x,y,z\right)$.

\bigskip

\begin{equation}
{\rm{E}}_{380}\left(a, b, c; x,y,z\right)=\sum\limits_{m,n,p=0}^\infty{(a)_{2p-m}(b)_{m-n}(c)_{n-p}}\frac{x^m}{m!}\frac{y^n}{n!}\frac{z^p}{p!},
\end{equation}

region of convergence:
$$ \left\{ r<\infty,\,\,\,s<\infty,\,\,\,t<\frac{1}{4}
\right\}.
$$

System of partial differential equations:

$
 \left\{\begin{aligned}
&{xu_{xx} -2zu_{xz} +\left(1-a +x\right)u_{x} -yu_{y} +bu=0,} \\
& {yu_{yy} -xu_{xy} +(1-b +y)u_{y} -zu_{z} + c u=0,} \\
& { z(1+4z)u_{zz} +x^{2} u_{xx} -4xzu_{xz} -yu_{yz}} {-2a xu_{x} +\left[1-c +\left(4a+6\right)z\right]u_{z} +a (1+a )u=0,} \end{aligned}\right.
 $

where $u\equiv \,\,{\rm{E}}_{380}\left(a, b, c; x,y,z\right)$.

\bigskip

\begin{equation}
{\rm{E}}_{381}\left(a, b, c; x,y,z\right)=\sum\limits_{m,n,p=0}^\infty{(a)_{m}(b)_{2p-m-n}(c)_{m+n-p}}\frac{x^m}{m!}\frac{y^n}{n!}\frac{z^p}{p!},
\end{equation}

region of convergence:
$$ \left\{ \left\{r^2t+r<1\right\}=\left\{r<\frac{\sqrt{1+4t}-1}{2t}\right\},\,\,\,s<\infty
\right\},
$$

System of partial differential equations:

$
 \left\{\begin{aligned}
&{x(1+x)u_{xx} }+(1+x)yu_{xy} -(2+x)zu_{xz}{+\left[1-b +\left(1+a+c\right)x\right]u_{x} +ayu_{y} -azu_{z} +ac u=0,} \\
& {yu_{yy} +xu_{xy} -2zu_{yz}+xu_{x} +(1-b +y)u_{y}  -zu_{z} + c u=0,} \\
& {z(1+4z)u_{zz}+ x^{2} u_{xx} +y^{2} u_{yy}  +2xyu_{xy}} \\
&\,\,\,\,\,\,\,\,\, -x(1+4z)u_{xz} -y(1+4z)u_{xz}{-2b xu_{x} -2b yu_{y} +\left[1-c +\left(4b+6\right)z\right]u_{z} +b(1+b)u=0,} \end{aligned}\right.
$

where $u\equiv \,\,{\rm{E}}_{381}\left(a, b, c; x,y,z\right)$.

\bigskip

\begin{equation}
{\rm{E}}_{382}\left(a, b, c; x,y,z\right)=\sum\limits_{m,n,p=0}^\infty{(a)_{m+n}(b)_{2p-m-n}(c)_{m-p}}\frac{x^m}{m!}\frac{y^n}{n!}\frac{z^p}{p!},
\end{equation}

region of convergence:
$$ \left\{ \left\{r^2t+r<1\right\}=\left\{r<\frac{\sqrt{1+4t}-1}{2t}\right\},\,\,\,s<\infty
\right\}.
$$

System of partial differential equations:

$
 \left\{\begin{aligned}
& {x(1+x)u_{xx} +(1+x)yu_{xy} }\\& \,\,\,\,\,\,\,\,\, -yzu_{yz} -(x+2)zu_{xz} {+\left[1-b +\left(a+c+1 \right)x\right]u_{x} +c yu_{y} -azu_{z} +a c u=0,} \\
& {yu_{yy} +xu_{xy} -2zu_{yz} +xu_{x} +(1-b +y)u_{y} +au=0,} \\
& { z(1+4z)u_{zz}+ x^{2} u_{xx} +y^{2} u_{yy} +2xyu_{xy} } \\
&\,\,\,\,\,\,\,\,\,-x(1+4z)u_{xz} -4yzu_{xz} {-2b xu_{x} -2b yu_{y} +\left[1-c +\left(4b+6\right)z\right]u_{z} +b (1+b)u=0,} \end{aligned}\right.
$

where $u\equiv \,\,{\rm{E}}_{382}\left(a, b, c; x,y,z\right)$.

\bigskip

\begin{equation}
{\rm{E}}_{383}\left(a, b, c; x,y,z\right)=\sum\limits_{m,n,p=0}^\infty{(a)_{n}(b)_{2p-m-n}(c)_{m-p}}\frac{x^m}{m!}\frac{y^n}{n!}\frac{z^p}{p!},
\end{equation}

region of convergence:
$$ \left\{ r<\infty,\,\,\,s<\infty,\,\,\,t<\frac{1}{4}
\right\}.
$$

System of partial differential equations:

$
\left\{\begin{aligned}
 & {xu_{xx} +yu_{xy} -2zu_{xz} +(1-b+x)u_{x} -zu_{z} + cu=0,} \\
 & {yu_{yy} +xu_{xy} -2zu_{yz} +(1-b +y)u_{y} +au=0,} \\
 & { z(1+4z)u_{zz}+x^{2} u_{xx} +y^{2} u_{yy} } +2xyu_{xy}\\
 &\,\,\,\,\,\,\,\,\, -x(1+4z)u_{xz} -4yzu_{xz} {-2b xu_{x} -2b yu_{y} +\left[1-c+\left(4b+6\right)z\right]u_{z} +b (1+b )u=0,} \end{aligned}\right.
 $

where $u\equiv \,\,{\rm{E}}_{383}\left(a, b, c; x,y,z\right)$.

\bigskip

\begin{equation}
{\rm{E}}_{384}\left(a,b, c; x,y,z\right)=\sum\limits_{m,n,p=0}^\infty{(a)_{m+n}(b)_{2p-m}(c)_{m-n-p}}\frac{x^m}{m!}\frac{y^n}{n!}\frac{z^p}{p!},
\end{equation}

region of convergence:
$$ \left\{ \left\{r^2t+r<1\right\}=\left\{r<\frac{\sqrt{1+4t}-1}{2t}\right\},\,\,\,s<\infty
\right\}.
$$

System of partial differential equations:

$
\left\{\begin{aligned}
&{x(1+x)u_{xx} -y^{2} u_{yy}  } \\
&\,\,\,\,\,\,\,\,\, -(x+2)zu_{xz} -yzu_{yz} {+\left[1-b +\left(a+c+1 \right)x\right]u_{x} -\left(a-c+1 \right)yu_{y} -azu_{z} +ac u=0,} \\
& {yu_{yy} -xu_{xy} +zu_{yz} +xu_{x} +(1-c +y)u_{y} +au=0,} \\
& {z(1+4z)u_{zz}+x^{2} u_{xx} -x(1+4z)u_{xz} +yu_{xz}} {-2b xu_{x} +\left[1-c +\left(4b +6\right)z\right]u_{z} +b (1+b )u=0,} \end{aligned}\right.
$

where $u\equiv \,\,{\rm{E}}_{384}\left(a,b, c; x,y,z\right)$.

\bigskip

\begin{equation}
{\rm{E}}_{385}\left(a, b, c; x,y,z\right)=\sum\limits_{m,n,p=0}^\infty{(a)_{2m+n-p}(b)_{n-m}(c)_{p-n}}\frac{x^m}{m!}\frac{y^n}{n!}\frac{z^p}{p!},
\end{equation}

region of convergence:
$$ \left\{ r<\frac{1}{4}\wedge s<\min\left\{\Psi_1(r), \Psi_2(r)\right\}=s<1 \wedge r<\min\left\{\Theta_1(s), \Theta_2(s)\right\},\,\,\,t<\infty
\right\}.
$$

System of partial differential equations:

$
 \left\{\begin{aligned}
&{x(1+4x)u_{xx} +y^{2} u_{yy} +z^{2} u_{zz} -(1-4x)yu_{xy}  } \\
&\,\,\,\,\,\,\,\,\,-4xzu_{xz} -2yzu_{yz} {+\left[1-a +\left(4a +6\right)x\right]u_{x} +2\left(a+1 \right)yu_{y} -2a zu_{z} +a (a +1)u=0,} \\
& {y(1+y)u_{yy}-2x^2u_{xx} +xyu_{xy} } +xzu_{xz} \\
 &\,\,\,\,\,\,\,\,\, -(y+1)zu_{yz}{-\left(a-2b+2\right)xu_{x} +\left[1-c +\left(a+b +1\right)y\right]u_{y} -bzu_{z} +ab u=0,} \\
& {zu_{zz} -2xu_{xz} -yu_{yz} -yu_{y} +(1-a +z)u_{z} +c u=0,} \end{aligned}\right.
$

where $u\equiv \,\,{\rm{E}}_{385}\left(a, b, c; x,y,z\right)$.

\bigskip

\begin{equation}
{\rm{E}}_{386}\left(a, b, c; x,y,z\right)=\sum\limits_{m,n,p=0}^\infty{(a)_{2m+n-p}(b)_{p-m}(c)_{p-n}}\frac{x^m}{m!}\frac{y^n}{n!}\frac{z^p}{p!},
\end{equation}

region of convergence:
$$ \left\{ \left\{rt^2+t<1\right\}=\left\{t<\frac{\sqrt{1+4r}-1}{2r}\right\},\,\,\,s<\infty
\right\}.
$$

System of partial differential equations:

$  \left\{\begin{aligned}
&{x(1+4x)u_{xx} +y^{2} u_{yy} +z^{2} u_{zz} +4xyu_{xy} }-(1+4x)zu_{xz} \\
&\,\,\,\,\,\,\,\,\,  -2yzu_{yz} {+\left[1-b +\left(4a +6\right)x \right]u_{x} +2\left(a+1 \right)yu_{y} -2a zu_{z} +a (a +1)u=0,} \\
& {yu_{yy} -zu_{yz} +2xu_{x} +(1-c +y)u_{y} -zu_{z} +a u=0,} \\
& {z(1+z)u_{zz} +xyu_{xy} } -x(2+z)u_{xz} -y(z+1)u_{xz}\\& \,\,\,\,\,\,\,\,\,{-c xu_{x} -b yu_{y} +\left[1-a +\left(b +c+1\right)z\right]u_{z} +b c u=0,} \end{aligned}\right.
$

where $u\equiv \,\,{\rm{E}}_{386}\left(a, b, c; x,y,z\right)$.

\bigskip

\begin{equation}
{\rm{E}}_{387}\left(a, b, c; x,y,z\right)=\sum\limits_{m,n,p=0}^\infty{(a)_{n}(b)_{p-m-n}(c)_{2m+n-p}}\frac{x^m}{m!}\frac{y^n}{n!}\frac{z^p}{p!},
\end{equation}

region of convergence:
$$ \left\{ \left\{r<\frac{1}{4} \wedge s<\frac{1}{2}+\frac{1}{2}\sqrt{1-4r}\right\}\cup\left\{s\leq\frac{1}{2}\right\},\,\,\,t<\infty
\right\}.
$$

System of partial differential equations:

$
  \left\{\begin{aligned}
 &{x(1+4x)u_{xx} +y^{2} u_{yy} +z^{2} u_{zz} +(1+4x)yu_{xy} } -(1+4x)zu_{xz}-2yzu_{yz}\\& \,\,\,\,\,\,\,\,\,  {+\left[1-b +\left(4c+6\right)x\right]u_{x} +2\left(1+c\right)yu_{y} -2c zu_{z} +c(c +1)u=0,} \\
 & y(1+y)u_{yy} +x(1+2y)u_{xy} -(1+y)zu_{yz} \\& \,\,\,\,\,\,\,\,\, {+2axu_{x} +\left[1-b+\left(a+c +1\right)y\right]u_{y} -azu_{z} +ac u=0,} \\
 & {zu_{zz} -2xu_{xz} -yu_{yz} -xu_{x} -yu_{y} +(1-c +z)u_{z} +bu=0,} \end{aligned}\right.
 $

where $u\equiv \,\,{\rm{E}}_{387}\left(a, b, c; x,y,z\right)$.

\bigskip

\begin{equation}
{\rm{E}}_{388}\left(a, b, c; x,y,z\right)=\sum\limits_{m,n,p=0}^\infty{(a)_{p}(b)_{p-m-n}(c)_{2m+n-p}}\frac{x^m}{m!}\frac{y^n}{n!}\frac{z^p}{p!},
\end{equation}

region of convergence:
$$ \left\{ t(1+2\sqrt{r})<1,\,\,\,s<\infty
\right\}.
$$

System of partial differential equations:

$
 \left\{\begin{aligned}
 &{x(1+4x)u_{xx} +y^{2} u_{yy} +z^{2} u_{zz} +(1+4x)yu_{xy}   } \\
 &\,\,\,\,\,\,\,\,\,-(1+4x)zu_{xz}-2yzu_{yz} {+\left[1-b +\left(4c +6\right)x\right]u_{x} +2\left(1+c \right)yu_{y} -2c zu_{z} +c(c +1)u=0,} \\
 & {yu_{yy} +xu_{xy} -zu_{yz} +2xu_{x} +(1-b +y)u_{y} -zu_{z} + c u=0,} \\
 & {z(1+z)u_{zz} -x(2+z)u_{xz} -y(z+1)u_{xz} -axu_{x}}{-ayu_{y} +\left[1- c+\left(a+b+1\right)z\right]u_{z} +ab u=0,} \end{aligned}\right.
  $

where $u\equiv \,\,{\rm{E}}_{388}\left(a, b, c; x,y,z\right)$.

\bigskip

\begin{equation}
{\rm{E}}_{389}\left(a, b; x,y,z\right)=\sum\limits_{m,n,p=0}^\infty{(a)_{p-m-n}(b)_{2m+n-p}}\frac{x^m}{m!}\frac{y^n}{n!}\frac{z^p}{p!},
\end{equation}

region of convergence:
$$ \left\{ r<\frac{1}{4},\,\,\,s<\infty,\,\,\,t<\infty
\right\}.
$$

System of partial differential equations:

$ \left\{\begin{aligned}
&{x(1+4x)u_{xx} +y^{2} u_{yy} +z^{2} u_{zz} +(1+4x)yu_{xy}  }  -(1+4x)zu_{xz} -2yzu_{yz}\\& \,\,\,\,\,\,\,\,\,{+\left[1-a +\left(4b+6\right)x\right]u_{x} +2\left(1+b \right)yu_{y} -2bzu_{z} +b(b+1)u=0,} \\
& {yu_{yy} +xu_{xy} -zu_{yz} +2xu_{x} +(1-a +y)u_{y} -zu_{z} +b u=0,} \\
& {zu_{zz} -2xu_{xz} -yu_{yz} -xu_{x} -yu_{y} +(1-b +z)u_{z} +a u=0,} \end{aligned}\right.
$

where $u\equiv \,\,{\rm{E}}_{389}\left(a, b; x,y,z\right)$.

\bigskip

\begin{equation}
{\rm{E}}_{390}\left(a, b, c; x,y,z\right)=\sum\limits_{m,n,p=0}^\infty{(a)_{2m+n}(b)_{p-m-n}(c)_{n-p}}\frac{x^m}{m!}\frac{y^n}{n!}\frac{z^p}{p!},
\end{equation}

region of convergence:
$$ \left\{ \left\{r<\frac{1}{4} \wedge s<\frac{1}{2}+\frac{1}{2}\sqrt{1-4r}\right\}\cup\left\{s\leq\frac{1}{2}\right\},\,\,\,t<\infty
\right\}.
$$

System of partial differential equations:

$
 \left\{\begin{aligned}
&{x(1+4x)u_{xx} +y^{2} u_{yy} +(1+4x)yu_{xy} -zu_{xz} }\\& \,\,\,\,\,\,\,\,\,{+\left[1-b +(4a+6)x \right]u_{x} +2(a+1)yu_{y} +a(a+1)u=0,} \\
& {y(1+y)u_{yy} +x(1+2y)u_{xy}  }-(1+y)zu_{yz} -2xzu_{xz}\\& \,\,\,\,\,\,\,\,\,{+2c xu_{x} +\left[1-b +\left(a+c +1\right)y\right]u_{y} -azu_{z} +ac u=0,} \\
& {zu_{zz} -yu_{yz} -xu_{x} -yu_{y} +(1-c +z)u_{z} +b u=0,} \end{aligned}\right.
$

where $u\equiv \,\,{\rm{E}}_{390}\left(a, b, c; x,y,z\right)$.

\bigskip

\begin{equation}
{\rm{E}}_{391}\left(a, b, c; x,y,z\right)=\sum\limits_{m,n,p=0}^\infty{(a)_{2m+n}(b)_{n-m-p}(c)_{p-n}}\frac{x^m}{m!}\frac{y^n}{n!}\frac{z^p}{p!},
\end{equation}

region of convergence:
$$ \left\{ r<\frac{1}{4}\wedge s<\min\left\{\Psi_1(r), \Psi_2(r)\right\}=s<1 \wedge r<\min\left\{\Theta_1(s), \Theta_2(s)\right\},\,\,\,t<\infty
\right\}.
$$

System of partial differential equations:

$
\left\{\begin{aligned}
 &{x(1+4x)u_{xx} +y^{2} u_{yy}  } -(1-4x)yu_{xy} +zu_{xz}\\& \,\,\,\,\,\,\,\,\, {+\left[1-b +(4a+6)x \right]u_{x} +2(1+a)yu_{y} +a(a+1)u=0,} \\
 & {y(1+y)u_{yy} -2x^{2} u_{xx} +xyu_{xy}  } -2xzu_{xz}-(1+y)zu_{yz} \\& \,\,\,\,\,\,\,\,\, {-\left(a-2b +2\right)xu_{x} +\left[1-c +\left(a+b +1\right)y\right]u_{y} -azu_{z} +ab u=0,} \\
 & {zu_{zz} +xu_{xz} -yu_{yz} -yu_{y} +(1-b +z)u_{z} +c u=0,} \end{aligned}\right.
 $

where $u\equiv \,\,{\rm{E}}_{391}\left(a, b, c; x,y,z\right)$.

\bigskip

\begin{equation}
{\rm{E}}_{392}\left(a, b, c; x,y,z\right)=\sum\limits_{m,n,p=0}^\infty{(a)_{2n-m}(b)_{2p-n}(c)_{m-p}}\frac{x^m}{m!}\frac{y^n}{n!}\frac{z^p}{p!},
\end{equation}

region of convergence:
$$ \left\{ r<\infty,\,\,\,4s(1+2\sqrt{t})<1
\right\}.
$$

System of partial differential equations:

$
\left\{\begin{aligned}
& {xu_{xx} -2yu_{xy} +\left(1-a +x\right)u_{x} -zu_{z} +c u=0,} \\
& {y(1+4y)u_{yy} +x^{2} u_{xx} -4xyu_{xy} -2zu_{yz}} {-2a xu_{x} +\left[1-b +\left(4a +6\right)y\right]u_{y} +a \left(a +1\right)u=0,} \\
& {z(1+4z)u_{zz} +y^{2} u_{yy} -xu_{xz} -4yzu_{yz}} {-2b yu_{y} +\left[1-c +\left(4b +6\right)z\right]u_{z} +b (b +1)u=0,} \end{aligned}\right.
$

where $u\equiv \,\,{\rm{E}}_{392}\left(a, b, c; x,y,z\right)$.

\bigskip

\begin{equation}
{\rm{E}}_{393}\left(a, b, c; x,y,z\right)=\sum\limits_{m,n,p=0}^\infty{(a)_{m}(b)_{2n-p}(c)_{2p-m-n}}\frac{x^m}{m!}\frac{y^n}{n!}\frac{z^p}{p!},
\end{equation}

region of convergence:
$$ \left\{ r<\infty,\,\,\,Z_1\cap Z_2
\right\},
$$

$$
 Z_1=\left\{s<\Phi_1(t)\right\}=\left\{s<\frac{1}{4}\wedge t<\Phi_2(s)\right\},
$$

$$
Z_2=\left\{t<\Phi_1(s)\right\}=\left\{t<\frac{1}{4}\wedge s<\Phi_2(t)\right\}.
$$

System of partial differential equations:

$
\left\{\begin{aligned}
 &{xu_{xx} +yu_{xy} -2zu_{xz} +\left(1-c +x\right)u_{x} +au=0;} \\
 & {y(1+4y)u_{yy} +z^{2} u_{zz} +xu_{xy} -2(1+2y)zu_{yz} }\\& \,\,\,\,\,\,\,\,\, {+\left[1-c +\left(4b +6\right)y\right]u_{y} -2b zu_{z} +b \left(b +1\right)u=0,} \\
 & {z(1+4z)u_{zz}+x^{2} u_{xx} +y^{2} u_{yy}  +2xyu_{xy} } \\
&\,\,\,\,\,\,\,\,\,-4xzu_{xz} -2y(1+2z)u_{xz} {-2c xu_{x} -2c yu_{y} +\left[1-b+\left(4c +6\right)z\right]u_{z} +c(1+c )u=0,} \end{aligned}\right.
 $

where $u\equiv \,\,{\rm{E}}_{393}\left(a, b, c; x,y,z\right)$.

\bigskip

\begin{equation}
{\rm{E}}_{394}\left( a, b;  x,y,z\right)=\sum\limits_{m,n,p=0}^\infty{(a)_{2n-p}(b)_{2p-m-n}}\frac{x^m}{m!}\frac{y^n}{n!}\frac{z^p}{p!},
\end{equation}

region of convergence:
$$ \left\{ r<\infty,\,\,\,Z_1\cap Z_2
\right\},
$$

$$
Z_1=\left\{s<\Phi_1(t)\right\}=\left\{s<\frac{1}{4}\wedge t<\Phi_2(s)\right\},
$$

$$
Z_2=\left\{t<\Phi_1(s)\right\}=\left\{t<\frac{1}{4}\wedge s<\Phi_2(t)\right\}.
$$

System of partial differential equations:

$
 \left\{\begin{aligned}
&{xu_{xx} +yu_{xy} -2zu_{xz} +(1-b )u_{x} +u=0,} \\
& {y(1+4y)u_{yy} +z^{2} u_{zz} +xu_{xy} -2(1+2y)zu_{yz} } \\& \,\,\,\,\,\,\,\,\,{+\left[1-b +\left(4a +6\right)y\right]u_{y} -2a zu_{z} +a \left(a+1\right)u=0,} \\
& {z(1+4z)u_{zz} + x^{2} u_{xx} +y^{2} u_{yy} +2xyu_{xy}   } \\
&\,\,\,\,\,\,\,\,\,\,\, -4xzu_{xz}-2y(1+2z)u_{xz}{-2b xu_{x} -2b yu_{y} +\left[1-a +\left(4b +6\right)z\right]u_{z} +b (1+b)u=0,} \end{aligned}\right. $

where $u\equiv \,\,{\rm{E}}_{394}\left( a, b;  x,y,z\right)$.

\bigskip

\begin{equation} \label{e395}
{\rm{E}}_{395}\left(a, b; x,y,z\right)=\sum\limits_{m,n,p=0}^\infty{(a)_{2m+n-p}(b)_{2p-m-n}}\frac{x^m}{m!}\frac{y^n}{n!}\frac{z^p}{p!},
\end{equation}

region of convergence:
$$ \left\{ Z_1\cap Z_2,\,\,\,s<\infty
\right\},
$$

$$  Z_1=\left\{r<\Phi_1(t)\right\}=\left\{r<\frac{1}{4}\wedge t<\Phi_2(r)\right\}
,
$$

$$  Z_2=\left\{t<\Phi_1(r)\right\}=\left\{t<\frac{1}{4}\wedge r<\Phi_2(t)\right\}
.
$$

System of partial differential equations:

$
 \left\{\begin{aligned}
 &{x(1+4x)u_{xx} +y^{2} u_{yy} +z^{2} u_{zz} +(1+4x)yu_{xy} } -2(1+2x)zu_{xz} -2yzu_{yz} \\& \,\,\,\,\,\,\,\,\,{+\left[1-b +\left(4a +6\right)x\right]u_{x} +2\left(1+a \right)yu_{y} -2a zu_{z} +a \left(a+1\right)u=0,} \\
 & {yu_{yy} +xu_{xy} -2zu_{yz} +2xu_{x} +(1-b +y)u_{y} -zu_{z} +a u=0,} \\
 & { z(1+4z)u_{zz}+ x^{2} u_{xx} +y^{2} u_{yy} +2xyu_{xy} } -2x(1+2z)u_{xz} -y(1+4z)u_{xz}\\& \,\,\,\,\,\,\,\,\, {-2b xu_{x} -2b yu_{y} +\left[1-a +\left(4b +6\right)z\right]u_{z} +b (1+b )u=0,} \end{aligned}\right.
 $

where $u\equiv \,\,{\rm{E}}_{395}\left(a, b; x,y,z\right)$.

\section{Example of finding particular solutions near the origin}  \label{S9}

As an example, we consider  the system
\begin{equation}  \label{e10as}
\left\{
{\begin{array}{*{20}{l}}
  {x\left( {1 - x} \right){u_{xx}} - xy{u_{xy}} + \left[ {{c_1} - \left( {{a_1} + {a_3} + 1} \right)x} \right]{u_x} - {a_3}y{u_y} - {a_1}{a_3}u = 0,} \\
  y\left( {1 - y} \right){u_{yy}} - xy{u_{xy}} - xz{u_{xz}} - yz{u_{yz}} \\\,\,\,\,\,\,\,\,\,+ \left[ {{c_2} - \left( {{a_1} + {a_2} + 1} \right)y} \right]{u_y}
  - {a_2}x{u_x} - {a_1}z{u_z} - {a_1}{a_2}u = 0,  \\
  {z\left( {1 - z} \right){u_{zz}} - yz{u_{yz}} + \left[ {{c_3} - \left( {{a_2} + {a_4} + 1} \right)z} \right]{u_z} - {a_4}y{u_y} - {a_2}{a_4}u = 0,}
\end{array}} \right.
\end{equation}
associated with the function

\begin{equation}  \label{e10af}
{F_{10a}}\left( {{a_1},{a_2},{a_3},{a_4};{c_1},{c_2},{c_3};x,y,z} \right) \hfill \\
  = \sum\limits_{m,n,p = 0}^\infty  {} \frac{{{{\left( {{a_1}} \right)}_{m + n}}{{\left( {{a_2}} \right)}_{n + p}}{{\left( {{a_3}} \right)}_m}{{\left( {{a_4}} \right)}_p}}}{{{{\left( {{c_1}} \right)}_m}{{\left( {{c_2}} \right)}_n}{{\left( {{c_3}} \right)}_p}}}\frac{x^m}{m!}\frac{y^n}{n!}\frac{z^p}{p!}.
\end{equation}

We look for particular solutions of the system  (700) near the origin in the form

\begin{equation} \label{eq826}
u\left( {x,y,z} \right) = {x^\tau }{y^\nu }{z^\lambda }w\left( {x,y,z} \right).
\end{equation}

We calculate the derivatives

\begin{equation} \label{eq 827}
\begin{gathered}
  {u_x} = \tau {x^{\tau  - 1}}{y^\nu }{z^\lambda }w + {x^\tau }{y^\nu }{z^\lambda }{w_x}, \hfill \\
  {u_y} = \nu {x^\tau }{y^{\nu  - 1}}{z^\lambda }w + {x^\tau }{y^\nu }{z^\lambda }{w_y}, \hfill \\
  {u_z} = \lambda {x^\tau }{y^\nu }{z^{\lambda  - 1}}w + {x^\tau }{y^\nu }{z^\lambda }{w_z}, \hfill \\
\end{gathered}
\end{equation}

\begin{equation} \label{eq 828}
\begin{gathered}
  {u_{xy}} = \tau \nu {x^{\tau  - 1}}{y^{\nu  - 1}}{z^\lambda }w + \tau {x^{\tau  - 1}}{y^\nu }{z^\lambda }{w_y} + \nu {x^\tau }{y^{\nu  - 1}}{z^\lambda }{w_x} + {x^\tau }{y^\nu }{z^\lambda }{w_{xy}}, \hfill \\
  {u_{xz}} = \tau \lambda {x^{\tau  - 1}}{y^\nu }{z^{\lambda  - 1}}w + \tau {x^{\tau  - 1}}{y^\nu }{z^\lambda }{w_z} + \lambda {x^\tau }{y^\nu }{z^{\lambda  - 1}}{w_x} + {x^\tau }{y^\nu }{z^\lambda }{w_{xz}}, \hfill \\
  {u_{yz}} = \nu \lambda {x^\tau }{y^{\nu  - 1}}{z^{\lambda  - 1}}w + \nu {x^\tau }{y^{\nu  - 1}}{z^\lambda }{w_z} + \lambda {x^\tau }{y^\nu }{z^{\lambda  - 1}}{w_y} + {x^\tau }{y^\nu }{z^\lambda }{w_{yz}}, \hfill \\
\end{gathered}
\end{equation}

\begin{equation} \label{eq 829}
\begin{gathered}
  {u_{xx}} = \tau \left( {\tau  - 1} \right){x^{\tau  - 2}}{y^\nu }{z^\lambda }w + 2\tau {x^{\tau  - 1}}{y^\nu }{z^\lambda }{w_x} + {x^\tau }{y^\nu }{z^\lambda }{w_{xx}}, \hfill \\
  {u_{yy}} = \nu \left( {\nu  - 1} \right){x^\tau }{y^{\nu  - 2}}{z^\lambda }w + 2\nu {x^\tau }{y^{\nu  - 1}}{z^\lambda }{w_y} + {x^\tau }{y^\nu }{z^\lambda }{w_{yy}}, \hfill \\
  {u_{zz}} = \lambda \left( {\lambda  - 1} \right){x^\tau }{y^\nu }{z^{\lambda  - 2}}w + 2\lambda {x^\tau }{y^\nu }{z^{\lambda  - 1}}{w_z} + {x^\tau }{y^\nu }{z^\lambda }{w_{zz}}. \hfill \\
\end{gathered}
\end{equation}

Substituting (703) -- (705) into system (700) and dividing by ${x^\tau} {y^\nu} {z^\lambda}$, we have

\begin{equation} \label{eq 830}
\left\{ {\begin{array}{*{20}{l}}
  \begin{gathered}
  x\left( {1 - x} \right){w_{xx}} - xy{w_{xy}} + \left\{ {2\tau  + {c_1} - \left[ {\left( {\tau  + \nu  + {a_1}} \right) + \left( {\tau  + {a_3}} \right) + 1} \right]x} \right\}{w_x} \hfill \\
 \,\,\,\,\,\,\,\,\,\,\,\,\,\, - \left( {\tau  + {a_3}} \right)y{w_y}   - \left[ { - \tau \left( {\tau  - 1 + {c_1}} \right){x^{ - 1}} + \left( {\tau  + \nu  + {a_1}} \right)\left( {\tau  + {a_3}} \right)} \right]w = 0, \hfill \\
\end{gathered}  \\
  \begin{gathered}
  y\left( {1 - y} \right){w_{yy}} - xy{w_{xy}} - xz{w_{xz}} - yz{w_{yz}} - \left( {\nu  + \lambda  + {a_2}} \right)x{w_x} \hfill \\
  \,\,\,\,\,\,\,\,\,\,\,\,\,\, + \left\{ {2\nu  + {c_2} - \left[ {\left( {\tau  + \nu  + {a_1}} \right) + \left( {\nu  + \lambda  + {a_2}} \right) + 1} \right]y} \right\}{w_y} - \left( {\tau  + \nu  + {a_1}} \right)z{w_z} \hfill \\
  \,\,\,\,\,\,\,\,\,\,\,\,\,\, - \left[ { - \nu \left( {\nu  - 1 + {c_2}} \right){y^{ - 1}} + \left( {\tau  + \nu  + {a_1}} \right)\left( {\nu  + \lambda  + {a_2}} \right)} \right]w = 0, \hfill \\
\end{gathered}  \\
  \begin{gathered}
  z\left( {1 - z} \right){w_{zz}} - yz{w_{yz}} - \left( {\lambda  + {a_4}} \right)y{w_y} + \left\{ {2\lambda  + {c_3} - \left[ {\left( {\nu  + \lambda  + {a_2}} \right) + \left( {\lambda  + {a_4}} \right) + 1} \right]z} \right\}{w_z} \hfill \\
  \,\,\,\,\,\,\,\,\,\,\,\,\,\, + \left[ {\lambda \left( {\lambda  - 1 + {c_3}} \right){z^{ - 1}} - \left( {\nu  + \lambda  + {a_2}} \right)\left( {\lambda  + {a_4}} \right)} \right]w = 0. \hfill \\
\end{gathered}
\end{array}} \right.
\end{equation}

If the following algebraic system of equations
\begin{equation} \label{eq 831}
\left\{ {\begin{array}{*{20}{l}}
  {\tau \left( {\tau  - 1 + {c_1}} \right) = 0,} \\
  {\nu \left( {\nu  - 1 + {c_2}} \right) = 0,} \\
  {\lambda \left( {\lambda  - 1 + {c_3}} \right) = 0,}
\end{array}} \right.
\end{equation}
has a solution, then the system (706) is simplified:

\begin{equation} \label{eq 8300000}
\left\{ {\begin{array}{*{20}{l}}
  \begin{gathered}
  x\left( {1 - x} \right){w_{xx}} - xy{w_{xy}} + \left\{ {2\tau  + {c_1} - \left[ {\left( {\tau  + \nu  + {a_1}} \right) + \left( {\tau  + {a_3}} \right) + 1} \right]x} \right\}{w_x} \hfill \\
 \,\,\,\,\,\,\,\,\,\,\,\,\,\, - \left( {\tau  + {a_3}} \right)y{w_y}   - \left( {\tau  + \nu  + {a_1}} \right)\left( {\tau  + {a_3}} \right)w = 0, \hfill \\
\end{gathered}  \\
  \begin{gathered}
  y\left( {1 - y} \right){w_{yy}} - xy{w_{xy}} - xz{w_{xz}} - yz{w_{yz}} - \left( {\nu  + \lambda  + {a_2}} \right)x{w_x} \hfill \\
  \,\,\,\,\,\,\,\,\,\,\,\,\,\, + \left\{ {2\nu  + {c_2} - \left[ {\left( {\tau  + \nu  + {a_1}} \right) + \left( {\nu  + \lambda  + {a_2}} \right) + 1} \right]y} \right\}{w_y} - \left( {\tau  + \nu  + {a_1}} \right)z{w_z} \hfill \\
  \,\,\,\,\,\,\,\,\,\,\,\,\,\, -  \left( {\tau  + \nu  + {a_1}} \right)\left( {\nu  + \lambda  + {a_2}} \right)w = 0, \hfill \\
\end{gathered}  \\
  \begin{gathered}
  z\left( {1 - z} \right){w_{zz}} - yz{w_{yz}} - \left( {\lambda  + {a_4}} \right)y{w_y} + \left\{ {2\lambda  + {c_3} - \left[ {\left( {\nu  + \lambda  + {a_2}} \right) + \left( {\lambda  + {a_4}} \right) + 1} \right]z} \right\}{w_z} \hfill \\
  \,\,\,\,\,\,\,\,\,\,\,\,\,\, - \left( {\nu  + \lambda  + {a_2}} \right)\left( {\lambda  + {a_4}} \right) w = 0. \hfill \\
\end{gathered}
\end{array}} \right.
\end{equation}

It is obvious that the system (707) has 8 solutions:

\[
\left. 1 \right) \,\,\,\,\tau_1=0, \,\,\,\,\, \nu_1=0, \,\,\,\,\,\lambda_1=0;
\]

\[
\left. 2 \right) \,\,\,\,\tau_2=1-c_1, \,\,\,\,\, \nu_2=0, \,\,\,\,\,\lambda_2=0;
\]

\[
\left. 3 \right) \,\,\,\,\tau_3=0, \,\,\,\,\, \nu_3=1-c_2, \,\,\,\,\,\lambda_3=0;
\]

\[
\left. 4 \right) \,\,\,\,\tau_4=0, \,\,\,\,\, \nu_4=0, \,\,\,\,\,\lambda_4=1-c_3;
\]

\[
\left. 5 \right) \,\,\,\,\tau_5=1-c_1, \,\,\,\,\, \nu_5=1-c_2, \,\,\,\,\,\lambda_5=0;
\]

\[
\left. 6 \right) \,\,\,\,\tau_6=0, \,\,\,\,\, \nu_6=1-c_2, \,\,\,\,\,\lambda_6=1-c_3;
\]

\[
\left. 7 \right) \,\,\,\,\tau_7=1-c_1, \,\,\,\,\, \nu_7=0, \,\,\,\,\,\lambda_7=1-c_3;
\]

\[
\left. 8 \right)\,  \tau_8=1-c_1, \nu_8=1-c_2, \lambda_8=1-c_3.
\]

Now  substituting these solutions into the system (708), by virtue of (702), we obtain 8 particular solutions of the system  (700) near the origin

$
{u_1} = {F_{10a}}\left( {{a_1},{a_2},{a_3},{a_4};{c_1},{c_2},{c_3};x,y,z} \right),
$

$
{u_2} = {x^{1 - {c_1}}}{F_{10a}}\left( {1 - {c_1} + {a_1},{a_2},1 - {c_1} + {a_3},{a_4};2 - {c_1},{c_2},{c_3};x,y,z} \right),
$

$
{u_3} = {y^{1 - {c_2}}}{F_{10a}}\left( {1 - {c_2} + {a_1},1 - {c_2} + {a_2},{a_3},{a_4};2 - {c_1},{c_2},{c_3};x,y,z} \right),
$

$
{u_4} = {z^{1 - {c_3}}}{F_{10a}}\left( {{a_1},1 - {c_3} + {a_2},{a_3},1 - {c_3} + {a_4};{c_1},{c_2},2 - {c_3};x,y,z} \right),
$

$
{u_5} = {x^{1 - {c_1}}}{y^{1 - {c_2}}}{F_{10a}}\left( {2 - {c_1} - {c_2} + {a_1},1 - {c_2} + {a_2},1 - {c_1} + {a_3},{a_4};2 - {c_1},2 - {c_2},{c_3};x,y,z} \right),
$

$
{u_6} = {y^{1 - {c_2}}}{z^{1 - {c_3}}}
{F_{10a}}\left( {1 - {c_2} + {a_1},2 - {c_2} - {c_3} + {a_2},{a_3},1 - {c_3} + {a_4};{c_1},2 - {c_2},2 - {c_3};x,y,z} \right),
$

$
{u_7} = {x^{1 - {c_1}}}{z^{1 - {c_3}}}{F_{10a}}\left( {1 - {c_1} + {a_1},1 - {c_3} + {a_2},1 - {c_1} + {a_3},1 - {c_3} + {a_4};2 - {c_1},{c_2},2 - {c_3};x,y,z} \right),
$

$
\begin{aligned}
  {u_8}& = {x^{1 - {c_1}}}{y^{1 - {c_2}}}{z^{1 - {c_3}}} \times\\& \times {F_{10a}}\left( {2 - {c_1} - {c_2} + {a_1},2 - {c_2} - {c_3} + {a_2},1 - {c_1} + {a_3},1 - {c_3} + {a_4};}
 {2 - {c_1},2 - {c_2},2 - {c_3};x,y,z} \right).
 \end{aligned}
$

\bigskip

\bigskip

\textbf{Funding:} The research   has been initiated by the Intercontinental Research Center "Analysis and PDE" (Ghent University, Belgium), supported by the FWO Odysseus 1 grant G.0H94.18N: Analysis and partial differential equations and by the Methusalem programme of the Ghent University Special Research Fund (BOF) (Grant number 01M01021).

\end{document}